\documentclass[12pt]{article}
\usepackage[dvips]{graphicx}
\font\teneufm=eufm10
\font\seveneufm=eufm7
\font\fiveeufm=eufm5
\newfam\eufmfam
\textfont\eufmfam=\teneufm
\scriptfont\eufmfam=\seveneufm
\scriptscriptfont\eufmfam=\fiveeufm
\def\eufm#1{{\fam\eufmfam\relax#1}}

\newcommand\beq[1]{ \begin{equation}\label{#1} }
\newcommand{\eeq}{ \end{equation} }

\newcommand\beqa[1]{ \begin{eqnarray} \label{#1}}
\newcommand{\eeqa}{ \end{eqnarray} }
\newcommand{\beqano}{ \begin{eqnarray*} }
\newcommand{\eeqano}{ \end{eqnarray*} }
\newcommand\arr[1]{\left\{\begin{array}{l}#1\end{array}\right.}
\renewcommand{\theequation}{\arabic{section}.\arabic{equation}}

\newtheorem{theorem}{Theorem}[section]
\newtheorem{definition}{Definition}[section]
\newtheorem{proposition}{Proposition}[section]
\newtheorem{lemma}{Lemma}[section]
\newtheorem{sublemma}{Sublemma}[section]
\newtheorem{remark}{Remark}[section]
\newtheorem{notationalremark}{Notational Remark}[section]
\newtheorem{corollary}{Corollary}[section]
\newtheorem{assumption}{Assumption}[section]
\newtheorem{claim}{Claim}[section]
\newtheorem{tools}{$\negsp\negsp$}[subsection]

\newcommand\thm[1]{ \begin{theorem}\label{#1}}
\newcommand\thmtwo[2]{ \begin{theorem}[#1]\label{#2}}
\newcommand\ethm{ \end{theorem} }
\newcommand\dfn[1]{ \begin{definition}\label{#1} \rm}
\newcommand\dfntwo[2]{ \begin{definition}[#1]\label{#2} \rm}
\newcommand\edfn{ \end{definition} }
\newcommand\pro[1]{ \begin{proposition}\label{#1}}
\newcommand\protwo[2]{ \begin{proposition}[#1]\label{#2}}
\newcommand\epro{ \end{proposition} }
\newcommand\lem[1]{ \begin{lemma}\label{#1}}
\newcommand\lemtwo[2]{ \begin{lemma}[#1]\label{#2}}
\newcommand\elem{ \end{lemma} }
\newcommand\sublem[1]{ \begin{sublemma}\label{#1}}
\newcommand\sublemtwo[2]{ \begin{sublemma}[#1]\label{#2}}
\newcommand\esublem{ \end{sublemma} }
\newcommand\rem[1]{ \begin{remark}\label{#1} \rm}
\newcommand\erem{ \end{remark} }
\newcommand\notrem[1]{ \begin{notationalremark}\label{#1} \rm}
\newcommand\enotrem{ \end{notationalremark} }
\newcommand\cor[1]{ \begin{corollary}\label{#1}}
\newcommand\cortwo[2]{ \begin{corollary}[#1]\label{#2}}
\newcommand\ecor{ \end{corollary} }
\newcommand\asmp[1]{ \begin{assumption}\label{#1}}
\newcommand\asmptwo[2]{ \begin{assumption}[#1]\label{#2}}
\newcommand\easmp{ \end{assumption} }
\newcommand\clm[1]{ \begin{claim}\label{#1}}
\newcommand\eclm{ \end{claim} }

\expandafter\chardef\csname pre amssym.def
at\endcsname=\the\catcode`\@
\catcode`\@=11
\def\undefine#1{\let#1\undefined}
\def\newsymbol#1#2#3#4#5{\let\next@\relax
 \ifnum#2=\@ne\let\next@\msafam@\else
 \ifnum#2=\tw@\let\next@\msbfam@\fi\fi
 \mathchardef#1="#3\next@#4#5}
\def\mathhexbox@#1#2#3{\relax
 \ifmmode\mathpalette{}{\m@th\mathchar"#1#2#3}%
 \else\leavevmode\hbox{$\m@th\mathchar"#1#2#3$}\fi}
\def\hexnumber@#1{\ifcase#1 0\or 1\or 2\or 3\or 4\or 5\or 6\or 7\or
8\or
 9\or A\or B\or C\or D\or E\or F\fi}
\ifcase\@ptsize
 \font\tenmsb=msbm10
 \font\sevenmsb=msbm7
 \font\fivemsb=msbm5
\or
 \font\tenmsb=msbm10 scaled \magstephalf
 \font\sevenmsb=msbm7 scaled \magstephalf
 \font\fivemsb=msbm5  scaled \magstephalf
\or
 \font\tenmsb=msbm10 scaled \magstep1
 \font\sevenmsb=msbm7 scaled \magstep1
 \font\fivemsb=msbm5 scaled \magstep1
\fi
\newfam\msbfam
\textfont\msbfam=\tenmsb
\scriptfont\msbfam=\sevenmsb
\scriptscriptfont\msbfam=\fivemsb
\edef\msbfam@{\hexnumber@\msbfam}
\def\Bbb#1{\fam\msbfam\relax#1}
\def\widehat#1{\setboxz@h{$\m@th#1$}%
 \ifdim\wdz@>\tw@ em\mathaccent"0\msbfam@5B{#1}%
 \else\mathaccent"0362{#1}\fi}
\def\widetilde#1{\setboxz@h{$\m@th#1$}%
 \ifdim\wdz@>\tw@ em\mathaccent"0\msbfam@5D{#1}%
 \else\mathaccent"0365{#1}\fi}

\def\RIfM@{\relax\ifmmode}
\def\nonmatherr@#1{\errmessage{\string#1\space allowed only in math mode}}
\def\Bbb{\RIfM@\expandafter\Bbb@\else
 \expandafter\nonmatherr@\expandafter\Bbb\fi}
\def\Bbb@#1{{\Bbb@@{#1}}}
\def\Bbb@@#1{\fam\msbfam\relax#1}
\def\setboxz@h{\setbox\z@\hbox}
\def\wdz@{\wd\z@}
\catcode`\@=\csname pre amssym.def at\endcsname
\newcommand{\ie}{{\it i.e.\  }}

\newcommand{\nl}{{\smallskip\noindent}}
\newcommand{\noi}{{\noindent}}
\newcommand{\pagina}{{\vfill\eject}}
\newcommand{\negsp}{\hspace{-.09truecm}}  

\newcommand{\dst}{\displaystyle}
\newcommand\ovl[1]{ \overline {#1} }

\newcommand{\torus}{ {\Bbb T}   }
\renewcommand{\natural}{ {\Bbb N}   }
\newcommand{\real}{ {\Bbb R}   }
\newcommand{\integer}{ {\Bbb Z}   }
\newcommand{\complex}{ {\Bbb C}   }
\newcommand{\rational}{ {\Bbb Q}  }

\renewcommand{\a }{ {\alpha}   }
\renewcommand{\b}{ {\beta}   }
\newcommand{\g}{ {\gamma}   }
\newcommand{\G}{ {\Gamma}   }
\renewcommand{\d}{ {\delta}   }
\newcommand{\D}{ {\Delta}   }
\newcommand{\e }{ {\varepsilon}   }

\renewcommand{\k}{ {\kappa}   }
\renewcommand{\l}{ {\lambda}   }
\renewcommand{\L}{ {\Lambda}   }
\newcommand{\m}{ {\mu}   }
\newcommand{\n}{ {\nu}   }

\newcommand{\p}{ {\pi}   }
\renewcommand{\P}{ {\Pi}   }
\renewcommand{\r}{ {\rho}   }
\newcommand{\s}{ {\sigma}   }

\renewcommand{\t}{ {\tau}   }

\renewcommand{\o}{ {\omega}   }
\renewcommand{\O}{ {\Omega}   }

\newcommand{\const}{{\, \rm const\, }}
\renewcommand{\Im}{{\, \rm Im\, }}
\renewcommand{\Re}{{\, \rm Re\, }}

\newcommand{\cA}{ {\cal A} }
\newcommand{\cB}{ {\cal B} }
\newcommand{\cE}{ {\cal E} }
\newcommand{\cT}{ {\cal T} }
\newcommand{\cR}{ {\cal R} }
\newcommand{\cH}{ {\cal H} }
\newcommand{\cK}{ {\cal K} }
\newcommand{\cC}{ {\cal C} }
\newcommand{\cD}{ {\cal D} }
\newcommand{\cF}{ {\cal F} }

\newcommand{\cM}{ {\cal M} }
\newcommand{\cP}{ {\cal P} }
\newcommand{\cI}{ {\cal I} }
\newcommand{\cJ}{ {\cal J} }

\newcommand\meas{{\rm meas}}

\newcommand\id{{\, \rm id \,}}
\begin{document}

\title{\bf On the Kolmogorov Set for Many--Body Problems\thanks{{\bf Acknowledgments.} I am deeply grateful to Prof. Luigi Chierchia for introducing me to Celestial Mechanics and teaching me perturbative KAM techniques,  suggesting the subject of this thesis,
inspiring  and next following step--by--step all the details  with care and { much} ({\sl much}, {\sl much} $\cdots$) patience. Without his precious help and encouragement, I could  never have written it. I am also indebted with Luca Biasco, who taught me  techniques for the construction and the proof of the convergence of KAM algorithms and with Jacques F\'ejoz, for  many enlightening discussions and his critical review. Finally, I wish to thank Massimiliano Berti for his interest.}
}
\author{}

\date{}

\maketitle

\vglue0.7 cm



\vglue7cm

\begin{flushleft}
Candidate\\
Dr. Gabriella Pinzari
\end{flushleft}
\begin{flushright}
Supervisor\\
Prof. Luigi Chierchia
\end{flushright}

\nl

\maketitle

\newpage

\vglue2 cm
\begin{flushright}
\large{\it a mia Madre}
\end{flushright}
\newpage

{\footnotesize
\tableofcontents}

\pagina
{\footnotesize
\section{Introduction}
The {\sl many--body problem} is the the study of the motion of $1+N$ point masses $m_0$, $\cdots$, $m_N$ interacting through gravity  only, hence, whose  coordinates $v_0$, $\cdots$, $v_N\in\real^{\rm d}$ (where ${\rm d}=2,\ 3$) obey to the Newton's equations 
 \beq{Newton equations}
\dst\ddot v_i=-\sum_{j\neq i}m_j\frac{v_i-v_j}{|v_i-v_j|^3}\quad \textrm{for}\quad 0\leq i\leq N\ .
\eeq 
As usual, the ``dot'' denotes the derivative with respect to  time  and $|\cdot|$ is the {Euclidean norm} of $\real^{\rm d}$.
In the planetary problem,  one mass, $m_0$ (the ``Sun''), is much greater than the others (the ``planets''). So, it is customary to  introduce a small parameter  $\m$ and take
\beqa{small masses}
m_0=\bar m_0\quad \textrm{and}\quad m_i=\m \bar m_i\quad \textrm{for}\quad 1\leq i\leq N\ .
\eeqa
The Hamiltonian formulation of (\ref{Newton equations}) is
\beqano
\arr{ 
\dst\dot u_i=-\partial_{v_i}\hat {\cal H}_{\rm plt}(\m;u,v)\\
\dst\dot v_i=\partial_{u_i}\hat {\cal H}_{\rm plt}(\m;u,v)\\
}
\eeqano
with
\beqa{non lin red}
\hat {\cal H}_{\rm plt}(\m;u,v):=\sum_{0\leq i\leq N}\frac{|u_i|^2}{2m_i}-\sum_{0\leq i<j\leq N}\frac{m_im_j}{|v_i-v_j|}\ .
\eeqa
(masses as in (\ref{small masses})) which is regular (real--analytic
\footnote{We recall that a  function $f:A\to \real$, where $A$ is an open, connected, bounded subset of $\real^n$, is real--analytic in $x_0\in A$ if there exist $\{a_k(x_0)\}_{k\in \natural^n}$ and an open neighborhood of $x_0$, $U(x_0)$, such that the series $\sum_{k\in \natural^n}a_k(x_0)(x-x_0)^k$ converges uniformly to $f$, for any $x\in U(x_0)$; $f$ is said to be real--analitic in $A$ if it is real--analytic in any point $x_0\in A$.}
) over the  domain (``phase space'')
$$\hat{\cal C}_{\rm cl}=\Big\{(u,v):=\Big((u_0,\cdots,u_N),(v_0,\cdots,v_N)\Big)\in \real^{\rm dN}\times \real^{\rm dN}:\quad v_i\neq v_j\quad \textrm{for}\quad 0\leq i<j\leq N\Big\}$$
Besides the energy, $\hat{\cal H}_{\rm plt}$ has, as integrals of the motion (\ie, conserved quantities along its trajectories), the three components of the {\sl linear momentum}
\beqa{lin mom}
\hat{\rm Q}=\sum_{0\leq i\leq N}u_i
\eeqa 
and, for ${\rm d}=3$, also the three components
of the {\sl angular momentum}
\beqa{ang mom intr}
\hat{\rm C}=\sum_{0\leq i\leq N} v_i\times u_i
\eeqa
(which are related to the {\sl translation} and {\sl rotation} invariance of $\hat {\cal H}_{\rm plt}$, respectively), where ``$\times$'' denotes the standard vector product of $\real^3$. Hence, the number of degrees of freedom of the system
\footnote{The number of degrees of freedom of an Hamiltonian system is defined as one half of the dimension of the phase space of the Hamiltonian. In the case of (\ref{non lin red}), the number of degrees of freedom is ${\rm d}(1+N)$}
of can be lowered.
\vskip.1in
\noi
The {\sl linear momentum reduction} is usually performed as follows. 
\vskip.1in
\noi
Consider the invariant manifold with dimension $2{\rm d}N$
\beqano
\cM_{\rm lin}:=\left\{u,\ v\in \hat\cC_{\rm cl}:\quad \sum_{0\leq i\leq N}u_i=\sum_{0\leq i\leq N}m_iv_i=0\right\}\ .
\eeqano
Let
\beqano
\cC_{\rm cl}:=\left\{(y,x)\in \real^{{\rm d}N}\times\real^{{\rm d}N}:\ 0\neq x_i\neq x_j\ \textrm{for}\ 1\leq i<j\leq N\right\}
\eeqano
be the ``collisionless phase space'' and define the embedding
\beqano
{\phi}_{\rm helio}:\ (y,x)\in\cC_{\rm cl}\subset \real^{{\rm d}N}\times \real^{{\rm d}N}&\to& (u,v)\in{\cal M}_{\rm lin}\nonumber\\
\eeqano
which acts as
\beqano\arr{\dst
u_0:=-\sum_{1\leq j\leq N}y_j\\
v_0:=-\left(\dst\sum_{0\leq j\leq N}m_j\right)^{-1}\dst\sum_{1\leq i\leq N}m_ix_i\\
u_i:=y_i\\
v_i:=x_i-\left(\dst\sum_{0\leq j\leq N}m_j\right)^{-1}\dst\sum_{1\leq j\leq N}m_jx_j
}
\eeqano
Then, it is not difficult to see that the evolution in time for the ``relative momenta--coordinates'' pairs  $(y,x)$ is governed by the Hamilton equations of
$$\widetilde {\cal H}_{\rm plt}(\m;y,x):=\hat{\cal H}_{\rm plt}\circ \phi_{\rm helio}(\m;y,x)\ .$$
A suitable rescaling of variables and Hamiltonian (which does not change the motion equations)
\beqano
{\cal H}_{\rm plt}(\m;y,x):=\m^{-1}\widetilde{\cal H}_{\rm plt}(\m;\m y,x)
\eeqano 
brings finally to 
\beqa{heliocentric hamiltonian}
{\cal H}_{\rm plt}(\m;y,x)&=&\m^{-1}\widetilde{\cal H}_{\rm plt}(\m;\m y,x)\nonumber\\
&=&\sum_{1\leq i\leq N}\left(\frac{|y_i|^2}{2\tilde m_i}-\frac{\tilde m_i\hat m_i}{|x_i|}\right)+\m\,\sum_{1\leq i<j\leq N}\left(\frac{y_i\cdot y_j}{\bar m_0}-\frac{\bar m_i\bar m_j}{|x_i-x_j|}\right)\ ,
\eeqa
where  $\hat m_i$, $\tilde m_i$ 
are the ``reduced masses''
\beqano
\hat m_i:=\bar m_0+\m\,\bar m_i\ ,\quad \tilde m_i:=\frac{\bar m_0\bar m_i}{\bar m_0+\m\,\bar m_i}
\eeqano
and $u\cdot v$
denotes the usual inner product of two vectors $u$, $v$ of $\real^{\rm d}$.

\vskip.1in
\noindent
Notice that the  angular momentum (in heliocentric variables)
\beq{angular momentum}
C:=\sum_{1\leq i\leq N}x_i\times y_i
\eeq
is still conserved along the ${\cal H}_{\rm plt}$--trajectories, which is still rotation invariant.

\vskip.1in
\noindent
When $\m=0$ (``integrable limit'' ), the ${\cal H}_{\rm plt}$--evolution is the resulting of $N$ independent Keplerian motions for the coordinates $x_1$, $\cdots$, $x_N$: each of them, accordingly to the {\sl Law of Equal Areas} 
\footnote{{\sl Law of Equal Areas}: the area spanned  by $x_i(t)$ on the ellipse  $\cE_i$ is given by $$S_i(t)=S_i(0)+\frac{a_ib_i}{2}\,n_it$$ where $n_i$, defined by $n_i^2a_i^3=\hat m_i$ is the {\sl mean motion} and $a_i$, $b_i=a_i\sqrt{1-e_i^2}$ are the semi--axes of $\cE_i$.}
draws an  ellipse in the space, whose position and shape depends only on the initial data $(\bar y_i,\bar x_i)$, all the ellipses possessing  a common focus (the Sun). The total motion is thus -- for $\m=0$ -- quasi--periodic with $N$ frequencies, provided each {\sl two -- body energy} \beq{Two Body}h_i:=\frac{|y_i|^2}{2\tilde m_i}-\frac{\tilde m_i\hat m_i}{|x_i|}\eeq is negative.

\vskip.1in
\noi
The (analytic, $C^{\infty}$) continuation, for $\m>0$, of the  quasi--periodic motions (with $N$--frequencies) of (\ref{Two Body}) with  quasi--periodic motions with more frequencies has been investigated by several authors. Of great interest is the case of ``maximal'' continuation, which consists in looking for tori with the maximum number ${\rm f}$ of frequecies possible, \ie, (analytic, $C^{\infty}$)  invariant manifolds for ${\cal H}_{\rm plt}$, diffeomorphic to the standard torus $\torus^{\rm f}$ where the angular coordinate evolves with linear low in time.
\vskip.1in
\noindent
The pioneering work on this subject is the one by Arnol'd \cite{Arn63}, who, in the framework of the KAM theory, stated the {\sl existence of a positive measure (``Cantor'') set of initial data giving rise to bounded motions}. He proved his statement  only in the case of the plane three -- body problem
(${\rm d}=N=2$) and, for the general case (spatial $(1+N)$--body problem), gave only some indications on how to extend the result.

\vskip.1in
\noi
It has been noticed in \cite{Fej04} that such indications contain a flaw. Nonetheless, Arnol'd's  proof of existence of quasi--periodic motions in the planar three body problem, is based on a refined KAM theorem -- constructed in the framework of  real--analytic functions 
and called by himself {\sl Fundamental Theorem} (quoted below), which could overcome the strong ``degeneracy'' of the problem. To explain this point, and for future use, we need a bit of preparation.
\vskip.1in
\noindent
Let us start by considering the planar case and let us introduce the  {planar Delaunay--Poincar\'e variables} as follows:
\beqa{planar poinc}
\arr{
\L_i=\tilde m_i\sqrt{\hat m_i\,a_i}\\
\l_i=\ell_i+g_i\\
\eta_i=\sqrt[4]{\hat m_i\,a_i}\sqrt{2\tilde m_i(1-\sqrt{1-e_i^2})}\cos{g_i}\\
\xi_i=-\sqrt[4]{\hat m_i\,a_i}\sqrt{2\tilde m_i(1-\sqrt{1-e_i^2})}\sin{g_i}\\
}\qquad \qquad 1\leq i\leq N
\eeqa
where,  denoting by $\cE_i$ the ``osculating'' ellipse spanned by the  solution of the {\sl two -- body} differential problem
\beqano
\arr{
\ddot v=-\hat m_i\frac{v}{|v|^3}\ ,\quad v\in \real^2\\
\Big(\tilde m_i\dot v(0),\ v(0)\Big)=\Big(y_i,\ x_i\Big)
}
\eeqano
  $a_i$, $e_i$, $g_i$, $\ell_i$, are  the {\sl semimajor axis}, the {\sl eccentricity}, the {\sl argument of perihelion} of  $\cE_i$ and the {\sl mean anomaly}  of $x_i$ on $\cE_i$
(assume that $(y,x)$ varies in a region of ${\cal C}_{\rm cl}$ for which each  $(y_i,\ x_i)$ gives rise to an ellipse, \ie, with $0<e_i<1$). It is a classical result (see \cite{Che88}, \cite{Lask88}, \cite{BCV03}) that the map $$\phi_{\rm DP}^{-1}:\quad \Big(y,x\Big)\to\Big(\L,\l,\eta,\xi\Big)\ ,$$ with $\L=(\L_1,\cdots,\L_N)$, $\cdots$ as in (\ref{planar poinc}) is real--analytic $1$:$1$ and symplectic  on a suitable open neighborhood of $\{\bar\L\}\times\torus^N\times\{0\}$.

\vskip.1in
\noi
When expressed in planar  Delaunay--Poincar\'e coordinates, the integrable limit of   ${\cal H}_{\rm plt}$ becomes
$$h_{\rm plt}:=\sum_{1\leq i\leq N}h_i\circ\phi_{\rm DP}=-\sum_{1\leq i\leq N}\frac{\tilde m_i^3\hat m_i^2}{2\L_i^2}$$ a function of the actions $\L=(\L_1$, $\cdots$, $\L_N)$, only -- a fact  usually called ``proper degeneracy'', which prevents the use of {standard} KAM (Kolmogorov, Arnol'd, Moser) theory
\footnote{The theory, on the persistence, under suitable assumptions, of quasi-periodic motions for nearly--integrable Hamiltonian systems, developed in the late 60's by  Moser, (1962, \cite{Mos62}), Arnol'd (1963, \cite{Arn63Kolm}) on the basis of an early paper (1954) by Kolmogorov \cite{Kolm54}. For a review --and a complete proof--of the original Kolmogorov's  Theorem, see \cite{Ch08}. For related references, see also \cite{CellCh07}, \cite{Ch09}, \cite{Sa04}.}
in order to construct {\sl maximal} tori.

\vskip.1in
\noi
It is also well--known, since Laplace, that the ``secular perturbation'' of the planetary problem, \ie, the mean
$$\bar f_{\rm plt}:=\frac{1}{(2\p)^N}\int_{\torus^N}\sum_{1\leq i<j\leq N}\left(\frac{y_i\cdot y_j}{\bar m_0}-\frac{\bar m_i\bar m_j}{|x_i-x_j|}\right)\circ\phi_{\rm DP}\,d\l$$
has an {\sl elliptic} equilibrium point at $z:=(\eta,\xi)=0$, for any $\L$, \ie, it has an equilibrium point there and can be symplectically  put into the form
\beq{bar f plt}\bar f_{\rm plt}\circ \phi_{\rm diag}=\ovl{f_{\rm plt}\circ \phi_{\rm diag}}=\bar f_0(\tilde\L)+\sum_{1\leq i\leq N}\O_i(\tilde\L)\frac{\tilde\eta_i+\tilde\xi_i^2}{2}+{\rm O}_4\ ,\eeq 
where $\O_i$ are usually called {\sl Birkhoff invariants of the first order}.

\noi
We recall here the Theorem by Arnol'd.
\vskip.1in
\noi
{\scshape fundamental theorem}\ (Arnol'd, 1963, \cite{Arn63})\ {\sl Assume that
\begin{itemize}
\item[{\rm ({\scshape ft}$_0$)}] 
$
\cH(I,\varphi,p,q)=h(I)+\e\,f(I,\varphi,p,q)
$
is real--analytic on $U(r_0):=\overline\cI\times \torus^{\bar n}\times B^{2\hat n}_{r_0}(0)$, with $\cI$  an open, bounded, connected subset of $\real^{\bar n}$;
\item[{\rm ({\scshape ft}$_1$)}] $h$ is a diffeomorphism of an open neighborhood of $\ovl\cI$, with non degenerate Jacobian $\partial\o=\partial^2 h$ on such neighborhood;
\item[{\rm ({\scshape ft}$_2$)}] the mean perturbation $\bar f(I,p,q):=\frac{1}{(2\p)^{\bar n}}\int_{\torus^{\bar n}}f(I,\varphi,p,q)d\varphi$ has the form\footnote{By Birkhoff theory, a sufficient condition for (\ref{6normalform}) is that  $\bar f$ has an elliptic equilibrium point at the origin, with non resonant Birkhoff invariants of the first order $\O$, that is, \beq{6 non res}|\O(I)\cdot k|\neq 0\quad \textrm{for any}\quad I\in \bar \cI, \quad k=(k_1,\cdots,k_m)\in \integer^{\hat n}:\ \sum_{1\leq i\leq \hat n}|k_i|\leq 6\ .\eeq} 
\beqa{6normalform}
\bar f&=&f_0(I)+\sum_{1\leq i\leq \hat n}\O_i(I)J_i+\frac{1}{2}\sum_{1\leq i,j\leq \hat n}A_{ij}(I)J_iJ_j+\frac{1}{6}\sum_{1\leq i,j,k\leq \hat n}B_{ijk}(I)\,J_iJ_jJ_k\nonumber\\
&+&o_6
\eeqa
with $\dst J_i:=\frac{p_i^2+q_i^2}{2}$ and $o_6/|(p,q)|^6\to 0$;
\item[{\rm ({\scshape ft}$_3$)}] $A$ is {\sl non singular} on $\cI$, \ie,
\beqano
{\rm det}A(I)\neq 0\quad \textrm{for any}\quad I\in \bar \cI\ .
\eeqano
\end{itemize}
Then, for any $\k>0$, there exists $r_*$,  such that, for any
\beq{Arnold bounds}
0<r<r_*\qquad\textrm{and}\qquad 0<\e<r^8
\eeq
an $\cH$--invariant ` set $F(r)\subset U(r)=\cI\times \torus^n\times B^{2\hat n}_{r}$ may be found, with

\beq{k}\meas\Big(U(r)\setminus F(r)\Big)<\k\,\meas\Big(F(r)\Big)\eeq
consisting of $n(:=\bar n+\hat n)$--dimensional tori where the $\cH$--flow is $\vartheta\to \vartheta+\n\,t$.
}

\vskip.1in
\noindent
Arnol'd's estimate for the tori density in phase space is
\footnote{This estimate is not explicitely quoted into the statement of {\scshape ft}, but can be deduced as follows. Using  the original Arnol'd's notations  $(\epsilon,\m,n_0,n_1):=(r^2,\e,\bar n,\hat n)$, in the course of the proof, we find the condition $\d:=\epsilon^{1/T}<C\,\k$ with $T=16(n+4):=16(n_0+n_1+4)$: see on page  $144$, eq. $(4.2.5)$ with $\d^{(3)}$ defined below and p. $145$, eq. $(4.2.7)$.}

\beq{Arnold measure}
\meas\big(F(r)\big)>(1-c\,r^{1/(8(n+4))})\meas\Big(U(r)\Big)
\eeq
\vskip.1in
\noindent
As told before, Arnol'd applied his theorem to the planar three--body problem, checking, in particular, assumptions {(\scshape ft$_3$)},{(\scshape ft$_4$)} (in fact, {(\scshape ft$_2$)} is a consequence of {(\scshape ft$_3$)} and Birkhoff Theory, in view of (\ref{bar f plt})).

\vskip.1in
\noindent
The {\sl spatial} three body problem  (${\rm d}=3$, $N=2$) was solved, in $1995$, by  P. Robutel  (\cite{Rob95}; see also \cite{LaskRob95}). After performing the Jacobi, or {\sl nodes} (angular momentum) {\sl reduction}, he checked the assumptions of the {\scshape ft}, proving, so, the existence of (maximal) tori with $4$ frequencies.

\vskip.1in
\noi
The first complete proof of the existence of a positive measure set of quasi--periodic motions was given only in 2004, by J. F\'ejoz (\cite{Fej04},{\scshape theor\'eme} $60$), who, completing the investigations of M.Herman \cite{Her98}, in  the framework of a refined   $C^{\infty}$ KAM theory, stated the existence of a positive measure set of initial data giving rise to quasi--periodic motions with $3N-1$ frequencies, with their density going to $1$
as $\m\to 0$.

\vskip.1in
\noindent
Another  recent proof of Arnol'd's statement, but in the real--analytic framework of 2001 R\"ussmann theory \cite{Rus01},   may be found in \cite{ChPu08} (see also \cite{Pus06}). The real--analytic framework appears more natural for the many--body problem (\ref{Newton equations}), which is formulated using real--analytic functions.

\vskip.1in
\noindent
Both the proofs presented in \cite{Fej04}, \cite{ChPu08} are based on the check of   ``weak'' conditions on the first invariants $\O$ of $\bar f$   (suitable non--planarity conditions, sometimes called {\sl Arnol'd-- Pyartli, R\"ussmann} conditions, respectively), which, however, cannot be applied directly to $\cH_{\rm plt}$, due to the presence of two  ``secular resonances''. 

\vskip.1in
\noi
Letting, in fact,  ${\cal H}_{\rm plt}$ in {\sl spatial} Delaunay--Poincar\'e variables (definition \ref{Delaunay variables space}), the frequencies $\O$ correspond to $2N$ frequencies (related to the motions of perihelia and ascending nodes, respectively) $\s=(\s_1,\cdots,\s_N)$, $\zeta=(\zeta_1,\cdots,\zeta_N)$ which are found to verify the  (unique, \cite{Fej04} ) linear relations (up to linear combinations)
\beq{secular resonances}\zeta_N=0\ ,\quad \sum_{1\leq k\leq N}\s_k+\sum_{1\leq k\leq N-1}\zeta_k=0\ .\eeq  The former  relation in (\ref{secular resonances}) is known since Laplace; the latter  was firstly noticed by M. Herman, so it is usually called {\sl Herman's resonance} (for an interesting investigation on the Herman's resonance, see \cite{AbdAlb01}).
Owing to such secular resonances (in particular, the Herman's resonance), both the non--planarity conditions required by the KAM theories used in \cite{Fej04}, \cite{ChPu08} are violated by $\cH_{\rm plt}$. In order to overcome this problem, in \cite{Fej04} a modified Hamiltonian is introduced, next considered on the symplectic manifold of vertical angular momentum; in \cite{ChPu08}, the phase space is extended by adding an extra degree of freedom.  

\vskip.1in
\noindent
Notice  that the former  relation  in (\ref{secular resonances}) is actually a resonance of (low) order $1$, and also prevents the direct application of {\scshape ft}, making (\ref{6 non res}) false; the Herman's resonance is of higher order, $2N-1$, so, it violates  (\ref{6 non res})  only for  $N=2$, $3$.
\vskip.1in
\noindent
A direct attack to the problem, in the sense specified by {\scshape ft}, using a good set of coordinates which performs the {\sl angular momentum reduction},  has never been attempted. We outline that such a strategy  would lead also  to a more precise insight into the properties of the quasi--periodic motions (tori measure,  frequencies, $\cdots$).
\vskip.1in
\noindent
The problems which this thesis addresses  are the following.
\begin{itemize}
\item[({\scshape p}$_1$)]
(Section \ref{KAM and degeneracy}) Construction of a {\scshape ft}--like KAM Theorem ({\scshape theorem} $1$ below, for a simplified verson) in the real-- analytic class for properly degenerate systems, in order to obtain a fine measure estimate for the ``invariant set'' (roughly speaking, the set of  the KAM quasi--periodic trajectories) of a properly degenerate $\cH$, nearly an equilibrium point so as to
\item[({\scshape p}$_2$)]
(Section $3$) establish the existence of maximal quasi--periodic motions and estimate the measure of  Kolmogorov's invariant set for the plane planetary problem;
\item[({\scshape p}$_3$)]
(Section $4$) reduction of the angular momentum in the spatial planetary problem which leads to 
\item[({\scshape p}$_4$)]
(Section $5$) a proof of existence of KAM tori with $3N-1$ Diophantine frequencies (via a {\sl partial reduction} of the angular momentum) and measure of the invariant set;
\item[({\scshape p}$_5$)]
(Section $6$) a direct proof  of existence of $(3N-2)$--dimensional KAM tori via analytic theories of \cite{ChPu08} (full reduction).
\end{itemize}

\noi
We briefly explain our results. 
\vskip.1in
\noi
As for {\scshape p}$_1$, we prove {\scshape theorem} $1$ below (for a more general statement, see Theorem \ref{more general degenerate KAM}), which may be viewed as a refinement of {\scshape ft}: compare the bound on $\e$ (\ref{new bounds}) and the estimate for the tori measure (\ref{new measure}) for the invariant set with (\ref{Arnold bounds}), (\ref{Arnold measure}).

\vskip.1in
\noi
{\scshape theorem} $1$. {\sl Assume {\rm({\scshape ft}$_0$)}, {\rm({\scshape ft}$_1$)}, {\rm({\scshape ft}$_3$)} as in {\scshape ft} and 
\beqano
&&\textrm{\rm({\scshape ft}$_2^{'}$)}\ \bar f=f_0(I)+\sum_{1\leq i\leq m}\O_i(I)J_i+\frac{1}{2}\sum_{1\leq i, j\leq N} A_{ij}(I)J_iJ_j+o_4\nonumber\\
&&\textrm{where}\quad |\O(I)\cdot k|\neq 0\quad \textrm{for any}\quad I\in \bar \cI, \quad k=(k_1,\cdots,k_m)\in \integer^m:\ \sum_{1\leq i\leq m}|k_i|\leq 4
\eeqano
where $o_4/|(p,q)|^4\to 0$. Then,
there exist $r_*$, $0<c<1<C$, $b>0$ such that, for any 
\beq{new bounds}
0<r<r_*\quad \textrm{and}\quad 0<\e<c(\log{r^{-1}})^{-2b}
\eeq
an invariant set $\cK(\e,r)\subset \cI\times \torus^{\bar n}\times B^{2\hat n}_{C(\log{r^{-1}})^{-1}}(0)$ (``Kolmogorov set'') with measure
\beq{new measure}
\meas\Big(\cK(\e,r)\Big)\geq \Big(1-C\e^{1/2}(\log{r^{-1}}\Big)^{b}-Cr^{1/2}\Big)\meas \Big(U(r)\Big)\ .
\eeq
consisting of $n=\bar n+\hat n$--dimensional invariant tori, with ($\e\,r^{5/2},\t$)--Diophantine frequencies $\n$, where the motion is analytically conjugated to $\vartheta\to \vartheta+\n\,t$.
}

\vskip.1in
\noindent
The proof  of Theorem \ref{more general degenerate KAM} is made in two steps.
\begin{itemize}
\item[{\scshape (s$_1$)}] On one side, proof of a quantitative isofrequencial   KAM theorem particularly well suited for {properly degenerate quasi integrable Hamiltonians in action--angle variables}  \beq{action angle degenerate}H(J,\psi)=\bar{\rm h}(\bar J)+\hat{\rm h}(J)+{\rm f}(J,\psi)\eeq \ie, with integrable part which splits into the sum of two terms: $\bar{\rm h}$  (thought dominant), which depends only on a part of the action variables $\bar J$ and $\hat{\rm h}$ (thought small with respect to $\bar{\rm h}$) which depends on {\sl all} the actions $J$. 
The peculiarity of this theorem is of  choosing (the idea goes back to Arnol'd) two different scales for the tori frequencies to be kept fixed.
\item[{\scshape (s$_2$)}] 
On the other side, we reduce
the properly degenerate Hamiltonian ${\cal H}(I,\varphi,p,q)$ to the form (\ref{action angle degenerate}), with $\bar{\rm h}$  of order $1$, $\hat{\rm h}$  of order $\e r^2$ and the perturbation small ($\e r^{5/2}$). The reduction  is based on a non standard averaging theory, developed by Biasco et al. \cite{BCV03}, and Birkhoff Theory.
\end{itemize}
\vskip.1in
\noi
As a second step, we apply {\scshape th$1$} to the plane $(1+N)$--Body Problem. In order to do that, we compute explicitely the Birkhoff invariants of order $1$ and $2$,  expanding the perturbation of $\cH$ in plane Delaunay variables $(\L,\l,\eta,\xi)$, up to order $4$, after suitable diagonalization and Birkhoff Theory. 
If $\bar f_{\rm pl}$ denotes the mean perturbation of the plane problem in Delaunay variables, the Hessian matrix $\partial^2\bar f_{\rm pl}$ has the form
$$\left(
\begin{array}{lrr}
{\cal F}(\L)&0\\
0&{\cal F}(\L)
\end{array}
\right)$$
with ${\cal F}(\L)$ a symmetric $N\times N$ matrix. The first Birkhoff invariants are thus the eigenvalues of ${\cal F}(\L)$.
We introduce a small parameter, the maximum semimajor axes ratio $\d$, letting
$$\frac{a_i}{a_{i+1}}=\hat\a_i\,\d\ .$$ For small $\d$, the asymptotics of ${\cal F}(\L)$ is
$${\cal F}(\L)\approx \d^n\left(
\begin{array}{lrrrr}
f_1&O(\d^{n_{12}})&\cdots &O(\d^{n_{1k}})&\cdots\\
&f_2\d^{n_{22}}&\cdots &O(\d^{n_{2k}})&\cdots\\
&&\ddots&\vdots&\\
&&&f_k\d^{n_{kk}}&\cdots\\
&&&&\ddots
\end{array}
\right)$$
with $f_i$ with order $1$ in $\d$, $n_{hk}$ positive integers verifying
$$n_{kk}<n_{k+1,k+1}\quad \textrm{and}\quad n_{h-1,k},\ n_{h,k+1}>n_{hk}\ .$$
The eigenvalues of ${\cal F}$ are thus $\O_k=f_k\d^{n_kk}$ up to higher orders, and are thus non resonant. The Hamiltonian can be put in Birkhoff normal form
$$\bar f_0(\L)+\O(\L)\cdot J+\frac{1}{2}J\cdot { A}(\L)J+\cdots\quad \textrm{where}\quad J_i=\frac{\eta_i^2+\xi_i^2}{2}$$
and the Birkhoff invariants of order $2$ are the eigenvalues of the symmetric matrix ${\cal A}(\L)$. We finally prove that ${\cal A}(\L)$ as the asymptotics
$${A}(\L)\approx \d^p\left(
\begin{array}{lrrrrr}
\a_{11}&\a_{12}&O(\d^{p_{13}})&\cdots &O(\d^{p_{1k}})&\cdots\\
\a_{21}&\a_{22}&O(\d^{p_{23}})&\cdots&O(\d^{p_{2k}})&\cdots\\
&&\a_{33}\d^{p_{33}}&\cdots&O(\d^{p_{3k}})&\cdots\\
&&&\ddots&\vdots&\\
&&&&\a_{kk}\d^{p_{kk}}&\cdots\\
&&&&&\ddots\\
\end{array}
\right)$$
with
$\a_{ij}$ with order $1$ in $\d$, $p_{k+1,k+1}>p_{kk}$, $\a_{11}\a_{22}-\a_{12}\a_{21}\neq 0$, $\a_{kk}\neq 0$. This allows us checking that
$$\textrm{det}\,A=\d^{Np}(\a_{11}\a_{22}-\a_{12}\a_{21})\prod_{3\leq k\leq N}\a_{kk}+o(\d^{Np})$$
is nonvanishing, for small $\d$, concluding the proof.


\vskip.1in
\noindent
The extension of the  previous proof to the spatial problem in Delaunay variables is forbidden, by the presence of the above mentioned secular resonances, closely related to the rotation invariance of ${\cal H}_{\rm plt}$. A reduction of the number of degrees of freedom is however possible, with the use of the {\sl Deprit variables}  \cite{BOI82}, \cite{Dep83}.

\vskip.1in
\noi
The remarkable property of this new set of variables is to have, among their conjugated momenta, two coordinates of the angular momentum: the modulus $G$ and the third component $C_{\rm z}$. Their conjugated angles will be then cyclic variables. In particular, the conjugated angle $\zeta$ of $C_{\rm z}$ has the meaning of the ascending node longitude of the total angular momentum $C$, \ie, its third component, so it is an integral of the motion, too. 

\vskip.1in
\noi
When expressed in these new variables, only one external parameter will appear (the modulus $G$) for the reduced problem. This new set of variables can be regularized in a similar way to Poincar\'e's one for the regularization of the Delaunay variables. The new regularized variables $(\L,\l,\eta,\xi,p,q)$ for the reduced problem are of dimension $2N+2N+2(N-2)=2(3N-2)$. The variables $(\L,\l,\eta,\xi)$, of dimension $4N$, play  the same role as the Poincar\'e variables in the plane problem. The variables $(p,q)$ are related to the couples (inclinations, nodes): only $N-2$ couples may be chosen as independent, having fixed the modulus $G$ of the angular momentum and its verical component $C_{\rm z}$.  As consequence of the reduction, the  D'Alembert symmetries, existing in the plane problem, are broken; the origin of the new secular coordinates $z=(\eta,\xi,p,q)$ is no longer an equilibrium point and an {\sl elliptic singularity} appears, that is a singularity for the perturbation over the manifold 
$$G=\sum_{1\leq i\leq N}\left(\L_i-\frac{\eta_i^2+\xi_i^2}{2}\right)-\sum_{1\leq i\leq N-2}\frac{p_i^2+q_i^2}{2}$$
owing to which,
the configuration with {\sl all} zero eccentricities and  inclinations (which corresponds to $z=0$ and $G=\sum\L_i$) is not allowed. As a consequence, motions arbitrarily close to cocircular and coplanar trajectories cannot be considered, a fact already known in the case of  the three body problem.

\vskip.1in
\noi
Nonetheless Deprit's reduction has a {\sl partial reduction} ({\sl partial reductions} were also studied in \cite{Mal02}) naturally associated, which  consists in using only $C_{\rm z}$ as generalized momentum, and not $G$ also, making a further symplectic change of variables $(G,g)\leftrightarrow(p_{N-1},q_{N-1})$, where $g$ is the (cyclic) variable associated to $G$. In this way, a further inclination is treated as independent, but the number of degrees of freedom is enhanced from $3N-2$ to $3N-1$, having lost the cyclic variable $g$. Differently from what happens using Delaunay variables, Deprit's {\sl partial reduction} leaves the mean perturbation regular and even around the secular origin, which is thus an equilibrium point corresponding to zero eccentricities and mutual inclinations. Thus,  Deprit's partial reduction allows us to consider a larger region of the phase space than in the case of full reduction, even if at the price of one degree of freedom more.

\vskip.1in
\noi
In Section $5$, we show that the set of partially reduced Deprit variables provides a natural proof of the existence of $(3N-1)$--dimensional KAM tori via {\scshape th1}, at least for $N\geq 3$ planets, that is, from the four body problem on. In these variables, the planetary Hamiltonian (\ref{heliocentric hamiltonian}) 
$${\cal H}_{\rm plt,pr}=h_{\rm plt}(\L)+\m f_{\rm plt,pr}(\L,\l,\eta,\xi,p,q)\ ,$$
where $h_{\rm plt}$ is the usual Kepler's integrable part, satisfies the following. The
 ``secular perturbation'' $\dst \bar f_{\rm plt,pr}:=(2\p)^{-N}\int_{\torus^N}f_{\rm plt,pr}d\l$ is even and regular around the secular origin, as said before,
 and has the form
\beqa{Fejnot}\bar f_{\rm plt,pr}=\bar f^0_{\rm plt,pr}+{\cal Q}^*_h\cdot\frac{\eta^2+\xi^2}{2}+{\cal Q}^*_v\cdot\frac{p^2+q^2}{2}+{\eufm F}(\eta,\xi,p,q)+\cdots\eeqa
with ${\cal Q}^*_h$, ${\cal Q}^*_v$  suitable quadratic forms acting only on the ``horizontal'', ``vertical'' variables, respectively,
\footnote{Following \cite{Fej04}'s notation, in (\ref{Fejnot}), the dot ``.''  denotes contraction of indices:${\cal Q}\cdot\eta^2:=\sum_{i,j}{\cal Q}_{ij}\eta_i\eta_j$ if $\eta=(\eta_1,\cdots,)$.} 
${\eufm F}$  a suitable quartic form, all depending parametrically on $\L$ as well as $\bar f^0_{\rm plt,pr}$ and verifying the following. The respective sets of eigenvalues $s=(s_1$, $\cdots$, $s_N$), $z=(z_1$, $\cdots$, $z_{N-1}$) of ${\cal Q}^*_h$, ${\cal Q}^*_v$ (together also with the ``mean motions'' ${\rm n}=({\rm n}_1,\cdots,{\rm n}_N):=\partial_{\L}h_{\rm plt}$)
are found to satisfy the {\sl Herman's resonance} 
$$\sum_{1\leq i\leq N}s_i+\sum_{1\leq i\leq N-1}z_i=0\ ,$$
and {\sl only} that (Proposition \ref{true non resonance}). Since this resonance is of order $N+(N-1)=2N-1$, it does not violate the condition {\scshape ft}$_3^{'}$  of {\scshape th$1$} when $N\geq 3$, and the normal form of {\scshape ft}$_2^{'}$ can be constructed. This normal form turns out to be non degenerate, \ie, it also satisfies the second order non--degeneracy condition {\scshape ft$_4$} (Propositin \ref{spatial non degeneracy}). Both the proofs of non degeneracy (of first and second order) are inductive and  are developed with similar techniques as in \cite{Fej04}.  Then,  invoking {\scshape th}$1$, we can state the existence of $(3N-1)$--dimensional KAM tori for the planetary problem and thus estimate the density of the invariant set (Theorem \ref{3N-1KAMtori and meas})
$$1-\m^{1/2}(\log{\e^{-1}})^{b}-\e^{1/2}$$
into a ball with volume $\e^{2(2N-1)}$, where $\e$ is an upper bound for eccentricities and inclinations. Notice that the partial reduction generates an extra dimension for the KAM tori, relatated to the rotation invariance of ${\cal H}_{\rm plt,pr}$.

\vskip.1in
\noi
Nicely,  Deprit's partial reduction, for $N\geq 3$, makes us appear the spatial problem as the natural extension of the plane problem:  the set of  Birkhoff invariants of order $1$ has the planar one as subset, and the same happens at order $2$: when comparing the two matrices (planar and spatial) of Birkhoff invariants of order $2$, the planar one is a submatrix of the spatial one. This fact cannot be observed in the $3$--body problem as treated in \cite{Rob95}, because, there, the full reduction is made, and the coinclination of the two planets is expressed as a function of the eccentricities.

\vskip.1in
\noi
In Section $6$, we look at the {\sl full} reduction, that is, we use also the modulus of the angular momentum $G$ as generalized momentum.  As said before, this makes us gain a cyclic variable, the angle $g$ conjugate to $G$, lowering the number of  degrees of feedom to $3N-2$, but also causes, for $N\geq 3$, an {\sl elliptic singularity} and lack of symmetries (facts already known in Poincar\'e--Delaunay variables, trying to do a partial reduction, \ie, to eliminate one inclination with the use of the integral $C_{\rm z}$). A consequence of the lack of symmetries is that, for $N\geq 3$, the secular origin is no longer an equilibrium point for the mean perturbation.  Nonetheless, in the range of small eccentricities and inclinations, it is possible to find a new equilibrium point, which is ``small'', \ie, consists of almost circular and coplanar orbits, in the region of phase space which is considered, such that after suitable re--centering around it and  symplectic diagonalization of the quadratic part,  the planetary Hamiltonian is finally put into the form
$${\cal H}_{\rm plt}=h_{\rm plt}+\m f_{\rm plt}$$
where the mean $\bar f_{\rm plt}:=(2\p)^{-N}\int f_{\rm plt}$ becomes
$$\bar f_{\rm plt}=\hat f^0_{\rm plt}(\L,G)+\sum_{1\leq i\leq N}{\rm s}_i(\L,G)\frac{\eta_i^2+\xi_i^2}{2}+\sum_{1\leq i\leq N-2}{\rm z}_i(\L,G)\frac{p_i^2+q_i^2}{2}+O(3)$$ 
and the first Birkhoff invariants ${\rm s}=({\rm s}_1$, $\cdots$, ${\rm s}_N$), ${\rm z}=({\rm z}_1$, $\cdots$, ${\rm z}_{N-2}$), together with the mean motions ${\rm n}=({\rm n}_1,\cdots,{\rm n}_N)$, {\sl do not satisfy any linear relation}. Then, applying the real--analytic first order theory developed in \cite{ChPu08}, we can state the existence of $(3N-2)$--dimensional KAM tori for the planetary problem (Theorem \ref{3N-2 KAM tori}).
}
\pagina
\section{Properly Degenerate KAM Theory}\label{KAM and degeneracy}
\setcounter{equation}{0}
We recall  some basic notations and definitions.

\vskip.1in
\noi
Let $\bar n$, $\hat n\in \natural$, $n:=\bar n+\hat n$, $\cI$ a bounded, connected subset of $\real^{\bar n}$, $\torus^{\bar n}:=\real^{\bar n}/2\pi\integer^{\bar n}$  the usual real ``flat'' $\bar n$--torus and $B^{\hat n}_r(x)$  the   real open  $\hat n$--ball in $\real^{\hat  n}$  with radius $r$, centered at $x$.

\vskip.1in
\noi
In order that a compact set ${\cal T}\subset V$ is called  a $(\g,\hat \g;\t)$--{\sl Lagrangian torus with frequency} $\nu$ for a given
\beqano
\dst\cH(I,\varphi,p,q)=h(I)+f(I,\varphi,p,q)
\eeqano 
assumed to be real--analytic on the phase space $V:=\cI\times \torus^{\bar n}\times B^{2\hat n}_r(0)$,
we require that
\begin{itemize}
\item[({\scshape t}$_1$)]
there exists a real--analytic embedding onto
\beqano
\phi=(\phi_I,\phi_{\varphi},\phi_p,\phi_q):\ \vartheta=(\bar\vartheta,\hat\vartheta) \in\torus^{\bar n}\times\torus^{\hat n} \to \phi(\vartheta)\in {\cal T}
\eeqano
(and, hence, $2\p$ --periodic in each variable) 
given by
\footnote{For shortness, in (\ref{KAMtorus}), the symbol $\sqrt{2r}\cos{\psi}$ denotes the ${\hat n}$--vector \beqano(\sqrt{2r_1}\cos{\psi_1},\cdots,\sqrt{2r_{\hat n}}\cos{\psi_{\hat n}})\ ,\eeqano if $r=(r_1,\cdots,r_{\hat n})\in \real^{\hat n}$, $\psi=(\psi_1,\cdots,\psi_{\hat n})\in \torus^{\hat n}$, and similarly for $\sqrt{2r}\sin{\psi}$.}
\beq{KAMtorus}
\arr{
\phi_I(\vartheta)=\bar v(\vartheta)\\
\phi_{\varphi}(\vartheta)=\bar \vartheta+\bar u(\vartheta)\\
\phi_p(\vartheta)=p_0(\vartheta)+\sqrt{2\hat v(\vartheta)}\,\cos{(\hat \vartheta+\hat u(\vartheta))}\\
\phi_{q}(\vartheta)=q_0(\vartheta)+\sqrt{2\hat v(\vartheta)}\,\sin{(\hat \vartheta+\hat v(\vartheta))}\ ,
}
\eeq
such that
\item[({\scshape t}$_2$)]
the map
\beqano
\vartheta\to (\bar\vartheta+\bar u(\vartheta),\ \hat\vartheta+\hat u(\vartheta))
\eeqano
is a diffeomorphism of $\torus^{\bar n}\times\torus^{\hat n}$;
\item[({\scshape t}$_3$)]
${\cal T}$ is invariant under the 
$\cH$--flow $ \phi_t^\cH$, which, on $\cT$, acts as  a translation by $\n$, i.e.,
\beqano 
\phi^{-1}\circ\phi_t^\cH\circ \phi:\quad \vartheta \to\vartheta+\nu\,t \ ,\quad \forall\ \vartheta\in\torus^{\bar n}\times\torus^{\hat n}\ ;\eeqano

\item[({\scshape t}$_4$)] the torus {\sl frequency} $\nu\in\real^{\bar n+{\hat n}}$   belongs to the generalized Diophantine set \beqano\cD^{\bar n,{\hat n}}:=\bigcup_{\g,\hat \g>0 \atop \t>\bar n+{\hat n}} \cD^{\bar n,{\hat n}}_{\g,\hat \g;\t}\ ,\eeqano
where 
\beqano{\cal D}^{\bar n,{\hat n}}_{\g,\hat \g;\t}:=\bigcap_{0\neq k\in \integer^{\bar n}\times \integer^{{\hat n}}}{\cal D}^{\bar n,{\hat n}}_{\g,\hat \g;\t}(k)\eeqano
with
\beqa{two gammas}{\cal D}^{\bar n,{\hat n}}_{\g,\hat \g;\t}(k)&:=&\left\{\nu\in \real^{\bar n+{\hat n}}:\quad |\n\cdot k|\geq \frac{\g}{|k|^{\t}}\quad \textrm{if}\quad k=(\bar k,\hat k)\quad \textrm{with}\ \bar k\in \integer^{\bar n}\setminus\{0\}\ ,\right.\nonumber\\
&&\left. |\n\cdot k|\geq \frac{\hat\g}{|k|^{\t}}\quad \textrm{if}\quad \bar k= 0\in \integer^{\bar n}\right\}
\eeqa
\end{itemize}
We shall say that the embedding $\phi$ as in {$(T_1)$}$\div${$(T_3)$} {\sl realizes the Lagrangian torus} $\cT$.

\vskip.1in
\noi
We are now ready to quote the following refined version of Arnol'd's Fundamental Theorem \cite{Arn63}. For a simpler formulation, see Remark \ref{degenerate KAM}.

\vskip.1in
\noi
We assume that
\begin{itemize}
\item[{\rm({\scshape d}$_0$)}] 
\beq{prop deg ham}
\cH(\e;I,\varphi,p,q)=h(I)+\e\,f(I,\varphi,p,q)
\eeq
is real--analytic on $V(\bar r):=\ovl\cI\times\torus^{\bar n}\times B^{2\hat n}_{\bar r}(0)$, where $\cI$ is an open, bounded, connected subset of $\real^{\bar n}$;
\item[{\rm({\scshape d}$_1$)}] $\partial h$ is a diffeomorphism of an open neighborhood of $\ovl\cI$, with non singular Jacobian on such neighborhood;
\item[{\rm({\scshape d}$_2$)}] the mean perturbation: $$\bar f(I,p,q):=\frac{1}{(2\p)^{\bar n}}\int_{\torus^{\bar n}}f(I,\varphi,p,q)d\varphi$$ has the form
\beq{Birkh norm form}\bar f(I,p,q)=f_0(I)+\sum_{1\leq i\leq \hat n}\O_i(I)\frac{p_i^2+q_i^2}{2}+\frac{1}{2}\sum_{1\leq i\leq j\leq \hat n}A_{ij}(I)\frac{p_i^2+q_i^2}{2}\frac{p_j^2+q_j^2}{2}+o_4\eeq
where $o_4/|(p,q)|^4\to 0$ as $(p,q)\to 0$, where
\item[{\rm({\scshape d}$_3$)}] $\O=(\O_1,\cdots,\O_{\hat n})$ verifies
$$\min_{0<|k|\leq 4}\min_{\ovl\cI}|\O(I)\cdot k|>0\qquad\qquad \textrm{(``}4\textrm{{\sl --non resonance}'')}$$
and
\item[{\rm({\scshape d}$_4$)}] $A=(A_{ij})_{1\leq i,j\leq N}$ {non singular} on $\ovl\cI$, \ie,
$$\min_{\ovl\cI}|{\rm det}\,A|>0\qquad\qquad \textrm{(``{\sl non degeneracy}'')}$$
\end{itemize}
\begin{theorem}\label{more general degenerate KAM} Let $\bar n$, $\hat n\in \natural$, $n:=\bar n+\hat n$, $\bar\t>\bar n$, $\t>n$ and assume {\scshape d$_0\div$\scshape d$_4$} above.
Then, there exists $r_*$, $\g_*$, $\g^*$, $C_*>0$ such that, for any $0<r<r_*$ and $\bar\g$, $\g$, $\hat\g$ in the range
\beqa{how small gammas}
\arr{
\g_{*}\max\{\sqrt{\e}(\log{r^{-1}})^{\bar\t+1},\ 
\sqrt[3]{\e r}(\log{r^{-1}})^{\bar\t+1},\ r^2(\log{r^{-1}})^{\bar\t+1}\}<\bar\g<\g^*\\
\g_{*} r^{5/2}<\g<\g^*\\
\g_{*}r^{5/2}(\log_+{(r^5/\g^2)^{-1}})^{\t+1}<\hat\g<\g^*r^2\\
}
\eeqa
an invariant set ${\eufm K}(\e,r,\bar\g,\g,\hat\g)\subset V(\bar r)$ may be found, with measure
\beq{fund measure}\meas\Big({\eufm K}(\e,r,\bar\g,\g,\hat\g)\Big)\geq \left[1-C_*\left(\bar\g+\g+\frac{\hat\g}{r^2}+r^{\hat n/2}\right)\right]\meas\Big(V(r)\Big)
\eeq
consisting of Lagrangian tori $\{{\eufm T}_\n(\e,r,\bar\g,\g,\hat\g)\}_{\n}$ with generalized $(\bar n$, $\hat n$, $\g$, $\e\hat\g$, $\t)$--Dio-phantine frequencies $\n$.
\end{theorem}
\begin{remark}\label{rem on FT}\rm From the proof of Theorem \ref{more general degenerate KAM}, the following amplifications follow.
\begin{itemize}
\item[{\rm (D$_1$)}] A Cantor set ${\eufm J}_*(\e,r,\bar\g,\g,\hat\g)\subseteq \cI\times B^{\hat n}_{r^2}(0)$  and a bi--Lipschitz homeomorphism (onto) 
\beqano
&& \varpi_*(\e,r,\bar\g,\g,\hat\g;\cdot)=(\bar\varpi_*(\e,r,\bar\g,\g,\hat\g;\cdot),\e\hat\varpi_*(\e,r,\bar\g,\g,\hat\g;\cdot):\nonumber\\
&& {\eufm J}_*(\e,r,\bar\g,\g,\hat\g)\to {\eufm O}_*(\e,r,\bar\g,\g,\hat\g)\subset\cD^{\bar n,\hat n}_{(\g,\e\hat\g,\t)}
\eeqano
such that
\item[{\rm (D$_2$)}]
for any $\n$ $\in$ ${\eufm O}_*(\e,r,\bar\g,\g,\hat\g)$, the embedding
\beqano
{\eufm F}(\e,r,\bar\g,\g,\hat\g,\n;\cdot)&=&({\eufm F}_I(\e,r,\bar\g,\g,\hat\g,\n;\cdot),{\eufm F}_{\varphi}(\e,r,\bar\g,\g,\hat\g,\n;\cdot),\nonumber\\
& &{\eufm F}_p(\e,r,\bar\g,\g,\hat\g,\n;\cdot),{\eufm F}_q(\e,r,\bar\g,\g,\hat\g,\n;\cdot))\eeqano which realizes ${\eufm T}_\n(\e,r,\bar\g,\g,\hat\g)$ is given by
\beqa{eq of tori}
\left\{
\begin{array}{lrr}
{\eufm F}_I(\e,r,\bar\g,\g,\hat\g,\n;\vartheta)=\bar {\eufm j}_*(\e,r,\bar\g,\g,\hat\g,\n)+\bar {\eufm u}(\e,r,\bar\g,\g,\hat\g,\n;\vartheta)\\
{\eufm F}_{\varphi}(\e,r,\bar\g,\g,\hat\g,\n;\vartheta)=\bar{\vartheta}+\bar {\eufm v}(\e,r,\bar\g,\g,\hat\g,\n;\vartheta)\\
{\eufm F}_p(\e,r,\bar\g,\g,\hat\g,\n;\vartheta)={\eufm p}_0(\e,r,\bar\g,\g,\hat\g,\n;\vartheta)\\
\qquad\qquad\qquad\qquad\ +\sqrt{2\hat {\eufm j}_{*}(\e,r,\bar\g,\g,\hat\g,\n)+2\hat{{\eufm u}}(\e,r,\bar\g,\g,\hat\g,\n;\vartheta)}\\
\qquad\qquad\qquad\qquad\ \times\cos{[\hat \vartheta+\hat {\eufm v}_{}(\e,r,\bar\g,\g,\hat\g,\n;\vartheta)]}\\
{\eufm F}_q(\e,r,\bar\g,\g,\hat\g,\n;\vartheta)={\eufm q}_0(\e,r,\bar\g,\g,\hat\g,\n;\vartheta)\\
\qquad\qquad\qquad\qquad\ +\sqrt{2\hat {\eufm j}_{*}(\e,r,\bar\g,\g,\hat\g,\n)+2\hat {\eufm u}(\e,r,\bar\g,\g,\hat\g,\n;\vartheta)}
\\
\qquad\qquad\qquad\qquad\ \times\sin{[\hat \vartheta_{}+\hat {\eufm v}(\e,r,\bar\g,\g,\hat\g,\n;\vartheta)]}
\end{array}
\right.\nonumber\\
\eeqa
where ${\eufm j}_*(\e,r,\bar\g,\g,\hat\g;\n)$ is the $\varpi_*(\e,r,\bar\g,\g,\hat\g;\cdot)$--preimage of $\n$. 
\end{itemize}
Furtheremore, the unperturbed frequencies $\partial h$ may be chosen  $(\bar\g,\bar\t)$--Diophantine on ${\eufm J}_*$$(\e$, $r$, $\bar\g$, $\g$, $\hat\g)$ and the following bounds hold, for $\varpi_*(\e,r,\bar\g,\g,\hat\g,\cdot)$, 
and ${\eufm F}(\e,r,\bar\g,\g,\hat\g,\cdot,\cdot)$:
\beqa{bounds for varpi}
 \sup_{{\eufm J}_*(\e,r,\bar\g,\g,\hat\g)}|\bar\varpi_*-\partial h|&\leq& C_*\min\left\{\frac{\g}{(\log_+{(r^5/\g^2)^{-1}})^{\t+1}}\ ,\ \frac{\hat\g}{(\log_+{(r^5/\g^2)^{-1}})^{\t+1}}\ ,\right.\nonumber\\
&& \left.\frac{\bar\g}{(\log{r^{-1}})^{\bar\t+1}}\right\}+C_*\e\nonumber\\
 \sup_{{\eufm J}_*(\e,r,\bar\g,\g,\hat\g)}|\hat\varpi_*-\O|&\leq& C_*\min\left\{\frac{\g}{(\log_+{(r^5/\g^2)^{-1}})^{\t+1}}\ ,\ \frac{\hat\g}{(\log_+{(r^5/\g^2)^{-1}})^{\t+1}}\ ,\right.\nonumber\\
&& \left.\frac{\bar\g}{(\log{r^{-1}})^{\bar\t+1}}\right\}+C_*\e\frac{(\log{r^{-1}})^{2\bar\t+1}}{\bar\g^2}+C_* r^2\\
\nonumber\\
\label{bounds for FF}
\sup_{{\eufm O}_*(\e,r,\bar\g,\g,\hat\g)\times\torus^n} |\bar {\eufm u}| &\leq& C_*\frac{\hat\g}{\g}\min\left\{\frac{\g}{(\log_+{(r^5/\g^2)^{-1}})^{\t+1}}\ ,\ \frac{\hat\g}{(\log_+{(r^5/\g^2)^{-1}})^{\t+1}}\ ,\right.\nonumber\\
&&\left.\frac{\bar\g}{(\log{r^{-1}})^{\bar\t+1}}\right\}+ C_*\e\frac{(\log{r^{-1}})^{\bar\t}}{\bar\g}\nonumber\\ 
\sup_{{\eufm O}_*(\e,r,\bar\g,\g,\hat\g)\times\torus^n} |\hat{\eufm u}|&\leq& C_*\min\left\{\frac{\g}{(\log_+{(r^5/\g^2)^{-1}})^{\t+1}}\ ,\ \frac{\hat\g}{(\log_+{(r^5/\g^2)^{-1}})^{\t+1}}\ ,\right.\nonumber\\
&&\left.\frac{\bar\g}{(\log{r^{-1}})^{\bar\t+1}}\ ,\quad r^{5/2}\right\}\nonumber\\ 
\sup_{{\eufm O}_*(\e,r,\bar\g,\g,\hat\g)\times\torus^n} |\bar{\eufm v}| &\leq&C_* \max\left\{r^5\frac{(\log_+{(r^5/\g^2)^{-1}})^{2(\t+1)}}{\g^2}\ ,\ \frac{r^5(\log_+{(r^5/\g^2)^{-1}})^{2(\t+1)}}{\hat\g^2}\ ,\right.\nonumber\\ 
&&\left.  \frac{r^2(\log{r^{-1}})^{\bar\t+1}}{\bar\g}\ ,\ \frac{\e(\log{r^{-1}})^{2\bar\t+1}}{\bar\g^2}\ ,\ \frac{\e r(\log{r^{-1}})^{3\t+2}}{\bar\g^3}\ ,\ r\right\}\nonumber\\ 
\sup_{{\eufm O}_*(\e,r,\bar\g,\g,\hat\g)\times\torus^n}|\hat{\eufm v}| &\leq&C_* \max\left\{r^5\frac{(\log_+{(r^5/\g^2)^{-1}})^{2(\t+1)}}{\g^2}\ ,\ \frac{r^5(\log_+{(r^5/\g^2)^{-1}})^{2(\t+1)}}{\hat\g^2}\ ,\right.\nonumber\\ 
&&\left.  \frac{r^5(\log{r^{-1}})^{2(\bar\t+1)}}{\bar\g^2}\ ,\ r\right\}\nonumber\\ 
\sup_{{\eufm O}_*(\e,r,\bar\g,\g,\hat\g)\times\torus^n}|{\eufm p}_0|&\leq& C_*\frac{\sqrt{\e}(\log{r^{-1}})^{\bar\t+1}}{\bar\g}\ ,\nonumber\\
\sup_{{\eufm O}_*(\e,r,\bar\g,\g,\hat\g)\times\torus^n}|{\eufm q}_0|&\leq& C_*\frac{\sqrt{\e}(\log{r^{-1}})^{\bar\t+1}}{\bar\g}
\eeqa
\end{remark}
{\bf On the proof of Theorem \ref{more general degenerate KAM} and Remark \ref{rem on FT}.}\ The strategy for the proof  of Theorem \ref{more general degenerate KAM} is the following. Firstly, we construct an isofrequencial   KAM theorem (Theorem \ref{two scales KAM},  Section \ref{A two times scale KAM theorem}) which is well suited to {\sl properly degenerate quasi integrable Hamiltonians in action--angle variables}  $$H(J,\psi)={\rm h}(J)+{\rm f}(J,\psi)$$ \ie, with non degenerate integrable part  ($\partial^2{\rm h}\neq 0$)
which splits as
$${\rm h}(J)=\bar{\rm h}(\bar J)+\hat{\rm h}(J)\qquad J=(\bar J,\hat J)$$
where $\bar{\rm h}$  (thought dominant) depends only on a part of the action variables $\bar J$ and $\hat{\rm h}$ (thought small) depends on all the actions.
The peculiarity of this theorem is the one of  choosing (the idea goes back to Arnol'd \cite{Arn63}) generalized $(\g,\e\hat\g;\t)$--Diophantine frequencies (definition) and its smallness condition is
\beq{small cond two freq}c_0\frac{F\max\{M,\hat M^{-1},N\}}{(\min\{\g/M,\e\hat\g/\hat M,R\})^2}<1\eeq
where $c_0$ is a universal constant, $F$, $M$, $N$, $\hat M$, $R$, are a measure of
\footnote{As usual, if ${\rm h}$ is real--analytic on $\cI_\r$, the symbol $\partial\,{\rm h}$ denotes its gradient $(\partial_{I_1}{\rm h}, \cdots,\partial_{I_p}{\rm h})$; the Hessian $\partial^2{\rm h}$ is  the $p\times p$ matrix with entries $\dst \partial^2_{I_iI_j}{\rm h}$ (where $i$ is the row, $j$ the coloumn).}
${\rm f}$, $\partial^2{\rm h}$, $(\partial^2{\rm h})^{-1}$, $\partial^2\hat{\rm h}$ and the strenghth of $\hat J$, respectively. Next, we reduce the properly degenerate Hamiltonian  (\ref{prop deg ham}) to the form
\beq{good form}{\rm H}={\rm h}_0(\bar J)+\e{\rm h}_1(J)+\e r^5{\rm f}(J,\psi):={\rm h}(J)+\e r^5{\rm f}(J,\psi)\ .\eeq
where ${\rm h}_0$ is $\e$--close to the unperturbed $h$ of $\cH$.
This reduction  is based on averaging theory, developed by Biasco et al., \cite{BCV03} (which carries itself a smallness condition $\bar\g>\g_*\sqrt{\e}(\log{r^{-1}})^{\bar\t+1}$ for the diophantine constant of the unperturbed frequencies $\partial{ h}$), Birkhoff Theory (Appendix \ref{app:Birkhoff Normal Form}) and use of symplectic polar coordinates. Taking then $F=C\,\e r^5$, $M=C$, $N=C\e^{-1}$, $\hat M=C\e$, $R=Cr^2$ (this is due to the use of symplectic polar coordinates), the smallness condition (\ref{small cond two freq}) essentially becomes
\beq{small cond two freq appl}c_1\frac{r^5}{\min\{\g^2,\hat\g^2,r^4\}}<1\eeq
and it will be guaranteed as soon as $r$ is small and $\g$, $\hat\g$ are chosen not smaller that $\g_*r^{5/2}$. The condition $\hat\g<\g^*r^2$ is necessary to find a not empty  $\partial {\rm h}$--pre image of $\cD^{\bar n,\hat n}_{\g,\e\hat\g}$. Observe the cancellation of $\e$ from (\ref{small cond two freq}) to (\ref{small cond two freq appl}), which makes us take $r$ as perturbative parameter.
\begin{remark}\rm({\sl Proof of} {\scshape theorem} $1$ {\sl and other details})\label{degenerate KAM}
The formulation we have chosen for Theorem \ref{more general degenerate KAM} is very general. Even if the parameters $\bar\g$, $\g$, $\hat\g$, in principle, might assume any value in the ranges (\ref{how small gammas}), nonetheless, in order to get the tori density as large as possible, the gamma--constants $\bar\g$, $\g$, $\hat\g$ (which are related to the amount of irrationality of the the unperturbed frequencies $\partial h$ and tori frequencies $\n$) should be taken as small a  possible. The choice
\beqa{good gammas}
\arr{
\bar\g= \g_0\max\left\{\sqrt{\e}(\log{r^{-1}})^{\bar\t+1}\ ,\ \sqrt[3]{\e r}(\log{r^{-1}})^{\bar\t+1},\ r^2(\log{r^{-1}})^{\bar\t+1}\right\}\\
\g=\hat\g =\g_0 r^{5/2}
}
\eeqa
(with a fixed  $\g_0>\g_*(\log{\g_*})^{\t+1}$) leads to an invariant set \beq{cal K}\cK(\e,r):={\eufm K}|_{(\bar\g,\g,\hat\g)=(\ref{good gammas})}\eeq
with  density just the one announced in the Introduction (as $\e r\leq \max\{\e^2,r^2\}$):
\beqano
1-\g_0\left(\bar\g+\g+\frac{\hat\g}{r^2}+r^{\hat n/2}\right)&\geq& 1-\hat C\Big(\sqrt{\e}(\log{r^{-1}})^{\bar\t+1}+\sqrt[3]{\e r} (\log{r^{-1}})^{\bar\t+1}\nonumber\\
&+&r^2(\log{r^{-1}})^{\bar\t+1}+r^{5/2}\nonumber\\
&+&\sqrt{r}+r^{\hat n/2}\Big)\nonumber\\
&\geq&1-\tilde C\Big(\sqrt{\e}(\log{r^{-1}})^{\bar\t+1+b}+\sqrt{r}\Big)
\eeqano
Furthermore, denoting by
\beqano
\arr{
\cJ_*(\e,r):={{\eufm J}_*}|_{(\bar\g,\g,\hat\g)=(\ref{good gammas})}\\
\\
\o_*(\e,r,\cdot)=(\bar\o_*(\e,r,\cdot),\e\hat\o_*(\e,r,\cdot)):=
\varpi=(\bar\varpi,\e\hat\varpi)|_{(\bar\g,\g,\hat\g)=(\ref{good gammas})}\\
\\
{\cal F}(\e,r,\cdot,\cdot)=({\cal F}_I(\e,r,\cdot,\cdot),{\cal F}_\varphi(\e,r,\cdot,\cdot),{\cal F}_p(\e,r,\cdot,\cdot),{\cal F}_q(\e,r,\cdot,\cdot)):=\\
{\eufm F}=({\eufm F}_I,{\eufm F}_\varphi,{\eufm F}_p,{\eufm F}_q)|_{(\bar\g,\g,\hat\g)=(\ref{good gammas})}\\
}
\eeqano
with
\beqano
\left\{
\begin{array}{lrr}
{\cal F}_I(\e,r,\n;\vartheta)=\bar j_*(\e,r,\n)+\bar u(\e,r,\n;\vartheta)\\
{\cal F}_{\varphi}(\e,r,\n;\vartheta)=\bar{\vartheta}+\bar v(\e,r,\n;\vartheta)\\
{\cal F}_p(\e,r,\n;\vartheta)=p_0(\e,r,\n;\vartheta)+\sqrt{2\hat j_{*}(\e,r,\n)+2\hat{u}(\e,r,\n;\vartheta)}\\
\qquad\qquad\qquad\ \times\cos{[\hat \vartheta+\hat v_{}(\e,r,\n;\vartheta)]}\\
{\cal F}_q(\e,r,\n;\vartheta)=q_0(\e,r,\n;\vartheta)+\sqrt{2\hat j_{*}(\e,r,\n)+2\hat u(\e,r,\n;\vartheta)}
\\
\qquad\qquad\qquad\ \times\sin{[\hat \vartheta_{}+\hat v(\e,r,\n;\vartheta)]}
\end{array}
\right.\nonumber\\
\eeqano
then, the bounds (\ref{bounds for varpi})$\div$(\ref{bounds for FF}) imply
\beqa{new bounds in r e}
 \sup_{\cJ_*(\e,r)}|\bar\o_*-\partial h|&\leq& C_*\max\left\{\sqrt{\e},\ \sqrt[3]{\e r},\ r^{5/2}\right\}\nonumber\\
 \sup_{\cJ_*(\e,r)}|\hat\o_*-\O|&\leq& \max\left\{\sqrt{\e},\ \sqrt[3]{\e r},\ (\log{r^{-1}})^{-1}\right\}\nonumber\\
 \sup_{\O_*(\e,r)\times\torus^n} |\bar u| &\leq& C_*\e r^{5/2}\nonumber\\
 \sup_{\O_*(\e,r)\times\torus^n} |\hat u|&\leq& C_*r^{5/2}\nonumber\\
 \sup_{\O_*(\e,r)\times\torus^n}|p_0|,\  \sup_{\O_*(\e,r)\times\torus^n}|q_0|&\leq& C_*(\log{r^{-1}})^{-1}
\eeqa 
\end{remark}
\begin{remark}[A physical comment]\rm
The tori in $\cK(\e,r)$ may be thought of the  analytic continuation of the tori
\beqa{where tori are}
\arr{
I=\bar I\\
\varphi\in \torus^n\\
(p_i-p_{0i})^2+(q_i-q_{0i})^2=\bar R_i^2\qquad C r^{5/2}\leq R_i^2\leq cr^2\\
}
\eeqa
being crossed by the system with frequencies $\sqrt{\e}+(\log{r^{-1}})^{-1}$--close to $(\partial h,0)$. Observe, in particular, that $p_0$, $q_0$, $\bar u$ go to $0$ with $\e$, for any fixed $0<r<r_*$. 
\vskip.1in
\noi
Notice  that if \beq{how small eps}\e<\const r^2(\log{r^{-1}})^{-2(\bar\t+1)}\ ,\eeq if $\bar\g$ is not too small, for instance, $\bar\g=\const \sqrt{\e}(\log{r^{-1}})^{-\bar\t+1}/r$ we than have $$|(p_0,q_0)|\leq\const\sqrt{\e}(\log{r^{-1}})^{2(\bar\t+1))}\bar\g^{-2}=\const\, r\ ,$$ \ie, the $(p,q)$ variables into a ball of radius $r$ around the origin (comapare also (\ref{where tori are})). This agrees with the result obtained in  Arnol'd's {\scshape fundamental theorem}. When (\ref{how small eps}) is no longer satisfied, we generally have a set of invariant tori for which the $(p,q)$--variables can stay away from the origin as far as $(\log{r^{-1}})^{-1}$: compare (\ref{new bounds in r e}) above.
\end{remark}
\subsection{A Two Times Scale KAM Theorem}\label{A two times scale KAM theorem}
In order to state Theorem \ref{two scales KAM} below, we introduce some  useful notations and definitions.
\begin{itemize}
\item[(i)]The $r$--neighorhood $\cI_r$ of $\cI\subset \real^p$ compact and the $s$--neighborhood $\torus^{p}_s$ of $\torus^p$ are defined as:
\beqano
\cI_r:=\bigcup_{I\in \cI}D^p_r(I)\ ,\quad \torus^p_s:=\{\varphi=(\varphi_1,\cdots,\varphi_p)\in \torus_{\complex}^p:\ |\Im{\varphi_i}|<s)\}
\eeqano
where $\torus_{\complex}^p:=\complex^p/2\p\integer^p$ and
\beqano
D^p_r(I)=\{I'\in \complex^p:\ |I'-I|<r\}
\eeqano
is the ususal complex open $p$--ball, where $\complex^p$ is equipped with the standard Euclidean norm: $\dst|(I_1,\cdots,I_p)|:=\sqrt{\sum_{1\leq i\leq p}|I_i|^2}$.
\item[(ii)]
Real--analytic functions $\dst {\rm f}:\ \cP\subset \real^{2p}\to \real$ on compact sets $\cP=\cI\times \torus^{p}\subset\real^{2p}$ ($\dst\cI\subset \real^{p}$ compact) are identified with their analytic extensions $\dst\bar {\rm f}:\ \cP_{r,s}\subset \complex^{2p}\to \complex$ over a suitable $(r,s)$--neighborhood $\cP_{r,s}=\cI_r\times \torus^{p}_s$ of their real domain. 
\item[(iii)] The ``sup--Fourier'' norm $\|{\rm f}\|_{r,s}$ of a real--analytic function ${\rm f}$ on $\cP_{r,s}$ is 
\beqano
\|{\rm f}\|_{r,s}:=\sum_{k\in \integer^p}\sup_{I\in\cI_r}|{\rm f}_k(I)|e^{|k|s}
\eeqano
where $|k|$ is the $1$--norm
\beqano
|(k_1,\cdots,k_p)|:=\sum_{1\leq i\leq p}|k_i|
\eeqano
and
\beqano
{\rm f}_k(I):=\frac{1}{(2\p)^p}\int_{\torus^p}{\rm f}(I,\varphi)e^{ik\cdot\varphi}d\varphi
\eeqano
is its $k^{th}$ Fourier coefficient.
\item[(iv)] If $A$ is a $n\times n$ matrix and $1\leq p$, $q\leq n$ the symbol $A^{[p,q]}$ denotes the $p\times q$ sub--matrix  of the {\sl first} $p$ rows and the {\sl last} $q$ coloumns of $A$, \ie, the matrix with elements $$A^{[p,q]}_{ij}=A_{i,n-q+j}\quad 1\leq i\leq p\ ,\quad 1\leq j\leq q\ .$$
Conversely, $A^{[p,q]}_*$ denotes the $p\times q$ sub--matrix of $A$ of the {\sl last} $p$ rows and {\sl first} $q$ coloumns.
\item[(v)] We recall that
$f: \cI\subset \real^p\to \real^p$ is {\sl Lipschitz} if 
$${\cal L}_{\|\cdot\|}(f):=\sup_{I\neq I'\in \cI}\frac{\|f(I)-f(I')\|}{\|I-I'\|}<+\infty$$
(with $\|\cdot\|$ a somewhat norm  of $\real^p$) . For a Lipschitz function $f$, we denote by ${\cal L}(f)$ the number
$${\cal L}(f):=\sup_{I\neq I'\in \cI}\frac{|f(I)-f(I')|}{|I-I'|}$$
(with respect to the Euclidean norm) and call it {\sl Lipschitz constant} for $f$.
We define the $\r$--{\sl Lipschitz norm} of $f$ on $\cI$
$$\|f\|^{\rm Lip}_{\r,{\cI}}:=\r^{-1}\sup_{{\cI}}|f|+{\cal L}(f)\ .$$
\item[(vi)]
 $f$ is called {\sl bi--Lipschitz} if $f$ is Lipschitz, injective, with Lipschitz inverse, or, equivalently, if there exist $0<{\cal L}_-(f)\leq{\cal L}_+(f)$, called {\sl Lipschitz constants}, such that
$${\cal L}_-(f)|I-I'|\leq|f(I)-f(I')|\leq {\cal L}_+(f)|I-I'|\quad \textrm{for all}\quad I,I'\in \cI\ ,$$
where
$${\cal L}_+(f)={\cal L}(f)\ ,\quad {\cal L}_-(f)=\frac{1}{{\cal L}(f^{-1})}\ .$$
\end{itemize}
\begin{theorem}\label{two scales KAM}
Let $\bar n$, $\hat n\in \natural$, $n:=\bar n+\hat n$, $\t>n$, $\g\geq\hat\g>0$, $0<2s\leq \bar s<1$, $\bar\cI\subset \real^{\bar n}$, $\hat\cI\subset\real^{\hat n}$  $\cI:=\bar\cI\times \hat\cI$ such that
\beqano
{\rm H}
(J,\psi)={\rm h}(J)+{\rm f}(J,\psi)
\eeqano 
real--analytic on $\cI_\r\times \torus_{\bar s+s}^{n}$. Assume that $\o:=\partial {\rm h}$
is a diffeomorphism of $\cI_\r$ and the Hessian matrix $U:=\partial^2{\rm h}$ is non singular on $\cI_\r$.
Let
\beqano
&&  \dst M\geq\sup_{{\cal I}_{\r}}\|U\|\\
&&  \dst \hat M\geq\sup_{{\cal I}_{\r}}\|U^{[n,\hat n]}\|\\
&&  \dst N\geq\sup_{{\cal I}_{\r}}\|T\|\\
&&  \dst \bar N\geq\sup_{{\cal I}_{\r}}\|T^{[\bar n,n]}\|\\
&&  \dst \hat N\geq\sup_{{\cal I}_{\r}}\|T_*^{[\hat n,n]}\|\\
&&  \dst F\geq\|{\rm f}\|_{\r,s}
\eeqano
where $T:=U^{-1}$, define
\beqano
&&  \dst c:=\max\left\{{2^{11}n},\ \frac{2}{3}(12)^{\t+1}\right\}\\
&&  \dst K:=\frac{6}{s}\ \log_+{\left(\frac{FM^2\,L}{\gamma^2}\right)^{-1}}\quad \textrm{where}\quad \log_+(a):=\max\{1,\log{a}\}\\ 
&&  \dst \tilde{\r}:=\min\left\{\frac{\g}{3MK^{\t+1}}\ ,\ \frac{\hat\g}{3\hat MK^{\t+1}}\ ,\  
 \r\right\}\\
&&  \dst L:=\max\ \{N\ , M^{-1},\hat M^{-1}\}\\
\eeqano
and assume the ``perturbation'' $\rm f$ so ``small'' that the following ``KAM condition'' holds
\beq{KAM cond}
cE:=c\frac{FL}{\tilde{\r}^2}<1\ .
\eeq
\begin{itemize}
\item[(i)]
Then, for any frequency $\nu\in\O_*:=\o(\cI)\cap\cD^{\bar n,\hat n}_{\g,\hat\g,\t}$, there exists a unique  Lagrangian KAM torus ${\rm T}_\n\subset \Re({\cal I}_{34\tilde{\r}E})\times \torus^{n}$ for $H$ with frequency $\nu$, such that the following holds. 
There exists a ``Cantor'' set ${\cal I}_{*}\subset \Re ({\cal I}_{32\tilde{\r}E})$  and a   bi--Lipschitz (onto) {homeomorphism} 
\beq{omega*}
\omega_*=(\bar\o_*,\hat\o_*):\ {\cal I}_{*}\to \O_*
\eeq
satisfying
\beqa{***}
&&\label{cantor close to predioph}\sup_{{\cal O}_*}|\bar \o_*^{-1}-\bar\o^{-1}|\leq 2^5\frac{\bar N}{N}\tilde\r\,E\ ,\quad \sup_{{\cal O}_*}|\hat \o_*^{-1}-\hat\o^{-1}|\leq 2^5\frac{\hat N}{N}\tilde\r\,E\\
&&\sup_{{\cal I}_*}|\bar\omega_{*}-\bar\omega|\leq{2^5}M\tilde{\r}E\  ,\quad \sup_{{\cal I}_*}|\hat\omega_{*}-\hat\omega|\leq{2^5}\hat M\tilde{\r}E
\label{frequency map}\\
&& \label{Lip Norm}
\|\o_*^{-1}\circ \o-\textrm{\rm id}\|^{\rm Lip}_{\tilde{\r},{\cal I}_{\g,\hat \g,\t}}\leq 2^{11}E\ ,
\qquad
{\cal I}_{\g,\hat \g,\t}:=\o^{-1}({\cal D}^{n,\hat n}_{\g,\hat \g,\t})\cap{\cal I}\ .
\eeqa
such that ${\rm T}_\n$ is
realized by the real--analytic embedding $\phi_\n=(\phi_{\n I},\phi_{\n \varphi})$ given by
\begin{eqnarray}\label{perturbed torus}
\left\{
\begin{array}{lrr}
\phi_{\n I}(\vartheta)=I_{*}(\nu)+v(\nu,\vartheta)&&\\
\phi_{\n \varphi}(\vartheta)=\vartheta+u(\nu,\vartheta)&&\\
\end{array}
\right. \quad \vartheta\in \torus^n\ ,
\end{eqnarray}
where $I_*(\n):=\o_*^{-1}(\n)$ and $v=(\bar v,\hat v)$, $u=(\bar u,\hat u)$ are bounded as
\begin{eqnarray}\label{alpha beta}
|\bar v(\nu,\vartheta)|\leq2\,\frac{\hat\g}{\g}E\,\tilde{\r}\ ,\quad |\hat v(\nu,\vartheta)|\leq2\,E\,\tilde{\r}\ ,\quad |u(\nu,\vartheta)|\leq2\,E\,s
\end{eqnarray}
\item[(ii)] The measure of the invariant set ${\rm K}=\phi_{\O_*}(\torus^n)$ satisfies
\beqa{tori measure}
\meas\Big(\cI\times\torus^n\setminus{\rm K}\Big)&\leq&\Big(1+(1+2^7E)^{2n}\Big)\Big(1+(1+2^{10}E)^{n}\Big)\meas\Big(\cI\setminus \cI_{\g,\hat\g,\t}\times \torus^n\Big)\nonumber\\
&+&\Big(1+(1+2^7E)^{2n}\Big)(1+2^{10}E)^n\meas\Big(\cI_{\r_1}\setminus\cI\times\torus^n\Big)\nonumber\\
&+&(1+2^7E)^{2n}\meas\Big(\cI_{\r_2}\setminus \cI\times\torus^n\Big)\ .
\eeqa
where
$\r_1=2^6E\tilde\r/(1-2^{10}E)$, $\r_2=4E\tilde\r/(1-2^7E)$.
\end{itemize}
\end{theorem}

\vskip.1in
\noindent
\subsubsection{Construction of the Approximating Sequences}
The proof of Theorem \ref{two scales KAM} is obtained by infinite iterations of real--analytic symplectomorphisms, converging over a Cantor set. 
 Each iteration is based on Lemma \ref{average} below.
 
\vskip.1in
\noi
Let $\bar n$, $\hat n$, $n=\bar n+\hat n\in \natural$, ${\cal I}\subset \real^{\bar n}\times\real^{\hat n}$. Following P\"oschel \cite{POSCH93}, we define the $\cP$--norm on $\cI_\r\times \torus^{n}_s$ as
$$|\Big((\bar I,\hat I),\varphi)\Big)|_{\cP}:=\max\{|\bar I|_1,\ |\hat I|_1,\ |\varphi|_{\infty}\}\ ,$$
where, as usual,
\beqano
|\bar I|_1:=\sum_{1\leq i\leq \bar n}|\bar I_i|\ ,\quad |\hat I|_1:=\sum_{1\leq i\leq \hat n}|\hat I_i|\ ,\quad |\varphi|_{\infty}:=\max_{1\leq i\leq n}|\varphi_i|
\eeqano
if $\bar I=(\bar I_1,\cdots,\bar I_{\bar n})$, $\hat I=(\hat I_1,\cdots,\hat I_{\hat n})$, $\varphi=(\varphi_1,\cdots,\varphi_{n})$.
We also introduce the matrices  $W^{p}_{\a,\b}$, $W_{\a,\b}$ ($1\leq p\leq 2n$) whch are defined as the $(2p)\times(2p)$ diagonal matrices  $$W^{p}_{\a,\b}:=(\a^{-1} \id_{p},\b^{-1} \id_{2n-p})\ ,\quad W_{\a,\b}:=W^{n}_{\a,\b}=(\a^{-1} \id_{n},\b^{-1} \id_{n})$$
 where
 $\id_p$ is the identity matrix with order $p$. 
\begin{lemma}[Averaging Theorem]\label{average}
Let ${\rm H}(I,\varphi)={\rm h}(I)+{\rm f}(I,\varphi)$ real--analytic on\\ $\cP_{r,\bar s+s}$ $:={{\cal I}}_{r}$ $\times \torus^n_{\bar s+s}$. 
Assume that $\o:=\partial\,{\rm h}$ verifies
\begin{eqnarray}\label{non resonance in avarage thm}
|\o(I)\cdot k|\geq \left\{\begin{array}{lrr}
\bar\a\quad \textrm{for}\quad k=(\bar k,\hat k)\in \integer^{\bar n}\times \integer^{\hat n}\setminus \L\quad \bar k\neq 0\ ,\quad |k|\leq K \\
\hat\a \quad \textrm{for}\quad k=(0,\hat k)\in \{0\}\times \integer^{\hat n}\setminus \L\quad 0<|\hat k|\leq K\\
\end{array}
\right.
\end{eqnarray}
\noi
where $\L\subseteq \integer^n$.
If 
\begin{eqnarray}\label{tilde rho alpha}
\|{\rm f}\|_{r,\bar s+s}\leq \frac{1}{2^7}\,\frac{\a\,r}{K}\ ,
\end{eqnarray}
where $\a:=\min\{\bar\a,\hat\a\}$ and $Ks\geq 6$, then, there exists a real--analytic, symplectic coordinate transformation $$\Psi:\ \cP_{r/2,\bar s+s/6}\to \cP_{r,\bar s+s}$$
such that
$${\rm H}\circ \Psi(I,\varphi)={\rm h}(I)+g(I,\varphi)+{\rm f}_*(I,\varphi)$$
with $g(I,\varphi)=\sum_{k\in \L}g_k(I)e^{ik\cdot\varphi}$ is $\L$--completely resonant and
$$\|g-{\rm f}_0\|_{r/2,\bar s+s/6}\leq \frac{2^5K}{\a\,r}\|{\rm f}\|_{r,\bar s+s}^2\ ,\quad \|{\rm f}_*\|_{r/2,\bar s+s/6}\leq e^{-Ks/6}\|{\rm f}\|_{r,\bar s+s}$$
where ${\rm f}_0:=P_{\L}T_K\,{\rm f}$.
Moreover, the following bounds hold, uniformly on $\cP_{r/2,\bar s+s/6}$
\begin{eqnarray}\label{identity bounds}
|W^{\bar n}_{\a/\bar\a,1}W_{r,s}(\Psi-\id)|_{\cal P}\leq \frac{2\,K\|{\rm f}\|_{r,\bar s+s}}{\a\,r}\ ,
\end{eqnarray}
and
\begin{eqnarray}\label{derivatives0}
\|W^{\bar n}_{\a/\bar\a,1}(W_{r,s}D\Psi W^{-1}_{r,s}-\id_{2n})\|_{\cal P}\leq \frac{2^6K}{\a\,r}\,\|{\rm f}\|_{r,\bar s+s}\ .
\end{eqnarray}
where $\|\cdot\|_{\cal P}$ denotes the operatorial norm induced by $|\cdot|_{\cal P}$
\footnote{I.e., $\|A\|_{\cal P}:=\sup_{z=(I,\varphi),\ |z|_{\cal P}=1}|A z|_{\cal P}$}
.
\end{lemma}
Lemma \ref{average} is a useful remake of the {\sl Normal Form Lemma} of \cite{POSCH93}. For sake of completeness, its proof may be found in Appendix \ref{app:average}.
\vskip.1in
\noi
%
\begin{lemma}[Iterative Lemma]\label{itera}
Let  $0<\hat\g\leq\g$, $\bar\cI\subset \real^{\bar n}$, $\hat\cI\subset\real^{\hat n}$ compact sets, put $\cI:=\bar\cI\times \hat\cI$ and let
\beqano
{\rm H}
(J,\psi)={\rm h}(J)+{\rm f}(J,\psi)
\eeqano 
real--analytic on $\cI_\r\times \torus_{\bar s+s}^{n}$, with $\bar s>0$, $0<s<1$. Assume that $\o:=\partial {\rm h}$ is a diffeomorphism of $\cI_\r$ with Jacobian matrix $U:=\partial^2{\rm h}$ non singular on $\cI_\r$ and $\o(\cI)\subseteq \cD^{\bar n,\hat n}_{\g,\hat\g,\t}$.
Let
\beqano
&&  \dst M\geq\sup_{{\cal I}_{\r}}\|U\|\\
&&  \dst \hat M\geq\sup_{{\cal I}_{\r}}\|U^{[n,\hat n]}\|\\
&&  \dst N\geq\sup_{{\cal I}_{\r}}\|T\|\\
&&  \dst \bar N\geq\sup_{{\cal I}_{\r}}\|T^{[\bar n,n]}\|\\
&&  \dst \hat N\geq\sup_{{\cal I}_{\r}}\|T_*^{\hat n,n}\|\\
&&  \dst F\geq\|{\rm f}\|_{\r,\bar s+s}
\eeqano
where $T:=U^{-1}$, define
\beqano
&&  \dst K:=\frac{6}{s}\ \log_+{\left(\frac{FM^2\,L}{\gamma^2}\right)^{-1}}\\ 
&&  \dst \tilde{\r}:=\min\left\{\frac{\g}{3MK^{\t+1}}\ ,\ \frac{\hat\g}{3\hat MK^{\t+1}}\ ,\  
 \r\right\}\\
&&  \dst L:=\max\ \{N\ , M^{-1},\hat M^{-1}\}\\
\eeqano
and assume 
\beq{KAM condnew}
2^8nE:=2^8n\frac{FL}{\tilde{\r}^2}\leq1\ .
\eeq
Then, a set ${\cal I}_+\subset\cI_{\tilde\r/32}$,  two  numbers $\r_+>0$, $0<s_+<1$ and a symplectic analytic transformation 
\begin{eqnarray*}
&\Psi&:\ {\cal I}_{+\rho_+}\times \torus^n_{\bar s+s_+}\to {\cal I}_{{\rho}}\times \torus^n_{\bar s+s}\ ,\nonumber\\
& &(I_+,\varphi_+)\to (I,\varphi)=\Psi(I_+,\varphi_+) 
\end{eqnarray*}
may be found, putting ${\rm H}$ into the form
\begin{eqnarray}\label{H+}
{\rm H}_+(I_+,\varphi_+):={\rm H}\circ \Psi_+(I_+,\varphi_+)={\rm h}_+(I_+)+{\rm f}_+(I_+,\varphi_+)\ ,
\end{eqnarray}
 where $\o_+:=\partial {\rm h_+}$  is a diffeomorphism of $\cI_{+\r_+}$ with Jacobian matrix $U_+:=\partial^2{\rm h_+}$ non singular on $\cI_{+\r+}$ such that $\o_+(\cI_+)=\o(\cI)$.
There also exist suitable constants
\beqano
&&  \dst M_+\geq\sup_{\cI_{+\r+}}\|U_+\|\label{M+}\\
&&  \dst \hat M_+\geq\sup_{\cI_{+\r+}}\|U_+^{[n,\hat n]}\|\\
&&  \dst N_+\geq\sup_{\cI_{+\r+}}\|T_+\|\\
&&  \dst \bar N_+\geq\sup_{\cI_{+\r+}}\|T_+^{[\bar n,n]}\|\\
&&  \dst \hat N_+\geq\sup_{\cI_{+\r+}}\|T_*^{[\hat n,n]}\|\\
&&  \dst F_+\geq\|{\rm f}_+\|_{\r_+,\bar s+s_+}\label{FF+}
\eeqano
where $T_+:=U_+^{-1}$, such that, defining
\beqano
&&  \dst K_+:=\frac{6}{s_+}\ \log_+{\left(\frac{F_+M_+^2\,L_+}{\gamma^2}\right)^{-1}}\\ 
&&  \dst \tilde{\r}_+:=\min\left\{\frac{\g}{3M_+K_+^{\t+1}}\ ,\ \frac{\hat\g}{3\hat M_+K_+^{\t+1}}\ ,\  
 \r_+\right\}\\
&&  \dst L_+:=\max\ \{N_+\ , M_+^{-1},\hat M_+^{-1}\}\\
\eeqano
then, 
\beq{KAM cond+}
E_+:=\frac{F_+L_+}{\tilde{\r}_+^2}\leq {E^2}\ .
\eeq
More precisely:
\begin{itemize}
\item[(i)] The numbers $\r_+$, $s_+$ and the constants $M_+$, $\cdots$ may be taken as
\beqano
\arr{\r_+=\frac{\tilde\r}{8}\\
s_+=\frac{s}{6}\\
M_+=2M\\	\hat M_+=2\hat M\\ N_+=2N\\ \bar N_+=2\bar N\\ \hat N_+=2\hat N\\ F_+=\frac{F^2LM^2}{\g^2}}\label{F+}
\eeqano
so as to satisfy
$$\frac{\tilde{\r}_+}{\tilde{\r}}\geq\frac{1}{2}\left(\frac{1}{12}\right)^{\t+1}\ ;$$
\item[(ii)] the set $\cI_+$ may be obtained as $\cI_+=l_+(\cI)$, where 
\item[(iii)]$\dst l_+:\ \cI\to \cI_+$ is an injective ``isofrequency map'', \ie, uniquely defined on $\cI$ by
$$\omega_+\circ l_+=\omega\quad \textrm{on}\quad \cI$$
which satisfies 
\begin{eqnarray}\label{lip}
{N}/{\bar N}\sup_{{\cal I}}|\bar l_+-\id|,\ {N}/{\hat N}\sup_{{\cal I}}|\hat l_+-\id|\leq 2^4E \tilde{\r}\ ,\quad {\cal L}_{\cal P}(l_+-\id)\leq2^8E\ ;
\end{eqnarray}
\item[(iv)] the map $\Psi$ satisfies
\begin{eqnarray}\label{Psi close to identity}
|W^{\bar n}_{\hat\g/\g,1}W_{\tilde \r,s}(\Psi_+-\textrm{\rm id})|_{\cal P}\leq E
\end{eqnarray}
and
\begin{eqnarray}\label{derivatives}
\left\|W^{\bar n}_{\hat\g/\g,1}(W_{\tilde \r,s}D\Psi_+\,W_{\tilde\r,s}^{-1}-I_{2p})
\right\|_{\cal P}\leq {2^5\,E}\ .
\end{eqnarray} 
\end{itemize}
\elem 
{\bf Proof.}\ We proceed by steps.
\vskip.1in
\noindent
{{$\underline{\textrm{\sl Claim $0$.}}$}:  {\it $\rm H$ is put into the form
$${\rm H}_+(I_+,\varphi_+):={\rm H}\circ \Psi_{\rm av}(I_+,\varphi_+)={\rm h}(I_+)+{\rm g}(I_+)+{\rm f}_+(I_+,\varphi_+):={\rm h}_+(I_+)+{\rm f}_+(I_+,\varphi_+)$$
where
\beq{g minus f0}
\sup_{\cI_{\tilde\r/2}}|{\rm g}-\langle {\rm f}\rangle_{\varphi}|\leq \frac{F}{16}
\eeq
and 
\beq{new F}
\|{\rm f}_+\|_{\tilde\r/2,\bar s+s/6}\leq \frac{F^2M^2L}{\g^2}=:F_+\ ,
\eeq
by means of a real--analytic symplectomorphism $\Psi_{\rm av}$ defined on
$\cI_{\tilde\r/2}\times\torus^{n}_{\bar s+s/6}$ and verifying
\beqa{bounds for Psiav}
&& |W^{n}_{\hat\g/\g,1}W_{r,s}(\Psi_{\rm av}-I)|_{\cal P}\leq E\ ,\nonumber\\
&& \|W^{n}_{\hat\g/\g,1}(W_{r,s}D\Psi_{\rm av} W^{-1}_{r,s}-I)\|_{\cal P}
\leq 2^5E\ .
\eeqa
}

\vskip.1in
\noi
{{$\underline{\textrm{\sl Proof.}}$} We are going to construct the transformation $\Psi_{\rm av}$ by means of application of the Averaging Theorem (Lemma \ref{average}) to ${\rm H}$, on the domain $\cI_{\tilde\r}\times \torus^n_{\bar s+s/6}$, with the trivial resonant lattice $\L=\{0\}\in \integer^n$. We first verify the ``non resonance'' assumption (\ref{non resonance in avarage thm}) out of $\L=\{0\}$.
By assumption, for any $I\in \cI$, $\o(I)\in D^{\bar n,\hat n}_{\g,\hat\g,\t}$, which means
\beqano
|\o(I)\cdot k|\geq \arr{
\frac{\g}{|k|^\t}\quad \textrm{if}\quad k=(\bar k,\hat k)\quad \textrm{with}\quad \bar k\neq 0\\
\frac{\hat\g}{|k|^\t}\quad \textrm{if}\quad k=(0,\hat k)\quad \textrm{with}\quad \hat k\neq 0\\
}
\eeqano
for any $k\in \integer^{\bar n}\times \integer^{\hat n}\setminus\{0\}$. In particular, for $\dst 0<|k|\leq K$ and $I\in \cI$,
\beqano
|\o(I)\cdot k|\geq \arr{
\frac{\g}{K^\t}\quad \textrm{if}\quad k=(\bar k,\hat k)\quad \textrm{with}\quad \bar k\neq 0\\
\frac{\hat\g}{K^\t}\quad \textrm{if}\quad k=(0,\hat k)\quad \textrm{with}\quad \hat k\neq 0\\
}
\eeqano
Let, now, $\dst I\in \cI_{\tilde\r}$. By definition, there exist $I_0\in \cI$ such that $|I-I_0|<\tilde\r$ and we find (recall that we have prefixed  the $1$--norm in $\cI$)
\beqano
|\o(I)-\o(I_0)|_{\infty}&:=&\max_{1\leq i\leq p}|\o_i(I)-\o_i(I_0)|\nonumber\\
&\leq&\sup_{\cI_\r}\|\partial\o\||I-I_0|\nonumber\\
&\leq&M\tilde\r
\eeqano 
Hence, if $\dst k=(\bar k,\hat k)$, with $\bar k\neq 0$ and $0<|k|\leq K$
\beqano
|\omega(I) \cdot k|&\geq&|\o(I_0)\cdot k|-|(\o(I)-\o(I_0))\cdot k|\nonumber\\
&\geq& \frac{\g}{K^\t}-|\o(I)-\o(I_0)|_{\infty}|k|\nonumber\\
&\geq& \frac{\g}{K^\t}-MK\tilde \r\nonumber\\
&\geq&\frac{2}{3}\,\frac{\g}{K^\t}=:\bar\a
\eeqano
having used $\dst\tilde\r\leq \frac{\g}{3MK^{\t+1}}$. Similarly, taking $k=(0,\hat k)$, with $0<|k|=|\hat k|\leq K$ and using $\dst\tilde\r\leq \frac{\hat\g}{3\hat MK^{\t+1}}$, we find that
$$|\omega(I) \cdot k|\geq \frac{2}{3}\,\frac{\hat\g}{K^{\t}}=:\hat\a\quad \textrm{for any}\quad I\in \cI_{\tilde\r}\ ,\quad 0<|k|\leq K\ .$$
But, then, the ``smallness condition'' (\ref{tilde rho alpha}) with $r=\tilde\r$ is also verified in $\cI_{\tilde\r}$, because, using
$$\a=\min\left\{\frac{2}{3}\,\frac{\g}{K^\t},\ \frac{2}{3}\,\frac{\hat\g}{K^\t}\right\}\geq 2K\min\{M,\hat M\}\tilde\r\geq 2K\frac{\tilde\r}{L}\ ,$$
we find that (\ref{KAM condnew}) is stronger:
$$2^7\frac{\|f\|_{\tilde\r,\bar s+s}K}{\a\tilde\r}\leq 2^6 \frac{FL}{\tilde\r^2}<1\ .$$
Therefore, Lemma \ref{average} applies (as also, trivially, $Ks=6\log_+(FM^2L/\g^2)\geq 6$), and $\rm H$ is put into the form
$${\rm H}_+:={\rm H}\circ \Psi_{\rm av}={\rm h}+{\rm g}+{\rm f}_+:={\rm h}_++{\rm f}_+$$
by means of a real--analytic symplectomorphism $\Psi_{\rm av}$ defined on
$\cI_{\tilde\r/2}\times\torus^{n}_{\bar s+s/6}$,
where ${\rm g}$ is a $\{0\}$--completely resonant  (which means that $\rm g$ is a function of $I$ only) real--analytic function, suitably close to ${\rm f}_0=P_{\{0\}}T_K\,{\rm f}=\langle f\rangle_{\varphi}$:
\beqano
\sup_{\cI_{\tilde\r/2}}|{\rm g}-{\rm f}_0|=\|{\rm g}-{\rm f}_0\|_{\tilde\r/2,\bar s+s/6}\leq \frac{2^5K}{\a\,\tilde\r}\|{\rm f}\|_{\tilde\r,\bar s+s}^2\leq 2^4\frac{F^2L}{\tilde\r^2}\leq \frac{F}{16}
\eeqano
and ${\rm f}_+$ is ``small'':
\beqano
\|{\rm f}_+\|_{\tilde\r/2,\bar s+s/6}\leq e^{-Ks/6}\|{\rm f}\|_{\tilde\r,\bar s+s}\leq e^{-Ks/6}F=\frac{F^2M^2L}{\g^2}\ ,
\eeqano
The bounds (\ref{bounds for Psiav}) are an easy consequence of (\ref{identity bounds}), (\ref{derivatives0}) (recall $\hat\g<\g$):
\begin{eqnarray*}
|W^{\bar n}_{\hat\g/\g,1}W_{r,s}(\Psi_{\rm av}-I)|_{\cal P}=|W^{\bar n}_{\a/\bar\a,1}W_{r,s}(\Psi_{\rm av}-I)|_{\cal P}\leq \frac{2\,K\|{\rm f}\|_{\tilde\r,\bar s+s}}{\a\,\tilde\r}\leq \frac{FL}{\tilde\r^2}=E\ ,
\end{eqnarray*}
and
\begin{eqnarray*}
\|W^{\bar n}_{\hat\g/\g,1}(W_{r,s}D\Psi_{\rm av} W^{-1}_{r,s}-I)\|_{\cal P}&=&\|W^{\bar n}_{\a/\bar\a,1}(W_{r,s}D\Psi_{\rm av} W^{-1}_{r,s}-I)\|_{\cal P}\nonumber\\
\leq \frac{2^6K}{\a\,\tilde\r}\,\|{\rm f}\|_{r,\bar s+s}
\leq 2^5E\ ,
\end{eqnarray*} {\sl Claim} $0$ is thus proved.

\vskip.1in
\noindent
{{$\underline{\textrm{\sl Claim $1$}}$}: {\it the Jacobian matrix $U_+:=\partial^2{\rm h}_+$ is non singular in $\cI_{\tilde\r/4}$ and satisfies
\beqa{bounds on M}
&&  \dst M_+:=2M\geq\sup_{I\in{\cal I}_{\r/4}}\|U_+\|\\
&&  \dst \hat M_+:=2\hat M\geq\sup_{I\in{\cal I}_{\r/4}}\|U_+^{[n,\hat n]}\|\\
&&  \dst N_+:=2N\geq\sup_{I\in{\cal I}_{\r/4}}\|T_+\|\\
&&  \dst \bar N_+:=2\bar N\geq\sup_{I\in{\cal I}_{\r/4}}\|T_+^{[\bar n,n]}\|\\
&&  \dst \hat N_+:=2\hat N\geq\sup_{I\in{\cal I}_{\r/4}}\|T_{+*}^{[\hat n,n]}\|\label{bounds on N}
\eeqa
where $T_+:=U_+^{-1}$.}
\vskip.1in
\noindent
{$\underline{\textrm{\sl Proof.}}:$} The bound (\ref{g minus f0}) gives 
\begin{eqnarray}\label{newhamlemma2}
\sup_{{\cal I}_{\tilde{\rho}/2}}|g|\leq \sup_{{\cal I}_{\tilde{\rho}/2}}|f_0|+\sup_{{\cal I}_{\tilde{\rho}/2}}|g-f_0|\leq \frac{17}{16}\,F\ ,
\end{eqnarray}
whence, applying twice the P\"oschel's General Cauchy Inequality (see Appendix \ref{The General Cauchy Inequality}),
\begin{eqnarray}\label{estimates on g}
\sup_{{\cal I}_{\tilde{\r}/4}}\|\partial^2\,g\|\leq \frac{\sup_{{\cal I}_{3\tilde{\r}/8}}\|\partial g\|}{\tilde{\r}/8}\leq\frac{17/16\,F}{(\tilde{\r}/8)^2}\leq 2^7\frac{F}{\tilde{\r}^2}<\frac{1}{2L}\leq \frac{1}{2}\min\left\{M, \hat M,\frac{1}{N}\right\}
\end{eqnarray}
which suddenly implies 
$$\sup_{{\cI}_{\tilde\r/4}}\|U_+\|=\sup_{{\cI}_{\tilde\r/4}}\|\partial^2{\rm h}_+\|=\sup_{{\cI}_{\tilde\r/4}}\|\partial^2{\rm h}+\partial^2{\rm g}\|\leq \sup_{{\cI}_{\tilde\r/4}}\|\partial^2{\rm h}\|+\sup_{{\cI}_{\tilde\r/4}}\|\partial^2{\rm g}\|\leq M+M=2M$$
Similarly, one finds
$$\sup_{{\cI}_{\tilde\r/4}}\|U_+^{[n,\hat n]}\|\leq 2\hat M$$
But (\ref{estimates on g}) also implies
\begin{eqnarray}\label{non degeneracy}
\sup_{{\cal I}_{\tilde{\r}/4}}\|\partial^2\,{\rm g}(\partial^2\,{\rm h})^{-1}\|
&\leq& \sup_{{\cal I}_{\tilde{\r}/4}}\|\partial^2g\|\sup_{{\cal I}_{\tilde{\r}/4}}\|(\partial^2\,h)^{-1}\|\nonumber\\
&\leq& \sup_{{\cal I}_{\tilde{\r}/4}}\|\partial^2\,g\|N\nonumber\\
&\leq&\frac{1}{2}
\end{eqnarray}
so, the matrix $$\id+\partial^2\,{\rm g}(\partial^2\,{\rm h})^{-1}$$ is non singular on ${\cal I}_{\tilde{\r}/4}$ with
$$\left\|\Big(\id+\partial^2\,{\rm g}(\partial^2\,{\rm h})^{-1}\Big)^{-1}\right\|\leq 2\ .$$
This implies that
$$\partial^2\,h_+=\partial^2{\rm h}+\partial^2{\rm g}=\Big(\id+\partial^2\,{\rm g}(\partial^2\,{\rm h})^{-1}\Big)\partial^2\,{\rm h}$$
is non singular on ${\cal I}_{\tilde{\r}/4}$, with
\begin{eqnarray}
\sup_{{\cal I}_{\tilde{\r}/4}}\|T_+\|=\sup_{{\cal I}_{\tilde{\r}/4}}\|(\partial^2\,h_+)^{-1}\|&=&\sup_{{\cal I}_{\tilde{\r}/4}}\left\|(\partial^2{\rm h})^{-1}\Big(\id+\partial^2\,{\rm g}(\partial^2\,{\rm h})^{-1}\Big)^{-1}\right\|\leq 2N\ .
\end{eqnarray}
Similarly,
\begin{eqnarray*}
\sup_{{\cal I}_{\tilde{\r}/4}}\|T^{[\bar n,n]}\|
&\leq& 2\,\bar N\ ,\qquad \sup_{{\cal I}_{\tilde{\r}/4}}\|T_*^{[\hat n,n]}\|
\leq 2\,\hat N\ ,
\end{eqnarray*}

\vskip.1in
\noindent
{{$\underline{\textrm{\sl Claim $2$}}$}:{\it The new frequency $\o_+:=\partial{\rm h}_+$ is a diffeomorphism of $\cI_{3\tilde\r/16}$}.} 
\vskip.1in
\noindent
{\sl $\underline{\textrm{Proof.}}$} We want to prove that, if $I_+$, $I_+'\in \cI_{3\tilde\r/16}$ verify $\o_+(I_+)=\o_+(I_+')$, then, $I_+=I_+'$.\\
Let $I_+$, $I_+'\in \cI_{\tilde\r/8}$ verify
$$\o(I'_+)+\partial{\rm g}(I_+')=\o_+(I_+)=\o_+(I_+')=\o(I_+)+\partial{\rm g}(I_+)\ .$$
Then,
$$|\o(I'_+)-\o(I_+)|=|\partial{\rm g}(I_+')-\partial{\rm g}(I_+)|\leq 2\sup_{\cI_{\tilde\r/4}}|\partial {\rm g}|\leq 2^4\frac{F}{\tilde\r}\ ,$$
hence,
$$|I_+-I_+'|=|\o^{-1}(\o(I_+))-\o^{-1}(\o(I_+'))|\leq \sup\|\partial(\o^{-1})\||\o(I'_+)-\o(I_+)|\leq 2^4n\frac{FN}{\tilde\r}\leq\frac{\tilde\r}{16}$$
The previous inequality implies that the segment $s(I_+,I_+')$ from $I_+$ to $\cI_+'$ lies interely in $\cI_{\tilde\r/4}$. So, let $\t$ the curve from $\o(I_+)$ to $\o(I_+')$ defined as $\t:=\o(s(I_+,I_+'))$; let $F:=\partial{\rm g}\circ\o^{-1}$ and observe $\sup_{\o(\cI_{\tilde\r/4})}\|\partial F\|=\sup_{\cI_{\tilde\r/4}}\|\partial^2{\rm g}(\partial {\rm h})^{-1}\|\leq 1/2$. Then,
\beqano
0&=&|\o_+(I_+)-\o_+(I_+')|\nonumber\\
&=&|\o(I_+)-\o(I_+')+F(\o(I_+))-F(\o(I_+'))|\nonumber\\
&\geq&|\o(I_+)-\o(I_+')|-\left|\int_{\t}\partial F(\zeta)\cdot d\zeta\right|\nonumber\\
&\geq&\frac{1}{2}|\o(I_+)-\o(I_+')|
\eeqano
which implies $\o(I_+)=\o(I_+')$, hence, $I_+=I_+'$.

\vskip.1in
\noindent
{{$\underline{\textrm{\sl Claim $3$}}$}: {\it The new frequency $\o_+$ maps $\cI_{\tilde\r/32}$ over $\o(\cI)$, \ie, $\o_+(\cI_{\tilde\r/32})\supseteq \o(\cI)$}.
\vskip.1in
\noindent
{{$\underline{\textrm{\sl Proof.}}$}: we prove that, for any $I_0\in \cI$, $\o_+(B^n_{\tilde\r/32}(I_0))\supseteq \o(I_0)$. If $|I-I_0|=r<\tilde\r/32$, then,
\beqano
|\o_+(I)-\o_+(I_0)|&=&|\o(I)-\o(I_0)+\partial{\rm g}\circ\o^{-1}(\o(I))-\partial{\rm g}\circ\o^{-1}(\o(I_0))|\nonumber\\
&\geq& \left(1-\|\partial^2{\rm g}(\partial{\rm h})^{-1}\|\right)|\o(I)-\o(I_0)|\nonumber\\
&\geq&\frac{1}{2}|\o(I)-\o(I_0)|\nonumber\\
&\geq&\frac{1}{2N}|I-I_0|\nonumber\\
&=&\frac{r}{2N}
\eeqano
hence, $\o_+(B^n_{\tilde\r/32}(I_0))\supseteq B^n_{\tilde\r/(64N)}(\o_+(I_0))$. We prove that
$$\o(I_0)\in B^n_{r/(64N)}(\o_+(I_0))\ .$$ Using the KAM condition (\ref{KAM cond}) and General Cauchy Estimate (for $\sup_{\cI}|\partial{\rm g}|$)
$$|\o_+(I_0)-\o(I_0)|=|\partial{\rm g}(I_0)|\leq \sup_{\cI}|\partial{\rm g}|<\frac{17}{16}\frac{F}{\tilde\r/2}<\frac{\tilde\r}{2^6N}$$
which concludes the proof of the Claim.

\vskip.1in
\noindent
{{$\underline{\textrm{\sl Claim $4$}}$} {\it For any $I\in \cI$, equation $\o_+(I')=\o(I)$ has a unique solution $I':=l_+(I)=\o_+^{-1}\circ \o(I)\in \cI_{\tilde\r/32}$ satisfying (\ref{lip})
}.
\vskip.1in
\noindent
$\underline{\textrm{\sl Proof.}}$\ Existence and uniqueness of  the solution $I'=l_+(I)$ of $\o_+(I')=\o(I)$  are consequences of {\sl claims} $2$, $3$. We prove (\ref{lip}). Let $0<r<\bar r\leq3\tilde\r/16$, with $\bar r$ so small that $\o(\cI_r)\subseteq \o_+(3\cI_{\tilde\r/16})$. For $I\in \cI_r$, we find, as $\o_+$ is a diffeomorphism of $\cI_r$ and General Cauchy Inequality,
\begin{eqnarray*}
|\ovl{\o_+^{-1}}(\o(I))-\ovl{\o_+^{-1}}(\o_+(I))|\leq 2\bar N|\o_+(I)-\o(I)|\leq2\bar N\sup_{{\cI}_r}|\partial{\rm g}|\leq\frac{2\bar N(17F/16)}{3\tilde\r/16-r}
\end{eqnarray*}
Hence, due to the arbitrariness of $r$, for $I\in \cI=\bigcap_{0<r<\bar r}\cI_r$,
\begin{eqnarray*}
|\bar l_+(I)-\bar I|\leq\frac{34\bar NF}{3\tilde\r}\leq 2^4\frac{\bar NF}{\tilde\r}
\end{eqnarray*}
Similarly,
\begin{eqnarray*}
|\hat l_+(I)-\hat I|\leq2^4 \frac{\hat NF}{\tilde\r}
\end{eqnarray*}
Also, using
$$Dl_+=D[(\o+\partial g)^{-1}\circ \o]=[\textrm{id}_n+(\partial^2{\rm h})^{-1}\partial^2{\rm g}]^{-1}$$
we find
\beqano
\sup_{\cI_{r}}\|D[(\o+\partial g)^{-1}\circ \o]-\textrm{id}_n\|&=&\sup_{\cI_r}\|[\textrm{id}_n+(\partial^2{\rm h})^{-1}\partial^2{\rm g}]^{-1}-\textrm{id}_n\|\nonumber\\
&\leq&\frac{\sup_{\cI_r}\|(\partial^2{\rm h})^{-1}\partial^2{\rm g}\|}{1-\sup_{\cI_r}\|(\partial^2{\rm h})^{-1}\partial^2{\rm g}\|}\nonumber\\
&\leq&2\sup_{\cI}\|(\partial^2{\rm h})^{-1}\partial^2{\rm g}\|\nonumber\\
&\leq& 2^8E\ .
\eeqano
which says
$${\cal L}_{\cP}(\o_+^{-1}\circ\o-\id)\leq 2^{8}E\quad \textrm{on}\quad \cI_r$$
and implies, due to the arbitrariness of $r$,
$${\cal L}_{\cP}(l_+-\id)\leq 2^{8}E\ .$$

\vskip.1in
\noindent
{{$\underline{\textrm{\sl Conclusion.}}$}: Let $\cI_+:=l_+(\cI)$, $\r_+:=\tilde\r/8$, $s_+:=s/6$. By {\sl claim $4$}, $\cI_+\subset \cI_{\tilde\r/32}$, so, the following inclusions hold:
$$\cI_{+\r_+}\subset \cI_{5\tilde\r/32}\subset \cI_{3\tilde\r/16}\subset\cI_{\tilde\r/4}$$
Hence, $\Psi_+:=\Psi_{\rm av}|_{\cI_{+\r_+}\times \torus^n_{\bar s+s_+}}$ is well defined and has the desired properties, upon recognizing that
\beqano
E_+:=\frac{F_+M_+^2L_+}{\tilde\r_+^2}\leq E^2\ ,\quad \tilde\r_+\geq \frac{1}{2}\left(\frac{1}{12}\right)^{\t+1}\tilde \r\ .
\eeqano
In order to do that, we first prove
\begin{eqnarray}\label{new K}10\,K\leq K_+\leq 12\,K\ .
\end{eqnarray}
We find
$$\frac{F\,LM^2}{\g^2}\leq \frac{F\, LM^2}{\g^2}\left[\log_+{\left(\frac{F\, LM^2}{\g^2}\right)^{-1}}\right]^{2(\t+1)}\leq \frac{1}{9}\left(\frac{s}{6}\right)^{2(\t+1)}\frac{F\,L}{\tilde \r^2}=\leq\frac{2^8\,E}{9\cdot6^2\cdot2^8}\leq\frac{1}{9\cdot6^2\cdot2^8}$$
which gives
\begin{eqnarray}\label{estimates on x}
2^{16}< x:=\left(\frac{F LM^2}{\g^2}\right)^{-1}\ ,
\end{eqnarray}
namely
$$\frac{108}{s}\log{2}<\frac{9}{8}\,K<2\,K$$
implying immediately (\ref{new K}), after using
$$K_+=\frac{6}{s_+}\log_+{(x^2/8)}=\frac{6}{s_+}\log{(x^2/8)}=12\,K-\frac{108}{s}\log{2}\ .$$ 
Now, using
$$\g\geq 3\,MK^{\t+1}\tilde \r\ ,\quad \hat \g\geq 3\,MK^{\t+1}\tilde \r\ ,\quad \frac{K_+}{K}\leq 12\ ,\quad K\geq \frac{96\log{2}}{s}\ ,\quad 0<s\leq 1\ ,$$
we find
\begin{eqnarray*}E_+&\leq&\max\left\{128\,\left(\frac{F\,LM}{\g\tilde \r}\right)^2\ ,\quad 648\left(\frac{FLM}{\g^2}\right)^2K_+^{2(\t+1)}\ ,\quad 648\left(\frac{FLM\hat M}{\g\hat\g}\right)^2K_+^{2(\t+1)}\right\}\nonumber\\
&\leq&\max\left\{\frac{128}{9K^{2(\t+1)}}\,\left(\frac{F\,L}{\tilde \r^2}\right)^2\ ,\quad 8\left(\frac{FL}{\tilde \r^2}\right)^2\frac{K_+^{2(\t+1)}}{K^{4(\t+1)}}\ ,\quad 8\left(\frac{FL}{\tilde \r^2}\right)^2\frac{K_+^{2(\t+1)}}{K^{2(\t+1)}K^{2(\t+1)}}\right\}\nonumber\\
&\leq&\max\left\{\frac{128}{9}\left(\frac{s}{96\log{2}}\right)^{2(\t+1)}\ ,\quad 8\left(\frac{s}{8\log{2}}\right)^{2(\t+1)}\right\}\,E^2\nonumber\\
&\leq& E^2\ .
\end{eqnarray*}
Similarly,
\begin{eqnarray*}
\tilde{\r}_{+}&=&\min \left\{\r_{+}\ ,\ \frac{\g}{3M_+K_{+}^{\t+1}}\ ,\ \frac{\hat\g}{3\hat M_+K_{+}^{\t+1}}\right\}\nonumber\\
&=&\min \left\{\frac{\r}{8}\ ,\ \frac{\g}{24MK^{\t+1}}\ ,\ \frac{\hat\g}{24MK^{\t+1}}\ ,\ \frac{\g}{6MK_{+}^{\t+1}}\ ,\ \frac{\hat\g}{6MK_{+}^{\t+1}}\right\}\nonumber\\
&=&\min \left\{\frac{\r}{8}\ ,\ \frac{\g}{6MK_{+}^{\t+1}}\ ,\ \frac{\hat\g}{6MK_{+}^{\t+1}}\right\}\nonumber\\
&\geq&\frac{1}{2}\left(\frac{1}{12}\right)^{\t+1}\min\left\{\r\ ,\ \frac{\g}{3MK^{\t+1}}\ ,\quad \frac{\hat\g}{3MK^{\t+1}}\right\}\nonumber\\
&=&\frac{1}{2}\left(\frac{1}{12}\right)^{\t+1}\tilde \r\nonumber\\
\end{eqnarray*}
(use $(M_+,\hat M_+)=2\,(M,\hat M)$, $12\,K\geq K_+\geq 10\,K$).
This concludes the proof. 
\subsubsection{Lemmas on Measure}
We recall the following classical results on Lipschitz functions and measure theory, referring to \cite{EG} for their proofs. 
\vskip.1in
\noindent
\begin{lemma}\label{extension}{\bf (Kirszbraun Theorem)}\ Assume ${A}\subset \real^n$, and let $f:A\to \real^m$ be Lipschitz. There exists a Lipschitz function $\bar f:\ \real^n\to \real^m$, such that
\vskip.1in
\noindent
$i)\quad\bar f=f\quad {\rm on}\quad A\ ;$
\vskip.1in
\noindent
$ii)\quad{\cal L}({\bar f})={\cal L}(f)\ .$ 
\end{lemma}
\begin{lemma}\label{measure}
Let $A\subset \real^n$ Lebesgue--measurable, $f:\ A\to \real^n$ Lipschitz (bi--Lipschitz). Then, 
\begin{eqnarray*}
\meas((f(A))&\leq& {\cal L}(f)^n\ \meas(A)\nonumber\\
{\cal L}_{-}(f)^n\ \meas(A)&\leq&\meas((f(A))\leq {\cal L}_{+}(f)^n\ \meas(A))\ .
\end{eqnarray*}
\end{lemma}
\vskip.1in
\noindent
\subsubsection{Proof of Theorem \ref{two scales KAM}}
Here also, we proceed by steps.
$\underline{\textrm{\sl Claim}\ 0.}$ {\sl (``construction of the sequences'')}: {\it For each $1\leq j\in \natural$, $\rm H$ is analitically conjugated to
$${\rm H}_j={\rm H}\circ \Phi_j={\rm h}_j+{\rm f}_j$$ real--analytic on $\cP_j=\cI_{j\r_j}\times\torus^n_{\bar s+s_j}$, where: 
\begin{itemize}
\item[(i)]$s_j=s/6^j$ and, letting
\beqano
\arr{M_j=2^jM\\
\hat M_j=2^j\hat M\\
N_j=2^j N\\
L_j=\max\left\{M_j^{-1},\ \hat M_j^{-1},\ N_j\right\}\\
F_j=\frac{F_{j-1}^2L_{j-1}M_j^2}{\g^2}\\
K_j:=\frac{6}{s_j}\log_+\left(\frac{F_{j}L_{j}M_j^2}{\g^2}\right)\\
\tilde\r_j:=\min\left\{\frac{\g}{3M_jK_j^{\t+1}},\ \frac{\hat\g}{3\hat M_jK_j^{\t+1}},\ \r_{j}\right\}}
\eeqano
then, $\r_j=\tilde\r_{j-1}/8$; 
\item[(ii)] $\cI_j\subseteq {\cI_{j-1}}_{\tilde\r_{j-1}/32}$ is obtained as $\cI_j=l_j(\cI_{j-1})$, where $l_j$ is a Lipschitz homeomorphism satisfying
\begin{eqnarray}\label{estimates for lj}
&&\o_j\circ l_j=\o_{j-1}\quad \textrm{on}\quad \cI_{j-1}\nonumber\\
&&\max\left\{\frac{N}{\bar N}\sup_{{\cal I}_j}|\bar l_{j+1}-\id_{\bar n}|,\ \frac{N}{\hat N}\sup_{{\cal I}_j}|\hat l_{j+1}-\id_{\hat n}|\right\}\leq 2^4\tilde{\r}_jE^{2^j}\\\label{estimates for lj1}
&&{\cal L}_{\cP}(l_{j+1}-\id)\leq 2^8E^{2^j}:=\m_j
\end{eqnarray}
\item[(iii)] $\o_j:=\partial{\rm h}_j$ is a diffeomorphism of $\cI_{j\r_j}$ with non singular Jacobian $U_j:=\partial^2{\rm h}_j$ such that $\o_j(\cI_j)=\cD^{\bar n,\hat n}_{\g,\hat\g,\t}\cap\o(\cI)$;
\item[(iv)] ${\rm f}_j$ satisfies
$$\|{\rm f}_j\|_{\cI_{j\r_j},\bar s+s_j}\leq F_j$$
\item[(v)] The real--analytic symplectomorphism $\Phi_j$ is obtained as $\Phi_j=\Psi_1\circ\cdots\circ\Psi_j$, where $\Psi_k:\cP_k\to\cP_{k-1}$ ($k\geq 1$) verifies
\begin{eqnarray}\label{Psij close to id}
&&\sup_{\cP_k}|W^{\bar n}_{\hat\g/\g,1}W_{\tilde\r,s}(\Psi_{k}-\textrm{\rm id})|_{\cal P}\leq \left(\frac{1}{6}\right)^{k-1}E^{2^{k-1}}\\
\label{deriv1}
&&\sup_{\cP_k}\left\|W_{\tilde\r, s}D\Psi_{k}\,W^{-1}_{\tilde\r, s}-\id_{2n}
\right\|_{\cal P}\leq2^5\,\left[\frac{(12)^{\t+1}}{3}\right]^{{k-1}}{E^{2^{k-1}}}=:\zeta_k\ .
\end{eqnarray}
\end{itemize}}
\noi
{\sl Proof.}\ Starting with 
$${\rm H}_{0}:={\rm H}={\rm h}+{\rm f}\ ,\quad \cP_0:={\cI}_{0\r}\times \torus^n_{\bar s+s}\ ,$$
where $\cI_0:=\{I\in\cI:\ \o(I)\in \cD^{\bar n,\hat n}_{\g,\hat\g,\t}\}$
and, labeling by ``$0$'' the quantities relatively to ${\rm H}_0$, apply (inductively) the Iterative Lemma (Lemma \ref{itera}) to
$${\rm H}_{j}={\rm h}_j+{\rm f}_j\ ,\quad \cP_j:={\cI}_{j\r_j}\times \torus^n_{\bar s+s_j}\ ,\quad j\geq 0$$
and label by ``$j+1$'' the ``$+$''--quantities of the Iterative Lemma.
Next, observe  that (\ref{deriv1}) is a consequence of
$$\sup_{\cP_k}|W^{\bar n}_{\hat\g/\g,1}W_{k-1}(\Psi_{k}-\textrm{\rm id})|_{\cal P}\leq E^{2^{k-1}}$$
$$\sup_{\cP_k}\left\|W_{k-1}D\Psi_{k}\,W_{k-1}^{-1}-\id_{2n}
\right\|_{\cal P}\leq 2^5E^{2^{k-1}}$$
where $W_{k-1}:=W_{\tilde\r_{k-1},s_{k-1}}$,
using
$$\|W_{\tilde\r, s}W_{k-1}^{-1}\|_{\cal P}=\max\left\{\frac{\tilde{\r}_{k-1}}{\tilde{\r}}\ ,\ \frac{s_{k-1}}{s}\right\}=\left(\frac{1}{6}\right)^{k-1}$$ $$\|W_{k-1}W_{\tilde\r, s}^{-1}\|_{\cal P}=\max\left\{\frac{\tilde{\r}}{\tilde{\r}_{k-1}}\ ,\ \frac{s}{s_{k-1}}\right\}\leq [2(12)^{\t+1}]^{k-1}\ .$$
 
\vskip.1in
\noindent
$\underline{\textrm{\sl Claim}\ 1.}$ {\sl(``construction of $\cI_*$'')}: {\it The sequence of Lipschitz homeomorphisms on ${\cal I}_0$
\begin{eqnarray}\label{ellj}
\ell_j:=l_j\circ l_{j-1}\circ \cdots \circ l_1
\end{eqnarray}
converges uniformly to a bi--Lipschitz homeomorphism $\ell=(\bar \ell,\hat\ell)$ satisfying
\begin{eqnarray}\label{ell close to id1}
\frac{N}{\bar N}\sup_{\cI_0}|\bar\ell-\id|\ ,\quad \frac{N}{\hat N}\sup_{\cI_0}|\hat\ell-\id|\leq {2^5}\tilde{\r}E
\end{eqnarray}
\beq{Lip const for ell}
{\cal L}_{\cP}(\ell-\id)\leq 2^{10}E
\eeq
Furthermore, the following holds:
\beqa{ellj-ell}
\sup_{\cI_0}|\ell_j-\ell|\leq 2^5\tilde\r E^{2^j}
\eeqa
and 
$$\cI_*:=\ell(\cI_0)\subseteq {\cI_{0}}_{32E\tilde\r}\bigcap\Big(\cap_j{\cal I}_{j\r_j}\Big)$$}
\vskip.1\in
\noindent
{\sl Proof.}\ Using (\ref{estimates for lj}), the inequality
\begin{eqnarray}\label{iterations}
\sup_{\cI_0}|\ell_i-\ell_j|&\leq& \sum_{l=j}^{i-1}\sup_{\cI_0}|\ell_{l+1}-\ell_l|\leq 2^4\tilde{\r}_j\,\sum_{l=j}^{i-1}\left(\frac{1}{8}\right)^{l-j}E^{2^l}
\end{eqnarray}
proves the uniform convergence of $\ell_j$. 
Letting, in (\ref{iterations}), $i\to +\infty$, we find \footnote{$2^l-2^j\geq l-j$ for any $l\geq j$} 
\begin{eqnarray}\label{ellj - ell}
\sup_{\cI_0}|\ell_j-\ell|
&\leq& 2^4\tilde{\r}_j\,\sum_{l=j}^{+\infty}\left(\frac{1}{8}\right)^{l-j}E^{2^l}\nonumber\\
&\leq& 2^4\tilde{\r}_jE^{2^j}\,\sum_{l=j}^{+\infty}\left(\frac{1}{8}\right)^{l-j}E^{2^l-2^j}\nonumber\\
&\leq&2^4\tilde{\r}_jE^{2^j}\sum_0^{+\infty}\left(\frac{E}{8}\right)^k\nonumber\\
&\leq&2^5{\r_j}E^{2^j}\nonumber\\
&<&\r_j
\end{eqnarray}
which also holds for $j=0$, with the convention $\ell_0:=\id$. Eq. (\ref{ellj - ell}) implies
$${\cal I}_*\subseteq \bigcap_j{\cal I}_{j\r_j}\ .$$
In particular, (\ref{ellj - ell}) with $j=0$ gives
\begin{eqnarray}\label{ell close to id}
\sup_{\cI}|\ell-\textrm{\rm id}|\leq {2^5}\tilde{\r}E\ .
\end{eqnarray}
which also implies
\begin{eqnarray}\label{ell(I0)}
{\cal I}_*\subset {{\cal I}_0}_{{32}\tilde{\r}E}\ .
\end{eqnarray}
With similar techniques, but using (\ref{ell close to id}), one proves (\ref{ell close to id1}). We prove that $\ell$ is injective on $\cI_0$. If $I$, $I'\in {\cal I}_0$ are such that $\ell(I)=\ell(I')=I_*$, then,
by (\ref{ellj - ell})
\begin{eqnarray*}
|\o(I)-\o(I')|&=&|\o_j(\ell_j(I))-\o_j(\ell_j(I'))|\nonumber\\
&\leq&M_j\,|\ell_j(I)-\ell_j(I')|\nonumber\\
&\leq&M_j\,\left(|\ell_j(I)-\ell(I)|+|\ell_j(I')-\ell(I')|\right)\nonumber\\
&\leq&2M\,2^j{\r}_jE^j
\end{eqnarray*}
which gives $\o(I)=\o(I')$ (as  the r.h.s goes to $0$ as $j\to \infty$), hence, $I=I'$. We prove (\ref{Lip const for ell}). The estimates (\ref{estimates for lj1}) for ${\cal L}_{\cal P}(l_{j+1}-\id)$ give
\footnote{Write 
$$i_{j+1}:=\ell_{j+1}-\id=(l_{j+1}-\id)\circ(\id+i_j)+i_j$$
to find
\begin{eqnarray*}
{\cal L}(i_{j+1})&\leq& {\cal L}(l_{j+1}-\id)\Big(1+{\cal L}(i_j)\Big)+{\cal L}(i_j)\nonumber\\
&=&{\cal L}(i_j)\Big({\cal L}(l_{j+1}-\id)+1\Big)+{\cal L}(l_{j+1}-\id)
\end{eqnarray*}
Iterating the above formula, we find
\begin{eqnarray*}
{\cal L}(i_{j+1})&\leq& {\cal L}(l_{j+1}-\id)+(1+{\cal L}(l_{j+1}-\id)){\cal L}(l_{j}-\id)+\cdots\nonumber\\
&+&(1+{\cal L}(l_{j+1}-\id))\cdots(1+{\cal L}(l_{2}-\id)){\cal L}(l_{1}-\id)\nonumber\\
&=&\prod_{k=1}^{j+1}(1+{\cal L}(l_{k}-\id))-1\ .
\end{eqnarray*}
}
\beqa{Lip for ellj}{\cal L}_{\cal P}(\ell_{j+1}-\textrm{id})\leq \prod_{l=1}^{j}(1+\mu_l)-1\leq \prod_{l=1}^{+\infty}(1+\mu_l)-1\eeqa
where the infinite productory $\prod_{l=1}^{+\infty}(1+\mu_l)$ converges, being bounded by
\beqa{prod conv}
\prod_{l=0}^{+\infty}(1+\mu_l)&=& \exp\left[\sum_l\log{(1+\mu_l)}\right]\nonumber\\
&\leq&\exp\left[\sum_l\mu_l\right]\nonumber\\
&=&\exp\left[2^8E\sum_l(E)^l\right]\nonumber\\
&\leq&\exp\left[2^9E\right]\nonumber\\
&\leq&1+2^{10}E\nonumber\\
\eeqa
having used the elementary estimate $e^x\leq 1+2x$ for $0\leq x\leq 1$. In follows from (\ref{Lip for ellj}), (\ref{prod conv})
$${\cal L}_{\cP}(\ell-\id)\leq\limsup_j {\cal L}_{\cal P}(\ell_{j+1}-\textrm{id})\leq 2^{10}E\ .$$ 
\vskip.1in
\noindent
$\underline{\textrm{\sl Claim}\ 2 .}$ {\sl(``definition and bounds for $\o_*$'')}: {\it The bi--Lipschitz homeomorphism defined on $\cI_*$ by $\o_*=(\bar\o_*,\hat\o_*):=\o\circ \ell^{-1}$ is onto on $\cD_{\g,\hat\g,\t}^{\bar n,\hat n}\cap\o(\cI)$ and is subject to the following bounds
\beqano
\sup_{\cI_*}|\bar\o_*-\bar\o|\leq 2^5M\tilde\r E\qquad \qquad&& \sup_{\cI_*}|\ovl{\o_*^{-1}}-\ovl{\o^{-1}}|\leq 2^5\frac{\bar N}{N}\tilde\r E\\
\sup_{\cI_*}|\hat\o_*-\hat\o|\leq 2^5\hat M\tilde\r E\qquad\qquad &&\sup_{\cI_*}|\widehat{\o_*^{-1}}-\widehat{\o^{-1}}|\leq 2^5\frac{\hat N}{N}\tilde\r E\\
\eeqano
}

\vskip.1in
\noi
{\sl Proof.}\ Trivially, $\n\in \cD^{\bar n,\hat n}_{\g,\hat\g,\t}$, then, $I_*:=\ell\circ\o^{-1}(\n)$ is its (unique) preimage, with
$$|\overline{\o_*^{-1}}(\n)-\overline{\o^{-1}}(\n)|=|\bar\ell\circ\o^{-1}(\n)-\overline{\o^{-1}}(\n)|\leq \sup_{\cI_0}|\bar\ell-\id|\leq \frac{\bar N}{N}{2^5\tilde\r E}\ .$$
Similarly,
$$|\widehat{\o_*^{-1}}(\n)-\widehat{\o^{-1}}(\n)|\leq \frac{\hat N}{N}{2^5\tilde\r E}\ .$$
Using (\ref{ell close to id1}) (recall $\bar N$, $\hat N\leq N$), we find
\begin{eqnarray*}
|\bar\o_*(I_*)-\bar\o(I_*)|&=&|\bar\o(I)-\bar\o(\ell(I)|\nonumber\\
&\leq&M\,|I-\ell(I)|\nonumber\\
&\leq&2^5M\tilde{\r}E\ ,
\end{eqnarray*}
The proof of
$$|\hat\o_*(I_*)-\hat\o(I_*)|\leq2^5\hat M\tilde{\r}E$$
is similar.
\vskip.1in
\noindent
$\underline{\textrm{\sl Claim}\ 3.}$ {\sl(``construction of $\Phi$'')}: {\it The sequence of real--analytic symplectomorphisms, 
defined on $\cP_j={\cI_j}_{\r_j}\times \torus^n_{\bar s+s_j}$,
$$\Phi_j:=\Psi_1 \circ \cdots \circ \Psi_j$$ converges uniformly on $\cP_*=\cI_*\times \torus^n_{\bar s}$,  to  an $\cI_*$--family of real--analytic embeddings $$I_*\in\cI_*\to \Phi(I_*,\cdot):\ \torus^n\to \Re(\cI_\r)\times\torus^n$$
where $\Phi(I_*,\vartheta)=(\Phi_I(I_*,\vartheta),\Phi_{\varphi}(I_*,\vartheta))$ is given by
\beqano
\arr{
\Phi_I(I_*,\vartheta)=I_*+a(I_*,\vartheta)\\
\Phi_{\varphi}(I_*,\vartheta)=\vartheta+b(I_*,\vartheta)
}
\eeqano
with $a=(\bar a,\hat a)$, $b$ verifying
\beqa{bounds on a b}\sup_{\cI_*\times\torus^n}|\bar a|\leq2\,\frac{\hat\g}{\g}E\,\tilde{\r}\ ,\quad \sup_{\cI_*\times\torus^n}|\hat a|\leq2\,E\,\tilde{\r}\ ,\quad \sup_{\cI_*\times\torus^n}|b|\leq2\,E\,s\ 
\eeqa
and $\vartheta\to \vartheta+b(I_*,\vartheta)$ a diffeomorphism of $\torus^n$, for any $I\in \cI_*$.
Furthermore, the rescaled map
$$\check{\Phi}:=W_{\tilde\r, s}\Phi\circ W_{\tilde\r, s}^{-1}:\quad \tilde{\r}^{-1}{\cal I}_*\times s^{-1}\torus^n_{\bar s}\to r^{-1}{\cal I}_{\r}\times s^{-1}\torus^n_{\bar s+s}$$
is bi--Lipschitz, with
$${\cal L}_{\cP}(\check\Phi-\id)\leq 2^7E\ .$$
}

\vskip.1in
\noindent
{\sl Proof.}\ The bound (\ref{deriv1}) implies that rescaled maps 
$$\check{\Phi}_j:=W_{\tilde\r, s}\Phi_j\circ W_{\tilde\r, s}^{-1}:\quad \tilde{\r}^{-1}{\cal I}_{j\r_j}\times s^{-1}\torus^n_{\bar s+s_j}\to r^{-1}{\cal I}_{\r}\times s^{-1}\torus^n_{\bar s+s}\ ,$$
 are bi--Lipschitz on $\tilde{\r}^{-1}{\cal I}_{j\r_j}\times s^{-1}\torus^n_{\bar s+s_j}$ (hence, on $\tilde{\r}^{-1}{\cal I}_{*}\times s^{-1}\torus^n$), with
\begin{eqnarray}\label{Phi is Lipschitz2}
{\cal L}_{\cal P}(\check{\Phi}_j-\rm id)\leq {2^7E}\leq \frac{1}{2}\ ,
\end{eqnarray}
because 
\begin{eqnarray}\label{Phi is Lipschitz}
{\cal L}_{\cal P}(\check{\Phi}_j-\rm id)&\leq&\sup_{\cP_j}\|W_{\tilde\r, s}D\Phi_jW_{\tilde\r, s}^{-1}-\id_{2n}\|_{\cal P}\nonumber\\
&=&\|W_{\tilde\r, s}D\Psi_1\cdots D\Psi_jW_{\tilde\r, s}^{-1}-\id_{2n}\|_{\cal P}\nonumber\\
&\leq& \prod_0^{j-1}(1+\varsigma_k)-1\nonumber\\
&\leq& \prod_0^{+\infty}(1+\varsigma_k)-1\nonumber\\
&\leq&2^7E
\end{eqnarray}
having used
\begin{eqnarray}\label{Phi is Lipschitz1}
\prod_{k=0}^{N}(1+\varsigma_k)
&=&\exp\left[\sum_{k=0}^{N}\log{(1+\varsigma_k)}\right]\nonumber\\
&\leq& \exp\left[\sum_{k=0}^{N}\varsigma_k\right]\nonumber\\
&=&\exp\left[2^5\,\sum_{k=0}^{N}\left[\frac{(12)^{\t+1}}{3}\right]^{k}{E^{2^k}}\right]\nonumber\\
&\leq&\exp\left[2^5\,\sum_{k=0}^{+\infty}\left[\frac{(12)^{\t+1}}{3}\right]^{k}{E^{2^k}}\right]\nonumber\\
&\leq&\exp\left[2^6E\right]\nonumber\\
&\leq&1+2^7E\nonumber\\
\end{eqnarray}
Hence, the uniform convergence of $\Phi_j$ on $\cP_*$ easily follows: if $i+1\leq j$,
\begin{eqnarray}\label{convergence of Phij}
\sup_{\cI_*\times\torus^n_{\bar s}}|W_{\tilde\r, s}(\Phi_j-\Phi_i)|_{\cal P}&=&\sup_{\cI_*\times\torus^n_{\bar s}}|W_{\tilde\r, s}\Phi_i(\Psi_{i+1}\circ \cdots \circ \Psi_j)-W_{\tilde\r, s}\Phi_i|_{\cal P}\nonumber\\
&=&\sup_{\cI_*\times\torus^n_{\bar s}}|\check{\Phi}_i(W_{\tilde\r, s}\Psi_{i+1}\circ \cdots \circ \Psi_j)-\check{\Phi}_i(W_{\tilde\r, s})|_{\cal P}\nonumber\\
&\leq&{\cal L}_{\cal P+}(\check{\Phi}_i)\sup_{\cI_*\times\torus^n_{\bar s}}|W_{\tilde\r, s}(\Psi_{i+1}\circ \cdots \circ \Psi_j-\id)|_{\cal P}\nonumber\\
&\leq&\left(1+{2^7E}\right)\,\sum_{l=i}^{j-1}\left(\frac{1}{6}\right)^{l}\,E^{2^l}
\end{eqnarray}
having used (\ref{Psij close to id}), for which
\begin{eqnarray}\label{Phi is close to identity}
\sup_{\cI_*\times\torus^n_{\bar s}}|W^{\bar n}_{\hat\g/\g,1}W_{\tilde\r, s}(\Psi_{i+1}\circ \cdots \circ \Psi_j-\textrm{id})|_{\cal P}\leq\sum_{l=i}^{j-1}\left(\frac{1}{6}\right)^{l}\,E^{2^l}\qquad i+1\leq j
\end{eqnarray}
Denote then by $\Phi$ the uniform limit of $\Phi_j$ on $\cP_*$. Taking, in (\ref{Phi is close to identity}), $i=0$ and letting $j\to \infty$, we find 
\begin{eqnarray}\label{W phi-W is small}
& &\max\left\{\frac{\g}{\hat\g\tilde\r}\sup_{\cI_*\times\torus^n_{\bar s}}|\Phi_{\bar I}-\id_{\bar n}|,\ \frac{1}{\tilde\r}\sup_{\cI_*\times\torus^n_{\bar s}}|\Phi_{\hat I}-\id_{\hat n}|,\ \frac{1}{s}\sup_{\cI_*\times\torus^n_{\bar s}}|\Phi_{\varphi}-\id_{n}|\right\}\nonumber\\
&=&
|W^{n}_{\hat\g/\g,1}W_{\tilde \r,s}(\Phi-\textrm{id})|_{\cP}\nonumber\\
&\leq&\sum_{l=0}^{+\infty}\left(\frac{1}{6}\right)^{l}\,E^{2^l}\nonumber\\
&\leq& \frac{E}{1-\frac{E}{6}}<2E\ ,
\end{eqnarray}
which clearly implies (\ref{bounds on a b}). But (\ref{W phi-W is small}) also implies that, for any fixed $I_*$, the analytic map $$\vartheta\to \Phi_{\varphi}(I_*,\vartheta)=\vartheta+b(I_*,\vartheta)$$ is a diffeomorphism of $\torus^n$. In fact, by General Cauchy Inequality,
$$\sup_{\cI_*\times\torus^n}|\partial_{\vartheta}b|\leq \frac{2s}{\bar s}E<1\qquad \textrm{as}\quad \bar s\geq 2s\quad \textrm{and}\quad E<1\ .$$
Finally, by (\ref{Phi is Lipschitz2}), the rescaled map 
$$\check{\Phi}:=W_{\tilde\r, s}\Phi\circ W_{\tilde\r, s}^{-1}:\quad \tilde{\r}^{-1}{\cal I}_*\times s^{-1}\torus^n_{\bar s}\to r^{-1}{\cal I}_{\r}\times s^{-1}\torus^n_{\bar s+s}$$
is bi--Lipschitz, with
$${\cal L}_{\cP}(\check\Phi-\id)\leq \limsup_j{\cal L}_{\cP}(\check\Phi_j-\id)\leq 2^7E\leq \frac{1}{2}$$

\vskip.1in
\noindent
$\underline{\textrm{\sl Claim}\ 4.}$ {\it $\Phi(I_*,\torus^n)$ is a Lagrangian torus with frequency $\o_*(I_*)$ for ${\rm H}$.}:
\vskip.1in
\noindent
{\sl Proof.}\ We have to prove that the H--flow $\phi_t^{\rm H}\Big(\Phi(I_*,\vartheta)\Big)$ of a point $\Phi(I_*,\vartheta)\in \Phi(I_*,\torus^n)$ evolves as 
$$\phi_t^{\rm H}\Big(\Phi(I_*,\vartheta)\Big)=\Phi(I_*,\vartheta+\o_*(I_*)t)\qquad t\geq 0\ .$$
We split $\dst |\phi_t^{\rm H}\Big(\Phi(I_*,\vartheta)\Big)-\Phi(I_*,\vartheta+\o_*(I_*)t)|$ as
\begin{eqnarray*}
|\phi_t^{\rm H}\Big( \Phi(I_*,\vartheta)\Big)-\Phi(I_*,\vartheta+\o_*(I_*)\,t)|&\leq&|\phi_t^{\rm H}\Big( \Phi(I_*,\vartheta)\Big)-\phi_t^{\rm H}\Big( \Phi_j(I_*,\vartheta)\Big)|\nonumber\\
&+&|\phi_t^{\rm H}\Big( \Phi_j(I_*,\vartheta)\Big)-\Phi_j(I_*,\vartheta+\o_*(I_*)t)|\nonumber\\
&+&|\Phi_j(I_*,\vartheta+\o_*(I_*)t)-\Phi(I_*,\vartheta+\o_*(I_*)t)|
\end{eqnarray*}
Due to the uniform convergence of $\Phi_j$ to $\Phi$ on $\cP_*$ and continuous dependence on the initial data, both 
$|\Phi_j(I_*,\vartheta+\o_*(I_*)t)-\Phi(I_*,\vartheta+\o_*(I_*)t)|$ and $|\phi_t^{\rm H}\Big( \Phi(I_*,\vartheta)\Big)-\phi_t^{\rm H}\Big( \Phi_j(I_*,\vartheta)\Big)|$ go to $0$
as $j\to \infty$. We prove that $|\phi_t^{\rm H}\Big( \Phi_j(I_*,\vartheta)\Big)-\Phi_j(I_*,\vartheta+\o_*(I_*)t)|$ also goes to $0$, which will conclude the proof. As for the canonicity of $\Phi_j$ on $\cP_j$, 
$$\phi_t^{\rm H}\Big( \Phi_j(I_*,\vartheta)\Big)=\Phi_j\Big(\phi_t^{{\rm H}_j}(I_*,\vartheta)\Big)$$
with $${\rm H}_j={\rm H}\circ\Phi_j={\rm h}_j+{\rm f}_j\ .$$ So, using the Lipschitz property for the rescaled maps $\check\Phi_j=W^{-1}_{\tilde\r,\s}\Phi\circ W_{\tilde\r,\s}$ on $\tilde\r^{-1}\cI_*\times s^{-1}\torus^n$, with ${\cal L}_{\cP}(\check\Phi_j)\leq 1+2^7E$, we find
\beqano
|\phi_t^{\rm H}\Big( \Phi_j(I_*,\vartheta)\Big)-\Phi_j(I_*,\vartheta+\o_*(I_*)t)|&=&|\Phi_j\Big(\phi_t^{\rm H}(I_*,\vartheta)\Big)-\Phi_j(I_*,\vartheta+\o_*(I_*)t)|\nonumber\\
&=&|W_{\tilde\r,\s}\Big(\check\Phi_j\circ W^{-1}_{\tilde\r,\s}\Big(\phi_t^{{\rm H}_j}(I_*,\vartheta)\Big)\nonumber\\
&-&\check\Phi_j\circ W^{-1}_{\tilde\r,\s}(I_*,\vartheta+\o_*(I_*)t)\Big)|\nonumber\\
&\leq&\|W_{\tilde\r,\s}\|(1+2^7E)\nonumber\\
&\times&|W^{-1}_{\tilde\r,\s}\Big(\phi_t^{{\rm H}_j}(I_*,\vartheta)\Big)-W^{-1}_{\tilde\r,\s}(I_*,\vartheta+\o_*(I_*)t)\Big)|\nonumber\\
&\leq&\|W_{\tilde\r,\s}\|(1+2^7E)\|W^{-1}_{\tilde\r,\s}\|\nonumber\\
&\times&|\phi_t^{{\rm H}_j}(I_*,\vartheta)-(I_*,\vartheta+\o_*(I_*)t)|\nonumber\\
&=&\|W_{\tilde\r,\s}\|(1+2^7E)\|W^{-1}_{\tilde\r,\s}\|\nonumber\\
&\times&|(I_j(t)-I_*,\varphi_j(t)-(\vartheta+\o_*(I_*)t))|
\eeqano
where we have let $\Big(I_j(t),\varphi_j(t)\Big)=\phi_t^{{\rm H}_j}(I_*,\vartheta)$  the ${\rm H}_j$--evolution  of $(I_*,\vartheta))$. Representing $\Big(I_j(t),\varphi_j(t)\Big)$ as
\beqano\arr{I_j(t)=I_*-\int_0^t\partial_{\varphi}{\rm H}_j(I_j(\t),\varphi_j(\t))d\t\\
\varphi_j(t)=\vartheta+\int_0^t\partial_{I}{\rm H}_j(I_j(\t),\varphi_j(\t))d\t
}\eeqano
we find, by General Cauchy Inequality,
\begin{eqnarray}\label{Ij-I*}
|I_j(t)-I_*|&=&\left|\int_0^t\partial_{\varphi}{\rm H}_j(I_j(\t),\varphi_j(\t))d\t\right|\nonumber\\
&=&\left|\int_0^t\partial_{\varphi}f_j(I_j(\t),\varphi_j(\t))d\t\right|\nonumber\\
&\leq&\int_0^t\left|\partial_{\varphi}f_j(I_j(\t),\varphi_j(\t))\right|d\t\nonumber\\
&\leq&\frac{F_j}{s_j}t\to 0\quad \textrm{as}\quad j\to+\infty
\end{eqnarray}
Now, letting $I_0:=\ell^{-1}(I_*)\in \cI_0$ and using $\o_j(\ell_j(I_0))=\o_*(I_*)$, we find
\beqa{thetaj-theta-nu t}
\left|\varphi_j(t)-\vartheta-\o_*(I_*)t\right|&=&\left|\int_0^t\partial_{I}{\rm H}_j(I_j(\t),\varphi_j(\t))d\t-\o_*(I_*)t\right|\nonumber\\
&=&\left|\int_0^t\Big(\partial_{I}{\rm H}_j(I_j(\t),\varphi_j(\t))-\o_*(I_*)\Big)d\t\right|\nonumber\\
&=&\left|\int_0^t\Big(\o_j(I_j(\t))+\partial_{I}{\rm f}_j(I_j(\t),\varphi_j(\t))-\o_*(I_*)\Big)d\t\right|\nonumber\\
&\leq&\int_0^t\left|\o_j(I_j(\t))+\partial_{I}{\rm f}_j(I_j(\t),\varphi_j(\t))-\o_*(I_*)\right|d\t\nonumber\\
&\leq&\int_0^t\left|\o_j(I_j(\t))-\o_j(I_*)\right|d\t+\int_0^t\left|\o_j(I_*)-\o_j(\ell_j(I_0))\right|d\t\nonumber\\
&+&\int_0^t\left|\partial_{I}{\rm f}_j(I_j(\t),\varphi_j(\t))\right|d\t\nonumber\\
\eeqa
where
$$\int_0^t\left|\o_j(I_j(\t))-\o_j(I_*)\right|d\t\leq M_j\sup_{[0,t]}|I_j(\t)-I_*|t\leq \frac{M_jF_j}{s_j}t\to 0$$
by (\ref{Ij-I*});
$$\int_0^t\left|\o_j(I_*)-\o_j(\ell_j(I_0))\right|d\t=\left|\o_j(\ell(I_0))-\o_j(\ell_j(I_0))\right|t\leq M_j\sup_{\cI_0}|\ell_j-\ell|t\leq 2^5\tilde\r M_j E^{2^j}t\to 0$$
by (\ref{ellj-ell}), and, finally, by General Cauchy Inequality,
$$\int_0^t\left|\partial_{I}{\rm f}_j(I_j(\t),\varphi_j(\t))\right|d\t\leq \frac{F_j}{\r_j}t\to 0\ .$$
\vskip.1in
\noindent
$\underline{\textrm{\sl Claim}\ 5}$. {\it The ``invariant'' set ${\rm K}:=\Phi(\cI_*\times \torus^n)$ satisfies the measure estimate
\beqano
\meas\Big(\cI\times\torus^n\setminus{\rm K}\Big)&\leq&\Big(1+(1+2^7E)^{2n}\Big)\Big(1+(1+2^{10}E)^{n}\Big)\meas\Big(\cI\setminus \cI_0\times \torus^n\Big)\nonumber\\
&+&\Big(1+(1+2^7E)^{2n}\Big)(1+2^{10}E)^n\meas\Big(\cI_{\r_1}\setminus\cI\times\torus^n\Big)\nonumber\\
&+&(1+2^7E)^{2n}\meas\Big(\cI_{\r_2}\setminus \cI\times\torus^n\Big)\ .
\eeqano
where
$\r_1=2^6E\tilde\r/(1-2^{10}E)$, $\r_2=4E\tilde\r/(1-2^7E)$.}

\vskip.1in
\noi
{\sl Proof.}\ Let $\r_1=2^6E\tilde\r/(1-2^{10}E)$. Extend the Lipschitz function  $\dst \ell-\id:\ \cI_0\to \cI_*$ to a Lipschitz function $\ell_{\rm e}-\id$ on $\cI_{0r}$, with the same Lipschitz constant ${\cal L}(\ell_{\rm e})(\ell-\id)\leq 2^{10}E$ (this is made possible thanks to Lemma \ref{extension}). Then, $\ell_{\rm e}$ is a bi--Lipschitz extension (hence, injective) of $\ell$ on $\cI_{0r}$, with lower Lipschitz constant ${\cal L}_-(\ell_{\rm e})\geq 1-2^{10}E$. This implies that $\ell_{\rm e}$ sends a ball with radius $\r_1$ centered at $I_0\in \cI_0$ over a ball with radius $(1-2^{10}E)\r_1=2^6E\tilde\r>2^5E\tilde\r\geq |\ell(I_0)-I_0|$ centered at $\ell(I_0)$, so as to conclude
$$\ell_{\rm e}(\cI_{0\r_1})\supset \cI_0\ .$$
Then,
\beqano
\meas\Big(\cI_0\setminus \cI_*\Big)&\leq&\meas\Big(\ell_{\rm e}(\cI_{0\r_1})\setminus \cI_*\Big)\nonumber\\
&=&\meas\Big(\ell_{\rm e}(\cI_{0\r_1})\setminus \ell_{\rm e}(\cI_0)\Big)\nonumber\\
&\leq&\meas\Big(\ell_{\rm e}(\cI_{0\r_1}\setminus \cI_0)\Big)\nonumber\\
&\leq&{\cal L}(\ell)^{n}\meas\Big(\cI_{0\r_1}\setminus \cI_0\Big)\nonumber\\
&\leq&{\cal L}(\ell)^{n}\meas\Big(\cI_{\r_1}\setminus \cI_0\Big)\nonumber\\
&\leq&{\cal L}(\ell)^{n}\Big(\meas\Big(\cI_{\r_1}\setminus \cI\Big)+\meas\Big(\cI\setminus \cI_0\Big)\Big)\nonumber\\
\eeqano
and, finally,
\beqa{est I-I*}
\meas\Big(\cI\setminus \cI_*\Big)&\leq&\meas\Big(\cI\setminus \cI_0\Big)+\meas\Big(\cI_0\setminus \cI_*\Big)\nonumber\\
&\leq&\Big(1+{\cal L}(\ell)^{n}\Big)\meas\Big(\cI\setminus \cI_0\Big)+{\cal L}(\ell)^{n}\meas\Big(\cI_{\r_1}\setminus \cI\Big)\nonumber\\
&\leq&\Big(1+(1+2^{10}E)^{n}\Big)\meas\Big(\cI\setminus \cI_0\Big)+(1+2^{10}E)^{n}\meas\Big(\cI_{\r_1}\setminus \cI\Big)\nonumber\\
\eeqa
On turn, using 
$$\arr{\sup_{\tilde\r^{-1}\cI_{*}\times s^{-1}\torus^n}|\check\Phi-\id|\leq 2E\\ {\cal L}(\check\Phi-\id)\leq 2^7E}$$
we can repeat the above argument, with $$(\check\Phi,\check{\rm K}=\check\Phi(\tilde\r^{-1}\cI_*\times s^{-1}\torus^n),\tilde\r^{-1}\cI_*\times s^{-1}\torus^n,2n,{\cal L}(\check\Phi),\check r=2^2E/(1-2^7E))$$ replacing $(\ell,\cI_*,\cI_0,n,{\cal L}(\ell),r)$ and we find
\beqano
\meas\Big(\tilde\r^{-1}\cI\times s^{-1}\torus^n\setminus\check{\rm K}\Big)&\leq&\Big(1+(1+2^7E)^{2n}\Big)\meas\Big(\tilde\r^{-1}\cI\setminus (\tilde\r^{-1}\cI_*)\times s^{-1}\torus^n\Big)\nonumber\\
&+&(1+2^7E)^{2n}\meas\Big((\tilde\r^{-1}\cI)_{\check r}\setminus (\tilde\r^{-1}\cI)\times s^{-1}\torus^n\Big)\ .
\eeqano
Hence, rescaling the variables,
\beqano
\meas\Big(\cI\times\torus^n\setminus{\rm K}\Big)&\leq&\Big(1+(1+2^7E)^{2n}\Big)\meas\Big(\cI\setminus \cI_*\times \torus^n\Big)\nonumber\\
&+&(1+2^7E)^{2n}\meas\Big(\cI_{\r_2}\setminus \cI\times\torus^n\Big)\qquad \r_2=\check r\tilde\r
\eeqano

\vskip.1in
\noi
Finally, taking into accunt (\ref{est I-I*}), we find the result.
\vskip.1in
\noi
$\underline{\textrm{{\sl Conclusion of the proof of Theorem \ref{two scales KAM}.}}}$\ Take
$$\phi_\n=\Phi(\o_*^{-1}(\n),\cdot)$$
and recognize that
$$\ell=\o_*^{-1}\circ \o\quad \textrm{on}\quad \cI_0=\cI_{\g,\hat\g,\t}$$
has Lipschitz norm bounded as in (\ref{Lip Norm}), by (\ref{ell close to id1}) and (\ref{Lip const for ell}).

\noi
\subsubsection{Nondegenerate KAM Theorem via Theorem \ref{two scales KAM}}\label{nondegenerate KAM}
Taking, in Theorem \ref{two scales KAM}, 
$$\g=\hat\g\ ,\quad \hat M=M\ ,\quad \bar N=\hat N=N$$
gives a standard (nondegenerate, isofrequencial) KAM Theorem:
\begin{theorem}\label{nondegenerate KAM Thm}
Let $n\in \natural$, $\t>n$, $\g>0$,  $\cI\subset\real^n$ compact and let
\beqano
{\rm H}
(J,\psi)={\rm h}(J)+{\rm f}(J,\psi)
\eeqano 
real--analytic on $\cI_\r\times \torus_s^{n}$, where $\o:=\partial {\rm h}$ is a diffeomorphism of $\cI_\r$ with Jacobian matrix $U:=\partial^2{\rm h}$ non singular on $\cI_\r$.
Let
\beqano
&&  \dst M\geq\sup_{{\cal I}_{\r}}\|U\|\\
&&  \dst N\geq\sup_{{\cal I}_{\r}}\|T\|\\
&&  \dst F\geq\|{\rm f}\|_{\r,s}
\eeqano
where $T:=U^{-1}$, define
\beqano
&&  \dst c:=\max\left\{{2^{11}n},\ \frac{2}{3}(12)^{\t+1}\right\}\\
&&  \dst K:=\frac{6}{s}\ \log_+{\left(\frac{FM^2\,L}{\gamma^2}\right)^{-1}}\quad \textrm{where}\quad \log_+(a):=\max\{1,\log{a}\}\\ 
&&  \dst \tilde{\r}:=\min\left\{\frac{\g}{3MK^{\t+1}}\ ,\ 
 \r\right\}\\
&&  \dst L:=\max\ \{N\ , M^{-1}\}\\
\eeqano
and assume the ``perturbation'' $\rm f$ so ``small'' that 
\beq{KAM cond non deg}
cE:=c\frac{FL}{\tilde{\r}^2}<1\ .
\eeq
\begin{itemize}
\item[(i)]
Then, for any frequency $\nu\in\O_*:=\o(\cI)\cap\cD^{n}_{\g,\t}$, there exists a unique  Lagrangian KAM torus ${\rm T}_\n\subset \Re({\cal I}_{34\tilde{\r}E})\times \torus^{n}$ for $H$ with frequency $\nu$, such that the following holds. 
There exists a ``Cantor'' set ${\cal I}_{*}\subset \Re ({\cal I}_{32\tilde{\r}E})$  and a   bi--Lipschitz (onto) {homeomorphism} 
\beqano
\omega_*:\ {\cal I}_{*}\to \O_*
\eeqano
satisfying
\beqano
&&\sup_{{\cal O}_*}| \o_*^{-1}-\o^{-1}|\leq 2^5\tilde\r\,E\ ,\quad \sup_{{\cal I}_*}|\omega_{*}-\omega|\leq{2^5}M\tilde{\r}E
\\
&&
\|\o_*^{-1}\circ \o-\textrm{\rm id}\|^{\rm Lip}_{\tilde{\r},{\cal I}_{\g,\t}}\leq 2^{10}E\ ,
\qquad
{\cal I}_{\g,\t}:=\o^{-1}({\cal D}^{n,\hat n}_{\g,\t})\cap{\cal I}\ .
\eeqano
such that ${\rm T}_\n$ is
realized by the real--analytic embedding $\phi_\n=(\phi_{\n I},\phi_{\n \varphi})$ given by
\begin{eqnarray*}
\left\{
\begin{array}{lrr}
\phi_{\n I}(\vartheta)=I_{*}(\nu)+v(\nu,\vartheta)&&\\
\phi_{\n \varphi}(\vartheta)=\vartheta+u(\nu,\vartheta)&&\\
\end{array}
\right. \quad \vartheta\in \torus^n\ ,
\end{eqnarray*}
where $I_*(\n):=\o_*^{-1}(\n)$ and $v$, $u$ are bounded as
\begin{eqnarray*}
|v(\nu,\vartheta)|\leq2\,E\,\tilde{\r}\ ,\quad |u(\nu,\vartheta)|\leq2\,E\,s
\end{eqnarray*}
\item[(ii)] The measure of the invariant set ${\rm K}=\phi_{\O_*}(\torus^n)$ satisfies
\beqano
\meas\Big(\cI\times\torus^n\setminus{\rm K}\Big)&\leq&\Big(1+(1+2^7E)^{2n}\Big)\Big(1+(1+2^{10}E)^{n}\Big)\meas\Big(\cI\setminus \cI_{\g,\t}\times \torus^n\Big)\nonumber\\
&+&\Big(1+(1+2^7E)^{2n}\Big)(1+2^{10}E)^n\meas\Big(\cI_{\r_1}\setminus\cI\times\torus^n\Big)\nonumber\\
&+&(1+2^7E)^{2n}\meas\Big(\cI_{\r_2}\setminus \cI\times\torus^n\Big)\ .
\eeqano
with
$\r_1=2^6E\tilde\r/(1-2^{10}E)$, $\r_2=4E\tilde\r/(1-2^7E)$.
\end{itemize}
\end{theorem}

\subsection{Properly Degenerate KAM Theory (Proof of Theorem \ref{more general degenerate KAM})}\label{Properly degenerate KAM Theorem}
The aim of this section is to prove Theorem \ref{more general degenerate KAM} as an application of Theorem \ref{two scales KAM}. 
\vskip.1in
\noi 
We quote a refined averaging theory by Biasco et al (\cite{BCV03}, `` full version'' on p. $110$), in which statement the ``sup--Fourier''  norm  of
$$f(I,\varphi,p,q)=\sum_kf_k(I,p,q)e^{ik\cdot\varphi}$$ 
(with $f_k$ the Fourier coefficient of $f$), real--analytic on $U_{\bar{r}}\times \torus^n_{\bar s}\times E_{r_p}\times F_{r_q}$, with $U\subset\real^n$, $E$, $F\subset\real^m$, is defined as
$$\|f\|_{\bar r,\bar s,r_p,r_q}:=\sum_k\sup_{U_{{\bar r}}\times E_{r_p}\times F_{r_q}}|f_k|e^{|k|\bar s}$$
In order to avoid confusion with other parameters here introduced, we denote by $a$, $\bar\e$, $\bar r$, $\bar s$, $\bar d$ the parameters $\a$, $\e$, $r$, $s$, $d$ of \cite{BCV03}, but we do not change the name of the dimensions $(n,m)$ (which correspond to ``our'' $(\bar n,\hat n)$) thereby used, letting the reader be aware not to confuse \cite{BCV03}'s $n$ (which corresponds to ``our'' $\bar n$) with ``our'' $n=\bar n+\hat n$.
\begin{proposition}[Fast Averaging Theorem]\label{average BCV}
Let $H:=h(I)+f(I,\varphi,p,q)$ be a real-analytic Hamiltonian on $U_{{\bar r}}\times \torus^n_{\bar s}\times E_{r_p}\times F_{r_q}$. Denoting  $\omega:=\partial_I\,h$ and $c_m:=e(1+em)/2$, suppose that
\begin{eqnarray*}
|\omega(I)\cdot k|\geq a\ ,\quad \textrm{for all}\quad I\in U_{\bar r}\ ,\quad k\in \integer^n\ ,\quad k\notin \L\ ,\quad |k|\leq K\ ,
\end{eqnarray*} 
where $\L$ is a $\integer^n$--module, $K\bar s \geq 6$, and
\begin{eqnarray}\label{smallness BCV}
\bar\e:=\|f\|_{\bar r,\bar s,r_p,r_q}<\frac{a \bar d}{2^7c_mK\bar s}\ , \quad \bar d=\min\{\bar r\bar s,\ r_p\,r_q\}\ .
\end{eqnarray} 
Then, there exists a real-analytic, symplectic transformation
\begin{eqnarray}\label{transformation}
\Psi:\quad U_{\bar r/2}\times \torus^n_{s/6}\times E_{r_p/2}\times F_{r_q/2}&\to&U_{\bar r}\times \torus^n_{\bar s}\times E_{r_p}\times F_{r_q}\nonumber\\
(I',\varphi',p',q')&\to&(I,\varphi,p,q)=\Psi(I',\varphi',p',q')
\end{eqnarray}
such that
\begin{eqnarray*}
H_{*}:=H\circ \Psi=h+g+f_{*}\ ,
\end{eqnarray*}
with $g$ in normal form:
\begin{eqnarray}\label{complete resonance}
g=\sum_{k\in \L}g_k(I',p',q')\,e^{ik\cdot\varphi'}\ .
\end{eqnarray}
Moreover, when the projection of $\Psi(I',\varphi',p',q')$ onto the $I$--variables is denoted by\\ $I(I',\varphi',p',q')$, etc,
\begin{eqnarray}\label{results}
&&\|g-P_{\L}T_Kf\|_{\bar r/2,\bar s/6,r_p/2,r_q/2}\leq \frac{12}{11}\frac{2^7c_m\,\bar\e^2}{a d}\leq \frac{\bar\e}{4}\ , \nonumber\\
&&\|f_{*}\|_{\bar r/2,\bar s/6,r_p/2,r_q/2}\leq e^{-K\bar s/6}\frac{2^9c_m\bar\e^2}{a\,d}\leq e^{-K\bar s/6}\,\bar\e\ ,\nonumber\\
&&\max\left\{\bar s\ |I(I',\varphi',p',q')-I'|\ , \quad \bar r\ |\varphi(I',\varphi',p',q')-\varphi'|\ ,\right.\nonumber\\
&&\left.\quad r_q\ |p(I',\varphi',p',q')-p'|\ , \quad r_p\ |q(I',\varphi',p',q')-q'|\right\}\ \leq \frac{9\bar\e}{a}\ .
\end{eqnarray}
\end{proposition}

\vskip.1in
\noi
We are ready to begin  our proof.

\vskip.1in
\noi
Let 
$$\cH(I,\varphi,p,q)=h(I)+\e\,f(I,\varphi,p,q)$$
a real--analytic on $\cI_{\r_0}\times\torus^{\bar n}_{s_0}\times \bar B(r_0)_{r_0}\times \bar B(r_0)_{r_0}$ Hamiltonian,
where
$\bar B(r)=\ovl{B^{\hat n}_r(0)}$  and $r_0\leq\bar r/2$.
We assume that $\r_0$, $r_0$ are so small that $\bar f$ preserves, on $\cI_{\r_0}\times\torus^{\bar n}_{s_0}\times \bar B(r_0)_{r_0}\times \bar B(r_0)_{r_0}$, the form
\beqano
\bar f&=&f_0(I)+\sum_{1\leq i\leq \hat n}\O_i(I)\frac{p_i^2+q_i^2}{2}+\frac{1}{2}\sum_{1\leq i,j\leq \hat n}A_{i,j}(I)\frac{p_i^2+q_i^2}{2}\frac{p_j^2+q_j^2}{2}+ o_4
\eeqano
and the {\sl non--resonance}, {\sl non--degeneracy} assumptions, \ie, 
\beqa{non res non deg}
\arr{
\min_{0<|k|\leq 4}\inf_{\cI_{\r_0}}|\O\cdot k|>0\ ,\\
\\ 
\inf_{\cI_{\r_0}}|\textrm{\rm det}\,A|>0\ \label{non deg}.
}
\eeqa
We proceed in $6$ steps. 
\vskip.1in
\noi
$\underline{\textrm{\sl Step}\ 1.}$ (``{\sl fast averaging}'')\ {\it There exist $0<c_{\rm av}<1<C_{\rm av}$, $r_{\rm av}>0$ $\g_{\rm  av}$ such that, for any $0<r\leq r_{\rm av}$, $\bar\g>\g_{\rm av}\sqrt{\e}(\log{r^{-1}})^{\bar\t+1}$, $\cH$ is put into the form
\begin{eqnarray}\label{form of H}
{\cal H}'(\e,r,\bar\g;I',\varphi',p',q')&=&\cH\circ\Psi_{\rm av}(\e,r,\bar\g;I',\varphi',p',q')\nonumber\\
&=&h(I')+\e\,g'(\e,r,\bar\g;I',p',q')+\e r^5\,f'(\e,r,\bar\g;I',\varphi',p',q')\nonumber\\
\end{eqnarray}
by means of a real--analytic symplectomorphism
\begin{eqnarray*}
\phi_{\rm av}:\ \bar\cI_{\bar\r(r)/2}\times \torus^{\bar n}_{s_0/6}\times {\bar B}(r_0)_{{r}_0/2}\times {\bar B}(r_0)_{{r_0}_0/2}&\to&\bar\cI_{\bar \r(r)}\times \torus^{\bar n}_{s_0}\times{\bar B}(r_0)_{{r}_0}\times {\bar B}(r_0)_{{r}_0}\nonumber\\
\end{eqnarray*}
where:
\begin{itemize} 
\item[(i)] $\bar\cI$, $\bar\r(r)$ are given by
\beqano
&& \bar\cI=\cI_{\bar\g,\bar\t}:=\Big\{I\in \cI:\ \o_0(I):=\partial h(I)\in \cD^{\bar n}_{\bar\g,\bar\t}\Big\}\nonumber\\
&& \bar\r(r):=\min\left\{c_{\rm av}\frac{\bar\g}{\Big(\log{r^{-1}}\Big)^{\bar\t+1}},\ \r_0\right\}
\eeqano
\item[(ii)]
$g'$, $f'$ satisfy 
\begin{eqnarray*}
&&\sup_{\bar\cI_{\bar\r(r)/2}\times \bar B(r_0)_{r_0/2}\times \bar B(r_0)_{r_0/2}}|g'-\bar f|\leq C_{\rm av}\frac{\e(\log{r^{-1}})^{2\bar\t+1}}{\bar\g^2}\ , \nonumber\\
&&\|f'\|_{\bar\cI_{\bar\r(r)/2}\times \torus^{\bar n}_{s_0/6}\times {\bar B}(r_0)_{{r}_0/2}\times {\bar B}(r_0)_{{r_0}_0/2}}\leq \|f\|_{\r_0,s_0,r_0}\ ,\nonumber\\
\end{eqnarray*}
($\bar f=\langle f \rangle_{\varphi}$)
\item[(iii)] The projections $I(I',\varphi',p',q')$, $\cdots$, of $\phi_{\rm av}$ over the $I$, $\cdots$ variables satisfy
\beqa{average estimates}
\arr{
|I(I',\varphi',p',q')-I'|\leq C_{\rm av}\frac{\e(\log{r^{-1}})^{\bar\t}}{\bar\g}\\
|\varphi(I',\varphi',p',q')-\varphi'|\leq C_{\rm av}\frac{\e(\log{r^{-1}})^{2\bar\t+1}}{\bar\g^2}\\
|p(I',\varphi',p',q')-p'|\leq C_{\rm av}\frac{\e(\log{r^{-1}})^{\bar\t}}{\bar\g}\\
|q(I',\varphi',p',q')-q'|\leq C_{\rm av}\frac{\e(\log{r^{-1}})^{\bar\t}}{\bar\g}\\
}
\eeqa
\end{itemize}
}
\vskip.1in
\noi
{\sl Proof.}\ We apply Proposition \ref{average BCV} to our case, taking
\beqano
&& \L=\{0\}\ ,\quad U=\bar\cI\ ,\quad E=F=\bar B(r_0)\ ,\quad \bar\e=\e\|f\|_{\r_0,s_0,r_0}\nonumber\\
&& \bar r=\bar\r(r):=\min\left\{\frac{1}{2\bar M}\left(\frac{s_0}{30}\right)^{\bar\t+1}\frac{\bar\g}{\Big(\log{r^{-1}}\Big)^{-(\bar\t+1)}},\ \r_0\right\}\ ,\quad \bar s=s_0\nonumber\\
&& r_p=r_q=r_0\ ,\quad n=\bar n\ ,\quad m=\hat n\ ,\nonumber\\
&& K=\bar K(r):=\frac{30}{s_0}\log{r^{-1}}\ ,\quad a=\bar\a(r):={\bar\g}/({2\bar K(r)^{\bar\t}})
\eeqano
Observe, in particular, that  $K$ has been chosen such in a way to get a new perturbation $f_*$ of order
$$\|f_*\|_{\r_0/2,s_0/6,r_0/2}\leq\e e^{-\bar K(r)s_0/6}\|f\|_{\r_0,s_0,r_0}=\e r^5 \|f\|_{\r_0,s_0,r_0}\ ,$$  
so, we will put $f_*:=\e r^5 f'$.
We check, then,
\begin{itemize}
\item[(i)] $\bar K(r)s_0=30\log{r^{-1}}\geq 6\qquad$  (for $0<r<e^{-1/5}$);
\item[(ii)] If $I\in \bar\cI_{\bar\r(r)}$ and $I_0\in \bar\cI$ is such that $|I-I_0|<\bar\r(r)$, then, for $0<|\bar k|\leq \bar K(r)$, 
\beqano
\Big|\o_0(I)\cdot \bar k\Big|&\geq& \Big|\o(I_0)\cdot \bar k\Big|-\Big|\Big(\o(I)-\o(I_0)\Big)\cdot \bar k\Big|\nonumber\\
&\geq& \frac{\bar\g}{\bar K(r)^{\bar\t}}-\bar M\bar K(r)\bar\r(r)\nonumber\\
&=&\frac{\bar\g}{2\bar K(r)^{\bar\t}}=\bar\a(r)\ ;
\eeqano
where
$$\bar M=\sup_{\cI_{\r_0}}\|\partial^2\,h\|\ .$$
\item[(iii)] 
with $$\bar d=d(r)=\min\left\{\bar\r(r) s_0,\ r_0^2\right\}=\min\left\{\frac{s_0}{2\bar M}\left(\frac{s_0}{30}\right)^{\bar\t+1}\frac{\bar\g}{\Big(\log{r^{-1}}\Big)^{\bar\t+1}},\ s_0\r_0,\ r_0^2\right\}$$
using $\e\leq \sqrt{\e}$ (as $0\leq\e\leq 1$), we find a suitable constant $\tilde\g$ depending on $\bar\t$, $s_0$, $\r_0$, $\bar M$, for which
\beqano
\bar\e\frac{2^7c_mK\bar s}{a \bar d}&=&\e\|f\|_{\r_0,s_0,r_0}\frac{2^7c_{\hat n}\bar K(r)s_0}{\bar\a(r) d(r)}\nonumber\\
&=&{2^8 c_{\hat n}s_0\|f\|_{\r_0,s_0,r_0}}\left(\frac{30}{s_0}\right)^{\bar\t+1}\frac{\e(\log{r^{-1}})^{\bar\t+1}}{\bar\g}\nonumber\\
&\times&\max\left\{\frac{1}{r_0^2},\ \frac{1}{s_0\r_0},\ \frac{2\bar M}{s_0}\left(\frac{30}{s_0}\right)^{\bar\t+1}\frac{\Big(\log{r^{-1}}\Big)^{\bar\t+1}}{\bar\g}\right\}\nonumber\\
&\leq&\tilde\g\max\left\{\frac{\sqrt{\e}\Big(\log{r^{-1}}\Big)^{\bar\t+1}}{\bar\g},\ \left(\frac{\sqrt{\e}\Big(\log{r^{-1}}\Big)^{\bar\t+1}}{\bar\g}\right)^2\right\}\nonumber\\
&<&1
\eeqano
provided 
$$\bar\g>\tilde\g\sqrt{\e}\Big(\log{r^{-1}}\Big)^{\bar\t+1}$$
\end{itemize}
Hence, Proposition \ref{average BCV}, applies. Let, then, $g$, $f_*$, $\Psi$ as in claimed there, and put $\e g':=g$, $\e r^5 f':=f_*$, $\phi_{\rm av}=:=\Psi$. By the definition of $\tilde\r$, we find, then,
\beqano
|\varphi(I',\varphi',p',q')-\varphi'|&\leq& {18\|f\|_{\r_0,s_0,r_0}}\left(\frac{30}{s_0}\right)^{\bar\t}\frac{\e(\log{r^{-1}})^{\bar\t}}{\tilde\r\bar\g}\nonumber\\
&\leq&\hat C\max\left\{\frac{\e(\log{r^{-1}})^{2\bar\t+1}}{\bar\g^2},\ \frac{\e(\log{r^{-1}})^{\bar\t}}{\bar\g}\right\}\nonumber\\
&=&\hat C \frac{\e(\log{r^{-1}})^{2\bar\t+1}}{\bar\g^2}
\eeqano
as soon as $\sqrt{\e}(\log{r^{-1}})^{-1}<1$
Then, $g'$, $f'$, $\Psi_{\rm av}$ satisfy the claim.
\vskip.1in
\noi
$\underline{\textrm{\sl Step}\ 2.}$ (``{\sl preparation to Birkhoff Theory}'')\ {\it There exist $0<r_{\rm T}<1<C_{\rm T}$, $\g_{\rm T}>0$ such that, for any $0<r<r_{\rm T}/(8)$,  and $$\bar\g>\g_{\rm T}\max\{\sqrt{\e}(\log{r^{-1}})^{\bar\t+1},\ r^2(\log{r^{-1}})^{\bar\t+1},\ \sqrt[3]{\e r}(\log{r^{-1}})^{\bar\t+1}\}\ ,$$ 
there exists a real--analytic symplectomorphism $$\phi_{\rm T}=(I'(\e,r,\bar\g;\cdot),\varphi'(\e,r,\bar\g;\cdot),p'(\e,r,\bar\g;\cdot),p'(\e,r,\bar\g;\cdot))$$ on
$$
\bar\cI_{\bar\r(r)/4}\times \torus^{\bar n}_{s_0/12}\times B({}r)_{{}r}\times B({}r)_{{}r}
$$ $\cH'$ which puts $\cH'$ into the form
\beqa{H''old}\cH''(I'',\varphi'',p'',q'')&=&\cH'\circ \phi_{\rm T}(I'',\varphi'',p'',q'')\nonumber\\
&=&h(I'')+\e\,g''(I'',p'',q'')+\e r^5 f''(I'',\varphi'',p'',q'')\ ,
\eeqa
where:
\begin{itemize}
\item[(i)]\ $g''$ has an equilibrium point at $(p'',q'')=0$ with Hessian in $(p'',q'')=0$ satisfying
\beqano
\sup_{{\bar\cI}_{\bar\r(r)/4}}\|\partial^2 g''|_0-\textrm{\rm diag}(\O_1,\cdots \O_{\hat n},\O_1,\cdots \O_{\hat n})\|\leq C_{\rm T}\frac{\e(\log{r^{-1}})^{2\bar\t+1}}{\bar\g^2}\ ;
\eeqano
\item[(iii)] $f$ satisfies
$$\|f''\|_{\bar\cI_{\bar\r(r)/4}\times \torus^{\bar n}_{s_0/12}\times B({}r)_{{}r}\times B({}r)_{{}r}}\leq C_{\rm T}$$
\item[(iv)] the following bounds hold for $\phi_{\rm T}$, uniformly on $\bar\cI_{\bar\r(r)/4}\times \torus^{\bar n}_{s_0/12}\times B({}r)_{{}r}\times B({}r)_{{}r}$,
\begin{eqnarray}\label{id2}
\arr{I'(\e,r,\bar\g;\cdot)=I''\\
|\varphi'(\e,r,\bar\g;I'',\varphi'',p'',q'')-\varphi''|\leq C_{\rm T}\max\left\{\frac{r^2(\log{r^{-1}})^{\bar\t+1}}{\bar\g},\frac{\e r(\log{r^{-1}})^{3\bar\t+2}}{\bar\g^3}\right\}\\
|p'(\e,r,\bar\g;I'',\varphi'',p'',q'')-p''| \leq\frac{C_{\rm T}}{\bar\g^2}\,\e(\log{r^{-1}})^{2\bar\t+1}\\
|q'(\e,r,\bar\g;I'',\varphi'',p'',q'')-q''|\leq\frac{C_{\rm T}}{\bar\g^2}\,\e(\log{r^{-1}})^{2\bar\t+1}\\
}\nonumber\\
\end{eqnarray}
\end{itemize}
}

\vskip.1in
\noi
{\sl Proof. }\ Write $g'$ as 
$$g'(I',p',q')=\bar f(I',p',q')+\tilde g(I',p',q')$$
where
$$\bar f(I',p',q'):=f_0(I')+f_2(I',p',q')+f_4(I',p',q')=f_0(I')+\sum_{1\leq i\leq \hat n}\O_i(I')\frac{{p'_i}^2+{q'_i}^2}{2}+f_4$$
where, by assumption, $f_4$ is a power series in $(p',q')$ starting with 
\beq{ff11}f_4(I',p',q')=\frac{1}{2}\sum_{1\leq i,j\leq \hat n}A_{i,j}(I)\frac{{p_i'}^2+{q_i'}^2}{2}\frac{{p_j'}^2+{q_j'}^2}{2}+\cdots
\eeq
and, by the previous step, 
\beq{gg11}\|\tilde g\|_{\bar\r(r)/2, s_0/6,r_0/2}\leq C\frac{\e(\log{r^{-1}})^{2\bar\t+1}}{\bar\g^2}\leq C(\log{r^{-1}})^{-1}
\eeq
Then, $F(I',p',q'):=\partial_{(p',q')} g'(I',p',q')$ splits as
$$F(I',p',q')=F_0(I',p',q')+F_1(I',p',q')$$
with
$$F_0:=\left
(\begin{array}{lrr}
\O_1(I'){p_1'}\\
\vdots\\
\O_{\hat n}(I')p'_{\hat n}\\
\O_1(I'){q_1'}\\
\vdots\\
\O_{\hat n}(I')q'_{\hat n}\\
\end{array}
\right)\quad \textrm{and}\quad F_1:=\partial_{(p',q')}(f_4+\tilde g)$$
where $F_0$ is a diffeomorphism of $\complex^{\hat n}$ sending $0$ to $0$ and det$\partial F_0=\O_1(I)^2\cdots\O_{\hat n}(I)^2\neq 0$ on $\bar\cI_{\bar\r(r)/4}$, thanks to the non resonance condition (first in (\ref{non res non deg})). Furthermore, by (\ref{ff11}) and Cauchy Estimates on(\ref{gg11}), the following bound holds for $F_1$, uniformly on $\bar\cI_{\bar\r(r)/4}$:
$$\sup_{B^{2\hat n}_{r_0/4}(0)}|F_1|\leq C\max\left\{r_0^3,\ \frac{\e(\log{r^{-1}})^{2\bar\t+1}}{\bar\g^2r_0/4}\right\}\ .$$
This implies, by Cauchy estimates, that
\beqano
&& \max\left\{\sup_{B^{2\hat n}_{r_0/4}(0)}\|\partial F_1\|\sup_{B^{2\hat n}_{r_0/4}(0)}\|(\partial F_0)^{-1}\|,\ \frac{|\sup_{B^{2\hat n}_{r_0/4}(0)}\|(\partial F_0)^{-1}\|\sup_{B^{2\hat n}_{r_0/4}(0)}|F_1|}{r_0/8}\right\}\nonumber\\
&& \leq C\max\left\{r_0^2,\ \frac{\e(\log{r^{-1}})^{2\bar\t+1}}{\bar\g^2(r_0/4)^2}\right\}\nonumber\\
&& \leq \frac{1}{2}
\eeqano
as soon as $r_0$ is small enough and
$$\bar\g>\tilde\g{\sqrt{\e}(\log{r^{-1}})^{\bar\t+1/2}}.$$
Then, by the Quantitative Implicit Function Theorem (Appendix \ref{Quantitative Implicit Function Theorem}), for any $I'\in \cI_{\bar\r(r)/2}$, we find an equilibrium point $(p_{\rm e}(\e,r,\bar\g;I'), q_{\rm e}(\e,r,\bar\g;I'))$, with 
\begin{eqnarray}\label{small peq}
\sup_{\cI_{\bar\r(r)/2}}|(p_{\rm e}(\e,r,\bar\g;I'),q_{\rm e}(\e,r,\bar\g;I'))|\leq \frac{C}{\bar\g^2}\e(\log{r^{-1}})^{2\bar\t+1}\leq \frac{r_0}{8}
\end{eqnarray} for $g$, \ie, satisfying
$$F(\e,r,\bar\g;I',(p_{\rm e}(\e,r,\bar\g;I'), q_{\rm e}(\e,r,\bar\g;I'))=0\quad \textrm{for all}\quad I'\in \cI_{\bar\r(r)/2}\ .$$
Define, on $$\bar\cI_{\bar\r(r)/4}\times \torus^{\bar n}_{s_0/12}\times B( r)_{ r}\times B( r)_{ r}\quad \textrm{where}\quad  r<\frac{r_0}{8}$$ the transformation $\phi_{\rm T}(\e,r,\bar\g;\cdot):=(I'(\e,r,\bar\g,\cdot),\varphi'(\e,r,\bar\g,\cdot),p'(\e,r,\bar\g,\cdot),q'(\e,r,\bar\g,\cdot))$  by means of 
\footnote{$\phi_{\rm T}$ is generated by $S_{\rm T}=I''\cdot\varphi'+(p''+p_{\rm e}(\e,r,\bar\g;I'')\cdot(q'-q_{\rm e}(\e,r,\bar\g;I''))$}
\begin{eqnarray}\label{T}
\left\{
\begin{array}{lrr}
I'=I''\\
\varphi'=\varphi''-\partial_{I''}\Big(p''+p_{\rm e}(\e,r,\bar\g;I'')\Big)\cdot\Big(q'-q_{\rm e}(\e,r,\bar\g;I'')\Big)\\
p'=p_e(\e,r,\bar\g;I'')+p''\\
q'=q_e(\e,r,\bar\g;I'')+q''\\
\end{array}
\right.
\end{eqnarray}
Put
$$D(\e,r,\bar\g):=\{(I'',p'',q'):\{I''\in\bar\cI_{\bar\r(r)/4}\ ,\quad p''\in B( r)_{ r}\ ,\quad q'-q_{\rm e}(\e,r,\bar\g;I'')\in B( r)_{ r} \}$$
Then,  by Cauchy estimate,
\footnote{\label{ab le a2 b2}Use $a\,b\leq \max\{a^2,\ b^2\}$, for any $a$, $b>0$.}
we find 
\beqa{phi'-phi''}
|\varphi'(\e,r,\bar\g;I'',p'',q'')-\varphi''|&=&\Big|\partial_{I''}\Big(p''+p_{\rm e}(\e,r,\bar\g;I'')\Big)\cdot\Big(q'-q_{\rm e}(\e,r,\bar\g;I'')\Big)|\Big|_{q'=q_{\rm e}+q''}\nonumber\\
&\leq& \sup_{D(\e,r,\bar\g)}\Big|\partial_{I''}\Big(p''+p_{\rm e}(\e,r,\bar\g;I'')\Big)\cdot\Big(q'-q_{\rm e}(\e,r,\bar\g;I'')\Big)|\nonumber\\
&\leq&\tilde C\frac{\sup_{D(\e,r,\bar\g)}\Big|\Big(p''+p_{\rm e}(\e,r,\bar\g;I'')\Big)\cdot\Big(q'-q_{\rm e}(\e,r,\bar\g;I'')\Big)|}{\bar\r(r)/4}\nonumber\\
&\leq& C_{\rm T}\frac{\left(r+\frac{\e(\log{r^{-1}})^{2\bar\t+1}}{\bar\g^2}\right)r}{\frac{\bar\g}{(\log{r^{-1}})^{\bar\t+1}}}\nonumber\\
&\leq&C_{\rm T}\max\left\{\frac{r^2(\log{r^{-1}})^{\bar\t+1}}{\bar\g}\ ,\quad \frac{\e r(\log{r^{-1}})^{3\bar\t+2}}{\bar\g^3}\right\}\nonumber\\
&\leq&\frac{s_0}{12}
\eeqa
for a suitably small $\g_{\rm T}$. By (\ref{small peq}) and (\ref{phi'-phi''}), $\phi_{\rm T}$ is well put.
By construction, $g'':=g'\circ \phi_{\rm T}$ has an equilibrium point at $(p'',q'')=0$. Furthermore, by the splitting
$$g''=\bar f\circ \phi_{\rm T}+\tilde g\circ \phi_{\rm T}=f_0+f_2\circ \phi_{\rm T}+f_4\circ \phi_{\rm T}+\tilde g\circ \phi_{\rm T}$$
as $\partial^2(f_0)=0$ and $\partial^2(f_2\circ \phi_{\rm T})=\partial^2f_2=\textrm{\rm diag}(\O_1,\cdots \O_{\hat n},\O_1,\cdots \O_{\hat n})$, we find
\beqano
\sup_{{\bar\cI}_{\bar\r(r)/4}}\|\partial^2 g''|_0-\textrm{\rm diag}(\O_1,\cdots \O_{\hat n},\O_1,\cdots \O_{\hat n})\|&=&\sup_{{\bar\cI}_{\bar\r(r)/4}}\|\partial^2(f_4\circ \phi_{\rm T}+\tilde g\circ \phi_{\rm T})|_0\|\nonumber\\
&\leq& C_{\rm T}\frac{\e(\log{r^{-1}})^{2\bar\t+1}}{\bar\g^2}
\eeqano
and the claim is proved.
\vskip.1in
\noi
$\underline{\textrm{\sl Step}\ 3.}$ (``{\sl Birkhoff Theory}'')\ {\it There exists $0<r_{\rm B}<1<C_{\rm B}$ $\g_{\rm B}>0$ such that, for any $0<r<r_{\rm B}$and $$\bar\g>\g_{\rm B}\max\{\sqrt{\e}(\log{r^{-1}})^{\bar\t+1},\ r^2(\log{r^{-1}})^{\bar\t+1},\ \sqrt[3]{\e r}(\log{r^{-1}})^{\bar\t+1}\}\ ,$$   $\cH''$ is put into the form
\beqano\cH'''(I''',\varphi''',p''',q''')&=&\cH''\circ \phi_{\rm B}(I''',\varphi''',p''',q''')\nonumber\\
&=&h(I''')+\e\,g'''(I''',p''',q''')+\e r^5 f'''(I''',\varphi''',p''',q''')\ ,
\eeqano
where
\beqa{g'''}
g'''&=&\tilde f_0(\e,r,\bar\g;I''')+\sum_{1\leq i\leq \hat n}\tilde\O_i(\e,r,\bar\g;I''')\frac{{p_i'''}^2+{q_i'''}^2}{2}\nonumber\\
&+&\frac{1}{2}\sum_{1\leq i,j\leq \hat n}\tilde A_{i,j}(\e,r,\bar\g;I''')\frac{{p_i'''}^2+{q_i'''}^2}{2}\frac{{{p'''_j}}^2+{{q'''_j}}^2}{2}
\eeqa
with $$\tilde f_0(\e,r,\bar\g,\cdot)\ ,\quad \tilde\O(\e,r,\bar\g,\cdot)=(\tilde\O_1(\e,r,\bar\g,\cdot),\cdots,\tilde\O_{\hat n}(\e,r,\bar\g,\cdot)\ ,\quad \tilde A(\e,r,\bar\g,\cdot)=(\tilde A_{i,j}(\e,r,\bar\g,\cdot))$$
$\bar\g^{-2}\e\,(\log{r^{-1}})^{2\bar\t+1}$--close to $f_0$ $\O$, $A$, respectively and 
$$
\|f'''\|_{\bar\cI_{\bar\r(r)/8}\times \torus^{\bar n}_{s_0/24}\times B( r/2)_{ r/2}\times B( r/2)_{ r/2}}\leq C_{\rm B}
$$ 
The change of coordinates $\phi_{\rm B}=(I''(\e,r,\bar\g;\cdot)$, $\varphi''(\e,r,\bar\g;\cdot)$, $p''(\e,r,\bar\g;\cdot)$, $p''(\e,r,\bar\g;\cdot))$  may be chosen real--analytic on 
$$
\bar\cI_{\bar\r(r)/8}\times \torus^{\bar n}_{s_0/24}\times B( r/2)_{ r/2}\times B( r/2)_{ r/2}
$$
and the following bounds hold, uniformly on $\bar\cI_{\bar\r(r)/8}\times \torus^{\bar n}_{s_0/24}\times B( r/2)_{ r/2}\times B( r/2)_{ r/2}$:
\begin{eqnarray}\label{Birkhoff estimates}
\arr{I''(\e,r,\bar\g;\cdot)=I'''\\
|\varphi''(\e,r,\bar\g;I'',\varphi'',p'',q'')-\varphi'''|\leq C_{\rm B}\frac{\e r^2(\log{r^{-1}})^{3\t+2}}{\bar\g^3}\\
|p''(\e,r,\bar\g;I'',\varphi'',p'',q'')-p'''|\leq C_{\rm B}\frac{\e r(\log{r^{-1}})^{2\bar\t+1}}{\bar\g^2}\\
|q''(\e,r,\bar\g;I'',\varphi'',p'',q'')-q'''|\leq C_{\rm B}\frac{\e r(\log{r^{-1}})^{2\bar\t+1}}{\bar\g^2}\\
}\nonumber\\
\end{eqnarray}}

\vskip.1in
\noi
{\sl Proof.}\ For small values of the number $\bar\g^{-2}\e(\log{r^{-1}})^{2\bar\t+1}$, the eigenvalues of $\partial^2 g''|_0$ are purely imaginary, $\bar\g^{-2}\e(\log{r^{-1}})^{2\bar\t+1}$--close to $(\O_1,\cdots,\O_{\hat n}, \O_1,\cdots,\O_{\hat n})$
\footnote{We can always assume that such eigenvalues are pairwise equal, otherwise we perform the change of variables
$$I''={\rm I}''\ ,\quad \varphi''=\Phi''_j-\frac{\partial_{{\rm I}_j}t_j}{t_j}P''\cdot Q''\ ,\quad p''_j=t_j P''_j\ ,\quad q''_j=\frac{1}{t_j}Q''_j$$
with 
$$t_j=\sqrt[4]{\frac{\tilde\O_{q_j}}{\tilde\O_{p_j}}}=1+O\left(\frac{\e(\log{r^{-1}})^{2\bar\t+1}}{\bar\g^2}\right)$$
if {\rm diag}$(\partial^2 g''|_0)=(\tilde\O_{p_1},\cdots,\tilde\O_{p_{\hat n}},\tilde\O_{q_1},\cdots,\tilde\O_{q_{\hat n}})$. Such a transformation, generated by $S_{\rm eig}={\rm I}''\cdot \varphi''+\sum_j t_jP''_jq''_j$ does not change the final estimate (\ref{Birkhoff estimates}).
}
hence, $4$--non resonant on $\bar\cI_{\bar\r(r)/4}$. Then after a suitable ``symplectic diagonalization'' 
$$\phi_{\rm D}=(I''(\e,r,\bar\g;\cdot),\varphi''(\e,r,\bar\g;\cdot),p''(\e,r,\bar\g;\cdot),q''(\e,r,\bar\g;\cdot))$$
such that $I''(\e,r,\bar\g;\cdot)=I_{\rm D}$, $(p''-p_{\rm D},q''-q_{\rm D})$ is $\e r\bar\g^{-2}(\log{r^{-1}})^{2\bar\t+1}$--close to the identity, $(\varphi''-\varphi_{\rm D}$ $\e r^2\bar\g^{-2}(\log{r^{-1}})^{2\bar\t+1}$--close to the identity\footnote{
We may take $\phi_{\rm D}$ as generated by $S_{\rm D}=I_{\rm D}\cdot\varphi''+s_{\rm D}(I_{\rm D},p_{\rm D},q'')$, where $s_{\rm D}(I_{\rm D},p_{\rm D},q'')$ is a polynomial of degree $2$ in $(p_{\rm D},q'')$, the coefficients of which are of order $\e\bar\g^{-2}(\log{r^{-1}})^{2\bar\t+1}$. 
}
which sends $g''$ to 
$$g_{\rm D}=g''_0+\sum_{1\leq i\leq \hat n}\tilde\O(\e,r,\bar\g;I_{\rm D})\frac{p_{\rm D i}^2+q_{\rm D i}^2}{2}+o_2\ ,$$
we may apply  Birkhoff Theory (Appendix \ref{app:Birkhoff Normal Form}), putting $g_{\rm D}$ into Birkhoff Normal Form
$$g''\circ\phi_{\rm B}=g'''+\tilde o_4$$
where $g'''$ is as in (\ref{g'''}) and $|\tilde o_4|\leq |(p''',q''')|^5\leq C_{\rm B}\,r^5$,
 by means of a real--analytic symplectomorphism $$\phi_{\rm B}(\e,r,\bar\g;\cdot)=(I_{\rm D}(\e,r,\bar\g;\cdot),\varphi_{\rm D}(\e,r,\bar\g;\cdot),p_{\rm D}(\e,r,\bar\g;\cdot),q_{\rm D}(\e,r,\bar\g;\cdot)$$ such that $I_{\rm D}=I'''$, $|p_{\rm D}(\e,r,\bar\g,\cdot)-p'''|$, $|q''(\e,r,\bar\g,\cdot)-q'''|$ is $\bar\g^{-2}\e r^2(\log{r^{-1}})^{2\bar\t+1}$--close to the identity  and $|\varphi_{\rm D}(\e,r,\bar\g,\cdot)-\varphi'''|$ is ${\e r^3(\log{r^{-1}})^{3\t+2}}/{\bar\g^3}$--close to the identity
\footnote{
We may obtain (see Appendix \ref{app:Birkhoff Normal Form} for details) $\phi_{\rm B}$ in $2$ steps (which reduce the diagonalized $g_{\rm D}$ in Birkhoff Normal form of order $3$, $4$), the first of which generated by $\tilde I\cdot\varphi_{\rm D}+\sum_{|\a|+|\b|=3}\tilde s_{\a,\b}(\tilde I)\tilde p^{\a}q_{\rm D}^{\b}$, with $\tilde s^{i,j}$ $\e\bar\g^{-2}(\log{r^{-1}})^{2\bar\t+1}$--close to $0$; the second one by $I'''\cdot\tilde\varphi+\sum_{|\a|+|\b|=4}\check s_{\a,\b}( I'''){p'''}^{\a}\tilde q^{\b}$ and $\check s^{i,j}$ $\e\bar\g^{-2}(\log{r^{-1}})^{2\bar\t+1}$--close to $0$.  Apply then Cauchy estimate.
}.

\vskip.1in
\noindent
In the following step, we introduce the symplectic polar coordinates. In order to do that, we must stay away from the singularities of these coordinates at $(p'''_i,q'''_i)=0$. So, following \cite{POSCH82}, we introduce a minimum radius $r_m$ for  $(p''',q''')$ and later on we will estimate the measure of the descarted zone.\\
For $0<r_1<r_2$, denote 
$$A^{p}({r_1,r_2}):=\{x\in \real^{p}:\ r_1\leq|x|\leq r_2\}$$
the real closed anulus with radii  $r_1$, $r_2$.

\vskip.1in
\noindent
{$\underline{\textrm{\sl Step\ 4.}}$ ({\sl ``the symplectic polar coordinates''}) }\ {\it 
There exist $C_{\rm pc}$, $s>0$ such that, for any fixed $r_m>0$, 
the symplectic (``polar coordinates'') transformation $\phi_{\rm pc}$ defined on the domain
$$\bar\cI_{\bar\r(r)/8}\times A^{\hat n}\left({r_m^2},r^2\right)_{\hat\r(r)}\times \torus_s^{\bar n}\times \torus_s^{\hat n}\ ,\quad \hat\r(r):=\min\left\{\frac{r_m^2}{2},\ r^2\right\}$$
by
\begin{eqnarray}\label{polar}
\phi_{\rm pc}:\ \left\{
\begin{array}{lrr}
I'''=\bar J\\
{\varphi'''}=\bar \psi\\
p'''=\sqrt{2\,\hat J}\cos{\hat\psi}\\
q'''=\sqrt{2\,\hat J}\sin{\hat\psi}\\
\end{array}\right.
\end{eqnarray}
puts
 $\cH'''(\e,r,\bar\g,\cdot)$  into the form
\beqano
{\rm H}\Big(\e,r,\bar\g;(\bar J,\hat J), (\bar\psi,\hat\psi)\Big)&=&\cH'''\circ\phi_{\rm pc}\nonumber\\
&=&{h}(\bar J)+\e{\rm h}_1(\e,r,\bar\g;(\bar J,\hat J))+\e r^5{\rm f}\Big(\e,r,\bar\g;(\bar J,\hat J),(\bar\psi,\hat\psi)\Big)
\eeqano
where 
$$
\arr{
{\rm h}_1(\e,r,\bar\g;(\bar J,\hat J))=\tilde f_0(\e,r,\bar\g;\bar J)+\tilde\O(\e,r,\bar\g;\bar J)\cdot \hat J+\frac{1}{2}\hat J\cdot \tilde A(\e,r,\bar\g;\bar J)\hat J)
\\
\|{\rm f}\|_{\bar\cI_{\bar\r(r)/8}\times A(r_m^2,cr^2)_{\hat\r(r)}\times \torus_s^{\bar n}\times \torus^{\hat n}_s}\leq C_{\rm pc}
}
$$
}

\vskip.1in
\noi
{\sl Proof.}\ Obvious.
\begin{remark}\label{reduction remark}\rm
Denote by $$\phi_{\rm red}(\e,r,\bar\g;\cdot,\cdot)=(\phi_{\rm red,I}(\e,r,\bar\g;\cdot,\cdot),\phi_{\rm red,\varphi}(\e,r,\bar\g;\cdot,\cdot),\phi_{\rm red,p}(\e,r,\bar\g;\cdot,\cdot),\phi_{\rm red,q}(\e,r,\bar\g;\cdot,\cdot))$$ the composition of the real--analytic symplectomorphisms described in steps $1\div 4$. Then, by the estimates (\ref{average estimates}), (\ref{id2}), (\ref{Birkhoff estimates}), we may let
\beqano
\arr{
\phi_{\rm red,I}\Big(\e,r,\bar\g;(\bar J,\hat J),(\bar\psi,\hat\psi)\Big)=\bar J+a \Big(\e,r,\bar\g;(\bar J,\hat J),(\bar\psi,\hat\psi)\Big)\\
\phi_{\rm red,\varphi}\Big(\e,r,\bar\g;(\bar J,\hat J),(\bar\psi,\hat\psi)\Big)=\bar\psi+ b\Big(\e,r,\bar\g;(\bar J,\hat J),(\bar\psi,\hat\psi)\Big)\\
\phi_{\rm red,p}\Big(\e,r,\bar\g;(\bar J,\hat J),(\bar\psi,\hat\psi)\Big)=\sqrt{2\hat J}\cos{\hat\psi}+u\Big(\e,r,\bar\g;(\bar J,\hat J),(\bar\psi,\hat\psi)\Big)\\
\phi_{\rm red,q}\Big(\e,r,\bar\g;(\bar J,\hat J),(\bar\psi,\hat\psi)\Big)=\sqrt{2\hat J}\sin{\hat\psi}+v\Big(\e,r,\bar\g;(\bar J,\hat J),(\bar\psi,\hat\psi)\Big)\\
}
\eeqano
where the functions $a(\e,r,\bar\g,\cdot,\cdot)$, $b(\e,r,\bar\g,\cdot,\cdot)$, $u(\e,r,\bar\g,\cdot,\cdot)$, $v(\e,r,\bar\g,\cdot,\cdot)$ satisfy, uniformly on $\bar\cI_{\bar\r(r)/8}\times A(r_m^2,r^2)_{\hat\r(r)}\times \torus_s^{\bar n}\times \torus_s^{\hat n}$,
\beqano
\arr{
|a(\e,r,\bar\g,\cdot,\cdot)|\leq C\frac{\e(\log{r^{-1}})^{\bar\t}}{\bar\g}\\
|b(\e,r,\bar\g,\cdot,\cdot)|\leq C\max\left\{\frac{\e(\log{r^{-1}})^{2\bar\t+1}}{\bar\g^2}\ ,\quad \frac{r^2(\log{r^{-1}})^{2\bar\t+1}}{\bar\g}\ ,\quad \frac{\e r(\log{r^{-1}})^{3\bar\t+2}}{\bar\g^3}\right\}\\
|u|\leq C\max\left\{\frac{\e(\log{r^{-1}})^{\bar\t}}{\bar\g}\ ,\quad \frac{\e(\log{r^{-1}})^{2\bar\t+1}}{\bar\g^2}\right\}\\
|v|\leq C\max\left\{\frac{\e(\log{r^{-1}})^{\bar\t}}{\bar\g}\ ,\quad \frac{\e(\log{r^{-1}})^{2\bar\t+1}}{\bar\g^2}\right\}\\
}
\eeqano
\end{remark}

\vskip.1in
\noi
$\underline{\textrm{\sl Step}\ 5.}$ (``{\sl KAM}'')\ {\it Let 
$$\bar\cJ:=\bar\cI\times A^{\hat n}\left({r_m^2},r^2\right)$$
\beq{good rho}\r:=\min\left\{\frac{\bar\r(r)}{16},\ \hat\r(r)\right\}=\min\left\{{r^2}\ ,\quad \frac{r_m^2}{2}\ ,\quad  c_*\frac{\bar\g}{(\log{r^{-1}})^{\bar\t+1}}\ ,\quad \r_0\right\}\eeq
There exists $r_{\rm KAM}$, $\g_{\rm KAM}$, $c_{\rm KAM}$ such that, 
for any $0< r<r_{\rm KAM}$, $r_m\geq r_{\rm KAM}^{-1} r^{5/4}$
and
$$\arr{\g>\g_{\rm KAM} r^{5/2}\\
\hat\g_\e>\g_{\rm KAM}\e r^{5/2}(\log_+{(r^5/\g^2)^{-1}})^{\t+1}\\
\bar\g>\g_{\rm KAM}\max\{\sqrt{\e}(\log{r^{-1}})^{\bar\t+1},\ 
\sqrt[3]{\e r}(\log{r^{-1}})^{\bar\t+1},\ r^2(\log{r^{-1}})^{\bar\t+1}\}\\
}
$$
then the Hamiltonian ${\rm H}$ satisfies the assumptions of Theorem \ref{two scales KAM} on the domain $\bar\cJ_\r\times \torus^{n}_s$.

}

\vskip.1in
\noi
$\underline{\textrm{\sl Proof.}}$

\vskip.1in
\noi 
$\underline{\textrm{\sl Claim $1$:}}$ {\sl $\o=\partial{\rm h}$ is a diffeomorphism of $\bar\cJ_\r$.}

\vskip.1in
\noi
$\underline{\textrm{\sl proof of claim $1$:}}$ Due to the analyticity assumptions, we have only to prove the injectivity for $\o$. We prove that, for any $\n=(\bar\n,\e\hat\n)\in \o(\e,r,\bar\g;\cI_{\r})$, equation $\o(\e,r,\bar\g;(\bar J,\hat J))=\n$
has at most one solution on $(\bar J,\hat J)\in\cI_\r$. Let, then, 
\beqano
\arr{\bar\partial\,h+\e(\bar\partial \tilde f_0+\bar\partial\tilde\O\cdot \hat J+\frac{1}{2}\hat J\cdot \bar\partial\tilde A\hat J)=\bar\n\\
\tilde\O+\tilde A\hat J=\hat\n
}
\eeqano
where $\bar\partial$, $\hat\partial$ denote, respectively, differentiation with respect to $\bar J$, $\hat J$.
For any fixed $\bar J\in \bar\cI$, the map
$$\hat J\to \tilde\O+\tilde A\hat J$$ is injective (as $\tilde A(\e,r,\bar\g;\bar J)$ is nonsingular): $\bar J$, we find a unique $$\hat J=\hat J_0(\e,r,\bar\g,\hat\n,\bar J):=\tilde A(\e,r,\bar\g,\bar J)^{-1}(\hat\n-\tilde\O(\e,r,\bar\g,\bar J)$$ solving the second equation.
Replacing this value into the equation for the first components, we find an equation of the kind
\beqano
&& \o_0(\bar J)+\o_1(\e,r,\bar\g,\hat\n;\bar J))=\bar\n
\eeqano
where $\o_0=\bar\partial h$ is well defined and analytic up to $\bar\cI_{\r_0}$, hence, with $\|(\bar\partial\o_0)^{-1}\|$ uniformly bounded on $\bar\cI_\r$ by a suitable constant $N_0$ (which does not depend on $(\e,r,\bar\g)$) and
\beqano
\o_1&=&\e\left(\partial \tilde f_0+\bar\partial\tilde\O(\bar J)\cdot \hat J_0(\e,r,\bar\g,\hat\n,\bar J)+\frac{1}{2}\hat J_0(\e,r,\bar\g,\hat\n,\bar J)\cdot\bar\partial\tilde A \hat J_0(\e,r,\bar\g,\hat\n,\bar J)\right)\nonumber\\
\eeqano well defined and analytic up to $\bar\cI_{\bar\r(r)/8}$, hence, by Cauchy estimates,
$$\sup_{\bar\cI_{\bar\r(r)/16}} \|\bar\partial\o_1\|\leq C'\frac{\e}{\bar\r(r)^2}\leq C\frac{\e(\log{r^{-1}})^{2(\bar\t+1)}}{\bar\g^2}\leq \frac{1}{2N_0}$$
(recall that $\tilde f_0$, $\tilde\O$, $\tilde A$ are $\e(\log{r^{-1}})^{2\bar\t+1}/\bar\g^2$--close to $f_0$, $\O$, $A$). This proves the claim.

\vskip.1in
\noi 
$\underline{\textrm{\sl Claim $2$:}}$ {\sl $\o=\partial{\rm h}$ has non singular Jacobian on $\bar\cJ_\r$ and $T:=(\partial^2{\rm h})^{-1}$ is bounded by
$$\sup_{{\bar\cJ}_{\r}}\|T\|\leq C\e^{-1}$$}

\vskip.1in
\noi
$\underline{\textrm{\sl proof of claim $2$:}}$ We have
\beqa{full hessian}
\partial\o&=&\partial^2{\rm h}\nonumber\\
&=&\left(
\begin{array}{lrr}
\partial^2 h(\bar J)& \e B&\\
\e B^T& \e \tilde A(\e,r,\bar\g,\bar J)&
\end{array}
\right)
\eeqa
where $B=(B_{i,j})$ is the $\bar n\times \hat n$ matrix with elements
$$B_{i,j}=\partial_{\bar J_i}\tilde\O_j+(\partial_{\bar J_i}\tilde A \hat J)_j$$
bounded in norm by (by Cauchy estimates)
\beq{BB}
\sup_{\cI_\r}\|B\|\leq \sup_{\cI_{\bar\r(r)}}\|B\|\leq\frac{C'}{\bar\r(r)}\leq  C\frac{\log{r^{-1}})^{\bar\t+1}}{\bar\g}
\eeq
and, for small ${\bar\g^{-2}}{\e(\log{r^{-1}})^{2\bar\t+1}}$ and 
$$\|\tilde A^{-1}\|\leq 2\|A^{-1}\|$$ (say).
We prove that the matrix 
$$\cM:=\left(
\begin{array}{lrr}
\partial^2 h(\bar J)& \e B&\\
B^T& \tilde A(\e,r,\bar\g,\bar J)&
\end{array}
\right)$$
is non singular
which will imply the claim, as, by (\ref{full hessian}),
$${\rm det}(\partial^2{\rm h})=\e^{\hat n}{\rm det}\cM\ .$$
We split $\cM$ as
$$\cM=\cM_0+\cM_1\ ,$$
where
$$\cM_0=\left(
\begin{array}{lrr}
\partial^2 h(\bar J)& 0&\\
B^T& \tilde A(\e,r,\bar\g,\bar J)&
\end{array}
\right)\ ,\qquad \cM_1=\left(
\begin{array}{lrr}
0& \e B&\\
0& 0&
\end{array}
\right)$$
The matrix $\cM_0$ is non singular (because $\partial^2 h$, $\tilde A$ are so), with inverse matrix
$$\cM_0^{-1}=\left(
\begin{array}{lrr}
(\partial^2 h(\bar J))^{-1}& 0&\\
-(\tilde A(\e,r,\bar\g,\bar J))^{-1}B^T(\partial^2 h(\bar J))^{-1}& (\tilde A(\e,r,\bar\g,\bar J))^{-1}&
\end{array}
\right)$$
Furthermore, the matrix
$$\d\cM:=\cM_0^{-1}\cM_1=\e\left(
\begin{array}{lrr}
0&(\partial^2 h(\bar J))^{-1}B&\\
0& -(\tilde A(\e,r,\bar\g,\bar J))^{-1}B^T(\partial^2 h(\bar J))^{-1}B&
\end{array}
\right)$$
has norm bounded by
\beqano
\sup_{\bar\cJ_\r}\|\d\cM\|&\leq& 2\e\max\Big\{\sup_{\bar\cJ_\r}\|(\partial^2 h)^{-1}\|\sup_{\bar\cJ_\r}\|B\|,\nonumber\\
&&\sup_{\bar\cJ_\r}\|\tilde A^{-1}\|\sup_{\bar\cJ_\r}\|\partial^2 h^{-1}\|(\sup_{\bar\cJ_\r}\|B\|)^2\Big\}\nonumber\\
&\leq&C\frac{\e(\log{r^{-1}})^{2(\bar\t+1)}}{\bar\g^2}\nonumber\\
&\leq& \frac{1}{2}
\eeqano
provided $\bar\g^{-2}{\e(\log{r^{-1}})^{2(\bar\t+1)}}$ is small.
This makes the matrix
$$\cM=\cM_0(\id+\d\cM)$$
invertible. Finally, by (\ref{full hessian}), it is clear that 
$$\|T\|\leq \e^{-1}\|\cM^{-1}\|\leq C\e^{-1}\ .$$

\vskip.1in
\noi
$\underline{\textrm{\sl Claim $3$:}}$ (``{\sl check of the KAM condition}''){\sl There exist $r_{KAM}$, $\g_{\rm KAM}$, $c_{KAM}>0$ such that, for any $0<r<r_{KAM}$ and $r_m\geq r_{KAM}^{-1}r^{5/4}$ and  
\beqa{all gammas}
\arr{\g>\g_{\rm KAM} r^{5/2}\\
\hat\g_\e>\g_{\rm KAM}\e r^{5/2}(\log_+{(r^5/\g^2)^{-1}})^{\t+1}\\
\bar\g>\g_{\rm KAM}\max\{\sqrt{\e}(\log{r^{-1}})^{\bar\t+1},\sqrt[3]{\e r}(\log{r^{-1}})^{\bar\t+1},\ r^2(\log{r^{-1}})^{\bar\t+1}\}\\
}
\eeqa
there exist 
\beqano
&&  \dst M\geq\sup_{{\cal J}_{\r}}\|U\|\ ,\quad \hat M\geq\sup_{{\cal J}_{\r}}\|U^{[n,\hat n]}\|\ ,\quad  N\geq \sup_{{\cal J}_{\r}}\|T\|\ ,\quad F\geq\|\e r^{5}{ f}(\e,r,\bar\g,\cdot,\cdot)\|_{\r,s}
\eeqano
such that, 
letting
\beqano
&&  c:=\max\left\{{2^{11}n},\ \frac{2}{3}(12)^{\t+1}\right\}\ ,\ \dst \tilde{\r}:=\min\left\{\frac{\g}{3MK^{\t+1}}\ ,\ \frac{\hat\g_\e}{3\hat MK^{\t+1}}\ ,\  
 \r\right\}\ ,\nonumber\\
&&  L:=\max\ \{N\ , M^{-1},\hat M^{-1}\}\ ,\quad K:=\frac{6}{s}\ \log_+{\left(\frac{FM^2\,L}{\gamma^2}\right)^{-1}}\\
\eeqano
where $\r$ is as in (\ref{good rho}), then,
\beq
cE:=c\frac{FL}{\tilde{\r}^2}<1\ .
\eeq}

\vskip.1in
\noi
$\underline{\textrm{\sl proof of claim $3$:}}$
To apply Theorem \ref{two scales KAM} to our case, we let, for suitable $c_+>1>c_-$, $$F=\e\,r^5 c_+\ ,\quad M=c_+\ ,\quad \hat M=\e\,c_+\ ,\quad N=c_+\e^{-1}\ ,$$
so that $$L=c_+\e^{-1}\ ,\quad K\leq c_+(\log{(r^5/\g^2)^{-1}})\ .$$
Taking $\hat\g_\e=\e\hat\g$, we find
\beqa{tilderho fund}
\tilde\r&=&\min\left\{\frac{\g}{3MK^{\t+1}}\ ,\ \frac{\hat\g_\e}{3\hat MK^{\t+1}}\ ,\  
 \r\right\}\nonumber\\
 &\geq& c_-\,\min\left\{\frac{\g}{(\log_+{(r^5/\g^2)^{-1}})^{\t+1}}\ ,\quad \frac{\hat\g}{(\log_+{(r^5/\g^2)^{-1}})^{\t+1}}\ ,\right.\nonumber\\
 & & \left.\frac{\bar\g}{(\log{r^{-1}})^{\bar\t+1}}\ ,\quad { \frac{r_m^2}{2}}\ ,
 \quad r^2\ ,\quad \r_0\right\}
\eeqa
and hence, the smallness or ``KAM'' condition
\beqano
cE&=&c\frac{FL}{\tilde\r^2}\nonumber\\
&\leq& c_+ \max\left\{r^5\frac{(\log_+{(r^5/\g^2)^{-1}})^{2(\t+1)}}{\g^2}\ ,\quad \frac{r^5(\log_+{(r^5/\g^2)^{-1}})^{2(\t+1)}}{\hat\g^2}\ ,\right.\nonumber\\
&&\left. \frac{r^5(\log{r^{-1}})^{2(\bar\t+1)}}{\bar\g^2}\ ,\quad r\ ,\quad {\frac{r^5}{r_m^4}}\ ,\quad 
\frac{r^5}{\r_0}\right\}\nonumber\\
&<&1
\eeqano
is fulfilled whenever we choose $r<r_{\rm KAM}$, then $\bar\g$, $\g$, $\hat\g_\e$  as in (\ref{all gammas}) and finally $r_m$ not less than $r_{\rm KAM}^{-1}r^{5/4}$, for a suitable small $r_{\rm KAM}$.

\vskip.1in
\noi
$\underline{\textrm{\sl Conclusion of the proof.}}$

\vskip.1in
\noi
Define
\beqano
\arr{
{\eufm J}_*(\e,r,\bar\g,\g,\hat\g)\subset \bar\cJ=\bar\cI\times A^{\hat n}\left({r_m^2},r^2\right)\\ 
\varpi_*(\e,r,\bar\g,\g,\hat\g,\cdot):\quad {\eufm J}_*(\e,r,\bar\g,\g,\hat\g)\to {\eufm O}_*(\e,r,\bar\g,\g,\hat\g)\subset\cD^{\bar n,\hat n}_{\g,\hat\g_\e,\t}\\
\phi(\e,r,\bar\g,\g,\hat\g,\cdot,\cdot):\quad {\eufm O}_*\times \torus^n\to \bar\cJ\times \torus^n\\
{\rm K}_*(\e,r,\bar\g,\g,\hat\g)\subset \bar\cJ\times\torus^n\\
}
\eeqano
as the Cantor set, the Lipschitz homeomorphism onto, the tori embedding and the invariant set which are obtained by  Theorem \ref{two scales KAM}. Define
 ${\eufm F}(\e,r,\bar\g,\g,\hat\g,\cdot,\cdot)$ as the image of $\phi(\e,r,\bar\g,\g,\hat\g,\cdot,\cdot)$ under $\phi_{\rm red}(\e,r,\bar\g;\cdot)$; ${\eufm U}(\e,r,\bar\g)\supset{\eufm K}(\e,r,\bar\g,\g,\hat\g)$,  the images of $\bar\cJ\times \torus^n\supset{\rm K}_*(\e,r,\bar\g,\g,\hat\g)$,  under $\phi_{\rm red}(\e,r,\bar\g;\cdot)$. 

\vskip.1in
\noi
$\underline{\textrm{\sl Proof of (\ref{bounds for varpi}).}}$\ Put $$\bar\varpi_*(\e,r,\bar\g,\g,\hat\g,\cdot):=\Big(\bar\varpi_*(\e,r,\bar\g,\g,\hat\g,\cdot),\ \e\hat\varpi_*(\e,r,\bar\g,\g,\hat\g,\cdot)\Big)\ .$$
By (\ref{frequency map}), 
\beqano
|\bar\varpi_*-\partial h|&\leq&|\bar\varpi_*-\partial{\rm h}|+C'\e\nonumber\\
&\leq&C''\tilde\r+C'\e\nonumber\\
&\leq&C\min\left\{\frac{\g}{(\log_+{(r^5/\g^2)^{-1}})^{\t+1}}\ ,\quad \frac{\hat\g}{(\log_+{(r^5/\g^2)^{-1}})^{\t+1}}\ ,\quad\frac{\bar\g}{(\log{r^{-1}})^{\bar\t+1}}\right\}\nonumber\\
&+&C\e
\eeqano
having used (\ref{tilderho fund}), for which
\beqano
\tilde\r&\leq& \tilde C\,\min\left\{\frac{\g}{(\log_+{(r^5/\g^2)^{-1}})^{\t+1}}\ ,\quad \frac{\hat\g}{(\log_+{(r^5/\g^2)^{-1}})^{\t+1}}\ ,\quad\frac{\bar\g}{(\log{r^{-1}})^{\bar\t+1}}\ ,\quad r^{5/2}\right\}
\eeqano
(recall $r_m^2=\const r^{5/2}$).
Similarly,
\beqano
|\e\hat\varpi_*-\e\O|&\leq&|\e\hat\varpi_*-\e\tilde\O-\e \tilde A\hat J|+|\e\tilde\O-\e\O|+|\e (\tilde A-A)\hat J|+|\e A\hat J|\nonumber\\
&\leq&C\e\min\left\{\frac{\g}{(\log_+{(r^5/\g^2)^{-1}})^{\t+1}}\ ,\quad \frac{\hat\g}{(\log_+{(r^5/\g^2)^{-1}})^{\t+1}}\ ,\quad\frac{\bar\g}{(\log{r^{-1}})^{\bar\t+1}}\right\}\nonumber\\
&+&C\e^2\frac{(\log{r^{-1}})^{2\bar\t+1}}{\bar\g^2}+C\e^2 r^2\frac{(\log{r^{-1}})^{2\bar\t+1}}{\bar\g^2}+C\e r^2
\eeqano

\vskip.1in
\noi
$\underline{\textrm{\sl Proof of (\ref{bounds for FF}).}}$\ By Theorem \ref{two scales KAM}, we find
\beqano
&& \phi(\e,r,\bar\g,\g,\hat\g,\cdot,\cdot)\nonumber\\
&& =\Big(\phi_{\bar J}(\e,r,\bar\g,\g,\hat\g,\cdot,\cdot),\phi_{\hat J}(\e,r,\bar\g,\g,\hat\g,\cdot,\cdot),\phi_{\bar\psi}(\e,r,\bar\g,\g,\hat\g,\cdot,\cdot),\phi_{\hat\psi}(\e,r,\bar\g,\g,\hat\g,\cdot,\cdot)\Big)\nonumber\\
&& =\arr{
\phi_{\bar J}(\e,r,\bar\g,\g,\hat\g,\cdot,\cdot)=\bar{\eufm j}_*(\e,r,\bar\g,\g,\hat\g;\n)+\bar U(\e,r,\bar\g,\g,\hat\g;(\bar\vartheta,\hat\vartheta))\\
\phi_{\hat J}(\e,r,\bar\g,\g,\hat\g,\cdot,\cdot)=\hat{\eufm j}_*(\e,r,\bar\g,\g,\hat\g;\n)+\hat U(\e,r,\bar\g,\g,\hat\g;(\bar\vartheta,\hat\vartheta))\\
\phi_{\bar \psi}(\e,r,\bar\g,\g,\hat\g,\cdot,\cdot)=\bar\vartheta+\bar V(\e,r,\bar\g,\g,\hat\g;(\bar\vartheta,\hat\vartheta))\\
\phi_{\hat \psi}(\e,r,\bar\g,\g,\hat\g,\cdot,\cdot)=\hat\vartheta+\hat V(\e,r,\bar\g,\g,\hat\g;(\bar\vartheta,\hat\vartheta))\\
}
\eeqano
with
$$
\arr{
|\bar U|\leq \check C\frac{\hat\g_\e}{\g}\tilde\r\leq C\frac{\hat\g_\e}{\g}\min\left\{\frac{\g}{(\log_+{(r^5/\g^2)^{-1}})^{\t+1}}\ ,\quad \frac{\hat\g}{(\log_+{(r^5/\g^2)^{-1}})^{\t+1}}\ ,\quad\frac{\bar\g}{(\log{r^{-1}})^{\bar\t+1}}\ ,\quad r^{5/2}\right\}\\
|\hat U|\leq \check C\tilde\r\leq C\min\left\{\frac{\g}{(\log_+{(r^5/\g^2)^{-1}})^{\t+1}}\ ,\quad \frac{\hat\g}{(\log_+{(r^5/\g^2)^{-1}})^{\t+1}}\ ,\quad\frac{\bar\g}{(\log{r^{-1}})^{\bar\t+1}}\ ,\quad r^{5/2}\right\}\\
|\bar V|, |\hat V|\leq \max\left\{r^5\frac{(\log_+{(r^5/\g^2)^{-1}})^{2(\t+1)}}{\g^2}\ ,\quad \frac{r^5(\log_+{(r^5/\g^2)^{-1}})^{2(\t+1)}}{\hat\g^2}\ ,\quad \frac{r^5(\log{r^{-1}})^{2(\bar\t+1)}}{\bar\g^2}\ ,\quad r\right\}\\
}
$$
Hence, recalling Remark \ref{reduction remark}, the estimates (\ref{bounds for FF}) follow for ${\eufm F}(\e,\bar\g,\g,\hat\g;\cdot,\cdot)$

\vskip.1in
\noi
$\underline{\textrm{\sl Step $6$: proof of (\ref{fund measure}).}}$\  Let $$\cJ:=\cI\times A^{\hat n}\left(r_m^2,\ r^2\right)\supset \bar\cI\times A^{\hat n}\left(r_m^2 ,\ r^2\right)=\bar\cJ\ ,\quad r_m\geq r_* r^{5/4}\ .$$ 
We first prove that 
$$\meas\Big(\cJ\times\torus^n\setminus{\rm K}_*(\e,r,\bar\g,\g,\hat\g)\Big)\leq C_*\Big(\bar\g+\hat\g_\e+\frac{\hat\g}{r^2}\Big)\meas\Big(\cJ\times\torus^n\Big)\ .$$
The set $${\eufm K}(\e,r,\bar\g,\g,\hat\g)$$ is measure--equivalent to $${\rm K}_*(\e,r,\bar\g,\g,\hat\g)$$  as $\phi_{\rm red}(\e,r,\bar\g;\cdot)$ is, in particular, a real symplectomorphism, hence, area--preserving. The density of $\bar\cJ\times \torus^n\setminus{\rm K}_*(\e,r,\bar\g,\g,\hat\g)$ in $\bar\cJ\times \torus^n$ is estimated by (\ref{tori measure}) of Theorem \ref{two scales KAM}:
\beqa{meas of invariant tori}
\meas\Big(\bar\cJ\times\torus^n\setminus{\rm K}_*(\e,r,\bar\g,\g,\hat\g)\Big)&\leq& \bar c\Big(\meas(\bar\cJ\setminus\bar\cJ_{\g,\hat\g_\e,\t}\times\torus^n)+\meas(\bar\cJ_{\r_{\rm max}}\setminus\bar\cJ\times\torus^n)\Big)\nonumber\\
\eeqa
where (recall (\ref{tilderho fund}))
$$\r_{\max}:=\max\{\r_1,\ \r_2\}\leq 2^{-5}\tilde\r\leq C\min\{\g,\hat\g\}\ .$$
The second term is easily bounded by
\footnote{We are using that $\cI$--being an open and bounded set of $\real^{\bar n}$--satisfies the following: there exists $D=D(\cI)>0$, $\bar\r=\bar\r(\cI)$ such that, for any $0<\r<\bar\r$,
\beqa{D property}\meas(\cI_\r\setminus\cI)\leq \frac{\r}{D(\cI)}\meas\cI\ .\eeqa
Then, $\cJ$ is a product $\cJ={\eufm A}\times{\eufm B}$ where both ${\eufm A}=\cI$ and ${\eufm B}=A^{\hat n}\left(r_m^2,r^2\right)$ have the property (\ref{D property}), with $D({\eufm A})=D_0$, $D({\eufm B})=D_0\,r^2$. So, using
$$({\eufm A}\times{\eufm B})_\r\setminus({\eufm A}\times{\eufm B})=({\eufm A}_{\r}\setminus{\eufm A})\times{\eufm B}\bigcup{\eufm A}\times({\eufm B}_{\r}\setminus{\eufm B})\bigcup ({\eufm A}_{\r}\setminus{\eufm A})\times ({\eufm B}_{\r}\setminus{\eufm B})$$
we find
\beqano
\meas\Big(\cJ_\r\setminus\cJ\Big)&\leq&\meas({\eufm A}_{\r}\setminus{\eufm A})\times{\eufm B}+\meas{\eufm A}\times({\eufm B}_{\r}\setminus{\eufm B})+\meas ({\eufm A}_{\r}\setminus{\eufm A})\times ({\eufm B}_{\r}\setminus{\eufm B})\nonumber\\
&\leq&C\left(\r+\frac{\r}{r^2}+\frac{\r^2}{r^2}\right)\meas(\cJ)
\eeqano
}
\beqano
\meas(\bar\cJ_{\r_{\rm max}}\setminus\bar\cJ\times\torus^n)&\leq&\meas(\cJ_{\r_{\rm max}}\setminus\bar\cJ\times\torus^n)\nonumber\\
&\leq&\meas(\cJ_{\r_{\rm max}}\setminus\cJ\times\torus^n)+\meas(\cJ\setminus\bar\cJ\times\torus^n)\nonumber\\
&\leq&\check c\max\left\{\frac{\hat\g}{r^2},\ \g\right\}\meas(\cJ\times\torus^n)+\meas(\cJ\setminus\bar\cJ\times\torus^n)
\eeqano
Inserting this bound into (\ref{meas of invariant tori}), we get
\beqa{measure estimates1}
\meas\Big(\bar\cJ\times\torus^n\setminus{\rm K}_*(\e,r,\bar\g,\g,\hat\g)\Big)&\leq&\bar c\left(\meas(\bar\cJ\setminus\bar\cJ_{\g,\hat\g_\e,\t}\times\torus^n)\right.\nonumber\\
&+&\left.\meas(\cJ\setminus\bar\cJ\times\torus^n)\right.\nonumber\\
&+&\left.\max\left\{\frac{\hat\g}{r^2},\ \g\right\}\meas(\cJ\times\torus^n)\right)\nonumber\\
\eeqa
Now, recalling that
$$\bar\cJ=\bar\cI\times A^{\hat n}\left(r_m^2,\ r^2\right)\quad \textrm{with}\quad \bar\cI=\cI_{\bar\g,\bar\t}=\{I\in \cI:\ \o(I)\in \cD^{\bar n}_{\bar\g,\bar\t}\}$$
we find that first two terms inside the parentheses of (\ref{measure estimates1}) are similar, and they are simultaneouly estimated by the Lemma \ref{predioph measure lemma} below. 
\begin{lemma}\label{predioph measure lemma}
Let $\hat n$, $\bar n\in \natural$, $\t>n:=\bar n+\hat n$, $1<\a<2$, $0<\hat r<1$, $\bar\cI$ compact, , $\hat\cI:=A_{\hat r}:=A^{\hat n}\Big(\hat r^\a,\hat r\Big)$
 $$\o=(\bar\o,\hat\o):\quad \cI:=\bar\cI\times\hat\cI\to \real^{\bar n}\times\real^{\hat n}$$
a diffeomorphism of an open neighborhood of $\cI$, with $\hat\o$ of the form
$$\hat\o(\bar I,\hat I)=\hat\o_0(\bar I)+A(\bar I)\hat I$$
where $\bar I\to A(\bar I)$ is non singular on $\bar\cI$. Let
$$\bar R>\max_{\cI}|\bar\o|\ ,\quad \cA>\max_{\bar\cI}\|A\|\ ,\quad c(n,\t):=\sum_{0\neq k\in \integer^n}\frac{1}{|k|^\t}\ ,$$
and denote
$$\cR^{\bar n,\hat n}_{g,\hat g,\t}:=\Big\{I=(\bar I,\hat I)\in \cI:\quad \o(I)\notin \cD^{\bar n,\hat n}_{g,\hat g,\t}\Big\}\ .$$
Then, there exists a suitable integer number $p$ such that
$$\meas\Big(\cR^{\bar n,\hat n}_{g,\hat g,\t}\Big)\leq \left(\bar cg+\hat c\,\frac{\hat g}{\hat r}\right)\meas\Big(\cI\Big)$$
where
$$\arr{\bar c:=\sup_{\cI}\|\o^{-1}\|^n\frac{\bar R^{\bar n-1}}{\meas(\bar\cI)}\frac{\cA^{\hat n}}{\meas(A_1)}c(n,\t)p\\
\hat c=\frac{\sup_{\bar\cI}\|A^{-1}\|^{\hat n}(\cA)^{\hat n-1}}{\meas(A_1)}c(\hat n,\t)
}
$$
\end{lemma}
For continuity reasons, the proof of Lemma \ref{predioph measure lemma} is postponed at the end of the actual one.

\vskip.1in
\noi
We may then take
\beqa{measure estimates}
\arr{
\meas(\bar\cJ\setminus\bar\cJ_{\g,\hat\g_\e,\t}\times\torus^n)\leq c_*\left(\g+\frac{\hat\g}{r^2}\right)\meas\Big(\cJ\times\torus^n\Big)\\
\meas(\cJ\setminus\bar\cJ\times\torus^n)\leq c_*\bar\g\meas\Big(\cJ\times\torus^n\Big)
}
\eeqa
with $c_*$ independent of $r$. Hence, using (\ref{measure estimates}) into (\ref{measure estimates1}), we finally find
\beqa{tori measure in the anulus}
\meas\Big(\cJ\times\torus^n\setminus{\rm K}_*(\e,r,\bar\g,\g,\hat\g)\Big)&\leq& \meas\Big(\cJ\setminus\bar\cJ\times\torus^n\Big)\nonumber\\
&+&\meas\Big(\bar\cJ\times\torus^n\setminus{\rm K}_*(\e,r,\bar\g,\g,\hat\g)\Big)\nonumber\\
&\leq&C_*\left(\bar\g+\g+\frac{\hat\g}{r^2}\right)\meas\Big(\cJ\times\torus^n\Big)
\eeqa
which is quite what we meant to prove.

\vskip.1in
\noi
Having now the masure estimate (\ref{tori measure in the anulus}) and using
$$\textrm{\rm meas}\Big(V(r)\setminus\cJ\times\torus^n\Big)\leq C_*\left(\frac{r_m}{r}\right)^{2\hat n}\textrm{\rm meas}\Big(V(r)\Big)\leq C_*r^{\hat n/2}\textrm{\rm meas}\Big(V(r)\Big)$$
(eventually with a different $C_*$)
we easily infer (\ref{fund measure}).

\vskip.1in
\noi
{\bf Proof of Lemma \ref{predioph measure lemma}.}\ The first part of the proof uses a compactness argument.
Let $$\bar R>\max_\cI|\bar\o|\ ,\quad \cA>\max_{\bar \cI}\|A\|\ ,\quad \hat R:=\cA\hat r>\max_\cI|\hat\o-\hat\o_0|$$ 
so that
$${\cal U}:=\Big\{B^{\bar n}_{\bar R}(0)\times B^{\hat n}_{\hat R}(\o_0(\bar I))\ ,\quad \bar I\in \bar\cI\Big\}$$
is an open covering of $\o(\cI)$, which is compact, as continuous image of a compact. Then, there exists a finite number of $\bar I_1$, $\cdots$, $\bar I_p\in \bar\cI$ such that
$$\bar{\cal U}:=\bigcup_{1\leq i\leq p}U_i\ ,\quad U_i:=B^{\bar n}_{\bar R}(0)\times B^{\hat n}_{\hat R}(\o_0(\bar I_i))$$
covers $\o(\cI)$. Now,  the ``resonant set'' $\cR^{\bar n,\hat n}_{g,\hat g,\t}$ in $\cI$ is
\beqa{resonant set}
\cR^{\bar n,\hat n}_{g,\hat g,\t}&=&\bigcup_{k=(\bar k,\hat k)\in \integer^{\bar n}\times \integer^{\hat n},\bar k\neq 0}\left\{I:\ |\o(I)\cdot k|\leq\frac{g}{|k|^\t}\right\}\bigcup_{0\neq\hat k\in \integer^{\hat n}}\left\{I:\ |\hat\o(I)\cdot \hat k|\leq\frac{\hat g}{|\hat k|^\t}\right\}\nonumber\\
\eeqa
The measure of the first set in (\ref{resonant set}) is bounded by
\beqa{first set}
&& \meas\left(\bigcup_{k\in \integer^{n},\bar k\neq 0}\left\{I:\ |\o(I)\cdot k|\leq\frac{g}{|k|^\t}\right\}\right)\nonumber\\
&& \leq \sup_{\cI}\|\o^{-1}\|^n\meas\left(\bigcup_{k\in \integer^{n},\bar k\neq 0}\left\{x\in \o(\cI):\ |x\cdot k|\leq\frac{g}{|k|^\t}\right\}\right)\nonumber\\
&& \leq  \sup_{\cI}\|\o^{-1}\|^n\meas\left(\bigcup_{k\in \integer^{n},\bar k\neq 0}\bigcup_{i=1}^p\left\{x\in U_i:\ |x\cdot k|\leq\frac{g}{|k|^\t}\right\}\right)\nonumber\\
&& \leq  \sup_{\cI}\|\o^{-1}\|^n\sum_{k\in \integer^{n},\bar k\neq 0}\sum_{1\leq i\leq p}\meas\left(\left\{x\in U_i:\ |x\cdot k|\leq\frac{g}{|k|^\t}\right\}\right)\nonumber\\
&& =  \sup_{\cI}\|\o^{-1}\|^n\sum_{k\in \integer^{n},\bar k\neq 0}\sum_{1\leq i\leq p}\int_{\cB_k^i}d\bar xd\hat x\nonumber\\
\eeqa
where
$$\cB_k^i:=\left\{x'=(\bar x',\hat x')\in U_i=B^{\bar n}_{\bar R}(0)\times B^{\hat n}_{\hat R}(\o_0(\bar I_i)):\ |\bar x'\cdot\bar k+\hat x'\cdot\hat k|\leq\frac{g}{|k|^\t}\right\}$$
Now, as $\bar k\neq 0$, we certainly find $1\leq j\leq\bar n$ with $|\bar k_j|\geq 1$. Perform, then, the change of variables
$$z_m=\bar x_m\quad \textrm{for}\quad 1\leq m\leq \bar n\ ,\ m\neq j\ ,\quad z_j=\bar x\cdot \bar k+\hat x\cdot \hat k\ ,\quad z_{m}=\hat x_{m-\hat n}\quad \textrm{for}\quad \bar n+1\leq m\leq n$$
Then, letting
\beqano
\tilde \cB^i_k&:=&\left\{z'=(z'_1,\cdots,z'_n):\right.\nonumber\\
&&\left. \left((z'_1,\cdots,z'_{j-1},\frac{1}{k_j}\left(z'_j-\sum_{m\neq j}\bar z'_m\bar k_m-\sum_{\bar n+1\leq m\leq n}z'_m\hat k_{m-\bar n}\right),z'_{j+1},\cdots,z'_{\bar n}\right)\right.\nonumber\\
&&\in B^{\bar n}_{\bar R}(0)\ ,\nonumber\\
&&\left.|z'_j|\leq \frac{g}{|k|^\t}\ ,\ (z'_{\bar n+1},\cdots,z'_{n})\in B^{\hat n}_{\hat R}(\o_0(\bar I_i))\right\}\supseteq\nonumber\\
&\supseteq&\left\{z'=(z_1',\cdots,z_n'):\ |z'_m|\leq \bar R\right.\nonumber\\
&& \textrm{for}\ 1\leq m\neq j\leq \bar n,\ |z'_j|\leq \frac{g}{|k|^\t},\nonumber\\
&&\left. |z'_m-\o_0(\bar I_i)|\leq\hat R\ \textrm{for}\ m=\bar n+1,\cdots,n\right\}\nonumber\\
&=:&\cC_k^i\ .
\eeqano
we find
$$\int_{\cB_k^i}d\bar x d\hat x=\frac{1}{|k_j|}\int_{\tilde \cB_k^i}dz\leq \frac{1}{|k_j|}\int_{\cC_k^i}dz\leq \int_{\cC_k^i}dz=\bar R^{\bar n-1}\hat R^{\hat n}\frac{g}{|k|^\t}$$
Hence, inserting this expression into (\ref{first set}), we find
\beqano\meas\left(\bigcup_{k\in \integer^{n},\bar k\neq 0}\left\{I:\ |\o(I)\cdot k|\leq\frac{g}{|k|^\t}\right\}\right)&\leq& \|\o^{-1}\|^n\bar R^{\bar n-1}\hat R^{\hat n}p\, g\,\sum_{k\in \integer^{n},\bar k\neq 0}\frac{1}{|k|^\t}\nonumber\\
&\leq& \sup_{\cI}\|\o^{-1}\|^n\bar R^{\bar n-1}\hat R^{\hat n}c(n,\t)g\nonumber\\
&=& \sup_{\cI}\|\o^{-1}\|^n\bar R^{\bar n-1}\cA^{\hat n}\hat r^{\hat n}c(n,\t)g\nonumber\\
&=&\bar c g\meas(\cI)
\eeqano
We now  estimate  the measure of the second set in (\ref{resonant set}). By Fubini's Theorem, 
we find
\beqa{second integral}
&& \meas\left(\bigcup_{0\neq\hat k\in \integer^{\hat n}}\left\{I:\ |\hat\o(I)\cdot \hat k|\leq\frac{\hat g}{|\hat k|^\t}\right\}\right)=\int_{\bar\cI}d\bar I\int_{\bigcup_{\hat k\neq 0}\hat\cB_k(\bar I)}d\hat I
\eeqa
where
$$\hat\cB_k(\bar I)=\bigcup_{0\neq\hat k\in \integer^{\hat n}}\left\{\hat I\in B^{\hat n}_{\hat r}(0):\ |(\hat\o_0(\bar I)+A(\bar I)\hat I)\cdot \hat k|\leq\frac{\hat g}{|\hat k|^\t}\right\}$$
Perform, in the inner integral, the change of variable
$$\hat x=\hat\o_0(\bar I)+A(\bar I)\hat I$$
and let
\beqano
\tilde\cC_k(\bar I)&:=&\left\{x\in \real^{\hat n}:\ A(\bar I)^{-1}(\hat x-\o_0(\bar I))\in B^{\hat n}_{\hat r}(0)\ ,\quad |\hat x\cdot \hat k|\leq\frac{\hat g}{|\hat k|^\t}\right\}\nonumber\\
&\subseteq&\left\{\hat x\in \real^{\hat n}:\ \hat x\in B^{\hat n}_{\cA\hat r}(\o_0(\bar I))\ ,\quad |\hat x\cdot \hat k|\leq\frac{\hat g}{|\hat k|^\t}\right\}\nonumber\\
&=:&\hat \cC_k(\bar I)
\eeqano
Then, proceeding as done for the first part of the proof (\ie, with a suitable change of variable, for which $z'_j=\hat x\cdot \hat k$ if $\hat k_j\neq 0$), 
we find
\beqano
\int_{\bigcup_{\hat k\neq 0}\hat\cB_k(\bar I)}d\hat I&\leq& \sup_{\bar\cI}\|A^{-1}\|^{\hat n}\int_{\bigcup_{\hat k\neq 0}\tilde\cC_k(\bar I)}d\hat I\nonumber\\
&\leq& \sup_{\bar\cI}\|A^{-1}\|^{\hat n}\int_{\bigcup_{\hat k\neq 0}\hat\cC_k(\bar I)}d\hat x\nonumber\\
&\leq& \sup_{\bar\cI}\|A^{-1}\|^{\hat n}(\cA\hat r)^{\hat n-1}{\hat g}\sum_{\hat k\neq 0}\frac{1}{|k|^\t}
\eeqano
Hence, inserting this value into (\ref{second integral}), we find
\beqano
\meas\left(\bigcup_{0\neq\hat k\in \integer^{\hat n}}\left\{I:\ |\hat\o(I)\cdot \hat k|\leq\frac{\hat g}{|\hat k|^\t}\right\}\right)&\leq& \meas(\bar\cI)\sup_{\bar\cI}\|A^{-1}\|^{\hat n}(\cA\hat r)^{\hat n-1}{\hat g}c(\hat n,\t)\nonumber\\
&=&\frac{\sup_{\bar\cI}\|A^{-1}\|^{\hat n}(\cA)^{\hat n-1}}{\meas(A_1)}c(\hat n,\t)\frac{\hat g}{\hat r}\meas(\cI)\nonumber\\
&=&\hat c \frac{\hat g}{\hat r}\meas(\cI)
\eeqano
because, for small $\hat r$,
$$\hat r^{\hat n-1}\leq \tilde c\,\frac{\textrm{\rm meas}\Big(\hat{\cal I}\Big)}{\hat r}$$
since $\a>1$.
\newpage
\section{Kolmogorov's Set in the Plane Planetary Problem}\label{Plane Problem}
\setcounter{equation}{0}
Let us consider the motion of a system of $1+N$ masses ${m}_0$, $\cdots$, ${m}_N$ moving in $\real^{\rm d}$ (but, soon, we will take ${\rm d}=2$) under the only influence of gravity. As customary, we restrict to the ``planetary case'': one mass, ${m}_0$ (the ``Sun'') is much greater than ${m}_1$, $\cdots$, ${m}_N$ (``the planets''), namely, we take
\beq{masses}{m}_0=\bar m_0\ ,\quad {m}_1=\m \bar m_1\ ,\quad \cdots\ ,\quad {m}_N=\m \bar m_N\quad (\m\ll1) \ .\eeq
The motion equations are
\beq{Newt eq}\ddot v_i=-\sum_{0\leq j\leq N,\ j\neq i}{m}_j\frac{v_i-v_j}{|v_i-v_j|^3}\ ,\quad \  \quad 0\leq i\leq N\ .\eeq
In the Hamiltonian formalism, equations (\ref{Newt eq}) are equivalent to the study of the Hamiltonian
\beq{not lin red}\hat {\cal H}_{\rm plt}(u,v)=\sum_{0\leq i\leq N}\frac{|u_i|^2}{2{m}_i}-\sum_{0\leq i<j\leq N}\frac{{m}_i{m}_j}{|v_i-v_j|}\qquad (u_i={m}_i \dot v_i)\eeq
on the phase space
\beqano
\hat\cC_{\rm cl,d}:=\Big\{(u,v)=\Big((u_0,\cdots,u_{N}),(v_0,\cdots,v_{N})\Big)\in (\real^{\rm d})^{1+N}\times (\real^{\rm d})^{1+N}:\nonumber\\
v_i\neq v_j\quad \textrm{for}\quad i\neq j\Big\}\ .
\eeqano 
The number of degrees of freedom of (\ref{not lin red}) may be reduced (from ${\rm d}(1+N)$ to ${\rm d}N$) as follows. 
On the (invariant) symplectic manifold  with dimension ${\rm d}N$
$$\cM_{\rm lin}=\left\{u=(u_{0},u_1,\cdots,u_N), v=(v_{0},v_1,\cdots v_N)\in \hat\cC_{\rm cl,d}:\ \sum_{0\leq i\leq N}u_i=\sum_{0\leq i\leq N}{m}_iv_i=0\right\}$$
we introduce the relative coordinates
\beqa{rel coord}
\arr{\tilde x_i=v_i-v_0\\
\tilde y_i=u_i
}\qquad \textrm{for}\qquad 1\leq i\leq N\ .
\eeqa
The motion of ${m}_0$  is then recovered by 
\beqa{coord of Sun}
\arr{
u_{0}=-\sum_{1\leq i\leq N}\tilde y_i\\
v_0=-\frac{\sum_{1\leq i\leq N}{m}_i\tilde x_i}{\sum_{0\leq i\leq N}{m_i}}
}
\eeqa
as it results by requiring
$$\arr{0=\sum_{0\leq i\leq N}u_i=u_{0}+\sum_{1\leq i\leq N}\tilde y_i\\
0=\sum_{0\leq i\leq N}{m}_i v_i={m}_0 v_0+\sum_{1\leq i\leq N}{m}_i(\tilde x_i+v_0)\ .
}
$$
The parametrization of ${\cM}_{\rm lin}$ (\ref{rel coord})$\div$(\ref{coord of Sun}), which expresses a point $(u,v)\in \cM_{\rm lin}$ in terms of the coordinates $(\tilde y,\tilde x)$, is a homogeneous symplectic embedding, \ie, it preserves the Liouville $1$--form:
\beqano
\sum_{0\leq i\leq N} u_i\,dv_i&=&u_0\,dv_0+\sum_{1\leq i\leq N}u_i dv_i\nonumber\\
&=&u_0\,dv_0+\sum_{1\leq i\leq N}\tilde y_i d(\tilde x_i+v_0)\nonumber\\
&=&\left(u_0+\sum_{1\leq i\leq N}\tilde y_i\right)dv_0+\sum_{1\leq i\leq N}\tilde y_i d\tilde x_i\nonumber\\
&=&\sum_{1\leq i\leq N}\tilde y_i d\tilde x_i
\eeqano
(because $u_{0}$ just coincides with $-\sum_{1\leq i\leq N}\tilde y_i$).\\
Then, Hamiltonian (\ref{not lin red}), in terms of the relative coordinates $(\tilde y,\tilde x)$, with the masses (\ref{masses}), becomes
$$\tilde {\cal H}_{\rm plt}(\m;\tilde y,\tilde x):=\hat {\cal H}_{\rm plt}\circ\phi_{\rm lin}=\sum_{1\leq i\leq N}\left(\frac{|\tilde y_i|^2}{2\tilde m_i\m}-\frac{\hat m_i\tilde m_i\m}{|\tilde x_i|}\right)+\m\sum_{1\leq i<j\leq N}\left(\frac{\tilde y_i\cdot \tilde y_j}{\bar m_0\m}-\frac{\bar m_i\bar m_j\m}{|\tilde x_i-\tilde x_j|}\right)\ ,$$
where \beq{reduced masses}\hat m_i=\bar m_0+\m \bar m_i\ ,\qquad \tilde m_i=\frac{\bar m_0\bar m_i}{\bar m_0+\m \bar m_i}\eeq are the {\sl reduced masses}. 
A rescaling of the variables
\beqa{rescale}
\arr{
\tilde y=\m y\\
\tilde x=x
}
\eeqa
joined with the rescaling of the Hamiltonian $${\cal H}_{\rm plt}(\m;y,x):=\m^{-1}\tilde {\cal H}_{\rm plt}(\m;\m y,x)$$ 
(which does not change Hamiltonian form of the equations of the motion) brings to the Hamiltonian
\beqa{planar 1+N BP}
{\cal H}_{\rm plt}(\m;y,x)=\sum_{1\leq i\leq N}\left(\frac{| y_i|^2}{2\tilde m_i}-\frac{\hat m_i\tilde m_i}{| x_i|}\right)+\m\sum_{1\leq i<j\leq N}\left(\frac{ y_i\cdot  y_j}{\bar m_0}-\frac{\bar m_i\bar m_j}{| x_i- x_j|}\right)
\eeqa
with $(y,x)$ varying in the ``collisionless'' domain 
\beqa{no coll}\cC_{\rm cl,d}:=\Big\{(y,x)=\Big((y_1,\cdots,y_N),(x_1,\cdots,x_N)\Big)\in \real^{{\rm d}N}\times \real^{{\rm d}N}:\nonumber\\
 x_i\neq  x_j\neq 0\quad\forall\ 1\leq i<j\leq N\Big\}\eeqa
 When $\m=0$, the Hamiltonian ${\cal H}_{\rm plt}$ (\ref{planar 1+N BP}) splits into the sum of $N$ Two--Body (integrable) Hamiltonians describing each the interction of a fictictious mass $\tilde m_i$ with a fixed star with mass $\hat m_i$.

 \vskip.1in
 \noi
 The aim of this section is the proof of Theorem \ref{plane problem} below.
\begin{theorem}\label{plane problem}
Consider the  evolution in time of the coordinates  of $1+N$ ``planetary'' masses (\ref{masses})
moving on the plane undergoing Newtonian attraction. Let $a_i$, $e_i$ denote the semimajor axis and the eccentricity of the Keplerian ellipse arising from the two -- body interaction of a fictictous mass   $\tilde m_i$ with a fixed star $\hat m_i$ in correspondence of the initial datum $(\bar y_i,\bar x_i)$ for the coordinates $(y_i,x_i)$ described in (\ref{rel coord})$\div$(\ref{rescale}), where $\tilde m_i$, $\hat m_i$ are as in (\ref{reduced masses}). Then, there exist $b$, $c$, $C$, $\d_*>0$ such that, for any $0<\d<\d_*$, a parameter $\e_*=\e_*(\d)$ may be found such that, for any $$0<\e<\e_*\quad \textrm{and}\quad 0<\m<(\log{\e^{-1}})^{-2b}$$
in the  set of $(\bar y,\bar x)=\Big((\bar y_1,\cdots,\bar y_N), (\bar x_1,\cdots,\bar x_N)\Big)\in (\real^{2})^N\times (\real^{2})^N$ such that
$$a_i=\hat a_i\d^{N-i}\ ,
\quad \textrm{where}\quad \underline a\leq {\hat a}_i\leq \ovl a$$
 there exists a positive Lebesgue measure set $\cK$ (``Kolmogorov set''), satisfying 
 $$c\,\e^{2N}>\meas\,\cK>c\,\Big(1-C(\sqrt{\e}+\sqrt{\m}(\log{\e^{-1}})^{b})\Big)\e^{2N}\ ,$$
  formed by the union of invariant tori of dimension $2N$ on which the  $\cH_{\rm plt}$--flow is linear in time, with Diophantine frequency.
 Furthermore, the eccentricities on the invariant tori are bounded by $c(\log\e^{-1})^{-1}$.
\end{theorem}
\subsection{The Plane Delaunay--Poincar\'e Map}
A good set of action--angle variables for the plane problem (\ref{planar 1+N BP})$\div$(\ref{no coll}) is the set of {\sl Delaunay variables} $(L$,{\bf G},$\ell$,{\bf g}),  $L=(L_1,\cdots,L_N)$, {\bf G}=$(G_1,\cdots,G_N)$, $\cdots$, 
with
$$0<G_i<L_i\ ,\qquad \ell_i,\ g_i\in \torus$$
defined as 
\beqa{Delaunay variables}
\arr{
L_i=\tilde m_i\sqrt{\hat m_ia_i}\\
G_i=|x_i\times y_i|=\sqrt{1-e_i^2}L_i\\
}\qquad 
\arr{
\ell_i=\frac{2{\cal A}_i}{a_i^2\sqrt{1-e_i^2}}\\
g_i=\textrm{argument of}\quad P_i\\
}
\eeqa
where, on the ($\hat m_i,\tilde m_i;y_i,x_i$)--``osculating'' ellipse
\footnote{
{\it I.e.}, the ellipse arising from the  initial datum $(\dot x(0),x(0))=(y_i/\tilde m_i,x_i)$ for  planet in  Newtonian interaction with  a star with mass $\hat m_i$.}  the quantities $a_i$, $e_i$, $P_i$, ${\cal A}_i$ denote, respectively, the semimajor axis, eccentricity, perihelion, the area of the elliptic sector from $P_i$ to $x_i$.
In terms of the Delaunay variables, the linear momenta $y_i$ and the positions $x_i$ are recovered by
\beqa{momentum position}
\arr{x_i=x_i^{\rm D}:=a_i{\cal R}_{\rm z}(g_i)\left(
\begin{array}{lrr}
\cos{u_i}-e_i\\
\sqrt{1-e_i^2}\sin{u_i}
\end{array}
\right)\\ y_i=y_i^{\rm D}:=\tilde m_in_i\partial_{\ell_i} x_i} \ n_i^2a_i^3:=\hat m_i\ ,\ (n_i:\ \textrm{\rm ``mean motion''})\eeqa
where $a_i$, $e_i$ are thought as functions of $L_i$, $G_i$,\ie, \beqa{ai ei}a_i=\frac{1}{\hat m_i}\left(\frac{L_i}{\tilde m_i}\right)^2\ ,\quad e_i=\sqrt{1-\left(\frac{G_i}{L_i}\right)^2}\eeqa $u_i$ solves the {\sl Kepler's Equation}
$$u_i-e_i\sin{u_i}=\ell_i$$
and ${\cal R}_{\rm z}(g)$ denotes a rotation of $g$ in the plane:
$${\cal R}_{\rm z}(g)=\left(
\begin{array}{lrr}
\cos{g}&-\sin{g}\\
\sin{g}&\cos{g}
\end{array}
\right)$$
The plane Delaunay variables (\ref{Delaunay variables}) are well defined whenever
$$e_i\neq 0\quad \textrm{for}\quad 1\leq i\leq N\ .$$
A  ``regularization'', due to Poincar\'e, allows the system reaching also zero eccentricities. It  is achieved with the symplectic change of variables
\beqa{reg Poi}\phi_{\rm P}^{-1}:\quad \arr{\L_i=L_i\\
\l_i=l_i+g_i
}\qquad \arr{\eta_i=\sqrt{2(L_i-G_i)}\cos{g_i}\\
\xi_i=-\sqrt{2(L_i-G_i)}\sin{g_i}}\eeqa
The regularized variables (\ref{reg Poi}) are usually called {\sl Poincar\'e variables} and, in terms of them, 
 (\ref{momentum position}) become
\footnote{As usual, $\L=(\L_1,\cdots,\L_N)$, $\cdots$.}
\beqa{DP2}\arr{x_i=\hat x_i:=x_{\rm DP}(\hat m_i,\tilde m_i;\L_i,\l_i,\eta_i,\xi_i)=(x^{\rm DP}_1, x^{\rm DP}_2)\\
\\
y_i=\hat y_i:=y_{\rm DP}(\hat m_i,\tilde m_i;\L_i,\l_i,\eta_i,\xi_i)=\frac{\hat m^2\tilde m^4}{\L^3}\,\partial_{\l}x_{\rm DP}}\eeqa
where 
\beqa{DP4}
\left\{
\begin{array}{lrr}
x^{\rm DP}_1=\frac{1}{ {\hat m}}\left(\frac{\L}{{\tilde m}}\right)^2\left[\cos{(\hat\zeta+\l)}-\frac{{\xi}}{2{\L}}\left({\eta}\sin{(\hat\zeta+\l)}+{\xi}\cos{(\hat\zeta+\l)}\right)-\frac{{\eta}}{\sqrt{\L}}\sqrt{1-\frac{{\eta}^2+{\xi}^2}{4{\L}}}\right]\\
\\
x^{\rm DP}_2=\frac{1}{ {\hat m}}\left(\frac{\L}{{\tilde m}}\right)^2\left[\sin{(\hat\zeta+\l)}-\frac{{\eta}}{2{\L}}\left({\eta}\sin{(\hat\zeta+\l)}+{\xi}\cos{(\hat\zeta+\l)}\right)+\frac{{\xi}}{\sqrt{\L}}\sqrt{1-\frac{{\eta}^2+{\xi}^2}{4{\L}}}\right]\\
\end{array}
\right.\nonumber\\
\eeqa
and $\hat \zeta $ solves the {\sl regularized Kepler Equation}
\beqano
\zeta=\frac{1}{\sqrt{\L}}\sqrt{1-\frac{{\eta}^2+{\xi}^2}{4\L}}\,\left({{\eta}}\sin{(\zeta+\l)}+{{\xi}}\cos{(\zeta+\l)}\right)\ .
\eeqano
The semimajor axes $a_i$ and the eccentricities $e_i$ (\ref{ai ei}) become
$$a_i=\frac{1}{\hat m_i}\left(\frac{\L_i}{\tilde m_i}\right)^2\ ,\quad e_i=\sqrt{1-\left(1-\frac{\eta_i^2+\xi_i^2}{2\L_i}\right)^2}$$
(zero eccentricities correspond to $(\eta_i,\xi_i)=0$). In order to avoid collisions we let
\footnote{$\real_+:=(0,+\infty)$.}
\beqa{DP5}\arr{\L\in{\cal A}_\e\\
 \l\in \torus^N\\ (\hat\eta,\hat\xi):=\left(\left(\frac{\eta_1}{\sqrt{\L_1}},\cdots,\frac{\eta_N}{\sqrt{\L_N}}\right),\left(\frac{\xi_1}{\sqrt{\L_1}},\cdots,\frac{\xi_N}{\sqrt{\L_N}}\right)\right)\in{\cal E}_\e}
\eeqa
where, for $0<\e<1$,
\beqa{DP6}
\arr{{\cal A}_\e:=\Big\{\L\in \real_+^N:\ a_i(1+\e)<a_{i+1}(1-\e)\quad 1\leq i\leq N\Big\}\\
{\cal E}_\e:=\Big\{(\hat\eta,\hat\xi)\in \real^N\times \real^N:\ \frac{\hat\eta_i^2+\hat\xi_i^2}{2}\leq 1-\sqrt{1-\e^2}\quad 1\leq i\leq N\Big\}\\}
\eeqa
\begin{proposition}[Delaunay--Poincar\'e] For any $0<\e<1$, in the domain (\ref{DP5}) $\div$ (\ref{DP6}), equations (\ref{DP2}) $\div$ (\ref{DP4}) well define a real--analytic symplectomorphism
\beqa{DP1}\phi_{\rm DP}:\quad \Big(\L,\l,\eta,\xi\Big)\to \Big(\hat y,\hat x\Big)=\Big((\hat y_1,\cdots \hat y_N),(\hat x_1,\cdots \hat x_N)\Big)\eeqa
namely, a $1:1$ onto, real--analytic and   symplectic map, with respect to the standard $2$--form $$\sum_{1\leq i\leq N}\Big(d\L_i\wedge d\l_i+d\eta_i\wedge d\xi_i\Big)\ ,$$
usually called {\sl plane Delaunay--Poincar\'e map}, which carries  Hamiltonian of the Plane $(1+N)$--Body Problem, \ie, the Hamiltonian (\ref{planar 1+N BP}) defined on the domain ${\cal C}_{\rm cl,2}$ (\ref{no coll}) to
\beq{calHNBP}{\cal H}_{\rm plt}:=-\sum_{1\leq i\leq N}\frac{\tilde m_i^3\hat m_i^2}{2\L_i^2}+\m\sum_{1\leq i<j\leq N}\left(\frac{ \hat y_i\cdot\hat y_j}{\bar m_0}-\frac{\bar m_i\bar m_j}{|\hat x_i- \hat x_j|}\right)\eeq
\end{proposition}
For a self--contained proof of this Proposition--not easy to be found in literature, see, \cite{Che88}, \cite{Lask88}, and also \cite{BCV03}.

\vskip.2in
\noi
{\bf Sketch of the proof of Theorem \ref{plane problem}.}\ The proof of Theorem \ref{plane problem} consists in applying Theorem \ref{more general degenerate KAM} (in the simplified version of Remark \ref{degenerate KAM}) to the properly degenerate Hamiltonian ${\cal H}_{\rm plt}$  of the plane $(1+N)$--Body problem, espressed in Delaunay--Poincar\'e variables (\ref{calHNBP}). We have thus to check all the assumptions thereby involved. We do this in the following  steps.

\vskip.1in
\noi
\begin{itemize}
\item[(i)] Let
$$f:=\sum_{1\leq i<j\leq N}\left(\frac{ \hat y_i\cdot\hat y_j}{\bar m_0}-\frac{\bar m_i\bar m_j}{|\hat x_i- \hat x_j|}\right)$$
denote the ``perturbation'' of ${\cal H}_{\rm plt}$ and
$$\bar f:=\frac{1}{(2\p)^N}\int_{\torus^N}f\,d\l$$
its  mean (also called ``secular'' perturbation) with respect to the ``fast'' angles $\l$.\\
Due to the {\sl D'Alembert relations} (Lemma \ref{simmetries remark}), $\bar f$ is even in $(\eta,\xi)$, hence, it has an equilibrium point at the origin of the ``secular'' coordinates, \ie, for $z:=(\eta,\xi)=0$ (Laplace); its quadratic and quartic parts have the form
\beq{quadratic}\frac{1}{2}\eta\cdot {\cal F}(\L)\eta+\frac{1}{2}\xi\cdot {\cal F}(\L)\xi\eeq
and
$$\sum_{1\leq i,j,k,l\leq N}q_{i,j,k,l}(\L)\,\eta_i\eta_j\eta_k\eta_l+r_{i,j,k,l}(\L)\,\eta_i\eta_j\xi_k\xi_l+q_{i,j,k,l}(\L)\,\xi_i\xi_j\xi_k\xi_l$$
respectively.
In particular, since ${\cal F}(\L)$ is a symmetric, in view of (\ref{quadratic}), $z=0$ in an {\sl elliptic} equilibrium point, and the eigenvalues $\O=(\O_1,\cdots,\O_N)$  of ${\cal F}(\L)$ have the meaning of the Birkhoff invariants with order $1$. Both the entries of ${\cal F}(\L)$ and the tensors ${\cal Q}:=(q_{ijkl})$, ${\cal R}:=(r_{ijkl})$ can be expressed in terms of the Laplace coefficients (Lemmas \ref{exp to 2}, \ref{order four}). 
\item[(ii)] The diagonalization (by  means of a unitary matrix $U(\L)$) of ${\cal F}(\L)$ is required, in order to check the  $4$--non resonance of the Birkhoff invariants  with order $1$ $\O=(\O_1,\cdots,\O_N)$ for $\bar f$. We prove a technical Lemma (Lemma \ref{diagonalization N bodies}) which implies  that the lowest asymptotics for $\O(\L)$ is just the one of the diagonal elements for ${\cal F}$, and so it is given by  
\beqano
\O_i(\L)&=&\d^{(9-3N)/2}\nonumber\\
&\times&\arr{
-\frac{3}{4}\frac{\hat a_1^2}{\hat a_2^3}\frac{m_1\bar m_2}{\tilde m_1\sqrt{\hat m_1\hat a_1}}+O(\d^2)\qquad \textrm{for}\quad i=1\\
\\
-\frac{3}{4}\frac{\hat a_{i-1}^2}{\hat a_i^3}\frac{\bar m_i\bar m_{i-1}}{\tilde m_i\sqrt{\hat m_i\hat a_i}}\d^{(3i-5)/2}+O(\d^{(3i-3)/2)})\quad \textrm{for}\quad 2\leq i\leq N-1\\
\\
-\frac{3}{4}\frac{\hat a_{N-1}^2}{\hat a_N^3}\frac{\bar m_N\bar m_{N-1}}{\tilde m_N\sqrt{\hat m_N\hat a_N}}\d^{(3N-5)/2}+O(\d^{(3N-2)/2})\quad \textrm{for}\quad i=N
}
\eeqano
which immediately implies non resonance up to any finite order, for small $\d$ (Corollary \ref{non res cor}).
We also compute the lowest $\d$--asymptotics for the entries of $U(\L)$ (Lemma \ref{order 2 asymptotics}). 
\item[(iii)] We diagonalize the quadratic part of $\bar f$ with the symplectic transformation
$$\phi_{\rm diag}:\quad \eta=U(\tilde\L)\tilde\eta\quad \xi=U(\tilde\L)\tilde\xi\quad \L=\tilde\L\quad \l=\tilde\l+\varphi(\tilde\L,\tilde\eta,\tilde\xi)$$
(where $\varphi$ a suitable shift of $\l$ which does not change the mean), hence, we put $\bar f$ into the form
\beqano
\tilde f&:=&\bar f\circ\phi_{\rm diag}=\bar f_0(\tilde\L)+\sum_{1\leq i\leq N}\O_i(\tilde\L)\frac{\tilde\eta_i^2+\tilde\xi_i^2}{2}\nonumber\\
&+&\sum_{1\leq i,j,k,l\leq N}\tilde q_{i,j,k,l}(\tilde\L)\,\tilde\eta_i\tilde\eta_j\tilde\eta_k\tilde\eta_l+\tilde r_{i,j,k,l}(\tilde\L)\,\tilde\eta_i\tilde\eta_j\tilde\xi_k\tilde\xi_l+\tilde q_{i,j,k,l}(\tilde\L)\,\tilde\xi_i\tilde\xi_j\tilde\xi_k\tilde\xi_l\nonumber\\
&+&{\rm o}_4
\eeqano
\item[(iv)] We compute the $\d$--asymptotics for $\tilde q_{ijkl}$, $\tilde r_{ijkl}$ (which involves those of $ r(\L)$, $s(\L)$, $U(\L)$) and hence the $\d$--asymptotics for the entries of the Birkhoff invariants with order $2$, which are the entries of the symmetric matrix $A(\bar\L)$ of the Birkhoff Normal form of $\tilde f$ 
\beqano
\bar f_0(\bar\L)+\O(\bar\L)\cdot J+\frac{1}{2}J\cdot A(\bar\L)J+{\rm o}(|J|^{3})\quad J_i:=\frac{\bar\eta_i^2+\bar\xi_i^2}{2}
\eeqano
obtained by projections of the entries  $\tilde q_{ijkl}$, $\tilde r_{ijkl}$ (Lemma \ref{order 2 invariants}). We check (Lemma \ref{Aij asymptotics}) that $A(\bar\L)$ has the form
$${A}(\L)\approx\d^p\left(
\begin{array}{lrrrrr}
\a_{11}&\a_{12}&O(\d^{p_{13}})&\cdots &O(\d^{p_{1k}})&\cdots\\
\a_{21}&\a_{22}&O(\d^{p_{23}})&\cdots&O(\d^{p_{2k}})&\cdots\\
&&\a_{33}\d^{p_{33}}&\cdots&O(\d^{p_{3k}})&\cdots\\
&&&\ddots&\vdots&\\
&&&&\a_{kk}\d^{p_{kk}}&\cdots\\
&&&&&\ddots\\
\end{array}
\right)$$
\end{itemize}
where $\a_{11}\a_{22}-\a_{12}\a_{21}\neq 0$ and $p_{k+1,k+1}>p_{kk}$,  $\a_{kk}\neq 0$ and that that this implies (Lemma \ref{non degeneracy lem})
$$\textrm{\rm det}\,A=(\a_{11}\a_{22}-\a_{12}\a_{21})\d^q+{\rm o}(\d^q)\neq 0\ ,$$
concluding the proof.

\subsection{Non Resonance and non Degeneracy for the Plane Planetary Problem}

\subsubsection{Expansion of the Hamiltonian}
The perturbation $f$ of the Plane Planetary Problem is composed of two terms
$$f_{\rm p}:=-\sum_{1\leq i<j\leq N}\frac{\bar m_i \bar m_j}{|\hat x_i-\hat x_j|}\ ,\quad f_{\rm s}:=\frac{1}{\bar m_0}\sum_{1\leq i<j\leq N}\hat y_i\cdot \hat y_j$$
usually refered as {\sl principal} and {\sl secondary part}, respectively. We are interested to the {\sl secular perturbation}, \ie, the mean, over $\l\in \torus$, of $f$, to which only $f_{\rm d}$ contributes:
\begin{lemma}\label{no indirect}
The secondary part of the perturbation
has zero mean.
\end{lemma}
{\bf Proof.}\ In fact, as $\dst \hat y_i$ is defined as $$\dst \hat y_i:=\frac{\hat m_i^2\tilde m_i^4}{\L_i^3}\partial_{\l_i}\,\hat x_i$$ and does not depend on the variables $\l_j$ with $j\neq i$,
\beqano
\frac{1}{(2\p)^N}\int_{\torus^N}\hat y_i\cdot \hat y_j\,d\l&=&\frac{1}{(2\p)^2}\int_{\torus^2}\hat y_i\cdot\hat y_j\,d\l_id\l_j\nonumber\\
&=&\frac{1}{(2\p)^2}\frac{\hat m_i^2\tilde m_i^4}{\L_i^3}\frac{\hat m_j^2\tilde m_j^4}{\L_j^3}\nonumber\\
&\times&\int_{\torus^2}\partial_{\l_i}\,\hat x_i\cdot\partial_{\l_j}\,\hat x_i\,d\l_id\l_j\nonumber\\
&=&\frac{1}{(2\p)^2}\frac{\hat m_i^2\tilde m_i^4}{\L_i^3}\frac{\hat m_j^2\tilde m_j^4}{\L_j^3}\int_{\torus}\partial_{\l_i}\,\hat x_id\l_i\nonumber\\
&\cdot&\int_{\torus}\partial_{\l_j}\,\hat x_i\,d\l_j\nonumber\\
&=&0\ .
\eeqano
\vskip.1in
\noindent
Notice that that  $\hat x_i$--component of the Plane Delaunay Poincar\'e   depends on $\L_i$, $\eta_i$, $\xi_i$ as a function of  $$a_i=\frac{1}{\hat m_i}\left(\frac{\L_i}{\tilde m_i}\right)^2\ ,\quad (\hat\eta_i, \hat\xi_i)=\left(\frac{\eta_i}{\sqrt{\L_i}},\ \frac{\xi_i}{\sqrt{\L_i}}\right)\ ,$$ only. 
\begin{lemma}[D'Alembert relations]\label{simmetries remark}
Let
\beq{Tay exp}\sum_{(j_1,j_2,i_1,i_2)}a_{j_1,j_2,i_1,i_2}(a_i,a_j)\hat{\eta_i}^{j_1}\,\hat{\eta_j}^{j_2}\,\hat{\xi_i}^{i_1}\,\hat{\xi_j}^{i_2}
\eeq
the Taylor expansion of
\beqa{gij}
\frac{1}{(2\p)^2}\int_{\torus^2}\frac{d\l_id\l_j}{|\hat x_i-\hat x_j|}\nonumber\\
\eeqa
Then, 
\begin{eqnarray}\label{simmetries of coeff}
i)&&\quad a_{j_2j_1i_2i_1}({a_i},{a_j})=a_{j_1j_2i_1i_2}({a_j},{a_i})\ ;\nonumber\\
ii)&&\quad a_{j_1j_2i_1i_2}({a_i},{a_j})=0\quad \textrm{if}\quad j_1+j_2\quad \textrm{is odd}\ ;\nonumber\\
iii)&&\quad a_{j_1j_2i_1i_2}({a_i},{a_j})=0\quad \textrm{if}\quad i_1+i_2\quad \textrm{is odd}\ ;\nonumber\\
iv)&&\quad a_{i_1i_2j_1j_2}({a_i},{a_j})=a_{j_1j_2i_1i_2}({a_i},{a_j})\ .
\end{eqnarray}
\end{lemma}
{\bf Proof.}\ Item $(i)$ is trivial.  Items $(ii)\div (iv)$ are related to the following symmeries of the plane Delaunay--Poincar\'e map
\begin{eqnarray*}
\hat x_i(\hat m,\tilde m;\L,\p-\l,-\eta,\xi)&=&{\cal R}_{x=0}\hat x_i(\hat m,\tilde m;\L,\l,\eta,\xi)\ ,\nonumber\\
\hat x_i(\hat m,\tilde m;\L,-\l,\eta,-\xi)&=&{\cal R}_{y=0}\hat x_i(\hat m,\tilde m;\L,\l,\eta,\xi)\ ,\nonumber\\
\hat x_i(\hat m,\tilde m;\L,\p/2-\l,\xi,\eta)&=&{\cal R}_{x=y}\hat x_i(\hat m,\tilde m;\L,\l,\eta,\xi)\nonumber\\
\end{eqnarray*}
where ${\cal R}_{y=0}$, ${\cal R}_{x=0}$, ${\cal R}_{x=y}$, denote the reflections in the plane with respect to the axes $x$, $y$, $x=y$, axes.
\vskip.1in

\vskip.1in
\noi
\begin{remark}\rm
Then, the secular perturbation $\bar f$
 contains only polynomials $f_{2j}(\L,\eta,\xi)$ with even degree $2j$:
\begin{eqnarray}\label{developing f}
\bar f(\L,\eta,\xi)=f_{0}(\L)+f_{2}(\L,\eta,\xi)+f_{4}(\L,\eta,\xi)+\cdots\ ;
\end{eqnarray}
where each $f_{2i}$ is an even function of $\eta$, $\xi$ separately. In particular,
\end{remark}
\begin{corollary}[Laplace]
The point $(\eta,\xi)=0$ is an equilibrium point for $\bar f$, for all $\L$.
\end{corollary}
\begin{remark}\rm The computation of the $0$--term $f_0(\L)$ in (\ref{developing f}) is trivial. When $(\eta,\xi)=0$, $\hat x_i$ reduces to
$$\hat x_i|_{(\eta,\xi)=0}=a_i(\cos\l_i,\sin\l_i)\ .$$
Hence,
$$|\hat x_i|_{(\eta,\xi)=0}-\hat x_j|_{(\eta,\xi)=0}|=\sqrt{a_i^2+a_j^2-2a_ia_j\cos{(\l_i-\l_j)}}$$
and, finally
\beqano
f_0&=&-\sum_{1\leq i<j\leq N}\bar m_i\bar m_j\frac{1}{(2\p)^2}\int_{\torus^2}\frac{d\l_id\l_j}{|\hat x_i|_{(\eta,\xi)=0}-\hat x_j|_{(\eta,\xi)=0}|}\nonumber\\
&=&-\sum_{1\leq i<j\leq N}{\bar m_i\bar m_j}\frac{1}{2\p}\int_{\torus}\frac{dt}{\sqrt{a_i^2+a_j^2-2a_ia_j\cos{(\l_i-\l_j)}}}\nonumber\\
&=&-\sum_{1\leq i<j\leq N}\frac{\bar m_i\bar m_j}{a_j}b_{1/2,0}(a_i/a_j)
\eeqano
where $b_{s,k}(\a)$ is the $(s, k)$--Laplace coefficient, defined as the $k^{th}$ Fourier coefficient of the function $t\to\Big[1+\a-2\a\cos{(\l_i-\l_j)}\Big]^{-s}$:
$$b_{s,k}(\a)=\frac{1}{2\p}\int_{\torus}\frac{\cos{kt}}{\Big[1+\a-2\a\cos{(\l_i-\l_j)}\Big]^s}\ .$$
Regularity properties and expansions (in $\a$) of the Laplace Coefficients are briefly discussed in Appendix \ref{Laplace coefficients}. 
\end{remark}

\vskip.1in
\noi
Next two lemmas are devoted to the computation of $f_2$, $f_4$ of (\ref{developing f}).
\begin{lemma}\label{exp to 2}
The polynomial with order $2$ in the expansion  of $\bar f$ is
$$f_2(\L,\eta,\xi)=\frac{1}{2}\eta\cdot\cF(\L)\eta+\frac{1}{2}\xi\cdot \cF(\L)\xi$$
where
$$f_{ij}(\L)=
\arr{-2\,\frac{\bar m_i}{\L_i}\,\left[\sum_{k\neq i}\bar m_k\,a_{2000}(a_i,a_k)\right]\quad \textrm{for}\quad i=j\\
-\bar m_i\bar m_j\frac{a_{1100}(a_i,a_j)}{\sqrt{\L_i\L_j}}\quad \textrm{for}\quad i<j\\
-\bar m_i\bar m_j\frac{a_{1100}(a_j,a_i)}{\sqrt{\L_i\L_j}}\quad \textrm{for}\quad i>j
}
$$
with
\beqa{a2000a1100}
a_{2000}({a},{b})&=&\frac{{a}}{8{b}^2}\,\left[-7{a}/{b}\,b_{5/2,0}({a}/{b})\right.\nonumber\\
&&+\left.4(1+{a}^2/{b}^2)\,b_{5/2,1}({a}/{b})-{a}/{b}\,b_{5/2,2}({a}/{b})\right]\nonumber\\
\nonumber\\
a_{1100}({a},{b})&=&\frac{{a}}{8{b}^2}\,\left[-17\,{a}/{b}\,b_{5/2,1}({a}/{b})\right.\nonumber\\
&&+\left.8(1+{a}^2/{b}^2)\,b_{5/2,2}({a}/{b})+{a}/{b}\,\,b_{5/2,3}({a}/{b})\right]\nonumber\\
\eeqa
\end{lemma}
{\bf Proof.}\ Using the symmetries (\ref{simmetries of coeff}) outlined in Corollary \ref{simmetries remark}, the non vanishing terms with order $2$ appearing in the expansion (\ref{Tay exp}) of $g_{ij}$, are only {\sl six}, and they are individuated by only {\sl two} independent coefficients, say $a_{2000}$ and $a_{1100}$:
\begin{eqnarray}\label{single order 2}
&&{a_{2000}(a_i,a_j)}\,{\hat{\eta}_i^2}+{a_{1100}(a_i,a_j)}\,{\hat\eta_i\hat{\eta}_j}+{a_{2000}(a_i,a_j)}\,\hat{\eta}_j^2\nonumber\\
&&+{a_{2000}(a_i,a_j)}\,{\hat{\xi}_i^2}+{a_{1100}(a_i,a_j)}\,{\hat\xi_i\hat{\xi}_j}+{a_{2000}(a_i,a_j)}\,\hat{\xi}_j^2\nonumber\\
\end{eqnarray}
Thus, multiplying by $-\bar m_i\bar m_j$ and summing over all $1\leq i<j\leq N$ and then symmetrizing the sum, we write $f_2$ as
\beqano
f_2=\frac{1}{2}\eta\cdot\cF(\L)\eta+\frac{1}{2}\xi\cdot \cF(\L)\xi
\eeqano
where the matrix $\cF(\L)=(f_{i,j}(\L)$ has elements
$$f_{ij}(\L)=
\arr{-2\,\frac{\bar m_i}{\L_i}\,\left[\sum_{1\leq k< i}\bar m_k\,a_{2000}(a_k,a_i)+\sum_{i<k\leq N}\bar m_k\,a_{2000}(a_i,a_k)\right]\quad \textrm{for}\quad i=j\\
-\bar m_i\bar m_j\frac{a_{1100}(a_i,a_j)}{\sqrt{\L_i\L_j}}\quad \textrm{for}\quad i<j\\
-\bar m_i\bar m_j\frac{a_{1100}(a_j,a_i)}{\sqrt{\L_i\L_j}}\quad \textrm{for}\quad i>j
}
$$
The coefficients $a_{2000}(a,b)$, $a_{1100}(a,b)$ coincide with the expressions $(ab/8)\cI(a,b)$,\\ $(ab/8)\cJ(a,b)$ computed in \cite{BCV04}, which, written in terms of the Laplace Coefficients are just (\ref{a2000a1100}). The result then follows taking into account the symmetry of the coefficient $a_{2000}$ ($a_{2000}(a,b)=a_{2000}(b,a)$).
\begin{lemma}\label{order four}
The polynomial with order $4$ in the expansion of $\bar f$ is
$$f_4=\sum_{1\leq i,j,k,l\leq N}q_{ijkl}(\L)\Big(\eta_i\eta_j\eta_k\eta_l+\xi_i\xi_j\xi_k\xi_l\Big)+\sum_{1\leq i,j,k,l\leq N}r_{ijkl}(\L)\eta_i\eta_j\xi_k\xi_l$$
where
\beqa{r and q}
q_{i,j,k,l}(\L)&:=&\left\{
\begin{array}{lrr}
-\frac{\bar m_i}{\L_i^2}\sum_{h:h\neq i}{\bar m_h}\,{a_{4000}(a_i,a_h)}\quad \textrm{for}\quad i=j=k=l\\
\\
-{\bar m_i\,\bar m_l}\frac{a_{3100}(a_i,a_l)}{\L_i\,\sqrt{\L_i\,\L_l}}\quad \textrm{for}\quad i=j=k\neq l\\
\\
-{\bar m_i\,\bar m_l}\frac{a_{2200}(a_i,a_l)}{\L_i\,\L_l}\quad \textrm{for}\quad i=j< k= l\\
\\
0\quad \textrm{otherwise;}
\end{array}
\right.\nonumber\\
r_{i,j,k,l}(\L)&:=&\left\{
\begin{array}{lrr}
-\frac{\bar m_i}{\L_i^2}\sum_{h:h \neq i}{\bar m_h}{a_{2020}(a_i,a_h)}\quad \textrm{for}\quad i=j=k=l\\
\\
-{\bar m_k\,\bar m_i}\frac{a_{0220}(a_k,a_i)}{\L_k\,\L_i}\quad \textrm{for}\quad i=j\neq k=l\\
\\
-{\bar m_i\,\bar m_j}\frac{a_{1120}(a_i,a_j)}{\L_i\,\sqrt{\L_i\,\L_j}}\quad \textrm{for}\quad i=k=l\neq j\\
\\
-{\bar m_i\,\bar m_l}\frac{a_{1120}(a_i,a_l)}{\L_i\,\sqrt{\L_i\,\L_l}}\quad \textrm{for}\quad i=j=k\neq l\\
\\
-{\bar m_i\,\bar m_j}\frac{a_{1111}(a_i,a_j)}{\L_i\,\L_j}\quad \textrm{for}\quad i=k< j=l\\
\\
0\quad \textrm{otherwise.}
\end{array}
\right.
\eeqa
where
\beqa{a4000 and others}
a_{4000}({a},{b})&=&\frac{{a}}{512{b}^2}\,\left[(-60({a}/{b})^5+4311({a}/{b})^3\right.\nonumber\\
&-&300({a}/{b}))\,b_{9/2,0}({a}/{b})+8(7({a}/{b})^{6}\nonumber\\
&-&252({a}/{b})^4-222({a}/{b})^2+7)\,b_{9/2,1}({a}/{b})\nonumber\\
&+&4(75({a}/{b})^{5}-503({a}/{b})^{3}+135({a}/{b}))\,b_{9/2,2}({a}/{b})\nonumber\\
&+&24(23({a}/{b})^{4}+13({a}/{b})^{2})\,b_{9/2,3}({a}/{b})\nonumber\\
&+&\left.37({a}/{b})^{3}\,b_{9/2,4}({a}/{b})\right]\nonumber\\
\nonumber\\
a_{3100}({a},{b})&=&-\frac{{a}}{256{b}^2}\left[(-744({a}/{b})^5+2014({a}/{b})^3\right.\nonumber\\
&-&864({a}/{b}))\,b_{9/2,1}({a}/{b})+8(28({a}/{b})^6\nonumber\\
&-&321({a}/{b})^4-321({a}/{b})^2+28)\,\,b_{9/2,2}({a}/{b})\nonumber\\
&+&(552({a}/{b})^5+423({a}/{b})^3+672({a}/{b}))\,b_{9/2,3}({a}/{b})\nonumber\\
&+&(1146({a}/{b})^4+1266({a}/{b})^2)\,b_{9/2,0}({a}/{b})\nonumber\\
&+&6(29({a}/{b})^4+9({a}/{b})^2)\,b_{9/2,4}({a}/{b})\nonumber\\
&-&\left.5({a}/{b})^3\,b_{9/2,5}({a}/{b})\right]\nonumber\\
\nonumber\\
a_{2200}({a},{b})&=&\frac{{a}}{512{b}^2}\,\left[(-324({a}/{b})^5+10584({a}/{b})^3-324({a}/{b}))\,b_{9/2,0}({a}/{b})\right.\nonumber\\
&+&8(17({a}/{b})^6-300({a}/{b})^4-300({a}/{b})^2+17)\,b_{9/2,1}({a}/{b})\nonumber\\
&-&(1272({a}/{b})^5+6337({a}/{b})^3+1272({a}/{b}))\,b_{9/2,2}({a}/{b})\nonumber\\
&+&(648({a}/{b})^6+396({a}/{b})^4+396({a}/{b})^2\nonumber\\
&+&648)\,b_{9/2,3}({a}/{b})+(348({a}/{b})^5\nonumber\\
&+&800({a}/{b})^3+348({a}/{b}))\,b_{9/2,4}({a}/{b})\nonumber\\
&+&(-60{a}/{b})^4-60({a}/{b})^2)\,b_{9/2,5}({a}/{b})\nonumber\\
&+&\left.9({a}/{b})^3\,b_{9/2,6}({a}/{b})\right]\nonumber\\
\nonumber\\
a_{2020}({a},{b})&=&\frac{{a}}{256{b}^2}\,[8(7({a}/{b})^6-252({a}/{b})^4-222({a}/{b})^2\nonumber\\
&+&7)\,b_{9/2,1}({a}/{b})+(-60({a}/{b})^5+4311({a}/{b})^3\nonumber\\
&-&300({a}/{b}))\,b_{9/2,0}({a}/{b})+4(75({a}/{b})^5-503({a}/{b})^3\nonumber\\
&+&135({a}/{b}))\,b_{9/2,2}({a}/{b})+24(23({a}/{b})^4\nonumber\\
&+&\left.13({a}/{b})^2)\,b_{9/2,3}({a}/{b})+37({a}/{b})^4\,b_{9/2,4}({a}/{b})\right]\nonumber\\
\nonumber\\
a_{1120}({a},{b})&=&-\frac{{a}}{256{b}^2}\,\left[(-744({a}/{b})^5+2014({a}/{b})^3\right.\nonumber\\
&-&864({a}/{b}))\,b_{9/2,1}({a}/{b})+8(28({a}/{b})^6-321({a}/{b})^4\nonumber\\
&-&321({a}/{b})^2+28)\,b_{9/2,2}({a}/{b})\nonumber\\
&+&(552({a}/{b})^5+423({a}/{b})^3+672({a}/{b}))\,b_{9/2,3}({a}/{b})\nonumber\\
&+&(1146({a}/{b})^4+1266({a}/{b})^2)\,b_{9/2,0}({a}/{b})\nonumber\\
&+&6(29({a}/{b})^4+9({a}/{b})^2)\,b_{9/2,4}({a}/{b})\nonumber\\
&-&\left.5({a}/{b})^3\,b_{9/2,5}({a}/{b})\right]\nonumber\\
\nonumber\\
a_{0220}({a},{b})&=&-\frac{3{a}}{512{b}^2}\,\left[(84({a}/{b})^5-8832({a}/{b})^3\right.\nonumber\\
&+&84({a}/{b}))\,b_{9/2,0}({a}/{b})-8(5({a}/{b})^6\nonumber\\
&-&652({a}/{b})^4-652({a}/{b})^2+5)\,b_{9/2,1}({a}/{b})\nonumber\\
&-&5(328({a}/{b})^5-561({a}/{b})^3+328({a}/{b}))\,b_{9/2,2}({a}/{b})\nonumber\\
&+&(216({a}/{b})^6-1020({a}/{b})^4\nonumber\\
&-&1020({a}/{b})^2+216)\,b_{9/2,3}({a}/{b})\nonumber\\
&+&(116({a}/{b})^5+200({a}/{b})^3\nonumber\\
&+&116({a}/{b}))\,b_{9/2,4}({a}/{b})\nonumber\\
&-&(20({a}/{b})^4+20({a}/{b})^2)\,b_{9/2,5}({a}/{b})\nonumber\\
&+&\left.3({a}/{b})^3\,b_{9/2,6}({a}/{b})\right]\nonumber\\
\nonumber\\
a_{1111}({a},{b})&=&\frac{{a}}{128{b}^2}\,\left[(-36({a}/{b})^5-7956({a}/{b})^3\right.\nonumber\\
&-&36({a}/{b}))\,b_{9/2,0}({a}/{b})+8(({a}/{b})^6\nonumber\\
&+&828({a}/{b})^4+828({a}/{b})^2+1)\,b_{9/2,1}({a}/{b})\nonumber\\
&+&(-3096({a}/{b})^5+1039({a}/{b})^3\nonumber\\
&-&3096({a}/{b}))\,b_{9/2,2}({a}/{b})+(648({a}/{b})^6\nonumber\\
&-&1332({a}/{b})^4-1332({a}/{b})^2+648)\,b_{9/2,3}({a}/{b})\nonumber\\
&+&(348({a}/{b})^5+700({a}/{b})^3\nonumber\\
&+&348({a}/{b}))\,b_{9/2,4}({a}/{b})-60(({a}/{b})^4\nonumber\\
&+&\left.({a}/{b})^2)\,b_{9/2,5}({a}/{b})+9({a}/{b})^3\,b_{9/2,6}({a}/{b})\right]
\eeqa
\end{lemma}
{\bf Proof.}\ As in the proof of the previous Lemma, we use the symmetries (\ref{simmetries of coeff}) outlined in Corollary \ref{simmetries remark}. We find, in the fourth order of the function $g_{ij}$, only $19$ (among the $35$ possible ones) non vanishing monomials with degree $4$, wich are individuated by $7$ independent coefficients, say $a_{4000}$, $a_{3100}$, $a_{2200}$, $a_{2020}$, $a_{1120}$, $a_{0220}$ and $a_{1111}$:
\begin{eqnarray*}
&&{a_{4000}(a_i,a_j)}\,\hat\eta_i^4+{a_{3100}(a_i,a_j)}\,\hat\eta_i^3\hat\eta_j+{a_{2200}(a_i,a_j)}\,\hat\eta_i^2\hat\eta_j^2+{a_{3100}(a_j,a_i)}\,\hat\eta_i\hat\eta_j^3\nonumber\\
&+&{a_{4000}(a_j,a_i)}\,\hat\eta_j^4+{a_{2020}(a_i,a_j)}\,\hat\eta_i^2{\hat\xi_i}^2+{a_{1120}(a_i,a_j)}\,\hat\eta_i\hat\eta_j{\hat\xi_i}^2+{a_{0220}(a_i,a_j)}\,\hat\eta_j^2{\hat\xi_i}^2\nonumber\\
&+&{a_{1120}(a_i,a_j)}\,\hat\eta_i^2{\hat\xi_i}{\hat\xi_j}+{a_{1111}(a_i,a_j)}\,\hat\eta_i\hat\eta_j{\hat\xi_i}{\hat\xi_j}+{a_{1120}(a_j,a_i)}\,\hat\eta_j^2{\hat\xi_i}{\hat\xi_j}+{a_{0220}(a_j,a_i)}\,\hat\eta_i^2{\hat\xi_j}^2\nonumber\\
&+&{a_{1120}(a_j,a_i)}\,\hat\eta_i\hat\eta_j{\hat\xi_j}^2+{a_{2020}(a_j,a_i)}\,\hat\eta_j^2{\hat\xi_j}^2+{a_{4000}(a_i,a_j)}\,{\hat\xi_i}^4+{a_{3100}(a_i,a_j)}\,{\hat\xi_i}^3{\hat\xi_j}\nonumber\\
&+&{a_{2200}(a_i,a_j)}\,{\hat\xi_i}^2{\hat\xi_j}^2+{a_{3100}(a_j,a_i)}\,{\hat\xi_i}{\hat\xi_j}^3+{a_{4000}(a_j,a_i)}\,{\hat\xi_j}^4\ ,\nonumber\\
\end{eqnarray*}
Thus, multiplying by $-\bar m_i\bar m_j$, and summing over all $1\leq i<j\leq N$, we find (\ref{r and q}).

\vskip.2in
\noi
We perform now the computation
of the $7$ coefficients (here, $a_1$, $a_2$, $\cdots$ are used as ``dummy'' variables)
\begin{eqnarray}\label{needed derivatives}
{a_{4000}(a_1,a_2)}&=&\left.\frac{1}{24}\,\partial_{\hat{\eta}_1^4}\hat g(a_1,a_2,\hat{\eta}_1,\hat{\eta}_2,\hat{\xi}_1,\hat{\xi}_2)\right|_{0}\nonumber\\
a_{3100}(a_1,a_2)&=&\left.\frac{1}{6}\,\partial_{\hat{\eta}_1^3\hat{\eta}_2}\hat g(a_1,a_2,\hat{\eta}_1,\hat{\eta}_2,\hat{\xi}_1,\hat{\xi}_2)\right|_{0}\nonumber\\
a_{2200}(a_1,a_2)&=&\left.\frac{1}{4}\,\partial_{\hat{\eta}_1^2\hat{\eta}_2^2}\hat g(a_1,a_2,\hat{\eta}_1,\hat{\eta}_2,\hat{\xi}_1,\hat{\xi}_2)\right|_{0}\nonumber\\
a_{2020}(a_1,a_2)&=&\left.\frac{1}{4}\,\partial_{\hat{\eta}_1^2\hat{\xi}_1^2}\hat g(a_1,a_2,\hat{\eta}_1,\hat{\eta}_2,\hat{\xi}_1,\hat{\xi}_2)\right|_{0}\nonumber\\
a_{0220}(a_1,a_2)&=&\left.\frac{1}{4}\,\partial_{\hat{\eta}_2^2\hat{\xi}_1^2}\hat g(a_1,a_2,\hat{\eta}_1,\hat{\eta}_2,\hat{\xi}_1,\hat{\xi}_2)\right|_{0}\nonumber\\
a_{1120}(a_1,a_2)&=&\left.\frac{1}{2}\,\partial_{\hat{\eta}_1\hat{\eta}_2\hat{\xi}_1^2}\hat g(a_1,a_2,\hat{\eta}_1,\hat{\eta}_2,\hat{\xi}_1,\hat{\xi}_2)\right|_{0}\nonumber\\
a_{1111}(a_1,a_2)&=&\left.\partial_{\hat{\eta}_1\hat{\eta}_2\hat{\xi}_1\hat{\xi}_2}\hat g(a_1,a_2,\hat{\eta}_1,\hat{\eta}_2,\hat{\xi}_1,\hat{\xi}_2)\right|_{0}\nonumber\\
\end{eqnarray}
where $|_0$ stands for $|_{(\hat{\eta}_1,\hat{\eta}_2,\hat{\xi}_1,\hat{\xi}_2)=0}$ and
\begin{eqnarray*}
&&\hat g(a_1,a_2,\hat{\eta}_1,\hat{\eta}_2,\hat{\xi}_1,\hat{\xi}_2)=\frac{1}{(2\pi)^2}\,\int_{[0,2\pi]^2}\,\frac{d\l_1\,d\l_2}{{\hat d(a_1,a_2,\l_1,\l_2,\hat{\eta}_1,\hat{\eta}_2,\hat{\xi}_1,\hat{\xi}_2)}}\ .
\end{eqnarray*}
We write
$$\hat d(a_1,a_2,\l_1,\l_2,\hat{\eta}_1,\hat{\eta}_2,\hat{\xi}_1,\hat{\xi}_2)^2=\sqrt{\hat d_2(a_1,a_2,\l_1,\l_2,\hat{\eta}_1,\hat{\eta}_2,\hat{\xi}_1,\hat{\xi}_2)}$$
where
\begin{eqnarray}\label{square distance}
{\hat d_2(a_1,a_2,\l_1,\l_2,\hat{\eta}_1,\hat{\eta}_2,\hat{\xi}_1,\hat{\xi}_2)}&:=&|\hat x(a_1,\l_1,\hat{\eta}_1,\hat{\xi}_1)-\hat x(a_2,\l_2,\hat{\eta}_2,\hat{\xi}_2)|^2\nonumber\\
&=&\hat x_1(a_1,\l_1,\hat{\eta}_1,\hat{\xi}_1)^2+\hat x_2(a_1,\l_1,\hat{\eta}_1,\hat{\xi}_1)^2\nonumber\\
&+&\hat x_1(a_2,\l_2,\hat{\eta}_2,\hat{\xi}_2)^2+\hat x_2(a_2,\l_2,\hat{\eta}_2,\hat{\xi}_2)^2\nonumber\\
&-&2\,\hat x_1(a_1,\l_1,\hat{\eta}_1,\hat{\xi}_1)\,\hat x_1(a_2,\l_2,\hat{\eta}_2,\hat{\xi}_2)\nonumber\\
&-&2\,\hat x_2(a_1,\l_1,\hat{\eta}_1,\hat{\xi}_1)\,\hat x_2(a_2,\l_2,\hat{\eta}_2,\hat{\xi}_2)\nonumber\\
\end{eqnarray}
where, for short, $$(a,\l,\hat\eta,\hat\xi)\to \hat x(a,\l,\hat\eta,\hat\xi)=\Big(\hat x_1(a,\l,\hat\eta,\hat\xi),\hat x_2(a,\l,\hat\eta,\hat\xi)\Big)$$
denotes the Delaunay--Poincar\'e map for $N=1$. 
Then, by usual calculus rules
\begin{eqnarray*}
\partial_{\zeta_1\zeta_2\zeta_3\zeta_4}\frac{1}{\hat d}&=&(16\,\hat d_2^{9/2})^{-1}\nonumber\\
&\times&\Big[105\,\partial_{\zeta_1}\hat d_2\,\partial_{\zeta_2}\hat d_2\,\partial_{\zeta_3}\hat d_2\,\partial_{\zeta_4}\hat d_2-30\,\hat d_2\,\partial_{\zeta_1\zeta_2}\hat d_2\,\partial_{\zeta_3}\hat d_2\,\partial_{\zeta_4}\hat d_2\nonumber\\
&-&30\,\hat d_2\,\partial_{\zeta_1\zeta_3}\hat d_2\,\partial_{\zeta_2}\hat d_2\,\partial_{\zeta_4}\hat d_2-30\,\hat d_2\,\partial_{\zeta_1\zeta_4}\hat d_2\,\partial_{\zeta_2}\hat d_2\,\partial_{\zeta_3}\hat d_2\nonumber\\
&-&30\,\hat d_2\,\partial_{\zeta_2\zeta_3}\hat d_2\,\partial_{\zeta_1}\hat d_2\,\partial_{\zeta_4}\hat d_2-30\,\hat d_2\,\hat d_2\,\partial_{\zeta_2\zeta_4}\hat d_2\,\partial_{\zeta_1}\hat d_2\,\partial_{\zeta_3}\nonumber\\
&+&-30\,\hat d_2\,\partial_{\zeta_3\zeta_4}\hat d_2\,\partial_{\zeta_1}\hat d_2\,\partial_{\zeta_2}\hat d_2+12\,\hat d_2^2\,\partial_{\zeta_2\zeta_3\zeta_4}\hat d_2\,\partial_{\zeta_1}\hat d_2\nonumber\\
&+&12\,\hat d_2^2\,\partial_{\zeta_2}\hat d_2\,\partial_{\zeta_1\zeta_3\zeta_4}\hat d_2+12\,\hat d_2^2\,\partial_{\zeta_3}\hat d_2\,\partial_{\zeta_1\zeta_2\zeta_4}\hat d_2\nonumber\\
&+&12\,\hat d_2^2\,\partial_{\zeta_4}\hat d_2\,\partial_{\zeta_1\zeta_2\zeta_3}\hat d_2+12\,\hat d_2^2\,\partial_{\zeta_1\zeta_2}\hat d_2\,\partial_{\zeta_3\zeta_4}\hat d_2\nonumber\\
&+&12\,\hat d_2^2\,\partial_{\zeta_1\zeta_3}\hat d_2\partial_{\zeta_2\zeta_4}\hat d_2+12\,\hat d_2^2\,\partial_{\zeta_1\zeta_4}\hat d_2\,\partial_{\zeta_2\zeta_3}\hat d_2\nonumber\\
&-&8\,\hat d_2^3\,\partial_{\zeta_1\zeta_2\zeta_3\zeta_4}\hat d_2\Big]\nonumber\\
\end{eqnarray*} 
We find, then, 
\begin{eqnarray}\label{derivatives od inversdist}
\left.\partial_{\hat{\eta}_1^4}\frac{1}{\hat d}\right|_0&=&\Big(16\,\hat d_2^{9/2}\Big)^{-1}\nonumber\\
&\times&\Big[105\,(\partial_{\hat{\eta}_1}\hat d_2)^4-180\,\hat d_2(\partial_{\hat{\eta}_1}\hat d_2)^2\,\partial_{\hat{\eta}_1^2}\hat d_2\nonumber\\
&+&48\,\hat d_2^2\,\partial_{\hat{\eta}_1}\hat d_2\,\partial_{\hat{\eta}_1^3}\hat d_2+36\,\hat d_2^2\,(\partial_{\hat{\eta}_1^2}\hat d_2)^2\nonumber\\
&-&\left.8\,\hat d_2^3\,\partial_{\hat{\eta}_1^4}\hat d_2\right.\Big]_0\nonumber\\
\nonumber\\
\left.\partial_{\hat{\eta}_1^3\hat{\eta}_2}\frac{1}{\hat d}\right|_0&=&\Big(16\,\hat d_2^{9/2}\Big)^{-1}\nonumber\\
&\times&\Big[105\,\partial_{\hat{\eta}_2}\hat d_2\,(\partial_{\hat{\eta}_1}\hat d_2)^3-90\,\hat d_2\,(\partial_{\hat{\eta}_1}\hat d_2)^2\,\partial_{\hat{\eta}_1\hat{\eta}_2}\hat d_2\nonumber\\
&-&90\,\hat d_2\,\partial_{\hat{\eta}_2}\hat d_2\,\partial_{\hat{\eta}_1}\hat d_2\partial_{\hat{\eta}_1^2}\hat d_2\nonumber\\
&+&36\,\hat d_2^2\partial_{\hat{\eta}_1}\hat d_2\partial_{\hat{\eta}_1^2\hat{\eta}_2}\hat d_2+12\,\hat d_2^2\partial_{\hat{\eta}_2}\hat d_2\,\partial_{\hat{\eta}_1^3}\hat d_2\nonumber\\
&+&36\,\hat d_2^2\,\partial_{\hat{\eta}_1\hat{\eta}_2}\hat d_2\partial_{\hat{\eta}_1^2}\hat d_2-8\,\hat d_2^3\,\partial_{\hat{\eta}_1^3\hat{\eta}_2}\hat d_2\Big]_0\nonumber\\
\nonumber\\
\left.\partial_{\hat{\eta}_1^2\hat{\eta}_2^2}\frac{1}{\hat d}\right|_0&=&\Big(16\,\hat d_2^{9/2}\Big)^{-1}\nonumber\\&\times&\Big[105\,(\partial_{\hat{\eta}_2}\hat d_2)^2\,(\partial_{\hat{\eta}_1}\hat d_2)^2-30\,\hat d_2\,\partial_{\hat{\eta}_2^2}\hat d_2\,(\partial_{\hat{\eta}_1}\hat d_2)^2\nonumber\\
&-&120\,\hat d_2\,\partial_{\hat{\eta}_1}\hat d_2\,\partial_{\hat{\eta}_2}\hat d_2\,\partial_{\hat{\eta}_1\hat{\eta}_2}\hat d_2-30\,\hat d_2\,(\partial_{\hat{\eta}_2}\hat d_2)^2\,\partial_{\hat{\eta}_1^2}\hat d_2\nonumber\\
&+&24\,\hat d_2^2\,\partial_{\hat{\eta}_2}\hat d_2\,\partial_{\hat{\eta}_1^2\hat{\eta}_2}\hat d_2+24\,\hat d_2^2\,\partial_{\hat{\eta}_1}\hat d_2\partial_{\hat{\eta}_1\hat{\eta}_2^2}\hat d_2\nonumber\\
&+&\left.12\,\hat d_2^2\,\partial_{\hat{\eta}_2^2}\hat d_2\,\partial_{\hat{\eta}_1^2}\hat d_2+24\,\hat d_2^2\,(\partial_{\hat{\eta}_1\hat{\eta}_2}\hat d_2)^2\right.\nonumber\\
&-&\left.8\,\hat d_2^3\,\partial_{\hat{\eta}_1^2\hat{\eta}_2^2}\hat d_2\right.\Big]_0\nonumber\\
\nonumber\\
\left.\partial_{\hat{\eta}_1^2\hat{\xi}_1^2}\frac{1}{\hat d}\right|_0&=&\Big(16\,\hat d_2^{9/2}\Big)^{-1}\nonumber\\&\times&\Big[105\,(\partial_{\hat{\xi}_1}\hat d_2)^2\,(\partial_{\hat{\eta}_1}\hat d_2)^2-30\,\hat d_2\,\partial_{\hat{\xi}_1^2}\hat d_2\,(\partial_{\hat{\eta}_1}\hat d_2)^2\nonumber\\
&-&120\,\hat d_2\,\partial_{\hat{\xi}_1}\hat d_2\,\partial_{\hat{\eta}_1}\hat d_2\,\partial_{\hat{\eta}_1\hat{\xi}_1}\hat d_2-30\,\hat d_2(\partial_{\hat{\xi}_1}\hat d_2)^2\,\partial_{\hat{\eta}_1^2}\hat d_2\nonumber\\
&+&24\,\hat d_2^2\partial_{\hat{\xi}_1}\hat d_2\,\partial_{\hat{\eta}_1^2\hat{\xi}_1}\hat d_2+24\,\hat d_2^2\partial_{\hat{\eta}_1}\hat d_2\,\partial_{\hat{\eta}_1\hat{\xi}_1^2}\hat d_2\nonumber\\
&+&\left.12\,\hat d_2^2\,\partial_{\hat{\xi}_1^2}\hat d_2\,\partial_{\hat{\eta}_1^2}\hat d_2+24\,\hat d_2^2\,(\partial_{\hat{\eta}_1\hat{\xi}_1}\hat d_2)^2-8\,\hat d_2^3\,\partial_{\hat{\eta}_1^2\hat{\xi}_1^2}\hat d_2\right.\Big]_0\nonumber\\
\nonumber\\
\left.\partial_{\hat{\eta}_2^2\hat{\xi}_1^2}\frac{1}{\hat d}\right|_0&=&\Big(16\,\hat d_2^{9/2}\Big)^{-1}\nonumber\\&\times&\Big[105\,(\partial_{\hat{\xi}_1}\hat d_2)^2\,(\partial_{\hat{\eta}_2}\hat d_2)^2-30\,\hat d_2\,\partial_{\hat{\xi}_1^2}\hat d_2\,(\partial_{\hat{\eta}_2}\hat d_2)^2\nonumber\\
&-&120\,\hat d_2\,\partial_{\hat{\xi}_1}\hat d_2\,\partial_{\hat{\eta}_2}\hat d_2\,\partial_{\hat{\xi}_1\hat{\eta}_2}\hat d_2-30\,\hat d_2(\partial_{\hat{\xi}_1}\hat d_2)^2\,\partial_{\hat{\eta}_2^2}\hat d_2\nonumber\\
&+&24\,\hat d_2^2\,\partial_{\hat{\xi}_1}\hat d_2\,\partial_{\hat{\eta}_2^2\hat{\xi}_1}\hat d_2+24\,\hat d_2^2\,\partial_{\hat{\eta}_2}\hat d_2\,\partial_{\hat{\eta}_2\hat{\xi}_1^2}\hat d_2\nonumber\\
&+&\left.12\,\hat d_2^2\,\partial_{\hat{\xi}_1^2}\hat d_2\,\partial_{\hat{\eta}_2^2}\hat d_2+24\,\hat d_2^2\,(\partial_{\hat{\xi}_1\hat{\eta}_2}\hat d_2)^2-8\,\hat d_2^3\,\partial_{\hat{\eta}_2^2\hat{\xi}_1^2}\hat d_2\right.\Big]_0\nonumber\\
\nonumber\\
\left.\partial_{\hat{\eta}_1\hat{\eta}_2\hat{\xi}_1^2}\frac{1}{\hat d}\right|_0&=&\Big(16\,\hat d_2^{9/2}\Big)^{-1}\nonumber\\&\times&\Big[105\,\partial_{\hat{\eta}_1}\hat d_2\,\partial_{\hat{\eta}_2}\hat d_2\,(\partial_{\hat{\xi}_1}\hat d_2)^2-30\,\hat d_2\,\partial_{\hat{\eta}_1}\hat d_2\,\partial_{\hat{\eta}_2}\hat d_2\,\partial_{\hat{\xi}_1^2}\hat d_2\nonumber\\
&-&60\,\hat d_2\,\partial_{\hat{\xi}_1}\hat d_2\,\partial_{\hat{\eta}_2\hat{\xi}_1}\hat d_2\,\partial_{\hat{\eta}_1}\hat d_2-60\,\hat d_2\,\partial_{\hat{\xi}_1}\hat d_2\,\partial_{\hat{\eta}_2}\hat d_2\,\partial_{\hat{\eta}_1\hat{\xi}_1}\hat d_2\nonumber\\
&-&30\,\hat d_2\,(\partial_{\hat{\xi}_1}\hat d_2)^2\,\partial_{\hat{\eta}_1\hat{\eta}_2}\hat d_2+12\,\hat d_2^2\,\partial_{\hat{\eta}_2\hat{\xi}_1^2}\hat d_2\,\partial_{\hat{\eta}_1}\hat d_2+24\,\hat d_2^2\,\partial_{\hat{\xi}_1}\hat d_2\,\partial_{\hat{\eta}_1\hat{\eta}_2\hat{\xi}_1}\hat d_2\nonumber\\
&+&12\,\hat d_2^2\partial_{\hat{\eta}_2}\hat d_2\,\partial_{\hat{\eta}_1\hat{\xi}_1^2}\hat d_2+24\,\hat d_2^2\,\partial_{\hat{\eta}_2\hat{\xi}_1}\hat d_2\,\partial_{\hat{\eta}_1\hat{\xi}_1}\hat d_2\nonumber\\
&+&\left.12\,\hat d_2^2\,\partial_{\hat{\xi}_1^2}\hat d_2\,\partial_{\hat{\eta}_1\hat{\eta}_2}\hat d_2-8\,\hat d_2^3\partial_{\hat{\eta}_1\hat{\eta}_2\hat{\xi}_1^2}\hat d_2\right.\Big]_0\nonumber\\
\nonumber\\
\left.\partial_{\hat{\eta}_1\hat{\eta}_2\hat{\xi}_1\hat{\xi}_2}\frac{1}{{\hat d}}\right|_0&=&\Big(16\,\hat d_2^{9/2}\Big)^{-1}\nonumber\\&\times&\Big[105\,\partial_{\hat{\eta}_1}\hat d_2\,\partial_{\hat{\eta}_2}\hat d_2\,\partial_{\hat{\xi}_1}\hat d_2\,\partial_{\hat{\xi}_2}\hat d_2-30\,\hat d_2\,\partial_{\hat{\eta}_1\hat{\eta}_2}\hat d_2\,\partial_{\hat{\xi}_1}\hat d_2\,\partial_{\hat{\xi}_2}\hat d_2\nonumber\\
&-&30\,\hat d_2\,\partial_{\hat{\eta}_1\hat{\xi}_1}\hat d_2\,\partial_{\hat{\eta}_2}\hat d_2\,\partial_{\hat{\xi}_2}\hat d_2-30\,\hat d_2\,\partial_{\hat{\eta}_1\hat{\xi}_2}\hat d_2\,\partial_{\hat{\eta}_2}\hat d_2\,\partial_{\hat{\xi}_1}\hat d_2\nonumber\\
&-&30\,\hat d_2\,\partial_{\hat{\eta}_2\hat{\xi}_1}\hat d_2\,\partial_{\hat{\eta}_1}\hat d_2\,\partial_{\hat{\xi}_2}\hat d_2-30\,\hat d_2\,\partial_{\hat{\eta}_2\hat{\xi}_2}\hat d_2\,\partial_{\hat{\eta}_1}\hat d_2\,\partial_{\hat{\xi}_1}\nonumber\\
&-&30\,\hat d_2\,\partial_{\hat{\xi}_1\hat{\xi}_2}\hat d_2\,\partial_{\hat{\eta}_1}\hat d_2\,\partial_{\hat{\eta}_2}\hat d_2+12\,\hat d_2^2\,\partial_{\hat{\eta}_2\hat{\xi}_1\hat{\xi}_2}\hat d_2\,\partial_{\hat{\eta}_1}\hat d_2+12\,\hat d_2^2\,\partial_{\hat{\eta}_2}\hat d_2\,\partial_{\hat{\eta}_1\hat{\xi}_1\hat{\xi}_2}\hat d_2\nonumber\\
&+&12\,\hat d_2^2\,\partial_{\hat{\xi}_1}\hat d_2\,\partial_{\hat{\eta}_1\hat{\eta}_2\hat{\xi}_2}\hat d_2+12\,\hat d_2^2\,\partial_{\hat{\xi}_2}\hat d_2\,\partial_{\hat{\eta}_1\hat{\eta}_2\hat{\xi}_1}\hat d_2\nonumber\\
&+&12\,\hat d_2^2\,\partial_{\hat{\eta}_1\hat{\eta}_2}\hat d_2\,\partial_{\hat{\xi}_1\hat{\xi}_2}\hat d_2+\left.12\,\hat d_2^2\,\partial_{\hat{\eta}_1\hat{\xi}_1}\hat d_2\partial_{\hat{\eta}_2\hat{\xi}_2}\hat d_2\right.\nonumber\\
&+&\left.12\,\hat d_2^2\,\partial_{\hat{\eta}_1\hat{\xi}_2}\hat d_2\,\partial_{\hat{\eta}_2\hat{\xi}_1}\hat d_2-8\,\hat d_2^3\,\partial_{\hat{\eta}_1\hat{\eta}_2\hat{\xi}_1\hat{\xi}_2}\hat d_2\right.\Big]_0\nonumber\\
\end{eqnarray}
where
$$\left.\hat d_2\right|_0=a_1^2+a_2^2-2a_1a_2\cos{(\l_1-\l_2)}$$
We start by computing the following derivatives in terms of the derivatives of the coordinates
\begin{eqnarray}\label{d2expansion}
\left.\partial_{\hat{\hat{\eta}}_1}\hat d_2\right|_0&=&\left.2\,[(\hat x_{11}-\hat x_{12})\,\partial_{\hat{\hat{\eta}}_1}\hat x_{11}+(\hat x_{21}-\hat x_{22})\,\partial_{\hat{\eta}_1}\hat x_{21}]\right|_0\nonumber\\ \nonumber\\
\left.\partial_{\hat{\eta}_1^2}\hat d_2\right|_0&=&\left.2\,[(\partial_{\hat{\eta}_1}\hat x_{11})^2+(\partial_{\hat{\eta}_1}\hat x_{21})^2+(\hat x_{11}-\hat x_{12})\,\partial_{\hat{\eta}_1^2}\hat x_{11}+(\hat x_{21}-\hat x_{22})\,\partial_{\hat{\eta}_1^2}\hat x_{21}]\right|_0\nonumber\\ \nonumber\\
\left.\partial_{\hat{\eta}_1\hat{\eta}_2}\hat d_2\right|_0&=&\left.-2\,[\partial_{\hat{\eta}_1}\hat x_{11}\,\partial_{\hat{\eta}_2}\hat x_{12}+\partial_{\hat{\eta}_1}\hat x_{21}\,\partial_{\hat{\eta}_2}\hat x_{22}]\right|_0\nonumber\\ \nonumber\\
\left.\partial_{\hat{\eta}_1\hat{\hat{\xi}}_1}\hat d_2\right|_0&=&2\,[\partial_{\hat{\hat{\xi}}_1}\hat x_{11}\,\partial_{\hat{\eta}_1}\hat x_{11}+\partial_{\hat{\xi}_1}\hat x_{21}\,\partial_{\hat{\eta}_1}\hat x_{21}+(\hat x_{11}-\hat x_{12})\,\partial_{\hat{\eta}_1\hat{\xi}_1}\hat x_{11}\nonumber\\
&+&(\hat x_{21}-\hat x_{22})\,\partial_{\hat{\eta}_1\hat{\xi}_1}\hat x_{21}]\nonumber\\ \nonumber\\
\left.\partial_{\hat{\eta}_2\hat{\xi}_1}\hat d_2\right|_0&=&\left.-2\,[\partial_{\hat{\xi}_1}\hat x_{11}\,\partial_{\hat{\eta}_2}\hat x_{12}+\partial_{\hat{\xi}_1}\hat x_{21}\,\partial_{\hat{\eta}_2}\hat x_{22}]\right|_0\nonumber\\ \nonumber\\
\left.\partial_{\hat{\eta}_1^3}\hat d_2\right|_0&=&\left.2\,[3\,\partial_{\hat{\eta}_1}\hat x_{11}\,\partial_{\hat{\eta}_1^2}\hat x_{11}+3\,\partial_{\hat{\eta}_1}\hat x_{21}\,\partial_{\hat{\eta}_1^2}\hat x_{21}\right.\nonumber\\
&+&\left.(\hat x_{11}-\hat x_{12})\,\partial_{\hat{\eta}_1^3}\hat x_{11}+(\hat x_{21}-\hat x_{22})\,\partial_{\hat{\eta}_1^3}\hat x_{21}]\right|_0\nonumber\\ \nonumber\\
\left.\partial_{\hat{\eta}_1^2\hat{\eta}_2}\hat d_2\right|_0&=&\left.-2\,[\partial_{\hat{\eta}_2}\hat x_{12}\,\partial_{\hat{\eta}_1^2}\hat x_{11}+\partial_{\hat{\eta}_2}\hat x_{22}\,\partial_{\hat{\eta}_1^2}\hat x_{21}]\right|_0\nonumber\\ \nonumber\\
\left.\partial_{\hat{\eta}_1\hat{\xi}_1^2}\hat d_2\right|_0&=&2\,[\partial_{\hat{\xi}_1^2}\hat x_{11}\,\partial_{\hat{\eta}_1}\hat x_{11}+\partial_{\hat{\xi}_1^2}\hat x_{21}\,\partial_{\hat{\eta}_1}\hat x_{21}+2\,\partial_{\hat{\xi}_1}\hat x_{11}\,\partial_{\hat{\xi}_1\hat{\eta}_1}\hat x_{11}+2\,\partial_{\hat{\xi}_1}\hat x_{21}\,\partial_{\hat{\xi}_1\hat{\eta}_1}\hat x_{21}\nonumber\\
&+&\left.(\hat x_{11}-\hat x_{12})\,\partial_{\hat{\eta}_1\hat{\xi}_1^2}\hat x_{11}+(\hat x_{21}-\hat x_{22})\,\partial_{\hat{\eta}_1\hat{\xi}_1^2}\hat x_{21}]\right|_0\nonumber\\ \nonumber\\
\left.\partial_{\hat{\eta}_2\hat{\xi}_1^2}\hat d_2\right|_0&=&\left.-2\,[\partial_{\hat{\xi}_1^2}\hat x_{11}\,\partial_{\hat{\eta}_2}\hat x_{12}+\partial_{\hat{\xi}_1^2}\hat x_{21}\,\partial_{\hat{\eta}_2}\hat x_{22}]\right|_0\nonumber\\ \nonumber\\
\left.\partial_{\hat{\eta}_1\hat{\eta}_2\hat{\xi}_1}\hat d_2\right|_0&=&\left.-2\,[\partial_{\hat{\eta}_2}\hat x_{12}\,\partial_{\hat{\eta}_1,\hat{\xi}_1}\hat x_{11}+\partial_{\hat{\eta}_2}\hat x_{22}\,\partial_{\hat{\eta}_1,\hat{\xi}_1}\hat x_{21}]\right|_0\nonumber\\ \nonumber\\
\left.\partial_{\hat{\eta}_1^4}\hat d_2\right|_0&=&2\,[3\,(\partial_{\hat{\eta}_1^2}\hat x_{11})^2+3\,(\partial_{\hat{\eta}_1^2}\hat x_{21})^2+4\,\partial_{\hat{\eta}_1}\hat x_{11}\,\partial_{\hat{\eta}_1^3}\hat x_{11}+4\,\partial_{\hat{\eta}_1}\hat x_{21}\,\partial_{\hat{\eta}_1^3}\hat x_{21}\nonumber\\
&+&\left.(\hat x_{11}-\hat x_{12})\,\partial_{\hat{\eta}_1^4}\hat x_{11}+(\hat x_{21}-\hat x_{22})\,\partial_{\hat{\eta}_1^4}\hat x_{21}]\right|_0\nonumber\\ \nonumber\\
\left.\partial_{\hat{\eta}_1^3\hat{\eta}_2}\hat d_2\right|_0&=&-2\,[\partial_{\hat{\eta}_1^3}\hat x_{11}\,\partial_{\hat{\eta}_2}\hat x_{12}+\partial_{\hat{\eta}_1^3}\hat x_{21}\,\partial_{\hat{\eta}_2}\hat x_{22}]\nonumber\\ \nonumber\\
\left.\partial_{\hat{\eta}_1^2\hat{\eta}_2^2}\hat d_2\right|_0&=&\left.-2\,[\partial_{\hat{\eta}_1^2}\hat x_{11}\,\partial_{\hat{\eta}_2^2}\hat x_{12}+\partial_{\hat{\eta}_1^2}\hat x_{21}\,\partial_{\hat{\eta}_2^2}\hat x_{22}]\right|_0\nonumber\\\nonumber\\
\left.\partial_{\hat{\eta}_1^2\hat{\xi}_1^2}\hat d_2\right|_0&=&2\,[2\,(\partial_{\hat{\eta}_1\hat{\xi}_1}\hat x_{11})^2+2\,(\partial_{\hat{\eta}_1\hat{\xi}_1}\hat x_{21})^2+2\,\partial_{\hat{\eta}_1}\hat x_{11}\,\partial_{\hat{\eta}_1\hat{\xi}_1^2}\hat x_{11}+2\,\partial_{\hat{\eta}_1}\hat x_{21}\,\partial_{\hat{\eta}_1\hat{\xi}_1^2}\hat x_{21}\nonumber\\
&+&\partial_{\hat{\xi}_1^2}\hat x_{11}\,\partial_{\hat{\eta}_1^2}\hat x_{11}+\partial_{\hat{\xi}_1^2}\hat x_{21}\,\partial_{\hat{\eta}_1^2}\hat x_{21}+2\,\partial_{\hat{\xi}_1}\hat x_{11}\,\partial_{\hat{\eta}_1^2\hat{\xi}_1}\hat x_{11}+2\,\partial_{\hat{\xi}_1}\hat x_{21}\,\partial_{\hat{\eta}_1^2\hat{\xi}_1}\hat x_{21}\nonumber\\
&+&\left.(\hat x_{11}-\hat x_{12})\,\partial_{\hat{\eta}_1^2\hat{\xi}_1^2}\hat x_{11}+(\hat x_{21}-\hat x_{22})\,\partial_{\hat{\eta}_1^2\hat{\xi}_1^2}\hat x_{21}\right|_0\nonumber\\\nonumber\\
\left.\partial_{\hat{\eta}_2^2\hat{\xi}_1^2}\hat d_2\right|_0&=&\left.-2\,[\partial_{\hat{\xi}_1^2}\hat x_{11}\,\partial_{\hat{\eta}_2^2}\hat x_{12}+\partial_{\hat{\xi}_1^2}\hat x_{21}\,\partial_{\hat{\eta}_2^2}\hat x_{22}]\right|_0\nonumber\\\nonumber\\
\left.\partial_{\hat{\eta}_1\hat{\eta}_2\hat{\xi}_1^2}\hat d_2\right|_0&=&\left.-2\,[\partial_{\hat{\eta}_2}\hat x_{12}\,\partial_{\hat{\eta}_1\hat{\xi}_1^2}\hat x_{11}+\partial_{\hat{\eta}_2}\hat x_{22}\,\partial_{\hat{\eta}_1\hat{\xi}_1^2}\hat x_{21}]\right|_0\nonumber\\ \nonumber\\
\left.\partial_{\hat{\eta}_1\hat{\eta}_2\hat{\xi}_1\hat{\xi}_2}\hat d_2\right|_0&=&\left.-2\,[\partial_{\hat{\eta}_1\hat{\xi}_1}\hat x_{11}\,\partial_{\hat{\eta}_2\hat{\xi}_2}\hat x_{12}+\partial_{\hat{\eta}_1\hat{\xi}_1}\hat x_{21}\,\partial_{\hat{\eta}_2\hat{\xi}_2}\hat x_{22}]\right|_0\nonumber\\ \nonumber\\
\end{eqnarray}
where we have denoted, for shortness,
\begin{eqnarray}\label{notation}
\hat x_{ij}:=\hat x_i(a_j,\l_j,\hat{\eta}_j,\hat{\xi}_j)\ ,\qquad i,\ j=1,\, 2\ .
\end{eqnarray}
Recalling the relation
\begin{eqnarray}\label{second coordinate}
&&\hat x_2(a,\l,\hat{\eta},\hat{\xi})=-\hat x_1(a,\l+\pi/2,\hat{\xi},-\hat{\eta})\ ,
\end{eqnarray}
the computation is then reduced to the one of the involved derivatives of the coordinates $q_{1}$, which we recall, is defined as
\begin{eqnarray}\label{poincare}
\hat x_1(a,\l,\hat{\eta},\hat{\xi})=a\,\left[\cos{(\hat{\zeta}+\l)}-\frac{\hat{\xi}}{2}\left(\hat{\eta}\sin{(\hat{\zeta}+\l)}+\hat{\xi}\cos{(\hat{\zeta}+\l)}\right)-{\hat{\eta}}\sqrt{1-\frac{\hat{\eta}^2+\hat{\xi}^2}{4}}\right]
\nonumber\\
\end{eqnarray}
where  $\hat{\zeta}=\hat{\zeta}(\l,\hat{\eta},\hat{\xi})$ is implicitely defined as the solution of
\begin{eqnarray}\label{bar zeta}
\hat{\zeta}=\sqrt{1-\frac{\hat{\eta}^2+\hat{\xi}^2}{4}}\,\left[(\hat{\eta}\cos{\l}-\hat{\xi}\sin{\l})\sin{\hat{\zeta}}+(\hat{\eta}\sin{\l}+\hat{\xi}\cos{\l})\cos{\hat{\zeta}}\right]\ .
\end{eqnarray} 
To this end, put
\begin{eqnarray}\label{st}
s(\l,\hat{\eta},\hat{\xi}):=\sqrt{1-\frac{\hat{\eta}^2+\hat{\xi}^2}{4}}\,(\hat{\eta}\cos{\l}-\hat{\xi}\sin{\l})\nonumber\\
t(\l,\hat{\eta},\hat{\xi}):=\sqrt{1-\frac{\hat{\eta}^2+\hat{\xi}^2}{4}}\,(\hat{\eta}\sin{\l}+\hat{\xi}\cos{\l})\nonumber\\
\end{eqnarray}
Write, from (\ref{bar zeta}), $\hat{\zeta}(\l,\hat{\eta},\hat{\xi})=\tilde{\zeta}(s(\l,\hat{\eta},\hat{\xi}),t(\l,\hat{\eta},\hat{\xi}))$, where $\tilde{\zeta}(s,t)$ is the solution of
\begin{eqnarray}\label{bar zeta1}
{\zeta}=s\,\sin{{\zeta}}+t\,\cos{{\zeta}}\ .
\end{eqnarray}
The $4$--order expansion of $\tilde{\zeta}$ around $(s,t)=0$ gives 
\begin{eqnarray}\label{bar zeta2}
\tilde{\zeta}(s,t)=t+s\,t+s^2\,t-\frac{1}{2}\,t^3+s^3\,t-\frac{5}{3}\,s\,t^3+O_5(s,t)
\end{eqnarray}
so that, inserting the previous expression into the $4$--order developing of 
\begin{eqnarray*}
\cos{(\l+\tilde{\zeta}(s,t))}&=&\cos{\l}-(\sin{\l})\,\tilde{\zeta}(s,t)-\frac{1}{2}(\cos{\l})\,\tilde{\zeta}(s,t)^2+\frac{1}{6}(\sin{\l})\,\tilde{\zeta}(s,t)^3\nonumber\\
&+&\frac{1}{24}(\cos{\l})\,\tilde{\zeta}(s,t)^4+O_5(\tilde{\zeta}(s,t))\nonumber\\
\end{eqnarray*}
we find
\begin{eqnarray}\label{CosAndSin}
\cos{(\l+\tilde{\zeta}(s,t))}&=&\cos{\l}-(\sin{\l})\,t-(\sin{\l})\,s\,t-\frac{1}{2}(\cos{\l})\,t^2-(\sin{\l})\,s^2\,t\nonumber\\
&-&(\cos{\l})\,s\,t^2+\frac{2}{3}\,(\sin{\l})\,t^3-(\sin{\l})\,s^3\,t-\frac{3}{2}\,(\cos{\l})\,s^2\,t^2\nonumber\\
&+&\frac{13}{6}\,(\sin{\l})\,s\,t^3+\frac{13}{24}\,(\cos{\l})\,t^4+O_5(s,t)\nonumber\\
\nonumber\\
\end{eqnarray}
On the other hand, taking into account the following expansions around $(\eta,\xi)=0$
\begin{eqnarray}\label{last term}
\hat{\eta}\,\sqrt{1-\frac{\hat{\eta}^2+\hat{\xi}^2}{4}}&=&\hat{\eta}-\frac{\hat{\eta}^3+\hat{\eta}\,\hat{\xi}^2}{8}+O_5(\hat{\eta},\hat{\xi})\nonumber\\
\hat{\xi}\,\sqrt{1-\frac{\hat{\eta}^2+\hat{\xi}^2}{4}}&=&\hat{\xi}-\frac{\hat{\eta}^2\,\hat{\xi}+\hat{\xi}^3}{8}+O_5(\hat{\eta},\hat{\xi})\nonumber\\
\end{eqnarray}
we find, by muliplying the first by $\cos{\l}$ ($\sin{\l}$) and the second by $-\sin{\l}$ ($\cos{\l}$), and, then, taking the sum
\begin{eqnarray}\label{st1}
s(\l,\hat{\eta},{\hat{\xi}})&=&\sqrt{1-\frac{\hat{\eta}^2+{\hat{\xi}}^2}{4}}\,(\hat{\eta}\cos{\l}-{\hat{\xi}}\sin{\l})\nonumber\\
&=&{(\cos{\l})\,\hat{\eta}-(\sin{\l})\,\hat{\xi}}-\frac{(\cos{\l})\,\hat{\eta}^3-(\sin{\l})\,\hat{\eta}^2\,\hat{\xi}+(\cos{\l})\,\hat{\eta}\,\hat{\xi}^2-(\sin{\l})\,\hat{\xi}^3}{8}\nonumber\\
&+&O_5(\hat{\eta},\hat{\xi})\nonumber\\
\nonumber\\
\left(t(\L,\hat{\eta},\l,\hat{\xi})\right.&=&\sqrt{1-\frac{\hat{\eta}^2+{\hat{\xi}}^2}{4}}\,(\hat{\eta}\sin{\l}+{\hat{\xi}}\cos{\l})\nonumber\\
&=&{(\sin{\l})\,\hat{\eta}+(\cos{\l})\,\hat{\xi}}-\frac{(\sin{\l})\,\hat{\eta}^3+(\cos{\l})\,\hat{\eta}^2\,\hat{\xi}+(\sin{\l})\,\hat{\eta}\,\hat{\xi}^2+(\cos{\l})\,\hat{\xi}^3}{8}\nonumber\\
&+&O_5(\hat{\eta},\hat{\xi})\ ,
\end{eqnarray}
respectively). Consequentely, substituting the (\ref{st1}), into (\ref{CosAndSin}), we obtain
\begin{eqnarray*}
\cos{(\l+\hat{\zeta}(\l,\hat{\xi},\hat{\eta}))}&=&\cos{(\l+\tilde{\zeta}(\l,s(\l,\hat{\xi},\hat{\eta}),t(\l,\hat{\xi},\hat{\eta})))}\nonumber\\
&=&\cos{\l}+\frac{\cos{2\l}-1}{2}\,\hat{\eta}-\frac{\sin{2\l}}{2}\,\hat{\xi}\nonumber\\
&-&3\,\frac{\cos{\l}-\cos{3\l}}{8}\,\hat{\eta}^2+\frac{\sin{\l}-3\,\sin{3\l}}{4}\,\hat{\eta}\,\hat{\xi}-\frac{\cos{\l}+3\,\cos{3\l}}{8}\,\hat{\xi}^2\nonumber\\
&+&\frac{3-19\,\cos{2\l}+16\,\cos{4\l}}{48}\,\hat{\eta}^3+\frac{9\sin{2\l}-16\,\sin{4\l}}{16}\,\hat{\eta}^2\,\hat{\xi}
\nonumber\\
&+&\frac{1-\cos{2\l}-16\,\cos{4\l}}{16}\,\hat{\eta}\,\hat{\xi}^2+\frac{11\,\sin{2\l}+16\,\sin{4\l}}{48}\,\hat{\xi}^3\nonumber\\
&+&\frac{46\,\cos{\l}-171\,\cos{3\l}+125\,\cos{5\l}}{384}\,\hat{\eta}^4\nonumber\\
&-&\frac{8\,\sin{\l}-99\,\sin{3\l}+125\sin{5\l}}{96}\,\hat{\eta}^3\,\hat{\xi}\nonumber\\
&+&\frac{10\,\cos{\l}+27\,\cos{3\l}-125\,\cos{5\l}}{64}\,\hat{\eta}^2\hat{\xi}^2\nonumber\\
&-&\frac{8\,\sin{\l}-45\,\sin{3\l}-125\sin{5\l}}{96}\,\hat{\eta}\,\hat{\xi}^3\nonumber\\
&+&\frac{14\,\cos{\l}+117\,\cos{3\l}+125\,\cos{5\l}}{384}\,\hat{\xi}^4+O_5(\hat{\eta},\hat{\xi})\nonumber\\
\end{eqnarray*}
The expansion of $\sin{(\l+\hat{\zeta}(\l,\hat{\eta},\hat{\xi}))}$ is obtained by the previous one, using
$$\sin{(\l+\hat{\zeta}(\l,\hat{\eta},\hat{\xi}))}=-\cos{(\l+\pi/2+\hat{\zeta}(\l+\pi/2,\hat{\xi},-\hat{\eta}))}$$
 and gives
\begin{eqnarray*}
\sin{(\l+\hat{\zeta}(\l,\hat{\eta},\hat{\xi}))}&=&\sin{\l}+\frac{\cos{2\l}+1}{2}\,\hat{\xi}+\frac{\sin{2\l}}{2}\,\hat{\eta}\nonumber\\
&-&3\,\frac{\sin{\l}+\sin{3\l}}{8}\,\hat{\xi}^2+\frac{\cos{\l}+3\,\cos{3\l}}{4}\,\hat{\eta}\,\hat{\xi}+\frac{-\sin{\l}+3\,\sin{3\l}}{8}\,\hat{\eta}^2\nonumber\\
&-&\frac{3+19\,\cos{2\l}+16\,\cos{4\l}}{48}\,\hat{\xi}^3-\frac{9\sin{2\l}+16\,\sin{4\l}}{16}\,\hat{\xi}^2\,\hat{\eta}
\nonumber\\
&-&\frac{1+\cos{2\l}-16\,\cos{4\l}}{16}\,\hat{\xi}\,\hat{\eta}^2+\frac{-11\,\sin{2\l}+16\,\sin{4\l}}{48}\,\hat{\eta}^3\nonumber\\
&+&\frac{46\,\sin{\l}+171\,\sin{3\l}+125\,\sin{5\l}}{384}\,\hat{\xi}^4\nonumber\\
&-&\frac{8\,\cos{\l}+99\,\cos{3\l}+125\cos{5\l}}{96}\,\hat{\xi}^3\,\hat{\eta}\nonumber\\
&-&\frac{-10\,\sin{\l}+27\,\sin{3\l}+125\,\sin{5\l}}{64}\,\hat{\xi}^2\hat{\eta}^2\nonumber\\
&-&\frac{8\,\cos{\l}+45\,\cos{3\l}-125\cos{5\l}}{96}\,\hat{\xi}\,\hat{\eta}^3\nonumber\\
&-&\frac{-14\,\sin{\l}+117\,\sin{3\l}-125\,\sin{5\l}}{384}\,\hat{\eta}^4+O_5(\hat{\eta},\hat{\xi})\ .\nonumber\\
\end{eqnarray*}
Both the previous expansions (of $\cos{(\l+\hat{\zeta}(\l,\hat{\eta},\hat{\xi})))}$ and $\sin{(\l+\hat{\zeta}(\l,\hat{\eta},\hat{\xi}))}$), together with the first line in (\ref{last term}), are, next, inserted into the expression of $\hat x_1$ given in (\ref{poincare}), and one obtains the expansion of $\hat x_1$ to the fourth order:
\begin{eqnarray}\label{q1expansion}
\hat q_1(a,\l,\hat{\eta},\hat{\xi})&=&a\,\left[\cos{\l}+\frac{-3+\cos{2\l}}{2}\,\hat{\eta}-\frac{\sin{2\l}}{2}\,\hat{\xi}\right.\nonumber\\
&-&3\,\frac{\cos{\l}-\cos{3\l}}{8}\,\hat{\eta}^2-\frac{\sin{\l}+3\,\sin{3\l}}{4}\,\hat{\eta}\,\hat{\xi}-\frac{5\,\cos{\l}+3\,\cos{3\l}}{8}\,\hat{\xi}^2\nonumber\\
&+&\frac{9-19\,\cos{2\l}+16\,\cos{4\l}}{48}\,\hat{\eta}^3+\frac{5\,\sin{2\l}-16\,\sin{4\l}}{16}\,\hat{\eta}^2\,\hat{\xi}\nonumber\\
&+&\frac{3-9\,\cos{2\l}-16\cos{4\l}}{16}\,\hat{\eta}\,\hat{\xi}^2+\frac{23\,\sin{2\l}+16\sin{4\l}}{48}\,\hat{\xi}^3\nonumber\\
&+&\frac{46\,\cos{\l}-171\,\cos{3\l}+125\,\cos{5\l}}{384}\,\hat{\eta}^4\nonumber\\
&+&\frac{-2\,\sin{\l}+81\,\sin{3\l}-125\,\sin{5\l}}{96}\,\hat{\eta}^3\,\hat{\xi}\nonumber\\
&+&\frac{14\,\cos{\l}-9\,\cos{3\l}-125\,\cos{5\l}}{64}\,\hat{\eta}^2\,\hat{\xi}^2\nonumber\\
&+&\frac{-2\,\sin{\l}+99\,\sin{3\l}+125\,\sin{5\l}}{96}\,\hat{\eta}\,\hat{\xi}^3\nonumber\\
&+&\left.\frac{38\,\cos{\l}+189\,\cos{3\l}+125\,\cos{5\l}}{384}\,\hat{\xi}^4\right]+O_5(\hat{\eta},\hat{\xi})
\end{eqnarray}
As outlined before, the expansion for $\hat x_2$ is obtained by
$$\hat x_2(a,\l,\hat{\eta},\hat{\xi})=-\hat x_1(a,\l+\pi/2,\hat{\xi},-\hat{\eta})\ ,$$
and the result is
\begin{eqnarray}\label{q2expansion}
\hat x_2(a,\l,\hat{\eta},\hat{\xi})&=&a\,\left[\sin{\l}+\frac{3+\cos{2\l}}{2}\,\hat{\xi}+\frac{\sin{2\l}}{2}\,\hat{\eta}\right.\nonumber\\
&-&3\,\frac{\sin{\l}+\sin{3\l}}{8}\,\hat{\xi}^2-\frac{\cos{\l}-3\,\cos{3\l}}{4}\,\hat{\eta}\,\hat{\xi}-\frac{5\,\sin{\l}-3\,\sin{3\l}}{8}\,\hat{\eta}^2\nonumber\\
&-&\frac{9+19\,\cos{2\l}+16\,\cos{4\l}}{48}\,\hat{\xi}^3-\frac{5\,\sin{2\l}+16\,\sin{4\l}}{16}\,\hat{\xi}^2\,\hat{\eta}\nonumber\\
&-&\frac{3+9\,\cos{2\l}-16\cos{4\l}}{16}\,\hat{\xi}\,\hat{\eta}^2+\frac{-23\,\sin{2\l}+16\sin{4\l}}{48}\,\hat{\eta}^3\nonumber\\
&+&\frac{46\,\sin{\l}+171\,\sin{3\l}+125\,\sin{5\l}}{384}\,\hat{\xi}^4\nonumber\\
&+&\frac{-2\,\cos{\l}-81\,\cos{3\l}-125\,\cos{5\l}}{96}\,\hat{\xi}^3\,\hat{\eta}\nonumber\\
&+&\frac{14\,\sin{\l}+9\,\sin{3\l}-125\,\sin{5\l}}{64}\,\hat{\xi}^2\,\hat{\eta}^2\nonumber\\
&+&\frac{-2\,\cos{\l}-99\,\cos{3\l}+125\,\cos{5\l}}{96}\,\hat{\xi}\,\hat{\eta}^3\nonumber\\
&+&\left.\frac{38\,\sin{\l}-189\,\sin{3\l}+125\,\sin{5\l}}{384}\,\hat{\eta}^4\right]+O_5(\hat{\eta},\hat{\xi})
\end{eqnarray}
Expansions (\ref{q1expansion}) and (\ref{q2expansion}) for $\hat x_1$ and  $\hat x_2$ give at sight their desired derivatives in $(\eta,\xi)=0$, which are substituted into (\ref{d2expansion}), giving
\begin{eqnarray*}
\left.\hat{d}_2\right|_0&=&a_1^2+a_2^2-2\,a_1\,a_2\,\cos{(\l_1-\l_2)}\nonumber\\ \nonumber\\
\left.\partial_{\hat{\eta}_1}\hat{d}_2\right|_0&=&-{a_1}\,\left[2\,a_1\,\cos{\l_1}+a_2\,(\cos{(2\l_1-\l_2)}-3\,\cos{\l_2})\right]\nonumber\\ \nonumber\\
\left.\partial_{\hat{\eta}_1^2}\hat{d}_2\right|_0&=&-\frac{a_1}{2}\,\left[-6\,a_1+2\,a_1\,\cos{(2\l_1)}-4\,a_2\,\cos{(\l_1-\l_2)}+3\,a_2\cos{(3\l_1-\l_2)}\right.\nonumber\\
&+&\left.a_2\,\cos{(\l_1+\l_2)}\right]\nonumber\\ \nonumber\\
\left.\partial_{\hat{\eta}_1\hat{\eta}_2}\hat{d}_2\right|_0&=&-\frac{a_1\,a_2}{2}\,\left[9-3\,\cos{(2\l_1)}+\cos{(2\l_1-2\l_2)}-3\,\cos{(2\l_2)}\right]\nonumber\\ \nonumber\\
\left.\partial_{\hat{\eta}_1\hat{\xi}_1}\hat{d}_2\right|_0&=&\frac{a_1}{2}\,\left[2\,a_1\,\sin{(2\l_1)}+a_2\,(3\,\sin{(3\l_1-\l_2})+\sin{(\l_1+\l_2)}\right]\nonumber\\ \nonumber\\
\left.\partial_{\hat{\eta}_2\hat{\xi}_1}\hat{d}_2\right|_0&=&\frac{a_1\,a_2}{2}\,\left[-3\,\sin{(2\l_1)}+\sin{(2\l_1-2\l_2)}-3\,\sin{2\l_2}\right]\nonumber\\ \nonumber\\
\left.\partial_{\hat{\eta}_1^3}\hat{d}_2\right|_0&=&-\frac{a_1}{4}\,\left[-12\,a_1\,\cos{\l_1}+6\,a_1\cos{3\l_1}+a_2\,(-21\,\cos{(2\l_1-\l_2)}\right.\nonumber\\
&+&\left.16\,\cos{(4\l_1-\l_2)}+9\,\cos{\l_2}+2\,\cos{(2\l_1+\l_2)})\right]\nonumber\\ \nonumber\\
\left.\partial_{\hat{\eta}_1^2\hat{\eta}_2}\hat{d}_2\right|_0&=&\frac{a_1\,a_2}{4}\,\left[-9\,\cos{\l_1}+9\cos{(3\l_1)}+4\,\cos{(\l_1-2\l_2)}-3\,\cos{(3\l_1-2\l_2)}\right.\nonumber\\
&-&\left.\cos{(\l_1+2\l_2)}\right]\nonumber\\ \nonumber\\
\left.\partial_{\hat{\eta}_1\hat{\xi}_1^2}\hat{d}_2\right|_0&=&\frac{a_1}{4}\,\left[4\,a_1\cos{\l_1}+6\,a_1\,\cos{(3\l_1)}+7\,a_2(\cos{(2\l_1-\l_2)}\right.\nonumber\\
&+&\left.16\,\cos{(4\l_1-\l_2)}-3\,\cos{\l_2}+2\,\cos{(2\l_1+\l_2)})\right]\nonumber\\ \nonumber\\
\left.\partial_{\hat{\eta}_2\hat{\xi}_1^2}\hat{d}_2\right|_0&=&\frac{a_1\,a_2}{4}\,\left.[-15\,\cos{\l_1}-9\,\cos{(3\l_1)}+4\,\cos{(\l_1-2\l_2)}\right.\nonumber\\
&+&\left.3\,\cos{(3\l_1-2\l_2)}\right]+\left.\cos{(\l_1+2\l_2)}]\right]\nonumber\\ \nonumber\\
\left.\partial_{\hat{\eta}_1\hat{\eta}_2\hat{\xi}_1}\hat{d}_2\right|_0&=&\frac{a_1\,a_2}{4}\,\left[-3\,\sin{\l_1}-9\,\sin{(3\l_1)}+3\,\sin{(3\l_1-2\l_2)}+\sin{(\l_1+2\l_2)}\right]\nonumber\\ \nonumber\\
\left.\partial_{\hat{\eta}_1^4}\hat{d}_2\right|_0&=&\frac{a_1}{8}\,\left[-72\,a_1+56\,a_1\,\cos{(2\l_1)}-32\,a_1\,\cos{(4\l_1)}-42\,a_2\,\cos{(\l_1-\l_2)}\right.\nonumber\\
&+&\left.180\,a_2\,\cos{(3\l_1-\l_2)}-125\,a_2\,\cos{(5\l_1-\l_2)}-4\,a_2\,\cos{(\l_1+\l_2)}\right.\nonumber\\
&-&\left.9\,a_2\,\cos{(3\l_1+\l_2)}\right]\nonumber\\ \nonumber\\
\left.\partial_{\hat{\eta}_1^3\hat{\eta}_2}\hat{d}_2\right|_0&=&\frac{a_1\,a_2}{8}\,\left[27\,-57\,\cos{(2\l_1)}+48\,\cos{(4\l_1)}+21\,\cos{(2\l_1-2\l_2)}\right.\nonumber\\
&-&\left.16\,\cos{(4\l_1-2\l_2)}-9\,\cos{(2\l_2)}-2\,\cos{(2\l_1+2\l_2)}\right]\nonumber\\ \nonumber\\
\left.\partial_{\hat{\eta}_1^2\hat{\eta}_2^2}\hat{d}_2\right|_0&=&\frac{a_1\,a_2}{8}\,\left[12\,\cos{(\l_1-3\l_2)}-9\,\cos{(3\l_1-3\l_2)}-17\,\cos{(\l_1-\l_2)}\right.\nonumber\\
&+&\left.12\,\cos{(3\l_1-\l_2)}+8\,\cos{(\l_1+\l_2)}-3\,\cos{(3\l_1+\l_2)}\right.\nonumber\\
&-&\left.3\,\cos{(\l_1+3\l_2)}\right]\nonumber\\\nonumber\\
\left.\partial_{\hat{\eta}_1^2\hat{\xi}_1^2}\hat{d}_2\right|_0&=&\frac{a_1}{8}\,\left[-24\,a_1+32\,a_1\,\cos{4\l_1}-14\,a_2\,\cos{(\l_1-\l_2)}\right.\nonumber\\
&+&\left.125\,a_2\,\cos{(5\l_1-\l_2)}+9\,a_2\,\cos{(3\l_1+\l_2)}\right]\nonumber\\\nonumber\\
\left.\partial_{\hat{\eta}_2^2\hat{\xi}_1^2}\hat{d}_2\right|_0&=&\frac{3\,a_1\,a_2}{8}\,\left[4\,\cos{(\l_1-3\l_2)}+3\,\cos{(3\l_1-3\l_2)}-5\,\cos{(\l_1-\l_2)}\right.\nonumber\\
&-&4\left.\,\cos{(3\l_1-\l_2)}+\cos{(3\l_1+\l_2)}+\cos{(\l_1+3\l_2)}\right]\nonumber\\\nonumber\\
\left.\partial_{\hat{\eta}_1\hat{\eta}_2\hat{\xi}_1^2}\hat{d}_2\right|_0&=&\frac{a_1\,a_2}{8}\,\left[9-27\,\cos{(2\l_1)}-48\,\cos{(4\l_1)}\right.\nonumber\\
&+&7\,\cos{(2\l_1-2\l_2)}+16\,\cos{(4\l_1-2\l_2)}-3\,\cos{(2\l_2)}\nonumber\\
&+&\left.2\,\cos{(2\l_1+2\l_2)}\right]\nonumber\\ \nonumber\\
\left.\partial_{\hat{\eta}_1\hat{\eta}_2\hat{\xi}_1\hat{\xi}_2}\hat{d}_2\right|_0&=&\frac{a_1\,a_2}{8}\,\left[-9\,\cos{(3\l_1-3\l_2)}-\cos{(\l_1-\l_2)}\right.\nonumber\\
&+&\left.3\,\cos{(3\l_1+\l_2)}+3\,\cos{(\l_1+3\l_2)}\right]\nonumber\\ \nonumber\\
\end{eqnarray*}
The remaining derivatives involved in (\ref{d2expansion}) are found by symmetry (as in the proof of Lemma \ref{simmetries remark}):
\begin{eqnarray*}
\left.\partial_{\hat{\eta}_2}\hat{d}_2\right|_0&=&{a_2}\,[3\,a_1\,\cos{\l_1}-a_1\,\cos{(\l_1-2\l_2)}-2\,a_2\,\cos{\l_2}]\nonumber\\\nonumber\\
\left.\partial_{\hat{\xi}_1}\hat{d}_2\right|_0&=&{a_1}\,[-3\,a_2\,\sin{\l_2}+a_2\,\sin{(2\l_1-\l_2)}+2\,a_1\,\sin{\l_1}]\nonumber\\\nonumber\\
\left.\partial_{\hat{\xi}_2}\hat{d}_2\right|_0&=&{a_2}\,[-3\,a_1\,\sin{\l_1}-a_1\,\sin{(\l_1-2\l_2)}+2\,a_2\,\sin{\l_2}]\nonumber\\\nonumber\\
\left.\partial_{\hat{\eta}_2^2}\hat{d}_2\right|_0&=&\frac{a_2}{2}\,\left[6\,a_2-3\,a_1\,\cos{(\l_1-3\l_2)}+4\,a_1\,\cos{(\l_1-\l_2)}\right.\nonumber\\
&-&\left.2\,a_2\,\cos{(2\l_2)}-a_1\,\cos{(\l_1+\l_2)}\right]\nonumber\\\nonumber\\
\left.\partial_{\hat{\eta}_1\hat{\xi}_2}\hat{d}_2\right|_0&=&\frac{a_1\,a_2}{2}\,\left[-3\,\sin{(2\l_2)}-\sin{(2\l_1-2\l_2)}-3\,\sin{2\l_1}\right] \nonumber\\\nonumber\\
\left.\partial_{\hat{\eta}_2\hat{\xi}_2}\hat{d}_2\right|_0&=&\frac{a_2}{2}\,\left[2\,a_2\,\sin{(2\l_2)}-3\,a_1\,\sin{(\l_1-3\l_2})+a_1\,\sin{(\l_2+\l_1)}\right]\nonumber\\ \nonumber\\
\left.\partial_{\hat{\xi}_1^2}\hat{d}_2\right|_0&=&\frac{a_1}{2}\,\left[6\,a_1+3\,a_2\,\cos{(3\l_1-\l_2)}+4\,a_2\,\cos{(\l_1-\l_2)}\right.\nonumber\\
&+&\left.2\,a_1\,\cos{(2\l_1)}+a_2\,\cos{(\l_1+\l_2)}\right]\nonumber\\\nonumber\\
\left.\partial_{\hat{\xi}_1\hat{\xi}_2}\hat{d}_2\right|_0&=&\frac{a_1\,a_2}{2}\,\left[-9-3\,\cos{(2\l_1)}-\cos{(2\l_1-2\l_2)}-3\,\cos{(2\l_2)}\right]\nonumber\\
\left.\partial_{\hat{\eta}_1\hat{\eta}_2^2}\hat{d}_2\right|_0&=&\frac{a_1\,a_2}{4}\,\left[-9\,\cos{\l_2}+9\cos{(3\l_2)}+4\,\cos{(2\l_1-\l_2)}-3\,\cos{(2\l_1-3\l_2)}\right.\nonumber\\
&-&\left.\cos{(2\l_1+\l_2)}\right]\nonumber\\ \nonumber\\
\left.\partial_{\hat{\eta}_1^2\hat{\xi}_1}\hat{d}_2\right|_0&=&\frac{a_1}{4}\,\left[-4\,a_1\sin{\l_1}+6\,a_1\,\sin{(3\l_1)}-7\,a_2(\sin{(2\l_1-\l_2)}\right.\nonumber\\
&+&16\,\sin{(4\l_1-\l_2)}+\left.3\,\sin{\l_2}+2\,\sin{(2\l_1+\l_2)})\right]\nonumber\\ \nonumber\\
\left.\partial_{\hat{\eta}_2\hat{\xi}_1\hat{\xi}_2}\hat{d}_2\right|_0&=&\frac{a_1\,a_2}{4}\,\left[3\,\cos{\l_2}-9\,\cos{(3\l_2)}-3\,\cos{(2\l_1-3\l_2)}+\cos{(2\l_1+\l_2)}\right]\nonumber\\ \nonumber\\
\left.\partial_{\hat{\eta}_1\hat{\xi}_1\hat{\xi}_2}\hat{d}_2\right|_0&=&\frac{a_1\,a_2}{4}\,\left[3\,\cos{\l_1}-9\,\cos{(3\l_1)}-3\,\cos{(3\l_1-2\l_2)}+\cos{(\l_1+2\l_2)}\right]\nonumber\\ \nonumber\\
\left.\partial_{\hat{\eta}_2^2\hat{\xi}_1}\hat{d}_2\right|_0&=&\frac{a_1\,a_2}{4}\,\left.[15\,\sin{\l_2}-9\,\sin{(3\l_2)}-4\,\sin{(2\l_1-\l_2)}+3\,\sin{(2\l_1-3\l_2)}\right]\nonumber\\
&+&\left.\sin{(2\l_1+\l_2)}]\right]\nonumber\\ \nonumber\\
\left.\partial_{\hat{\eta}_1\hat{\eta}_2\hat{\xi}_1}\hat{d}_2\right|_0&=&\frac{a_1\,a_2}{4}\,\left[-3\,\sin{\l_1}-9\,\sin{(3\l_1)}+3\,\sin{(3\l_1-2\l_2)}+\sin{(\l_1+2\l_2)}\right]\nonumber\\ \nonumber\\
\left.\partial_{\hat{\eta}_1\hat{\eta}_2\hat{\xi}_2}\hat{d}_2\right|_0&=&\frac{a_1\,a_2}{4}\,\left[-3\,\sin{\l_2}-9\,\sin{(3\l_2)}-3\,\sin{(2\l_1-3\l_2)}+\sin{(2\l_1+\l_2)}\right]\nonumber\\ \nonumber\\
\end{eqnarray*} 
The last step is to insert the previous list of derivatives into (\ref{derivatives od inversdist}). The result is
\begin{eqnarray*}
\left.\partial_{\hat{\eta}_1^4}\,\frac{1}{\hat{d}}\right|_0&=&\left[128\,(a_1^2+a_2^2-2\,a_1a_2\cos{(\l_1-\l_2)})^{9/2}\right]^{-1}\nonumber\\
&\times&\,\left\{a_1\,[-8\,a_1\,(224\,a_1^6-2964\,a_1^4\,a_2^2+1305\,a_1^2\,a_2^4-820\,a_2^6)\,\cos{2\l_1}\right.\nonumber\\
&+&(4096\,a_1^7-25080\,a_1^5\,a_2^2\,-20178\,a_1^3\,a_2^4-4952\,a_1\,a_2^6)\,\cos{(4\l_1)}\nonumber\\
&+&a_2\,(-360\,a_1^5\,a_2+25866\,a_1^3\,a_2^3-1800\,a_1\,a_2^5)-36\,a_1^3\,a_2^4\,\cos{(6\l_1-4\l_2)}\nonumber\\
&-&343\,a_1^3\,a_2^4\,\cos{(8\l_1-4\l_2)}+(11992\,a_1^4\,a_2^3-4488\,a_1^2\,a_2^5)\,\cos{(\l_1-3\l_2)}\nonumber\\
&+&(232\,a_1^4\,a_2^3+72\,a_1^2\,a_2^5)\,\cos{(5\l_1-3\l_2)}+(1676\,a_1^4\,a_2^3\nonumber\\
&+&1356\,a_1^2\,a_2^5)\,\cos{(7\l_1-3\l_2)}-5300\,a_1^3\,a_2^4\,\cos{(2\l_1-4\l_2)}+\nonumber\\
&+&(-4848\,a_1^5\,a_2^2-1280\,a_1^3\,a_2^4+144\,a_1\,a_2^6)\,\cos{(4\l_1-2\l_2)}\nonumber\\
&+&(-3492\,a_1^5\,a_2^2-6660\,a_1^3\,a_2^4-1908\,a_1\,a_2^6)\,\cos{(6\l_1-2\l_2)}\nonumber\\
&+&(336\,a_1^6\,a_2-12096\,a_1^4\,a_2^3-10656\,a_1^2\,a_2^5+336\,a_2^7)\,\cos{(\l_1-\l_2)}\nonumber\\
&+&(1800\,a_1^5\,a_2^2-12072\,a_1^3\,a_2^4+3240\,a_1\,a_2^6)\,\cos{(2\l_1-2\l_2)}\nonumber\\
&+&(3312\,a_1^4\,a_2^3+1872\,a_1^2\,a_2^5)\,\cos{(3\l_1-3\l_2)}+222\,a_1^3\,a_2^4\,\cos{(4\l_1-4\l_2)}\nonumber\\
&+&(-4896\,a_1^6\,a_2+14760\,a_1^4\,a_2^3-504\,a_1^2\,a_2^5-1440\,a_2^7)\,\cos{(3\l_1-\l_2)}\nonumber\\
&+&(6568\,a_1^6\,a_2+14148\,a_1^4\,a_2^3+9444\,a_1^2\,a_2^5\,+1000\,a_2^7)\,\cos{(5\l_1-\l_2)}\nonumber\\
&+&(-14256\,a_1^5\,a_2^2+33696\,a_1^3\,a_2^4-1584\,a_1\,a_2^6)\,\cos{(2\l_2)}\nonumber\\
&+&6561\,a_1^3\,a_2^4\,\cos{(4\l_2)}\nonumber\\
&+&(8096\,a_1^6\,a_2-42984\,a_1^4\,a_2^3\nonumber\\
&-&5448\,a_1^2\,a_2^5+32\,a_2^7)\,\cos{(\l_1+\l_2)}\nonumber\\
&+&(29436\,a_1^5\,a_2^2-18484\,a_1^3\,a_2^4-468\,a_1\,a_2^6)\,\cos{(2\l_1+2\l_2)}\nonumber\\
&+&(-17784\,a_1^6\,a_2+34236\,a_1^4\,a_2^3+10620\,a_1^2\,a_2^5+72\,a_2^7)\,\cos{(3\l_1+\l_2)}\nonumber\\
&+&\left.\left.(-22204\,a_1^4\,a_2^3+2340\,a_1^2\,a_2^5)\,\cos{(\l_1+3\l_2)}\right]\right\}\nonumber\\\nonumber\\\nonumber\\
\left.\partial_{\hat{\eta}_1^3\hat{\eta}_2}\,\frac{1}{\hat{d}}\right|_0&=&\left[128\,(a_1^2+a_2^2-2\,a_1a_2\cos{(\l_1-\l_2)})^{9/2}\right]^{-1}\nonumber\\
&\times&\,\left\{a_1\,a_2\,\left[-3438\,a_1^4\,a_2^2-3798\,a_1^2\,a_2^4\,\right.\right.\nonumber\\
&+&24\,a_1^2\,(76\,a_1^4-197\,a_1^2\,a_2^2+113\,a_2^4)\,\cos{(2\l_1)}\nonumber\\
&+&(-6144\,a_1^6+9462\,a_1^4\,a_2^2+2904\,a_1^2\,a_2^4)\,\cos{(4\l_1)}\nonumber\\
&-&343\,a_1^3\,a_2^3\,\cos{(\l_1-5\l_2)}\nonumber\\
&+&207\,a_1^3\,a_2^3\,\cos{(3\l_1-5\l_2)}+25\,a_1^3\,a_2^3\,\cos{(7\l_1-5\l_2)}\nonumber\\
&-&(162\,a_1^4\,a_2^2+18\,a_1^2\,a_2^4)\,\cos{(6\l_1-4\l_2)}\nonumber\\
&+&(72\,a_1^5\,a_2-18195\,a_1^3\,a_2^3+2736\,a_1\,a_2^5)\,\cos{(\l_1-3\l_2)}\nonumber\\
&+&(924\,a_1^5\,a_2+243\,a_1^3\,a_2^3-204\,a_1\,a_2^5)\,\cos{(5\l_1-3\l_2)}\nonumber\\
&-&75\,a_1^3\,a_2^3\cos{(7\l_1-3\l_2)}\nonumber\\
&+&(666\,a_1^4\,a_2^2+5274\,a_1^2\,a_2^4)\,\cos{(2\l_1-4\l_2)}\nonumber\\
&+&(2048\,a_1^6-228\,a_1^4\,a_2^2+1260\,a_1^2\,a_2^4+512\,a_2^6)\,\cos{(4\l_1-2\l_2)}\nonumber\\
&+&(486\,a_1^4\,a_2^2+270\,a_1^2\,a_2^4)\,\cos{(6\l_1-2\l_2)}\nonumber\\
&+&(2232\,a_1^5\,a_2-6042\,a_1^3\,a_2^3+2592\,a_1\,a_2^5)\,\cos{(\l_1-\l_2)}\nonumber\\
&+&(-672\,a_1^6+7704\,a_1^4\,a_2^2+7704\,a_1^2\,a_2^4-672\,a_2^6)\,\cos{(2\l_1-2\l_2)}\nonumber\\
&-&(1656\,a_1^5\,a_2+1269\,a_1^3\,a_2^3+2016\,a_2\,a_2^5)\,\cos{(3\l_1-3\l_2)}\nonumber\\
&-&(522\,a_1^4\,a_2^2+162\,a_1^2\,a_2^4)\,\cos{(4\l_1-4\l_2)}+15\,a_1^3\,a_2^3\,\cos{(5\l_1-5\l_2)}\nonumber\\
&-&(4656\,a_1^5\,a_2+4518\,a_1^3\,a_2^3+1656\,a_1\,a_2^5)\,\cos{(3\l_1-\l_2)}\nonumber\\
&-&(2772\,a_1^5\,a_2+1864\,a_1^3\,a_2^3+300\,a_1\,a_2^5)\,\cos{(5\l_1-\l_2)}\nonumber\\
&+&(16452\,a_1^4\,a_2^2-12588\,a_1^2\,a_2^4+288\,a_2^6)\,\cos{(2\l_2)}\nonumber\\
&+&(450\,a_1^4\,a_2^2-7350\,a_1^2\,a_2^4)\,\cos{(4\l_2)}\nonumber\\
&+&(-8820\,a_1^5\,a_2+21598\,a_1^3\,a_2^3-1068\,a_1\,a_2^5)\,\cos{(\l_1+\l_2)}\nonumber\\
&+&(64\,a_1^6-41694\,a_1^4\,a_2^2+4410\,a_2^2\,a_2^4+64\,a_2^6)\,\cos{(2\l_1+2\l_2)}\nonumber\\
&+&(26088\,a_1^5\,a_2-10854\,a_1^3\,a_2^3-360\,a_1\,a_2^5)\,\cos{(3\l_1+\l_2)}\nonumber\\
&+&\left.\left.(-276\,a_1^5\,a_2+29520\,a_1^3\,a_2^3-876\,a_1\,a_2^5)\,\cos{(\l_1+3\l_2)}\right]\right\}\nonumber\\\nonumber\\\nonumber\\
\left.\partial_{\hat{\eta}_1^2\hat{\eta}_2^2}\,\frac{1}{\hat{d}}\right|_0&=&\left[128\,(a_1^2+a_2^2-2\,a_1a_2\cos{(\l_1-\l_2)})^{9/2}\right]^{-1}\nonumber\\
&\times&\,\left\{-a_1\,a_2\,\left[324\,a_1^5\,a_2-10584\,a_1^3\,a_2^3+324\,a_1\,a_2^5\right.\right.\nonumber\\
&-&12\,a_1\,a_2\,(220\,a_1^4-1067\,a_1^2\,a_2^2+52\,a_2^4)\cos{(2\l_1)}\nonumber\\
&-&12\,a_1\,a_2\,(599\,a_1^4-277\,a_1^2\,_2^2-9\,a_2^4)\,\cos{(4\l_1)}-30\,a_1^3\,a_2^3\,\cos{(4\l_1-6\l_2)}\nonumber\\
&+&(264\,a_1^4\,a_2^2-1176\,a_1^2\,a_2^4)\,\cos{(\l_1-5\l_2)}\nonumber\\
&+&(36\,a_1^4\,a_2^2+564\,a_1^2\,a_2^4)\,\cos{(3\l_1-5\l_2)}\nonumber\\
&-&30\,a_1^3\,a_2^3\,\cos{(6\l_1-4\l_2)}\nonumber\\
&+&(96\,a_1^6-1920\,a_1^4\,a_2^2-14496\,a_1^2\,a_2^4+864\,a_2^6)\,\cos{(\l_1-3\l_2)}\nonumber\\
&+&75\,a_1^3\,a_2^3\,\cos{(2\l_1-6\l_2)}+(564\,a_1^4\,a_2^2+36\,a_1^2\,a_2^4)\,\cos{(5\l_1-3\l_2)}\nonumber\\
&+&(528\,a_1^5\,a_2+26\,a_1^3\,a_2^3+4272\,a_1\,a_2^5)\,\cos{(2\l_1-4\l_2)}\nonumber\\
&+&(4272\,a_1^5\,a_2+26\,a_1^3\,a_2^3+528\,a_1\,a_2^5)\,\cos{(4\l_1-2\l_2)}\nonumber\\
&+&75\,a_1^3\,a_2^3\,\cos{(6\l_1-2\l_2)}\nonumber\\
&+&(-136\,a_1^6+2400\,a_1^4\,a_2^2+2400\,a_1^2\,a_2^4-136\,a_2^6)\,\cos{(\l_1-\l_2)}\nonumber\\
&+&(1272\,a_1^5\,a_2+6337\,a_1^3\,a_2^3+1272\,a_1\,a_2^5)\,\cos{(2\l_1-2\l_2)}\nonumber\\
&+&(-648\,a_1^6-396\,a_1^4\,a_2^2-396\,a_1^2\,a_2^4-648\,a_2^6)\,\cos{(3\l_1-3\l_2)}\nonumber\\
&+&(-348\,a_1^5\,a_2-800\,a_1^3\,a_2^3-348\,a_1\,a_2^5)\,\cos{(4\l_1-4\l_2)}\nonumber\\
&+&(60\,a_1^4\,a_2^2+60\,a_1^2\,a_2^4)\,\cos{(5\l_1-5\l_2)}-9\,a_1^3\,a_2^3\,\cos{(6\l_1-6\l_2)}\nonumber\\
&+&(864\,a_1^6-14496\,a_1^4\,a_2^2-1920\,a_1^2\,a_2^4+96\,a_2^6)\,\cos{(3\l_1-\l_2)}\nonumber\\
&-&(1176\,a_1^4\,a_2^2+264\,a_1^2\,a_2^4)\,\cos{(5\l_1-\l_2)}\nonumber\\
&+&(-624\,a_1^5\,a_2+12804\,a_1^3\,a_2^3-2640\,a_1\,a_2^5)\,\cos{(2\l_2)}\nonumber\\
&+&(108\,a_1^5\,a_2+3324\,a_1^3\,a_2^3-7188\,a_1\,a_2^5)\,\cos{(4\l_2)}\nonumber\\
&+&(64\,a_1^6+456\,a_1^4\,a_2^2+456\,a_1^2\,a_2^4+64\,a_2^6)\,\cos{(\l_1+\l_2)}\nonumber\\
&+&(1224\,a_1^5\,a_2-45774\,a_1^3\,a_2^3+1224\,a_1\,a_2^5)\,\cos{(2\l_1+2\l_2)}\nonumber\\
&+&(-216\,a_1^6+29760\,a_1^4\,a_2^2-2736\,a_1^2\,a_2^4-24\,a_2^6)\,\cos{(3\l_1+\l_2)}\nonumber\\
&+&\left.\left.(-24\,a_1^6-2736\,a_1^4\,a_2^2+29760\,a_1^2\,a_2^4-216\,a_2^6)\,\cos{(\l_1+3\l_2)}\right]\right\}\nonumber\\\nonumber\\\nonumber\\
\left.\partial_{\hat{\eta}_1^2\hat{\xi}_1^2}\,\frac{1}{\hat{d}}\right|_0&=&\left[128\,(a_1^2+a_2^2-2\,a_1a_2\cos{(\l_1-\l_2)})^{9/2}\right]^{-1}\nonumber\\
&\times&\,\left\{-a_1\left[(4096\,a_1^7-25080\,a_1^5\,a_2^2-20178\,a_1^3\,a_2^4-4952\,a_1\,a_2^6)\,\cos{(4\l_1)}\right.\right.\nonumber\\
&+&120\,a_1^5\,a_2^2-8622\,a_1^3\,a_2^4+600\,a_1\,a_2^6-343\,a_1^3\,a_2^4\,\cos{(8\l_1-4\l_2)}\nonumber\\
&+&4\,a_1^2\,a_2^3\,(419\,a_1^2+339\,a_2^2)\,\cos{(7\l_1-3\l_2)}\nonumber\\
&+&(-3492\,a_1^5\,a_2^2-6660\,a_1^3\,a_2^4-1908\,a_1\,a_2^6)\,\cos{(6\l_1-2\l_2)}\nonumber\\
&+&(-112\,a_1^6\,a_2+4032\,a_1^4\,a_2^3+3552\,a_1^2\,a_2^5-112\,a_2^7)\,\cos{(\l_1-\l_2)}\nonumber\\
&+&(-600\,a_1^5\,a_2^2+4024\,a_1^3\,a_2^4-1080\,a_1\,a_2^6)\,\cos{(2\l_1-2\l_2)}\nonumber\\
&+&(-1104a_1^4\,a_2^3-624\,a_1^2\,a_2^5)\,\cos{(3\l_1-3\l_2)}-74\,a_1^3\,a_2^4\,{\cos{(4\l_1-4\l_2)}}\nonumber\\
&+&(6568\,a_1^6\,a_2+14148\,a_1^4\,a_2^3+9444\,a_1^2\,a_2^5+1000\,a_2^7))\,\cos{(5\l_1-\l_2)}\nonumber\\
&+&6561\,a_1^3\,a_2^4\,\cos{(4\l_2)}\nonumber\\
&+&(29436\,a_1^5\,a_2^2-18484\,a_1^3\,a_2^4-468\,a_1\,a_2^6)\,\cos{(2\l_1+2\l_2)}\nonumber\\
&+&(-17784\,a_1^6\,a_2+34236\,a_1^4\,a_2^3+10620\,a_1^2\,a_2^5+72\,a_2^7)\,\cos{(3\l_1+\l_2)}\nonumber\\
&+&\left.\left.(-22204\,a_1^4\,a_2^3+2340\,a_1^2\,a_2^5)\,\cos{(\l_1+3\l_2)}\right]\right\}\nonumber\\\nonumber\\\nonumber\\
\left.\partial_{\hat{\eta}_2^2\hat{\xi}_1^2}\,\frac{1}{\hat{d}}\right|_0&=&\left[128\,(a_1^2+a_2^2-2\,a_1a_2\cos{(\l_1-\l_2)})^{9/2}\right]^{-1}\nonumber\\
&\times&\,\left\{-3\,a_1\,a_2\,\left[84\,a_1^5\,a_2-8832\,a_1^3\,a_2^3+84\,a_1\,a_2^5\right.\right.\nonumber\\
&+&4\,a_1\,a_2\,(956\,a_1^4-3139\,a_1^2\,a_2^2+28\,a_2^4)\,\cos{(2\l_1)}\nonumber\\
&+&4\,a_1\,a_2\,(599\,a_1^4-277\,a_1^2\,a_2^2-9\,a_2^4)\,\cos{(4\l_1)}-10\,a_1^3\,a_2^3\,\cos{(4\l_1-6\l_2)}\nonumber\\
&+&(88\,a_1^4\,a_2^2+392\,a_1^2\,a_2^4)\,\cos{(\l_1-5\l_2)}\nonumber\\
&+&(12\,a_1^4\,a_2^2+188\,a_1^2\,a_2^4)\,\cos{(3\l_1-5\l_2)}\nonumber\\
&+&10\,a_1^3\,a_2^3\cos{(6\l_1-4\l_2)}\nonumber\\
&+&(32\,a_1^6-968\,a_1^4\,a_2^2-6920\,a_1^2\,a_2^4+288\,a_2^6)\,\cos{(\l_1-3\l_2)}\nonumber\\
&-&25\,a_1^3\,a_2^3\,\cos{(2\l_1-6\l_2)}-(188\,a_1^4\,a_2^2+12\,a_1^2\,a_2^4)\,\cos{(5\l_1-3\l_2)}\nonumber\\
&+&(176\,a_1^5\,a_2-158\,a_1^3\,a_2^3+1424\,a_1\,a_2^5)\,\cos{(2\l_1-4\l_2)}\nonumber\\
&+&(-1424\,a_1^5\,a_2+158\,a_1^3\,a_2^3-176\,a_1\,a_2^5)\,\cos{(4\l_1-2\l_2)}\nonumber\\
&-&25\,a_1^3\,a_2^3\,\cos{(6\l_1-2\l_2)}\nonumber\\
&+&(-40\,a_1^6+5216\,a_1^4\,a_2^2+5216\,a_1^2\,a_2^4-40\,a_2^6)\,\cos{(\l_1-\l_2)}\nonumber\\
&+&(-1640\,a_1^5\,a_2+2805\,a_1^3\,a_2^3-1640\,a_1\,a_2^5)\,\cos{(2\l_1-2\l_2)}\nonumber\\
&+&(216\,a_1^6-1020\,a_1^4\,a_2^2-1020\,a_1^2\,a_2^4+216\,a_2^6)\,\cos{(3\l_1-3\l_2)}\nonumber\\
&+&(116\,a_1^5\,a_2+200\,a_1^3\,a_2^3+116\,a_1\,a_2^5)\,\cos{(4\l_1-4\l_2)}\nonumber\\
&-&(20\,a_1^4\,a_2^2+20\,a_1^2\,a_2^4)\,\cos{(5\l_1-5\l_2)}+3\,a_1^3\,a_2^3\cos{(6\l_1-6\l_2)}\nonumber\\
&+&(-288\,a_1^6+6920\,a_1^4\,a_2^2+968\,a_1^2\,a_2^4-32\,a_2^6)\,\cos{(3\l_1-\l_2)}\nonumber\\
&+&(392\,a_1^4\,a_2^2+88\,a_1^2\,a_2^4)\,\cos{(5\l_1-\l_2)}\nonumber\\
&+&(-112\,a_1^5\,a_2+12556\,a_1^3\,a_2^3-2384\,a_1\,a_2^5)\,\cos{(2\l_2)}\nonumber\\
&+&(-36\,a_1^5\,a_2-1108\,a_1^3\,a_2^3+2396\,a_1\,a_2^5)\,\cos{(4\l_2)}\nonumber\\
&+&(-7056\,a_1^4\,a_2^2+7056\,a_1^2\,a_2^4)\,\cos{(\l_1+\l_2)}\nonumber\\
&+&(-408\,a_1^5\,a_2+15258\,a_1^3\,a_2^3-408\,a_1\,a_2^5)\,\cos{(2\l_1+2\l_2)}\nonumber\\
&+&(72\,a_1^6-9920\,a_1^4\,a_2^2+912\,a_1^2\,a_2^4+8\,a_2^6)\,\cos{(3\l_1+\l_2)}\nonumber\\
&+&\left.\left.(8a_1^6+912\,a_1^4\,a_2^2-9920\,a_1^2\,a_2^4+72\,a_2^6)\,\cos{(\l_1+3\l_2)}\right]\right\}\nonumber\\\nonumber\\\nonumber\\
\left.\partial_{\hat{\eta}_1\hat{\eta}_2\hat{\xi}_1^2}\,\frac{1}{\hat{d}}\right|_0&=&\left[128\,(a_1^2+a_2^2-2\,a_1a_2\cos{(\l_1-\l_2)})^{9/2}\right]^{-1}\nonumber\\
&\times&\,\,\left\{a_1\,a_2\,\left[-1146\,a_1^4\,a_2^2-1266\,a_1^2\,a_2^4\right.\right.\nonumber\\
&+&144\,a_1^2\,(6\,a_1^4-166\,a_1^2\,a_2^2-31\,a_2^4)\,\cos{(2\l_1)}\nonumber\\
&+&6\,a_1^2\,(1024\,a_1^4-1577\,a_1^2\,a-2^2-349\,a_2^4)\,\cos{(4\l_1)}\nonumber\\
&+&343\,a_1^3\,a_2^3\cos{(\l_1-5\l_2)}+69\,a_1^3\,a_2^3\,\cos{(3\l_1-5\l_2)}\nonumber\\
&-&25\,a_1^3a_2^3\,\cos{(7\l_1-5\l_2)}\nonumber\\
&+&(162\,a_1^4\,a_2^2+18\,a_1^2\,a_2^4)\,\cos{(6\l_1-4\l_2)}\nonumber\\
&+&(24\,a_1^5\,a_2-7437\,a_1^3\,a_2^3+912\,a_1\,a_2^5)\,\cos{(\l_1-3\l_2)}\nonumber\\
&+&(-924\,a_1^5\,a_2-303\,a_1^3\,a_2^3+204\,a_1\,a_2^5)\,\cos{(5\l_1-3\l_2)}\nonumber\\
&+&75\,a_1^3\,a_2^3\,\cos{(7\l_1-3\l_2)}\nonumber\\
&+&(222\,a_1^4\,a_2^2+1758\,a_1^2\,a_2^4)\,\cos{(2\l_1-4\l_2)}\nonumber\\
&+&(-2048\,a_1^6+2316\,a_1^4\,a_2^2-900\,a_1^2\,a_2^4-512\,a_2^6)\,\cos{(4\l_1-2\l_2)}\nonumber\\
&-&(486\,a_1^4\,a_2^2+270\,a_1^2\,a_2^4)\,\cos{(6\l_1-2\l_2)}\nonumber\\
&+&(744\,a_1^5\,a_2-2014\,a_1^3\,a_2^3+864\,a_1\,a_2^5)\,\cos{(\l_1-\l_2)}\nonumber\\
&+&(-224\,a_1^6+2568\,a_1^4\,a_2^2+2568\,a_1^2\,a_2^4-224\,a_2^6\,\cos{(2\l_1-2\l_2)}\nonumber\\
&+&(-552\,a_1^5\,a_2-423\,a_1^3\,a_2^3-672\,a_1\,a_2^5)\,\cos{(3\l_1-3\l_2)}\nonumber\\
&+&(-174\,a_1^4\,a_2^2-54\,a_1^2\,a_2^4)\,\cos{(4\l_1-4\l_2)}+5\,a_1^3\,a_2^3\,\cos{(5\l_1-5\l_2)}\nonumber\\
&+&(11280\,a_1^5\,a_2+294\,a_2^3\,a_2^3+2616\,a_1\,a_2^5)\,\cos{(3\l_1-\l_2)}\nonumber\\
&+&(2772\,a_1^5\,a_2+1864\,a_1^3\,a_2^3+300\,a_1\,a_2^5)\,\cos{(5\l_1-\l_2)}\nonumber\\
&+&(7284\,a_1^4\,a_2^2-6492\,a_1^2\,a_2^4+96\,a_2^6)\,\cos{(2\l_2)}\nonumber\\
&+&(-450\,a_1^4\,a_2^2+7350\,a_1^2\,a_2^4)\,\cos{(4\l_2)}\nonumber\\
&+&(-4044\,a_1^5\,a_2+23402\,a_1^3\,a_2^3-468\,a_1\,a_2^5)\,\cos{(\l_1+\l_2)}\nonumber\\
&+&(-64\,a_1^6+41694\,a_1^4\,a_2^2-4410\,a_1^2\,a_2^4-64\,a_2^6)\,\cos{(2\l_1+2\l_2)}\nonumber\\
&+&(-26088\,a_1^5\,a_2+10854\,a_1^3\,a_2^3+360\,a_1\,a_2^5)\,\cos{(3\l_1+\l_2)}\nonumber\\
&+&\left.\left.(276\,a_1^5\,a_2-29520\,a_1^3\,a_2^3+876\,a_1\,a_2^5)\,\cos{(\l_1+3\l_2)}\right]\right\}\nonumber\\\nonumber\\\nonumber\\
\left.\partial_{\hat{\eta}_1\hat{\eta}_2\hat{\xi}_1\hat{\xi}_2}\,\frac{1}{\hat{d}}\right|_0&=&\left[128\,(a_1^2+a_2^2-2\,a_1a_2\cos{(\l_1-\l_2)})^{9/2}\right]^{-1}\nonumber\\
&\times&\,\left\{a_1\,a_2\left[-36\,a_1^5\,a_2-7956\,a_1^3\,a_2^3-36\,a_1\,a_2^5\right.\right.\nonumber\\
&-&12\,a_1\,a_2\,(599\,a_1^4-277\,a_1^2\,a_2^2-9\,a_2^4)\,\cos{(4\l_1)}\nonumber\\
&-&24\,a_1^2\,a_2^2(11\,a_1^2+49\,a_2^2)\,\cos{(\l_1-5\l_2)}+75\,a_1^3\,a_2^3\,\cos{(2\l_1-6\l_2)}\nonumber\\
&+&75\,a_1^3\,a_2^3\,\cos{(6\l_1-2\l_2)}\nonumber\\
&+&(8\,a_1^6+6624\,a_1^4\,a_2^2+6624\,a_1^2\,a_2^4+8\,a_2^6)\,\cos{(\l_1-\l_2)}\nonumber\\
&+&(-3096\,a_1^5\,a_2+1039\,a_1^3\,a_2^3-3096\,a_1\,a_2^5)\,\cos{(2\l_1-2\l_2)}\nonumber\\
&+&(648\,a_1^6-1332\,a_1^4\,a_2^2-1332\,a_1^2\,a_2^4+648\,a_2^6)\,\cos{(3\l_1-3\l_2)}\nonumber\\
&+&(348\,a_1^5\,a_2+700\,a_1^3\,a_2^3+348\,a_1\,a_2^5)\,\cos{(4\l_1-4\l_2)}\nonumber\\
&-&(60\,a_1^4\,a_2^2+60\,a_1^2\,a_2^4)\,\cos{(5\l_1-5\l_2)}+9\,a_1^3\,a_2^3\,\cos{(6\l_1-6\l_2)}\nonumber\\
&-&(1176\,a_1^4\,a_2^2+264\,a_1^2\,a_2^4)\,\cos{(5\l_1-\l_2)}\nonumber\\
&+&(108\,a_1^5\,a_2+3324\,a_1^3\,a_2^3-7188\,a_1\,a_2^5)\,\cos{(4\l_2)}\nonumber\\
&+&(1224\,a_1^5\,a_2-45774\,a_1^3\,a_2^3+1224\,a_1\,a_2^5)\,\cos{(2\l_1+2\l_2)}\nonumber\\
&+&(-216\,a_1^6+29760\,a_1^4\,a_2^2-2736\,a_1^2\,a_2^4-24\,a_2^6)\,\cos{(3\l_1+\l_2)}\nonumber\\
&+&\left.\left.(-24\,a_1^6-2736\,a_1^4\,a_2^2+29760\,a_1^2\,a_2^4-216\,a_2^6)\,\cos{(\l_1+3\l_2)}\right]\right\}\ .
\end{eqnarray*}
The previous expressions have all the form of a finite sum
$$\sum_{k_1,k_2}{\cal P}_{k_1,k_2}(a_1,a_2)\,\frac{\cos{(k_1\,\l_1+k_2\l_2)}}{[a_1^2+a_2^2-2\,a_1\,a_2\,\cos{(\l_1-\l_2)}]^{9/2}}$$
with ${\cal P}_{k_1,k_2}(a_1,a_2)$ a homogeneous polynomial in $(a_1,a_2)$ with degree $8$. When taking the mean with respect to $(\l_1,\l_2)$, only the terms with $k_1+k_2=0$ survive. Taking also into account the normalizations factors in front of the r.h.s of (\ref{needed derivatives}), one finds the expressions (\ref{a4000 and others}), in terms of the Laplace coefficients $b_{9/2,k}(a_1/a_2)$.

\subsubsection{Proof of non Resonance}
In this paragraph, we prove

\begin{lemma}\label{order 2 asymptotics}
There exists $\d_*$ such that, for any $0<\d<\d_*$, the matrix $\cF(\L)$ defining the second order $f_2$ of the secular perturbation $\bar f$ is negative definite in $${\cal L}_{\underline\a,\ovl\a}(\m,\d)=\Big\{\L=(\L_1,\cdots,\L_N)\in{\real_+^N}:\ a_i:=\frac{1}{\hat m_i}\left(\frac{\L_i}{\tilde m_i}\right)^2\in [\underline\a,\ovl\a]\d^{N-i}\quad 1\leq i\leq N\Big\}\ .$$
Denoting, for such values of $\d$, by  $\O=(\O_1,\cdots,\O_N)$, $U(\L)=(u_{ij}(\L)$ the eigenvalues of $\cF(\L)$ and the unitary matrix through which $\cF(\L)$ is put in diagonal form
$$U(\L)^T\cF(\L)U(\L)=\textrm{\rm diag}[\O_1,\cdots,\O_N]\ ,\qquad U(\L)^TU(\L)=\id_N\ ,$$
then, $\O(\L)$ and $U(\L)$ satisfy the following asymptotics in $\d$:
\beqano
\O_i(\L)&=&\d^{(9-3N)/2}\times\arr{
-\frac{3}{4}\frac{\hat a_1^2}{\hat a_2^3}\frac{\bar m_1\bar m_2}{\tilde m_1\sqrt{\hat m_1\hat a_1}}+O(\d^2)\qquad \textrm{for}\quad i=1\\
\\
-\frac{3}{4}\frac{\hat a_{i-1}^2}{\hat a_i^3}\frac{\bar m_i\bar m_{i-1}}{\tilde m_i\sqrt{\hat m_i\hat a_i}}\d^{(3i-5)/2}+O(\d^{(3i-3)/2)})\quad \textrm{for}\quad 2\leq i\leq N-1\\
\\
-\frac{3}{4}\frac{\hat a_{N-1}^2}{\hat a_N^3}\frac{\bar m_N\bar m_{N-1}}{\tilde m_N\sqrt{\hat m_N\hat a_N}}\d^{(3N-5)/2}+O(\d^{(3N-2)/2})\quad \textrm{for}\quad i=N
}\nonumber\\
\nonumber\\
\nonumber\\
u_{ij}(\L)&=&\arr{
1
-\frac{25}{32}\left[\frac{\hat a_1\hat a_2^3}{\hat a_2^4}\sqrt{\frac{\tilde m_1}{\tilde m_2}}\sqrt[4]{\frac{\hat m_1\hat a_1}{\hat m_2\hat a_2}}\right]^2\d^{5/2}+o(\d^{5/2})\quad \\ \textrm{for}\quad 1\leq i=j\leq 2
\\
\\
1-\frac{25}{32}\left[\left(\frac{\bar m_i}{\bar m_{i-2}}\frac{a_{i-1}^6}{\hat a_i^4\hat a_{i-2}^2}\sqrt{\frac{\tilde m_{i-1}}{\tilde m_i}}\sqrt[4]{\frac{\hat m_{i-1}\hat a_{i-1}}{\hat m_i\hat a_i}}\right)^2\right.\\
+\left.\left(\frac{\bar m_{i+1}}{\bar m_{i-1}}\frac{a_i^6}{\hat a_{i-1}^4\hat a_{i-1}^2}\sqrt{\frac{\tilde m_i}{\tilde m_{i-1}}}\sqrt[4]{\frac{\hat m_{i}\hat a_i}{\hat m_{i-1}\hat a_{i-1}}}\right)^2 \right]\d^{9/2}+o(\d^{9/2})\quad \\ \textrm{for}\quad 2<i=j< N\\
\\
1-\frac{25}{32}\left(\frac{\bar m_N}{\bar m_{N-2}}\frac{a_N^6}{\hat a_N^4\hat a_{N-2}^2}\sqrt{\frac{\tilde m_{N-1}}{\tilde m_N}}\sqrt[4]{\frac{\hat m_{N-1}\hat a_{N-1}}{\hat m_N\hat a_N}}\right)^2\d^{9/2}+o(\d^{9/2})\quad \\ \textrm{for}\quad i=j=N\\
\\
\frac{5}{4}\frac{\bar m_j}{\bar m_2}\frac{\hat a_1\hat a_2^3}{\hat a_j^4}\sqrt{\frac{\tilde m_1}{\tilde m_2}}\sqrt[4]{\frac{\hat m_1\hat a_1}{\hat m_j\hat a_j}}\d^{(17j-29)/4}+o(\d^{(17j-29)/4})\quad \\ \textrm{for}\quad 1=i<j\leq N\\
\\
\frac{5}{4}\frac{\bar m_j}{\bar m_{i-1}}\frac{a_i^6}{\hat a_j^4\hat a_{i-1}^2}\sqrt{\frac{\tilde m_i}{\tilde m_j}}\sqrt[4]{\frac{\hat m_i\hat a_i}{\hat m_j\hat a_j}}\d^{(17j-17i-8)/4}+o(\d^{(17j-17i-8)/4})\quad \\ \textrm{for}\quad 1<i<j\leq N\\
\\
-\frac{5}{4}\frac{\bar m_i}{\bar m_2}\frac{\hat a_1\hat a_2^3}{\hat a_i^4}\sqrt{\frac{\tilde m_1}{\tilde m_i}}\sqrt[4]{\frac{\hat m_1\hat a_1}{\hat m_i\hat a_i}}\d^{(17i-29)/4}+o(\d^{(17i-29)/4})\quad \\ \textrm{for}\quad 1=j<i\leq N\\
\\
-\frac{5}{4}\frac{\bar m_i}{\bar m_{j-1}}\frac{a_j^6}{\hat a_i^4\hat a_{j-1}^2}\sqrt{\frac{\tilde m_j}{\tilde m_i}}\sqrt[4]{\frac{\hat m_j\hat a_j}{\hat m_i\hat a_i}}\d^{(17i-17j-8)/4}+o(\d^{(17i-17j-8)/4})\quad \\ \textrm{for}\quad 1<j<i\leq N\\
}
\eeqano
\end{lemma}
\begin{corollary}\label{non res cor}
There exists $\d^*>0$ such that, for any $0<\d\leq \d_*$, $\O(\L)$ is $4$--non resonant in ${\cal L}_{\underline\a,\ovl\a}(\m,\d)$:
$$\min_{1\leq |k|\leq 4}\min_{{\cal L}_{\underline\a,\ovl\a}(\m,\d)}|\O(\L)\cdot k|>0\ .$$
\end{corollary}
\vskip.1in
\noi
We will obtain Lemma \ref{order 2 asymptotics} as a consequence of Lemmas \ref{diagonalization N bodies} and \ref{asymp for F} below: apply Lemma \ref{diagonalization N bodies} to the matrix $\cF(\L)\d^{-(9-3N)/2}$, hence, with
\beqano
n_{ij}&=&\arr{
0\quad \textrm{for}\quad i=j=1\\
\frac{3i-5}{2}\quad \textrm{for}\quad 2\leq i=j\leq N\\
\frac{17j-11i-18}{4}\quad \textrm{for}\quad i<j\\
}\nonumber\\
\nonumber\\
a_{ij}(\d)&=&\left\{
\begin{array}{lrr}
-\frac{3}{4}\,\bar m_1\,\bar m_2\,\frac{(\hat a_1/\hat a_2)^2}{\hat a_2\,\tilde m_1\,\sqrt{\hat m_1\,\hat a_1}}+O(\d^{2})\ , \quad \textrm{for}\quad i=j=1\ ;\\
\\
-\frac{3}{4}\,\bar m_i\,\bar m_{i-1}\,\frac{(\hat a_{i-1}/\hat a_i)^2}{\hat a_i\,\tilde m_i\,\sqrt{\hat m_i\,\hat a_i}}\,\d^{(3i-5)/2}+O(\d^{(3i-3)/2})\\
\textrm{for}\quad 2\leq i=j\leq N-1\ ;\\
\\
-\frac{3}{4}\,\bar m_N\,\bar m_{N-1}\,\frac{(\hat a_{N-1}/\hat a_N)^2}{\hat a_N\,\tilde m_N\,\sqrt{\hat m_N\,\hat a_N}}\,\d^{(3N-5)/2}+O(\d^{(3N-1)/2})\\
\textrm{for}\quad i=j=N\ ,\\
\\
\frac{15}{16}\,\bar m_i\,\bar m_j\frac{(\hat a_i/\hat a_j)^3}{\hat a_j\,\sqrt{\tilde m_i\,\tilde m_j}\,\sqrt[4]{\hat m_i\,\hat m_j\,\hat a_i\,\hat a_j}}\,\d^{(17j-11i-18)/4}+O(\d^{(25j-19i-18)/4})\\
\textrm{for}\qquad i<j\\
\end{array}
\right.\nonumber\\
\eeqano
and next on comparing the remainder terms.
\begin{lemma}\label{diagonalization N bodies}
Let ${\cal A}=(a_{i,j})_{1\leq i,j\leq N}$ a real symmetric matrix with order $n$ and elements
\beq{aijneq0}a_{i,j}=\bar a_{i,j}(\d)\,\d^{n_{i,j}}\ ,\qquad \textrm{with}\qquad \bar a_{i,j}(0)\neq 0\eeq
and
\begin{eqnarray}\label{tilde F}
\arr{
n_{1,1}=0\ ;\quad n_{i,i}< n_{j,j}\quad \textrm{for}\ i<j\ \\
n_{i,j}<n_{i-1,j}\ ,\quad n_{i,j}<n_{i,j+1}\ \textrm{for}\ i<j+1\ .
}
\end{eqnarray}
Then, there exists $\bar\d$ such that, for any $0<\d<\bar\d$,  the eigenvalues  $\l_1$, $\cdots$, $\l_n$ of $\cA$ satisfy
$$|\l_i-a_{i,i}|\leq C\d^{m_i}\ ,\textrm{where}\qquad m_i=2\min\{n_{i-1,i},n_{i,i-1}\}-n_{i,i}$$
Furthermore,  the orthogonal matrix $V=\{v_{i,j}\}_{i,j=1,\cdots,N}$ which diagonalizes ${\cal A}$
\begin{eqnarray*}
 V^T\, V=\textrm{\rm id}_n\ ,\quad  V^T\, {\cal A}\, V=\textrm{\rm diag}(\l_1,\cdots,\l_n)
\end{eqnarray*}
satisfies
$$v_{i,j}=\d_{i,j}+\check v_{i,j}(\d)\,\d^{\n_{i,j}}$$ where
\begin{eqnarray}\label{check v and nu}
\n_{i,j}&=&\left\{
\begin{array}{lrr}
{n_{i,j}-n_{i,i}}\quad \textrm{for}\quad i<j\\
\\
2\, n_{12}\quad \textrm{for}\quad i=j=1\\
\\
2\,\min\{n_{i-1,i}-n_{i-1,i-1},n_{i+1,i}-n_{i,i}\}\quad\textrm{for}\quad 2\leq i=j\leq n-1 \\
\\
2(n_{n-1,n}-n_{n-1,n-1})\quad\textrm{for}\quad 2\leq i=j=n
\\
\\
{n_{i,j}-n_{j,j}}\quad \textrm{for}\quad i>j\ ,\\
\end{array}
\right.\nonumber\\
\nonumber\\
\check v_{i,j}(0)&=&\left\{
\begin{array}{lrr}
-\frac{\bar a_{i,j}(0)}{\bar a_{i,i}(0)}\quad \textrm{for}\quad i<j\ ,
\\
\\
\frac{\bar a_{i,j}(0)}{\bar a_{j,j}(0)}\quad \textrm{for}\quad i>j\ .\\
\\
-\frac{1}{2}\check v_{i+1,i}(0)^2 \quad\textrm{for}\quad  i=j=1\\
\\
-\frac{1}{2}\arr{
\check v_{i-1,i}(0)^2\qquad \textrm{if}\qquad n_{i-1,i}-n_{i-1,i-1}<n_{i+1,i}-n_{i,i}\\
\check v_{i+1,i}(0)^2+\check v_{i-1,i}(0)^2\qquad \textrm{if}\qquad n_{i-1,i}-n_{i-1,i-1}=n_{i+1,i}-n_{i,i}\\
\check v_{i+1,i}(0)^2\qquad \textrm{if}\qquad n_{i-1,i}-n_{i-1,i-1}>n_{i+1,i}-n_{i,i}\\
}\\
\textrm{for}\quad 2\leq i=j\leq n-1\\
\\
-\frac{1}{2}\check v_{n-1,n}(0)^2 \quad\textrm{for}\quad  i=j=n\\
\end{array}
\right.
\end{eqnarray}
if $\d_{i,j}$ is the Kronecker symbol.
\end{lemma}
This Lemma is purely technical and thus is proved in Appendix \ref{Proof of Lemma on diagonalization}.
\begin{lemma}\label{asymp for F}
The matrix $\cF(\L)$ satisfies the following asymptotics:
\begin{eqnarray*}
f_{ij}=\left\{
\begin{array}{lrr}
-\frac{3}{4}\,\bar m_1\,\bar m_2\,\frac{(\hat a_1/\hat a_2)^2}{\hat a_2\,\tilde m_1\,\sqrt{\hat m_1\,\hat a_1}}\,\d^{3(3-N)/2}+O(\d^{(13-3N)/2})\ , \  \textrm{for}\  i=j=1\ ;\\
\\
-\frac{3}{4}\,\bar m_i\,\bar m_{i-1}\,\frac{(\hat a_{i-1}/\hat a_i)^2}{\hat a_i\,\tilde m_i\,\sqrt{\hat m_i\,\hat a_i}}\,\d^{(3i+4-3N)/2}+O(\d^{(3i+6-3N)/2})\\ \textrm{for}\  2\leq i=j\leq N-1\ ;\\
\\
-\frac{3}{4}\,\bar m_N\,\bar m_{N-1}\,\frac{(\hat a_{N-1}/\hat a_N)^2}{\hat a_N\,\tilde m_N\,\sqrt{\hat m_N\,\hat a_N}}\,\d^{2}+O(\d^4)\ ,\  \textrm{for}\  i=j=N\ ,\\
\\
\frac{15}{16}\,\bar m_i\,\bar m_j\frac{(\hat a_i/\hat a_j)^3}{\hat a_j\,\sqrt{\tilde m_i\,\tilde m_j}\,\sqrt[4]{\hat m_i\,\hat m_j\,\hat a_i\,\hat a_j}}\,\d^{(17j-11i-6N)/4}+O(\d^{(25j-19i-6N)/4})\\ \textrm{for}\ i<j\\
\end{array}
\right.\nonumber\\
\end{eqnarray*}
\end{lemma}
{\bf Proof.}\ We start with computing the asymptotics for the diagonal elements, which, we write as 
\begin{eqnarray*}
{f}_{ii}=-2\,\frac{\bar m_i}{\L_i}\,\left[\sum_{k< i}\bar m_k\,a_{2000}(a_k,a_i)+\sum_{k> i}\bar m_k\,a_{2000}(a_i,a_k)\right]\qquad 1\leq i\leq N
\end{eqnarray*}
with
\beqa{order2newapprox}
a_{2000}({a},{b})&=&\frac{{a}}{8{b}^2}\,\left[-7{a}/{b}\,b_{5/2,0}({a}/{b})\right.\nonumber\\
&&+\left.4(1+{a}^2/{b}^2)\,b_{5/2,1}({a}/{b})-{a}/{b}\,b_{5/2,2}({a}/{b})\right]\nonumber\\
&=&\frac{3\,a^2}{8\,{b}^3}+O\left(\frac{a^4}{b^5}\right)
\eeqa
having used the following asymptotics for the involved Laplace coefficients (see Appendix \ref{Laplace coefficients}, Lemma \ref{asymp}):
\begin{eqnarray*}
b_{5/2,0}(\a)&=&1+O(\a^2)\nonumber\\
b_{5/2,1}(\a)&=&\frac{5}{2}\,\a+O(\a^3)\nonumber\\
b_{5/2,2}(\a)&=&O(\a^2)\nonumber\\
b_{5/2,3}(\a)&=&O(\a^3)\nonumber\\
\end{eqnarray*}
Then, letting $$a_i=\frac{1}{\hat m_i}\left(\frac{\L_i}{\tilde m_i}\right)^{2}=\hat a_i \d^{N-i}\ ,\quad \hat a_i\in [\underline\a,\ovl\a]\ ,$$
we find
\beqa{summations}
{f}_{ii}&=&-2\,\frac{\bar m_i}{\L_i}\,\left[\sum_{1\leq k< i}\bar m_k\,a_{2000}(a_k,a_i)+\sum_{i< k\leq N}\bar m_k\,a_{2000}(a_i,a_k)\right]\nonumber\\
&=&-2\frac{\bar m_i}{\tilde m_i\sqrt{\hat m_i\hat a_i \d^{N-i}}}\left[\sum_{1\leq k< i}\left(\frac{3}{8}\bar m_k\frac{\hat a_k^2 \d^{2N-2k}}{\hat a_i^3 \d^{3N-3i}}+O\left(\d^{5i-4k-N}\right)\right)\right.\nonumber\\
&+&\left.\sum_{i< k\leq N}\left(\frac{3}{8}\bar m_k\frac{\hat a_i^2 \d^{2N-2i}}{\hat a_k^3 \d^{3N-3k}}+O\left(\d^{5k-4i-N}\right)\right)\right]\nonumber\\
&=&-\frac{3}{4}\frac{\bar m_i}{\tilde m_i\sqrt{\hat m_i\hat a_i}}\left[\sum_{1\leq k< i}\left(\bar m_k\frac{\hat a_k^2}{\hat a_i^3}\d^{(7i-4k-3N)/{2}}+O(\d^{(11i-8k-3N)/2})\right)\right.\nonumber\\
&+&\left.\sum_{i<k\leq N}\left(\bar m_k\frac{\hat a_i^2}{\hat a_k^3}\d^{(6k-3i-3N)/{2}}+O(\d^{(10k-7i-3N)/2})\right)\right]
\eeqa 
So, when $i=1$, only the second summation appears and the dominant term is the one with $k=2$, namely, the term
$$-\frac{3}{4}\frac{\bar m_1\bar m_2}{\tilde m_1\sqrt{\hat m_1\hat a_1}}\frac{\hat a_1^2}{\hat a_2^3}\d^{3(3-N)/2}$$
and the dominant neglected term is of order $\d^{(13-3N)/2}$.\\
When $i=N$, only the first summation appears in (\ref{summations}), the dominant term is reached in the sum for $k=N-1$:
$$-\frac{3}{4}\frac{\bar m_{N-1}\bar m_N}{\tilde m_N\sqrt{\hat m_N\hat a_N}}\frac{\hat a_{N-1}^2}{\hat a_N^3}\d^{2}$$
and the dominant neglected term is of order $\d^4$.\\
When $2\leq i\leq N-1$ (for $N\geq 3$ planets), the first summation gives the lowest order term  for $k=i-1$
$$-\frac{3}{4}\frac{\bar m_{i-1}\bar m_i}{\tilde m_i\sqrt{\hat m_i\hat a_i}}\frac{\hat a_{i-1}^2}{\hat a_i^3}\d^{(3i+4-3N)/2}$$
which is dominant with respect to the dominant term of the second summand
$$-\frac{3}{4}\frac{\bar m_{i}\bar m_{i+1}}{\tilde m_i\sqrt{\hat m_i\hat a_i}}\frac{\hat a_{i}^2}{\hat a_{i+1}^3}\d^{(3i+6-3N)/2}$$ 
which, on turn, is the dominant neglected term. The evaluation of the off--diagonal elements of $\cF(\L)$ is easier. In fact, we find, for $i< j$,
\beqano
f_{ij}&=&-\frac{\bar m_i\bar m_j}{\sqrt{\L_i\L_j}}{a_{1100}(a_i,a_j)}\nonumber\\
&=&-\frac{\bar m_i\bar m_j}{\sqrt{\tilde m_i\tilde m_j\sqrt{\hat m_i\hat m_j\hat a_i\hat a_j\d^{2N-i-j}}}}\left[-\frac{15}{16}\frac{\hat a_i^3\d^{3N-3i}}{\hat a_j^4\d^{4N-4j}}+O\left(\frac{\d^{5N-5i}}{\d^{6N-6j}}\right)\right]\nonumber\\
&=&\frac{15}{16}\frac{\bar m_i\bar m_j}{\sqrt{\tilde m_i\tilde m_j\sqrt{\hat m_i\hat m_j\hat a_i\hat a_j}}}\frac{\hat a_i^3}{\hat a_j^4}\d^{(17j-11i-6N)/4}+O\left(\d^{(25j-19i-6N)/4}\right)
\eeqano
because
\beqano
a_{1100}({a},{b})&=&\frac{{a}}{8{b}^2}\,\left[-17\,{a}/{b}\,b_{5/2,1}({a}/{b})\right.\nonumber\\
&&+\left.8(1+{a}^2/{b}^2)\,b_{5/2,2}({a}/{b})+{a}/{b}\,\,b_{5/2,3}({a}/{b})\right]\nonumber\\
&=&-\frac{15}{16}\frac{a^3}{b^4}+O\left(\frac{a^5}{b^6}\right)
\eeqano
having used
\beqano
b_{5/2,1}(\a)&=&\frac{5}{2}\,\a+O(\a^3)\nonumber\\
b_{5/2,2}(\a)&=&\frac{35}{8}\,\a^2+O(\a^4)\nonumber\\
b_{5/2,3}(\a)&=&O(\a^3)\nonumber\\
\eeqano
\subsubsection{Proof of non Degeneracy}
The aim of this section is to prove the non degeneracy of the Plane Planetary Problem:
\begin{lemma}\label{non degeneracy lem}
There exists $\d^*>0$ such that, for any $0<\d<\d^*$ and $0<\m<1$, the matrix of the Birkhoff invariants with order $2$ for $\bar f$ is non singular on ${\cal L}_{\underline\a, \ovl\a}(\m,\d)$:
$$\inf_{{\cal L}_{\underline\a, \ovl\a}(\m,\d)}|{\rm det}\,A|>0\ .$$
\end{lemma}
We will obtain this result as a consequence of Lemma \ref{Aij asymptotics} below, which will be obtained by direct check.
\begin{lemma}\label{Aij asymptotics}
The matrix $A_{i,j}(\L)$ of the Birkhoff invariants with  order $2$  for $\bar f(\L,\cdot,\cdot)$ satisfies the following asymptotics
\begin{eqnarray}\label{order 2 asymp}
A_{i,j}=\d^{5-2N}\times\left\{
\begin{array}{lrr}
+\frac{3}{4}\frac{\bar m_1\,\bar m_2}{\,\tilde m_1^2\,\hat m_1}\,\frac{(\hat a_1/\hat a_2)^2}{\hat a_1\hat a_2}(1+\a_{11})\quad \textrm{for}\quad i=j=1\ ;\\
\\
-3\frac{\bar m_i\,\bar m_{i-1}}{\tilde m_i^2\,\hat m_i}\,\frac{(\hat a_{i-1}/\hat a_i)^2}{\hat a_i^2}\,\d^{2i-3}(1+\a_{ii})\quad \textrm{for}\quad i=j=2,\cdots,N\ ;\\
\\
-\frac{9}{4}\,\frac{\bar m_i\,\bar m_j}{\tilde m_i\,\tilde m_j\,\sqrt{{\hat m}_i\,{\hat m}_j\hat a_i\hat a_j}}\frac{\hat a_i^2}{\hat a_j^3}\,\d^{(7j-3i-10)/2}(1+\a_{ij})\quad \textrm{for}\quad i<j\\
\\
A_{j,i}\quad \textrm{for}\quad i>j\ .
\end{array}
\right.
\end{eqnarray}
where $\a_{ij}=O(\d)$.
\end{lemma}
We show here how Lemma \ref{non degeneracy lem} follows from Lemma \ref{Aij asymptotics}

\vskip.1in
\noi
{\bf Proof of Lemma \ref{non degeneracy lem}.}\ Put, for shortness,
\begin{eqnarray*}
\tilde A_{i,j}:=\left\{
\begin{array}{lrr}
+\frac{3}{4}\frac{\bar m_1\,\bar m_2}{\,\tilde m_1^2\,\hat m_1}\,\frac{(\hat a_1/\hat a_2)^2}{\hat a_1\hat a_2}(1+\a_{11})\quad \textrm{for}\quad i=j=1\ ;\\
\\
-3\frac{\bar m_i\,\bar m_{i-1}}{\tilde m_i^2\,\hat m_i}\,\frac{(\hat a_{i-1}/\hat a_i)^2}{\hat a_i^2}(1+\a_{ii})\quad \textrm{for}\quad i=j=2,\cdots,N\ ;\\
\\
-\frac{9}{4}\,\frac{\bar m_i\,\bar m_j}{\tilde m_i\,\tilde m_j\,\sqrt{{\hat m}_i\,{\hat m}_j\hat a_i\hat a_j}}\frac{\hat a_i^2}{\hat a_j^3}(1+\a_{ij})\quad \textrm{for}\quad i<j\\
\\
\tilde A_{j,i}\quad \textrm{for}\quad i>j\ ,
\end{array}
\right.
\end{eqnarray*}
In the case $N=2$ (Plane $3$--Body Problem), the matrix $A$ is
$$
A=\d\tilde A_2:=\d\left(
\begin{array}{lr}
\tilde A_{11}&\tilde A_{12}\d^{1/2}\\
\tilde A_{21}\d^{1/2}&\tilde A_{22}\d
\end{array}
\right)
$$
so, 
\beqano
\textrm{\rm det}\,A&=&\d^3(\tilde A_{11}\tilde A_{22}-\tilde A_{12}^2)\nonumber\\
&=&-\frac{117}{16}\frac{\bar m_1^2\,\bar m_2^2}{\tilde m_1^2\,\tilde m_2^2\,{{\hat m}_1\,{\hat m}_2}}\frac{\hat a_1^3}{\hat a_2^7}\d^3(1+D_2)
\eeqano
where $D_2=O(\d)$, hence, for a suitable $\d^*>0$, $\textrm{\rm det}\,A_2\neq 0$ on ${\cal L}_{\underline\a,\overline\a}(\m,\d)$ for $0<\d<\d^*$. The following claim concludes the proof.
\vskip.1in
\noi
$\underline{\textrm{\sl Claim}:}$ {\it For $N\geq 3$,
 the matrix $\tilde A:=\d^{-(5-2N)}A$ has determinant $$\textrm{\rm det}\tilde A=-\frac{117}{16}\frac{\bar m_1^2\,\bar m_2^2}{\tilde m_1^2\,\tilde m_2^2\,{{\hat m}_1\,{\hat m}_2}}\frac{\hat a_1^3}{\hat a_2^7}\prod_{3\leq i\leq N}\left(-3\frac{\bar m_i\,\bar m_{i-1}}{\tilde m_i^2\,\hat m_i}\,\frac{(\hat a_{i-1}/\hat a_i)^2}{\hat a_i^2}\d^{2i-3}\right)(1+D)$$
where $D=O(\d)$.}
\vskip.1in
\noi
$\underline{\textrm{\sl Proof.}}$
We prove, by induction, that any matrix  with order $n$  the form of $\tilde A_n$
has determinant $$\textrm{\rm det}\tilde A_n=-\frac{117}{16}\frac{\bar m_1^2\,\bar m_2^2}{\tilde m_1^2\,\tilde m_2^2\,{{\hat m}_1\,{\hat m}_2}}\frac{\hat a_1^3}{\hat a_2^7}\prod_{3\leq i\leq n}\left(-3\frac{\bar m_i\,\bar m_{i-1}}{\tilde m_i^2\,\hat m_i}\,\frac{(\hat a_{i-1}/\hat a_i)^2}{\hat a_i^2}\right)(1+D_n)$$
where $D_n=O(\d)$. The claim will be then obtained taking $n=N$, $D:=D_N$.
For $n=3$, the claim is true by direct computation, since
$$
\tilde A_3=\left(
\begin{array}{lrrr}
\tilde A_{11}&\tilde A_{12}\d^{1/2}&\tilde A_{13}\,\d^{4}\\
\tilde A_{21}\d^{1/2}&\tilde A_{22}\d &\tilde A_{23}\d^{5/2}\\
\tilde A_{31}\d^{4}&\tilde A_{23}\d^{5/2}&\tilde A_{33}\d^{3} 
\end{array}
\right)\ .
$$
Assume now that the claim is true for $n-1$ and let us prove it for $n$. 
We write 
$$
\tilde A_n=\left(
\begin{array}{lrrr}
 & & &\tilde A_{1,n}\d^{(7n-13)/2}\\
& \tilde A_{n-1}&&\vdots\\
& & &\tilde A_{n-1,n}\d^{(4n-7)/2}\\
\tilde A_{n,1}\d^{(7n-13)/2}&\cdots& \tilde A_{n,n-1}\d^{(4n-7)/2}&\tilde A_{nn}\d^{2n-3}
\end{array}
\right)
$$
where $\tilde A_{n-1}$ is the submatrix of $\tilde A$ composed of the first $n-1$ rows and coloumns, which, by the inductive hypothesis, has determinant
$$\textrm{\rm det}\tilde A_{n-1}=-\frac{117}{16}\frac{\bar m_1^2\,\bar m_2^2}{\tilde m_1^2\,\tilde m_2^2\,{{\hat m}_1\,{\hat m}_2}}\frac{\hat a_1^3}{\hat a_2^7}\prod_{3\leq i\leq n-1}\left(-3\frac{\bar m_i\,\bar m_{i-1}}{\tilde m_i^2\,\hat m_i}\,\frac{(\hat a_{i-1}/\hat a_i)^2}{\hat a_i^2}\d^{2i-3}\right)(1+D_{n-1})$$
Then, expanding the determinant of $\tilde A_n$ along the $n^{th}$ coloumn, we find
\beqa{dominant term}
\textrm{\rm det} \tilde A_n&=&\tilde A_{nn}\d^{2n-3}\textrm{\rm det}\,\tilde A_{n-1}+R_n\nonumber\\
&=&-\frac{117}{16}\frac{\bar m_1^2\,\bar m_2^2}{\tilde m_1^2\,\tilde m_2^2\,{{\hat m}_1\,{\hat m}_2}}\frac{\hat a_1^3}{\hat a_2^7}\prod_{3\leq i\leq n}\left(-3\frac{\bar m_i\,\bar m_{i-1}}{\tilde m_i^2\,\hat m_i}\,\frac{(\hat a_{i-1}/\hat a_i)^2}{\hat a_i^2}\d^{2i-3}\right)(1+D_{n-1})\nonumber\\
&+&R_n
\eeqa
where
\beq{remaindern}R_n=\sum_{1\leq i\leq n-1}(-1)^{n-i}\tilde A_{in}\d^{(7n-3i-10)/2)}\textrm{\rm det}\,\tilde{\cal A}_{in}\eeq
if $\tilde{\cal A}_{in}$  is the $(i,n)$--minor of $\tilde A_{n}$, hence, given by
$$\tilde{\cal A}_{in}=\left(
\begin{array}{c}
\tilde r_1\\
\vdots\\
\tilde r_{i-1}\\
\tilde r_{i+1}\\
\vdots\\
\tilde r_{n-1}\\
\hat r_n
\end{array}
\right)$$
where $$\tilde r_i:=[\tilde r_{i1},\cdots\tilde r_{i,n-1}]\qquad 1\leq i\leq n-1$$
is the $i^{th}$ row of $\tilde A_{n-1}$ and $$\hat r_n=[\hat r_{n1},\cdots,\hat r_{n,n-1}]$$ is the $n^{th}$ row of $\tilde A$, deprivated of its $n^{th}$ component. We prove that the  remainder terms appearing in the summation (\ref{remaindern}) are at leas $\d$ times the dominant term 
\beqa{dominant termnew}-\frac{117}{16}\frac{\bar m_1^2\,\bar m_2^2}{\tilde m_1^2\,\tilde m_2^2\,{{\hat m}_1\,{\hat m}_2}}\frac{\hat a_1^3}{\hat a_2^7}\prod_{3\leq i\leq n}\left(-3\frac{\bar m_i\,\bar m_{i-1}}{\tilde m_i^2\,\hat m_i}\,\frac{(\hat a_{i-1}/\hat a_i)^2}{\hat a_i^2}\d^{2i-3}\right)\eeqa in (\ref{dominant term}), which will conclude the proof. We distinguish $2$ cases.

\vskip.1in
\noi
{\bf 1.} $1\leq i\leq N-2$.

\vskip.1in
\noi
Write $\hat r_n=\d^{(-3n+10+3i)/2}(\tilde r_i-\hat r_i)$, where the $k^{th}$ component of
$\hat r_i$ is $\tilde r_{ik}(1-\r_{ik})$, with
$$\r_{ik}=\frac{\tilde A_{nk}}{\tilde A_{ik}}\times\arr{
\d^{10(n-i-1)/2}\quad \textrm{for}\quad 1\leq k<i\\
\d^{(10(n-i)-14)/2}\quad \textrm{for}\quad k=i\\
\d^{(10(n-k-1)/2}\quad \textrm{for}\quad i<k\leq n-1\\
}\ .$$
Then,
\beqa{induction}
\textrm{\rm det}\tilde{\cal A}_{in}&=&\textrm{\rm det}\,\left(
\begin{array}{c}
\tilde r_1\\
\vdots\\
\tilde r_{i-1}\\
\tilde r_{i+1}\\
\vdots\\
\tilde r_{n-1}\\
\hat r_n
\end{array}
\right)\nonumber\\
&=&\d^{(-3n+10+3i)/2}\textrm{\rm det}\,\left(
\begin{array}{c}
\tilde r_1\\
\vdots\\
\tilde r_{i-1}\\
\tilde r_{i+1}\\
\vdots\\
\tilde r_{n-1}\\
\tilde r_i-\hat r_i
\end{array}
\right)\nonumber\\
&=&\d^{(-3n+10+3i)/2}\textrm{\rm det}\,\left(
\begin{array}{c}
\tilde r_1\\
\vdots\\
\tilde r_{i-1}\\
\tilde r_{i+1}\\
\vdots\\
\tilde r_{n-1}\\
\tilde r_i
\end{array}
\right)-\d^{(-3n+10+3i)/2}\textrm{\rm det}\,\left(
\begin{array}{c}
\tilde r_1\\
\vdots\\
\tilde r_{i-1}\\
\tilde r_{i+1}\\
\vdots\\
\tilde r_{n-1}\\
\hat r_i
\end{array}
\right)
\eeqa
where both the matrices appearing in (\ref{induction}) may be rearranged (changing their determinants at most for a sign) such in a way to take the form of $\tilde A_{n-1}$. So, using the inductive hypothesys we see that each term $\d^{(7n-3i-10)/2)}\textrm{\rm det}\,\tilde{\cal A}_{in}$ is of order at least $\d^3$ times the dominant term
(\ref{dominant termnew}).

\vskip.1in
\noi
{\bf 2.} $i=n-1$. In this case, we write $$\hat r_n=(\tilde r_{n-1}-\hat r_{n-1})\d^{3/2}$$
where $\hat r_{n-1}$ has components $\tilde r_{n-1,k}(1-\r_{n-1,k})$ with
$$\r_{n-1,k}=\arr{
\frac{\tilde A_{n,k}}{\tilde A_{n-1,k}}\d^{2}\quad \textrm{for}\quad 1\leq k\leq n-2\\
\frac{\tilde A_{n,k}}{\tilde A_{n-1,n-1}}\quad \textrm{for}\quad k=n-1\\
}
$$
So, using, as in the previous case, the inductive hypoyhesys for
\beqano
\textrm{\rm det}\tilde{\cal A}_{n-1,n}&=&\textrm{\rm det}\,\left(
\begin{array}{c}
\tilde r_1\\
\vdots\\
\tilde r_{n-2}\\
\hat r_n
\end{array}
\right)\nonumber\\
&=&\d^{3/2}\textrm{\rm det}\,\left(
\begin{array}{c}
\tilde r_1\\
\vdots\\
\tilde r_{i-1}\\
\tilde r_{i+1}\\
\vdots\\
\tilde r_{n-1}\\
\tilde r_{n-1}-\hat r_{n-1}
\end{array}
\right)\nonumber\\
&=&\d^{3/2}\textrm{\rm det}\,\left(
\begin{array}{c}
\tilde r_1\\
\vdots\\
\tilde r_{i-1}\\
\tilde r_{i+1}\\
\vdots\\
\tilde r_{n-1}\\
\tilde r_{n-1}
\end{array}
\right)-\d^{3/2}\textrm{\rm det}\,\left(
\begin{array}{c}
\tilde r_1\\
\vdots\\
\tilde r_{i-1}\\
\tilde r_{i+1}\\
\vdots\\
\tilde r_{n-1}\\
\hat r_{n-1}
\end{array}
\right)
\eeqano
and multiplying by $\d^{(4n-7)/2}$, we find that the term with $i=n-1$ in (\ref{remaindern}) is at least $\d$ times the dominant term (\ref{dominant termnew}).
This completes the proof.
\vskip.1in
\noi
The following Lemmas are devoted to the proof of Lemma \ref{Aij asymptotics}.

\begin{lemma}\label{order 2 invariants}
Let $U(\L)=(u_{ij}(\L))$ the  unitary matrix which diagonalizes $\cF(\L)$:
$$U(\L)^T\cF(\L)U(\L)={\rm diag}(\O_1,\cdots,\O_N)\ ,\quad U(\L)^T U(\L)=\id_N$$
and let $q_{i,j,k,l}(\L)$, $r_{i,j,k,l}(\L)$  the coefficients of the order $4$--expansion of $\bar f$ in Delaunay--Poincar\'e variables:
\beqano
\bar f(\L,\eta,\xi)&=&f_0(\L)+\frac{1}{2}\left(\eta\cdot\cF(\L)\eta+\xi\cdot\cF(\L)\xi\right)+\sum_{1\leq i,j,k,l\leq N}q_{i,j,k,l}\Big(\eta_i\eta_j\eta_k\eta_l+\xi_i\xi_j\xi_k\xi_l\Big)\nonumber\\
&+&\sum_{1\leq i,j,k,l\leq N}r_{i,j,k,l}(\L)\eta_i\eta_j\xi_k\xi_l+o_4\ .
\eeqano
Define
\beqa{q tilde r tilde}
\left\{\begin{array}{lrr}
\tilde q_{i,j,k,l}({\L}):=\sum_{1\leq i',j',k',l'\leq N}\,q_{i',j',k',l'}({\L})\,u_{i',i}({\L})\,u_{j',j}({\L})\,u_{k',k}({\L})\,u_{l',l}({\L})\\
\\
\tilde r_{i,j,k,l}({\L}):=\sum_{1\leq i',j',k',l'\leq N}\,r_{i',j',k',l'}({\L})\,u_{i',i}({\L})\,u_{j',j}({\L})\,u_{k',k}\,u_{l',l}({\L})
\end{array}
\right.
\eeqa
Then, the Birkhoff invariants with order $2$ for $\bar f$, namely, the elements of the symmetric matrix $A_{i,j}(\L)$ defining the Birkhoff normal form for $\bar f$ of order $2$
\footnote{The variables $\L$, $p$, $q$ are thought ``dummy'' in eq (\ref{BBNNFF}): the existence  of the Birkhoff transformaton realizing (\ref{BBNNFF}) has  been proved in the previous section.}
\beq{BBNNFF}f_0(\L)+\sum_{1\leq i\leq N}\O_i(\L)\frac{p_i^2+q_i^2}{2}+\frac{1}{2}\sum_{1\leq i,j\leq N}A_{i,j}(\L)\frac{p_i^2+q_i^2}{2}\frac{p_j^2+q_j^2}{2}\ ,\eeq
are given by 
\begin{eqnarray}\label{Aij}
A_{i,j}(\L):=\left\{
\begin{array}{lrr}
6\,\tilde q_{i,i,i,i}+\tilde r_{i,i,i,i}\quad \textrm{for}\quad i=j\\
\\
2\tilde q_{i,i,j,j}+2\tilde q_{j,j,i,i}+2\tilde q_{i,j,i,j}+2\tilde q_{j,i,j,i}+2\tilde q_{i,j,j,i}+2\tilde q_{j,i,i,j}+\tilde r_{i,i,j,j}+\tilde r_{j,j,i,i}\\ \textrm{for}\quad i\neq j
\end{array}
\right.
\end{eqnarray}
\end{lemma}
{\bf Proof.} The transformation
$$
\phi_{\rm D}:\quad \arr{
\L=\tilde\L\\
\l_i=\tilde\l_i-\tilde\xi\cdot U(\tilde\L)^T \partial_{\L_i}U(\tilde\L)\tilde\eta\\
\eta=U(\tilde\L)\tilde\eta\\
\xi=U(\tilde\L)\tilde\xi\\
}
$$
diagonalizes the quadratic part of $\bar f$, sending it to
$$\tilde f_2=\sum_{1\leq i\leq N}\O_i(\tilde\L)\frac{\tilde\eta_i^2+\tilde\xi_i^2}{2}$$
and puts the order $4$ term $f_4$ to  
\beq{tilde f4}
\tilde f_4=\sum_{1\leq i,j,k,l\leq N}\tilde q _{i,j,k,l}\Big({\tilde\eta}_i{\tilde\eta}_j{\tilde\eta}_k{\tilde\eta}_l+{\tilde\xi}_i{\tilde\xi}_j{\tilde\xi}_k{\tilde\xi}_l\Big)+\sum_{1\leq i,j,k,l\leq N}\tilde r_{i,j,k,l}(\tilde\L){\tilde\eta}_i{\tilde\eta}_j\tilde\xi_k\tilde\xi_l
\eeq
with $(\tilde q_{ijkl}):=\tilde {\cal Q}$, $(\tilde r_{ijkl}):=\tilde {\cal R}$ as in (\ref{q tilde r tilde}).
Next, as outlined in Appendix \ref{app:Birkhoff Normal Form}, Remark \ref{Birkhoff algorithm},  the Birkhoff invariants with order $2$ may easily be found through the identification
\beqa{BNF}
\frac{1}{2}\,\sum_{i,j=1}^N\,A_{i,j}(\tilde{\L})\,\frac{p_i^2+q_i^2}{2}\,\frac{p_j^2+q_j^2}{2}&:=&\frac{1}{(2\p)^N}\int_{[0,2\p]^N}\tilde f_{4}(\tilde\L,\tilde\eta,\tilde\xi)|_{(\tilde\eta_h,\tilde\xi_h)=\sqrt{p_h^2+q_h^2}\,(\cos{\varphi}_h,\sin{\varphi}_h)}d\varphi\nonumber\\
\eeqa
But, replacing $\tilde f_4$ as in (\ref{tilde f4}), the r.h.s in (\ref{BNF}) becomes
\beqa{summtions}
&&  \frac{1}{(2\p)^N}\int_{[0,2\p]^N}\tilde f_{4}(\tilde\L,\tilde\eta,\tilde\xi)|_{(\tilde\eta_h,\tilde\xi_h)=\sqrt{p_h^2+q_h^2}\,(\cos{\varphi}_h,\sin{\varphi}_h)}d\varphi\nonumber\\
&&  =2\sum_{1\leq i,j,k,l\leq N}\tilde q_{i,j,k,l}I_{ijkl}\sqrt{(p_i^2+q_i^2)(p_j^2+q_j^2)(p_k^2+q_k^2)(p_l^2+q_l^2)}\nonumber\\
&&  +\sum_{1\leq i,j,k,l\leq N}\tilde r_{i,j,k,l}J_{ijkl}\sqrt{(p_i^2+q_i^2)(p_j^2+q_j^2)(p_k^2+q_k^2)(p_l^2+q_l^2)}
\eeqa
where we have let
\beqano
I_{ijkl}&:=&\frac{1}{(2\p)^N}\int_{\torus^N}\cos\varphi_i\cos\varphi_j\cos\varphi_k\cos\varphi_ld\varphi\nonumber\\
&=&\frac{1}{(2\p)^N}\int_{\torus^N}\sin\varphi_i\sin\varphi_j\sin\varphi_k\sin\varphi_ld\varphi\nonumber\\
J_{ijkl}&:=&\frac{1}{(2\p)^N}\int_{\torus^N}\cos\varphi_i\cos\varphi_j\sin\varphi_k\sin\varphi_ld\varphi
\nonumber\\
\eeqano

\vskip.1in
\noi
The elementary integral $I_{ijkl}$ does not vanish only when
 $i=j=k=l$, or $i=j\neq k=l$, or $i=k\neq j=l$, or $i=l\neq j=k$. In the first case, it gives
 $$I_{iiii}=\frac{1}{2\pi}\,\int_0^{2\pi}\,\cos^4{x}\,dx=\frac{3}{8}\ .$$ In the three remaining cases, it gives
 $$I_{iijj}=I_{ijij}=I_{ijji}=\left(\frac{1}{2\pi}\,\int_0^{2\pi}\,\cos^2{x}\,dx\right)^2=\frac{1}{4}\qquad i\neq j$$
So, the first summation in (\ref{summations})
\beqa{I1}
&& \sum_{1\leq i,j,k,l\leq N}\tilde q_{i,j,k,l}I_{ijkl}\sqrt{(p_i^2+q_i^2)(p_j^2+q_j^2)(p_k^2+q_k^2)(p_l^2+q_l^2)}\nonumber\\
&& =\frac{3}{8}\sum_{1\leq i\leq N}\tilde q_{i,i,i,i}(\tilde\L)(p_i^2+q_i^2)^2\nonumber\\
&& +\frac{1}{4}\sum_{1\leq i\neq j\leq N} (\tilde q_{i,i,j,j}(\tilde\L)+\tilde q_{i,j,i,j}(\tilde\L)+\tilde q_{i,j,j,i}(\tilde\L))(p_i^2+q_i^2)(p_j^2+q_j^2)
\eeqa

\vskip.1in
\noi
Besides, the integral $J_{ijkl}$ does not vanish only when $i=j=k=l$ or $i=j\neq k=l$. In the first case, it gives
$$J_{iiii}=\frac{1}{2\p}\int_0^{2\pi}\cos^2{x}\sin^2{x}dx=\frac{1}{8}\ ,$$
in the second case, it gives
$$J_{iijj}=\left(\frac{1}{2\pi}\,\int_0^{2\pi}\,\cos^2{x}\,dx\right)\left(\frac{1}{2\pi}\,\int_0^{2\pi}\,\sin^2{x}\,dx\right)=\frac{1}{4}\qquad i\neq j$$
Then, the second summation in (\ref{summtions}) is
\beqa{I2}
&& \sum_{1\leq i,j,k,l\leq N}\tilde r_{i,j,k,l}J_{ijkl}\sqrt{(p_i^2+q_i^2)(p_j^2+q_j^2)(p_k^2+q_k^2)(p_l^2+q_l^2)}\nonumber\\
&& \frac{1}{8}\sum_{1\leq i\leq N}\tilde r_{i,i,i,i}(\tilde\L)(p_i^2+q_i^2)^2+\frac{1}{4}\sum_{1\leq i\neq j\leq N}\tilde r_{i,i,j,j}(\tilde\L)(p_i^2+q_i^2)(p_j^2+q_j^2)\ .\eeqa
Finally, replacing  (\ref{I1}) and  (\ref{I2}) into (\ref{BNF}) and next simmetrizing the summation, we find the result.

\vskip.1in
\noi
Our next step is the computation of the  asymptotics for ${\cal Q}$, ${\cal R}$, which, together with the one for the diagonalization matrix  $U(\L)$ (Lemma \ref{order 2 asymptotics}), will give the one for the $A_{ij}$'s.
\begin{lemma}\label{Q and R asymptotics}
The $\d$-- asymptotics for the functions ${\cal Q}=(q_{ijkl})$, ${\cal R}=(r_{ijkl})$ defining the $4$--expansion of $\bar f$:
\beqano
\bar f&=&f_0(\L)+\frac{1}{2}\Big(\eta\cF(\L)\eta+\xi\cdot\cF(\L)\xi\Big)\nonumber\\
&+&\sum_{1\leq i,j,k,l\leq N}q_{i,j,k,l}(\L)\Big(\eta_i\eta_j\eta_k\eta_l+\xi_i\xi_j\xi_k\xi_l\Big)+\sum_{1\leq i,j,k,l\leq N}r_{i,j,k,l}(\L)\eta_i\eta_j\xi_k\xi_l+o_4
\eeqano
is 
\begin{eqnarray}\label{qiiiiasymp}
q_{i,i,i,i}&=&\left\{
\begin{array}{lrr}
+\frac{3\,\bar m_1\,\bar m_2}{32\,\tilde m_1^2\,\hat m_1}\,\frac{(\hat a_1/\hat a_2)^2}{\hat a_1\hat a_2}\,\d^{5-2N}+O(\d^{7-2N})\quad \textrm{for}\quad i=1\ ;\\
\\
-\frac{3\,\bar m_i\,\bar m_{i-1}}{8\,\tilde m_i^2\,\hat m_i}\,\frac{(\hat a_{i-1}/\hat a_i)^2}{\hat a_i^2}\,\d^{2i+2-2N}+O\left(\d^{2i+3-2N}\right)\quad \textrm{for}\quad 2\leq i\leq N-1\ ;\\
\\
-\frac{3\,\bar m_N\,\bar m_{N-1}}{8\,\tilde m_N^2\,\hat m_N}\,\frac{(\hat a_{N-1}/\hat a_N)^2}{\hat a_N^2}\,\d^{2}+O\left(\d^{4}\right)\quad \textrm{for}\quad i=N\ ;\\
\end{array}
\right.
\nonumber\\
\\
\nonumber\\
q_{i,i,i,j}&=&r_{i,j,i,i}=r_{i,i,i,j}\nonumber\\
&=&\left\{
\begin{array}{lrr}
+\frac{75}{128}\,\frac{\bar m_i \bar m_j}{\sqrt{\tilde m_i^3\tilde m_j\sqrt{\hat m_i^3\hat m_j\hat a_i^3\hat a_j}}}\frac{\hat a_i^3}{\hat a_j^4}\,\d^{(17j-9i-8N)/4}\\
\\
+O\left(\d^{(25j-17i-8N)/4}\right)\quad \textrm{for}\quad i<j\ ;\\
\\
+\frac{285}{128}\,\frac{\bar m_i \bar m_j}{\sqrt{\tilde m_i^3\tilde m_j\sqrt{\hat m_i^3\hat m_j\hat a_i^3\hat a_j}}}\frac{\hat a_j^3}{\hat a_i^4}\,\d^{(19i-11j-8N)/4}\  ;\\
\\
+O\left(\d^{(27i-19j-8N)/4}\right)\quad \textrm{for}\quad i>j\ ;\\
\end{array}
\right.
\nonumber\\
\nonumber\\
\nonumber\\
\label{qiijjasymp}q_{i,i,j,j}&=&\left\{\begin{array}{lrr}
-\frac{9}{16}\,\frac{\bar m_i\,\bar m_j}{\tilde m_i\,\tilde m_j\,\sqrt{{\hat m}_i\,{\hat m}_j\hat a_i\hat a_j}}\frac{\hat a_i^2}{\hat a_j^3}\,\d^{(7j-3i-4N)/2}+O\left(\d^{(11j-7i-4N)/2}\right)\ ,\quad \textrm{for}\quad i<j\ ;\\
\\
\equiv 0\quad \textrm{for}\quad i>j\ ;\\
\end{array}
\right.
\nonumber\\
\\
\nonumber\\
\label{riiiiasymp}r_{i,i,i,i}&=&\left\{
\begin{array}{lrr}
+\frac{3\,\bar m_1\,\bar m_2}{16\,\tilde m_1^2\,\hat m_1}\,\frac{(\hat a_1/\hat a_2)^2}{\hat a_1\hat a_2}\,\d^{5-2N}+O(\d^{7-2N})\quad \textrm{for}\quad i=1\ ;\\
\\
-\frac{3\,\bar m_i\,\bar m_{i-1}}{4\,\tilde m_i^2\,\hat m_i}\,\frac{(\hat a_{i-1}/\hat a_i)^2}{\hat a_i^2}\,\d^{2i+2-2N}+O\left(\d^{2i+3-2N}\right)\quad \textrm{for}\quad 2\leq i\leq N-1\ ;\\
\\
-\frac{3\,\bar m_N\,\bar m_{N-1}}{4\,\tilde m_N^2\,\hat m_N}\,\frac{(\hat a_{N-1}/\hat a_N)^2}{\hat a_N^2}\,\d^{2}+O\left(\d^{4}\right)\quad \textrm{for}\quad i=N\ ;\\
\end{array}
\right.
\nonumber\\
\\
\nonumber\\
\label{riijjasymp}r_{i,i,j,j}&=&\left\{
\begin{array}{lrr}-\frac{9}{16}\,\frac{\bar m_i\,\bar m_j}{\tilde m_i\,\tilde m_j\,\sqrt{{\hat m}_i\,{\hat m}_j\hat a_i\hat a_j}}\frac{\hat a_i^2}{\hat a_j}\,\d^{(7j-3i-4N)/2}+O\left(\d^{(11j-7i-4N)/2}\right)\quad \textrm{for}\quad i<j\ ;\\
\\
-\frac{9}{16}\,\frac{\bar m_i\,\bar m_j}{\tilde m_i\,\tilde m_j\,\sqrt{\hat m_i\hat mj\hat a_i\hat a_j}}\frac{\hat a_j^2}{\hat a_i^3}\,\d^{(7i-3j-4N)/2}+O\left(\d^{(11i-7j-4N)/2}\right)\quad \textrm{for}\quad i>j\ .\\
\end{array}
\right.
\nonumber\\
\\
\nonumber\\
\label{rijijasymp}r_{i,j,i,j}&=&\left\{\begin{array}{lrr}
-\frac{315\,\bar m_i\,\bar m_j}{64\tilde m_i\tilde m_j\sqrt{\hat m_i\hat m_j\hat a_i\hat a_j}}\,\frac{\hat a_i^4}{\hat a_j^5}\,\d^{(11j-7i-4N)/2}+O(\d^{(15j-11i-4N)/2})\quad \textrm{for}\quad i<j\ ;\\
\\
\equiv 0\quad \textrm{for}\quad i>j\ .\\
\end{array}
\right.
\end{eqnarray}
\end{lemma}
\vskip.1in
\noindent
{\bf Proof.}
\vskip.1in
\noi 
{\bf 1.} {\sl Expansion of} $\dst q_{i,i,i,i}=-\frac{\bar m_i}{\L_i^2}\sum_{h \neq i}{\bar m_h}{a_{4000}(a_i,a_h)}\ ,\qquad \L_i=\tilde m_i\sqrt{\hat m_i\,a_i}$.

\vskip.1in
\noi
The  coefficient $a_{4000}(a,b)$ may be written simultaneoulsy as
\begin{eqnarray*}
a_{4000}(a,b)&=&\frac{a}{512b^2}\,\left[(-60(a/b)^5+4311(a/b)^3\right.\nonumber\\
&-&300(a/b))\,b_{9/2,0}(a/b)+8(7(a/b)^{6}\nonumber\\
&-&252(a/b)^4-222(a/b)^2+7)\,b_{9/2,1}(a/b)\nonumber\\
&+&4(75(a/b)^{5}-503(a/b)^{3}+135(a/b))\,b_{9/2,2}(a/b)\nonumber\\
&+&24(23(a/b)^{4}+13(a/b)^{2})\,b_{9/2,3}(a/b)\nonumber\\
&+&\left.37(a/b)^{3}\,b_{9/2,4}(a/b)\right]\nonumber\\
&=&\frac{b}{512a^2}\,\left[(-60(b/a)+4311(b/a)^{3}\right.\nonumber\\
&-&300(b/a)^{5})\,b_{9/2,0}(b/a)\nonumber\\
&+&8(7-252(b/a)^{2}-222(b/a)^{4}+7(b/a)^6)\,b_{9/2,1}(b/a)\nonumber\\
&+&4(75(b/a)-503(b/a)^{3}+135(b/a)^{5})\,b_{9/2,2}(b/a)\nonumber\\
&+&24(23(b/a)^{2}+13(b/a)^{4})\,b_{9/2,3}(b/a)\nonumber\\
&+&\left.37(b/a)^{3}\,b_{9/2,4}(b/a)\right]\nonumber\\
\end{eqnarray*}
having used
$$b_{9/2,k}(1/\a)=\a^9\,b_{9/2,k}(\a)\ .$$
The involved Laplace coefficients satisfy the following asymptotics (see Appendix \ref{Laplace coefficients}), for small $\a$
\begin{eqnarray}\label{Laplace asymp}
\left\{
\begin{array}{lrr}
b_{9/2,0}(\a)=1+O(\a^2)\\
b_{9/2,1}(\a)=\frac{9}{2}\,\a+O(\a^3)\\
b_{9/2,2}(\a)=\frac{99}{8}\,\a^2+O(\a^4)\\
b_{9/2,k}(\a)=O(\a^k)\ ,\quad \textrm{for}\quad k\geq 3\ ,
\end{array}
\right.
\end{eqnarray}
and we find thus the asymptotics for $a_{4000}(a,b)$
\begin{eqnarray*}
a_{4000}(a,b)=\left\{
\begin{array}{lrr}-\frac{3}{32}\,\frac{a^2}{b^3}+O\left(\frac{a^4}{b^5}\right)\qquad \textrm{for small}\qquad a/b\\
\\
+\frac{3}{8}\,\frac{b^2}{a^3}+O\left(\frac{b^4}{a^5}\right)\qquad \textrm{for small}\qquad b/a\\
\end{array}
\right.
\end{eqnarray*}
So, replacing
\begin{eqnarray}\label{a b}
a=a_i=\hat a_i\d^{N-i}\ , \qquad b=a_h=\hat a_h\,\d^{N-h}
\end{eqnarray}
we obtain
\begin{eqnarray*}
a_{4000}(a_i,a_h)=\left\{
\begin{array}{lrr}-\frac{3}{32}\,\frac{\hat a_i^2}{\hat a_h^3}\,\d^{3h-2i-N}+O\left(\d^{5h-4i-N}\right)\quad \textrm{for}\quad i<h\\
\\
+\frac{3}{8}\,\frac{\hat a_h^2}{\hat a_i^3}\,\d^{3i-2h-N}+O\left(\d^{5i-4h-N}\right)\quad \textrm{for}\quad i>h\\
\end{array}
\right.
\end{eqnarray*}
Finally,
\begin{eqnarray*}
q_{i,i,i,i}&=&-\frac{\bar m_i}{\L_i^2}\sum_{h\neq i}{\bar m_h}{a_{4000}(a_i,a_h)}\nonumber\\
&=&-\frac{\bar m_i}{\tilde m_i^2\,\hat m_i\,\hat a_i}\nonumber\\
&\times&\left\{\frac{3}{8}\,\sum_{h<i}\,\left[\bar m_h\,\frac{\hat a_h^2}{\hat a_i^3}\,\d^{4i-2h-2N}+O\left(\d^{6i-4h-2N}\right)\right]\right.\nonumber\\
&-&\left.\frac{3}{32}\,\sum_{h>i}\,\left[\bar m_h\,\frac{\hat a_i^2}{\hat a_h^3}\,\d^{3h-i-2N}+O\left(\d^{5h-3i-2N}\right)\right]\right\}
\end{eqnarray*}
The lowest order term is reached for $h=2$ when $i=1$, for $h=i-1$ when $i\geq 2$. The first neglected powers of $\d$ are the ones coming from the remainder term with $h=2$, for $i=1$, from the dominant term with $h=i+1$, for $i=2,\cdots,N-1$, from the remainder  term with $h=N-1$ when $i=N$. The final result is  then (\ref{qiiiiasymp}).

\vskip.1in
\noi
{\bf 2.} {\sl Expansion of}
$\dst r_{iiii}=-\frac{\bar m_i}{\L_i^2}\sum_{h\neq i}{\bar m_h}{a_{2020}(a_i,a_h)}$.

\vskip.1in
\noi
>From the identity
$$a_{2020}(a,b)=2\,a_{4000}(a,b)$$
one finds $r_{i,i,i,i}=2\,q_{i,i,i,i}$, hence, (\ref{riiiiasymp}).

\vskip.1in
\noi
{\bf 3.} {\sl Expansion of}
$\dst 
q_{i,i,i,j}=r_{i,i,i,j}=r_{i,j,i,i}=-\frac{\bar m_i\,\bar m_j\,a_{3100}(a_i,a_j)}{\L_i\,\sqrt{\L_i\,\L_j}}=-\frac{\bar m_i\,\bar m_j\,a_{1120}(a_i,a_j)}{\L_i\,\sqrt{\L_i\,\L_j}}$ $ \textrm{for}\quad i\neq j$.

\vskip.1in
\noi
Proceding similarly to the expansion of $q_{iiii}$, we write the coefficients $a_{3100}$, $a_{1120}$ as
\begin{eqnarray*}
a_{3100}(a,b)=a_{1120}(a,b)=\left\{\begin{array}{lrr}
-\frac{a}{256b^2}\left[(-744(a/b)^5+2014(a/b)^3\right.\\
-864(a/b))\,b_{9/2,1}(a/b)+8(28(a/b)^6\\
-321(a/b)^4-321(a/b)^2+28)\,\,b_{9/2,2}(a/b)\\
+(552(a/b)^5+423(a/b)^3+672(a/b))\,b_{9/2,3}(a/b)\\
+(1146(a/b)^4+1266(a/b)^2)\,b_{9/2,0}(a/b)\\
+6(29(a/b)^4+9(a/b)^2)\,b_{9/2,4}(a/b)\\
-\left.5(a/b)^3\,b_{9/2,5}(a/b)\right]\\
\\
-\frac{b}{256a^2}\left[(-744(b/a)+2014(b/a)^3\right.\\
-864(b/a)^5)\,b_{9/2,1}(b/a)+8(28\\
-321(b/a)^2-321(b/a)^4+28(b/a)^6)\,\,b_{9/2,2}(b/a)\\
+(552(b/a)+423(b/a)^3+672(b/a)^5)\,b_{9/2,3}(b/a)\\
+(1146(b/a)^2+1266(b/a)^4)\,b_{9/2,0}(b/a)\\
+6(29(b/a)^2+9(b/a)^4)\,b_{9/2,4}(b/a)\\
-\left.5(b/a)^3\,b_{9/2,5}(b/a)\right]\\
\end{array}
\right.\ .
\end{eqnarray*}
The asymptotics for these coefficients is computed using the asymptotics (\ref{Laplace asymp})
for the involved Laplace cefficients
\begin{eqnarray*}
a_{3100}(a,b)=a_{1120}(a,b)=\left\{\begin{array}{lrr}
-\frac{75}{128}\,\frac{a^3}{b^4}+O\left(\frac{a^5}{b^6}\right)\qquad \textrm{for small}\quad a/b\\
\\
-\frac{285}{128}\,\frac{b^3}{a^4}+O\left(\frac{b^5}{a^6}\right)\qquad \textrm{for small}\quad b/a\\
\end{array}
\right.
\end{eqnarray*}
Hence, replacing $a_i=\hat a_i \d^{N-i}$, $a_j=\hat a_i \d^{N-j}$
\begin{eqnarray*}
a_{3100}(a_i,a_j)=a_{1120}(a_i,a_j)=\left\{
\begin{array}{lrr}
-\frac{75}{128}\,\frac{\hat a_i^3}{\hat a_j^4}\,\d^{4j-3i-N}+O\left(\d^{6j-5i-N}\right)\ \ \textrm{for}\ i<j\\
\\
-\frac{285}{128}\,\frac{\hat a_j^3}{\hat a_i^4}\,\d^{4i-3j-N}+O\left(\d^{6i-5j-N}\right)\ \textrm{for}\ i>j\\
\end{array}
\right.
\end{eqnarray*}
and finally, the result, for $q_{iiij}$, $r_{ijii}$, multiplying in front by
$$-\frac{\bar m_i\bar m_j}{\sqrt{\L_i^3\L_j}}=-\frac{\bar m_i \bar m_j}{\sqrt{\tilde m_i^3\tilde m_j\sqrt{\hat m_i^3\hat m_j\hat a_i^3\hat a_j}}}\d^{(3i+j-4N)/4}$$

\vskip.1in
\noi
{\bf 4.} {\sl Expansion of}
$$q_{iijj}=\arr{
-{\bar m_i\,\bar m_j}\frac{a_{2200}(a_i,a_j)}{\L_i\,\L_j}\quad \textrm{for}\quad i<j\\
0\quad \textrm{for}\quad i>j
}$$
In order to compute the coefficient $q_{i,i,j,j}$, we
start from the symmetric coefficient
\begin{eqnarray*}
a_{2200}(a,b)=a_{2200}(b,a)&=&\frac{a}{512b^2}\,\left[(-324(a/b)^5+10584(a/b)^3-324(a/b))\,b_{9/2,0}(a/b)\right.\nonumber\\
&+&8(17(a/b)^6-300(a/b)^4-300(a/b)^2+17)\,b_{9/2,1}(a/b)\nonumber\\
&-&(1272(a/b)^5+6337(a/b)^3+1272(a/b))\,b_{9/2,2}(a/b)\nonumber\\
&+&(648(a/b)^6+396(a/b)^4+396(a/b)^2\nonumber\\
&+&648)\,b_{9/2,3}(a/b)+(348(a/b)^5\nonumber\\
&+&800(a/b)^3+348(a/b))\,b_{9/2,4}(a/b)\nonumber\\
&+&(-60a/b)^4-60(a/b)^2)\,b_{9/2,5}(a/b)\nonumber\\
&+&\left.9(a/b)^3\,b_{9/2,6}(a/b)\right]\nonumber\\
\end{eqnarray*}
and we find
\begin{eqnarray*}
a_{2200}(a,b)=
\frac{9}{16}\,\frac{a^2}{b^3}+O\left(\frac{a^4}{b^5}\right)\qquad \textrm{for small}\quad a/b
\end{eqnarray*} 
which gives, for $a=a_i=\hat a\d^{N-i}$, $b=a_j=\hat a\d^{N-j}$ and $i<j$,
\begin{eqnarray*}
a_{2200}(a_i,a_j)=
\frac{9}{16}\,\frac{\hat a_i^2}{\hat a_j^3}\,\d^{3j-2i-N}+O\left(\d^{5j-4i-N}\right)\quad \textrm{for}\quad i<j
\end{eqnarray*}
hence, (\ref{qiijjasymp}) follows, after multiplying by the factor
\begin{eqnarray*}
-\frac{\bar m_i\,\bar m_j}{\L_i\,\L_j}=-\frac{\bar m_i\bar m_j}{\tilde m_i\tilde m_j\sqrt{\hat m_i\hat m_j\hat a_i\hat a_j\d^{2N-i-j}}}
\end{eqnarray*}

\vskip.1in
\noi
{\bf 5.} {\sl Expansion of}
$$r_{iijj}=-{\bar m_i\,\bar m_j}\frac{a_{0220}(a_j,a_i)}{\L_i\,\L_j}\qquad \textrm{for}\qquad i\neq j\ .$$

\vskip.1in
\noi
We expand the symmetric coefficient
\begin{eqnarray*}
a_{0220}(a,b)&=&-\frac{3a}{512b^2}\,\left[(84(a/b)^5-8832(a/b)^3\right.\nonumber\\
&+&84(a/b))\,b_{9/2,0}(a/b)-8(5(a/b)^6\nonumber\\
&-&652(a/b)^4-652(a/b)^2+5)\,b_{9/2,1}(a/b)\nonumber\\
&-&5(328(a/b)^5-561(a/b)^3+328(a/b))\,b_{9/2,2}(a/b)\nonumber\\
&+&(216(a/b)^6-1020(a/b)^4\nonumber\\
&-&1020(a/b)^2+216)\,b_{9/2,3}(a/b)\nonumber\\
&+&(116(a/b)^5+200(a/b)^3\nonumber\\
&+&116(a/b))\,b_{9/2,4}(a/b)\nonumber\\
&-&(20(a/b)^4+20(a/b)^2)\,b_{9/2,5}(a/b)\nonumber\\
&+&\left.3(a/b)^3\,b_{9/2,6}(a/b)\right]\ .\nonumber\\
\end{eqnarray*}
and we find
\begin{eqnarray*}
a_{0220}(a,b)=\left\{
\begin{array}{lrr}
+\frac{9}{16}\,\frac{a^2}{b^3}+O\left(\frac{a^4}{b^5}\right)\qquad\textrm{for small}\quad a/b\\
\\
+\frac{9}{16}\,\frac{b^2}{a^3}+O\left(\frac{b^4}{a^5}\right)\qquad\textrm{for small}\quad b/a\\
\end{array}
\right.
\end{eqnarray*}
which gives (\ref{riijjasymp}) for
\begin{eqnarray*}
r_{i,i,j,j}&=&-\frac{\bar m_i\,\bar m_j}{\L_i\,\L_j}a_{0220}(a_i,a_j)\nonumber\\
\end{eqnarray*}

\vskip.1in
\noi
{\bf 6.} {\sl Expansion of}
$$r_{ijij}=\arr{-{\bar m_i\,\bar m_j}\frac{a_{1111}(a_i,a_j)}{\L_i\,\L_j}\qquad \textrm{for}\quad i<j\\
0 \qquad \textrm{for}\quad i>j\\
}$$

\vskip.1in
\noi
We expand the symmetric coefficient
\begin{eqnarray*}
a_{1111}(a,b)&=&\frac{a}{128b^2}\,\left[(-36(a/b)^5-7956(a/b)^3\right.\nonumber\\
&-&36(a/b))\,b_{9/2,0}(a/b)+8((a/b)^6\nonumber\\
&+&828(a/b)^4+828(a/b)^2+1)\,b_{9/2,1}(a/b)\nonumber\\
&+&(-3096(a/b)^5+1039(a/b)^3\nonumber\\
&-&3096(a/b))\,b_{9/2,2}(a/b)+(648(a/b)^6\nonumber\\
&-&1332(a/b)^4-1332(a/b)^2+648)\,b_{9/2,3}(a/b)\nonumber\\
&+&(348(a/b)^5+700(a/b)^3\nonumber\\
&+&348(a/b))\,b_{9/2,4}(a/b)-60((a/b)^4\nonumber\\
&+&\left.(a/b)^2)\,b_{9/2,5}(a/b)+9(a/b)^3\,b_{9/2,6}(a/b)\right]
\end{eqnarray*}
The term of order ($a^2/b^3$) in the expansion of $a_{1111}$ vanishes, so, we shall go on in the asymptotics for the involved Laplace Coefficients:
\begin{eqnarray*}
\left\{
\begin{array}{lrr}
b_{9/2,0}(\a)=1+\frac{81}{4}\,\a^2+O(\a^4)\\
b_{9/2,1}(\a)=\frac{9}{2}\,\a+\frac{891}{16}\,\a^3+O(\a^5)\\
b_{9/2,2}(\a)=\frac{99}{8}\,\a^2+O(\a^4)\\
b_{9/2,3}(\a)=\frac{429}{16}\,\a^3+O(\a^5)\\
b_{9/2,k}(\a)=O(\a^k)\quad \textrm{for}\quad k\geq 4\ .
\end{array}
\right.
\end{eqnarray*}
We find, for small $a/b$,
\begin{eqnarray*}
a_{1111}(a,b)&=&\left(-7956-36\cdot \frac{81}{4}+8\cdot 828\cdot \frac{9}{2}+8\cdot \frac{891}{16}-3096 \cdot \frac{99}{8}+648 \cdot \frac{429}{16}\right)\nonumber\\
&\times&\frac{a^4}{128\,b^5}+O\left(\frac{a^6}{b^7}\right)\nonumber\\
&=&\frac{315}{64}\,\frac{a^4}{b^5}+O\left(\frac{a^6}{b^7}\right)\ ,
\end{eqnarray*}
consequentely, taking, for $i<j$, $a=a_i=\hat a_i\d^{N-i}$, $b=a_j=\hat a_j\d^{N-j}$, we find, for $r_{ijij}$ the expansion (\ref{rijijasymp}).
This completes the proof.
\vskip.1in
\noindent 
We are ready for the 
{\bf proof of Lemma \ref{Aij asymptotics}}. 
\vskip.1in
\noi
For any $i<j$, the functions $\tilde q_{iiii}$, $\cdots$ involved into equation (\ref{Aij}) of Lemma \ref{order 2 invariants} may be written as
\beqa{dominant terms on Aij}
\arr{
\tilde q_{i,i,i,i}=q_{iiii}(1+\k_i)\\
\tilde r_{i,i,i,i}=r_{iiii}(1+\r_i)\\
\tilde q_{iijj}=q_{iijj}(1+\k_{ij})\\
\tilde q_{jjii}=q_{iijj}\k_{ji}\\
\tilde q_{i,j,i,j}=q_{iijj}\hat\k_{ij}\\
\tilde q_{j,i,j,i}=q_{iijj}\hat\k_{ji}\\
\tilde q_{i,j,j,i}=q_{iijj}\check\k_{ij}\\
\tilde q_{jiij}=q_{iijj}\check\k_{ji}\\
\tilde r_{iijj}=r_{iijj}(1+\r_{ij})\\
\tilde r_{jjii}=r_{jjii}(1+\r_{ji})\\
}
\eeqa
where
\beqa{kidef}
\k_i&:=&(u_{i,i}^4-1)+q_{i,i,i,i}^{-1}\left[\sum_{k\neq i}q_{k,k,k,k}\,(u_{k,i})^4+\sum_{k\neq l}q_{k,k,k,l}\,(u_{k,i})^3\,u_{l,i}\right.\nonumber\\
&+&\left.\sum_{k< l}q_{k,k,l,l}\,(u_{k,i})^2\,(u_{l,i})^2\right]\nonumber\\
\r_i&:=&(u_{i,i}^4-1)+r_{i,i,i,i}^{-1}\left[\sum_{k\neq i}r_{k,k,k,k}\,(u_{k,i})^4+\sum_{k\neq l}r_{k,k,l,l}\,(u_{k,i})^2\,(u_{l,i})^2\right.\nonumber\\
&+&\sum_{k\neq l}r_{k,l,k,k}\,(u_{k,i})^3\,u_{l,i}
+\left.\sum_{k\neq l}r_{k,k,k,l}\,(u_{k,i})^3\,u_{l,i}+\sum_{k< l}r_{k,l,k,l}\,(u_{k,i})^2\,(u_{l,i})^2\right]\ ,\nonumber\\
\k_{i,j}&:=&(u_{i,i}^2\,u_{j,j}^2-1)+q_{i,i,j,j}^{-1}\left[\sum_{k< l,\ (k,l)\neq (i,\ j)}q_{k,k,l,l}\,(u_{k,i})^2\,(u_{l,j})^2\right.\nonumber\\
&+&\left.\sum_{k}q_{k,k,k,k}\,(u_{k,i})^2\,(u_{k,j})^2+\sum_{k\neq l}q_{k,k,k,l}\,(u_{k,i})^2\,u_{k,j}\,u_{l,j}\right]\nonumber\\
\k_{j,i}&:=&q_{i,i,j,j}^{-1}\left[\sum_{k< l}q_{k,k,l,l}\,(u_{k,j})^2\,(u_{l,i})^2+\sum_{k}q_{k,k,k,k}\,(u_{k,i})^2\,(u_{k,j})^2\right.\nonumber\\
&+&\left.\sum_{k\neq l}q_{k,k,k,l}\,(u_{k,j})^2\,u_{k,i}\,u_{l,i}\right]\nonumber\\
\hat\k_{ij}&:=&q_{iijj}^{-1}\left[\sum_{k}q_{k,k,k,k}\,(u_{k,i})^2\,(u_{k,j})^2+\sum_{k\neq l}q_{k,k,k,l}\,(u_{k,i})^2\,u_{k,j}\,u_{l,j}\right.\nonumber\\
&+&\left.\sum_{k< l}q_{k,k,l,l}\,u_{k,i}\,u_{k,j}\,u_{l,i}\,u_{l,j}\right]\nonumber\\
\hat\k_{ji}&:=&q_{iijj}^{-1}\left[\sum_{k}q_{k,k,k,k}\,(u_{k,i})^2\,(u_{k,j})^2+\sum_{k\neq l}q_{k,k,k,l}\,(u_{k,j})^2\,u_{k,i}\,u_{l,i}\right.\nonumber\\
&+&\left.\sum_{k< l}q_{k,k,l,l}\,u_{k,i}\,u_{k,j}\,u_{l,i}\,u_{l,j}\right]\nonumber\\
\check\k_{ij}&:=&q_{iijj}^{-1}\left[\sum_{k}q_{k,k,k,k}\,(u_{k,i})^2\,(u_{k,j})^2+\sum_{k\neq l}q_{k,k,k,l}\,u_{k,i}\,(u_{k,j})^2\,u_{l,i}\right.\nonumber\\
&+&\left.\sum_{k< l}q_{k,k,l,l}\,u_{k,i}\,u_{k,j}\,u_{l,j}\,u_{l,i}\right]\nonumber\\
\check\k_{ji}&:=&q_{iijj}^{-1}\left[\sum_{k}q_{k,k,k,k}\,(u_{k,i})^2\,(u_{k,j})^2+\sum_{k\neq l}q_{k,k,k,l}\,u_{k,j}\,(u_{k,i})^2\,u_{l,j}\right.\nonumber\\
&+&\left.\sum_{k< l}q_{k,k,l,l}\,u_{k,i}\,u_{k,j}\,u_{l,j}\,u_{l,i}\right]\nonumber\\
\r_{ij}&:=&(u_{i,i})^2\,(u_{j,j})^2-1\nonumber\\
&+&r_{iijj}^{-1}\left[\sum_{k\neq l,\ (k,l)\neq (i,j)}r_{k,k,l,l}\,(u_{k,i})^2\,(u_{l,j})^2+\sum_{k}r_{k,k,k,k}\,(u_{k,i})^2\,(u_{k,j})^2\right.\nonumber\\
&+&\left.\sum_{k\neq l}r_{k,l,k,k}\,u_{k,i}\,u_{l,i}\,(u_{k,j})^2+\sum_{k\neq l}r_{k,k,k,l}\,(u_{k,i})^2\,u_{k,j}\,u_{l,j}\right.\nonumber\\
&+&\left.\sum_{k< l}r_{k,l,k,l}\,u_{k,i}\,u_{l,i}\,u_{k,j}\,u_{l,j}\right]\nonumber\\
\r_{ji}&:=&(u_{i,i})^2\,(u_{j,j})^2-1\nonumber\\
&+&r_{jjii}^{-1}\left[\sum_{k\neq l,\ (k,l)\neq (j,i)}r_{k,k,l,l}\,(u_{k,j})^2\,(u_{l,i})^2+\sum_{k}r_{k,k,k,k}\,(u_{k,j})^2\,(u_{k,i})^2\right.\nonumber\\
&+&\left.\sum_{k\neq l}r_{k,l,k,k}\,u_{k,j}\,u_{l,j}\,(u_{k,i})^2+\sum_{k\neq l}r_{k,k,k,l}\,(u_{k,j})^2\,u_{k,i}\,u_{l,i}\right.\nonumber\\
&+&\left.\sum_{k< l}r_{k,l,k,l}\,u_{k,i}\,u_{l,i}\,u_{k,j}\,u_{l,j}\right]\nonumber\\
\eeqa
Let us prove, for istance, the first in (\ref{dominant terms on Aij}), with $\k_i$ as in (\ref{kidef}) (since the other ones equalities are similar).
Taking into account only the non vanishing components of ${\cal Q}=q_{i,j,k,l}$, we write $\tilde q_{i,i,i,i}$ as
\begin{eqnarray}\label{qiiii}
\tilde q_{i,i,i,i}&=&\sum_{i',j',k',l'}q_{i',j',k',l'}\,u_{i',i}\,u_{j',i}\,u_{k',i}\,u_{l',i}\nonumber\\
&=&\sum_{k}q_{k,k,k,k}\,(u_{k,i})^4+\sum_{k\neq l}q_{k,k,k,l}\,(u_{k,i})^3\,u_{l,i}+\sum_{k< l}q_{k,k,l,l}\,(u_{k,i})^2\,(u_{l,i})^2\nonumber\\
&=&q_{i,i,i,i}+\left[q_{i,i,i,i}\,(u_{i,i}^4-1)+\sum_{k\neq i}q_{k,k,k,k}\,(u_{k,i})^4+\sum_{k\neq l}q_{k,k,k,l}\,(u_{k,i})^3\,u_{l,i}\right.\nonumber\\
&+&\left.\sum_{k< l}q_{k,k,l,l}\,(u_{k,i})^2\,(u_{l,i})^2\right]\nonumber\\
&=&q_{i,i,i,i}(1+\k_i)\ .
\end{eqnarray}

\vskip.1in
\noi
$\underline{\textrm{\sl Claim}}$:\ {\it The functions $\k_i$, $\r_i$, $\k_{ij}$, $\cdots$ defined in (\ref{kidef}) are $O(\d)$.}
\vskip.1in
\noi
This claim follows by direct computation (throughout the asymptotics for ${\cal Q}$, ${\cal R}$, $U$ given in Lemmas \ref{diagonalization N bodies}, \ref{Q and R asymptotics}) of the order of the functions appearing in the summations defining the functions $\k_i$, $\r_i$, $\k_{ij}$, $\cdots$. The details of this computation, for sake of continuity, are postponed at the end of the present proof: an inspection of such orders, however,  shows that they never go under $O(\d)$.

\vskip.1in
\noi
To conclude the proof, we use  Lemma \ref{order 2 invariants}. Using equation (\ref{dominant terms on Aij}) into equation (\ref{Aij}) and using the asymptotics (\ref{qiiiiasymp}), (\ref{riiiiasymp}) for $q_{iiii}$, $r_{iiii}$, we find, for the diagonal elements of $A$, the asyptotics
\beqano
A_{i,i}&=&6\,\tilde q_{i,i,i,i}+\tilde r_{i,i,i,i}\nonumber\\
&=&6\,q_{i,i,i,i}(1+\k_i)+r_{i,i,i,i}(1+\r_i)\nonumber\\
&=&(6\,q_{i,i,i,i}+r_{i,i,i,i})(1+\a_{i,i})\nonumber\\
&=& (1+\a_i)\times\left\{
\begin{array}{lrr}
+\frac{3\,m_1\,m_2}{4\,\m_1^2\,M_1}\,\frac{(\hat a_1/\hat a_2)^2}{\hat a_1\hat a_2}\,\d^{5-2N}\quad \textrm{for}\quad i=1\ ;\\
\\
-\frac{3\,m_i\,m_{i-1}}{\m_i^2\,M_i}\,\frac{(\hat a_{i-1}/\hat a_i)^2}{\hat a_i^2}\,\d^{2i+2-2N}\quad \textrm{for}\quad 2\leq i\leq N\ 
\end{array}
\right.\nonumber\\
\eeqano
where
$$\a_{i,i}:=\frac{6\,q_{i,i,i,i}\k_i+r_{i,i,i,i}\r_i}{6\,q_{i,i,i,i}+r_{i,i,i,i}}$$
is an $O(\d)$. Similarly, taking into account the asymptotics (\ref{qiijjasymp}), (\ref{riijjasymp}) for $q_{iijj}$, $r_{iijj}$, for for the upper diagonal elements of $A$, we find
\beqano
A_{i,j}(\L)&=&
2\tilde q_{i,i,j,j}+2\tilde q_{j,j,i,i}+2\tilde q_{i,j,i,j}+2\tilde q_{j,i,j,i}+2\tilde q_{i,j,j,i}+2\tilde q_{j,i,i,j}+\tilde r_{i,i,j,j}+\tilde r_{j,j,i,i}\nonumber\\
&=& 2q_{iijj}(1+\k_{ij})+r_{i,i,j,j}(1+\r_{ij})+r_{j,j,i,i}(1+\r_{ji})\nonumber\\
&+&q_{iijj}(2\k_{ji}+2\hat k_{ij}+2\hat k_{ji}+2\check k_{ij}+2\check k_{ji})\nonumber\\
&=&(2q_{iijj}+2r_{iijj})(1+\a_{ij})\nonumber\\
&=&-\frac{9}{4}\,\frac{m_i\,m_j}{\m_i\,\m_j\,\sqrt{M_iMj\hat a_i\hat a_j}}\frac{\hat a_i^2}{\hat a_j^3}\,\d^{(7j-3i-4N)/2}(1+\a_{ij})\qquad (i<j)
\eeqano
(recall $r_{iijj}=r_{jjii}$)
with
$$\a_{ij}=\frac{(2\k_{ji}+2\hat k_{ij}+2\hat k_{ji}+2\check k_{ij}+2\check k_{ji})q_{iijj}+(\r_{ij}+\r_{ji})r_{iijj}}{2\k_{ij}+2q_{iijj}+2r_{iijj}}$$
an $O(\d)$ again.

\vskip.1in
\noi
$\underline{\textrm{\sl Proof of the Claim:}}$
Using the asymptotics for ${\cal Q}$, ${\cal R}$, $U$ given in Lemmas \ref{diagonalization N bodies}, \ref{Q and R asymptotics}, by direct check, that

\vskip.2in
\noi
{\bf 1.}\ for $1\leq i\leq N$,
\beqano
(u_{i,i})^4-1\ ,\ u_{i,i}^2u_{j,j}^2-1=O(\d^{5/2})\ ;\\ 
\eeqano

\vskip.2in
\noi
{\bf 2.}\ for $k\neq i$
\beqano
\frac{q_{kkkk}}{q_{iiii}}(u_{ki})^4\ ,\ \frac{r_{kkkk}}{r_{iiii}}(u_{ki})^4&=&\nonumber\\
O(\d^{19 k-32})\ &&\textrm{if}\ 1=i<k\\
O(\d^{19(k-i)-8})\ &&\textrm{if}\ 1<i<k\\
O(\d^{15 i-26})\ &&\textrm{if}\ 1=k<i\\
O(\d^{15(i-k)-8})\ &&\textrm{if}\ 1<k<i
\eeqano

\vskip.2in
\noi
{\bf 3.} for $k\neq l$,
\beqano
\frac{q_{kkkl}}{q_{iiii}}(u_{ki})^3(u_{li})\ ,\ \frac{r_{kkkl}}{r_{iiii}}(u_{ki})^3(u_{li})\ ,\ \frac{r_{klkk}}{r_{iiii}}(u_{ki})^3(u_{li})&=&\nonumber\\
O(\d^{7k-8})\ &&\textrm{if}\ 1=l<i=k\\
O(\d^{7(k-l)-4})\ &&\textrm{if}\ 1<l<i=k\\
O(\d^{(17l-29)/2})\ &&\textrm{if}\ i=k=1<l\\
O(\d^{(17(l-k)-8)/2})\ &&\textrm{if}\ 1<i=k<l\\
O(\d^{15l-26})\ &&\textrm{if}\ 1=k<l=i\\
O(\d^{15l-15k-8})\ &&\textrm{if}\ 1<k<l=i\\
O(\d^{(35k-59)/2})\ &&\textrm{if}\ 1=i=l<k\\
O(\d^{(35(k-l)-16)/2})\ &&\textrm{if}\ 1<i=l\\
O(\d^{(21k+17l-68)/2})\ &&\textrm{if}\ 1=i<k<l\\
O(\d^{(35k+3l-68)/2})\ &&\textrm{if}\ 1=i<l<k\\
O(\d^{(17(l-k)+38(k-i)-20)/2})\ &&\textrm{if}\ 1<i<k<l\\
O(\d^{(30(i-k)+17(l-i)-20)/2})\ &&\textrm{if}\ 1<k<i<l\\
O(\d^{15(i-l)+15(l-k)-10})\ &&\textrm{if}\ 1<k<l<i\\
O(\d^{15(i-k)+7(k-l)-10})\ &&\textrm{if}\ 1<l<k<i\\
O(\d^{(35(k-l)+38(l-i)-20)/2})\ &&\textrm{if}\ 1<i<l<k\\
O(\d^{10(k-i)+7(i-l)-10})\ &&\textrm{if}\ 1<l<i<k\\
\eeqano

\vskip.2in
\noi
{\bf 4.} for $k<l$,
\beqano
\frac{q_{kkll}}{q_{iiii}}(u_{ki})^2(u_{li})^2\ ,\ \frac{r_{kkll}}{r_{iiii}}(u_{ki})^2(u_{li})^2&=&\nonumber\\
O(\d^{12l-23})\ &&\textrm{if}\ 1=i=k<l\\
O(\d^{12(l-i)-6})\ &&\textrm{if}\ 1<i=k<l\\
O(\d^{10(i-k)-6})\ &&\textrm{if}\ 1<k<i=l\\
O(\d^{10i-18})\ &&\textrm{if}\ 1=k<i=l\\
O(\d^{12l+7k-34}\ &&\textrm{if}\ 1=i<k<l\\
O(\d^{12(l-k)+19(k-i)-10})\ &&\textrm{if}\ 1<i<k<l\\
O(\d^{12(l-i)+10i-22})\ &&\textrm{if}\ 1=k<i<l\\
O(\d^{12(l-i)+10(i-k)-10})\ &&\textrm{if}\ 1< k<i<l\\
O(\d^{15(i-l)+10l-22})\ &&\textrm{if}\ 1=k<l<i\\
O(\d^{15(i-l)+10(l-k)-10})\ &&\textrm{if}\ 1< k<l<i
\eeqano

\vskip.2in
\noi
{\bf 5.} for $k<l$
\beqano
\frac{r_{klkl}}{r_{iiii}}(u_{ki})^2(u_{li})^2&=&\nonumber\\
O(\d^{17l-23})\ &&\textrm{if}\ 1=i=k<l\\
O(\d^{14(l-i)-6})\ &&\textrm{if}\ 1<k=i<l\\
O(\d^{14i-22})\ &&\textrm{if}\ 1=k<i=l\\
O(\d^{14(i-k)-6})\ &&\textrm{if}\ 1<k<i=l\\
O(\d^{5k+14l-34})\ &&\textrm{if}\ 1=i<k<l\\
O(\d^{14(l-k)+17(k-i)-10})\ &&\textrm{if}\ 1<i<k<l\\
O(\d^{14(l-i)+12i-24})\ &&\textrm{if}\ 1=k<i<l\\
O(\d^{14(l-i)+12(i-k)-10})\ &&\textrm{if}\ 1<k<i<l\\
O(\d^{15(i-l)+12l-24})\ &&\textrm{if}\ 1=k<l<i
\eeqano

\vskip.1in
\noi
{\bf 6.} for $i<j$, $k<l$, $(i,j)\neq(k,l)$
\beqano
\frac{q_{kkll}}{q_{iijj}}(u_{k,i})^2(u_{l,j})^2,\ \frac{r_{kkll}}{r_{iijj}}(u_{k,i})^2(u_{l,j})^2&=&\nonumber\\
O(\d^{12(l-k)+19(k-j)+7j-17})\ && \textrm{if}\ 1=i<j<k<l\\
O(\d^{12(l-k)+19(k-j)+7(j-i)-8})\ &&\textrm{if}\ 1<i<j<k<l\\
O(\d^{12(l-j)+7j-17})\ &&\textrm{if}\ 1=i<k=j<l\\
O(\d^{12(l-j)+7(j-i)-8})\ &&\textrm{if}\ 1<i<k=j<l\\
O(\d^{12(l-j)+7k-17}))\ &&\textrm{if}\  1=i<k<j<l\\
O(\d^{12(l-j)+7(j-k)-8})\ &&\textrm{if}\ 1<i<k<j<l\\
O(\d^{7k-13})\ &&\textrm{if}\ 1=i<k<l=j\\
O(\d^{7(k-i)-4})\ &&\textrm{if}\ 1<i<k<l=j\\
O(\d^{5(j-l)+7k-17})\ &&\textrm{if}\ 1=i<k<l<j\\
O(\d^{5(j-l)+7(k-i)-8})\ &&\textrm{if}\ 1<i<k<l<j\\
O(\d^{12(l-j)-4})\ &&\textrm{if}\ 1=k=i<j<l\\
O(\d^{12(l-j)-4})\ &&\textrm{if}\ 1<k=i<j<l\\
O(\d^{5(j-l)-4})\ &&\textrm{if}\ 1=k=i<l<j\\
O(\d^{5(j-l)-4})\ &&\textrm{if}\ 1<k=i<l<j\\
O(\d^{12(l-j)+10i-20})\ &&\textrm{if}\ 1=k<i<j<l\\
O(\d^{12(l-j)+10(i-k)-8})\ &&\textrm{if}\ 1<k<i<j<l\\
O(\d^{10i-16})\ &&\textrm{if}\ 1=k<i<l=j\\
O(\d^{5(j-l)+10(i-k)-8})\ &&\textrm{if}\ 1<k<i<l<j\\
O(\d^{10(i-k)-4})\ &&\textrm{if}\ 1<k<i<l=j\\
O(\d^{5(j-l)+10i-20})\ &&\textrm{if}\ 1=k<i<l<j\\
O(\d^{5(j-i)+10i-20})\ &&\textrm{if}\ 1=k<l=i<j\\
O(\d^{5(j-i)+10(i-k)-8})\ &&\textrm{if}\ 1<k<l=i<j\\
O(\d^{5(j-i)+15(i-l)+10l-20})\ &&\textrm{if}\ 1=k<l<i<j\\
O(\d^{5(j-i)+15(j-l)+10(l-k)-8})\ &&\textrm{if}\ 1<k<l<i<j
\eeqano

\vskip.2in
\noi
{\bf 7.} for $i<j$, $k<l$, 
\beqano
\frac{q_{kkll}}{q_{iijj}}(u_{k,j})^2(u_{l,i})^2,\ \frac{r_{kkll}}{r_{iijj}}(u_{k,j})^2(u_{l,i})^2&=&\nonumber\\
O(\d^{12(l-k)+19(k-j)+7j-17})\ && \textrm{if}\ 1=i<j<k<l\\
O(\d^{12(l-k)+19(k-j)+7(j-i)-8})\ &&\textrm{if}\ 1<i<j<k<l\\
O(\d^{12(l-j)+7j-13})\ &&\textrm{if}\ 1=i<k=j<l\\
O(\d^{12(l-j)+7(j-i)-4})\ &&\textrm{if}\ 1<i<k=j<l\\
O(\d^{12(l-j)+17(j-k)+7k-17}))\ &&\textrm{if}\  1=i<k<j<l\\
O(\d^{12(l-j)+17(j-k)+7(k-i)-8})\ &&\textrm{if}\ 1<i<k<j<l\\
O(\d^{17(j-k)+7k-17})\ &&\textrm{if}\ 1=i<k<l=j\\
O(\d^{17(j-k)+7(k-i)-8})\ &&\textrm{if}\ 1<i<k<l=j\\
O(\d^{5(j-l)+17(l-k)+7k-17})\ &&\textrm{if}\ 1=i<k<l<j\\
O(\d^{5(j-l)+17(l-k)+7(k-i)-8})\ &&\textrm{if}\ 1<i<k<l<j\\
O(\d^{12(l-j)+17j-29})\ &&\textrm{if}\ 1=k=i<j<l\\
O(\d^{12(l-j)+17(j-i)-8})\ &&\textrm{if}\ 1<k=i<j<l\\
O(\d^{(17j-29)/2})\ &&\textrm{if}\ 1=k=i<j=l\\
O(\d^{17(j-i)-8})\ &&\textrm{if}\ 1<k=i<j=l\\
O(\d^{5j+12l-29})\ &&\textrm{if}\ 1=k=i<l<j\\
O(\d^{5(j-l)+17(l-i)-8})\ &&\textrm{if}\ 1<k=i<l<j\\
O(\d^{12(l-j)+17(l-i)+10i-20})\ &&\textrm{if}\ 1=k<i<j<l\\
O(\d^{12(l-j)+17(l-i)+10(i-k)-8})\ &&\textrm{if}\ 1<k<i<j<l\\
O(\d^{17(j-i)+10i-20})\ &&\textrm{if}\ 1=k<i<l=j\\
O(\d^{5(j-l)+17(l-i)+10(i-k)-8})\ &&\textrm{if}\ 1<k<i<l<j\\
O(\d^{17(j-i)+10(i-k)-8})\ &&\textrm{if}\ 1<k<i<l=j\\
O(\d^{5(j-l)+17(l-i)+10i-20})\ &&\textrm{if}\ 1=k<i<l<j\\
O(\d^{5j+5i-16})\ &&\textrm{if}\ 1=k<l=i<j\\
O(\d^{5(j-i)+10(i-k)-4})\ &&\textrm{if}\ 1<k<l=i<j\\
O(\d^{5(j-i)+15(i-l)+10l-20})\ &&\textrm{if}\ 1=k<l<i<j\\
O(\d^{5(j-i)+15(j-l)+10(l-k)-8})\ &&\textrm{if}\ 1<k<l<i<j
\eeqano

\vskip.2in
\noi
{\bf 8.} for $i<j$,
\beqano
\frac{q_{kkkk}}{q_{iijj}}(u_{ki})^2(u_{kj})^2,\ \frac{r_{kkkk}}{r_{iijj}}(u_{ki})^2(u_{kj})^2&=&\nonumber\\
O(\d^{5j-8})\ &&\textrm{if}\ 1=k=i<j\\
O(\d^{5(j-i)-2})\ &&\textrm{if}\ 1<k=i<j\\
O(\d^{7(j-i)-2})\ &&\textrm{if}\ k=j>i>1\\
O(\d^{7j-11})\ &&\textrm{if}\ k=j>i=1\\
O(\d^{10i+5j-25})\ &&\textrm{if}\ 1=k\neq i,j\\
O(\d^{12(k-j)+7(k-i)-6})\ &&\textrm{if}\ k>j,\ i>1\\
O(\d^{12(k-j)+7k-14})\ &&\textrm{if}\ k>j,\ i=1\\
O(\d^{7(k-i)+5(j-k)-6})\ &&\textrm{if}\ 1<i<k<j,\ k>1\\
O(\d^{7k+5(j-k)-15})\ &&\textrm{if}\ 1=i<k<j,\ k>1\\
O(\d^{10(i-k)+5(j-k)-6})\ &&\textrm{if}\  1<k<i
\eeqano

\vskip.2in
\noi
{\bf 9.} for $i<j$, $k\neq l$,
\beqano
\frac{q_{kkkl}}{q_{iijj}}(u_{ki})^2u_{kj}u_{lj},\ \frac{q_{kkkl}}{r_{iijj}}(u_{ki})^2u_{kj}u_{lj},\ \frac{q_{klkk}}{r_{iijj}}(u_{ki})^2u_{kj}u_{lj}&=&\nonumber\\
O(\d^{5j-8})\ &&\textrm{if}\ 1=k=i<j=l\\
O(\d^{5(j-k)-2})\ &&\textrm{if}\ 1<k=i<j=l\\
O(\d^{10i+5j-24})\ &&\textrm{if}\ 1=k<i<j=l\\
O(\d^{15(i-k)+5(j-i)-6}\ &&\textrm{if}\ 1<k<i<j=l\\
O(\d^{2k+5j-15})\ &&\textrm{if}\ 1=i<k<j=l\\
O(\d^{5(j-k)+7(k-i)-6})\ &&\textrm{if}\ 1<i<k<j=l\\
O(\d^{(35(k-j)+14j-30)/2})\ &&\textrm{if}\ 1=i<j,\ k>j=l\\
O(\d^{(35(k-j)+14(j-i)-12)/2})\ &&\textrm{if}\ 1<i<j,\ k>j=l\\
O(\d^{5j+2i-12})\ &&\textrm{if}\ 1=l<i=k<j\\
O(\d^{5(j-i)+7(i-l)-4})\ &&\textrm{if}\ 1<l<i=k<j\\
O(\d^{5j-10})\ &&\textrm{if}\ 1=i=k<l<j\\
O(\d^{5(j-i)-4})\ &&\textrm{if}\ 1<i=k<l<j\\
O(\d^{(17(l-j)+10j-20)/2})\ &&\textrm{if}\ 1=i=k<j<l\\
O(\d^{(17(l-j)+10(j-i)-8)/2})\ &&\textrm{if}\ 1<i=k<j<l\\
O(\d^{14(j-i)+7i-14})\ &&\textrm{if}\ 1=l<i<k=j\\
O(\d^{14(j-i)+7(i-l)-6})\ &&\textrm{if}\ 1<l<i<k=j\\
O(\d^{14j-23}\ &&\textrm{if}\ 1=l=i<k=j\\
O(\d^{(14(j-i)-6)}\ &&\textrm{if}\ 1<l=i<k=j\\
O(\d^{14(j-l)+7l-13})\ &&\textrm{if}\ 1=i<l<k=j\\
O(\d^{14(j-l)+7(l-i)-6}\ &&\textrm{if}\ 1<i<l<k=j\\
O(\d^{(17(l-j)+14j-30)/2}\ &&\textrm{if}\ 1=i<k=j<l\\
O(\d^{(17(l-j)+14(j-i)-12)/2}\ &&\textrm{if}\ 1<i<k=j<l\\
O(\d^{5(j-i)+15(i-k)+7k-16}\ &&\textrm{if}\ 1=l<k<i<j\\
O(\d^{5(j-i)+15(i-k)+7(k-l)-8})\ &&\textrm{if}\ 1<l<k<i<j\\
O(\d^{(17|l-j|+3(j-l)+20l+40i-104)/4})\ &&\textrm{if}\ 1=k<i<j,\ l\neq j,\ k\\
O(\d^{5(j-i)+15(i-k)+7k-16}\ &&\textrm{if}\ 1=l<k<i<j\\
O(\d^{(5(j-i)+15(i-k)+7(k-l)-8)}\ &&\textrm{if}\ 1<l<k<i<j\\
O(\d^{(17|l-j|+3(j-l)+20(l-k)+40(i-k)-32)/4}\ &&\textrm{if}\ 1<k<i<j,\ l>k,\ l\neq j\\
O(\d^{5j+9k-25})\ &&\textrm{if}\ 1=i=l<k<j\\
O(\d^{5(j-k)+13(k-l)+6l-17})\ &&\textrm{if}\ 1=i<l<k<j\\
O(\d^{(17|l-j|+3(j-l)+20l+8k-68)/4})\ &&\textrm{if}\ 1=i<k<j,\ l>k,\ l\neq j\\
O(\d^{5(j-k)+14(k-i)+7i-16})\ &&\textrm{if}\ 1=l<i<k<j\\
O(\d^{5(j-k)+14(k-i)+7(i-l)-8})\ &&\textrm{if}\ 1<l\leq i<k<j\\
O(\d^{5(j-k)+14(k-l)+7(l-i)-8})\ &&\textrm{if}\ 1<i<l<k<j\\
O(\d^{5(j-k)+7(k-i)-8})\ &&\textrm{if}\ 1<i<k<l<j\\
O(\d^{(17(l-j)+10(j-k)+14(k-i)-16)/2})\ &&\textrm{if}\ 1<i<k<j<l\\
O(\d^{9(k-j)+14j-25})\ &&\textrm{if}\ 1=i=l<j<k\\
O(\d^{(35(k-j)+28(j-l)+14l-34)/2})\ &&\textrm{if}\ 1=i<l\leq j<k\\
O(\d^{(35(k-l)+38(l-j)+14j-34)/2)})\ &&\textrm{if}\ 1=i<j<l<k\\
O(\d^{(17(l-k)+38(k-j)+14j-34)/2})\ &&\textrm{if}\ 1=i<j<k<l\\
O(\d^{(35(k-j)+28(j-i)+14i-32)/2})\ &&\textrm{if}\ 1=l<i<j<k\\
O(\d^{(35(k-j)+28(j-l)+14(l-i)-16)/2})\ &&\textrm{if}\ 1<i\leq l\leq j<k\\
O(\d^{(35(k-j)+28(j-i)+14(i-l)-16)/2})\ &&\textrm{if}\ 1<l<i<j<k\\
O(\d^{(35(k-l)+38(l-j)+14(j-i)-16)/2})\ &&\textrm{if}\ 1<i<j<l<k\\
O(\d^{(17(l-k)+38(k-j)+14(j-i)-16)/2})\ &&\textrm{if}\ 1<i<j<k<l
\eeqano
{\bf 10.} for $i<j$, $k\neq l$,
\beqano
\frac{q_{kkkl}}{q_{iijj}}(u_{kj})^2u_{ki}u_{li},\ \frac{r_{kkkl}}{r_{iijj}}(u_{kj})^2u_{ki}u_{li},\ \frac{r_{klkk}}{r_{iijj}}(u_{kj})^2u_{ki}u_{li}&=&\\
O(\d^{(27j-45)/2})\ && \textrm{if}\ 1=k=i<j=l\\
O(\d^{(27(j-k)-12)/2})\ && \textrm{if}\ 1<k=i<j=l\\
O(\d^{(27(j-i)+30i-52)/2})\ && \textrm{if}\ 1=k<i<j=l\\
O(\d^{(27(j-i)+30(i-k)-12)/2}\ && \textrm{if}\ 1<k<i<j=l\\
O(\d^{(27(j-k)+14k-34)/2})\ && \textrm{if}\ 1=i<k<j=l\\
O(\d^{(27(j-k)+14(k-i)-16)/2})\ && \textrm{if}\ 1<i<k<j=l\\
O(\d^{(35(k-j)+14j-34)/2})\ && \textrm{if}\ 1=i<j=l<k\\
O(\d^{(35(k-j)+14(j-i)-16)/2})\ && \textrm{if}\ 1<i<j=l<k\\
O(\d^{5j+2i-14})\ && \textrm{if}\ 1=l<i=k<j\\
O(\d^{5(j-i)+7(i-l)-6})\ && \textrm{if}\ 1<l<i=k<j\\
O(\d^{5j+17l-35})\ && \textrm{if}\ 1=i=k<l<j\\
O(\d^{(10(j-l)+27(j-i)-12)/2})\ && \textrm{if}\ 1<i=k<l<j\\
O(\d^{(17(l-j)+27j-45)/2})\ && \textrm{if}\ 1=i=k<j<l\\
O(\d^{(17(l-j)+27(j-i)-12)/2})\ && \textrm{if}\ 1<i=k<j<l\\
O(\d^{(11(j-i)+14i-24)/2})\ && \textrm{if}\ 1=l<i<k=j\\
O(\d^{11(j-i)+14(i-l)-8})\ && \textrm{if}\ 1<l<i<k=j\\
O(\d^{(11j-17)/2}\ && \textrm{if}\ 1=l=i<k=j\\
O(\d^{(11(j-i)-4)}\ && \textrm{if}\ 1<l=i<k=j\\
O(\d^{(11(j-l)+14l-22)/2})\ && \textrm{if}\ 1=i<l<k=j\\
O(\d^{(11(j-l)+14(l-i)-8)/2}\ && \textrm{if}\ 1<i<l<k=j\\
O(\d^{(17(l-j)+14j-26)/2}\ && \textrm{if}\ 1=i<k=j<l\\
O(\d^{(17(l-j)+14(j-i)-8)/2}\ && \textrm{if}\ 1<i<k=j<l\\
O(\d^{5(j-i)+15(i-k)+7k-16}\ && \textrm{if}\ 1=l<k<i<j\\
O(\d^{5(j-i)+15(i-k)+7(k-l)-8})\ && \textrm{if}\ 1<l<k<i<j\\
O(\d^{(17|l-i|+17(l-i)+40i+20j-104)/4})\ && \textrm{if}\ 1=k<i<j,\ l\neq j,\ k\\
O(\d^{5(j-i)+15(i-k)+7k-16}\ && \textrm{if}\ 1=l<k<i<j\\
O(\d^{(5(j-i)+15(i-k)+7(k-l)-8)}\ && \textrm{if}\ 1<l<k<i<j\\
O(\d^{(17|l-i|+17(l-i)+20(j-i)+60(i-k)-32)/4}\ && \textrm{if}\ 1<k<i<j,\ l>k,\ l\neq j\\
O(\d^{(10j+k-25)/2})\ && \textrm{if}\ 1=i=l<k<j\\
O(\d^{(9(k-l)+10(j-k)+12l-34)/2})\ && \textrm{if}\ 1=i<l<k<j\\
O(\d^{(34(l-k)+8k+20j-68)/4})\ && \textrm{if}\ 1=i<k<j,\ l>k,\ l\neq j\\
O(\d^{(10(j-k)+11(k-i)+14i-32)/2})\ && \textrm{if}\ 1=l<i<k<j\\
O(\d^{(10(j-k)+11(k-i)-12)/2})\ && \textrm{if}\ 1<l= i<k<j\\
O(\d^{(10(j-k)+11(k-i)+14(i-l)-16)2})\ && \textrm{if}\ 1<l< i<k<j\\
O(\d^{(10(j-k)+11(k-l)+14(l-i)-16)/2})\ && \textrm{if}\ 1<i<l<k<j\\
O(\d^{(10(j-l)+27(l-k)+14(k-i)-16)/2})\ && \textrm{if}\ 1<i<k<l<j\\
O(\d^{(17(l-j)+27(j-k)+14(k-i)-16)/2})\ && \textrm{if}\ 1<i<k<j<l\\
O(\d^{(18(k-j)+11j-25)/2})\ && \textrm{if}\ 1=i=l<j<k\\
O(\d^{(35(k-j)+14j-38)/2})\ && \textrm{if}\ 1=i<l=j<k\\
O(\d^{(35(k-j)+11(j-l)+14l-34)/2})\ && \textrm{if}\ 1=i<l<j<k\\
O(\d^{(35(k-l)+38(l-j)+14j-34)/2)})\ && \textrm{if}\ 1=i<j<l<k\\
O(\d^{(17(l-k)+38(k-j)+14j-34)/2})\ && \textrm{if}\ 1=i<j<k<l\\
O(\d^{(35(k-j)+11(j-i)+14i-32)/2})\ && \textrm{if}\ 1=l<i<j<k\\
O(\d^{(35(k-j)+11(j-i)-12)/2})\ && \textrm{if}\ 1<i=l\leq j<k\\
O(\d^{(35(k-j)+28(j-l)+14(l-i)-20)/2})\ && \textrm{if}\ 1<i<l=j<k\\
O(\d^{(35(k-j)+11(j-l)+14(l-i)-16)/2})\ && \textrm{if}\ 1<i\leq l\leq j<k\\
O(\d^{(35(k-j)+11(j-i)+14(i-l)-16)/2})\ && \textrm{if}\ 1<l<i<j<k\\
O(\d^{(35(k-l)+38(l-j)+14(j-i)-16)/2})\ && \textrm{if}\ 1<i<j<l<k\\
O(\d^{(17(l-k)+38(k-j)+14(j-i)-16)/2})\ && \textrm{if}\ 1<i<j<k<l
\eeqano

\vskip.1in
\noi
{\bf 11.} for $i<j$, $k<l$,
\beqano
\frac{q_{kkll}}{q_{iijj}}u_{ki}u_{kj}u_{li}u_{lj}\ ,\ \d^{-2(l-k)}\frac{r_{klkl}}{r_{iijj}}u_{ki}u_{kj}u_{li}u_{lj}&=&\nonumber\\
O(\d^{12(l-k)+19(k-j)+7j-17})\ && \textrm{if}\ 1=i<j<k<l\\
O(\d^{12(l-k)+19(k-j)+7(j-i)-8})\ &&\textrm{if}\ 1<i<j<k<l\\
O(\d^{12(l-j)+7j-15})\ &&\textrm{if}\ 1=i<k=j<l\\
O(\d^{12(l-j)+7(j-i)-6})\ &&\textrm{if}\ 1<i<k=j<l\\
O(\d^{(19(l-j)+12(j-k)+9k-34)/2}))\ &&\textrm{if}\  1=i<k<j<l\\
O(\d^{(24(l-j)+17(j-k)+14(k-i)-16)/2)})\ &&\textrm{if}\ 1<i<k<j<l\\
O(\d^{(17(j-k)+14k-30)/2})\ &&\textrm{if}\ 1=i<k<l=j\\
O(\d^{(17(j-k)+14(k-i)-12)/2})\ &&\textrm{if}\ 1<i<k<l=j\\
O(\d^{(10(j-l)+17(l-k)+14k-37)/2})\ &&\textrm{if}\ 1=i<k<l<j\\
O(\d^{(10(j-l)+17(l-k)+14(k-i)-16)/2})\ &&\textrm{if}\ 1<i<k<l<j\\
O(\d^{(24(l-j)+17j-33)/2})\ &&\textrm{if}\ 1=k=i<j<l\\
O(\d^{(24(l-j)+17(j-i)-12)/2})\ &&\textrm{if}\ 1<k=i<j<l\\
O(\d^{(17j-29)/2})\ &&\textrm{if}\ 1=k=i<j=l\\
O(\d^{(17(j-i)-8)/2})\ &&\textrm{if}\  1<k=i<j=l\\
O(\d^{(10j+7l-33)/2})\ &&\textrm{if}\ 1=k=i<l<j\\
O(\d^{(10(j-l)+17(l-i)-12)/2})\ &&\textrm{if}\ 1<k=i<l<j\\
O(\d^{(24(l-j)+17(l-i)+20i-40)/2})\ &&\textrm{if}\ 1=k<i<j<l\\
O(\d^{(24(l-j)+17(l-i)+20(i-k)-16)/2})\ &&\textrm{if}\ 1<k<i<j<l\\
O(\d^{(17j+3i-36)/2})\ &&\textrm{if}\ 1=k<i<l=j\\
O(\d^{(17(j-i)+20(i-k)-12)/2})\ &&\textrm{if}\ 1<k<i<l=j\\
O(\d^{(10j+7l+3i-40)/2})\ &&\textrm{if}\ 1=k<i<l<j\\
O(\d^{(10(j-l)+17(l-i)+3(i-k)-16)/2})\ &&\textrm{if}\ 1<k<i<l<j\\
O(\d^{(5j+5i-18)/2})\ &&\textrm{if}\ 1=k<l=i<j\\
O(\d^{5(j-i)+10(i-k)-6})\ &&\textrm{if}\ 1<k<l=i<j\\
O(\d^{5(j-i)+15(i-l)+5l-20})\ &&\textrm{if}\ 1=k<l<i<j\\
O(\d^{5(j-i)+15(j-l)+10(l-k)-8})\ &&\textrm{if}\ 1<k<l<i<j
\eeqano
\newpage
\section{Deprit Variables for the Spatial Planetary Problem}
\setcounter{equation}{0}
Consider the {\sl Spatial Planetary Problem}
\beqa{spatial problem}
{\cal H}_{\rm plt}(\m;y,x)=\sum_{1\leq i\leq N}\left(\frac{| y_i|^2}{2\tilde m_i}-\frac{\hat m_i\tilde m_i}{| x_i|}\right)+\m\sum_{1\leq i<j\leq N}\left(\frac{ y_i\cdot  y_j}{\bar m_0}-\frac{\bar m_i\bar m_j}{| x_i- x_j|}\right)
\eeqa
where
$(y,x)=\Big((y_1,\cdots, y_N),(x_1,\cdots, x_N)\Big)$ varies in the $6N$--dimensional collisionless phase space
$$\cC_{\rm cl,3}:=\left\{ y',\  x'\in \real^{3N}:\  x'_i\neq  x'_j\neq 0\quad\forall\ 1\leq i<j\leq N\right\}$$
and, as usual, 
$$\hat m_i=\bar m_0+\m\,\bar m_i\,\quad \tilde m_i=\frac{\bar m_0\bar m_i}{\bar m_0+\m\,\bar m_i}$$
are the reduced masses.

\vskip.1in
\noi
The system (\ref{spatial problem}) exhibits three integrals of the motion (besides the energy) related to its rotation invariance: the three components  of the  total angular momentum 
\beq{total angular momentum}
C=(C_{\rm x},C_{\rm y},C_{\rm z})=\sum_{1\leq i\leq N}x_i\times y_i\ .
\eeq
Hence, the number of degrees of freedom of (\ref{spatial problem}) can be furtherly reduced. Without performing such a reduction, any attempt of extending to the spatial case the strategy described in the previous section for the plane problem inevitably fails: two well known resonances, called {\sl secular resonances} (one of which with high order $2N-1$ and firstly noticed by M. Herman) appear, preventing the direct application of Theorem \ref{more general degenerate KAM}.

\vskip.1in
\noi
This section is devoted to the description of the reduction of the number of degrees of freedom of (\ref{spatial problem}), by means of a change of variables  essentially discovered, in the case of the Four Body Problem, by  Francoise Boigey \cite {BOI82} and then extended to the general case by A. Deprit (1926,2006), \cite{Dep83}. It may be viewed as a natural extension  of the {\sl Jacobi} or {\sl nodes reduction}, used in \cite{Rob95}, to prove the existence of quasi--periodic motions in the Three--Body Problem.

\vskip.1in
\noi
The three components $C_{\rm x}$, $C_{\rm y}$, $C_{\rm z}$ of the total angular momentum do not commute, but they verify the cyclic rules
\footnote{As usual, $\Big\{f,\ g\Big\}$ denotes the usual Poisson brackets of $f$, $g$:
$$\Big\{f,\ g\Big\}:=\sum_{1\leq i\leq N}\Big(\partial_{x_i}f\partial_{y_i}g-\partial_{y_i}f\partial_{x_i}g\Big)$$}
$$\Big\{C_{\rm x},\ C_{\rm y}\Big\}=C_{\rm z}\ ,\quad \Big\{C_{\rm y},\ C_{\rm z}\Big\}=C_{\rm x}\ ,\quad \Big\{C_{\rm z},\ C_{\rm x}\Big\}=C_{\rm y}\ .$$
However, as well known, starting with $C_{\rm x}$, $C_{\rm y}$, $C_{\rm z}$, it is possible to construct {\sl two} commuting integrals, for istance
$$C_{\rm z}\qquad \textrm{and}\qquad {\rm G}:=|C|=\sqrt{C_{\rm x}^2+C_{\rm y}^2+C_{\rm z}^2}\ .$$ We define then a system of (action--angle) symplectic coordinates, which are adapted to the reduction, since they  have $C_{\rm z}$ and ${\rm G}$ among their generalized momenta. The angle $\zeta$ conjugate to $C_{\rm z}$ is an integral of the motion too, implying that the Hamiltonian (\ref{spatial problem}), when expressed in such variables,  does not depend on the couple $(C_{\rm z},\zeta)$ and the angle ${\rm g}$ conjugate to ${\rm G}$. The constant value ${\rm G}={\rm G}_0$ will appear into the Hamiltonian as an ``external parameter'', meaning with this that the motion of the remaining $2(3N-2)$ variables will take place on a phase space parametrized by  ${\rm G}_0$. Owing to the rotation invariance of (\ref{spatial problem}), in particular, we find a set of symplectic variables on the manifold of dimension $2(3N-2)$
\beqano
\cM_{\rm vert,{\rm G}_0}:=\left\{y,\ x\in (\real^3)^N:\ C_{\rm x}=C_{\rm y}=C_{N{\rm x}}=0,\ C_{\rm z}={\rm G}_0\right\}\ ,
\eeqano
where $C_{N\rm x}$ denotes the first component of the $N^{th}$ angular momentum $C_N=x_N\times y_N$.  A further trivial integration will reconstruct the full motion on the full phase space.

\vskip.1in
\noi
Successively, we define a set of regularized variables (analogue to Poincar\'e's ones) on a larger domain, accordingly to the non--planarity condition.

\subsection{Angular Momentum Reduction}
Fix an orthonormal $3$--ple $({\rm k}_{\rm x},{\rm k}_{\rm y},{\rm k}_{\rm z})$ in $\real^3$. Denote by
$$C_i:=x_i\times y_i\qquad 1\leq i\leq N$$
the angular momentum of the `` body $i$'' and let
$$S_i:=\sum_{1\leq j\leq i}C_j\qquad 2\leq i\leq N$$
the sum of the first $i$ angular momenta, so that $S_N\equiv C$ coincides with the total angular momentum (\ref{total angular momentum}) of the system ($S_1$ is not defined because it coincides with $C_1$). Consider also, on the plane orthogonal to $C_i$, the $(\hat m_i,\tilde m_i)$--Keplerian motion evolving from $(y_i,x_i)$, which is defined as the solution of the differential problem
\beqa{two body problem}
\arr{
\ddot v=-\hat m_i\frac{v}{|v|^3},\ v\in \real^3\\
(\tilde m_i \dot v(0),v(0))=(y_i,x_i)
}
\eeqa
As well known,  the curve $t\to v(\hat m_i,\tilde m_i,y_i,x_i;t)$ solution of (\ref{two body problem}) draws in the space a conic section ${\cal E}_i:={\cal E}(\hat m_i,\tilde m_i,y_i,x_i)$ and we denote by $e_i:=e(\hat m_i,\tilde m_i,y_i,x_i)$ its eccentricity.

\vskip.1in
\noi
On the subset ${\cal C}_*$ of inital data $(y,x)\in(\real^{3})^N\times(\real^{3})^N$ for which
\beqa{assumptions}
\arr{
C_1\times C_2\neq 0\\
S_i\times C_{i+1}\neq 0\qquad 2\leq i\leq N-1\\
{\rm k}_{\rm z}\times C\neq 0\\
0<e_i<1\qquad 1\leq i\leq N
}
\eeqa
(in particular, each ${\cal E}_i$ is an ellipse), we define the  set of variables
\footnote{$\real_+:=(0,+\infty)\subset \real$}
\beqano
&& \Big((L,\G,\Psi),(\ell,\g,\psi)\Big)\nonumber\\
&& :=\Big((L_1,\cdots,L_N,\G_1,\cdots,\G_N,\Psi_1,\cdots,\Psi_N),(\ell_1,\cdots,\ell_N,\g_1,\cdots,\g_N,\psi_1,\cdots,\psi_N)\Big)\nonumber\\
&& \in \Big(\real_+^N\times \real_+^N\times (\real_+^{N-1}\times \real)\Big)\times (\torus^N)^3
\eeqano
as follows.
\vskip.1in
\noi
\begin{itemize}
\item[($\textrm{\scshape{d}}_1$)] For $1\leq i\leq N$,  if $a_i:=a(\hat m_i,\tilde m_i,y_i,x_i)$ is the semimajor axis of ${\cal E}_i$, then,
$$L_i:=\tilde m_i\sqrt{\hat m_ia_i}\ ;$$
\item[($\textrm{\scshape{d}}_2$)] if ${\cal A}_i:={\cal A}(\hat m_i,\tilde m_i,y_i,x_i)$ denotes the area spanned from the perihelion $P_i$ of ${\cal E}_i$ to $x_i$, then, the angle $\ell_i$ is the {\sl mean anomaly}
$$\ell_i:=2\,\frac{{\cal A}_i}{a_i^2\sqrt{1-e_i^2}}\ ;$$
\item[($\textrm{\scshape{d}}_3$)] the action $\G_i$ is 
$$\G_i:=|C_i|=L_i\sqrt{1-e_i^2}\ .$$
\item[($\textrm{\scshape{d}}_4$)] For $1\leq i\leq N-1$, the action $\Psi_i$ is
$$\Psi_i=|S_{i+1}|\ .$$
Notice that $\Psi_{N-1}={\rm G}=|C|$ is an integral of the motion.
\item[($\textrm{\scshape{d}}_5$)] The action $\Psi_N$ is
$$\Psi_N=C_{\rm z}$$
the third component of $C$. Also this variable is an integral of the motion.
\end{itemize}

\vskip.1in
\noi
Now, in order to define the conjugated angles $\g$, $\psi$, we introduce the following notations. Given $0\neq w\in \real^3$, we define the plane $\p_w$ orthogonal to $w$:
$$\p_w:=\{u\in \real^3:\ u\cdot w=0\}\ .$$
If $u$, $v$ are two non vanishing vectors in $\p_w$, we define
$${\rm k}_u:=\frac{u}{|u|}\ ,\qquad {\rm k}_w:=\frac{w}{|w|}\ ,\qquad {\rm k}_{uv}:={\rm k}_w\wedge {\rm k}_u$$
so that the triple $({\rm k}_u,{\rm k}_{uw}, {\rm k}_{w})$ is an orthonormal positively oriented basis
\footnote{I. e., the determinant of the matrix with coloumns the components of the oriented triple $({\rm k}_u,{\rm k}_{uw}, {\rm k}_{w})$  is positive (and in fact 1).}.
We then define the oriented angle seen from $w$ from $u$ to $v$, and denote it by $\a_w(u,v)$, as the angle $t+2\p\integer$, where $t$ is the unique number in $[0,2\p)$ such that
$$v=\cos{t}{\rm k}_u+\sin{t}{\rm k}_{wu}\ .$$

\vskip.1in
\noi
We can now define the angles $\g$, $\psi$. In view of assumptions (\ref{assumptions}), the following ``nodes'' are non vanishing
\beqano
n_i&:=&\arr{
C_2\times C_1\ ,\qquad i=1\\
S_i\times C_i\ ,\qquad 2\leq i\leq N 
}\nonumber\\
\bar n&:=&{\rm k}_{\rm z}\times C\ ,
\eeqano
hence, the following definitions are well put.

\vskip.1in
\noi
\begin{itemize}
\item[($\textrm{\scshape{d}}_6$)] 
For $1\leq i\leq N$, the angle $\g_i$ locates the perihelion $P_i$ of ${\cal E}_i$:
$$\g_i=\a_{C_i}(n_i,\ P_i)\ .$$
\item[($\textrm{\scshape{d}}_7$)] When $N\geq 3$, the angles $\psi_1$, $\cdots$, $\psi_{N-2}$ are
$$\psi_i=\a_{S_{i+1}}(n_{i+2},\ n_{i+1})\ ,\qquad 1\leq i\leq N-2\ .$$
\item[($\textrm{\scshape{d}}_8$)] The angle $\psi_{N-1}$ is 
$$\psi_{N-1}={\rm g}:=\a_C(\bar n,\ -n_N)\ .$$
\item[($\textrm{\scshape{d}}_9$)] The angle $\psi_N$ is the {longitude of the node} 
\footnote{The {\sl longitude of the node} of $v$ with respect to the orthonormal $3$--ple $({\rm e}_{\rm x}, {\rm e}_{\rm y}, {\rm e}_{\rm z})$ is defined as the angle $\a_{{\rm e}_{\rm z}}({\rm e}_{\rm x},{\rm e}_{\rm z}\times v)$.}
of $C$ with respect to $({\rm k}_{\rm x},{\rm k}_{\rm y},{\rm k}_{\rm z})$, namely,
$$\psi_N=\zeta:=\a_{{\rm k}_{\rm z}}({\rm k}_{\rm x},\bar n)\ .$$
Notice that this angle, together with the actions ${\rm G}$, $C_{\rm z}$, is the third component of the total angular momentum $C$.
\end{itemize}
\vskip.1in
\noi

\vskip.1in
\noi
The variables $\Big((L,\G,\Psi),(\ell,\g,\psi)\Big)$ defined via $\textrm{\scshape{d}}_1\div \textrm{\scshape{d}}_9$ will be referred as action--angle Deprit variables (or, simply, Deprit variables); the map
$$\Phi_*:\ {\cal C}_*\to \Big(\real_+^N\times \real_+^N\times (\real_+^{N-1}\times \real)\Big)\times (\torus^N)^3$$   
 which sends a point $(y,x)\in{\cal C}_*$ to the  Deprit variables  {\sl Deprit map};  their phase space  is denoted as ${\cal D}_*$. It corresponds to be the subset of $\Big((L,\G,\Psi),(\ell,\g,\psi)\Big)\in\Big(\real_+^N\times \real_+^N\times (\real_+^{N-1}\times \real)\Big)\times (\torus^N)^3$ defined by the inequalities
\beqa{D*}
\arr{
\G_i<L_i\qquad 1\leq i\leq N\\
|\G_1-\G_2|<\Psi_1<\G_1+\G_2\\
|\Psi_{i-1}-\G_{i+1}|<\Psi_{i}<\Psi_{i-1}+\G_{i+1}\qquad 2\leq i\leq N-1\\
|\Psi_N|<\Psi_{N-1}
}
\eeqa
In fact, by definition,
\beqano
\Phi_*({\cal C}_*)\subseteq{\cal D}_*
\eeqano
and we can prove
\begin{theorem}\label{invertible and symplectic}
The Deprit map $\Phi_*$ is a real--analytic symplectomorphism (symplectic diffeomorphism onto) of ${\cal C}_*$ onto ${\cal D}_*$. 
\end{theorem}
Real--analyticity follows immediately from the definition. To check injectivity and surjectivity, we shall exhibit its inverse transformation. The basis of the inversion formulae is to express the angular momenta $C_i$ ($1\leq i\leq N$), in terms of the Deprit variables : this is done in the following Lemma.
\begin{lemma}\label{reconstruction lem} The angular momenta $C_1$, $\cdots$, $C_N$  can be expressed in terms of the variables $(\G,\Psi,\psi)$ as follows. First, define $N-1$ orthonormal triples $({\rm e}_{\rm x}^i,{\rm e}^i_{\rm y},{\rm e}_{\rm z}^i)$, $2\leq i\leq N$ by letting
\begin{eqnarray}\label{Nth triple}
\arr{
{\rm e}_{\rm x}^N:=-\cos{\zeta} {\rm k}_{\rm x}-\sin{\zeta}{\rm k}_{\rm y}\\
{\rm e}_{\rm z}^N:=\sqrt{1-\left(\frac{C_{\rm z}}{{\rm G}}\right)^2}\sin{\zeta}{\rm k}_{\rm x}-\sqrt{1-\left(\frac{C_{\rm z}}{{\rm G}}\right)^2}\cos{\zeta}{\rm k}_{\rm y}+\frac{C_{\rm z}}{{\rm G}}{\rm k}_{\rm x}\\
{\rm e}_{\rm y}^N:={\rm e}_{\rm z}^N\times{\rm e}_{\rm x}^N=\frac{C_{\rm z}}{{\rm G}}\sin{\zeta}{\rm k}_{\rm x}-\frac{C_{\rm z}}{{\rm G}}\cos{\zeta}{\rm k}_{\rm y}-\sqrt{1-\left(\frac{C_{\rm z}}{{\rm G}}\right)^2}{\rm k}_{\rm z}
}
\end{eqnarray}
then (inductively), given $({\rm e}_{\rm x}^{i+1},{\rm e}^{i+1}_{\rm y},{\rm e}_{\rm z}^{i+1})$, for $3\leq i+1\leq N$, let
\beqa{invariable frames}
\arr{{\rm e}^{i}_{\rm z}:=-\frac{r_{i+1}}{\Psi_{i-1}}\,\sin{{\psi_{i}}}\,{\rm e}^{i+1}_{\rm x}+\frac{r_{i+1}}{\Psi_{i-1}}\,\cos{{\psi_{i}}}\,{\rm e}^{i+1}_{\rm y}+\frac{\tilde h_{i}}{\Psi_{i-1}}\,{\rm e}^{i+1}_{\rm z}\\
{\rm e}^{i}_{\rm x}:={\rm e}^{i+1}_{\rm x}\cos{\psi_i}+{\rm e}^{i+1}_{\rm y}\sin{\psi_i}\\
 {\rm e}^{i}_{\rm y}:={\rm e}^{i}_{\rm z}\times {\rm e}^{i}_{\rm x}=-\frac{\tilde h_i}{\Psi_{i-1}}\sin\psi_i {\rm e}_{\rm x}^{i+1}+\frac{\tilde h_i}{\Psi_{i-1}}\cos\psi_i {\rm e}_{\rm y}^{i+1}-\frac{r_{i+1}}{\Psi_{i-1}}{\rm e}_{\rm z}^{i+1}  
 }\ .
\eeqa
Then,
$$C_i=\arr{
r_i\sin\psi_{i-1}{\rm e}_{\rm x}^i-r_i\cos\psi_{i-1}{\rm e}_{\rm y}^i+h_i{\rm e}_{\rm z}^i\qquad 2\leq i\leq N\\
-r_2\sin\psi_{1}{\rm e}_{\rm x}^2+r_2\cos\psi_{1}{\rm e}_{\rm y}^2+h_1{\rm e}_{\rm z}^2\qquad i=1
}$$where, with the convention $\Psi_0:=\G_1$,
\begin{eqnarray*}
r_i&=&
\frac{\sqrt{(\Psi_{i-2}^2-(\G_{i}-\Psi_{i-1})^2)((\G_{i}+\Psi_{i-1})^2-\Psi_{i-2}^2)}}{2\Psi_{i-1}}\quad \textrm{for}\quad 2\leq i\leq N\nonumber\\
\nonumber\\
h_i&=&\left\{
\begin{array}{lrr}
\frac{\Psi_{1}^2-\G_{2}^2+\G_1^2}{2\Psi_{1}}\quad \textrm{for}\quad i=1\\
\\
\frac{\G_{i}^2+\Psi_{i-1}^2-\Psi_{i-2}^2}{2\Psi_{i-1}}\quad \textrm{for}\quad 2\leq i\leq N\quad 
\end{array}
\right.\nonumber\\
\tilde h_{i}&:=&\frac{{\Psi^2_{i}}+\Psi_{i-1}^2-\G_{i+1}^2}{2{\Psi_{i}}}\quad \textrm{for}\quad 2\leq i\leq N-1
\end{eqnarray*}
\end{lemma}
{\bf Proof.} By definition of $\zeta$, ${\rm G}$, ${C}_{\rm z}$, the components of $C$ are
$$C=\left(
\begin{array}{c}
\sqrt{{\rm G}^2-C_{\rm z}^2}\sin{\zeta}\\
-\sqrt{{\rm G}^2-C_{\rm z}^2}\cos{\zeta}\\
C_{\rm z}
\end{array}
\right)$$
Consider the orthonormal $3$--ple (\ref{Nth triple}), which has  ${\rm e}^N_{\rm z}$ in the direction of $C$, ${\rm e}^N_{\rm x}$ is in the direction of $-\bar n=-{\rm k}_{\rm z}\times C$. Then, the modulus, the third component  and the longitude of $C_N$ with respect to $({\rm e}_{\rm x}^N, {\rm e}_{\rm y}^N, {\rm e}_{\rm z}^N)$ are given, respectively, by
\begin{eqnarray}\label{elements of CN}
|C_N|&=&\G_N\nonumber\\
C_N\cdot {\rm e}_{\rm z}&=&\frac{C_N\cdot S_N}{G}\nonumber\\
&=&\frac{|C_N|^2+|S_N|^2-|S_{N}-C_N|^2}{2\Psi_{N-1}}\nonumber\\
&=&\frac{\Psi_{N-1}^2+\G_N^2-\Psi_{N-2}^2}{2\Psi_{N-1}}\nonumber\\
&=&h_{N-1}\nonumber\\
\a_{{\rm e}_{\rm z}}\left({\rm e}_{\rm x},\ {\rm e}_{\rm z}\times C_N\right)&=&\a_{S_N}\left(-\bar n,\ S_N\times C_N\right)\nonumber\\
&=&\psi_{N-1}
\end{eqnarray}
which is equivalent to
$$C_N=r_N\,\sin{\psi_{N-1}}\,{\rm e}_{\rm x}^N-r_N\,\cos{\psi_{N-1}}\,{\rm e}_{\rm y}^N+h_N\,{\rm e}_{\rm z}^N\ ,$$ 
with $$r_N:=\sqrt{\G_N^2-h_N^2}=\frac{\sqrt{(\Psi_{N-2}^2-(\G_N-\Psi_{N-1})^2)((\G_N+\Psi_{N-1})^2-\Psi_{N-2}^2)}}{2\Psi_{N-1}}\ .$$ 
Assume, now, that 
$$C_{i+1}=r_{i+1}\,\sin{{\psi_{i}}}\,{\rm e}_{\rm x}^{i+1}-r_{i+1}\,\cos{\psi_{i}}\,{\rm e}_{\rm y}^{i+1}+h_{i+1}\,{\rm e}_{\rm z}^{i+1}\qquad\textrm{for}\quad 3\leq i+1\leq N$$
Then,
\begin{eqnarray*}
S_{i}&:=&C_1+\cdots+C_i\nonumber\\
&=&S_{i+1}-C_{i+1}\nonumber\\
&=&\Psi_{i}{\rm e}^{i+1}_{\rm z}-\left(r_{i+1}\,\sin{{\psi_{i}}}\,{\rm e}^{i+1}_{\rm x}-r_{i+1}\,\cos{{\psi_{i}}}\,{\rm e}^{i+1}_{\rm y}+h_{i+1}\,{\rm e}^{i+1}_{\rm z}\right)\nonumber\\
&=&-r_{i+1}\,\sin{{\psi_{i}}}\,{\rm e}^{i+1}_{\rm x}+r_{i+1}\,\cos{{\psi_{i}}}\,{\rm e}^{i+1}_{\rm y}+\tilde h_{i}\,{\rm e}^{i+1}_{\rm z}\ ,
\end{eqnarray*}
with
$$\tilde h_{i}:={\Psi_{i}}-h_{i+1}=\frac{{\Psi^2_{i}}+\Psi_{i-1}^2-\G_{i+1}^2}{2{\Psi_{i}}}\ .$$
Let, now, $({\rm e}^{i}_{\rm x}, {\rm e}^{i}_{\rm y}, {\rm e}^{i}_{\rm z})$ the orthonormal $3$--ple with ${\rm e}^{i}_{\rm z}$ in the direction of $S_i$, ${\rm e}^{i}_{\rm x}$ in the direction of $S_{i}\times C_{i+1}$:
$$\arr{{\rm e}^{i}_{\rm z}=\frac{S_{i}}{\Psi_{i-1}}=-\frac{r_{i+1}}{\Psi_{i-1}}\,\sin{{\psi_{i}}}\,{\rm e}^{i+1}_{\rm x}+\frac{r_{i+1}}{\Psi_{i-1}}\,\cos{{\psi_{i}}}\,{\rm e}^{i+1}_{\rm y}+\frac{\tilde h_{i}}{\Psi_{i-1}}\,{\rm e}^{i+1}_{\rm z}\\
{\rm e}^{i}_{\rm x}:=\frac{S_{i}\times C_{i+1}}{|S_{i}\times C_{i+1}|}={\rm e}^{i+1}_{\rm x}\cos{\psi_i}+{\rm e}^{i+1}_{\rm y}\sin{\psi_i}\\
 {\rm e}^{i}_{\rm y}:={\rm e}^{i}_{\rm z}\times {\rm e}^{i}_{\rm x}=-\frac{\tilde h_i}{\Psi_{i-1}}\sin\psi_i {\rm e}_{\rm x}^{i+1}+\frac{\tilde h_i}{\Psi_{i-1}}\cos\psi_i {\rm e}_{\rm y}^{i+1}-\frac{r_{i+1}}{\Psi_{i-1}}{\rm e}_{\rm z}^{i+1}
 }\ .
 $$
Repeating the argument in (\ref{elements of CN}), we find that the modulus, the third component  and the longitude of the node of $C_{i}$ with respect to  $({\rm e}_{\rm x}^{i},{\rm e}_{\rm y}^{i},{\rm e}_{\rm z}^{i})$ are given by
\begin{eqnarray*}
|C_{i}|&=&\G_{i}\nonumber\\
C_{i}\cdot {\rm e}^i_z&=&\frac{C_{i}\cdot S_{i}}{\Psi_{i-1}}\nonumber\\
&=&\frac{|C_{i}|^2+|S_{i}|^2-|S_{i}-C_{i}|^2}{2\Psi_{i-1}}\nonumber\\
&=&\frac{\G_{i}^2+\Psi_{i-1}^2-\Psi_{i-2}^2}{2\Psi_{i-1}}\nonumber\\
&=&h_{i}\nonumber\\
\a_{{\rm e}_{\rm z}^{i}}\left({\rm e}_{\rm x}^{i},\ {\rm e}_{\rm z}^{i}\times C_{i}\right)&=&\a_{S_{i}}\left({n_{i+1},\ n_{i}}\right)\nonumber\\
&=&\psi_{i-1}
\end{eqnarray*}
which is equivalent to
$$C_{i}=r_{i}\,\sin{{\psi_{i-1}}}\,{\rm e}^{i}_{\rm x}-r_{i}\,\cos{\psi_{i-1}}\,{\rm e}^{i}_{\rm y}+h_{i}\,{\rm e}^{i}_{\rm z}$$
hence,
$$n_i=S_i\times C_i=\Psi_{i-1}{\rm e}^{i}_{\rm z}\times C_i=\Psi_{i-1} r_i(\cos{{\psi_{i-1}}}\,{\rm e}^{i}_{\rm x}+\sin{\psi_{i-1}}\,{\rm e}^{i}_{\rm y})$$
where
$$r_{i}=\sqrt{\G_{i}^2-h_{i}^2}=\frac{\sqrt{(\Psi_{i-2}^2-(\G_{i}-\Psi_{i-1})^2)((\G_{i}+\Psi_{i-1})^2-\Psi_{i-2}^2)}}{2\Psi_{i-1}}\ .$$ 
with the convention $\Psi_0:=\G_1$. At the $N^{th}$ step, put
\begin{eqnarray*}
C_1&=&S_{2}-C_2=-r_2\,\sin{{\psi_{1}}}\,{\rm e}^2_{\rm x}+r_2\,\cos{{\psi_{1}}}\,{\rm e}^2_{\rm y}+\tilde h_{1}\,{\rm e}^2_{\rm z}\nonumber\\
\end{eqnarray*}
with
$$h_{1}:={\Psi_{1}}-h_2=\Psi_1-\frac{\G_{2}^2+\Psi_{1}^2-\G_1^2}{2\Psi_{1}}=\frac{\G_{1}^2+\Psi_{1}^2-\G_2^2}{2\Psi_{1}}\ .$$
This completes the proof of the Lemma.
\begin{remark}\rm 
The Deprit map may be seen as an ``unfolding'' of the Jacobi reduction of the nodes, available for $N=2$. To well understand this point, we write it in {\sl spatial} Delaunay variables $(L$,{\bf G}, $\Theta,\ell$,{\bf g},$\vartheta)$, with $L=(L_1,\cdots,L_N)$, {\bf G}=$(G_1,\cdots,G_N)$, $\cdots$, which, we recall,  
are defined
as 
\beqa{Delaunay variables space}
\arr{
L_i=\tilde m_i\sqrt{\hat m_ia_i}\\
G_i=|C_i|=\sqrt{1-e_i^2}L_i\\
\Theta_i=C_i\cdot{\rm k}_{\rm z}\\
}\qquad 
\arr{
\ell_i=\frac{2{\cal A}_i}{a_i\sqrt{1-e_i^2}}\\
g_i=\a_{C_i}(\bar n_i,P_i))\\
\vartheta_i=\a_{{\rm k}_{\rm z}}({\rm k}_{\rm x},\bar n_i)
}
\eeqa
and they are well defined whenever
$$\bar n_i:={\rm k}_{\rm z}\times C_i\neq 0\ ,\quad e_i\neq 0\quad \textrm{for}\quad 1\leq i\leq N\ .$$ 
The variables $L_i$, $\ell_i$, $G_i$ are then left unchanged. To find
the expressions of the  remaining Delaunay variables in terms of the Deprit variables, we use the expressions of the angular momenta $C_i$ of Lemma \ref{reconstruction lem}, in the case $N=2$:
\beqa{ang mom}
C_1&=&\left(
\begin{array}{ccc}
r\cos{\zeta}\sin {\rm g}+ r\frac{C_{\rm z}}{{\rm G}}\sin{\zeta}\cos{{\rm g}}+ h_1\frac{\sqrt{{\rm G}^2-C_{\rm z}^2}}{{\rm G}}\sin{\zeta}\\
r\sin{\zeta}\sin {\rm g}- r\frac{C_{\rm z}}{{\rm G}}\cos{\zeta}\cos{{\rm g}}- h_1\frac{\sqrt{{\rm G}^2-C_{\rm z}^2}}{{\rm G}}\cos{\zeta}\\
\frac{C_{\rm z}}{{\rm G}} h_1-r\frac{\sqrt{{\rm G}^2-C_{\rm z}^2}}{{\rm G}}\cos{{\rm g}}
\end{array}
\right)\nonumber\\
C_2&=&\left(
\begin{array}{ccc}
-r\cos{\zeta}\sin {\rm g}- r\frac{C_{\rm z}}{{\rm G}}\sin{\zeta}\cos{{\rm g}}+ h_2\frac{\sqrt{{\rm G}^2-C_{\rm z}^2}}{{\rm G}}\sin{\zeta}\\
-r\sin{\zeta}\sin {\rm g}+r\frac{C_{\rm z}}{{\rm G}}\cos{\zeta}\cos{{\rm g}}- h_2\frac{\sqrt{{\rm G}^2-C_{\rm z}^2}}{{\rm G}}\cos{\zeta}\\
\frac{C_{\rm z}}{{\rm G}} h_2+ r\frac{\sqrt{{\rm G}^2-C_{\rm z}^2}}{{\rm G}}\cos{{\rm g}}
\end{array}
\right)\nonumber\\
\eeqa
with
$$r=\frac{\sqrt{(\G_1^2-(\G_2-{\rm G})^2)((\G_2+{\rm G})^2-\G_1^2)}}{2{\rm G}}$$
$$h_1=\frac{{\rm G}^2+\G_1^2-\G_2^2}{2{\rm G}}\ ,\quad  h_2=\frac{{\rm G}^2+\G_2^2-\G_1^2}{2{\rm G}}\ .$$
\end{remark}
This allows us to find the nodes $\bar n$, $\bar n_i$, $n_i$, and hence, the Delaunay perihelia arguments $g_i$:
$$g_i=\a_{C_i}(\bar n_i,P_i)=\a_{C_i}(n_i,P_i)+\a_{C_i}(\bar n_i,n_i)=\g_i+\a_{C_i}(\bar n_i,n_i)\qquad (i=1,2)\ .$$
and the Delaunay nodes
$$\vartheta_i=\a_{{\rm k}_{\rm z}}({\rm k}_{\rm x},\bar n_i)=\a_{{\rm k}_{\rm z}}({\rm k}_{\rm x},\bar n)+\a_{{\rm k}_{\rm z}}(\bar n,\bar n_i)=\zeta+\a_{{\rm k}_{\rm z}}(\bar n,\bar n_i)\ .$$
Finally, identifying $\Theta_1$, $\Theta_2$ with the third components of $C_1$, $C_2$ in (\ref{ang mom}), we complete the inversion formulae of the Delaunay variables in therm of the $\cD$--variables with
\beqa{Thetai}\arr{
\Theta_1=\frac{C_{\rm z}}{2}+\frac{C_{\rm z}}{2{\rm G}^2}(\G_1^2-\G_2^2)
-\frac{\sqrt{({\rm G}^2-C_{\rm z}^2)(\G_1^1-(\G_2-{\rm G})^2)((\G_2+{\rm G})^2-\G_1^1)}}
{2{\rm G}^2}\cos{\rm g}\\
\Theta_2=\frac{C_{\rm z}}{2}-\frac{C_{\rm z}}{2{\rm G}^2}(\G_1^2-\G_2^2)
+\frac{\sqrt{({\rm G}^2-C_{\rm z}^2)(\G_1^1-(\G_2-{\rm G})^2)((\G_2+{\rm G})^2-\G_1^1)}}
{2{\rm G}^2}\cos{\rm g}\\
}\eeqa
However, due to the rotation invariance, the expression of the Hamiltonian is independent on the choice of the reference frame $({\rm k}_{\rm x}, {\rm k}_{\rm y}, {\rm k}_{\rm z})$. If we choose ${\rm k}_{\rm z}$ parallel to $C$ and ${\rm k}_{\rm x}$  parallel to $C\times C_1=n_1$,  we have $\bar n_i=n_i$, hence,
$$g_i=\g_i\qquad (i=1,2)\ .$$
Also, since $n_1=-n_2$, 
$$\arr{
\vartheta_1=\a_{{\rm k}_{\rm z}}({\rm k}_{\rm x},\bar n_1)=\a_{{\rm k}_{\rm z}}(n_1,n_1)=0\\
\vartheta_2=\a_{{\rm k}_{\rm z}}({\rm k}_{\rm x},\bar n_2)=\a_{{\rm k}_{\rm z}}(n_1,n_2)=\p
}\qquad \textrm{(``opposition of the nodes'')}
$$
Finally, when the total angular momentum $C$ is seen vertical, ${\rm G}=C_{\rm z}$, hence, (\ref{Thetai}) becomes
$$
\arr{
\Theta_1=\frac{{\rm G}}{2}+\frac{\G_1^2-\G_2^2}{2{\rm G}}\\
\Theta_2=\frac{{\rm G}}{2}-\frac{\G_1^2-\G_2^2}{2{\rm G}}\\
}
$$
The previous formulae 
(completed with the identity on $L_i$, $\ell_i$) are recognized as the classical formulae for the Jacobi's reduction of the nodes.
\vskip.1in
\noi
\begin{proposition}\label{inversion}
The Deprit map $\Phi_*$ is invertible on $\cD_*$ and its inverse $\Phi_*^{-1}$ is defined as follows. Let ${\cal R}_{\rm x}$, ${\cal R}_{\rm z}$ denote the elementary rotations
$${\cal R}_{\rm x}(\a)=\left(
\begin{array}{ccc}
1&0&0\\
0&\cos{\a}&-\sin\a\\
0&\sin\a&\cos\a\\
\end{array}
\right)\ ,\quad {\cal R}_{\rm z}(\a)=\left(
\begin{array}{ccc}
\cos{\a}&-\sin\a&0\\
\sin\a&\cos\a&0\\
0&0&1
\end{array}
\right)$$
${\rm I}_{\rm z}$  the reflection
$${\rm I}_{\rm z}=\left(
\begin{array}{ccc}
-1&0&0\\
0&-1&0\\
0&0&1\\
\end{array}
\right)\ ;$$
let $i$,  $i_j$, $\tilde i_j\in (0,\p)$  be defined by \beq{ijtildeij}\cos{i}=\frac{C_{\rm z}}{{\rm G}}\ ,\quad \cos{i_j}=\frac{h_j}{\G_j}\quad 1\leq j\leq N\ ,\quad\cos{\tilde i_j}=\frac{\tilde h_j}{\Psi_{j-1}}\quad 2\leq j\leq N-1\ ,\eeq
with $h_j$, $\tilde h_j$ as in Lemma \ref{reconstruction lem}
put
\beqa{new definition}
{\rm R}_j&:=&\arr{
{\cal R}_{\rm z}(\psi_j){\cal R}_{\rm x}(\tilde i_j)\quad 2\leq j\leq N-1\\
{\cal R}_{\rm z}(\zeta){\cal R}_{\rm x}(i)\quad j=N}\ ,\nonumber\\
{\rm S}_j&:=&\arr{{\cal R}_{\rm z}(\psi_1){\cal R}_{\rm x}(i_1)\qquad j=1\\
{\cal R}_{\rm z}(\psi_{j-1}){\rm I}_{\rm z}{\cal R}_{\rm x}(i_j)\qquad 2\leq j\leq N\\
}\nonumber\\
\eeqa
and define
\beqa{cal Ri}{\cal R}_i=\arr{
{\rm R}_{N}\cdots{\rm R}_2{\rm S_1}\qquad i=1\\
{\rm R}_{N}\cdots{\rm R}_{i}{\rm S}_{i}\qquad 2\leq i\leq N
}\eeqa
Denote by
$${\rm D}_i:\ (\bar L,\bar \ell,\bar\G,\bar\g)\to (Y_i(\bar L,\bar \ell,\bar\G,\bar\g),X_i(\bar L,\bar \ell,\bar\G,\bar\g))$$
the {\sl plane} $(\hat m_i,\tilde m_i)$--{\sl Delaunay map}, defined as the four dimensional map 
$${\rm D}_i:\ ((\bar L,\bar\G),(\bar\ell,\bar\g))\in \real_+^2\times\torus^2:\quad \bar\G<\bar L\to (Y_i,X_i)\in \real^3\times\real^3$$
given by
\beqano
X_i=\frac{1}{\hat m_i}\left(\frac{\bar L}{{\tilde m_i}}\right)^2{\cal R}_{\rm z}(\g_i)\left(
\begin{array}{c}
\cos{\bar u}-\sqrt{1-\left(\frac{\bar\G}{\bar L}\right)^2}\\
\\
\frac{\bar\G}{\bar L}\sin{\bar u}\\
\\
0
\end{array}
\right)\ ,\quad 
Y_i= \frac{\tilde m_i^4\hat m_i^2}{\bar L^3}\partial_{\bar \ell}X_i
\eeqano
with $\bar u:=u(\bar L,\bar\G,\bar\ell)$ the unique solution of the {\sl Kepler's Equation}
$$u-\sqrt{1-\left(\frac{\bar\G}{\bar L}\right)^2}\,\sin{u}=\bar\ell$$
Then,
\beqa{inversion formulae}
\arr{
y_i={\cal R}_i(\G,\Psi,\psi)Y_i(L_i,\ell_i,\G_i,\g_i)\\
x_i={\cal R}_i(\G,\Psi,\psi)X_i(L_i,\ell_i,\G_i,\g_i)
}\qquad 1\leq i\leq N
\eeqa
\end{proposition}
{\bf Proof}.\ By definition, the $(\hat m_i,\tilde m_i)$--plane Delaunay map 
\beqano
{\rm D}_i:\ (L_i,\ell_i,\G_i,\g_i)&\to&(Y_i(L_i,\ell_i,\G_i,\g_i), X_i(L_i,\ell_i,\G_i,\g_i)\nonumber\\
&:=&\Big((Y_{i\rm x},Y_{i\rm y},0),(X_{i\rm x},X_{i\rm y},0)\Big)
\eeqano
gives the coordinates of $y_i$, $x_i$ on the basis of the  ``orbital triples'', \ie the (orthonormal) triples $({\rm f}_{\rm x}^i,{\rm f}_{\rm y}^i,{\rm f}_{\rm z}^i)$, where ${\rm f}_{\rm x}^i$ is in the direction of $n_i$, ${\rm f}_{\rm z}^i$ in the direction of $C_i$ (and, hence, ${\rm f}_{\rm y}^i={\rm f}_{\rm z}^i\times{\rm f}_{\rm x}^i$)
\beqa{YX}
\arr{
y_i=Y_{i\rm x}{\rm f}_{\rm x}^i+Y_{i\rm y}{\rm f}_{\rm y}^i\\
x_i=X_{i\rm x}{\rm f}_{\rm x}^i+X_{i\rm y}{\rm f}_{\rm y}^i
}\qquad 1\leq i\leq N\ .
\eeqa
Having the expressions of $C_1$, $\cdots$, $C_N$ in terms of the Deprit variables  allows to find the  orbital triples $({\rm f}_{\rm x}^i,{\rm f}_{\rm y}^i,{\rm f}_{\rm z}^i)$:
\footnote{Use $\Psi_{i-1}e_{\rm z}^i=\sum_{1\leq j\leq i}C_i=S_i$, hence, for $2\leq i\leq N$, $n_i=S_i\times C_i=\Psi_{i-1}e^i_{\rm z}\times C_i=\Psi_{i-1}r_i(\cos{\psi_{i-1}}{\rm e}_{\rm x}^i+\sin{\psi_{i-1}}{\rm e}_{\rm y}^i)$. For $i=1$, recall that $n_1=-n_2$.}
\beqa{orbital frames}
&&\arr{{\rm f}_{\rm x}^1:=\frac{n_1}{|n_1|}=-\cos{\psi_{1}}{\rm e}_{\rm x}^2-\sin{\psi_{1}}{\rm e}_{\rm y}^2\\
{\rm f}_{\rm z}^1:=\frac{C_1}{\G_1}=-\frac{r_2}{\G_1}\sin\psi_{i-1}{\rm e}_{\rm x}^2+\frac{r_2}{\G_1}\cos\psi_{1}{\rm e}_{\rm y}^2+\frac{h_1}{\G_1}{\rm e}_{\rm z}^2\\
{\rm f}_{\rm y}^1:={\rm f}_{\rm z}^1\times {\rm f}_{\rm x}^1=\frac{h_1}{\G_1}\sin{\psi_{1}}{\rm e}_{\rm x}^2-\frac{h_1}{\G_1}\cos{\psi_{i-1}}{\rm e}_{\rm y}^2+\frac{r_2}{\G_1}{\rm e}_{\rm z}^2
}\nonumber\\
\nonumber\\
&&\arr{{\rm f}_{\rm x}^i:=\frac{n_i}{|n_i|}=\cos{\psi_{i-1}}{\rm e}_{\rm x}^i+\sin{\psi_{i-1}}{\rm e}_{\rm y}^i\\
{\rm f}_{\rm z}^i:=\frac{C_i}{\G_i}=\frac{r_i}{\G_i}\sin\psi_{i-1}{\rm e}_{\rm x}^i-\frac{r_i}{\G_i}\cos\psi_{i-1}{\rm e}_{\rm y}^i+\frac{h_i}{\G_i}{\rm e}_{\rm z}^i\\
{\rm f}_{\rm y}^i:={\rm f}_{\rm z}^i\times {\rm f}_{\rm x}^i=-\frac{h_i}{\G_i}\sin{\psi_{i-1}}{\rm e}_{\rm x}^i+\frac{h_i}{\G_i}\cos{\psi_{i-1}}{\rm e}_{\rm y}^i+\frac{r_i}{\G_i}{\rm e}_{\rm z}^i
}\qquad 2\leq i\leq N\nonumber\\
\eeqa
Put, for shortness,
$${\rm F}_i:=\left(
\begin{array}{c}
{\rm f}_{\rm x}^i\\
{\rm f}_{\rm y}^i\\
{\rm f}_{\rm z}^i
\end{array}
\right)\quad 1\leq i\leq N\ ,\quad {\rm E}_i:=\left(
\begin{array}{c}
{\rm e}_{\rm x}^i\\
{\rm e}_{\rm y}^i\\
{\rm e}_{\rm z}^i
\end{array}
\right)\quad 2\leq i\leq N\ ,\quad {\rm K}:=\left(
\begin{array}{c}
{\rm k}_{\rm x}\\
{\rm k}_{\rm y}\\
{\rm k}_{\rm z}
\end{array}
\right)$$
Then, equations (\ref{Nth triple}), (\ref{invariable frames}), (\ref{orbital frames}) may be written as
\beqa{contractions}
E_N=\hat{\rm R}_N^T{\rm K}\ ,\quad E_{i}=\hat{\rm R}_i^T{\rm E}_{i+1}\quad 2\leq i\leq N-1\ ,\quad {\rm F}_i=\arr{\hat{\rm S}_1^T{\rm E}_2\quad i=1\\
\hat{\rm S}_i^T{\rm E}_i\quad 2\leq i\leq N
}
\eeqa
where:
\beqano
\hat{\rm R}_j&:=&\arr{
{\rm I}_{\rm z}{\cal R}_{\rm z}(\psi_j){\cal R}_{\rm x}(\tilde i_j){\rm I}_{\rm z}\quad 2\leq j\leq N-1\\
{\cal R}_{\rm z}(\zeta){\cal R}_{\rm x}(i){\rm I}_{\rm z}\quad j=N}\ ,\nonumber\\
\hat{\rm S}_j&:=&\arr{{\rm I}_{\rm z}{\cal R}_{\rm z}(\psi_1){\cal R}_{\rm x}(i_1)\qquad j=1\\
{\cal R}_{\rm z}(\psi_{j-1}){\cal R}_{\rm x}(i_j)\qquad 2\leq j\leq N
}
\eeqano
Then, in view of (\ref{contractions}), we can write
\beqa{change of basis}
{\rm F}_i&=&{\cal R}_i^T{\rm K}\nonumber\\
\eeqa
with
\beqa{cal Ri old}
{\cal R}_i&:=&\arr{
\hat{\rm R}_N\cdots\hat{\rm R}_2\hat{\rm S}_1\quad i=1\\
\hat{\rm R}_N\cdots\hat{\rm R}_i\hat{\rm S}_i\quad 2\leq i\leq N
}
\eeqa
Equations (\ref{YX}), (\ref{change of basis}), with ${\cal R}_i$ as in (\ref{cal Ri old}), give  the inversion formulae 
\beqano
\arr{
y_i={\cal R}_i(\G,\Psi,\psi) Y_i(L_i,\ell_i,\G_i,\g_i)\\
x_i={\cal R}_i(\G,\Psi,\psi) X_i(L_i,\ell_i,\G_i,\g_i)\\
}\qquad 1\leq i\leq N
\eeqano
and this concludes the proof (the definitions   (\ref{new definition}), (\ref{cal Ri}) of ${\cal R}_i$ are a rewrite of (\ref{cal Ri old})).

\vskip.1in
\noi
The frames  ${\rm E}_i$ ($2\leq i\leq N$), ${\rm F}_i$ ($1\leq i\leq N$) introduced in the previous proof correspond to the frames  ${F}^*_{N-i}$, ${F}_{N-i+1}$  of the binary tree of kinetic frames arising from ${F}={\rm K}$ of \cite{Dep83}.

\vskip.1in
\noi
There remains to prove symplecticity of $\Phi_*$. This is done by induction. We write then explicitly the dependence on $N$ for the Deprit map, \ie, we write the Deprit map as $\Phi_*^N:\ {\cal C}_*^N\to {\cal D}_*^N$. The basis for the induction is  $N=2$.

\vskip.1in
\noi 
\begin{lemma}\label{Phi*2 is symplectic}
The $2$--Deprit map $\Phi_*^2:\ {\cal C}_*^2\to {\cal D}_*^2$  is symplectic.
\end{lemma}
The technique for the proof of Lemma \ref{Phi*2 is symplectic} is similar to the one presented in \cite{Dep83}, apart from the introduction of the { plane} Delaunay map.

\vskip.1in
\noi 
{\bf Proof.}\ In the case $N=2$, we have 
$${\cal R}_1={\rm R}{\rm S_1}={\cal R}_{\rm z}(\zeta){\cal R}_{\rm x}(i){\cal R}_{\rm z}({\rm g}){\cal R}_{\rm x}(i_1)\ ,\quad 
{\cal R}_2={\rm R}{\rm S_2}={\cal R}_{\rm z}(\zeta){\cal R}_{\rm x}(i){\cal R}_{\rm z}({\rm g}){\rm I}_{\rm z}{\cal R}_{\rm x}(i_2)$$
where  $i_1$, $i_2$, $i\in (0,\p)$ are  defined by $$\cos{i_1}=\frac{h_1}{\G_1}\ ,\quad \cos{i_2}=\frac{h_2}{\G_2}\ ,\quad \cos{i}=\frac{C_{\rm z}}{{\rm G}}\ .$$
Differentiating $x_i$ in (\ref{inversion formulae}), we find
$$dx_i={\cal R}_i(\G,\Psi,\psi) dX(L_i,\ell_i,\G_i,\g_i)+(d{\cal R}_i(\G,\Psi,\psi)) X(L_i,\ell_i,\G_i,\g_i)\qquad i=1,\ 2\ .$$
So, since ${\cal R}_i$ is unitary and the plane Delaunay map ${\rm D}_i$ is symplectic,
\beqano
y_i\cdot dx_i&=&Y_i(L_i,\ell_i,\G_i,\g_i)\cdot dX_i(L_i,\ell_i,\G_i,\g_i)+y_i\cdot (d{\cal R}_i(\G,\Psi,\psi)) {\cal R}_i(\G,\Psi,\psi)^Tx_i\nonumber\\
&=& L_id\ell_i+\G_id\g_i+y_i\cdot (d{\cal R}_i(\G,\Psi,\psi)) {\cal R}_i(\G,\Psi,\psi)^Tx_i
\eeqano
Thus, summing over $i=1$, $2$, 
\beqa{Cartan 1 form}
y\cdot dx&=& L\cdot d\ell+\G\cdot d\g\nonumber\\
&+&y_1\cdot (d{\cal R}_1(\G,\Psi,\psi)) {\cal R}_1(\G,\Psi,\psi)^Tx_1+y_2\cdot (d{\cal R}_2(\G,\Psi,\psi)) {\cal R}_2(\G,\Psi,\psi)^Tx_2\nonumber\\
\eeqa
Differantiating ${\cal R}_1={\cal R}_{\rm z}(\zeta){\cal R}_{\rm x}(i){\cal R}_{\rm z}({\rm g}){\cal R}_{\rm x}(i_1)$ and using, as well known,
$$(d{\cal R}_{\rm x}(\a)){\cal R}_{\rm x}(\a)^Tq={\rm k}_{\rm x}\times q\,d\a\ ,\quad (d{\cal R}_{\rm z}(\a)){\cal R}_{\rm z}(\a)^Tq={\rm k}_{\rm z}\times q\, d\a$$
we find
\beqano
(d{\cal R}_1(\G,\Psi,\psi)) {\cal R}_1(\G,\Psi,\psi)^Tx_1&=&{\rm k}_{\rm z}\times x_1\,d\zeta+{\rm e}_{\rm z}\times x_1\,d{\rm g}\nonumber\\
&-&{\rm e}_{\rm x}\times x_1\,di+{\rm f}_{\rm x}^1\times x_1 di_1
\eeqano
which gives, taking the scalar product with $y_1$,
\beqa{y1}
y_1\cdot (d{\cal R}_1(\G,\Psi,\psi)) {\cal R}_1(\G,\Psi,\psi)^Tx_1&=&C_1\cdot{\rm k}_{\rm z}d\zeta+ C_1\cdot {\rm e}_{\rm z}d{\rm g}\nonumber\\
&-& C_1\cdot {\rm e}_{\rm x}\,di+C_1\cdot{\rm f}_{\rm x}^1\,di_1\nonumber\\
&=&C_1\cdot{\rm k}_{\rm z}d\zeta+ C_1\cdot {\rm e}_{\rm z}d{\rm g}\nonumber\\
&-& C_1\cdot {\rm e}_{\rm x}\,di
\eeqa
since $C_1\cdot{\rm f}_{\rm x}^1=0$. Similarly, differentiating ${\cal R}_2={\cal R}_{\rm z}(\zeta){\cal R}_{\rm x}(i){\cal R}_{\rm z}({\rm g}){\rm I}_{\rm z}{\cal R}_{\rm x}(i_2)$, we find
\beqa{y2}
y_2\cdot (d{\cal R}_2(\G,\Psi,\psi)) {\cal R}_2(\G,\Psi,\psi)^Tx_1
&=&C_2\cdot{\rm k}_{\rm z}d\zeta+ C_2\cdot {\rm e}_{\rm z}d{\rm g}\nonumber\\
&-& C_2\cdot {\rm e}_{\rm x}\,di
\eeqa
The  sum of (\ref{y1}) and (\ref{y2}) gives then
\beqano
&& y_1\cdot (d{\cal R}_1(\G,\Psi,\psi)) {\cal R}_1(\G,\Psi,\psi)^Tx_1+y_2\cdot (d{\cal R}_2(\G,\Psi,\psi)) {\cal R}_2(\G,\Psi,\psi)^Tx_2\nonumber\\
&=& C\cdot{\rm k}_{\rm z}d\zeta+C\cdot{\rm e}_{\rm z}d{\rm g}-C\cdot{\rm e}_{\rm x}di\nonumber\\
&=&C_{\rm z}d\zeta+{\rm G}d{\rm g}
\eeqano
since
$$C=C_1+C_2={\rm G}{\rm e}_{\rm z}$$
(which also implies $C\cdot{\rm e}_{\rm x}=0$).
The proof is complete, in view of (\ref{Cartan 1 form}).

\vskip.1in
\noi
It remains to prove the inductive step. 
\begin{lemma}
Assume that the $N$--Deprit map $\Phi_*^N:\ {\cal C}_*^N\to {\cal D}_*^N$, is symplectic for a given $N\geq 2$. Then, the $(N+1)$--
Deprit map $\Phi_*^{N+1}:\ {\cal C}_*^{N+1}\to {\cal D}_*^{N+1}$ is symplectic.
\end{lemma}
{\bf Proof.}\ Without loss of generality, we shall restrict to the subset $\hat{\cal C}_*^{N+1}$ of ${\cal C}_*^{N+1}$ where also
\beq{good Delaunay}
\hat n_i:={\rm k}_{\rm z}\times S_i\neq 0\quad 1\leq i\leq N ,\quad  \bar n_i:={\rm k}_{\rm z}\times C_i\neq 0\quad 1\leq i\leq N+1
\eeq
and then we will recover the result by continuity. 
Under the assumption (\ref{good Delaunay}), we can view the Deprit map $\Phi_*^{N+1}$ in {\sl Delaunay variables}, namely, we shall write $\Phi_*^{N+1}=\hat\Phi_*^{N+1}\circ\Phi_D^{N+1}$, where 
$(L,{\textrm{\bf G}},\Theta,\ell,{\textrm{\bf g}},\vartheta)=\Phi_D^N(y,x)$ is the map  (\ref{Delaunay variables space}) which defines the Delaunay variables. 
Let $\hat{\cal D}^{N+1}_*:=\Phi_D^{N+1}(\hat{\cal C}^{N+1}_*)$. Then, $\hat\Phi_*^2$ is symplectic on $\hat{\cal D}^{2}_*$, by Lemma \ref{Phi*2 is symplectic};  $\hat\Phi_*^N$ is symplectic on $\hat{\cal D}^{N}_*$ by assumption. We equivalently prove that  $\hat\Phi_*^{N+1}$ is symplectic on $\hat{\cal D}^{N+1}_*$, which will conclude (since $\Phi_D^n$ is symplectic on $\hat{\cal C}_*^n$, for any $n$). Neglecting the variables  $(L,\ell)$ (on which $\hat\Phi_*^{N+1}$ acts as the identity), the map $\hat\Phi_*^{N+1}$ is described by equations 
\beqano
\arr{
\G_i=G_i\qquad 1\leq i\leq N+1\\
\Psi_i=|\sum_{1\leq j\leq i+1}C_j|\qquad 1\leq i\leq N-1\\
{\rm G}=|S_{N+1}|\\
C_{\rm z}=\sum_{1\leq j\leq N+1}\Theta_{j}\\
\g_i=g_i+\a_{C_i}(n_i,\bar n_i)\qquad 1\leq i\leq N+1\\
\psi_i=\a_{S_{i+1}}(n_{i+2},n_{i+1})\qquad 1\leq i\leq N-1\\
{\rm g}=\a_{S_{N+1}}(\bar n,-n_{N+1})\\
\zeta=\a_{{\rm k}_{\rm z}}({\rm k}_{\rm x},\bar n)
}
\eeqano
where the $C_i$'s, hence, $S_i=C_1+\cdots+C_{i+1}$ and the nodes $\bar n_i={\rm k}_{\rm z}\times C_i$, $n_i=S_i\times C_i$
, $\bar n={\rm k}_{\rm z}\times S_{N+1}$, are thought as functions of the Delaunay variables
\footnote{The angular momenta $C_i$'s, in terms of the Delaunay variables are
$$C_i=\left(\begin{array}{c}
\sqrt{G_i^2-\Theta_i^2}\sin{\vartheta_i}\\
-\sqrt{G_i^2-\Theta_i^2}\cos{\vartheta_i}\\
\Theta_i
\end{array}\right)\ .$$
This follows from the definitions of $(G_i,\Theta_i,\vartheta_i)$.
}
.\\
Let us introduce the following notations
\beqano
&z_i:=(L_i,G_i,\Theta_i,\ell_i,g_i,\vartheta_i)\ ,\quad &Z_i:=(L_i,\G_i,\Psi_i,\ell_i,\g_i,\psi_i)\nonumber\\
&\arr{z:=(z_1,\cdots,z_{N+1})\\
\hat z_i:=(z_1,\cdots,z_N)}\qquad &\arr{Z:=(Z_1,\cdots,Z_{N+1})\\
\hat Z_i:=(Z_1,\cdots,Z_N)}
\eeqano
Now, if $z\in\hat{\cal D}^{N+1}_*$, then, the point $\hat z$ lies in the domain of definition $\hat{\cal D}^{N}_*$ of $\hat\Phi_*^N$ and we can set
\footnote{In (\ref{intermediate}), let $\hat Z'=(Z'_1,\cdots,Z'_N)$, with $Z'_i=(L'_i,\G'_i,\Psi'_{i},\ell'_i,\g_i',\psi_i')$, for $1\leq i\leq N$.}
\beq{intermediate}\tilde\Phi_*^{N+1}(z)=\tilde\Phi_*^{N+1}(\hat z,z_N)=\hat\Phi_*^{N}(\hat z,z_{N+1})=:Z'=(\hat Z',z_{N+1})\eeq
\ie, $\tilde\Phi_*^{N+1}$ acts as $\hat\Phi_*^N$ on $\hat z$, while on the last block $z_{N+1}$ of the Delaunay variables acts as the identity. $\tilde\Phi_*^{N+1}$ is thus symplectic since  $\hat\Phi_*^N$ is, as already outlined.
Now, leaving the remaining variables unchanged, we apply  $\hat\Phi_*^2$ to the two blocks consisting, the former, to the block of variables $$z_N':=(L'_{N},{\rm G}'=\Psi'_{N-1},C_{\rm z}'=\Psi'_N,\ell_{N+1},{\rm g}'=\psi'_{N-1}, \zeta'=\psi'_N)$$  and, the latter, to the block of variables $z_{N+1}$  left unvaried by $\tilde\Phi_*^{N+1}$. We define, then,
$$\tilde\Phi_*^2:\ Z'\to Z$$ as follows:
$$
\arr{
Z_i=Z_i'\qquad 1\leq i\leq N-2\\
(L_{N-1},\G_{N-1},\G_N,\ell_{N-1},\g_{N-1},\g_N)=(L'_{N-1},\G'_{N-1},\G_N',\ell'_{N-1},\g'_{N-1},\g'_N)\\
((L_{N},\Psi_{N-1},\Psi_N,\ell_{N+1},\psi_{N-1}, \psi_N),(L_{N+1},\G_{N+1},\Psi_{N+1},\ell_{N+1},\g_{N+1},\Psi_{N+1}))=\hat\Psi_*^2(z'_N,z_{N+1})
}
$$
Also $\tilde\Phi_*^2$ is symplectic, because it is obtained lifting $\hat\Phi_*^2$ with the identity map, and, therefore, so is the composition \beq{composition}\tilde\Phi_*^2\circ \tilde\Phi_*^{N+1}\ .\eeq
The claim, now, follows upon recognizing that (\ref{composition}) reconstructs $\hat\Phi_*^{N+1}$:
\beq{check}\tilde\Phi_*^2\circ \tilde\Phi_*^{N+1}=\hat\Phi_*^{N+1}\ .\eeq
The key point while checking  (\ref{check}) is
\beqano
\psi_{N-1}|_{{\tilde\Phi_*^2\circ \tilde\Phi_*^{N+1}}}&=&
{\rm g}'+\a_{S_{N}}(S_{N+1}\times S_N,{\rm k}_{\rm z}\times S_N)\nonumber\\
&=&\a_{S_{N}}({\rm k}_{\rm z}\times S_N,-S_{N}\times C_N)+\a_{S_{N}}(S_{N+1}\times S_N,{\rm k}_{\rm z}\times S_N)\nonumber\\
&=&\a_{S_{N}}(S_{N+1}\times S_N,-S_{N}\times C_N)\nonumber\\
&=&\a_{S_{N}}(-S_{N+1}\times C_{N+1},-S_{N}\times C_N)\nonumber\\
&=&\a_{S_{N}}(S_{N+1}\times C_{N+1},S_{N}\times C_N)\nonumber\\
&=&\psi_{N-1}|_{\hat\Phi_*^{N+1}}\ ,
\eeqano
since, by definition, ${\rm g}'=\a_{S_{N}}({\rm k}_{\rm z}\times S_N,-S_{N}\times C_N)$ and $S_N=S_{N+1}-C_{N+1}$.

\subsection{Regularization}\label{Regularization}
The action--angle Deprit variables  discussed in the previous section become singular when some of (\ref{assumptions}) do not hold. In this paragraph, we discuss
 a Poincar\'e regularization. For $N=2$, we put
\beqano
&&\arr{G:={\rm G}\\
g:={\rm g}+\zeta\\
P:=\sqrt{2({\rm G}-C_{\rm z})}\cos{\zeta}\\
Q:=-\sqrt{2({\rm G}-C_{\rm z})}\sin{\zeta}\\
}\qquad \arr{\L_i=L_i\\
\l_i=\ell_i+\g_i}\nonumber\\
&&\arr{\eta_i=\sqrt{2(L_i-\G_i)}\cos{\g_i}\\
\xi_i=-\sqrt{2(L_i-\G_i)}\sin{\g_i}
\\}
\eeqano 
and we recover the ``unfolding'' of the Jacobi regularized coordinates. So, we discuss in detail only the case $N\geq 3$.

\vskip.1in
\noi 
Let then $N\geq 3$ and let ${\cal C}\supset{\cal C}_*$ the set of $(y,x)\in\real^{3N}\times\real^{3N}$ where
\beq{dot prod}
\arr{
e_i<1\\
\frac{C}{|C|}\cdot{\rm k}_{\rm z}\neq -1\\
\frac{C_1}{|C_1|}\cdot\frac{C_2}{|C_2|}\neq -1\\
\frac{C_i}{|C_i|}\cdot\frac{S_i}{|S_i|}\neq -1\qquad 2\leq i\leq N-1\\
\frac{C}{|C|}\cdot \frac{C_N}{|C_N|}\neq \pm 1\\
}\ .\eeq
\ie, the eccentricities are allowed to go to $0$; $C_1$ is allowed to go parallel to $C_2$; $C_i$ are allowed to go parallel to $S_i$ for $2\leq i\leq N-1$; $C$ is allowed to go parallel to ${\rm k}_{\rm z}$.

\vskip.1in
\noi
Other regularizations than (\ref{dot prod}) (relatively to different choices for the signs of the dot products in (\ref{dot prod})) might be discussed.

\vskip.1in
\noi
In order to regularize zero eccentricities and the first $N-1$ mutual inclinations, \ie, in order to define a new set of variables in a region of the phase space where
$$e_i=0\quad \textrm{for some}\quad 1\leq i\leq N\ ,\quad \textrm{or}\quad S_i\parallel C_{i}\quad \textrm{for some}\quad 2\leq i\leq N-1\ ,$$
we assume that the variables $\L$, ${\rm G}$ satisfy
\begin{eqnarray}\label{selfconsistency}
\left|\sum_{1\leq i\leq N-1}\L_i-\L_N\right|<{\rm G}<\sum_{1\leq i\leq N}\L_i\ .
\end{eqnarray}
This guarantees that the configuations of the phase space corresponding to (simultanously) zero eccentricities and first $N-1$ mutual inclinations might be reached by the system, being inner points of the phase space.
%
%
%
%
%
%
\begin{theorem}
When $\L$, ${\rm G}$ also satisfy (\ref{selfconsistency}), define the real--analytic symplectomorphism
\beqano
\Phi_{\rm r}:\ \Big((L,\G,\Psi),(\ell,\g,\psi)\Big)\in{\cal D_*}&\to& \Big((\L,\l),(\eta,\xi),(p,q),(G,g),(P,Q)\Big)\nonumber\\
&=&\Big((\L_1,\cdots,\L_N),(\l_1,\cdots,\l_N),(\eta_1,\cdots,\eta_N,\xi_1,\cdots,\xi_N)\nonumber\\  
&& (p_1,\cdots,p_{N-2},q_1,\cdots,q_{N-2}),(G,g),(P,Q)\Big)\nonumber\\
&&\in (\real_+^N\times\torus^N)\times (\real^N\times\real^N)\times(\real^{N-2}\times\real^{N-2})\times\nonumber\\
&&\times(\real_+\times\torus)\times(\real\times\real)                      
\eeqano
as follows. Let
\beqano
H_i&:=&L_i-\G_i\qquad 1\leq i\leq N\quad K_i:=\arr{\G_1-\Psi_1+\G_2\quad i=1\\
\Psi_{i-1}-\Psi_i+\G_{i+1}\quad 2\leq i\leq N-2}\nonumber\\
\k_i&:=&\sum_{i\leq j\leq N-2}\psi_j\quad 1\leq i\leq N-2\quad 
\hat \k_i:=\arr{\k_1\quad i=1\\
\k_{i-1}\quad 2\leq i\leq N-1\\
0\quad i=N}\nonumber\\
h_i&:=&\g_i+\hat\k_i
\eeqano
and put
\beqano
&&\arr{
\L_i:=L_i\\
\l_i:=\ell_i+h_i\\
\eta_i:=\sqrt{2H_i}\cos{h_i}\\
\xi_i:=-\sqrt{2H_i}\sin h_i\\
}1\leq i\leq N\ ,\quad\arr{
p_i:=\sqrt{2K_i}\cos{\k_i}\\
q_i:=-\sqrt{2K_i}\sin\k_i}1\leq i\leq N-2\nonumber\\
&&\arr{G:={\rm G}\\
g:={\rm g}+\zeta\\
P:=\sqrt{2({\rm G}-C_{\rm z})}\cos{\zeta}\\
Q:=-\sqrt{2({\rm G}-C_{\rm z})}\sin{\zeta}\\
}
\eeqano
Then, the map $\Phi_{\rm BD}:=\Phi_{\rm r}\circ\Phi_*$ extends to a real--analytic symplectomorphism on ${\cal C}$.
\end{theorem}
The variables $\Big((\L,\l)$ $(\eta,\xi)$, $(p,q)$, $(G,g)$, $(P,Q)\Big)$ will be referred as {\sl regularized} Deprit variables. Observe that, now, the role of cyclic variables for (\ref{spatial problem}) is played by $(P,Q,g)$.
\begin{remark}\rm The inverse $\phi_{\rm r}:=\Phi_{\rm r}^{-1}$  on 
${\cal D}_{\rm r*}:=\Phi_{\rm r}({\cal D}_*)$ is given by
\beqa{regularization}
&&\arr{L_i=\L_i\\
\ell_i=\l_i-h_i}\ ,\quad \arr{\G_i=\L_i-\frac{\eta_i^2+\xi_i^2}{2}\\
\g_i=h_i-\hat\k_i}\nonumber\\
&&\Psi_i=\arr{\sum_{1\leq j\leq i+1}\left(\L_j-\frac{\eta_j^2+\xi_j^2}{2}\right)-\sum_{1\leq j\leq i}\frac{p_j^2+q_j^2}{2}\quad 1\leq i\leq N-2\\
G\quad i=N-1\\
G-\frac{P^2+Q^2}{2}\quad i=N
}\nonumber\\
&&\psi_i=\arr{
\k_i-\k_{i+1}\quad 1\leq i\leq N-3\\
\k_{N-2}\quad i=N-2\\
g-\zeta(P,Q)\quad i=N-1\\
\zeta(P,Q)\quad i=N
}
\eeqa
where
$$h_i=\arg{(\eta_i,-\xi_i)}\ (1\leq i\leq N) ,\quad k_i=\arg(p_i,-q_i)\ (1\leq i\leq N-2)\ ,\quad \zeta(P,Q)=\arg(P,-Q)$$
(the previous expressions are well put on  ${\cal D}_{\rm r*}$). 
\end{remark}
\begin{remark}\rm The domain ${\cal D}_{\rm r*}$ is the set of $(\L,\l,\eta,\xi,p,q,G,g,P,Q)$  where $\L\in \real_+^N$, $(\l,g)\in \torus^N\times \torus$ and the {\sl functions} $\G$, $\Psi$ as in (\ref{regularization}) verify (\ref{D*}).
\end{remark}
\vskip.1in
\noi
{\bf Proof.}\ 
Put $${\cal D}_0:=\Big\{(P,Q)=0\Big\}\bigcup_{1\leq i\leq N}\Big\{(\eta_i,\xi_i)=0\Big\}\bigcup_{1\leq i\leq N-2}\Big\{(p_i,q_i)=0\Big\}$$ and observe that
$${\cal D}_{\rm r*}\cap{\cal D}_0=\emptyset$$
We prove that the map
\beqano
\Big((\L,\l),(\eta,\xi),(p,q),(G,g),(P,Q)\Big)\in{\cal D}_{\rm r*}&\to& (y,x)=\Big((y_1,\cdots,y_N),(x_1,\cdots,x_N)\Big)\nonumber\\
&=&\phi_*\circ\phi_{\rm r}\Big((\L,\l),(\eta,\xi),(p,q),(G,g),(P,Q)\Big)\nonumber\\
&:=&\phi_{\rm BD}\Big((\L,\l),(\eta,\xi),(p,q),(G,g),(P,Q)\Big)\nonumber\\
\eeqano
where $\phi_*:=\Phi_*^{-1}$, can be bijectively and analytically (hence, symplectically) extended to the domain
$${\cal D}_{\rm r}:={\cal D}_{\rm r*}\bigcup{\cal D}_0\ .$$
By Proposition \ref{inversion}, $\phi_{\rm fr}$ is defined by
\beqa{inver Phi*}
\arr{
x_i=\Big({\cal R}_i\circ{\phi_{\rm r}}\Big)X_i\circ{\phi_{\rm r}}\\
y_i=\Big({\cal R}_i\circ{\phi_{\rm r}}\Big)Y_i\circ{\phi_{\rm r}}
}\qquad 1\leq i\leq N
\eeqa
We explicitate the matrices ${\cal R}_i$ defined in Lemma \ref{inversion} as
\footnote{Clearly, in the second and third line of (\ref{rewrite Ri}), the productories do not appear when $N=3$.}
\beqa{rewrite Ri}
{\cal R}_i&=&{\cal R}_{\rm z}(\zeta){\cal R}_{\rm x}(i){\cal R}_{\rm z}({\rm g})\nonumber\\
&\times&\arr{
{\rm I}_{\rm z}{\cal R}_{\rm x}(i_N)\quad i=N\\
{\cal R}_{\rm x}(\tilde i_{N-1})\left(\prod_{j=2}^{N-i}{\cal R}_{z}(\psi_{N-j}){\cal R}_{x}(\tilde i_{N-j})\right){\cal R}_{\rm z}(\psi_{i-1}){\rm I}_{\rm z}{\cal R}_{\rm x}(i_i)\quad 2\leq i\leq N-1\\
{\cal R}_{\rm x}(\tilde i_{N-1})\left(\prod_{j=2}^{N-2}{\cal R}_{z}(\psi_{N-j}){\cal R}_{x}(\tilde i_{N-j})\right){\cal R}_{\rm z}(\psi_{1}){\cal R}_{\rm x}(i_1)\quad  i=1\\
}\nonumber\\
\eeqa
and we think (without changing their names) the ``inclinations'' $i$, $i_j$, $\tilde i_j$ expressed in regularized Deprit variables  ($i$ is a function of $P$, $Q$, $G$; $i_i$, $\tilde i_i$ are functions of $\L$, $(\eta,\xi)$, $(p,q)$ and $G$), then, using the expressions for $\psi_i$, with $1\leq i\leq N-2$, given (\ref{regularization}), we can rewrite ${\cal R}_i$ in terms of the regularized Deprit variables  as
\beqa{inver Phi}
{\cal R}_i\circ{\phi_{\rm r}}&=&{\eufm R}_0\times\arr{
{\rm I}_{\rm z}{\cal R}_{\rm x}(i_N)\quad i=N\\
{\cal R}_{\rm x}(\tilde i_{N-1})\left(\prod_{j=2}^{N-i}\tilde{\cal S}_{N-j}\right){\rm I}_{\rm z}{\cal S}_i{\cal R}_{\rm z}(\k_{i-1})\quad 2\leq i\leq N-1\\
{\cal R}_{\rm x}(\tilde i_{N-1})\left(\prod_{j=2}^{N-2}\tilde {\cal S}_{N-j}\right){\cal S}_1{\cal R}_{\rm z}(\k_1)\quad  i=1\\
}\nonumber\\
&=:&{\eufm R}_0{\eufm R}_i{\cal R}_{\rm z}(\hat \k_i)
\eeqa
where
\beqa{matrices1}
{\eufm R}_0&:=&{\cal R}_{\rm z}(\zeta(P,Q)){\cal R}_{\rm x}(i){\cal R}_{\rm z}(-\zeta(P,Q)+g)
\eeqa
and
\beqa{matrices}
\tilde {\cal S}_j&:=&{\cal R}_{\rm z}(\k_j){\cal R}_{\rm x}(\tilde i_j){\cal R}_{\rm z}(-\k_j)\quad 2\leq j\leq N-2\nonumber\\
{\cal S}_j&:=&\arr{{\cal R}_{\rm z}(\k_{1}){\cal R}_{\rm x}( i_1){\cal R}_{\rm z}(-\k_{1})\quad j=1\\
{\cal R}_{\rm z}(\k_{j-1}){\cal R}_{\rm x}( i_j){\cal R}_{\rm z}(-\k_{j-1})\quad 2\leq j\leq N-1}\nonumber\\
\eeqa
\begin{lemma}\label{Si are regular}
With the convention $(p_0,q_0):=(p_1,q_1)$, $\Psi_0:=\G_1$, on ${\cal D}_{\rm r*}$, the matrices ${\cal S}_j$'s, $\tilde{\cal S}_j$'s  have the following expressions:
\beqa{calSj}
{\cal S}_j&=&\left(
\begin{array}{ccc}
1-q_{j-1}^2{\eufm c}_j&-p_{j-1}q_{j-1}{\eufm c}_j&-q_{j-1}{\eufm s}_j\\
-p_{j-1}q_{j-1}{\eufm c}_j&1-p_{j-1}^2{\eufm c}_j&-p_{j-1}{\eufm s}_j\\
q_{j-1}{\eufm s}_j&p_{j-1}{\eufm s}_j&1-(p_{j-1}^2+q_{j-1}^2){\eufm c}_j
\end{array}
\right)\qquad (1\leq j\leq N-1)\nonumber\\
\tilde{\cal S}_j&=&\left(
\begin{array}{ccc}
1-q_{j}^2\tilde{\eufm c}_j&-p_{j}q_{j}\tilde{\eufm c}_j&-q_{j}\tilde{\eufm s}_j\\
-p_{j}q_{j}\tilde{\eufm c}_j&1-p_{j}^2\tilde{\eufm c}_j&-p_{j}\tilde{\eufm s}_j\\
q_{j}\tilde{\eufm s}_j&p_{j}\tilde{\eufm s}_j&1-(p_{j}^2+q_{j}^2)\tilde{\eufm c}_j
\end{array}
\right)\qquad (2\leq j\leq N-2)\nonumber\\
\eeqa
where
\beqa{eufm ci si}
&&\arr{{\eufm c}_1:=\frac{1-\cos{i_1}}{p_{1}^2+q_{1}^2}=\frac{\G_2-\G_1+\Psi_1}{4\G_1\Psi_1}\\
{\eufm s}_1:=\frac{\sin{i_1}}{\sqrt{p_{1}^2+q_{1}^2}}=\sqrt{{\eufm c}_1(2-(p_1^2+q_1^2){\eufm c}_1)}}\nonumber\\
\nonumber\\
&&\arr{
{\eufm c}_j:=\frac{1-\cos{i_j}}{p_{j-1}^2+q_{j-1}^2}=\frac{\Psi_{j-2}-\G_j+\Psi_{j-1}}{4\G_j\Psi_{j-1}}\\
{\eufm s}_j:=\frac{\sin{i_j}}{\sqrt{p_{j-1}^2+q_{j-1}^2}}=\sqrt{{\eufm c}_j(2-(p_{j-1}^2+q_{j-1}^2){\eufm c}_j)}}\quad (2\leq j\leq N-1)\nonumber\\
\nonumber\\
&&\arr{
\tilde{\eufm c}_j:=\frac{1-\cos\tilde{i_j}}{p_{j}^2+q_{j}^2}=\frac{\G_{j+1}+\Psi_j-\Psi_{j-1}}{4\Psi_j\Psi_{j-1}}\\
\tilde{\eufm s}_j:=\frac{\sin\tilde{i_j}}{\sqrt{p_{j}^2+q_{j}^2}}=\sqrt{\tilde{\eufm c}_j(2-(p_j^2+q_j^2)\tilde{\eufm c}_j)}}\quad (2\leq j\leq N-2)\nonumber\\
\eeqa
where $\G_i$, $\Psi_j$ are thought as functions of $(\L,\eta,\xi,p,q)$:
\beqa{GiPsii}\arr{\G_i=\L_i-\frac{\eta_i^2+\xi_i^2}{2}\qquad 1\leq i\leq N\\
\Psi_i=\sum_{1\leq i\leq i+1}\left(\L_j-\frac{\eta_j^2+\xi_j^2}{2}\right)-\sum_{1\leq j\leq i}\frac{p_j^2+q_j^2}{2}\qquad 1\leq i\leq N-2
}\eeqa
\end{lemma}
{\bf Proof.}\ Let us prove, for istance, the expression for ${\cal S}_1$, since the other ones are similar.
We have
\beqano
{\cal S}_1&=&{\cal R}_{\rm z}(\k_1){\cal R}_{\rm x}(i_1){\cal R}_{\rm z}(-k_1)\nonumber\\
&=&\left(
\begin{array}{ccc}
1-\sin^2k_1(1-\cos{i_1})&\sin{k_1}\cos{k_1}(1-\cos{i_1})&\sin{k_1}\sin{i_1}\\
\sin{k_1}\cos{k_1}(1-\cos{i_1})&1-\cos^2k_1(1-\cos{i_1})&-\cos{k_1}\sin{i_1}\\
-\sin{k_1}\sin{i_1}&\cos{k_1}\sin{i_1}&\cos{i_1}
\end{array}
\right)
\eeqano
Now, letting $\k_1=\arg{(p_1,-q_1)}$, hence,
$$
\arr{
\cos{\k_1}=\frac{p_1}{\sqrt{p_1^2+q_1^2}}\\
\sin{\k_1}=-\frac{q_1}{\sqrt{p_1^2+q_1^2}}
}
$$
and recognizing that
$${\eufm c}_1=\frac{\G_2-\G_1+\Psi_1}{4\G_1\Psi_1}=\frac{1-\cos{i_1}}{p_{1}^2+q_{1}^2}\ ,\quad {\eufm s}_1=\sqrt{{\eufm c}_1(2-(p_1^2+q_1^2){\eufm c}_1)}=\frac{\sin{i_1}}{\sqrt{p_{1}^2+q_{1}^2}}$$
the claim follows. Recall, in fact, the definition of $i_1$:
$${\cos{i_1}}:=\frac{\Psi_1^2+\G_1^2-\G_2^2}{2\Psi_1\G_1}\qquad i_1\in (0,\p)$$
which gives
\beqano
1-\cos{i_1}&=&\frac{\G_2^2-(\Psi_1-\G_1)^2}{2\Psi_1\G_1}\nonumber\\
&=&\frac{(\G_1+\G_2-\Psi_1)(\G_2-\G_1+\Psi_1)}{2\Psi_1\G_1}\nonumber\\
&=&\frac{\G_2-\G_1+\Psi_1}{4\Psi_1\G_1}(p_1^2+q_1^2)\nonumber\\
&=&{\eufm c}_1(p_1^2+q_1^2)
\eeqano
since
$$\G_1+\G_2-\Psi_1=\frac{p_1^2+q_1^2}{2}\ ,$$
hence (as $i_1\in (0,\p)$),
$$\sin{i_1}=\sqrt{1-\cos^2{i_1}}=\sqrt{1-\Big(1-{\eufm c}_1(p_1^2+q_1^2)\Big)^2}=\sqrt{{\eufm c}_1(2-(p_1^2+q_1^2){\eufm c}_1)}\sqrt{p_1^2+q_1^2}={\eufm s}_1\sqrt{p_1^2+q_1^2}$$
This completes the proof.

\vskip.1in
\noi
By the previous Lemma, the matrices $\tilde{\cal S}_j$, ${\cal S}_j$  can be analytically  extended to ${\cal D}_{\rm r}$. It is not difficult to prove that the same holds for the matrix ${\eufm R}_0$, hence, for ${\eufm R}_j$.
The transformation  $\Phi^{-1}$ is thus regularized on ${\cal D}_{\rm r}$, being given by
\beqano
\arr{y_i={\eufm R}_0{\eufm R}_i\hat y_i\\
x_i={\eufm R}_0{\eufm R}_i\hat x_i}
\eeqano
where $(\L_i,\l_i,\eta_i,\xi_i)\to (\hat y_i,\hat x_i$) is the embedding
\footnote{This means taking the third coordinate of $(\hat y_i,\hat x_i)$ equal to zero.}
in $\real^{3}\times \real^{3}$ of the $(\hat m_i,\tilde m_i)$--Plane Delaunay--Poincar\'e Map.
In fact, using the expressions for $L$, $\G$, $\ell$, $\g$ in (\ref{regularization}),
we find that 
$$(Y_i\circ{\phi_{\rm r}},X_i\circ{\phi_{\rm r}})=\Big({\cal R}_{\rm z}(-\hat \k_i)\hat y_i(\L_i,\l_i,\eta_i,\xi_i),\ {\cal R}_{\rm z}(-\hat \k_i)\hat x_i(\L_i,\l_i,\eta_i,\xi_i)\Big)$$
hence, use (\ref{inver Phi*}), (\ref{inver Phi}).
The bijectivity of this extension is trivial.

\vskip.1in
\noi
In particular, we have proven
\begin{proposition}\label{good variables}
Let ${\cal D}_{\rm r}$ is the set of $\Big((\L,\l),(\eta,\xi),(p,q),(G,g),(P,Q)\Big)\in (\real_+^N\times\torus^N)\times(\real^N\times\real^N)\times (\real^{N-2}\times\real^{N-2})\times (\real_+\times\torus)\times (\real\times \real)$ where the functions (\ref{GiPsii}) verify
$$\arr{\L_i>0\\
 0<\G_i\leq \L_i\ (1\leq i\leq N)\\
  |\Psi_{i-1}-\G_{i+1}|<\Psi_{i}\leq\Psi_{i-1}+\G_{i+1}\quad (1\leq i\leq N-2)\\ 
 |\Psi_{N-2}-\G_N|< G<\Psi_{N-2}+\G_N\\
  -G<C_{\rm z}=G-\frac{P^2+Q^2}{2}\leq G
  }
  $$
The (real--analytic and symplectic) ``full reduction'' change of variable 
\beqano
&& \phi_{\rm BD}:=\Phi_{\rm BD}^{-1}=\phi_*\circ\phi_{\rm r}:\nonumber\\
&& \Big((\L,\l),(\eta,\xi),(p,q),(G,g),(P,Q)\Big)\in{\cal D}_{\rm r}\to (y,x)=\Big((y_1,\cdots,y_N), (x_1,\cdots,x_N\Big)\in{\cal C}
 \eeqano
expressing the cartesian coordinates $(y,x)$ in terms of the regularized Deprit variables  is given by
\beqa{reg depr map}
\arr{y_i={\eufm R}_0{\eufm R}_i\hat y_i\\
x_i={\eufm R}_0{\eufm R}_i\hat x_i}\qquad 1\leq i\leq N
\eeqa
where ${\eufm R}_0={\eufm R}_0((P,Q),(G,g))$ is  (\ref{matrices1}), ${\eufm R}_i={\eufm R}_i(\L,(\eta,\xi),(p,q);G)$ is defined via (\ref{inver Phi}), \ie,
\beqa{eufm R i}
{\eufm R}_i&=&\arr{
{\cal R}_{\rm x}(\tilde i_{N-1})\left(\prod_{j=2}^{N-2}\tilde {\cal S}_{N-j}\right){\cal S}_1\quad  i=1\\
{\cal R}_{\rm x}(\tilde i_{N-1})\left(\prod_{j=2}^{N-i}\tilde{\cal S}_{N-j}\right){\rm I}_{\rm z}{\cal S}_i\quad 2\leq i\leq N-2\\
{\cal R}_{\rm x}(\tilde i_{N-1}){\rm I}_{\rm z}{\cal S}_{N-1}\quad i=N-1\\
{\rm I}_{\rm z}{\cal R}_{\rm x}(i_N)\quad i=N\\
}\nonumber\\
\eeqa
where $\tilde i_{N-1}$, $i_N\in (0,\p)$ are the inclinations (well defined and regular on ${\cal D}_{\rm r}$)
 $$\dst \tilde i_{N-1}=\cos^{-1}\left(\frac{G^2-\Psi_{N-2}^2-\G_N^2}{2G\Psi_{N-2}}\right)\ ,\quad \dst i_{N}=\cos^{-1}\left(\frac{G^2-\G_N^2-\Psi_{N-2}^2}{2G\G_N}\right)$$ 
with ${\cal S}_i$ ($1\leq i\leq N-1$), $\tilde{\cal S}_i$ ($2\leq i\leq N-2$)  as in Lemma \ref{Si are regular} 
and, finally, $(\L_i,\l_i,\eta_i,\xi_i)\to (\hat y_i,\hat x_i$) is the embedding
in $\real^{3}\times \real^{3}$ of the $(\hat m_i,\tilde m_i)$--Plane Delaunay--Poincar\'e Map.
\end{proposition}
We will refer to the map $\phi_{\rm BD}$ defined via Proposition \ref{good variables} as {\sl Regularized Deprit Map}.
\subsection{Partial Reduction}\label{Partial Reduction}
The regularized Deprit map discussed in Proposition \ref{good variables} becomes singular when 
\beq{elliptic singularity manifold}\sum_{1\leq i\leq N}\left(\L_i-\frac{\eta_i^2+\xi_i^2}{2}\right)-\sum_{1\leq i\leq N-2}\frac{p_i^2+q_i^2}{2}=G\ ,\eeq
which corresponds to the configuration with $S_{N-1}=\sum_{1\leq i\leq N-1}C_i$ parallel to $C_N$ (the two rotations ${\cal R}_{\rm x}(\tilde i_{N-1})$, ${\cal R}_{\rm x}(i_N)$ loss their regularity). Consider then the transformation,  denoted as $\phi_{\rm pr}$ (``{\rm pr}'' stands for ``partial reduction'') which acts as
\beqa{partial reduction}
&&  \arr{\dst
G=\sum_{1\leq k\leq N}\left(\L_k-\frac{\bar\eta_k^2+\bar\xi_k^2}{2}\right)-\sum_{1\leq k\leq N-1}\frac{\bar p_k^2+\bar q_k^2}{2}\\
\\
g=\arg{(\bar p_{N-1},-\bar q_{N-1})}\\
} \nonumber\\ \nonumber\\ 
&&   \l_i=\bar\l_i-g \quad 1\leq i\leq N\nonumber\\ 
&&  \left(\begin{array}{lr}
\eta_i\\
\xi_i
\end{array}
\right)={\cal R}_{\rm z}(g)
\left(\begin{array}{lr}
\bar \eta_i\\
\bar \xi_i
\end{array}
\right)\qquad \quad 1\leq i\leq N\nonumber\\
&&  \left(\begin{array}{lr}
p_j\\
q_j
\end{array}
\right)={\cal R}_{\rm z}(g)
\left(\begin{array}{lr}
\bar p_j\\
\bar q_j
\end{array}
\right)\quad 1\leq j\leq N-2
\eeqa 
leaving the remaining variables unvaried. It not difficult to check that $\phi_{\rm pr}$ is symplectic, since the Liouville $1$--form
remains unvaried:
$$\sum_{1\leq i\leq N}\L_id\l_i+\sum_{1\leq i\leq N}I_id\varphi_i+\sum_{1\leq i\leq N-2}J_id\psi_i+Gdg+PdQ$$ $$=\sum_{1\leq i\leq N}\L_id\bar\l_i+\sum_{1\leq i\leq N}I_id\bar\varphi_i+\sum_{1\leq i\leq N-2}J_id\bar \psi_i+J_{N-1}d\bar\psi_{N-1}+PdQ$$
where $(I_i,\varphi_i)$, $(J_j,\psi_j)$, $(I_i,\bar\varphi_i)$, $(J_k,\bar\psi_k)$ are the polar coordinates associated to $(\eta_i,\xi_i)$, $(p_j,q_j)$, $(\bar\eta_i,\bar\xi_i)$, $(\bar p_k,\bar q_k)$, with the indices $i$, $j$, $k$ running on $1\leq i\leq N$, $1\leq j\leq N-2$, $1\leq k\leq N-1$.

\noi
In terms of the variables $\Big((\L,\bar\l), (\bar\eta,\bar\xi), (\bar p,\bar q),(P,Q)\Big)$, the functions $\G=(\G_1,\cdots,\G_N)$, $\Psi=(\Psi_1,\cdots,\Psi_N)$
are
\beqa{new funct}\arr{\G_i=\L_i-\frac{\bar\eta_i^2+\bar\xi_i^2}{2}\qquad 1\leq i\leq N\\
\Psi_i=\sum_{1\leq i\leq i+1}\left(\L_j-\frac{\bar\eta_j^2+\bar\xi_j^2}{2}\right)-\sum_{1\leq j\leq i}\frac{\bar p_j^2+\bar q_j^2}{2}\qquad 1\leq i\leq N-1\quad (\Psi_{N-1}=G)\\
\Psi_N=C_{\rm z}=G-\frac{P^2+Q^2}{2}=\sum_{1\leq i\leq N}\left(\L_j-\frac{\bar\eta_j^2+\bar\xi_j^2}{2}\right)-\sum_{1\leq j\leq N-1}\frac{\bar p_j^2+\bar q_j^2}{2}-\frac{P^2+Q^2}{2}
}\eeqa
Denote now ${\cal D}_{\rm pr}$ the subset of \beqa{part red1}\Big((\L,\bar\l),(\bar\eta,\bar\xi),(\bar p,\bar q),(P,Q)\Big)\in\real_+^N\times\torus^N\times (\real^N\times\real^N)\times (\real^{N-1}\times\real^{N-1})\times\real\times\real\eeqa
where 
\beqa{part red2}\arr{\L_i>0\\
 0<\G_i\leq \L_i\ (1\leq i\leq N)\\
  |\Psi_{i-1}-\G_{i+1}|<\Psi_{i}\leq\Psi_{i-1}+\G_{i+1}\quad (1\leq i\leq N-1)\\ 
  -G<\Psi_N=C_{\rm z}\leq G
  }
\eeqa
(\ie, allow also $S_{N-1}\parallel C_N$, which corresponds to $(p_{N-1},q_{N-1})=0$). Then, the transormations $\phi_{\rm pr}$ regularizes, as the following Proposition claims, the proof of which is omitted.
\begin{proposition}\label{partial reduction prop}
Let ${\cal  D}_{\rm r}$ be defined via (\ref{part red1})$\div$(\ref{part red2}). The real--analytic and symplectic change of variable 
\beqano
&& \phi_{\rm  BD,pr}:=\phi_{\rm  BD}\circ\phi_{\rm pr}:\nonumber\\
&& \Big((\L,\bar\l),(\bar\eta,\bar\xi),(\bar p,\bar q),(P,Q)\Big)\in{\cal D}_{\rm pr}\to (y,x)=\Big((y_1,\cdots,y_N), (x_1,\cdots,x_N\Big)\in{\cal C}
 \eeqano
expressing the cartesian coordinates $(y,x)$ in terms of the {\sl partially reduced}, regularized Deprit variables $\Big((\L,\bar\l),(\bar\eta,\bar\xi),(\bar p,\bar q),(P,Q)\Big)$ is given by
\beqa{xi}
\arr{y_i={\eufm R}_0^{\rm pr}{\eufm R}_i^{\rm pr}\hat y_i\\
x_i={\eufm R}_0^{\rm pr}{\eufm R}_i^{\rm pr}\hat x_i}\qquad 1\leq i\leq N
\eeqa
where $(\L_i,\l_i,\eta_i,\xi_i)\to (\hat y_i,\hat x_i$) is the embedding
in $\real^{3}\times \real^{3}$ of the $(\hat m_i,\tilde m_i)$--Plane Delaunay--Poincar\'e Map
and 
$${\eufm R}_0^{\rm pr}={\cal R}_{\rm z}(\zeta(P,Q)){\cal R}_{\rm x}(i){\cal R}_{\rm z}(-\zeta(P,Q))\ ,$$ with 
$$i=\cos^{-1}\left(1-\frac{P^2+Q^2}{2\Psi_{N-1}}\right)\ ,\quad \zeta=\arg{(P,-Q)}$$
and
\beqano
{\eufm R}_i^{\rm pr}&=&\arr{
\left(\prod_{j=1}^{N-2}\tilde {\cal S}^{\rm pr}_{N-j}\right){\cal S}^{\rm pr}_1\quad  i=1\\
\left(\prod_{j=1}^{N-i}\tilde{\cal S}^{\rm pr}_{N-j}\right){\rm I}_{\rm z}{\cal S}^{\rm pr}_i\quad 2\leq i\leq N-1\\
{\rm I}_{\rm z}{\cal S}^{\rm pr}_N\quad i=N\\
}\nonumber\\
\eeqano
where the matrices ${\cal S}^{\rm pr}_j$'s, $\tilde{\cal S}^{\rm pr}_j$'s  have the following expressions:
\beqa{calSj part}
{\cal S}^{\rm pr}_j&=&\left(
\begin{array}{ccc}
1-\bar q_{j-1}^2{\eufm c}_j&-\bar p_{j-1}\bar q_{j-1}{\eufm c}_j&-\bar q_{j-1}{\eufm s}_j\\
-\bar p_{j-1}\bar q_{j-1}{\eufm c}_j&1-\bar p_{j-1}^2{\eufm c}_j&-\bar p_{j-1}{\eufm s}_j\\
\bar q_{j-1}{\eufm s}_j&\bar p_{j-1}{\eufm s}_j&1-(\bar p_{j-1}^2+\bar q_{j-1}^2){\eufm c}_j
\end{array}
\right)\qquad (1\leq j\leq N)\nonumber\\
\tilde{\cal S}^{\rm pr}_j&=&\left(
\begin{array}{ccc}
1-\bar q_{j}^2\tilde{\eufm c}_j&-\bar p_{j}\bar q_{j}\tilde{\eufm c}_j&-\bar q_{j}\tilde{\eufm s}_j\\
-\bar p_{j}\bar q_{j}\tilde{\eufm c}_j&1-\bar p_{j}^2\tilde{\eufm c}_j&-\bar p_{j}\tilde{\eufm s}_j\\
\bar q_{j}\tilde{\eufm s}_j&\bar p_{j}\tilde{\eufm s}_j&1-(\bar p_{j}^2+\bar q_{j}^2)\tilde{\eufm c}_j
\end{array}
\right)\qquad (2\leq j\leq N-1)\nonumber\\
\eeqa
where
\beqa{eufm ci si part}
&&\arr{{\eufm c}_1:=\frac{1-\cos{i_1}}{\bar p_{1}^2+\bar q_{1}^2}=\frac{\G_2-\G_1+\Psi_1}{4\G_1\Psi_1}\\
{\eufm s}_1:=\frac{\sin{i_1}}{\sqrt{\bar p_{1}^2+\bar q_{1}^2}}=\sqrt{{\eufm c}_1(2-(\bar p_1^2+\bar q_1^2){\eufm c}_1)}}\nonumber\\
\nonumber\\
&&\arr{
{\eufm c}_j:=\frac{1-\cos{i_j}}{\bar p_{j-1}^2+\bar q_{j-1}^2}=\frac{\Psi_{j-2}-\G_j+\Psi_{j-1}}{4\G_j\Psi_{j-1}}\\
{\eufm s}_j:=\frac{\sin{i_j}}{\sqrt{\bar p_{j-1}^2+\bar q_{j-1}^2}}=\sqrt{{\eufm c}_j(2-(\bar p_{j-1}^2+\bar q_{j-1}^2){\eufm c}_j)}}\quad (2\leq j\leq N)\nonumber\\
\nonumber\\
&&\arr{
\tilde{\eufm c}_j:=\frac{1-\cos\tilde{i_j}}{\bar p_{j}^2+\bar q_{j}^2}=\frac{\G_{j+1}+\Psi_j-\Psi_{j-1}}{4\Psi_j\Psi_{j-1}}\\
\tilde{\eufm s}_j:=\frac{\sin\tilde{i_j}}{\sqrt{\bar p_{j}^2+\bar q_{j}^2}}=\sqrt{\tilde{\eufm c}_j(2-(\bar p_j^2+\bar q_j^2)\tilde{\eufm c}_j)}}\quad (2\leq j\leq N-1)\nonumber\\
\eeqa
where, $\G_i$, $\Psi_j$ are thought as functions of $(\L,\bar\eta,\bar\xi,\bar p,\bar q)$:
\beqa{eufm R pr}\arr{\G_i=\L_i-\frac{\bar\eta_i^2+\bar\xi_i^2}{2}\qquad 1\leq i\leq N\\
\Psi_i=\sum_{1\leq i\leq i+1}\left(\L_j-\frac{\bar\eta_j^2+\bar\xi_j^2}{2}\right)-\sum_{1\leq j\leq i}\frac{\bar p_j^2+\bar q_j^2}{2}\qquad 1\leq i\leq N-1\quad (\Psi_{N-1}=G)
}\eeqa
\end{proposition}
\vskip.1in
\noi
For $1\leq i\leq N-1$ and $1\leq j\leq N-2$, the functions ${\eufm c}_i$, ${\eufm s}_i$, $\tilde{\eufm c}_j$, $\tilde{\eufm s}_j$ defined in (\ref{eufm ci si part}) coincide with the corresponding functions related to the full reduction (eq. (\ref{eufm ci si})).

\vskip.1in
\noi
We will refer to the maps $\phi_{\rm  BD}$, $\phi_{\rm BD,pr}$ as regularized {\sl full reduction} (or, simply, {\sl  reduction}), regularized {\sl partial reduction}  (or, simply, {\sl partial reduction}), respectively.
\newpage
\section{Kolmogorov's Set in the Space Planetary Problem I (Partial Reduction)}
\setcounter{equation}{0}
\subsection{Non--Degeneracy Conditions ($N\geq 3$)}
The construction of KAM tori for the spatial planetary problem with $N\geq 3$ planets via Theorem \ref{more general degenerate KAM} becomes quite natural and direct, with the use of  the Deprit variables.

\vskip.1in
\noi
In this section, we show  that, for $N\geq 3$, the set of  {\sl Deprit's partially reduced} variables  discussed in section \ref{Partial Reduction} (which, we recall, corresponds to the reduction of  $C_{\rm z}$) is a good set of coordinates in order to obtain KAM tori with $3N-1$ frequencies. The pregium of Deprit's partial reduction is that, differently from what happens  trying a partial reduction in Poincar\'e--Delaunay's variables,
\footnote{In Poincar\'e--Delaunay variables, the third component of the angular momentum is $$C_{\rm z}=\sum_{1\leq i\leq N}\left(\L_i-\frac{\eta_i^2+\xi_i^2}{2}-\frac{p_i^2+q_i^2}{2}\right)\ ,$$ 
which is quite the same of the expression of the modulus of the anguar momentum $G$ in partially reduced Deprit variables (first equation in (\ref{partial reduction})), apart for the dimension ($N$ instead than $N-1$) of the $(p,q)$ variables. 
Put then, similarly to (\ref{partial reduction}),
\beqano
&&  \arr{\dst
C_{\rm z}=\sum_{1\leq k\leq N}\left(\L_k-\frac{\eta_k^2+\xi_k^2}{2}\right)-\sum_{1\leq k\leq N}\frac{ p_k^2+ q_k^2}{2}\\
\\
\zeta=\arg{( p_{N},- q_{N})}\\
}\ ,\quad   \hat\l_i=\l_i-\zeta \quad 1\leq i\leq N\nonumber\\ 
&&  \left(\begin{array}{lr}
\hat\eta_i\\
\hat\xi_i
\end{array}
\right)={\cal R}_{\rm z}(\zeta)
\left(\begin{array}{lr}
 \eta_i\\
 \xi_i
\end{array}
\right)\qquad \quad 1\leq i\leq N\ ,\quad   \left(\begin{array}{lr}
\hat p_j\\
\hat q_j
\end{array}
\right)={\cal R}_{\rm z}(\zeta)
\left(\begin{array}{lr}
 p_j\\
 q_j
\end{array}
\right)\quad 1\leq j\leq N-1
\eeqano 
The variables $\Big((\L,\hat\l),(\hat\eta,\hat\xi),(\hat p,\hat q),(C_{\rm z},\zeta)\Big)$ realize a (Delaunay) partial reduction, however, singular, relatively to the configurations with vanishing $N^{th}$ inclination, \ie, when $$\frac{p_N^2+q_N^2}{2}=\sum_{1\leq i\leq N}\left(\L_i-\frac{\hat \eta_i^2+\hat \xi_i^2}{2}\right)-\sum_{1\leq i\leq N-1}\frac{\hat p_i^2+\hat q_i^2}{2}-C_{\rm z}=0\ .$$ This singularity is sometimes called ``elliptic singularity''.}
it leaves the secular perturbation regular and symmetric around the secular origin, which, as in the planar case, turns out to be an elliptic equilibrium point, corresponding to the configurations with all zero eccentricities and mutual inclinations.

\vskip.1in
\noi
As said before, the construction of the KAM tori is obtained as an application of Theorem \ref{more general degenerate KAM}, so, it is based on the check of the two non degeneracy conditions thereby involved:
\begin{itemize}
\item[(i)] {\sl check of $4$--non resonance for the Birkhoff invariants of the first order};
\item[(ii)] {\sl check of second order non degeneracy, \ie, proof of non singularity of the second order Birkhoff invariants} matrix. 
\end{itemize}
Both ${\rm (i)}$ and ${\rm (ii)}$ are proved by induction, in the range of well separated semimajor axes.
\vskip.1in
\noi
The restriction to $N\geq 3$ is due to the following. When the secular perturbation is put in Deprit partially reduced variables, as in the case of Poincar\'e--Delaunay variables, its quadratic part splits into the sum of a ``horizontal'' part ${\cal Q}^*_h$ and a `` vertical'' part ${\cal Q}^*_v$, of order $N$, $N-1$ respectively. Then, using partial reduction, a unique secular resonance (compare Proposition \ref{true non resonance} below) is exhibited by the respective eigenvalues $s=(s_1$, $\cdots$, $s_N)$, $z=(z_1,\cdots,z_{N-1})$ of ${\cal Q}^*_h$, ${\cal Q}^*_v$, once again the {\sl Herman's resonance}:
\beq{Hermann res}\sum_{1\leq i\leq N}s_i+\sum_{1\leq i\leq N-1}z_i=0\ .\eeq
This resonance is  of order $2N-1$, hence, it  prevents the construction of the Birkhoff normal form up to order $4$ only when $N=2$. 
When $N\geq 3$, the Herman's resonance (\ref{Hermann res}) is of high order ($2N-1\geq 5$), allowing us to use partial reduction for the construction (and proof of non--degeneracy) of the normal form.

\vskip.1in
\noi
The first step consists into the expansion of the secular perturbation up to order $2$. We denote
$${\cal H}_{\rm plt,pr}:={\cal H}_{\rm plt}\circ\phi_{\rm BD,pr}=h_{\rm plt}+\m f_{\rm plt,pr}$$
the planetary Hamiltonian function, put in regularized, partially reduced Deprit variables, 
where, as usual
$$h_{\rm plt}=-\sum_{1\leq i\leq N}\frac{\tilde m_i^3\hat m_i^2}{2\L_i^2}$$
is the Kepler's unperturbed integrable part.
\begin{lemma}\label{expansion to order two}
For $N\geq 2$, the mean  $\dst \bar f_{\rm  plt,pr}:=(2\p)^{-N}\int_{\torus^N}f_{\rm  plt,pr}d\bar\l$   is an even function of the ``secular fully regularized variables'' $\bar z:=(\bar\eta,\bar\xi,\bar p,\bar q)$ and its expansion around $\bar z=0$ is the following. Define:

\vskip.1in
\noi
the constants
\beq{C}\arr{
C_0(m,a):=-\sum_{1\leq j<k\leq N}\frac{m_jm_k}{a_j}b_{1/2,0}(a_k/a_j)\\
C_1(a_j,a_k):=-\frac{a_k}{2a_j^2}b_{3/2,1}(a_k/a_j)\\
C_2(a_j,a_k):=\frac{a_k}{a_j^2}b_{3/2,2}(a_k/a_j)
}
\eeq
where $\a\to b_{s,k}(\a)$ is the $(s,k)$--Laplace coefficient;

\vskip.1in
\noi
the quadratic forms
\beq{QhQv}
\arr{
{\cal Q}_h\cdot \bar\eta^2:=\sum_{1\leq j<k\leq N}m_jm_k\left(C_1(a_j,a_k)\left(\frac{\bar\eta_j^2}{\L_j}+\frac{\bar\eta_k^2}{\L_k}\right)+2C_2(a_j,a_k)\frac{\bar\eta_j\bar\eta_k}{\sqrt{\L_j\L_k}}\right)\\
\hat{\cal Q}_v\cdot \bar p^2:=-\sum_{1\leq j<k\leq N}m_jm_kC_1(a_j,a_k)\Big({\bar p_j}-{\bar p_k}\Big)^2
}
\eeq
the linear operator
$${\eufm L}:\quad \bar p\in\real^{N-1}\to {\eufm L}\bar p=\Big(\cdots,{\eufm L}\bar p_i,\cdots\Big)\in \real^{N}$$  which acts as
\beqa{Lop}\arr{{\eufm L}\bar p_{1}:=\ell_1\cdot\bar p:=c_1\,\bar p_1+{\sum}^*_{2\leq j\leq N-1}\tilde c_j\,\bar p_j\\
{\eufm L}\bar p_i:=\ell_i\cdot\bar p:=c_i\,\bar p_{i-1}+{\sum}^*_{i\leq j\leq N-1}\tilde c_j\,\bar p_j\quad 2\leq i\leq N-1\quad (\textrm{if}\quad N\geq 3)\\
{\eufm L}\bar p_{N}:=\ell_N\cdot\bar p:=c_{N}\,\bar p_{N-1}\\
}\ ,\eeqa
where the summands denoted $\Sigma^*$ does not appear when $N=2$, and
\beqa{ccc}\arr{\dst c_1:=-\sqrt{\frac{\L_2}{\L_1{\rm L}_2}}\\
\dst c_j:=\sqrt{\frac{{\rm L}_{j-1}}{{\rm L}_j\L_j}}\quad (2\leq j\leq N)\\
\dst \tilde c_j:=-\sqrt{\frac{\L_{j+1}}{{\rm L}_{j+1}{\rm L}_j}}\quad (2\leq j\leq N-1)\\
\dst {\rm L}_j:=\sum_{1\leq k\leq j}\L_k
}
\eeqa 
Then,
\beq{Q*h}\bar f_{\rm  plt,pr}=C_0(m,a)+{\cal Q}^*_h\cdot\frac{\bar\eta^2+\bar\xi^2}{2}+{\cal Q}^*_v\cdot\frac{\bar p^2+\bar q^2}{2}+\bar f^4_{\rm plt, pr}\ ,\eeq
where
\beq{Q*v}\arr{
{\cal Q}^*_h\cdot\eta^2:={\cal Q}_h\cdot\eta_*^2\ ,\quad \textrm{with}\quad\eta_*:=(-\bar\eta_1,\bar\eta_2,\cdots,\bar\eta_N)\\
Q^*_v\cdot \bar p^2:=\hat{\cal Q}_v\cdot ({\eufm L}\bar p)^2\\
\bar f^4_{\rm plt, pr}={\rm O(4)}\ .
}\eeq
\end{lemma}
The details of the computation of the expansion of $\bar f_{\rm  plt,pr}$ (up to order $4$, for future use) are in Section \ref{Expansion of the Hamiltonian space}.
\subsubsection{First Order Conditions}
\begin{proposition}\label{true non resonance}
For any $N\geq 2$, there exists an open set with full measure ${\cal U}\subset{\cal A}$ where the eigenvalues of ${\cal Q}^*_h$ and ${\cal Q}^*_v$  are paiwise distinct and verify the following. For any open, simply connected set ${\cal V}\subset {\cal U}$, they define $2N-1$ holomorphic functions $s_1$, $\cdots$, $s_N$, $z_1$, $\cdots$, $z_{N-1}$ which satisfy the only linear relation
\beq{only res}\sum_{1\leq i\leq N}s_i+\sum_{1\leq i\leq N-1}z_i=0\eeq
(up to an arbitrary multiplicative constant).
\end{proposition}
{\bf Proof.}\ Let us introduce matricial notations. Let ${\cal F}_h$, ${\cal F}_v$ denote the matrices (having both order $N$) associated to the quadratic forms ${\cal Q}_h$, ${\cal Q}_v$ (\ref{PoincDel1}) of the quadratic part of the secular peturbation in Delaunay--Poincar\'e variables and let ${\cal F}_h^*$, ${\cal G}_v^*$ denote the matrices associated to ${\cal Q}^*_h$, ${\cal Q}^*_v$ (having order $N$, $N-1$ respectively) of eq. (\ref{Q*v}):
$$\arr{{\cal Q}_h\cdot \eta^2=\eta\cdot {\cal F}_h\eta\\
{\cal Q}_v\cdot P^2=P\cdot {\cal F}_vP
}\qquad \arr{{\cal Q}^*_h\cdot \eta^2=\eta\cdot {\cal F}_h^*\eta\\
{\cal Q}^*_v\cdot p^2=p\cdot {\cal G}_v^*p
}$$
The matrices ${\cal F}_h^*$,  ${\cal F}_h$ are related by
\footnote{Using the symplectic variables $(\eta_*,\xi_*)$ (compare (\ref{Q*v})), rather than $(\eta,\xi)$, namely, with $(\eta^*_1,\xi^*_1)=-(\eta_1,\xi_1)$  related to the aphelion position of the first osculating ellipse, rather than the perihelion, would transform $\cF_h^*$ to $\cF_h$.  We do not use this change of variables here, because unnecessary, but in the next section, for the computation of the Birkhoff invariants of order $2$, in order to have simpler expressions.}
$${\cal F}_h^*={\rm I}{\cal F}_h{\rm I}$$ where ${\rm I}$ changes the sign of the first coordinate, so, they have the same eigenvalues. Then, in order to prove (\ref{only res}), in view, of (\ref{Fej only res}), we only need to prove that ${\cal F}_v$, ${\cal F}_v^*$ have the same trace. 
We have
\beq{cal G cal G*}{\cal F}_v^*=\ell^{\rm T}\hat{\cal F}_v^*\ell\ ,\quad {\cal F}_v=\ell_0^{\rm T}\hat{\cal F}_v\ell_0\ ,\eeq
where $\ell$  denotes  the $N\times (N-1)$ matrix associated to the linear operator (\ref{Lop})$\div$(\ref{ccc}), $\ell_0$  the diagonal matrix $$\ell_0=\textrm{\rm diag}\,(1/\sqrt{\L_1},\cdots,1/\sqrt{\L_N})$$ and $\hat{\cal F}_v=(\hat g_{ij})$  the $N\times N$ matrix of $\hat{\cal Q}_v$:
\beq{hatcalQ}\hat{\cal Q}_v\cdot \bar p^2:=-\sum_{1\leq j<k\leq N}m_jm_kC_1(a_j,a_k)\Big({\bar p_j}-{\bar p_k}\Big)^2:=\bar p\cdot\hat{\cal F}_v\bar p\ .\eeq
Equations (\ref{cal G cal G*}) give
\beqa{trG*}
\arr{
\textrm{\rm tr}\Big({\cal F}_v^*\Big)=\textrm{\rm tr}\Big(\ell^{\rm T}\hat{\cal F}_v\ell\Big)=\textrm{\rm tr}\Big({\cal L}\,\hat{\cal F}_v\Big)\\
\textrm{\rm tr}\Big({\cal F}_v\Big)=\textrm{\rm tr}\Big(\ell_0^{\rm T}\hat{\cal F}_v\ell_0\Big)=\textrm{\rm tr}\Big({\cal L}_0\,\hat{\cal F}_v\Big)
}
\eeqa
where ${\cal L}_0$ is the diagonal matrix $${\cal L}_0=\textrm{\rm diag}(1/\L_1,\cdots,1/\L_N)\ $$
and 
${\cal L}$ is the symmetric matrix with entries
$${\cal L}_{ij}=\ell_i\cdot\ell_j$$
if 
$$\ell_1=(c_1,\tilde c_2,\cdots,\tilde c_{N-1})\ ,\ \ell_2=(c_2,\tilde c_2,\cdots,\tilde c_{N-1})\ ,\cdots,\ell_i=(0,\cdots,0,c_{i},\tilde c_i,\cdots,\tilde c_{N-1})$$
is the $i^{\rm th}$ row of $\ell$, as in (\ref{Lop})$\div$(\ref{ccc}). We have
\beqa{telescopic sums}
\arr{\dst{\cal L}_{11}=|\ell_1|^2=c_1^2+\sum_{2\leq j\leq N-1}\tilde c_j^2\\
\dst{\cal L}_{ii}=|\ell_i|^2=c_1^2+\sum_{i\leq j\leq N-1}\tilde c_j^2\quad 1\leq j\leq N-1\\
\dst{\cal L}_{NN}=|\ell_N|^2=c_N^2\\
\dst{\cal L}_{ij}=\ell_i\cdot\ell_j=c_j\tilde c_{j-1}+\sum_{j\leq k\leq N-1}\tilde c_k^2\quad 1\leq i<j\leq N-1\\
\dst{\cal L}_{iN}=\ell_i\cdot\ell_N=c_N\tilde c_{N-1}\quad 1\leq i\leq N-1
}
\eeqa
Using  (\ref{ccc}), we can write
\beqano
\arr{
c_1^2=\frac{1}{\L_1}-\frac{1}{{\rm L}_2}\\
c_j^2=\frac{1}{\L_j}-\frac{1}{{\rm L}_j}\quad 2\leq j\leq N\\
\tilde c_j^2=\frac{1}{{\rm L}_j}-\frac{1}{{\rm L}_{j+1}}\quad 2\leq j\leq N-1\\
c_j\tilde c_{j-1}=-\frac{1}{{\rm L}_j}\quad 2\leq j\leq N\quad (\tilde c_1:=c_1)
}
\eeqano hence, inserting thee previous expressions into (\ref{telescopic sums}), by telescopic arguments,
$${\cal L}={\cal L}_0-\frac{{\rm E}}{{\rm L}_N}\ ,$$
where ${\rm E}$ has entries ${\rm E}_{ij}=1$ for any $i$, $j$.
Hence, in view of (\ref{trG*})
$$\textrm{\rm tr}\Big({\cal F}_v^*\Big)=\textrm{\rm tr}\left(\left({\cal L}_0-\frac{{\rm E}}{{\rm L}_N}\right)\hat{\cal F}_v\right)=\textrm{\rm tr}\left({\cal L}_0\hat{\cal F}_v\right)-\frac{1}{{{\rm L}_N}}\textrm{\rm tr}\left({{\rm E}}\hat{\cal F}_v\right)=\textrm{\rm tr}\left({\cal L}_0\hat{\cal F}_v\right)=\textrm{\rm tr}\left({\cal F}_v\right)$$
since, from (\ref{hatcalQ}), it is easy to check
$$\textrm{\rm tr}\Big({{\rm E}}\hat{\cal F}_v\Big)=\sum_{j,k}\hat g_{jk}=0\ .$$
We prove now uniqueness of (\ref{only res}), proceeding by induction on the number $N$ of planets. For $N=2$, it is a consequence of existence and the fact that the planar eigenvalues $s_1=\s_1$, $s_2=\s_2$ do not satisfy any linear condition, as proved in \cite{Fej04}. Assume now that uniqueness of (\ref{only res}) holds for $N-1$ and let  $$c=(c_1,\cdots,c_N)\in \real^N\ ,\quad g=(g_1,\cdots,g_{N-1})\in \real^{N-1}$$ such that 
$$c\cdot s+g\cdot z=0$$
where $s=(s_1,\cdots,s_N)$, $z=(z_1,\cdots,z_{N-1})$ are the eigenvalues of ${\cal F}_h:={\cal F}_h^N$, ${\cal F}_v:={\cal F}_v^{*N}$. We write explicitely  the entries of the matrices ${\cal F}_h^N$, ${\cal F}_v^{*N}$
\beqa{entries of cal Fh cal Fv}({\cal F}_h^N)_{ij}&=&\arr{
\frac{\bar m_i}{\L_i}\sum_{i<k\leq N}\bar m_k C_1(a_i,a_k)+\frac{\bar m_i}{\L_i}\sum_{1\leq k<i}\bar m_k C_1(a_k,a_i)\quad 1\leq i=j\leq N\\
\\
\frac{\bar m_i\bar m_j}{\sqrt{\L_i\L_j}} C_2(a_i,a_j)\quad 1\leq i<j\leq N
}\nonumber\\
\nonumber\\
\nonumber\\
({\cal F}_v^{*N})_{ij}&=&\arr{
-\tilde c_i^2\sum_{1\leq k\leq i}\sum_{l\geq i+2}\bar m_k\bar m_lC_1(a_k,a_l))-c_{i+1}^2\sum_{k\geq i+2}\bar m_{i+1}\bar m_k C_1(a_{i+1},a_k)\\
-\bar c_i^2\sum_{1\leq k\leq i}\bar m_k\bar m_{i+1}C_1(a_k,a_{i+1})\quad 1\leq i=j\leq N-2\\
\\
-\bar c_{N-1}^2\sum_{1\leq k\leq N-1}\bar m_k\bar m_NC_1(a_k,a_N)\quad i=j=N-1\\
\\
-\tilde c_i\tilde c_j\sum_{1\leq k\leq i,\ l\geq j+2}\bar m_k\bar m_lC_1(a_k,a_l)-\bar c_j\tilde c_i\sum_{1\leq k\leq i}\bar m_k\bar m_{j+1}C_1(a_k,a_{j+1})\\
-\bar c_jc_{i+1}\bar m_{i+1}\bar m_{j+1}C_1(a_{i+1},a_{j+1})-\tilde c_jc_{i+1}\sum_{l\geq j+2}\bar m_{i+1}\bar m_lC_1(a_{i+1},a_l)\\
\quad 1\leq i< j\leq N-2\\
\\
-\bar c_{N-1}\tilde c_i\sum_{1\leq k\leq i}\bar m_k\bar m_NC_1(a_k,a_N)
-\bar c_{N-1}c_{i+1}\bar m_{i+1}\bar m_N C_1(a_{i+1},a_N)\\
\quad 1\leq i< j=N-1\\
\\
}\nonumber\\
\eeqa
where $c_i$, $\tilde c_j$ are as in (\ref{ccc}) and
$$\bar c_k:=c_{k+1}-\tilde c_k=\sqrt{\frac{1}{\L_{k+1}}+\frac{1}{{\rm L}_k}}\quad 1\leq k\leq N-1\ .$$
It is easy to check  that, when $a_N\to \infty$,
$${\cal F}_h^N\to \left(
\begin{array}{lrr}
\hat{\cal F}_h^{N-1}&0\\
0&0
\end{array}
\right)\qquad {\cal F}_v^{*N}\to \left(
\begin{array}{lrr}
\hat{\cal F}_v^{*(N-1)}&0\\
0&0
\end{array}
\right)$$
and that when $a_1\to 0$,
$${\cal F}_h^N\to \left(
\begin{array}{lrr}
0&0\\
0&\check{\cal F}_h^{N-1}
\end{array}
\right)\qquad {\cal F}_v^N\to \left(
\begin{array}{lrr}
0&0\\
0&\check{\cal F}_v^{*(N-1)}
\end{array}
\right)$$
where  $\hat{\cal F}_{h}^{N-1}$, $\hat{\cal F}_{v}^{*(N-1)}$ denote the horizontal, vertical quadratic forms related to the ``first'' $N-1$ bodies and  $\check{\cal F}_{h}^{N-1}$, $\check{\cal F}_{v}^{*(N-1)}$, the horizontal, vertical quadratic forms related to the ``last'' $N-1$ bodies. Then, as $s$, $z$ are continuous function of the entries of their respective matrices, when $a_N\to \infty$,
$$s\to (\hat s,0)\ ,\quad z\to (\hat z,0)$$
and, when $a_1\to 0$,
$$s\to (0,\check s)\ ,\quad z\to (0,\check z)$$
where $\hat s$, $\hat z$ are the eigenvalues of  $\hat{\cal F}_{h}^{N-1}$, $\hat{\cal F}_{v}^{*(N-1)}$; $\check s$, $\check z$ are the eigenvalues of  $\check{\cal F}_{h}^{N-1}$, $\check{\cal F}_{v}^{*(N-1)}$. By the inductive hypoyhesis, the first limit implies $\dst c_1=\cdots=c_{N-1}=g_1=\cdots=g_{N-2}$; the second limit implies $\dst c_2=\cdots=c_{N}=g_2=\cdots=g_{N-1}$, hence, the thesis.
\subsubsection{Second Order Conditions (``Torsion'')}
By Proposition \ref{true non resonance} and Birkhoff theory, when $N\geq 3$, the secular perturbation $\bar f_{\rm plt,pr}$ of the planetary problem can be put in normal form up to order $4$. 

\noi
In this section, we prove the non--degeneracy of this normal form. 

\vskip.1in
\noi
\begin{itemize}
\item[(i)]
We call {\bf Birkhoff form} of a given polynomial of $2m$ variables $$(y,x)\in\real^m\times\real^m\to{\cal P}(y,x)\in \real$$  and even degree $p\geq 4$,  the polynomial
$${\cal P}_{\rm B}(y,x):=\frac{1}{(2\p)^m}\int_{\torus^m}{\cal P}\circ\phi_{\rm pc}(J,\varphi)d\varphi\quad \textrm{with}\quad J=\left(\cdots,\frac{y_i^2+x_i^2}{2},\cdots\right)$$ 
where $\phi_{\rm pc}$ is the usual symplectic polar coordinates map.
\item[(ii)] When ${\cal P}$ has degree $4$, we call {\bf Birkhoff matrix} associated to ${\cal P}$ the {\sl symmetric} matrix $A=(A_{ij})$ of  order $m$ for which
$${\cal P}_{\rm B}(y,x)=\frac{1}{2}\sum_{1\leq i,j\leq m}A_{ij}J_iJ_j\quad \textrm{where}\quad J_i:=\frac{y_i^2+x_i^2}{2}$$
\end{itemize}
\vskip.1in
\noi
Let $\r_h^*$, $\r_v^*$ denote the matrices which diagonalize the quadratic forms ${\cal Q}_h^*$, ${\cal Q}_v^*$, and let $\phi_{\rm diag}$ the symplectic transformation{\footnote{We neglect to write the action on the $\bar\l$--variables, which we do non need.}
\beqano
\phi_{\rm diag}:\quad &&\arr{\bar\eta=\r_h^*\tilde\eta\\ \bar\xi=\r_h^*\tilde\xi}\Leftrightarrow\arr{\eta_*=\r_h\tilde\eta\\ \xi_*=\r_h\tilde\xi}\ ,\nonumber\\ &&\arr{\bar p=\r_v^*\tilde p\\ \bar q=\r_v^*\tilde q}\ ,\quad \L=\tilde\L
\eeqano
where $\r_h={\rm I}\r_h^*$ is the matrix which diagonalizes the quadratic form ${\cal Q}_h={\rm I}{\cal Q}_h^*{\rm I}$ of the plane problem, as in \cite{Fej04}.
Then,  the secular perturbation $\bar f_{\rm plt,pr}$ is put into the form
\beqano
\bar f_{\rm diag}&:=&\bar f_{\rm  plt,pr}\circ\phi_{\rm diag}\nonumber\\
&=&C_0(m,a)+\sum_{1\leq i\leq N}s_i\frac{\tilde\eta_i^2+\tilde\xi_i^2}{2}+\sum_{1\leq i\leq N-1}z_i\frac{\tilde p_i^2+\tilde q_i^2}{2}+{\eufm F}\Big(\r_h\tilde\eta,\r_h\tilde\xi,\r_v^*\tilde p,\r_v^*\tilde q\Big)\nonumber\\
&+&O(6)
\eeqano
where ${\eufm F}$ is the polynomial of degree $4$ in $ z_*=(\eta_*,\xi_*, \bar p, \bar q)$ for which $$\dst \tilde f^4_{\rm plt,pr}(\bar\eta,\bar\xi,\bar p,\bar q)={\eufm F}(\eta_*,\xi_*,\bar p,\bar q)+O(6)\ .$$

\vskip.1in
\noi
We then have
\begin{proposition}\label{spatial non degeneracy}
For any $N\geq 2$, the Birkhoff matrix $A_{\rm plt}$, of order $2N-1$, associated to the polynomial $\dst {\eufm F}(\r_h\eta$, $\r_h\xi$, $\r_v^*p$, $\r_v^*q)$ is non singular, provided the semimajor axes $0<a_1<a_2<\cdots<a_N$ are well separated.
\end{proposition}
\begin{remark}\rm When $N\geq 3$, $A_{\rm plt}$ coincides with the matrix of the {\sl Birkhoff invariants of order $2$} of the planetary problem.
\end{remark}
In order to prove Proposition \ref{spatial non degeneracy}, we need the exact expression of  ${\eufm F}$, which  is computed in Section \ref{Expansion of the Hamiltonian space} and summarized in the following Lemma. 
\begin{lemma}\label{order four exp}
In the expansion for the  secular perturbation $\bar f_{\rm  plt,pr}$  around the secular origin $\bar z=(\bar\eta,\bar\xi,\bar p,\bar q)=0$ described in Lemma \ref{expansion to order two}, the term $\bar f^4_{\rm plt, pr}$ of eq. (\ref{Q*h}) begins with $\dst \bar f^4_{\rm plt, pr}={\eufm F}+O(6)$, where 
\beq{eufm F}{\eufm F}={\eufm F}_h+{\eufm F}_{hv}+{\eufm F}_v\eeq
and ${\eufm F}_h$, ${\eufm F}_{hv}$, ${\eufm F}_v$ are three polynomials of degree $4$, defined as follows. The ``horizontal'' part ${\eufm F}_h$ is 
\beq{horizquartic}{\eufm F}_h=q\cdot(\eta_*^4+\xi_*^4)+r\cdot\eta_*^2\xi_*^2\eeq
where, as in Lemma \ref{expansion to order two}, $\eta_*:=(-\eta_1,\eta_2,\cdots)$ and 
\beqa{rq}
q\cdot\eta^4&:=&\sum_{1\leq i\leq N}q_{iiii}\eta_i^4+\sum_{1\leq i<j\leq N}\Big(q_{iiij}\eta_i^3\eta_j+q_{iijj}\eta_i^2\eta_j^2+q_{jjji}\eta_i\eta_j^3\Big)\nonumber\\
r\cdot\eta^2q^2&:=&\sum_{1\leq i\leq N}r_{iiii}\eta_i^2\xi_i^2+\sum_{1\leq i<j\leq N}\Big(r_{iiij}\eta_i^2\xi_i\xi_j+r_{iijj}\eta_i^2\xi_j^2+r_{ijii}\eta_i\eta_i\xi_i^2+r_{ijij}\eta_i\eta_j\xi_i\xi_j\nonumber\\
&&+r_{jijj}\eta_i\eta_j\xi_j^2+r_{jjii}\eta_j^2\xi_i^2+r_{jjji}\eta_j^2\xi_i\xi_j\Big)\nonumber\\
\eeqa
if
$q=(q_{ijkl})$, $r=(r_{ijkl})$ are the $4$--indices tensors defining the quartic form of the secular perturbation of the plane problem, defined in  (\ref{r and q})$\div$(\ref{a4000 and others}). The ``vertical'' parts ${\eufm F}_{hv}$,  ${\eufm F}_v$ are
\beqa{vertical quadratic}
{\eufm F}_{hv}&:=&\frac{1}{2}\sum_{1\leq i<j\leq N}\bar m_i\bar m_j\Big({\cal Q}_{ij}^{11}\,({\eufm L}\bar q_j-{\eufm L}\bar q_i)^2+{\cal Q}_{ij}^{22}\,({\eufm L}\bar p_j-{\eufm L}\bar p_i)^2\nonumber\\
&+&({\cal Q}_{ij}^{12}+{\cal Q}_{ij}^{21})\Big({\eufm L}\bar p_j-{\eufm L}\bar p_i\Big)\Big({\eufm L}\bar q_j-{\eufm L}\bar q_i\Big)\nonumber\\
&+&({\cal Q}_{ij}^{21}-{\cal Q}_{ij}^{12})\,{\sum}^*_{1\leq h<k\leq j-1}\ell^{ij}_h\ell^{ij}_k(\bar p_h\bar q_k-\bar p_k\bar q_h)\Big)\nonumber\\
{\eufm F}_v&:=&\sum_{1\leq i<j\leq N}{\bar m_i\bar m_j}({\eufm Q}^{ij}_{11}+{\eufm Q}^{ij}_{22})C_1(a_i,a_j)\nonumber\\
&+&\sum_{1\leq i<j\leq N}{\bar m_i\bar m_j}\left(\Big({\eufm L}\bar p_i-{\eufm L}\bar p_j\Big)^2+\Big({\eufm L}\bar q_i-{\eufm L}\bar q_j\Big)^2\right)^2C_{13}(a_i,a_j)\nonumber\\
&+&\sum_{1\leq i<j\leq N}{\bar m_i\bar m_j}\left({\sum}^*_{1\leq h<k\leq j-1}\ell^{ij}_h\ell^{ij}_k(\bar p_{h}\bar q_k-\bar p_k\bar q_h)\right)^2C_{14}(a_i,a_j)\nonumber\\
\eeqa
where the  summands denoted ${\sum}^*$ do not appear  when $N=2$ or $j=2$ and where:
\begin{itemize}
\item[--]for $h,\ k=1,\ 2$, ${\cal Q}_{ij}^{hk}$ are the four quadratic forms acting on $(\eta^*_i,\eta^*_j,\xi^*_i,\xi^*_j)$ as in (\ref{Iij}), with $C_3(a_i,a_j)\div C_{12}(a_i,a_j)$ as in (\ref{Cij}) below;
\item[--]${\eufm L}$ is the linear operator from $\real^{N-1}$ to $\real^N$ defined  in eq. (\ref{Lop})$\div$(\ref{ccc}); 
\item[--]$\ell^{ij}_k:=\ell_{jk}-\ell_{ik}$ with $(\ell_{ik}):=\ell$  the matrix  associated to ${\eufm L}$;
\item[--]$C_{13}(a_i,a_j)$, $C_{14}(a_i,a_j)$ are as in (\ref{C13C14}) below and, finally, :
\item[--]${\eufm Q}^{ij}_{11}$, ${\eufm Q}^{ij}_{22}$ are the two first diagonal entries of the matrix of homogeneous polynomials of degree $4$ resulting from the productories (\ref{eufmSij})  below, when ${\cal S}_i$, $\cdots$ have the expansions defined in eq. (\ref{Sigma order 4})$\div$(\ref{R0}) below.
\end{itemize}
\end{lemma}

\vskip.1in
\noi
We outline that the main difference (and complexity), with respect to the same computation in Delaunay--Poincar\'e variables, is represented by the productory form  of the ``verticalizing'' matrices ${\eufm S}_{ij}$ (compare equation (\ref{eufmSij}) below) which describe the mutual orientation between the planes of the osculating orbits of the planets $i$ and $j$. We recall that, in turn, this productory structure is a consequence of the ``tree'' structure Deprit's {\sl kinetic frames}.

\vskip.1in
\noi
We are now ready for the Proof of Proposition \ref{spatial non degeneracy}. For simplicity of notations, during the proof, we  write
$$\eta=(\eta_1,\cdots,\eta_N)\ ,\quad \xi=(\xi_1,\cdots,\xi_N)\ ,\quad p=(p_1,\cdots,p_{N-1})\ ,\quad q=(q_1,\cdots,q_{N-1})$$
for
$$\eta_*=(\eta^*_1,\cdots,\eta^*_N)\ ,\quad \xi_*=(\xi^*_1,\cdots,\xi^*_N)\ ,\quad \bar p=(\bar p_1,\cdots,\bar p_{N-1})\ ,\quad \bar q=(\bar q_1,\cdots,\bar q_{N-1})$$
believing that no confusion arises with the full reduced variables, which are never used in this section. 

\vskip.1in
\noi
{\bf Proof of Proposition \ref{spatial non degeneracy}.}\
We proceed by induction on the number $N$ of planets, starting with $N=2$. We esplicitate  the dependence on $N$ marking ${\eufm F}^{N}$, $\r_h^N$, $\cdots$ the quantities ${\eufm F}$, $\r_h$, $\cdots$ relatively to $N$ planets. 
\vskip.1in
\noi
$\underline{\textrm{\sl Proof for $N=2$.}}$\ When $N=2$, the two matrices ${\cal F}_h^{2}$, ${\cal F}_v^{*2}$ of  the quadratic forms ${\cal Q}_h^2$, ${\cal Q}_v^{*2}$, have order $2$, $1$, respectively, so, their diagonalizing matrices $\r_h^{2}$, $\r_v^{*2}$ can be exactly computed: $\r_v^{*2}=\id$ is trivial and, diagonalizing
\beqano
{\cal F}_h^2&=&\bar m_1\bar m_2\left(
\begin{array}{lrr}
\frac{C_1(a_1,a_2)}{\L_1}&\frac{C_2(a_1,a_2)}{\sqrt{\L_1\L_2}}\\
\frac{C_2(a_1,a_2)}{\sqrt{\L_1,\L_2}}&\frac{C_1(a_1,a_2)}{\L_2}
\end{array}
\right)\nonumber\\
&=&\bar m_1\bar m_2\left(
\begin{array}{lrr}
-\frac{3}{4}\frac{a_1^2}{a_2^3\L_1}+O\left(\frac{a_1^4}{a_2^5\L_1}\right)&O\left(\frac{a_1^3}{a_2^4\sqrt{\L_1\L_2}}\right)\\
O\left(\frac{a_1^3}{a_2^4\sqrt{\L_1\L_2}}\right)&-\frac{3}{4}\frac{a_1^2}{a_2^3\L_2}+O\left(\frac{a_1^4}{a_2^5\L_2}\right)
\end{array}
\right)
\eeqano
we find, 
for $$a_1=O(1)\ ,\quad a_2\to \infty\ ,$$  
\beq{diagonalize Fh2}\r_h^2=\left(
\begin{array}{lrr}
\frac{1}{\sqrt{1+\e^2}}&\frac{\e}{\sqrt{1+\e^2}}\\
-\frac{\e}{\sqrt{1+\e^2}}&\frac{1}{\sqrt{1+\e^2}}
\end{array}
\right)=\id+\left(
\begin{array}{lcr}
O(a_2^{-5/2})&O(a_2^{-5/4})\\
O(a_2^{-5/4})&O(a_2^{-5/2})
\end{array}
\right)
\eeq
where, letting for shortness,
$$a:=\frac{C_1(a_1,a_2)}{\L_1}\ ,\quad b:=\frac{C_1(a_1,a_2)}{\L_2}\ ,\quad c:=\frac{C_2(a_1,a_2)}{\sqrt{\L_1,\L_2}}\ ,$$
then $\e$ denotes
$$\e:=\frac{|a-b|}{2c}\left(\sqrt{1+\frac{4c^2}{(a-b)^2}}-1\right)=O(a_2^{-5/4})\ .$$
\vskip.1in
\noi
Using these expressions, we compute now the asymptitics (in $a_2$) of the Birkhoff matrix associated to the three polynomials ${\eufm F}_h^2(\r_h\eta,\r_h\xi,p,q)$, ${\eufm F}_{hv}^2(\r_h\eta,\r_h\xi,p,q)$, ${\eufm F}_v^2(\r_h\eta,\r_h\xi,p,q)$.
The dominant contribute of ${\eufm F}_h^2(\r_h\eta,\r_h\xi,p,q)$  to the Birkhoff matrix is  (compare with Arnol'd's computation, \cite{Arn63})
\beq{horiz}\left(
\begin{array}{lcr}
\frac{3}{4}\frac{\bar m_1\bar m_2}{\L_1^2a_2}\frac{a_1^2}{a_2^2}&-\frac{9}{4}\frac{\bar m_1\bar m_2}{\L_1\L_2a_2}\frac{a_1^2}{a_2^2}&*\\
-\frac{9}{4}\frac{\bar m_1\bar m_2}{\L_1\L_2a_2}\frac{a_1^2}{a_2^2}&-3\frac{\bar m_1\bar m_2}{\L_2^2a_2}\frac{a_1^2}{a_2^2}&*\\
{*}&*&*\\
\end{array}
\right)
\eeq
We recall  for completeness  the computation, essentially done in the study of the plane problem, which leads to this result, since now we want to expand with respect to the semimajor axes, rather that their ratios. It is a consequence of the expression of the horizontal quartic form
\beqa{F2}
{\eufm F}_h^2&=&q_{1111}\eta_1^4+q_{1112}\eta_1^3\eta_2+q_{1122}\eta_1^2\eta_2^2+q_{2221}\eta_1\eta_2^3+q_{2222}\eta_2^4\nonumber\\
&+&q_{1111}\xi_1^4+q_{1112}\xi_1^3\xi_2+q_{1122}\xi_1^2\xi_2^2+q_{2221}\xi_1\xi_2^3+q_{2222}\xi_2^4\nonumber\\
&+&r_{1111}\eta_1^2\xi_1^2+r_{1112}\eta_1^2\xi_1\xi_2+r_{1122}\eta_1^2\xi_2^2\nonumber\\
&+&r_{1211}\eta_1\eta_1\xi_1^2+r_{1212}\eta_1\eta_2\xi_1\xi_2+r_{2122}\eta_1\eta_2\xi_2^2\nonumber\\
&+&r_{2211}\eta_2^2\xi_1^2+r_{2221}\eta_2^2\xi_1\xi_2+r_{2222}\eta_2^2\xi_2^2\nonumber\\
\eeqa
with the following expansions  (based on the expansions of the Laplace coefficients)  of the involved entries of the tensors $q$, $r$
\beqa{neglect}
&& q_{1111}=\frac{3}{32}\frac{\bar m_1\bar m_2}{\L_1^2}\frac{1}{a_2}\left(\frac{a_1^2}{a_2^2}+O\left(\frac{a_1^4}{a_2^4}\right)\right)\nonumber\\
&& q_{2222}=-\frac{3}{8}\frac{\bar m_1\bar m_2}{\L_2^2}\frac{1}{a_2}\left(\frac{a_1^2}{a_2^2}+O\left(\frac{a_1^4}{a_2^4}\right)\right)\nonumber\\
&& q_{1122}=-\frac{9}{16}\frac{\bar m_1\bar m_2}{\L_1\L_2}\frac{1}{a_2}\left(\frac{a_1^2}{a_2^2}+O\left(\frac{a_1^4}{a_2^4}\right)\right)\nonumber\\
&& r_{1111}=\frac{3}{16}\frac{\bar m_1\bar m_2}{\L_1^2}\frac{1}{a_2}\left(\frac{a_1^2}{a_2^2}+O\left(\frac{a_1^4}{a_2^4}\right)\right)\nonumber\\
&& r_{2222}=-\frac{3}{4}\frac{\bar m_1\bar m_2}{\L_2^2}\frac{1}{a_2}\left(\frac{a_1^2}{a_2^2}+O\left(\frac{a_1^4}{a_2^4}\right)\right)\nonumber\\
&& r_{1122}=-\frac{9}{16}\frac{\bar m_1\bar m_2}{\L_1\L_2}\frac{1}{a_2}\left(\frac{a_1^2}{a_2^2}+O\left(\frac{a_1^4}{a_2^4}\right)\right)\nonumber\\
&& r_{2211}=-\frac{9}{16}\frac{\bar m_1\bar m_2}{\L_1\L_2}\frac{1}{a_2}\left(\frac{a_1^2}{a_2^2}+O\left(\frac{a_1^4}{a_2^4}\right)\right)\nonumber\\
&& q_{1112},\ r_{1211},\ r_{1112}=O\left(\frac{a_1^3}{a_2^4\L_1\sqrt{\L_1\L_2}}\right)\nonumber\\
&& q_{2221},\ r_{2122},\ r_{2221}=O\left(\frac{a_1^3}{a_2^4\L_2\sqrt{\L_1\L_2}}\right)\nonumber\\
&& r_{1212}=O\left(\frac{a_1^4}{a_2^5\L_1\L_2}\right)
\eeqa

\vskip.1in
\noi
We compute now the Birkhoff matrix associated to ${\eufm F}_{hv}^2(\r_h\eta,\r_h\xi,p,q)$. Replacing
\beqa{cal Q}
{\eufm L}p_2-{\eufm L}p_1&=&(\ell_2-\ell_1)p=\bar c_1p=\sqrt{\left(\frac{1}{\L_1}+\frac{1}{\L_2}\right)}p\nonumber\\
{\cal Q}_{12}^{11}\cdot(\eta_1,\eta_2,\xi_1,\xi_2)&=&C_3(a_1,a_2)\frac{\eta_1^2}{\L_1}+C_4(a_1,a_2)\frac{\eta_1\eta_2}{\sqrt{\L_1\L_2}}+C_5(a_1,a_2)\frac{\eta_2^2}{\L_2}\nonumber\\
&+&C_6(a_1,a_2)\frac{\xi_1^2}{\L_1}+C_7(a_1,a_2)\frac{\xi_1\xi_2}{\sqrt{\L_1\L_2}}+C_8(a_1,a_2)\frac{\xi_2^2}{\L_2}\nonumber\\
{\cal Q}_{12}^{22}\cdot(\eta_1,\eta_2,\xi_1,\xi_2)&=&C_6(a_1,a_2)\frac{\eta_1^2}{\L_1}+C_7(a_1,a_2)\frac{\eta_1\eta_2}{\sqrt{\L_1\L_2}}+C_8(a_1,a_2)\frac{\eta_2^2}{\L_2}\nonumber\\
&+&C_3(a_1,a_2)\frac{\xi_1^2}{\L_1}+C_4(a_1,a_2)\frac{\xi_1\xi_2}{\sqrt{\L_1\L_2}}+C_5(a_1,a_2)\frac{\xi_2^2}{\L_2}\nonumber\\
\eeqa
the polynomial ${\eufm F}_{hv}^2(\eta,\xi,p,q)$ (compare eq. (\ref{vertical quadratic}))
reduces to
\beqa{Fhvcontr}
{\eufm F}_{hv}^2&=&\frac{1}{2}\bar m_1\bar m_2\Big({\cal Q}_{12}^{11}\,({\eufm L}q_2-{\eufm L}q_1)^2+{\cal Q}_{12}^{22}\,({\eufm L}p_2-{\eufm L}p_1)^2\Big)+\hat{\eufm F}_{hv}^2\nonumber\\
&=&\frac{1}{2}\bar m_1\bar m_2\left(\frac{1}{\L_1}+\frac{1}{\L_2}\right)\nonumber\\
&\times&\left[C_3(a_1,a_2)\frac{\eta_1^2q^2+\xi_1^2p^2}{\L_1}+C_4(a_1,a_2)\frac{\eta_1\eta_2q^2+\xi_1\xi_2p^2}{\sqrt{\L_1\L_2}}+C_5(a_1,a_2)\frac{\eta_2^2q^2+\xi_2^2p^2}{\L_2}\right.\nonumber\\
&+&C_6(a_1,a_2)\frac{\xi_1^2q^2+\eta_1^2p^2}{\L_1}+C_7(a_1,a_2)\frac{\xi_1\xi_2q^2+\eta_1\eta_2p^2}{\sqrt{\L_1\L_2}}+C_8(a_1,a_2)\frac{\xi_2^2q^2+\eta_2^2p^2}{\L_2}\nonumber\\
&+&\hat{\eufm F}_{hv}^2
\eeqa
(referring to (\ref{ccc}), (\ref{Cij}), (\ref{C13C14}) for the definition of $C_1$, $C_3\div C_8$, $C_{13}$)
where 
\beqano
\hat{\eufm F}_{hv}^2&:=&\frac{1}{2}\bar m_1\bar m_2\Big({\cal Q}_{12}^{12}+{\cal Q}_{12}^{21})\Big({\eufm L}p_2-{\eufm L}p_2\Big)\Big({\eufm L}q_j-{\eufm L}q_i\Big)\nonumber\\
\eeqano
This term gives no contribute to the Birkhoff matrix, since, when computed in $(\r_h^2\eta, \r_h^2\xi,p,q)$, it contains only monomials of the form $\eta_i\xi_jpq$, hence, with vanishing Birkhoff form.

\vskip.1in
\noi
Then, from  (\ref{Fhvcontr}), ${\eufm F}_{hv}^2(\r_h\eta,\r_h\xi,p,q)$ generates on the entries with place $(1,3)$, $(2,3)$, $(3,1)$, $(3,2)$ of the Birkhoff matrix the dominant terms (in $a_1$) 
\beq{horiz vert}\left(
\begin{array}{lcr}
{*}&*&\frac{9}{4}\frac{\bar m_1\bar m_2}{\L_1^2a_2}\frac{a_1^2}{a_2^2}\\
{*}&{*}&\frac{9}{4}\frac{\bar m_1\bar m_2}{\L_1\L_2a_2}\frac{a_1^2}{a_2^2}\\
\frac{9}{4}\frac{\bar m_1\bar m_2}{\L_1^2a_2}\frac{a_1^2}{a_2^2}&\frac{9}{4}\frac{\bar m_1\bar m_2}{\L_1\L_2a_2}\frac{a_1^2}{a_2^2}&*\\
\end{array}
\right)
\eeq
The result is found  taking into account the diagonalizing matrices as in (\ref{diagonalize Fh2})  and using the  expansions  for the $C_3\div C_8$
\beqano
&& C_3=\frac{3}{a_2}\left(\frac{a_1^2}{a_1^2}+O\left(\frac{a_1^4}{a_2^4}\right)\right)\nonumber\\
&& C_5=\frac{9}{8a_2}\left(\frac{a_1^2}{a_1^2}+O\left(\frac{a_1^4}{a_2^4}\right)\right)\nonumber\\
&& C_6=-\frac{3}{4a_2}\left(\frac{a_1^2}{a_1^2}+O\left(\frac{a_1^4}{a_2^4}\right)\right)\nonumber\\
&& C_8=\frac{9}{8a_2}\left(\frac{a_1^2}{a_1^2}+O\left(\frac{a_1^4}{a_2^4}\right)\right)\nonumber\\
&& C_4,\ C_7=O\left(\frac{a_1^3}{a_2^4}\right)\nonumber\\
\eeqano

\vskip.1in
\noi
We finally compute now the Birkhoff matrix associated to   (compare eq. (\ref{vertical quadratic}))
\beqano
{\eufm F}_v^2(\eta,\xi,p,q)&=&{\bar m_1\bar m_2}({\eufm Q}^{12}_{11}+{\eufm Q}^{12}_{22})C_1(a_1,a_2)\nonumber\\
&+&\bar m_1\bar m_2\left(\Big({\eufm L}p_2-{\eufm L}p_1\Big)^2+\Big({\eufm L}q_2-{\eufm L}q_1\Big)^2\right)^2C_{13}(a_1,a_2)\Big)\ .\nonumber\\
\eeqano
We recall that in the case $N=2$, ${\eufm S}_{12}=\bar{\cal S}_1$, so,   ${\eufm Q}_{12}^{11}$, ${\eufm Q}_{12}^{22}$ coincide with the order $4$ terms of the antries at places $(1,1)$, $(2,2)$ of $\bar{\cal S}_1$, which are (compare (\ref{Sigma order 4})$\div$(\ref{R0}) below)
\beqano{\eufm Q}_{12}^{11}&=&-q^2\bar{\cal C}_1=-q^2\frac{2\L_{2}^2(\frac{\eta_1^2+\xi_1^2}{2})+2\L_1^2\frac{\eta_2^2+\xi_2^2}{2}-\L_1\L_{2}\frac{p^2+q^2}{2}}{4\L_1^2\L_{2}^2}\nonumber\\
{\eufm Q}_{12}^{22}&=&-p^2\bar{\cal C}_1=-p^2\frac{2\L_{2}^2(\frac{\eta_1^2+\xi_1^2}{2})+2\L_1^2\frac{\eta_2^2+\xi_2^2}{2}-\L_1\L_{2}\frac{p^2+q^2}{2}}{4\L_1^2\L_{2}^2}\nonumber\\
\eeqano
Hence, 
\beqa{Fv}
{\eufm F}_v^2(\eta,\xi,p,q)&=&{\bar m_1\bar m_2}({\eufm Q}^{12}_{11}+{\eufm Q}^{12}_{22})C_1(a_1,a_2)\nonumber\\
&+&\bar m_1\bar m_2\left(\Big({\eufm L}p_2-{\eufm L}p_1\Big)^2+\Big({\eufm L}q_2-{\eufm L}q_1\Big)^2\right)^2C_{13}(a_1,a_2)\Big)\nonumber\\
&=&-\frac{\bar m_1\bar m_2}{2}C_1(a_1,a_2)\nonumber\\
&\times&\left(\frac{2}{\L_1^2}\frac{\eta_1^2+\xi_1^2}{2}\frac{p^2+q^2}{2}+\frac{2}{\L_2^2}\frac{\eta_2^2+\xi_2^2}{2}\frac{p^2+q^2}{2}-\frac{1}{\L_1\L_2}\left(\frac{p^2+q^2}{2}\right)^2\right)\nonumber\\
&+&\bar m_1\bar m_2\left(\frac{1}{\L_1}+\frac{1}{\L_2}\right)^2(p^2+q^2)^2C_{13}(a_1,a_2)
\eeqa
Since (\ref{Fv}), letting $(\eta,\xi,p,q)\to(\r_h\eta,\r_h\xi,p,q)$ with $\r_h$ as in (\ref{diagonalize Fh2}) and using the expansions $$C_1=-\frac{3}{4}\frac{a_1^2}{a_2^3}+O\left(\frac{a_1^4}{a_2^5}\right)\ ,\quad C_{13}=-\frac{3}{32}\frac{a_1^2}{a_2^3}+O\left(\frac{a_1^4}{a_2^5}\right)\ ,$$
we find the contribute of the Birkhoff matrix associated to ${\eufm F}_{v}^2(\r_h\eta,\r_h\xi,p,q)$
\beq{vert}\left(
\begin{array}{lcr}
{*}&*&\frac{3}{4}\frac{\bar m_1\bar m_2}{\L_1^2a_2}\frac{a_1^2}{a_2^2}\\
{*}&{*}&O\left(\frac{1}{\L_2^2a_2}\frac{a_1^2}{a_2^2}\right)\\
\frac{3}{4}\frac{\bar m_1\bar m_2}{\L_1^2a_2}\frac{a_1^2}{a_2^2}&O\left(\frac{1}{\L_2^2a_2}\frac{a_1^2}{a_2^2}\right)&-\frac{3}{4}\frac{\bar m_1\bar m_2}{\L_1^2a_2}\frac{a_1^2}{a_2^2}\\
\end{array}
\right)
\eeq
Finally, on count of (\ref{horiz}), (\ref{horiz vert}), (\ref{vert}), we find that the Birkhoff matrix associated to ${\eufm F}^2(\r_h\eta,\r_h\xi,p,q)$ is
\beqano
\frac{\bar m_1\bar m_2}{a_2}\frac{a_1^2}{a_2^2}\left(
\begin{array}{lcr}
\frac{3}{4\L_1^2}&-\frac{9}{4\L_1\L_2}&\frac{3}{\L_1^2}\\
-\frac{9}{4\L_1\L_2}&-\frac{3}{\L_2^2}&\frac{9}{4\L_1\L_2}\\
\frac{3}{\L_1^2}&\frac{9}{4\L_1\L_2}&-\frac{3}{4\L_1^2}\\
\end{array}
\right)(1+o(1))
\eeqano
hence, it is non singular, having determinant
$$-\left(\frac{\bar m_1\bar m_2}{a_2}\frac{a_1^2}{a_2^2}\right)^3\frac{27}{16\L_1^4\L_2^2}(1+o(1))\neq 0\ ,$$
and the basis of induction is proved.
\vskip.1in
\noi
For the proof of the inductive step, we need the following result, due to J. Fej\'oz, to whom we refer for the proof.
\begin{lemma}\label{Fejoz lemma}({\sl J. Fej\'oz, \cite{Fej04}, {\scshape corollaire} $72$ })\\
Let $\d_1$, $\cdots$, $\d_{n-1}\in \real$, $\d_n=0$ such that $\dst \s:=\min_{1\leq j\neq k\leq n}|\d_j-\d_k|\neq 0$, $\hat D\in {\rm Matr}_{n-1}(\real)$ a symmetric with eigenvalues $\d_1$, $\cdots$, $\d_{n-1}$, $D$ the symmetric matrix
$$D=\left(\begin{array}{lrr}
\hat D&0\\
0&0
\end{array}
\right)
$$
and $A_\e\in{\rm Matr}_{n-1}(\real)$ a symmetric matrix with last coefficient 
$$(A_\e)_{nn}=c_1+c_2\e^\b\ ,\quad c_1,\ c_2\in \real\ ,\quad 0\leq \b<2\ .$$
Then, when $\e\to 0$, the matrix $D+\e\,A$ has an eigenvalue
$$\s_n(\e)=\e(c_1+c_2\e^\b)+O(\e^2)\ .$$
Furthermore, if $\hat D$ is diagonal, $D+\e A$ is coniugated to a diagonal matrix through a matrix $\r\in SO_n(\real)$ verifying $\dst \r=I+O(\e)$.
\end{lemma}
We can apply the previous lemma taking ${\cal F}_h^{N-1}$, ${\cal F}_v^{*(N-1)}$ for $\hat D$ and ${\cal F}_h^N$, ${\cal F}^{*N}$ for $D+\e\,A$, with $\e=a_N^{-3}$. Observe in fact that both ${\cal F}_h^{N-1}$ ${\cal F}_v^{*(N-1)}$ verify the assumptions of the Lemma, since their respective eigenvalues do not satisfy any other linear relation than the Herman's resonance (Proposition \ref{true non resonance} above) and, furthermore,
$${\cal F}_h^N=\left(\begin{array}{lrr}
{\cal F}_h^{N-1}&0\\
0&0
\end{array}
\right)+O(a_N^{-7/2})\ ,\quad{\cal F}_v^{*N}=\left(\begin{array}{lrr}
{\cal F}_v^{*(N-1)}&0\\
0&0
\end{array}
\right)+O(a_N^{-3})\quad (a_N\to \infty)
$$
Then, the diagonalizing matrices  $\r_h^{N}$, $\r_v^{*N}$ at step $N$ are related to the matrices $\r_h^{N-1}$, $\r_h^{*(N-1)}$ of step $N-1$ by
\beqa{rhohNrhovN}\r_h^{N}=\left(\begin{array}{lrr}
\r_h^{N-1}&0\\
0&1
\end{array}
\right)+O(a_N^{-7/2})\ ,\quad \r_v^{*N}=\left(\begin{array}{lrr}
\r_v^{*(N-1)}&0\\
0&1
\end{array}
\right)+O(a_N^{-3})\ .\eeqa
This result will be used in the inductive step, which we are now ready to prove.

\vskip.1in
\noi
$\underline{\textrm{\sl Proof of the inductive step $(N-1)\to N$.}}$\ Assume that, when
$$a_1\ll\cdots\ll a_{N-2}\ll a_{N-1}\to \infty$$ the Birkhoff matrix $A_{\rm plt}^{N-1}$ associated to $${\eufm F}^{N-1}(\r_h^{N-1}\hat\eta,\r_h^{N-1}\hat\xi,\r_v^{*(N-1)}\hat p,\r_v^{*(N-1)}\hat q)$$ is non degenerate, and let us prove that, when 
$$a_1\ll\cdots\ll a_{N-1}\ll a_N\to \infty$$
 then, the matrix $A_{\rm plt}^N$ associated to $${\eufm F}^{N}(\r_h^{N}\eta,\r_h^{N}\xi,\r_v^{*N} p,\r_v^{*N} q)$$ is so.

\vskip.1in
\noi
Let ${\eufm F}:={\eufm F}^{N}$ as in (\ref{eufm F}) $\div$(\ref{vertical quadratic}) and let us split
\beq{slitting}{\eufm F}:={\eufm F}^{N}=\hat{\eufm F}^{N-1}+{\eufm F}'\eeq
where $$\dst  \hat{\eufm F}^{N-1}=\hat{\eufm F}^{N-1}_h+\hat{\eufm F}^{N-1}_{hv}+\hat{\eufm F}^{N-1}_v$$ is $4$--order polynomial associated to $N-1$ bodies, in the variables the variables $(\hat\eta,\hat\xi, \hat p,\hat q)$, when the variables $(\eta,\xi,p,q)$ related to $N$ bodies are written as
$$\eta:=(\hat\eta,\eta_N):=\Big((\eta_1,\cdots,\eta_{N-1}),\eta_N\Big)\ ,\quad p:=(\hat p,p_{N-1}):=\Big((p_1,\cdots,p_{N-2}),p_{N-1}\Big)\ ,\cdots$$
and
$${\eufm F}':={\eufm F}^{N}-\hat{\eufm F}^{N-1}={\eufm F}'_h+{\eufm F}'_{hv}+{\eufm F}'_v\ ,$$
with $\dst {\eufm F}'_h={\eufm F}_h-\hat{\eufm F}^{N-1}_h$, $\cdots$, and similarly for the definitions of ${\eufm F}'_{hv}$, ${\eufm F}'_v$.
By inspection of its coefficients, ${\eufm F}'_h$ is $\dst O(a_N^{-3})$, 
so, making  use of  (\ref{rhohNrhovN}), it is not difficult to see that (\ref{slitting}) implies that the first approximation for the Birkhoff matrix  $A_{\rm plt}^N$ (suitably rearranged) is
\beq{Birkh matr}A_{\rm plt}^N=\left(
\begin{array}{lrr}
A_{\rm plt}^{N-1}+O(a_N^{-3})&O(a_N^{-3})\\
O(a_N^{-3})&\hat A+O(a_N^{-6})
\end{array}
\right)\eeq
where $\hat A=O(a_N^{-3})$ is the square matrix of order $2$ associated to the quadratic form in  the couple of variables $$\dst \left(\frac{\eta_N^2+\xi_N^2}{2},\ \frac{p_{N-1}^2+q_{N-1}^2}{2}\right)$$ appearing the Birkhoff polynomial of $${\eufm F}'\Big((\r_h^{N-1}\hat\eta,\eta_N),(\r_h^{N-1}\hat\xi,\xi_N)(\r_v^{*(N-1)}\hat p,p_{N-1})(\r_v^{*(N-1)}\hat q,q_{N-1})\Big)$$ By (\ref{Birkh matr}) and the inductive hypothesis, we only need to prove that $\hat A$ is non singular.  As for the proof of the iductive basis, this is done by direct computation. We start with the horizontal  terms 
\beqano
{\eufm F}'_h&=&{\eufm F}_h^N-\hat{\eufm F}_h^{N-1}\nonumber\\
&=&q_{NNNN}\eta_N^4+q_{NNNN}\xi_N^4+r_{NNNN}\eta_N^2\xi_N^2\nonumber\\
&+&\sum_{1\leq i\leq N-1}\Big(q_{iiiN}\eta_i^3\eta_N+q_{iiNN}\eta_i^2\eta_N^2+q_{NNNi}\eta_i\eta_N^3\nonumber\\
&+&q_{iiiN}\xi_i^3\xi_N+q_{iiNN}\xi_i^2\xi_N^2+q_{NNNi}\xi_i\xi_N^3\nonumber\\
&+&r_{iiiN}\eta_i^2\xi_i\xi_N+r_{iiNN}\eta_i^2\xi_N^2\nonumber\\
&+&r_{iNii}\eta_i\eta_i\xi_i^2+r_{iNiN}\eta_i\eta_N\xi_i\xi_N+r_{NiNN}\eta_i\eta_N\xi_N^2\nonumber\\
&+&r_{NNii}\eta_N^2\xi_i^2+r_{NNNi}\eta_N^2\xi_i\xi_N\Big)\nonumber\\
\eeqano
By the previous discussion, we have to pick the monomials in $\eta_N$, $\xi_N$ only, namely
$$q_{NNNN}\eta_N^4+q_{NNNN}\xi_N^4+r_{NNNN}\eta_N^2\xi_N^2$$
and then we have to compute the related Birkhoff polynomial. Proceeding as done in the previous step, we use the expansions
$$q_{NNNN}=\sum_{1\leq i\leq N-1}\left(-\frac{3}{8}\frac{\bar m_i\bar m_N}{\L_N^2a_N}\frac{a_i^2}{a_N^2}+O\left(\frac{a_i^2}{a_N^3}\right)\right)=-\frac{3}{8}\frac{\bar m_{N-1}\bar m_N}{\L_N^2a_N}\frac{a_{N-1}^2}{a_N^2}+O\left(\frac{a_{N-1}^2}{a_N^3}\right)$$ 
and, similarly,
$$r_{NNNN}=\sum_{1\leq i\leq N-1}\left(-\frac{3}{16}\frac{\bar m_i\bar m_N}{\L_N^2a_N}\frac{a_i^2}{a_N^2}+O\left(\frac{a_i^2}{a_N^3}\right)\right)=-\frac{3}{16}\frac{\bar m_{N-1}\bar m_N}{\L_N^2a_N}\frac{a_{N-1}^2}{a_N^2}+O\left(\frac{a_{N-1}^2}{a_N^3}\right)$$
At this point the computation is just the same we have seen for $N=2$ and, at the place $(1,1)$ of $\hat A$ we will find
$$\hat A=\left(
\begin{array}{lrr}
-3\frac{\bar m_{N-1}\bar m_N}{\L_N^2a_N}\frac{a_{N-1}^2}{a_N^2}&{*}\\
{*}&{*}\\
\end{array}
\right)$$
Also the analysis of the vertical part ${\eufm F}'_{hv}$ 
will give, as dominant terms on the off diagonal entries, a similar result as in the cae $N=2$ 
\beq{off diag}\hat A=\left(
\begin{array}{lrr}
{*}&\frac{9}{4}\frac{\bar m_{N-1}\bar m_N}{\L_{N-1}\L_Na_N}\frac{a_{N-1}^2}{a_N^2}\\
\frac{9}{4}\frac{\bar m_{N-1}\bar m_N}{\L_{N-1}\L_Na_N}\frac{a_{N-1}^2}{a_N^2}&{*}\\
\end{array}
\right)
\eeq
This result follows computing the  Birkhoff matrix relatively to the the monomials with $p_{N-1}$, $q_{N-1}$, $\eta_N$, $\xi_N$ of
\beqano
{\eufm F}'_{hv}&:=&{\eufm F}_{hv}^N-\hat{\eufm F}_{hv}^{N-1}\nonumber\\
&=&\frac{1}{2}\sum_{1\leq i\leq N-1}\bar m_i\bar m_N\Big({\cal Q}_{iN}^{11}\,({\eufm L} q_N-{\eufm L} q_i)^2+{\cal Q}_{iN}^{22}\,({\eufm L} p_N-{\eufm L} p_i)^2\Big)\nonumber\\
&+&\hat{\eufm F}_{hv}^N
\eeqano
(where $\hat{\eufm F}^N$ is a suitable polynomial which gives no contribute to the Birkhoff matrix), noticing that ${\eufm L} q_N-{\eufm L} q_i$ has the form 
$${\eufm L} q_N-{\eufm L} q_i=\bar c_{N-1}p_{N-1}+\hat c_i\cdot\hat p\ ,\quad \hat p=(p_1,\cdots,p_{N-2})$$
with
$$\bar c_{N-1}=\sqrt{\frac{1}{{\rm L}_{N-1}}+\frac{1}{\L_N}}=\sqrt{\frac{1}{\L_{N-1}}}+O\left(\frac{1}{\L_N}\right)+O\left(\L_{N-2}\right)$$
and ${\cal Q}_{iN}^{kk}$ as in (\ref{cal Q}), with $i$, $N$ replacing $1$, $2$ respectively.

\vskip.1in
\noi
The last step consists in evaluating the contribute to $\hat A$ of
\beqa{FvN}
{\eufm F}_v'&:=&{\eufm F}_{v}^N-\hat{\eufm F}_{v}^{N-1}\nonumber\\\nonumber\\
&=&\sum_{1\leq i\leq N-1}{\bar m_i\bar m_N}({\eufm Q}^{iN}_{11}+{\eufm Q}^{iN}_{22})C_1(a_i,a_N)\nonumber\\
&+&\sum_{1\leq i\leq N-1}{\bar m_i\bar m_N}\left(\Big({\eufm L} p_i-{\eufm L} p_N\Big)^2+\Big({\eufm L} q_i-{\eufm L} q_N\Big)^2\right)^2C_{13}(a_i,a_N)\nonumber\\
&+&\hat{\eufm F}_v^N\nonumber\\
\eeqa
where
$$\hat{\eufm F}_v^N:=\sum_{1\leq i\leq N-1}{\bar m_i\bar m_N}\left({\sum}^*_{1\leq h<k\leq N-1}\ell^{iN}_h\ell^{iN}_k( p_{h} q_k- p_k q_h)\right)^2C_{14}(a_i,a_N)$$
gives only a negligible contribute on $\hat A$.
Notice that matrices ${\cal S}_i$, $\tilde{\cal S}_j$ appearing in the the productories (\ref{eufmSij}) do not involve the variables $p_{N-1}$, $q_{N-1}$, $\eta_N$, $\xi_N$, so,  the monomials of degree $4$ in ${\eufm Q}^{iN}_{11}$, ${\eufm Q}^{iN}_{22}$ involving only $p_{N-1}$, $q_{N-1}$, $\eta_N$, $\xi_N$coincide with the corresponding monomials of degree  $4$ of $\bar{\cal S}_{N-1}$, on the entries $(1,1)$, $(2,2)$, which are
$$-q_{N-1}^2\frac{2{\rm L}_{N-1}^2\t_{N}-{\rm L}_{N-1}\L_{N}\r_{N-1}}{4{\rm L}_{N-1}^2\L_{N}^2}\ ,\ -p_{N-1}^2\frac{2{\rm L}_{N-1}^2\t_{N}-{\rm L}_{N-1}\L_{N}\r_{N-1}}{4{\rm L}_{N-1}^2\L_{N}^2}$$
Thus, the term
$$\sum_{1\leq i\leq N-1}{\bar m_i\bar m_N}({\eufm Q}^{iN}_{11}+{\eufm Q}^{iN}_{22})C_1(a_i,a_N)$$
gives 
$$\left(
\begin{array}{lrr}
{*}&O\left(\frac{1}{\L_{N}^2a_{N}}\frac{a_{N-1}^2}{a_{N}^2}\right)\\
O\left(\frac{1}{\L_{N}^2a_{N}}\frac{a_{N-1}^2}{a_{N}^2}\right)&O\left(\frac{1}{\L_{N}\L_{N-1}a_{N}}\frac{a_{N-1}^2}{a_{N}^2}\right)\\
\end{array}
\right)$$
As in the case $N=2$, the off--diagonal terms above can be neglected with respect to the ones appearing in (\ref{off diag}), while the diagonal term with place  $(2,2)$ can be neglected with respect to the corresponding  term
$$\hat A=\left(
\begin{array}{lrr}
{*}&{*}\\
{*}&-\frac{3}{4}\frac{\bar m_{N-1}\bar m_{N}}{\L_{N-1}^2a_{N}}\frac{a_{N-1}^2}{a_{N}^2}\\
\end{array}
\right)$$
generated by the second line in (\ref{FvN}). On count of the previous computations, the final result extends the one found for $N=2$, giving
$$\hat A=\left(
\begin{array}{lrr}
-3\frac{\bar m_{N-1}\bar m_N}{\L_N^2a_N}\frac{a_{N-1}^2}{a_N^2}&\frac{9}{4}\frac{\bar m_{N-1}\bar m_N}{\L_{N-1}\L_Na_N}\frac{a_{N-1}^2}{a_N^2}\\
\frac{9}{4}\frac{\bar m_{N-1}\bar m_N}{\L_{N-1}\L_Na_N}\frac{a_{N-1}^2}{a_N^2}&-\frac{3}{4}\frac{\bar m_{N-1}\bar m_{N}}{\L_{N-1}^2a_{N}}\frac{a_{N-1}^2}{a_{N}^2}\\
\end{array}
\right)$$
so, 
$$\textrm{\rm det}\hat A=-\frac{45}{16}\left(\frac{\bar m_{N-1}\bar m_N}{\L_{N-1}\L_Na_N}\frac{a_{N-1}^2}{a_N^2}\right)^2\neq 0$$
which finishes the proof.
\subsubsection{Expansion of the Hamiltonian}\label{Expansion of the Hamiltonian space}
\begin{lemma}
The ``secular perturbation'', \ie, the mean
$$\bar f_{\rm  plt,pr}:=\frac{1}{(2\p)^N}\int_{\torus^N}\,f_{\rm  plt,pr}\qquad f_{\rm  plt,pr}:=\sum_{1\leq i<j\leq N}\frac{y_i\cdot y_j}{\bar m_0}-\frac{\bar m_i \bar m_j}{|x_i-x_j|}$$
($y_i$, $x_i$ as in (\ref{xi})$\div$(\ref{eufm R pr})) coincides with the mean of the Newtonian potential:
\beq{Newtonian potential}\bar f_{\rm  plt,pr}=-\sum_{1\leq i<j\leq N}\frac{\bar m_i \bar m_j}{4\p^2}\int_{\torus^2}\frac{d\bar\l_id\bar\l_j}{|x_i-x_j|}\ .\eeq
\end{lemma}
The proof of this lemma is trivial and so it is omitted.

\begin{lemma}\label{bar f is even}
The secular perturbation $\bar f_{\rm  plt,pr}$ is an  even function of the ``secular variables'' $\bar z:=(\bar\eta,\bar\xi,\bar p,\bar q)$.
\end{lemma}
{\bf Proof.}\ Actually, it is even in $(\bar\eta,\bar\xi)$ and $(\bar p,\bar q)$ separately. In fact, letting $(\bar\eta,\bar\xi)\to -(\bar\eta,\bar\xi)$ and simultaneously $\bar\l\to\bar\l+\p$, the Plane Delaunay--Poincar\'e Map  changes for a sign  and the matrices ${\eufm R}_i^{\rm pr}$ do not change (they are even in $(\bar\eta,\bar\xi)$), so, taking the mean over $\bar\l$ of the Newtonian potential, we find that $\bar f_{\rm  plt,pr}$ is even in $(\bar\eta,\bar\xi)$. Observing that  $\bar f_{\rm plt,pr}$ depends on $(\bar p,\bar q)$ only through on the entries of ${\eufm R}_i^{\rm pr}$ with place $(1,1)$, $(1,2)$, $(2,1)$, $(2,2)$, which are even in $(\bar p,\bar q)$ (since they are products of the matrices ${\cal S}_i$, $\tilde{\cal S}_j$'s, the entries of which  with place $(1,1)$, $(1,2)$, $(2,1)$, $(2,2)$, $(3,3)$ are even in $(\bar p,\bar q)$, while the ones with place $(1,3)$, $(2,3)$, $(3,1)$, $(3,2)$ are odd), we also find that $\bar f_{\rm  plt,pr}$ is even in $(\bar p,\bar q)$.

\vskip.1in
\noi
We proceed with the expansion of the mean of the Newtonian potential (\ref{Newtonian potential}).

\vskip.1in
\noi
Unsing (\ref{eufm R pr}), we write the  Euclidean distance $|x_i-x_j|$ as 
\beqano
|x_i-x_j|=|{\eufm R}_i^{\rm pr}\hat x_i-{\eufm R}_j^{\rm pr}\hat x_j|=|{\rm I}_i\hat x_i-{\eufm S}_{ij}\hat x_j|=|{\rm I}_i\hat x_i-\hat{\eufm S}_{ij}\hat x_j|\ ,
\eeqano
where $\hat{\eufm S}$ denotes the submatrix of order $2$ of a give matrix ${\eufm S}$ of order $3$, 
$${\rm I}_i:=\arr{{\rm I}_{\rm z}\quad i=1\\
\id\quad i>1
}$$
changes the sign of $\hat x_1$ and ${\eufm S}_{ij}:={\eufm S}_i^T{\eufm S}_j$, with
\beqa{eufmSi}
\quad {\eufm S}_i&:=&\arr{
\left(\prod_{j=1}^{N-2}(\tilde{\cal S}^{\rm pr})_{N-j}^{\rm T}\right){{\cal S}_1^{\rm pr}}^T\quad  i=1\\
\left(\prod_{j=1}^{N-i}(\tilde{\cal S}^{\rm pr}_{N-j})^{\rm T}\right){\cal S}^{\rm pr}_i\quad 2\leq i\leq N-1\\
{\cal S}^{\rm pr}_N\quad i=N\\
}\nonumber\\
\eeqa
Changing,  into the integral (\ref{Newt pot}), the integration variable $\bar\l_1$ with $\bar\l_1+\p$ and making use of the relation $$\hat x_1(\L_1,\bar\l_1+\p,\bar\eta_1,\bar\xi_1)=-\hat x_1(\L_1,\bar\l_1,-\bar\eta_1,-\bar\xi_1)=:-\hat x_1^*(\L_1,\bar\l_1,\bar\eta_1,\bar\xi_1)=:-\hat x_1(\L_1,\bar\l_1,\bar\eta_1^*,\bar\xi_1^*)$$
where $\eta_*:=(-\eta_1,\eta_2,\cdots)$,
we can write the secular perturbation as
\beqa{Newt pot}
\bar f_{\rm  plt,pr}&=&-\sum_{1\leq i<j\leq N}\frac{\bar m_i\bar m_j}{4\p^2}\int_{\torus^2}\frac{d\bar\l_id\bar\l_j}{|x_i- x_j|}\nonumber\\
&=&-\sum_{1\leq i<j\leq N}\frac{\bar m_i\bar m_j}{4\p^2}\int_{\torus^2}\frac{d\bar\l_id\bar\l_j}{|{\rm I}_i\hat x_i-\hat{\eufm S}_{ij}\hat x_j|}\nonumber\\
&=&-\sum_{1\leq i<j\leq N}\frac{\bar m_i\bar m_j}{4\p^2}\int_{\torus^2}\frac{d\bar\l_id\bar\l_j}{|\hat x_i^*-\hat{\eufm S}_{ij}\hat x_j|}
\eeqa
where
$$\hat x_i^*=\arr{\hat x_1^*\quad i=1\\
\hat x_i\quad i>1\ .}$$
Let us now think that $\hat{\eufm S}_{ij}$ is expanded in powers of $(\bar\eta,\bar\xi,\bar p,\bar q)$, up to order four, that is, let us put
$$\hat{\eufm S}_{ij}=\id+\hat{\eufm S}_{ij}^2+\hat{\eufm S}_{ij}^4+{\rm O}(6)$$
and let us consequentely expand the square distance 
\beqano
D_{ij}&:=&|\hat x_i^*-\hat{\eufm S}_{ij}\hat x_j|^2\nonumber\\
&=&|\hat x_i^*|^2+|\hat x_j^*|^2-2\hat x_i^*\cdot\hat{\eufm S}_{ij}\hat x_j\nonumber\\
&=&|\hat x_i^*-\hat x_j|^2-2\,\hat x_i^*\cdot\hat{\eufm S}_{ij}^2\hat x_j-2\,\hat x_i^*\cdot\hat{\eufm S}_{ij}^4\hat x_j+{\rm O}(6)\nonumber\\
&=&D_0^{ij}+D_2^{ij}+D_4^{ij}+{\rm O}(6)
\eeqano
where
\beqa{Dijk}
\arr{
D_0^{ij}:=|\hat x_i^*-\hat x_j|^2\\
D_2^{ij}:=-2\,\hat x_i^*\cdot\hat{\eufm S}_{ij}^2\hat x_j\\
D_4^{ij}:=-2\,\hat x_i^*\cdot\hat{\eufm S}_{ij}^4\hat x_j
}
\eeqa
Using now the elementary expansion
\beqano
\frac{1}{\sqrt{D}}&=&\frac{1}{{D_0}^{1/2}}-\frac{D_2}{2D_0^{3/2}}-\frac{D_4}{2D_0^{3/2}}+\frac{3}{8}\frac{D_2^2}{D_0^{5/2}}+O(6)
\eeqano
when $D$ has the expansion
$$D=D_0+D_2+D_4+{\rm O}(6)$$
with $D_k:=D_k^{ij}$ as in (\ref{Dijk}), we write
\beqano
\frac{1}{|\hat x_i^*-\hat{\eufm S}_{ij}\hat x_j|}&=&\frac{1}{|\hat x_i^*-\hat x_j|}+\frac{\hat x_i^*\cdot\hat{\eufm S}_{ij}^2\hat x_j}{|\hat x_i^*-\hat x_j|^3}+\frac{\hat x_i^*\cdot\hat{\eufm S}_{ij}^4\hat x_j}{|\hat x_i^*-\hat x_j|^3}+\frac{3}{2}\frac{(\hat x_i^*\cdot\hat{\eufm S}_{ij}^2\hat x_j)^2}{|\hat x_i^*-\hat x_j|^5}+O(6)
\eeqano
When we multiply by $-\bar m_i\bar m_j$, sum over all $1\leq i<j\leq N$, take the mean on $(\bar\l_i,\bar\l_j)$, we can split $\bar f_{\rm  plt,pr}$ as
\beq{split bar fBD}\bar f_{\rm  plt,pr}=\bar f_{\rm pl}^*+\bar f_{\rm two}+\bar f_{\rm four}\eeq
where
\beqano
\bar f_{\rm pl}^*&:=&-\sum_{1\leq i<j\leq N}\frac{\bar m_i\bar m_j}{4\p^2}\int_{\torus^2}\frac{d\bar\l_id\bar\l_j}{|\hat x_i^*-\hat x_j|}\nonumber\\
\bar f_{\rm two}&:=&-\sum_{1\leq i<j\leq N}\frac{\bar m_i\bar m_j}{4\p^2}\int_{\torus^2}\frac{\hat x_i^*\cdot\hat{\eufm S}_{ij}^2\hat x_j}{|\hat x_i^*-\hat x_j|^3}d\bar\l_id\bar\l_j\nonumber\\
\bar f_{\rm four}&:=&-\sum_{1\leq i<j\leq N}\frac{\bar m_i\bar m_j}{4\p^2}\int_{\torus^2}\left(\frac{\hat x_i^*\cdot\hat{\eufm S}_{ij}^4\hat x_j}{|\hat x_i^*-\hat x_j|^3}+\frac{3}{2}\frac{(\hat x_i^*\cdot\hat{\eufm S}_{ij}^2\hat x_j)^2}{|\hat x_i^*-\hat x_j|^5}\right)d\bar\l_id\bar\l_j\nonumber\\
\eeqano
and we have now to expand the $\hat x$--coordinate of the plane Delaunay--Poincar\'e map $(\L$, $\bar\l$, $\bar\eta$, $\bar\xi)$ $\to$ $(\hat y$, $\hat x)$ in powers of $(\eta,\xi)$. In the following, we perform this expansion, collecting only the terms of order $2$, $4$. 
\begin{itemize}
\item[(i)] {\bf Expansion of $\bar f_{\rm pl}^*$.} The function  $$\bar f_{\rm pl}^*=-\sum_{1\leq i<j\leq N}\frac{\bar m_i\bar m_j}{4\p^2}\int_{\torus^2}\frac{d\bar\l_id\bar\l_j}{|\hat x_i^*-\hat x_j|}$$ coincides with the function $\bar f_{\rm pl}(\L,\eta_*,\xi_*)$, where $\bar f_{\rm pl}$ is the secular perturbation  of the plane problem. So,
recalling that the expansion of $\bar f_{\rm pl}$ is $$\bar f_{\rm pl}=-\sum_{1\leq i<j\leq N}\frac{\bar m_i\bar m_j}{4\p^2}\int_{\torus^2}\frac{d\l_id\l_j}{|\hat x_i-\hat x_j|}=C_0+\frac{1}{2}{\cal Q}_h\cdot(\eta^2+\xi^2)+{\eufm F}_h(\eta,\xi)+{\rm O}(6)$$
where  
$$
C_0:=\bar f_{0}:=\bar f_{\rm pl}|_{(\eta,\xi)=0}=-\sum_{1\leq i<j\leq N}\frac{\bar m_i\bar m_j}{a_j}b_{1/2,0}(a_i/a_j)\ ,$$ ${\cal Q}_h$ is the quadratic form associated to  the matrix ${\cal F}_h$ defined in Lemma \ref{exp to 2} and
\footnote{The entries of the matrix ${\cal F}_h$ defined in Lemma \ref{exp to 2} can be written in terms of the only Laplace coefficients $b_{3/2,0}(\a)$, $b_{3/2,1}(\a)$ as in (\ref{QhQv}).},  
${\eufm F}_h$ is the quartic form
$${\eufm F}_h=q\cdot(\eta^4+\xi^4)+r\cdot\eta^2\xi^2$$
where $q$, $r$ are the $4$--tensors of (\ref{r and q})$\div$(\ref{a4000 and others}), then,
\beqano\bar f_{\rm pl}^*=C_0+\frac{1}{2}{\cal Q}_h\cdot(\eta_*^2+\xi_*^2)+{\eufm F}_h(\eta_*,\xi_*)+{\rm O}(6)\ .\eeqano

\item[(ii)] {\bf Expansion of $\bar f_{\rm two}$} The function \beq{fmix}\bar f_{\rm two}=-\sum_{1\leq i<j\leq N}\frac{\bar m_i\bar m_j}{4\p^2}\int_{\torus^2}\frac{\hat x_i^*\cdot\hat{\eufm S}_{ij}^2\hat x_j}{|\hat x_i^*-\hat x_j|^3}d\bar\l_id\bar\l_j\eeq is of order $2$ in $\bar z$, so, it contributes to the $4$--expansion of $\bar f_{\rm  plt,pr}$ with terms of order $2$ and $4$, which we denote  $\bar f_{\rm two}|_2$, $\bar f_{\rm two}|_4$: 
\beqano\bar f_{\rm two}=\bar f_{\rm two}|_2+\bar f_{\rm two}|_4+{\rm O}(6)\ .\eeqano

\vskip.1in
\noi
Let us represent $\hat{\eufm S}_{ij}^2$ through its entries 
\beqa{repr Sij2 Sij4}\hat{\eufm S}_{ij}^2=\left(
\begin{array}{lrr}
{\eufm q}^{ij}_{11}&{\eufm q}^{ij}_{12}\\
{\eufm q}^{ij}_{21}&{\eufm q}^{ij}_{22}
\end{array}
\right)\eeqa
so as to  write the integrand function of (\ref{fmix}) as
\beq{exp fmix}\frac{\hat x_i^*\cdot\hat{\eufm S}_{ij}^2\hat x_j}{|\hat x_i^*-\hat x_j|^3}=\frac{\hat x_i^{*1}\hat x_j^1{\eufm q}^{ij}_{11}+\hat x_i^{*2}\hat x_j^2{\eufm q}^{ij}_{22}+\hat x_i^{*1}\hat x_j^2{\eufm q}^{ij}_{12}+\hat x_i^{*2}\hat x_j^1{\eufm q}^{ij}_{21}}{|\hat x_i^*-\hat x_j|^3}\ ,\eeq
 where $(\hat x_i^1,\hat x_i^2)$ are the components of $\hat x_i$.

\vskip.1in
\noi
Then, the term $f_{\rm two}|_2$ is found  replacing  $\hat x_i$ with its  $0$--approximation 
\beq{Order 0 Poinc}\hat x_i^{0}=a_i(\cos{\bar\l}_i,\sin{\bar\l_i})\eeq
into (\ref{fmix}) and next into (\ref{exp fmix}); this gives
\beq{fmix2}f_{\rm two}|_2=f_{\rm two}|_{\hat x=(\ref{Order 0 Poinc})}=-\sum_{1\leq i<j\leq N}{\bar m_i\bar m_j}({\eufm q}^{ij}_{11}+{\eufm q}^{ij}_{22})\frac{a_i}{2a_j^2}b_{3/2,1}(a_i/a_j)\eeq

\vskip.1in
\noi
For the computation of $f_{\rm two}|_4$,  
we denote as ${\cal Q}_{ij}^{hk}$ the quadratic forms acting on $(\eta_i,\eta_j,\xi_i,\xi_j)$ which realize the $\bar z$--expansion up to order $2$, 
 of the four integrals 
\beq{integrals}{\cal I}_{ij}^{hk}:=\frac{1}{4\p^2}\int_{\torus^2}\frac{\hat x_i^{*h}\hat x_j^k}{|\hat x_i^*-\hat x_j|^3}d\bar\l_id\bar\l_j\qquad h,\ k=1,\ 2\ ,\eeq \ie, we let
$${\cal I}_{ij}^{hk}=\arr{\frac{a_i}{2a_j^2}b_{3/2,1}(a_i/a_j)+{\cal Q}_{ij}^{hh}+\cdots\quad h=k=1,\ 2\\
{\cal Q}_{ij}^{hk}+\cdots\quad h\neq k=1,\ 2}\ .$$ 
Then, in view of (\ref{fmix})$\div$(\ref{exp fmix}), we have
\beq{fmix4}f_{\rm two}|_4=-\sum_{1\leq i<j\leq N}\bar m_i\bar m_j\sum_{1\leq h,k\leq 2}\,{\cal Q}_{ij}^{hk}\,{\eufm q}^{ij}_{hk}\eeq

\vskip.1in
\noi
The computation of the polynomials  ${\cal Q}_{ij}^{hk}$ is quite lenghty.   It is performed using, into (\ref{integrals}), the approximation of $\hat x_i$ up to order $2$
\footnote{which is
\beqa{2 Poinc}
\hat x_i^1&=&a_i\,\left(\cos\bar\l_i-\frac{3-\cos{2\bar\l_i}}{2}\hat\eta_i-\frac{\sin{2\bar\l_i}}{2}\hat\xi_i\right.\nonumber\\
&+&\left.\frac{-3\cos\bar\l_i+3\cos{3\bar\l_i}}{8}\hat\eta_i^2-\frac{\sin\bar\l_i+3\sin3\bar\l_i}{4}\hat\eta_i\hat\xi_i-\frac{5\cos\bar\l_i+3\cos3\bar\l_i}{8}\hat\xi_i^2\right)\nonumber\\
&+&{\rm O}(3)\nonumber\\
\nonumber\\
\hat x_i^2&=&a_i\,\left(\sin\bar\l_i+\frac{\sin{2\bar\l_i}}{2}\hat\eta_i+\frac{3+\cos{2\bar\l_i}}{2}\hat\xi_i+\right.\nonumber\\
&+&\frac{-5\sin\bar\l_i+3\sin3\bar\l_i}{8}\hat\eta_i^2\left.+\frac{-\cos\bar\l_i+3\cos3\bar\l_i}{4}\hat\eta_i\hat\xi_i-\frac{3\sin\bar\l_i+3\sin{3\bar\l_i}}{8}\hat\xi_i^2\right)\nonumber\\
&+&{\rm O}(3)\qquad \textrm{with}\quad (\hat\eta_i,\hat\xi_i):=\left(\frac{\bar\eta_i}{\sqrt{\L_i}},\frac{\bar\xi_i}{\sqrt{\L_i}}\right)\nonumber\\
\eeqa 
}
and next isolating the quadratic terms in $\bar z$, whose coefficients, as in the expansion of the secular perturbation of the plane problem, have the form of the mean over $(\bar\l_i$, $\bar\l_j)\in \torus^2$ of ratios of trigonometric polynomials in $\bar\l_i$, $\bar\l_j$, with ``Laplace'' denominators  $d_{ij}^s:=|\hat x_i^0-\hat x_j^0|^s$. The result is 
\beqa{Iij}
{\cal Q}_{ij}^{11}(\eta^*_i,\eta^*_j,\xi^*_i,\xi^*_j)&=&C_3(a_i,a_j)\frac{{\eta^*_i}^2}{\L_i}+C_4(a_i,a_j)\frac{\eta^*_i\eta^*_j}{\sqrt{\L_i\L_j}}+C_5(a_i,a_j)\frac{{\eta^*_j}^2}{\L_j}\nonumber\\
&+&C_6(a_i,a_j)\frac{{\xi^*_i}^2}{\L_i}+C_7(a_i,a_j)\frac{\xi^*_i\xi^*_j}{\sqrt{\L_i\L_j}}+C_8(a_i,a_j)\frac{{\xi^*_j}^2}{\L_j}\nonumber\\
{\cal Q}_{ij}^{22}(\eta^*_i,\eta^*_j,\xi^*_i,\xi^*_j)&=&C_6(a_i,a_j)\frac{{\eta^*_i}^2}{\L_i}+C_7(a_i,a_j)\frac{\eta^*_i\eta^*_j}{\sqrt{\L_i\L_j}}+C_8(a_i,a_j)\frac{{\eta^*}_j^2}{\L_j}\nonumber\\
&+&C_3(a_i,a_j)\frac{{\xi^*_i}^2}{\L_i}+C_4(a_i,a_j)\frac{\xi^*_i\xi^*_j}{\sqrt{\L_i\L_j}}+C_5(a_i,a_j)\frac{{\xi^*_j}^2}{\L_j}\nonumber\\
{\cal Q}_{ij}^{12}(\eta^*_i,\eta^*_j,\xi^*_i,\xi^*_j)&=&C_9(a_i,a_j)\frac{\eta^*_i\xi^*_i}{\L_i}+C_{10}(a_i,a_j)\frac{\eta^*_i\xi^*_j}{\sqrt{\L_i\L_j}}+C_{11}(a_i,a_j)\frac{\eta^*_j\xi^*_i}{\sqrt{\L_i\L_j}}\nonumber\\
&+&C_{12}(a_i,a_j)\frac{\eta^*_j\xi^*_j}{\L_j}\nonumber\\
{\cal Q}_{ij}^{21}(\eta^*_i,\eta^*_j,\xi^*_i,\xi^*_j)&=&C_9(a_i,a_j)\frac{\eta^*_i\xi^*_i}{\L_i}+C_{11}(a_i,a_j)\frac{\eta^*_i\xi^*_j}{\sqrt{\L_i\L_j}}+C_{10}(a_i,a_j)\frac{\eta^*_j\xi^*_i}{\sqrt{\L_i\L_j}}\nonumber\\
&+&C_{12}(a_i,a_j)\frac{\eta^*_j\xi^*_j}{\L_j}\nonumber\\
\eeqa
where $C_3(a_i,a_j)\div C_{12}(a_i,a_j)$ are defined, in terms of the Laplace coefficient, as
\beqa{Cij}
C_3(a_i,a_j)&:=& \frac{(a_i/a_j)^2}{32a_j}(57 (a_i/a_j)^2 +
    117)b_{7/2,0}(a_i/a_j)\nonumber\\
    &+&\frac{(a_i/a_j)}{64 a_j}(-12(a_i/a_j)^4-291(a_i/a_j)^2-12) b_{7/2,1}(a_i/a_j))\nonumber\\
    &+&\frac{(a_i/a_j)^2}{32 a_j}(15(a_i/a_j)^2-45)b_{7/2,2}(a_i/a_j))+27\frac{(a_i/a_j)^3}{64a_j}b_{7/2,3}(a_i/a_j)) \nonumber\\
C_4(a_i,a_j)&:=&277\frac{(a_i/a_j)^3}{32a_j}b_{7/2,0}(a_i/a_j))\nonumber\\
&+&\frac{(a_i/a_j)^2}{64 a_j}(-376-376(a_i/a_j)^2)b_{7/2,1}(a_i/a_j))\nonumber\\
    &+&\frac{(a_i/a_j)}{32 a_j}(16(a_i/a_j)^4+10(a_i/a_j)^2+16)b_{7/2,2}(a_i/a_j)\nonumber\\
    &+&56\frac{(a_i/a_j)^2}{64a_j}(a_i/a_j)^2+1)b_{7/2,3}(a_i/a_j)+\frac{(a_i/a_j)^3}{32a_j}b_{7/2,4}(a_i/a_j)\nonumber\\
C_5(a_i,a_j)&:=&\frac{(a_i/a_j)^2}{32a_j}(117 (a_i/a_j)^2 +57)b_{7/2,0}(a_i/a_j)\nonumber\\
&+&\frac{(a_i/a_j)}{64 a_j}(-12(a_i/a_j)^4-291(a_i/a_j)^2-12)b_{7/2,1}(a_i/a_j)\nonumber\\
     &+&\frac{(a_i/a_j)^2}{32 a_j}(-45(a_i/a_j)^2+15)b_{7/2,2}(a_i/a_j)+27\frac{(a_i/a_j)^3}{64a_j}b_{7/2,3}(a_i/a_j)\nonumber\\  
C_6(a_i,a_j)&:=&\frac{(a_i/a_j)^2}{32a_j}(71 (a_i/a_j)^2 + 11)b_{7/2,0}(a_i/a_j)\nonumber\\
&+&\frac{(a_i/a_j)}{64 a_j}(-20(a_i/a_j)^4+119(a_i/a_j)^2-20)b_{7/2,1}(a_i/a_j)\nonumber\\
&+&\frac{(a_i/a_j)^2}{32 a_j}(-79(a_i/a_j)^2-19)b_{7/2,2}(a_i/a_j)-47\frac{(a_i/a_j)^3}{64a_j}b_{7/2,3}(a_i/a_j)\nonumber\\
C_7(a_i,a_j)&:=&-215 \frac{(a_i/a_j)^3}{32a_j}b_{7/2,0}(a_i/a_j)\nonumber\\
&+&\frac{(a_i/a_j)^2}{64 a_j}(8+8(a_i/a_j)^2)b_{7/2,1}(a_i/a_j)\nonumber\\
    &+&\frac{(a_i/a_j)}{32 a_j}(16(a_i/a_j)^4+118(a_i/a_j)^2+16)b_{7/2,2}(a_i/a_j)\nonumber\\
    &+&56\frac{(a_i/a_j)^2}{64a_j}(a_i/a_j)^2+1)b_{7/2,3}(a_i/a_j)+\frac{(a_i/a_j)^3}{32a_j}b_{7/2,4}(a_i/a_j)\nonumber\\
C_8(a_i,a_j)&:=&\frac{(a_i/a_j)^2}{32a_j} (11 (a_i/a_j)^2+ 71)b_{7/2,0}(a_i/a_j)\nonumber\\
&+&\frac{(a_i/a_j)}{64 a_j}(-20(a_i/a_j)^4+119(a_i/a_j)^2-20)b_{7/2,1}(a_i/a_j)\nonumber\\
    &+&\frac{(a_i/a_j)^2}{32 a_j}(-19(a_i/a_j)^2+79)b_{7/2,2}(a_i/a_j)-47\frac{(a_i/a_j)^3}{64a_j}b_{7/2,3}(a_i/a_j)\nonumber\\  
C_9(a_i,a_j)&:=&\frac{(a_i/a_j)^2}{32a_j}(14(a_i/a_j)^2-106)b_{7/2,0}(a_i/a_j)\nonumber\\
&+&\frac{(a_i/a_j)}{32a_j}(-4(a_i/a_j)^4+205(a_i/a_j)^2-4)b_{7/2,1}(a_i/a_j)\nonumber\\
&+&\frac{(a_i/a_j)^2}{16a_j}(-47 (a_i/a_j)^2+13)b_{7/2,2}(a_i/a_j)\nonumber\\
&-&\frac{37}{32a_j}(a_i/a_j)^3b_{7/2,3}(a_i/a_j)\nonumber\\
C_{10}(a_i,a_j)&:=&-\frac{35}{32a_j}(a_i/a_j)^3b_{7/2,0}(a_i/a_j)\nonumber\\
&+&\frac{1}{32a_j}(a_i/a_j)^2(4(a_i/a_j)^2+4)b_{7/2,1}(a_i/a_j)\nonumber\\   
&+&\frac{(a_i/a_j)}{16a_j}(8(a_i/a_j)^4-31(a_i/a_j)^2+8)b_{7/2,2}(a_i/a_j)\nonumber\\
&+&\frac{(a_i/a_j)^2}{32a_j}(28(a_i/a_j)^2+28)b_{7/2,3}(a_i/a_j)\nonumber\\
&+&\frac{(a_i/a_j)^3}{32a_j}b_{7/2,4}(a_i/a_j)\nonumber\\
C_{11}(a_i,a_j)&:=&-\frac{457}{32a_j}(a_i/a_j)^3b_{7/2,0}(a_i/a_j)\nonumber\\
&+&\frac{1}{32a_j}(a_i/a_j)^2(188(a_i/a_j)^2+188)b_{7/2,1}(a_i/a_j)\nonumber\\   
&+&\frac{(a_i/a_j)}{16a_j}(-8(a_i/a_j)^4+85(a_i/a_j)^2-8)b_{7/2,2}(a_i/a_j)\nonumber\\
&+&\frac{(a_i/a_j)^2}{32a_j}(-28(a_i/a_j)^2-28)b_{7/2,3}(a_i/a_j)\nonumber\\
&-&\frac{(a_i/a_j)^3}{32a_j}b_{7/2,4}(a_i/a_j)\nonumber\\   
C_{12}(a_i,a_j)&:=&\frac{(a_i/a_j)^2}{32a_j}(14-106(a_i/a_j)^2)b_{7/2,0}(a_i/a_j)\nonumber\\
&+&\frac{(a_i/a_j)}{32a_j}(-4(a_i/a_j)^4+205(a_i/a_j)^2-4)b_{7/2,1}(a_i/a_j)\nonumber\\
&+&\frac{(a_i/a_j)^2}{16a_j}(-47+13(a_i/a_j)^2)b_{7/2,2}(a_i/a_j)\nonumber\\
&-&\frac{37}{32a_j}(a_i/a_j)^3b_{7/2,3}(a_i/a_j)\nonumber\\   
\eeqa
\item[(iii)]  {\bf Expansion of $\bar f_{\rm four}$.} The function \beq{ffour1}\bar f_{\rm four}=-\sum_{1\leq i<j\leq N}\frac{\bar m_i\bar m_j}{4\p^2}\int_{\torus^2}\left(\frac{\hat x_i^*\cdot\hat{\eufm S}_{ij}^4\hat x_j}{|\hat x_i^*-\hat x_j|^3}+\frac{3}{2}\frac{(\hat x_i^*\cdot\hat{\eufm S}_{ij}^2\hat x_j)^2}{|\hat x_i^*-\hat x_j|^5}\right)d\bar\l_id\bar\l_j\eeq  is of order $4$ in $\bar z=(\bar\eta,\bar\xi,\bar p,\bar q)$, so, it sufficies replace $\hat x_i$ with its $0$--order approximation (\ref{Order 0 Poinc}) to find 
\beqano\bar f_{\rm four}=\bar f_{\rm four}|_4+{\rm O}(6)\eeqano
As in the previous step, we represent $\hat{\eufm S}_{ij}^2$, $\hat{\eufm S}_{ij}^4$ through their entries:  $\hat{\eufm S}_{ij}^2$ as in eq. (\ref{repr Sij2 Sij4}) and $\hat{\eufm S}_{ij}^4$ as
$$\hat{\eufm S}_{ij}^4=\left(
\begin{array}{lrr}
{\eufm Q}^{ij}_{11}&{\eufm Q}^{ij}_{12}\\
{\eufm Q}^{ij}_{21}&{\eufm Q}^{ij}_{22}
\end{array}
\right)$$
Using (\ref{Order 0 Poinc}) , into (\ref{ffour1}), we find
\beqa{ffour}
\bar f_{\rm four}|_4&=&\bar f_{\rm four}|_{\hat x=(\ref{Order 0 Poinc})}\nonumber\\
&=&-\sum_{1\leq i<j\leq N}{\bar m_i\bar m_j}({\eufm Q}^{ij}_{11}+{\eufm Q}^{ij}_{22})\frac{a_i}{2a_j^2}b_{3/2,1}(a_i/a_j)\nonumber\\
&-&\frac{3}{8}\sum_{1\leq i<j\leq N}{\bar m_i\bar m_j}(({\eufm q}^{ij}_{11})^2+({\eufm q}^{ij}_{12})^2+({\eufm q}^{ij}_{21})^2+({\eufm q}^{ij}_{22})^2)\frac{a_i^2}{a_j^3}b_{5/2,0}(a_i/a_j)\nonumber\\
&-&\frac{3}{16}\sum_{1\leq i<j\leq N}{\bar m_i\bar m_j}(({\eufm q}^{ij}_{11}+{\eufm q}^{ij}_{22})^2+({\eufm q}^{ij}_{12}-{\eufm q}^{ij}_{21})^2)\frac{a_i^2}{a_j^3}b_{5/2,2}(a_i/a_j)\nonumber\\
\eeqa
\item[(iv)] {\bf Computation of $\hat{\eufm S}_{ij}^2$.}
The matrices ${\eufm S}_i$ are defined through equation (\ref{eufmSi}) as suitable products of  the matrices ${\cal S}_i^{\rm pr}$'s, $\tilde{\cal S}_j^{\rm pr}$'s, which (recall Proposition \ref{partial reduction prop}) are easily expanded up to order $2$ as
\beqa{exp cal S}
\arr{
{{\cal S}_1^{\rm pr}}^{\rm T}=\Sigma_{c_1}(\bar p_1,\bar q_1)\\
{{\cal S}_i^{\rm pr}}=\Sigma_{c_i}(\bar p_{i-1},\bar q_{i-1})\quad (2\leq i\leq N)\\
\tilde{{\cal S}_i^{\rm pr}}^{\rm T}=\Sigma_{\tilde c_i}(\bar p_i,\bar q_i)\qquad (2\leq i\leq N-1)\\
}\quad +{\rm O}(3)
\eeqa
where $c_i$'s, $\tilde c_j$'s are the constants (\ref{ccc}) which define the entries the matrix associated to the operator ${\eufm L}$ and $\Sigma_{\hat c}(\hat p,\hat q)$ denotes
\beqa{Sigma}\Sigma_{\hat c}(\hat p,\hat q)=\left(
\begin{array}{lcrr}
1-\frac{1}{2}\hat c^2\hat q^2&-\frac{1}{2}\hat c^2\hat p\hat q&-\hat c\hat q\\
-\frac{1}{2}\hat c^2\hat p\hat q&1-\frac{1}{2}\hat c^2\hat p^2&-\hat c\hat p\\
\hat c\hat q&\hat c\hat p&1-\frac{1}{2}\hat c^2(\hat p^2+\hat q^2)
\end{array}
\right)
\eeqa
Taking then  the products of ${\cal S}_i^{\rm pr}$, $\cdots$ as prescribed in (\ref{eufmSi}), we have
\beqano
{\eufm S}_i^{\rm pr}
&=&\Sigma_{\ell_{i,N-1}}(\bar p_{N-1},\bar q_{N-1})\cdots \Sigma_{\ell_{i,1}}(\bar p_{1},\bar q_{1})+{\rm O}(3)\nonumber\\
\eeqano
and hence,
\beqa{first write}
{\eufm S}_{ij}&=&{\eufm S}_i^T{\eufm S}_j=\Sigma_{-\ell_{i1}}(\bar p_{1},\bar q_{1})\cdots\Sigma_{-\ell_{i,N-1}}(\bar p_{N-1},\bar q_{N-1})\nonumber\\
&\times&\Sigma_{\ell_{j,N-1}}(\bar p_{N-1},\bar q_{N-1})\cdots\Sigma_{\ell_{j1}}(\bar p_{1},\bar q_{1})\nonumber\\
\eeqa where
$$\ell_i=(\ell_{i1}, \cdots,\ell_{i,{N-1}})$$
is the $i^{\rm th}$ row of the matrix $\ell$ associated to the operator ${\eufm L}$ (eq. \ref{Lop})$\div$(\ref{ccc}))
Using (\ref{eufmSi}), we can write
\beqa{eufmSij}
{\eufm S}_{ij}={\eufm S}_i^T{\eufm S}_j=\arr{
\bar{\cal S}^{\rm pr}_1\quad i=1,\ j=2\\
{\cal S}^{\rm pr}_1\bar{\cal S}^{\rm pr}_2\quad i=1,\ j=3\\
{\cal S}^{\rm pr}_1\tilde{\cal S}^{\rm pr}_2\cdots\tilde{\cal S}^{\rm pr}_{j-2}\bar{\cal S}^{\rm pr}_{j-1}\quad i=1,\ j\geq 4\\
{{\cal S}^{\rm pr}_i}^T\bar{\cal S}^{\rm pr}_{i}\quad 2\leq i=j-1\\
{{\cal S}_i^{\rm pr}}^T\tilde{\cal S}^{\rm pr}_i\cdots\tilde{\cal S}^{\rm pr}_{j-2}\bar{\cal S}_{j-1}^{\rm pr}\quad 2\leq i\leq j-2
}
\eeqa
having let
\beq{barcalSi}\bar{\cal S}^{\rm pr}_i:=\tilde{\cal S}^{\rm pr}_i{\cal S}^{\rm pr}_{i+1}=\Sigma_{\bar c_i}(p_i,q_i)+{\rm O}(3)\eeq
with
$$\bar c_i=c_{i+1}-\tilde c_i\qquad 1\leq i\leq N-1\quad (\tilde{\cal S}_1:={\cal S}_1\ ,\quad \tilde{c}_1:={c}_1)$$
which, together with (\ref{exp cal S}), implies (\ref{second write})$\div$(\ref{second write1}).
Hence, we may also  write the matrices
${\eufm S}_{ij}={\eufm S}_i^T{\eufm S}_j$ as
\beqa{second write}
{\eufm S}_{ij}&=&\prod_{1\leq k\leq j-1}\Sigma_{\ell^{ij}_k}(\bar p_{k},\bar q_k)+{\rm O}(3)
\eeqa
%
%
%
%
%
%
where
\beqa{second write1}\ell^{ij}_k=\ell_{jk}-\ell_{ik}
\eeqa
and then we use the following elementary Lemma.
\begin{lemma}
Let 
$$c=(c_1,\cdots,c_m)\ ,\quad p=(p_1,\cdots,p_m)\ ,\quad q=(q_1,\cdots,q_m)\in \real^m\ ,$$
let $\Sigma_{\hat c}(\hat p,\hat q)$ as in (\ref{Sigma}) and put
$$\P_c(p,q):=\Sigma_{c_1}(p_1,q_1)\cdots \Sigma_{c_m}(p_m,q_m)\ .$$
Then, the submatrix $\hat\P_c(p,q)$ of order $2$  of $\P_c(p,q)$
is 
$$\hat\P_c(p,q)=\id-\frac{1}{2}\left(
\begin{array}{lrrr}
(c\cdot q)^2&(c\cdot p)(c\cdot q)-\D_c(p,q)&\\
(c\cdot p)(c\cdot q)+\D_c(p,q)&(c\cdot p)^2&\\
\end{array}
\right)+{\rm O}(4)$$
where $\D_c(p,q)$ denotes
$$\D_c(p,q):=\sum_{1\leq i<j\leq m}c_ic_j(p_iq_j-p_jq_i)\quad (\textrm{when}\quad m\geq 2)\ .$$
\end{lemma}
In view of the previous Lemma and making use of (\ref{first write})$\div$(\ref{second write}), we find the following expressions for the entries $({\eufm q}^{ij}_{hk})_{hk}$ of ${\eufm S}_{ij}^2$
\beqa{order two vertical}
\arr{
{\eufm q}^{ij}_{11}=-\frac{1}{2}\Big({\eufm L}\bar q_j-{\eufm L}\bar q_i\Big)^2\\
{\eufm q}^{ij}_{12}=-\frac{1}{2}\Big({\eufm L}\bar q_j-{\eufm L}\bar q_i\Big)\Big({\eufm L}\bar p_j-{\eufm L}\bar p_i\Big)+\frac{1}{2}\sum_{1\leq h<k\leq j-1}\ell^{ij}_h\ell^{ij}_k(\bar p_h\bar q_k-\bar p_k\bar q_h)\\
{\eufm q}^{ij}_{21}=-\frac{1}{2}\Big({\eufm L}\bar q_j-{\eufm L}\bar q_i\Big)\Big({\eufm L}\bar p_j-{\eufm L}\bar p_i\Big)-\frac{1}{2}\sum_{1\leq h<k\leq j-1}\ell^{ij}_h\ell^{ij}_k(\bar p_h\bar q_k-\bar p_k\bar q_h)\\
{\eufm q}^{ij}_{22}=-\frac{1}{2}\Big({\eufm L}\bar p_j-{\eufm L}\bar p_i\Big)^2\ .\quad \\
}
\eeqa
\item[(v)] {\bf Computation of $\hat{\eufm S}_{ij}^4$.} 
Equation (\ref{eufmSij})  gives the expression of ${\eufm S}_{ij}$ in terms of the matrices ${\cal S}_i^{\rm pr}$, $\tilde{\cal S}_j^{\rm pr}$ (eq. (\ref{eufm ci si part})) and $\bar{\cal S}_k^{\rm pr}=\tilde{\cal S}_k^{\rm pr}{\cal S}_{k+1}^{\rm pr}$. Hence, the expression of ${\eufm S}_{ij}^4$ is uniquely determined by the expansion of these matrices  up to order $4$.

\vskip.1in
\noi
Let us first notice that we can write the matrices $\bar{\cal S}_k^{\rm pr}$ in the same form as ${\cal S}_i^{\rm pr}$, $\tilde{\cal S}_j^{\rm pr}$. In fact, recalling the definitions (eq. (\ref{matrices})) of $\tilde{\cal S}_k^{\rm pr}$, ${\cal S}_{k+1}^{\rm pr}$, we find, for $1\leq k\leq N-1$,
\beq{bar cal Sij as prod}\bar{\cal S}_k^{\rm pr}=\tilde{\cal S}_k^{\rm pr}{\cal S}_{k+1}^{\rm pr}={\cal R}_{\rm z}(\k_{k}){\cal R}_{\rm x}(\bar i_k){\cal R}_{\rm z}(-\k_{k})\qquad \k_{k}:=\arg{(p_k,-q_k)}\ ,\eeq
where $\bar i_k:=\tilde i_k+i_{k+1}$ has the meaning of the outern angle
\footnote{In fact, by the definition \ref{ijtildeij} of $\tilde i_k$, $i_{k+1}$  given in Proposition \ref{inversion} it is clear that $\tilde i_k$, $i_{k+1}$  have the meaning, respectively, of the inner angles corresponding to the couples of sides with lenghth ($\Psi_{k-1}$, $\Psi_{k}$) and ($\Psi_{k-1}$, $\G_{k+1}$).} 
of $\Psi_{k-1}$, $\G_{k+1}$ in the triangle with sides with lenghth $\Psi_{k-1}$, $\G_{k+1}$ $\Psi_k$, hence,
\beq{cosbari}
\cos{\bar i_k}=\frac{\Psi_k^2-\Psi_{k-1}^2-\G_{k+1}^2}{2\Psi_{k-1}\G_{k+1}}\qquad \bar i_k\in (0,\p)
\eeq
Equations (\ref{bar cal Sij as prod})$\div$(\ref{cosbari}) easily imply
\beqa{barcalS}
\bar{\cal S}_k^{\rm pr}&=&\left(
\begin{array}{ccc}
1-\bar q_{k}^2\bar{\eufm c}_k&-\bar p_{k}\bar q_{k}\bar{\eufm c}_k&-\bar q_{k}\bar{\eufm s}_k\\
-\bar p_{k}\bar q_{k}\bar{\eufm c}_k&1-\bar p_{k}^2\bar{\eufm c}_k&-\bar p_{k}\bar{\eufm s}_k\\
\bar q_{k}\bar{\eufm s}_k&\bar p_{k}\bar{\eufm s}_k&1-(\bar p_{k}^2+\bar q_{k}^2)\bar{\eufm c}_k
\end{array}
\right)\quad 1\leq k\leq N-2\nonumber\\
\eeqa
where
\beq{eufmbarc}\bar{\eufm c}_k=\frac{1-\cos{\bar i}_k}{\bar p_k^2+\bar q_k^2}=\frac{\Psi_{k-1}+\G_{k+1}+\Psi_k}{4\Psi_{k-1}\G_{k+1}}\ ,\quad \bar{\eufm s}_k=\frac{\sin{\bar i}_k}{\sqrt{\bar p_k^2+\bar q_k^2}}=\sqrt{\bar{\eufm c}_k\,\Big(2-(\bar p_k^2+\bar q_k^2)\bar{\eufm c}_k\Big)}\ .\eeq

\vskip.1in
\noi
We are now ready for the expansions of ${\cal S}_i^{\rm pr}$, $\tilde{\cal S}_j^{\rm pr}$, $\bar{\cal S}_k^{\rm pr}$.

\vskip.1in
\noi
Let us  observe that the functions ${\eufm c}_i$, ${\eufm s}_j$, $\bar{\eufm s}_k$ $\tilde{\eufm c}_j$, $\tilde{\eufm s}_j$, $\bar{\eufm c}_k$, $\bar{\eufm s}_k$ (eq. (\ref{eufm ci si})$\div$(\ref{eufmbarc})) are even in $\bar z$, so, the entries of ${\cal S}_i$, $\cdots$, with places $(1,1)$, $(1,2)$, $(2,1)$, $(2,2)$, $(3,3)$ are even (as functions of $\bar z$) and the ones with places $(1,3)$, $(2,3)$, $(3,1)$, $(3,2)$ are odd. Then, in order to obtain expansions of ${\cal S}_i^{\rm pr}$, $\tilde{\cal S}^{\rm pr}_j$, $\bar{\cal S}^{\rm pr}_k$ up to order $4$, it sufficies to expand the functions ${\eufm c}_i$, ${\eufm s}_i$, $\cdots$ up to order $2$. Making this operation leads to the fllowing expansions (which generalizes (\ref{exp cal S})$\div$(\ref{Sigma})). If 
$\dst \Sigma_{\hat c,\hat{\cal C},\hat S}(\hat p,\hat q)$ denotes the matrix
\beqa{Sigma order 4}
\Sigma_{\hat c,\hat{\cal C},\hat S}(\hat p,\hat q)=\left(
\begin{array}{lcrr}
1-\hat q^2\left(\frac{\hat c^2}{2}+\hat{\cal C}\right)&-\hat p\hat q\left(\frac{\hat c^2}{2}+\hat{\cal C}\right)&-\hat q\left(\hat c+\hat S\right)\\
-\hat p\hat q\left(\frac{\hat c^2}{2}+\hat{\cal C}\right)&1-\hat p^2\left(\frac{\hat c^2}{2}+\hat{\cal C}\right)&-\hat p\left(\hat c+\hat S\right)\\
\hat q\left(\hat c+\hat S\right)&\hat p\left(\hat c+\hat S\right)&1-(\hat p^2+\hat q^2)\left(\frac{\hat c^2}{2}+\hat{\cal C}\right)
\end{array}
\right)
\eeqa
then, \beqa{exp S}
\arr{
{{\cal S}^{\rm pr}_1}^{\rm T}=\Sigma_{c_1,{\cal C}_1,S_1}(p_1,q_1)\\
{\cal S}_i^{\rm pr}=\Sigma_{c_i,{\cal C}_i,S_i}(p_{i-1},q_{i-1})\quad (2\leq i\leq N)\\
{{\tilde{\cal S}}_i^{\rm prT}}=\Sigma_{\tilde c_i,\tilde{\cal C}_i,\tilde S_i}(p_i,q_i)\qquad (2\leq i\leq N-1)\\
{\bar{\cal S}_i}^{\rm pr}=\Sigma_{\bar c_i,\bar{\cal C}_i,\bar S_i}(p_i,q_i)\qquad (1\leq i\leq N-1)
}\quad +{\rm O}(5)
\eeqa
where, for $2\leq i\leq N$, $2\leq j\leq N-1$, $1\leq k\leq N-1$,
\beqa{cal C}
&&\arr{\dst c_1:=-\sqrt{\frac{\L_2}{\L_1{\rm L}_2}}\\
\dst {\cal C}_1:=\frac{2\L_2(2\L_1+\L_2)\t_1-2\L_1^2\t_2+\L_1(\L_2-\L_1)\r_1}{4\L_1^2{\rm L}_2^2}\\ 
\dst S_1:=\frac{1}{c_1}\left({\cal C}_1-\frac{c_1^4}{4}\r_{1}\right)}\nonumber\\
&&\arr{\dst c_i:=\sqrt{\frac{{\rm L}_{i-1}}{{\rm L}_i\L_i}}\\
\dst {\cal C}_i:=\frac{2{\rm L}_{i-1}({\rm L}_{i-1}+2\L_i)\t_i+\L_i({\rm L}_{i-1}-\L_i)\r_{i-1}-2\L_i^2({\rm T}_{i-1}+{\rm R}_{i-2})}{4{\rm L}_i^2\L_i^2}\\ \dst S_i:=\frac{1}{c_i}\left({\cal C}_i-\frac{c_i^4}{4}\r_{i-1}\right)} \nonumber\\
&&\arr{\dst \tilde c_j:=-\sqrt{\frac{\L_{j+1}}{{\rm L}_{j+1}{\rm L}_j}}\\
\dst \tilde{\cal C}_j:=\frac{-2{\rm L}_j^2\t_{j+1}+{\rm L}_j(\L_{j+1}-{\rm L}_j)\r_j+2\L_{j+1}(2{\rm L}_j+\L_{j+1})({\rm T}_j+{\rm R}_{j-1})}{4{\rm L}_j^2{\rm L}_{j+1}^2}\\
\dst \tilde S_j:=\frac{1}{\tilde c_j}\left(\tilde{\cal C}_j-\frac{c_j^4}{4}\r_{j}\right)}\nonumber\\
&&\arr{\dst \bar c_k:=c_{k+1}-\tilde c_k=\sqrt{\frac{{\rm L}_{k+1}}{{\rm L}_k\L_{k+1}}}\\
\dst \bar{\cal C}_k:=\frac{2\L_{k+1}^2({\rm T}_k+{\rm R}_{k-1})+2{\rm L}_k^2\t_{k+1}-{\rm L}_k\L_{k+1}\r_k}{4{\rm L}_k^2\L_{k+1}^2}\\
\dst \bar S_k:=\frac{1}{\bar c_k}\left(\bar{\cal C}_k-\frac{\bar c_k^4}{4}\r_{k}\right)}\nonumber\\
\eeqa
having let, for shortness,
$$\t_i:=\frac{\bar\eta_i^2+\bar\xi_i^2}{2}\ ,\quad \r_i:=\frac{\bar p_i^2+\bar q_i^2}{2}\ ,\quad {\rm T}_i:=\sum_{1\leq j\leq i}\t_j\ ,\quad {\rm R}_i:=\sum_{1\leq j\leq i}\r_j\ ,\quad {\rm L}_i:=\sum_{1\leq j\leq i}\L_j$$
with \beq{R0}{\rm R}_0:=0\eeq
\item[(vi)]
{\bf Conclusion of the computation.} 
In view of the expressions (\ref{order two vertical}) for the diagonal terms ${\eufm q^{ij}_{11}}$, ${\eufm q^{ij}_{22}}$ of $\hat{\eufm S}_{ij}^2$, the function $\bar f_{\rm two}|_2$ (eq. (\ref{fmix2})) becomes
\beqano
f_{\rm two}|_2&=&\frac{1}{2}\sum_{1\leq i<j\leq N}{\bar m_i\bar m_j}\Big(({\eufm L}\bar q_j-{\eufm L}\bar q_i)^2+({\eufm L}\bar p_j-{\eufm L}\bar p_i)^2\Big)\frac{a_i}{2a_j^2}b_{3/2,1}(a_i/a_j)\nonumber\\
&=:&{\cal Q}^*_v\cdot\frac{\bar p^2+\bar q^2}{2}
\eeqano
Similarly,  replacing diagonal and off--diagonal entries  ${\eufm q^{ij}_{hk}}$ as in (\ref{order two vertical}), the function $f_{\rm two}|_4$(eq. (\ref{fmix4})) becomes
\beqano
f_{\rm two}|_4&=&\frac{1}{2}\sum_{1\leq i<j\leq N}\bar m_i\bar m_j({\cal Q}_{ij}^{11}\,({\eufm L}\bar q_j-{\eufm L}\bar q_i)^2+{\cal Q}_{ij}^{22}\,({\eufm L}\bar p_j-{\eufm L}\bar p_i)^2\nonumber\\
&+&({\cal Q}_{ij}^{12}+{\cal Q}_{ij}^{21})\Big({\eufm L}\bar p_j-{\eufm L}\bar p_i\Big)\Big({\eufm L}\bar q_j-{\eufm L}\bar q_i\Big)\nonumber\\
&+&({\cal Q}_{ij}^{21}-{\cal Q}_{ij}^{12})\,\sum_{1\leq h<k\leq j-1}\ell^{ij}_h\ell^{ij}_k(\bar p_h\bar q_k-\bar p_k\bar q_h)\nonumber\\
&=:&{\eufm F}_{hv}(\eta_*,\xi_*,\bar p,\bar q)
\eeqano
and the function $\bar f_{\rm four}|_4$ 
(eq. (\ref{ffour}) becomes
\beqano
\bar f_{\rm four}|_4
&=&\sum_{1\leq i<j\leq N}{\bar m_i\bar m_j}({\eufm Q}^{ij}_{11}+{\eufm Q}^{ij}_{22})C_1(a_i,a_j)\nonumber\\
&+&\sum_{1\leq i<j\leq N}{\bar m_i\bar m_j}\left(\Big({\eufm L}\bar p_i-{\eufm L}\bar p_j\Big)^2+\Big({\eufm L}\bar q_i-{\eufm L}\bar q_j\Big)^2\right)^2C_{13}(a_i,a_j)\nonumber\\
&+&\sum_{1\leq i<j\leq N}{\bar m_i\bar m_j}\left(\sum_{1\leq h<k\leq j-1}\ell^{ij}_h\ell^{ij}_k(\bar p_{h}\bar q_k-\bar p_k\bar q_h)\right)^2C_{14}(a_i,a_j)\nonumber\\
&=:&{\eufm F}_{v}(\eta_*,\xi_*,\bar p,\bar q)
\eeqano
where
\beqa{C13C14}
C_{13}(a_i,a_j)&:=&-\frac{3}{64}\frac{a_i^2}{a_j^3}(2b_{5/2,0}(a_i/a_j)+b_{5/2,2}(a_i/a_j))\nonumber\\
C_{14}(a_i,a_j)&:=&-\frac{3}{16}\frac{a_i^2}{a_j^3}(b_{5/2,0}(a_i/a_j)+b_{5/2,2}(a_i/a_j))\nonumber\\
\eeqa

\vskip.1in
\noi
Collecting then the expansions 
\beqano
\bar f_{\rm pl}^*&=&C_0+{\cal Q}_h\cdot\frac{\eta_*^2+\xi_*^2}{2}+{\eufm F}_h(\eta_*,\xi_*)+{\rm O}(6)\quad \eta_*=(-\bar\eta_1,\bar\eta_2,\cdots)\nonumber\\
\bar f_{\rm two}&=&{\cal Q}^*_v\cdot\frac{\bar p^2+\bar q^2}{2}+{\eufm F}_{hv}(\eta_*,\xi_*,\bar p,\bar q)+{\rm O}(6)\nonumber\\
\bar f_{\rm four}&=&{\eufm F}_{v}(\eta_*,\xi_*,\bar p,\bar q)+{\rm O}(6)\ .
\eeqano
we  finally find 
\beqano
\bar{f}_{\rm plt,pr}&=&\bar f_{\rm pl}^*+\bar f_{\rm two}+\bar f_{\rm four}\nonumber\\
&=&C_0+{\cal Q}_h\cdot\frac{\eta_*^2+\xi_*^2}{2}+{\cal Q}^*_v\cdot\frac{\bar p^2+\bar q^2}{2}\nonumber\\
&+&{\eufm F}_h(\eta_*,\xi_*)+{\eufm F}_{hv}(\eta_*,\xi_*,\bar p,\bar q)+{\eufm F}_{v}((\eta_*,\xi_*,\bar p,\bar q)+{\rm O}(6)
\eeqano
\end{itemize}
\subsection{$(3N-1)$--Dimensional KAM Tori and Measure of the Kolmogorov's Set}
Having checked, for $N\geq 3$, the assumptions of non resonance up to order $4$ for the first Birkhoff invariants and non degeneracy for the second Birkhoff invariants, invoking Theorem \ref{more general degenerate KAM}, we can state the following  result  concerning existence of KAM tori of dimension $3N-1$ for the planetary $(1+N)$ body problem and measure estimates of the invariant set.
\begin{theorem}\label{3N-1KAMtori and meas}
Consider, in $\real^3$, a star with mass $\bar m_0$ and $N\ge 3$ planets with masses $\m\bar m_1$, $\cdots$, $\m\bar m_N$, interacting only through gravity. Let $a_i$ denote the instantaneous major semiaxis of the $i^{th}$ planet. Then, there exists $\d_*>0$, $\e_*>0$, $\m_*>0$, $b>0$, $c>0$, $C>0$ such that, if $0<a_i/a_{i+1}<\d_*$,  $0<\e<\e_*$ and $0<\m<\m_*$ and
$$\m<c(\log{\e^{-1}})^{-2b}\ ,$$
 there exists a positive Lebesgue measure set $\cK$ 
 such that
 \begin{itemize}
 \item[(i)] $\cK$ (``Kolmogorov set'') is formed by the union of invariant tori of dimension $3N-1$ on which the  $\cH_{\rm plt}$--flow is linear in time, with Diophantine frequency;
 \item[(ii)] the measure of $\cK$ satisfies
 $$c\,\e^{2(2N-1)}>\meas\,\cK>c\,\Big(1-C(\sqrt{\e}+\sqrt{\m}(\log{\e^{-1}})^{b})\Big)\e^{2(2N-1)}\ .$$
 \end{itemize}
 Furthermore, the eccentricities and the mutual inclinations on the invariant tori are bounded by $c(\log\e^{-1})^{-1}$.
\end{theorem}
\newpage

\section{Kolmogorov's Set in the Space Planetary Problem II (Total Reduction)}
\setcounter{equation}{0}
The proofs of existence of quasi--periodic motions for the planetary problem presented in  \cite{Fej04} and \cite{ChPu08} are based on the application of a ($C^{\infty}$, analytic, respectively) KAM theory based on ``weak'' non--degeneracy conditions, for  a given properly degenerate system, nearly an elliptic equilibrium point.
\vskip.1in
\noi
For istance, the proof in \cite{ChPu08}, in the real--analytic framework, is obtained as an application of Theorem \ref{Russmann Thm} below, based, on turn, on 2001 R\"ussmann Theory \cite{Rus01} (compare \cite{ChPu08}, {\bf Theorem 4}), where the following weak non--degeneracy condition is required 
\begin{definition}[R\"ussmann nondegeneracy condition]\rm A real--analytic function
$$\o:\quad y\in {\cal B}\subset \real^n\to \o(y)=(\o_1(y),\cdots,\o_m(y))\in \real^m$$
is called {\bf R--non degenerate} if ${\cal B}$ is a non--empty open connected set in $\real^n$ and if for any $c\in\real^m\setminus\{0\}$, the map
\beq{Russmann}y\to c\cdot\o(y)\neq 0\ \eeq
\end{definition}

\begin{theorem}\label{Russmann Thm}
Consider a  Hamiltonian function 
$${\cal H}_\m=h(I)+\m\,f(I,\varphi,p,q)$$
which assume to be real-analytic for
$$(I,\varphi,p,q)\in {\cal I}\times\torus^{\bar n}\times B^{2\hat n}_r(0):={\cal M}$$
with the mean perturbation $\dst \bar f:=(2\p)^{-\bar n}\int_{\torus^{\bar n}}fd\varphi$ of the form
\beq{first order theory}\bar f=\bar f_0(I)+\O(I)\cdot J+O(|J|^{3/2})\quad J=\left(\cdots,\frac{p_i^2+q_i^2}{2},\cdots\right)\ .\eeq
Assume also that the ``frequency map''
$$I\in{\cal I}\to (\partial\,h(I),\O(I))\in \real^{\bar n}\times \real^{\hat n}$$
is $R$--non degenerate. Then, if $\m$ is sufficiently small, there exists a positive measure set of phase space points belonging to real--analytic ${\cal H}_\m$--invariant tori which are close to $\torus^{\bar n}\times\{I_0\}\times_j\{p_j^2+q_j^2=\r_j\}$, with $\r_j=O(\m)$. Furtherore, the ${\cal H}_\m$--flow
on such tori is quasi--periodic with Diophantine frequencies.
\end{theorem}

\vskip.1in
\noi
Following \cite{Fej04} (who deals with Herman's $C^{\infty}$ KAM Theory
), the strategy of the proof in \cite{ChPu08} consists in applying the previous KAM Theory  to a suitably modified Hamiltonian function, which is obtained from the planetary Hamiltonian expressed in Delaunay--Poincar\'e variables  by adding a commuting Hamiltonian. The necessity of modifying the Hamiltonian function is that when the planetary Hamiltonian is put in Poincar\'e--Delaunay variables, the frequencies $\O$ of (\ref{first order theory}) correspond the $2N$--dimensional vector  of the {{\rm DP}--\sl secular frequencies}
$$(\s,\zeta)=\Big((\s_1,\cdots,\s_N),(\zeta_1,\cdots,\zeta_N)\Big)$$
defined, respectively, as the eigenvalues of the two ``horizontal'' and ``vertical'' quadratic forms 
\beqa{PoincDel1}
\arr{
{\cal Q}_h\cdot \eta^2:=\sum_{1\leq j<k\leq N}m_jm_k\left(C_1(a_j,a_k)\left(\frac{\eta_j^2}{\L_j}+\frac{\eta_k^2}{\L_k}\right)+2C_2(a_j,a_k)\frac{\eta_j\eta_k}{\sqrt{\L_j\L_k}}\right)\\
{\cal Q}_v\cdot p^2:=-\sum_{1\leq j<k\leq N}m_jm_kC_1(a_j,a_k)\left(\frac{p_j}{\L_j}-\frac{p_k}{\L_k}\right)^2
}
\eeqa
which are been proved \cite{Fej04} to satisfy, together with the mean motions ${\rm n}:=\partial h_{\rm plt}$   the {\sl only} two independent linear combinations, usually called {\sl secular resonances} 
\beq{Fej only res}
\sum_{1\leq i\leq N}(\s_i+\zeta_i)=0\ ,\quad \zeta_N=0\ \eeq
and the R\"ussmann condition (\ref{Russmann}) 
is clearly violated.
Adding a commuting Hamiltonian  makes the above 
non degeneracy condition (\ref{Russmann}) verified. The final result is reached with the use of an abstract argument: invariant ergodic tori for the modified Hamiltonian are recognized to be  invariant manifolds for the original Hamiltonian.

\vskip.1in
\noi
The use of the regularized (fully) reduced Deprit variables  provides a direct application of the KAM machinery of \cite{ChPu08} because {no secular resonance appears}. We recover then a result already found with a different technique in the 2007 revised version of the paper by J. F\'ejoz \cite{Fej04}.

\vskip.1in
\noi
\begin{theorem}\label{3N-2 KAM tori}
Consider, in $\real^3$ a star of mass $\bar m_0$ and $N\geq 2$ planets of mass $\m\bar m_1$, $\cdots$, $\m\bar m_N$, interacting only through gravity.  Let $a_i$ denote the instantaneous major semiaxis of the $i^{th}$ planet and let $\e$ be an upper bound of the instantaneous eccentricity and inclination of the planets. Then, there exists $\d_*>0$, $\e_*>0$ and $\m_*>0$ such that, if $a_i/a_{i+1}<\d_*$,  $0<\e<\e_*$ and $0<\m<\m_*$, there exists a positive measure set of phase space points whose time evolution lies on real--analytic, $3N-2$ dimensional invariant tori; the time evolution being quasi--periodic with $3N-2$  Diophantine frequencies. Furthermore, during the motion, eccentricities and inclinations are bounded by $C\sqrt{\m}$.
\end{theorem}
{\bf Proof.}\ It is a corollary of  Theorem \ref{Russmann Thm} above and Lemma \ref{non resonance} of the following section.
\subsection{R\"ussmann Non--Degeneracy and $(3N-2)$--Dimensional KAM Tori}\label{The Arnol'd Pyartli Condition in the Space}
{\bf Remarks on notations.} Referring especially to paragraphs \ref{Regularization}, \ref{Partial Reduction}, throughout all this section,
\begin{itemize}
\item[(i)]
We disregard the (cyclic, for ${\cal H}_{\rm plt}$) Deprit variables   $(P,Q)$, on which we will always think to lift the maps we will discuss, extending them through the identity map. Quite abusively, we do not change the name of the domains ${\cal D}_{\rm r}$,  ${\cal D}_{\rm pr}$ of the fully, partially reduced regularized Deprit variables.

\item[(ii)] We denote the set of {\sl fully reduced regularized Deprit variables} as
$$(\L,\l)\ ,\quad (\eta,\xi)\ ,\quad(p,q)\ ,\quad(G,g)$$
($g$ cyclic for ${\cal H}_{\rm plt}$), hence, in particular,  {\sl $p$, $q$ have dimension $N-2$}. The planetary Hamiltonian put in {\sl fully reduced regularized Deprit variables}
is denoted
$${\cal H}_{\rm plt,fr}=h_{\rm plt}+\m f_{\rm plt,fr}$$ where, as usual
$$h_{\rm plt}=-\sum_{1\leq i\leq N}\frac{\tilde m_i^3\hat m_i^2}{2\L_i^2}\ .$$

\item[(iii)] The set of {\sl partially reduced regularized Deprit variables} with
$$(\L,\bar\l)\ ,\quad (\bar\eta,\bar\xi)\ ,\quad(\bar p,\bar q)
$$
hence, with  {\sl with $\bar p$, $\bar q$ of dimension $(N-1)$}. The planetary Hamiltonian put in {\sl  partially reduced regularized Deprit variables}
is denoted
$${\cal H}_{\rm plt,pr}=h_{\rm plt}+\m f_{\rm plt,pr}$$
\end{itemize}
We start with the planetary Hamiltonian written in fully reduced variables
$${\cal H}_{\rm plt,fr}=h_{\rm plt}+\m f_{\rm plt,fr}\ .$$

\begin{lemma}\label{non resonance}
For a sufficiently small $\d_*$, in the set $\tilde{\cal D}_{\rm r}$ of $\dst \Big((\L,\l),(\eta,\xi),(p,q),(G,g)\Big)$ of $\dst (\real_+^N\times \torus^N)\times (\real^{N}\times \real^N)\times(\real^{N-2}\times \real^{N-2})\times(\real_+\times \torus)$
such that
\beqano
&& a(\L)\in{\cal A}\ ,\quad \d^2:=\dst \sum_{1\leq i\leq N}\L_i-G<\d_*^2\ ,\quad |(\eta,\xi,p,q)|_2<2\,\d\ ,
\eeqano
with ${\cal A}$ the  set of semimajor axes
$${\cal A}:=\Big\{a=(a_1,\cdots,a_N)\in \real^N:\quad 0<a_1<a_2<\cdots<a_N\Big\}\ ,$$
there exists a symplectic real--analytic  change of variable 
$$\phi:\quad \tilde{\cal D}_{\rm r}\to {\cal D}_{\rm r}$$
which leaves $G$, $\L$ unvaried and puts $\bar f_{\rm  plt,fr}$ into the form
\beqano
\bar f_{\rm plt}&:=&\bar f_{\rm plt,fr}\circ\phi=\ovl{f_{\rm plt,fr}\circ\phi}\nonumber\\
&=&\bar f_0(\L,G)+\sum_{1\leq i\leq N}{\rm s}_i(\L,G)\frac{\eta_i^2+\xi_i^2}{2}+\sum_{1\leq i\leq N-2}{\rm z}_i(\L,G)\frac{p_i^2+q_i^2}{2}+{\rm O}(3)
\eeqano
where,  for any fixed $G\in \real_+$, the ``secular frequencies'' ${\rm s}=({\rm s}_1$, $\cdots$, ${\rm s}_N)$, ${\rm z}=({\rm z}_1$, $\cdots$, ${\rm z}_{N-2})$ together with the mean motions ${\rm n}=({\rm n}_1,\cdots,{\rm n}_N):=\partial h_{\rm plt}$ do not satisfy any linear relation in any simply connected subset ${\cal V}_G$ of a suitable subset ${\cal U}_G$ with full measure of 
$${\cal A}_G:=\Big\{a(\L)\in{\cal A}\ ,\quad \dst \sum_{1\leq i\leq N}\L_i-G<\d_*^2\Big\}\ .$$ 
\end{lemma}

\vskip.1in
\noi
{\bf Proof.}\ We  discuss only the case $N\geq 3$, since the case $N=2$ is well understood.
\footnote{As already remarked, for $N=2$, the full Deprit reduction corresponds to the Jacobi reduction and the two spatial secular frequencies ${\rm s}_1$, ${\rm s}_2$ of the spatial three body problem  are manifestly related (see \cite{Rob95}) to the frequencies of the plane problem $\s_1$, $\s_2$ in Delaunay--Poincar\'e variables by
$${\rm s}_1=2\s_1+\s_2+{\rm O}\left(\d^2\right)\ ,\qquad {\rm s}_2=\s_1+2\s_2+{\rm O}\left(\d^2\right)\ .$$
Hence, ${\rm s}_1$, ${\rm s}_2$ have the desired property since $\s_1$, $\s_2$ have it, as proved in \cite{Fej04}. 
}

\vskip.1in
\noi
{\bf Step 1:}\ {\sl partial reduction (or full regularization).}\\  Let $\phi_{\rm pr}$  the map ``partial reduction map'' $\phi_{\rm pr}$ which acts as described in (eq. (\ref{partial reduction})) Section \ref{Partial Reduction}. This leads ${\cal H}_{\rm plt,fr}=h_{\rm plt}+\m f_{\rm plt,fr}$ to ${\cal H}_{\rm plt,pr}=h_{\rm plt}+\m f_{\rm plt,pr}$, where $\bar f_{\rm plt,pr}$ is as in Lemma \ref{expansion to order two}:

\beqano\bar f_{\rm  plt,pr}=C_0(m,a)+{\cal Q}^*_h\cdot\frac{\bar\eta^2+\bar\xi^2}{2}+{\cal Q}^*_v\cdot\frac{\bar p^2+\bar q^2}{2}+\bar f^4_{\rm plt, pr}\ .\eeqano
 
\vskip.2in
\noi
{\bf Step 2:}\ {\sl diagonalization of ${\cal Q}_{h}^*$, ${\cal Q}_{v}^*$}.\ Let  $\r^*_h$, $\r^*_v$ the unitary matrices which leave ${\cal Q}^*_h$, ${\cal Q}^*_v$ diagonal:
$${\cal Q}_{h}^*\cdot\tilde\eta^2=\sum_{1\leq i\leq N}s_i{\tilde\eta_i^2}\ ,\quad {\cal Q}_{v}^*\cdot\tilde p^2=\sum_{1\leq i\leq N-1}z_i{\tilde p_i^2}$$
where
$$\tilde\eta:={\r^*_h}^t\bar\eta\ ,\quad \tilde p:={\r^*_v}^t\bar p\ .$$ Then, the transformation $\phi_{\rm diag}$ which leaves $\L$ unvaried and acts on $\tilde p$, $\tilde q$, $\tilde\l$ as
\beq{diag}\phi_{\rm diag}:\quad\arr{\bar\eta:=\r^*_h\tilde\eta\\
\bar\xi:=\r^*_h\tilde\xi
}\ ,\quad \arr{\bar p:=\r^*_v\tilde p\\
\bar q:=\r^*_v\tilde q}\ ,\quad \bar\l=\tilde\l+\varphi(\L,\tilde p, \tilde q)\eeq  
(where $(\L,\bar p, \bar q)\to \varphi(\L,\bar p, \bar q)$ is a suitable shift which makes $\phi_{\rm diag}$ symplectic) puts $\bar f_{\rm  plt,pr}$ into the form
$$\bar f_{\rm diag}:=\bar f_{\rm  plt,pr}\circ\phi_{\rm diag}=\ovl{f_{\rm  plt,pr}\circ\phi_{\rm diag}}=C_0(m,a)+\sum_{1\leq i\leq N}s_i\frac{\tilde\eta_i^2+\tilde\xi_i^2}{2}+\sum_{1\leq i\leq N-1}z_i\frac{\tilde p_i^2+\tilde q_i^2}{2}+\tilde f^4_{\rm diag}\ .$$
where $\tilde f^4_{\rm diag}={\rm O}(4)$ and, as proved in Proposition \ref{true non resonance}, $s=(s_1,\cdots,s_N)$ and $z=(z_1,\cdots,z_N)$ do not satisfy any other linear condtion than the Herman's resonance. Notice that, since $\r^*_h$, $\r^*_v$ are unitary, the angular momentum $G$, in their terms has just the same expression as in full regularized variables:
$$G=\sum_{1\leq i\leq N}\left(\L_i-\frac{\tilde \eta_i^2+\tilde\xi_i^2}{2}\right)-\sum_{1\leq i\leq N-1}\frac{\tilde p_i^2+\tilde q_i^2}{2}$$
We are then ``justified'' if we do  not  change the name of the variable $G$ we introduce into the following step.

\vskip.2in
\noi
{\bf Step 3:}\ {\sl full reduction}.\ Apply now $\phi_{\rm fr}:=\phi_{\rm pr}^{-1}$, \ie, put
\beqano
&&  \arr{\dst
\tilde p_{N-1}=\sqrt{2\left(\sum_{1\leq i\leq N}\L_i-G-\sum_{1\leq i\leq N}\frac{\check \eta_i^2+\check \xi_i^2}{2}-\sum_{1\leq i\leq N-2}\frac{\check p_i^2+\check q_i^2}{2}\right)}\cos{\check g}\\
\tilde q_{N-1}=-\sqrt{2\left(G-\sum_{1\leq i\leq N}\L_i-\sum_{1\leq i\leq N}\frac{\check \eta_i^2+\check \xi_i^2}{2}-\sum_{1\leq i\leq N-2}\frac{\check p_i^2+\check q_i^2}{2}\right)}\sin{\check g}\\
} \nonumber\\ \nonumber\\ 
&&   \tilde\l_i=\check\l_i+\check g \nonumber\\ 
&&  \left(\begin{array}{lr}
\tilde\eta_i\\
\tilde\xi_i
\end{array}
\right)={\cal R}_{\rm z}(-\check g)
\left(\begin{array}{lr}
\check\eta_i\\
\check\xi_i
\end{array}
\right)\qquad \left(\begin{array}{lr}
\tilde p_j\\
\tilde q_j
\end{array}
\right)={\cal R}_{\rm z}(-\check g)
\left(\begin{array}{lr}
\check p_j\\
\check q_j
\end{array}
\right)
\eeqano
This carries $\bar f_{\rm diag}$ to
$$\bar f_{\rm fr}:=\bar f_{\rm diag}\circ\phi_{\rm fr}=\ovl{f_{\rm diag}\circ\phi_{\rm fr}}=\check C_0(m,a,G)+\sum_{1\leq i\leq N}{\rm s}^0_i\frac{\check\eta_i^2+\check\xi_i^2}{2}+\sum_{1\leq i\leq N-2}{\rm z}^0_i\frac{\check p_i^2+\check q_i^2}{2}+\check f\ ,$$
where
$$\arr{\check C_0(m,a,G):=C_0(m,a)+s_{N-1}\left(\sum_{1\leq i\leq N}\L_i-G\right)\\
\check f:= \tilde f^4_{\rm diag}\circ \phi_{\rm fr}\\
{\rm s}_i^0:=s_i-z_{N-1}\quad 1\leq i\leq N\\
{\rm z}_i^0:=z_i-z_{N-1}\quad 1\leq i\leq N-2\\
}$$
\begin{lemma}
The functions ${\rm s}_1^0$, $\cdots$, ${\rm s}_N^0$, ${\rm z}_1^0$, $\cdots$, ${\rm z}_{N-2}^0$ do not satisfy  any linear relation in any open, simply connected set ${\cal V}\subset{\cal U}$ for a suitable open set ${\cal U}\subset{\cal A}$ with full measure.
\end{lemma}
{\bf Proof.}\ Let ${\cal U}\subset {\cal A}$ the open set with full measure where Proposition \ref{true non resonance} holds and assume that, in some pont of some simply connected set $\bar{\cal V}\subset{\cal U}$, we had 
$$\sum_{1\leq i\leq N}c_i{\rm s}_i^0+\sum_{1\leq i\leq N-2}g_i{\rm z}_i^0=\sum_{1\leq i\leq N}c_i(s_i-z_{N-1})+\sum_{1\leq i\leq N-2}g_i(z_i-z_{N-1})=0\ ,
$$
with $(c_1$, $\cdots$, $c_N$, $g_1$, $\cdots$, $g_{N-2})\in \real^{2N-2}\setminus\{0\}$. Then, by Proposition \ref{true non resonance}
$$c_1=\cdots=c_N=g_1=\cdots=g_{N-2}=-\sum_{1\leq i\leq N}c_i-\sum_{1\leq i\leq N-2}g_i\ ,$$
which is a contraddiction.

\vskip.2in
\noi
{\bf Step 4:}\ {\sl recentering $f_{\rm fr}$ at its equilibrium point}.\ For small values of $\d^2:=\sum_{1\leq i\leq N}\L_i-G$ we find a new equilibrium point $\check z_{\rm eq}=(\check \eta_{\rm eq},\check \xi_{\rm eq},\check p_{\rm eq},\check q_{\rm eq})$ for $\bar f_{\rm fr}$, which is ${\rm O}(\d)$. 
Rescale, in fact, the variables as
$$\check \eta=2\d\,\hat p\ ,\quad \check \xi=2\d\,\hat q\ ,\quad \check p=2\d\,\hat p\ ,\quad \check q=2\d\,\hat q$$
and then discuss equation
$$\d^{-2}\partial_{\check z} \bar f_{\rm fr}=0\quad \textrm{where}\quad \check z=(\check \eta,\check \xi,\check p,\check q)\ .$$
by an Implicit Function Theorem argument.\\
Perform then the change of variable
$$\phi_{\rm eq}:\quad\check z=z_*+\check z_{\rm e}\quad \check g=g_*+\check\psi(\L,G)\ ,\quad \check\l_i=\l_i^*+\check\varphi_i(\L,G)$$
leaving the remaining variables unvaried, where $\check\psi(\L,G)$, $\check\varphi_i(\L,G)$ are suitable shifts which make $\phi_{\rm e}$ symplectic. \\
The result then follows after a suitable symplectic diagonalization of the Hessian matrix of $\dst \bar f_{\rm eq}:=\bar f_{\rm fr}\circ\phi_{\rm eq}=\ovl{f_{\rm fr}\circ\phi_{\rm eq}}$, which gives linear invariants ${\rm s}_1$, $\cdots$, ${\rm s}_N$, ${\rm z}_1$, $\cdots$, ${\rm z}_{N-2}$ $\d$--close to ${\rm s}_1^0$, $\cdots$, ${\rm s}_N^0$, ${\rm z}_1^0$, $\cdots$, ${\rm z}_{N-2}^0$, hence, with the desired property.
\appendix
%
\newpage
\section{Proof of the  Averaging Theorem (Lemma \ref{average})}
\setcounter{equation}{0}
\label{app:average}
\renewcommand{\theequation}{\ref{app:average}.\arabic{equation}}
\begin{lemma}\label{derivatives lemma}
Let $D$, $A_0$, $A_1$, $\cdots$, $A_N$ square complex matrices of order $n$, with $D$ non singular, such that 
$$\|D(A_i-\id_n)\|\leq \e_i\ ,\quad i=0,\cdots,N \ .$$
Then, if $\|\cdot\|$ is a norm on ${\rm Mat}_{n\times n}^{\complex}$ such that $\|A\,B\|\leq \|A\|\,\|B\|$, 
\begin{eqnarray*}
\|D(A_0A_1A_2\cdots A_N-\id_n)\|&\leq&\e_0+(1+\e_0\|D^{-1}\|)\e_1+\cdots\nonumber\\
&+&(1+\e_0\|D^{-1}\|)\cdots (1+\e_{N-1}\|D^{-1}\|)\e_N\ .
\end{eqnarray*}
\end{lemma}
{\bf Proof.}\ Let
$$T:=A_1\cdots A_N-\id_n\ .$$
Then, writing
$$D(A_0A_1\cdots A_N-\id_n)=D(A_0-\id_n)\,(\id_n+T)+DT\ ,$$
we find
\begin{eqnarray}\label{first iterate}
\|D(A_0A_1\cdots A_N-\id_n)\|&\leq& \|D(A_0-\id_n)\|(1+\|T\|)+\|DT\|\nonumber\\
&\leq&\e_0+(1+\e_0\|D^{-1}\|)\|DT\|\nonumber\\
&=&\e_0+(1+\e_0\|D^{-1}\|)\|D(A_1\cdots A_N-\id_n)\|
\end{eqnarray}
The Lemma is then proved after $N$ iterations of (\ref{first iterate}).
\vskip.1in
\noindent
\begin{lemma}[Iterative Lemma]\label{iterative lemma}
Let $\bar n+\hat n=n$, $0<2\r<r$ and $0<2\s<s$. Suppose that the Hamiltonian
$$H=h+g+f$$
is real--analytic on $\cP_{r,s}:=\cI_{r}\times \torus^n_s$ with $\o:=\partial\,h$ verifying
\begin{eqnarray}\label{complex non resonance}
|\o(I)\cdot k|\geq \left\{\begin{array}{lrr}
\bar\a\quad \textrm{for}\quad k=(\bar k,\hat k)\in \integer^{\bar n}\times \integer^{\hat n}\setminus \L\quad \bar k\neq 0\ ,\quad |k|_1\leq K \\
\hat\a \quad \textrm{for}\quad k=(0,\hat k)\in \{0\}_{\integer^{\bar n}}\times \integer^{\hat n}\setminus \L\quad |k|_1\leq K\\
\end{array}
\right.
\end{eqnarray}
for $I\in \cI_{r}$, and $f$ so small that
\begin{eqnarray}\label{smallness}
\|f\|_{r,s}<\frac{\a\r\s}{2}\ ,\quad \textrm{where}\quad \a:=\min\{\bar\a,\hat\a\}\ .
\end{eqnarray}
Then, there exists a real--analyitic, symplectic transformation
$$\Phi:\quad \cP_{r-2\r,s-2\s}\to \cP_{r,s}$$
which carries $H$ into
$$H_+:=H\circ\Phi(I,\varphi)=h+g_++f_+$$
with
$$g_+=g+P_\L T_K f$$
and
\begin{eqnarray}\label{f+}\|f_+\|_{r-2\r,s-2\s}\leq \left(1-\frac{\|f\|_{r,s}}{\a\r\s/2}\right)^{-1}\left(\frac{\|g\|_*+\|f\|_{r,s}}{\a\r\s}+e^{-K\s}\right)\,\|f\|_{r,s}\end{eqnarray}
where
$$\|g\|_*=\frac{1}{2}\,\sum_j\left(\frac{\r}{r_j-r+\r}+\frac{\s}{s_j-s+\s}\right)\|g_j\|_{r_j,s_j}\ ,$$
if $g=\sum_j\,g_j$, with terms $g_j$ bounded on $\cP_{r_j,s_j}\supset  \cP_{r,s}$. Moreover,
\begin{eqnarray}\label{close to id}
|W^{\bar n}_{\bar\a^{-1},\a^{-1}}W_{\r,\s}(\Phi-\textrm{\rm id})|_{\cal P}\ ,\qquad \frac{1}{2}\,\|W^{\bar n}_{\bar\a^{-1},\a^{-1}}(W_{\r,\s}\,D\Phi\,W^{-1}_{\r,\s}-\id_{2n})\|_{\cal P}\leq \frac{\|f\|_{r,s}}{\r\s}
\end{eqnarray}
uniformly on $\cP_{r-2\r,s-2\s}$.
\end{lemma}
Lemma \ref{iterative lemma} is a useful remake  of the {\sl Iterative Lemma} of  \cite{POSCH93}. We outline the differences. 
\begin{itemize}
\item[(i)]  
In \cite{POSCH93}, instead than (\ref{complex non resonance}), the {\sl real} nonresonance for $\o$ up to order $K$ and a ``smallness'' assumption for $r$: 
\begin{eqnarray}\label{real non resonance}
(a):\ |\o(I)\cdot k|\geq \a\ \textrm{for}\ I\in \cI\ , \ k\in \integer^n\setminus \L\ ,\ |k|_1\leq K\ ,\  (b):\ r\leq \frac{\a}{pMK}\nonumber\\
\end{eqnarray}
are required, where $p>1$ is a prefixed number.
But, in the proof, $(a)$ and $(b)$ are used only to prove 
$$|\o(I)\cdot k|\geq \frac{\a}{q}\quad \textrm{for}\quad I\in \cI_{r}\ , \quad k\in \integer^n\setminus \L\ ,\quad |k|_1\leq K$$
which is next needed. So, the Lemma remains true when the assumption
$$|\o(I)\cdot k|\geq \a\quad \textrm{for}\quad I\in \cI_{r}\ , \quad k\in \integer^n\setminus \L\ ,\quad |k|_1\leq K$$
replaces (\ref{real non resonance}) and 
$\a$ replaces $\a/q$ in all its occurences. But, (\ref{complex non resonance}) obviously implies
\begin{eqnarray*}
|\o(I)\cdot k|\geq \a:=\min\{\bar\a,\ \hat\a\}\quad \textrm{for}\quad I\in \cI_{r}\ , \quad k\in \integer^n\setminus \L\ ,\quad |k|_1\leq K\ .
\end{eqnarray*}
It now is enough observing that (\ref{smallness}), (\ref{f+}) are just the same of \cite{POSCH93}, with $\a$ replacing $\a/q$.
\item[(ii)] For what concernes (\ref{close to id}), in \cite{POSCH93}, taking into account  {\sl(i)}, we find
\begin{eqnarray}\label{closeness to id}
|W_{\r,\s}(\Phi-\textrm{\rm id})|_{\cal P^*}\leq \frac{\|f\|_{r,s}}{\a\r\s}
\end{eqnarray}
where
$$|(y,x)|_{\cP^*}:=\max\{|y|_1,\ |x|_{\infty}\}\ .$$
In particular, (\ref{closeness to id}) holds when
\footnote{$\P_{\hat I}$, $\cdots$ denotes the projection ove the $\hat I$, $\cdots$ variables.}
$$((0,\P_{\hat I}\Phi),\P_{\varphi}\Phi)$$ replaces $$\Phi=((\P_{\bar I}\Phi,\P_{\hat I}\Phi),\P_{\varphi}\Phi)\ .$$
In order to obtain an estimate for $|\r^{-1}(\P_{\bar I}\Phi-\textrm{\rm id}_{\bar n})|_1$, we recall that, in \cite{POSCH93}, 
$\Phi$ is constructed as the time $1$ map of the flow $X_{\phi}^t$ of a hamiltonian vectorfield $X_{\phi}=(-\partial_{\varphi}\phi,\partial_I\phi)$, where $\phi$ is the trigonometric polynomial
$$\phi(I,\varphi)=\sum_{k\in \integer^{\bar n}\times\integer^{\hat n}\setminus\L,\ |k|\leq K}\frac{f_k(I)}{i\o(I)\cdot k}\,e^{ik\cdot\varphi}\ .$$
But, taking into account the non resonance assumptions (\ref{complex non resonance}), we may split $\phi$ as
$$\phi(I,\varphi)=\bar \phi(I,\varphi)+\hat\phi(I,\hat\varphi)$$
where
$$\bar\phi(I,\varphi)=\sum_{k=(\bar k,\hat k)\in \integer^{\bar n}\times\integer^{\hat n}\setminus\L,\ |k|\leq K,\ \bar k\neq 0}\frac{f_k(I)}{i\o(I)\cdot k}\,e^{ik\cdot\varphi}$$
and $\hat\phi=\phi-\bar\phi$ does not depend from the $\bar\varphi$--variables.
Hence, the projection $\bar X_{\phi}:=\partial_{\bar{\varphi}}\phi=\partial_{\bar{\varphi}}\bar\phi$ of the vectorfield $X_{\phi}$ over $\complex^{\bar n}\times\{0\}_{\complex^{2n-\bar n}}$, by the General Cauchy Inequality, is bounded as
$$|\bar X_{\phi}|_1\leq \frac{\|f\|_{r,s}}{\bar\a\s}$$
uniformly on $V_{r-3\r/2,s-3\s/2}$, so, joining this result with (\ref{closeness to id}), we find
$$|W^{\bar n}_{\bar\a^{-1},\a^{-1}}W_{\r,\s}(\Phi-I_{2n})|_{\cal P}\leq \frac{\|f\|_{r,s}}{\r\s}\ .$$
The second equation in (\ref{close to id}), then, follows by the General Cauchy Inequality,
uniformly on $V_{r-2\r,s-2\s}$.
\end{itemize}
\vskip.1in
\noindent
{\bf Proof of Lemma \ref{average}.}
We distinguish two cases:
$${\rm (a)} \quad 48\,\frac{\|f\|_{r,\bar s+s}}{\a\,r\,s}<e^{-Ks/6}\ ,\qquad\qquad\quad \textrm{\rm (b)}\quad 48\,\frac{\|f\|_{r,\bar s+s}}{\a\,r}\geq e^{-Ks/6}\ .$$
The case (a) is trivial, because we apply the Iterative Lemma with $g\equiv 0$ and parameters
$$\r=\r':=\frac{r}{4}\ ,\qquad \s={\s'}:=\frac{s}{3}\ .$$
This is made possible, because (using $Ks\geq 6$)
$$\|f\|_{r,\bar s+s}< \frac{1}{2^7}\frac{\a\,r\,s}{Ks}\leq\frac{1}{6\cdot2^7}\a\,r\,s=\frac{\a\,{\r'}\,{\s'}}{2^6\,}<\frac{\a{\r'}{\s'}}{2}\ .$$
Letting, then,
$$r'=r-2{\r'}=\frac{r}{2}\ ,\quad s'=s-2\,{\s'}=\frac{s}{3}$$
we see that the map
$\Phi:\quad \cP_{r',\bar s+s'}\to \cP_{r,\bar s+s}$
verifies
\beqano
&&|W^{\bar n}_{\bar\a^{-1},\a^{-1}}W_{\r',\s'}(\Phi-\textrm{id})|_{\cal P}\ ,\nonumber\\
&& \frac{1}{2}\,\|W^{\bar n}_{\bar\a^{-1},\a^{-1}}(W_{\r',\s'}\,D\Phi\,{W}_{\r',\s'}^{-1}-\id_{2n})\|_{\cal P}\nonumber\\
&\leq& \frac{12}{rs}\,\|f\|_{r,\bar s+s}\nonumber\\
&\leq&\frac{2K}{r}\,\|f\|_{r,\bar s+s}\ ,
\eeqano
and 
$$H\circ \Phi=h+f_0+f_1$$
where $f_0=P_{\L}T_K\,f$. Furthermore,  
\begin{eqnarray*}
\|f_1\|_{r',\bar s+s'}&\leq&\left(1-\frac{\|f\|_{r,\bar s+s}}{\a{\r'}{\s'}/2}\right)^{-1}\left(\frac{\|f\|_{r,\bar s+s}}{\a{\r'}{\s'}}+e^{-Ks/3}\right)\,\|f\|_{r,\bar s+s}\nonumber\\
&\leq& (1-2^{-5})^{-1}\,\left(\frac{\|f\|_{r,\bar s+s}}{\a{\r'}{\s'}}+e^{-Ks/3}\right)\,\|f\|_{r,\bar s+s}\nonumber\\
&=&\frac{32}{31}\,\left(\frac{\|f\|_{r,\bar s+s}}{\a\,r\,s/(12)}+e^{-Ks/3}\right)\,\|f\|_{r,\bar s+s}\nonumber\\
&\leq&\frac{32}{31}\left(\frac{1}{4}+\frac{1}{e}\right)e^{-Ks/6}\,\|f\|_{r,\bar s+s}\nonumber\\
&<&e^{-Ks/6}\,\|f\|_{r,\bar s+s}\ .
\end{eqnarray*}
The lemma is proved with $\Psi=\Phi$. 

\nl
In case (b), we move in $N+1$ steps. First, we apply the Iterative Lemma with $g\equiv0$ and parameters
$$\r=\r_0:=\frac{r}{8}\ ,\qquad \s=\s_0:=\frac{s}{6}\ ,$$ 
thanks to the inequality
$$\|f\|_{r,\bar s+s}\leq \frac{1}{2^7}\frac{\a\,r\,s}{Ks}\leq\frac{1}{6\cdot2^7}\a\,r\,s=\frac{\a\,\r_0\,\s_0}{16}\ .$$
Letting, 
$$r_1=r-2\r_0=\frac{3}{4}\,r\ ,\quad s_1=s-2\s_0=\frac{2}{3}\,s$$
we find
$$\Phi_0:\quad \cP_{r_1,\bar s+s_1}\to \cP_{r,\bar s+s}$$
verifying
\beqano
&&|W^{\bar n}_{\bar\a^{-1},\a^{-1}}W_0(\Phi_0-\textrm{id})|_{\cal P}\ ,\nonumber\\
&&\frac{1}{2}\,\|W^{\bar n}_{\bar\a^{-1},\a^{-1}}(W_0\,D\Phi_0{W_0}^{-1}-\id_{2n})\|_{\cal P}\nonumber\\
&\leq&\frac{1}{\r_0\s_0}\,\|f\|_{r,\bar s+s}=\frac{48}{rs}\,\|f\|_{r,\bar s+s}\ ,
\eeqano
where $W_0:=W_{\r_0,\s_0}$, and
$$H_1:=H\circ \Phi_0=h+g_1+f_1$$
with 
\begin{eqnarray}\label{f1estimates}
\|f_1\|_{r_1,\bar s+s_1}&\leq& \left(1-\frac{\|f\|_{r,\bar s+s}}{\a\r_0\s_0/2}\right)^{-1}\left(\frac{\|f\|_{r,\bar s+s}}{\a\r_0\s_0}+e^{-K\s_0}\right)\,\|f\|_{r,\bar s+s}\nonumber\\
&\leq& \frac{8}{7}\,\left(\frac{\|f\|_{r,\bar s+s}}{\a\r_0\s_0}+e^{-K\s_0}\right)\,\|f\|_{r,\bar s+s}\nonumber\\
&=& \frac{8}{7}\,\left(\frac{\|f\|_{r,\bar s+s}}{\a\,r\,s/48}+e^{-Ks/6}\right)\,\|f\|_{r,\bar s+s}\nonumber\\
&\leq&\frac{16}{7}\,\frac{\|f\|_{r,\bar s+s}}{\a\,r\,s/48}\,\|f\|_{r,\bar s+s}\nonumber\\
&=&\frac{1}{7}\,\frac{\|f\|_{r,\bar s+s}}{\a\r_0\s_0/(16)}\,\|f\|_{r,\bar s+s}\nonumber\\
&\leq&\frac{\|f\|_{r,\bar s+s}}{7}\ .
\end{eqnarray}
Now, let $N$ an integer number. Our aim is to apply the Iterative Lemma $N$ times, each with parameters 
\begin{eqnarray}\label{rhos0}
\r={\r_N}:=\frac{r}{8N}\ ,\quad \s={\s_N}:=\frac{s}{4N}\ ,
\end{eqnarray}
so as to construct, at each step, a symplectic, analytic transformation
$$\Phi_i:\quad \cP_{r_{i+1},\bar s+s_{i+1}}\to \cP_{r_i,\bar s+s_i}\ ,\quad i=1,\cdots,N$$
where
$$r_i=r_1-2(i-1)\,\r_N\ ,\quad s_i=s_1-2(i-1)\,\s_N\ ,\quad i=1,\cdots,N+1\ ,$$
verifying

\beqano
&&|W^{\bar n}_{\bar\a^{-1},\a^{-1}}W_N(\Phi_i-\textrm{id})|_{\cal P}\ ,\nonumber\\
&& \frac{1}{2}\,\|W^{\bar n}_{\bar\a^{-1},\a^{-1}}(W_N\,D\Phi_i\,{W_N}^{-1}-\id_{2n})\|_{\cal P}\nonumber\\
&\leq&\frac{1}{\r_N\s_N}\,\|f_{i}\|_{r_{i},\bar s+s_{i}}\ ,
\eeqano
where
$W_N:=W_{\r_N,\s_N}$, and
$$H_{i+1}:=H\circ \Phi_0\circ \cdots \circ \Phi_i=h+g_{i+1}+f_{i+1}\ ,\quad i=1,\cdots,N$$
with
$$\|f_i\|_{r_i,\bar s+s_i}\leq \left(\frac{1}{4}\right)^{i-1}\|f_1\|_{r_1,\bar s+s_1}\leq \left(\frac{1}{4}\right)^{i}\|f\|_{r,\bar s+s}\ .$$
Consequentely, after $N$ steps, we will find 
$H_{N+1}=H\circ \Psi_N=h+g_{N+1}+f_{N+1}$ analytic on $\cP_{r/2,\bar s+s/6}$
$$f_{N+1}\leq \left(\frac{1}{4}\right)^{N+1}\|f\|_{r,\bar s+s}\leq e^{-Ks/6}\|f\|_{r,\bar s+s}$$
provided $N$ is sufficiently large:
$$N+1\geq\frac{Ks}{12\log{2}}\ .$$
Choice, thus, $N$ such in a way that
$$N<\frac{Ks}{12\log{2}}\leq N+1\ .$$
If $N=0$, the theorem is proved with $\Psi=\Phi_0$. Otherwise, if $N\geq 1$, we verify, by induction, inequalities
\begin{eqnarray}\label{inequalities}
\|f_{i}\|_{r_{i},\bar s+s_{i}}&\leq& \frac{\|f_{{i}-1}\|_{r_{{i}-1},\bar s+s_{{i}-1}}}{4}\ ,\quad \|f_{i}\|_{r_{i},\bar s+s_{i}}<\frac{\a\,r\,s}{2^{13}N^2}=\frac{\a{\r_N}{\s_N}}{2^8}\ ,\nonumber\\
\|g_{i}-g_{i-1}\|_{{r_{i},\bar s+s_{i}}}&\leq& \|f_{i-1}\|_{r_{i-1},\bar s+s_{i-1}}\ ,\quad i=1,\cdots,N
\end{eqnarray}
For $i=1$, we have
$$\|g_1-g_0\|_{r_1,\bar s+s_1}=\|g_1\|_{r_1,\bar s+s_1}=\|P_\L T_Kf\|_{r_1,\bar s+s_1}\leq \|f\|_{r_1,\bar s+s_1}$$
and, by (\ref{f1estimates}),
\begin{eqnarray}\label{f1est}
\|f_1\|_{r_1,\bar s+s_1}\leq 2^7\frac{\|f\|_{r,\bar s+s}^2}{\a\,r\,s}\leq\frac{\a\,r\,s}{2^7(Ks)^2}<\frac{\a\,r\,s}{2^{13}N^2}=\frac{\a{\r_N}{\s_N}}{2^8}\leq \frac{\|f\|_{r,\bar s+s}}{7}
\end{eqnarray}
provided
$${\|f\|_{r,\bar s+s}}\leq\frac{\a\,r}{2^7K}\ .$$
Assume, now, that (\ref{inequalities}) hold for a given $i<N$.
Then, Lemma \ref{iterative lemma} is applicable once again and we find $$\Phi_i:\ \cP_{r_{i+1},\bar s+s_{i+1}}\to \cP_{r_{i},\bar s+s_{i}}$$ 
such that $$H_{i+1}:=H_i\circ \Phi_i=h+g_{i+1}+f_{i+1}$$ with $g_{i+1}$ and $f_{i+1}$ verifying
$$\|g_{i+1}-g_{i}\|_{r_{i+1},\bar s+s_{i+1}}=\|P_\L T_Kf_{i}\|_{r_{i+1},\bar s+s_{i+1}}\leq \|f_{i}\|_{r_{i+1},\bar s+s_{i+1}}$$
and
\begin{eqnarray*}
&&\|f_{i+1}\|_{r_{i+1},\bar s+s_{i+1}}\leq\nonumber\\
&& \left(1-\frac{\|f_i\|_{r_i,\bar s+s_i}}{\a{\r_N}{\s_N}/2}\right)^{-1}\left(\frac{\|f_i\|_{r_i,\bar s+s_i}+\|g_i\|_*}{\a{\r_N}{\s_N}}+e^{-Ks/(4N)}\right)\,\|f_i\|_{r_i,\bar s+s_i}\leq\nonumber\\
&&\left(1-\frac{1}{2^7}\right)^{-1}\left(\frac{1}{2^8}+\frac{1}{2^4}+\frac{1}{8}\right)\,\|f_i\|_{r_i,\bar s+s_i}=\nonumber\\
&&\frac{49}{254}\,\|f_i\|_{r_i,\bar s+s_i}\leq\nonumber\\
&& \frac{\|f_i\|_{r_i,\bar s+s_i}}{4}
\end{eqnarray*}
where we have used $Ks>12\log{2}$ which implies $$e^{-Ks/(4N)}\leq \frac{1}{8}$$
and the telescopic expansion
$$g_{i}=g_1+\sum_{k=2}^{i}(g_k-g_{k-1})$$
with $g_1:=P_\L T_Kf$, $g_k-g_{k-1}=P_\L T_Kf_{k-1}$ bounded on $\cP_{r,\bar s+s}$, $\cP_{r_{k},\bar s+s_{k}}\supset \cP_{r_{i},\bar s+s_{i}}$, for $1\leq k\leq i$. This gives, according to the definition of $\|g_{i}\|_*$  relatively to this expansion,
\begin{eqnarray*}
\|g_{i}\|_*&=&\frac{1}{2}\left(\frac{\r_N}{r-r_i+\r_N}+\frac{\s_N}{\bar s+s-(\bar s+s_i)+\s_N}\right)\|g_1\|_{r,\bar s+s}\nonumber\\
&+&\frac{1}{2}\,\sum_{k=2}^{i}\left(\frac{\r_N}{r_k-r_i+\r_N}+\frac{\s_N}{(\bar s+s_k)-(\bar s+s_i)+\s_N}\right)\|g_k-g_{k-1}\|_{{r_{k},\bar s+s_{k}}}\nonumber\\
&=&\frac{1}{2}\left(\frac{\r_N}{r-r_i+\r_N}+\frac{\s_N}{s-s_i+\s_N}\right)\|g_1\|_{r,\bar s+s}\nonumber\\
&+&\frac{1}{2}\,\sum_{k=2}^{i}\left(\frac{\r_N}{r_k-r_i+\r_N}+\frac{\s_N}{s_k-s_i+\s_N}\right)\|g_k-g_{k-1}\|_{{r_{k},\bar s+s_{k}}}\nonumber\\
&\leq&\frac{1}{4}\left(\frac{\r_N}{\r_0}+\frac{\s_N}{\s_0}\right)\|g_1\|_{r,\bar s+s}+\frac{1}{2}\,\sum_{k=2}^i\|f_{k-1}\|_{{r_{k},\bar s+s_{k}}}\nonumber\\
&\leq&\frac{5}{8N}\|g_1\|_{r,\bar s+s}+\frac{2}{3}\,\|f\|_{{r,\bar s+s}}\nonumber\\
&\leq&\frac{\a\r_N\,\s_N}{2^5}+\frac{\a\r_N\,\s_N}{2^8}\nonumber\\
&\leq&\frac{\a\r_N\,\s_N}{2^4}
\end{eqnarray*}
(use the following inequalities: $$\|f_{k-1}\|_{{r_{k},\bar s+s_{k}}}\leq \|f_{k-1}\|_{{r_{k-1},\bar s+s_{k-1}}}\leq4^{-(k-1)}\|f\|_{r,\bar s+s}\ ,$$ 
$$\|g_1\|_{r,\bar s+s}\leq \|f\|_{r,\bar s+s}\leq \a\,r/(2^7K)\leq \a\,rs/(2^{10}N)= N\,\a \r_N\s_N/(2^5)\ ,$$ $$r_k-r_i\geq \r_N\ ,\quad s_k-s_i\geq \s_N\ ,$$ $$r-r_i+\r_N\geq r-r_i\geq r-r_1=r/4=2\r_0\ ,$$ $$s-s_i+\s_N\geq s-s_i\geq s-s_1=s/3=2\s_0)\ .$$
and (\ref{inequalities}) are proved for any $1\leq i\leq N$.
Let
$$H:=H_{N+1}=h+g+f_*\ ,\quad \left(g:=g_{N+1}\ ,\quad f_*:=f_{N+1}\right)\ .$$
Then, by construction,
$$\|f_*\|_{r/2,\bar s+s/6}=\|f_{N+1}\|_{r/2,\bar s+s/6}\leq e^{-Ks/6}\|f\|_{r,\bar s+s}$$
and, using $Ks\geq 8$, $\|f_{i}\|_{{r_{i},\bar s+s_{i}}}\leq 4^{-i}\|f\|_{{r,\bar s+s}}$ and (\ref{f1est}), 
\begin{eqnarray*}
\|g-P_\L T_Kf\|_{r/2,\bar s+s/6}&=&\|g_{N+1}-g_1\|_{r/2,\bar s+s/6}\nonumber\\
&\leq&\sum_{k=1}^N\|g_{k+1}-g_{k}\|_{r/2,\bar s+s/6}\nonumber\\
&\leq&\sum_{k=1}^N\|g_{k+1}-g_{k}\|_{r_k,\bar s+s_k}\nonumber\\
&\leq& \sum_{k=1}^N\|f_k\|_{r_k,\bar s+s_k}\nonumber\\
&\leq&\frac{4}{3}\|f_1\|_{r_1,\bar s+s_1}\nonumber\\
&\leq &2^8\frac{\|f\|_{r,\bar s+s}^2}{\a\,r\,s}\nonumber\\
&\leq& \frac{2^5K}{\a\,r}\|f\|_{r,\bar s+s}^2\ .
\end{eqnarray*}
Furthermore,  by the usual telescopic arguments,
\begin{eqnarray}\label{Pi I- I}
|\P_{\bar I}(\Psi_N-\textrm{id})|_1&\leq& \sum_{k=0}^N|\P_{\bar I}(\Phi_k-\textrm{id})|_1\nonumber\\
&=&|\P_{\bar I}(\Phi_0-\textrm{id})|_1+\sum_{k=1}^N|\P_{\bar I}(\Phi_k-\textrm{id})|_1\nonumber\\
&\leq& \frac{\|f\|_{r,\bar s+s}}{\bar\a\s_0}+\frac{1}{\bar\a\s_N}\sum_{k=1}^N\|f_k\|_{r_k,\bar s+s_k}\nonumber\\
&\leq&\frac{\|f\|_{r,\bar s+s}}{\bar\a\s_0}+\frac{2}{\bar\a\s_N}2^7\frac{\|f\|_{r,\bar s+s}^2}{\bar\a\,r\,s}\nonumber\\
&=&\frac{\|f\|_{r,\bar s+s}}{\bar\a\s_0}\left(1+2^8N\frac{\|f\|_{r,\bar s+s}}{\bar\a\,r\,s}\right)\nonumber\\
&\leq&\frac{\|f\|_{r,\bar s+s}}{\bar\a\s_0}\left(1+2^5K\frac{\|f\|_{r,\bar s+s}}{\bar\a\,r}\right)\nonumber\\
&\leq&\frac{2\,\|f\|_{r,\bar s+s}}{\bar\a\s_0}\nonumber\\
&\leq&\frac{2^4\,\|f\|_{r,\bar s+s}}{\bar\a\,s}\ . 
\end{eqnarray}
Similarly, one proves
\begin{eqnarray}\label{Pi phi- phi}
|\P_{ I}(\Psi_N-\textrm{id})|_1\leq \frac{2^4\,\|f\|_{r, s+s}}{\a\,s}\ ,\quad |\P_{\varphi}(\Psi_N-\textrm{id})|_{\infty}\leq\frac{2^4\,\|f\|_{r,\bar s+s}}{\a\,r}\ ,
\end{eqnarray}
which gives (\ref{identity bounds}). We prove now (\ref{derivatives0}). Writing 
$$W^{\bar n}_{\bar\a^{-1},\a^{-1}}(W_{r,s}D\Psi_0W^{-1}_{r,s}-\id_{2n})=W_{r,s}W_0^{-1}[W^{\bar n}_{\bar\a^{-1},\a^{-1}}(W_0D\Psi_0W_0^{-1}-\id_{2n})]W_0W_{r,s}^{-1}$$ and using 
$$\|W^{\bar n}_{\bar\a^{-1},\a^{-1}}(W_0D\Psi_0W_0^{-1}-\id_{2n})\|_{\cal P}\leq \frac{2\,}{\r_0\s_0}\,\|f\|_{r,\bar s+s}$$
$$\|W_{r,s}\,W_0^{-1}\|_{\cal P}=\frac{1}{6}\ ,\quad \|W_0W_{r,s}^{-1}\|_{\cal P}=8$$
we find 

\begin{eqnarray}\label{bound0}
\|W^{\bar n}_{\bar\a^{-1},\a^{-1}}(W_{r,s}D\Psi_0W_{r,s}^{-1}-\id_{2n})\|_{\cal P}&\leq& \|W_{r,s}W_0^{-1}\|_{\cal P}\|W_0W_{r,s}^{-1}\|_{\cal P}\nonumber\\
&\times&\|W^{\bar n}_{\bar\a^{-1},\a^{-1}}(W_0D\Psi_0W_0^{-1}-\id_{2n})\|_{\cal P}\nonumber\\
&=&\frac{4}{3}\,\|W^{\bar n}_{\bar\a^{-1},\a^{-1}}(W_0\Psi_0W_0^{-1}-\id_{2n})\|_{\cal P}\nonumber\\
&\leq&\frac{8}{3}\frac{1}{\r_0\s_0}\,\|f\|_{r,\bar s+s}\leq1\ .
\end{eqnarray}
Similarly, writing  $$W^{\bar n}_{\bar\a^{-1},\a^{-1}}(W_{r,s}D\Psi_iW_{r,s}^{-1}-\id_{2n})=W_{r,s}W_N^{-1}[W^{\bar n}_{\bar\a^{-1},\a^{-1}}(W_ND\Psi_iW_N^{-1}-\id_{2n})]W_NW_{r,s}^{-1}$$ and using
$$\|W^{\bar n}_{\bar\a^{-1},\a^{-1}}(W_ND\Phi_i{W_N}^{-1}-\id_{2n})\|_{\cal P}\leq \frac{2\,}{\r_N\s_N}\,\|f_{i}\|_{r_{i},\bar s+s_{i}}\ ,$$
$$\|W_{r,s}W_N^{-1}\|_{\cal P}=\frac{1}{4N}\ ,\quad \|W_NW_{r,s}^{-1}\|_{\cal P}=8N$$
we arrive at
\begin{eqnarray}\label{boundi}
\|W^{\bar n}_{\bar\a^{-1},\a^{-1}}(W_{r,s}D\Phi_i{W}_{r,s}^{-1}-\id_{2n})\|_{\cal P}\leq \frac{4\,}{\r_N\s_N}\,\|f_{i}\|_{r_{i},\bar s+s_{i}}\leq\a\ ,\quad i=1,\cdots,N\ .
\end{eqnarray}
Taking into account (\ref{bound0}), (\ref{boundi}) and Lemma \ref{derivatives lemma} (where $W^{\bar n}_{\bar\a^{-1},\a^{-1}}$ plays the role of the invertible matrix $D$, and $4\|f_i\|/\r_N\s_N$ the one of $\e_i$, for $i\neq 0$), we find, for 
$$D\Psi_N=D\Phi_0D\Phi_1\cdots D\Phi_N$$
the bound
\begin{eqnarray*}
\|W^{\bar n}_{\bar\a^{-1},\a^{-1}}(W_{r,s}\,D\Psi_N{W}_{r,s}^{-1}-\id_{2n})\|_{\cal P}&\leq&\frac{8}{3}\frac{\|f\|_{r,\bar s+s}}{\r_0\s_0}+\frac{8\,}{\r_N\s_N}\,\sum_{i=1}^N2^{i-1}\|f_{i}\|_{r_{i},\bar s+s_{i}}
\end{eqnarray*}
with 
$$\sum_{i=1}^N2^{i-1}\|f_{i}\|_{r_{i},\bar s+s_{i}}\leq 2\|f_1\|_{r_{1},\bar s+s_{1}}\leq 2^8\frac{\|f\|_{r,\bar s+s}^2}{\a\,r\,s}$$
so that
\begin{eqnarray*}
\|W^{\bar n}_{\bar\a^{-1},\a^{-1}}(W_{r,s}\,D\Psi_N{W}_{r,s}^{-1}-\id_{2n})\|_{\cal P}&\leq&\frac{2^5K}{r}\,\|f\|_{r,\bar s+s}+\left(\frac{2^7K\|f\|}{\a r}\right)\left(\frac{K\|f\|}{r}\right)\nonumber\\
&\leq&\frac{2^6K}{r}\,\|f\|_{r,\bar s+s}\ .
\end{eqnarray*}
This completes the proof.

\pagina
\section{Birkhoff Normal Form}\label{app:Birkhoff Normal Form}
\setcounter{equation}{0}
\renewcommand{\theequation}{\ref{app:Birkhoff Normal Form}.\arabic{equation}}
In this section, we discuss quantitatively  the reduction to Birkhoff Normal Form, for Hamiltonians possessing elliptic equilibrium points. For further references, see also \cite{HofZeh94}.
\begin{proposition}{\bf (Birkhoff Normal Form)}\label{Birkhoff Theorem}
Let $0<\theta<1$; let ${\cal D}\subseteq \complex^n$ such that $$\O=(\O_1,\cdots,\O_m):{\cal D}\to \complex^m$$ is $(\a,K)$ non resonant on $\cD$, with $K\geq 2$,
and let
\begin{eqnarray*}
&&f(I,\varphi,p,q)=\sum_{i=1}^m\ \Omega_i(I)\,\frac{p_i^2+q_i^2}{2}+o_2(p,q;I)\ , 
\end{eqnarray*}
where $o_2$ is real--analytic in ${\cal D}\times B^{2m}_{{r}}(0)$, and verifies
$$\lim_{(p,q)\to 0}\frac{o_2(p,q;I)}{|(p,q)|^2}=0\ \quad \textrm{for all}\quad I\in{\cal D}\ .$$
Then, there exists $r_K>0$ and a symplectic, analytic transformation
\begin{eqnarray*}
\pi&:&\quad {\cal D}\times \complex^n/(2\pi\integer^n)\times B^{2m}_{r_K}(0)\to {\cal D}\times \complex^n/(2\pi\integer^n)\times B^{2m}_{r}(0)
\nonumber\\
&&(J,\vartheta,P,Q)\to (I,\varphi,p,q)=\pi(J,\vartheta,P,Q)
\end{eqnarray*}
with $I$, $\varphi-\vartheta$, $p$, $q$ independent from $\vartheta$, 
which puts $f$ into Birkhoff normal form up to order $K$. Furthermore, the following holds. 
\vskip.1in
\noindent
$i)\quad$
The transformation $\p$ may be obtained as a product
$$\pi={\cal B}_2\circ \cdots\circ{\cal B}_K\quad ({\cal B}_2=\rm id)\ ,$$ 
where
\begin{eqnarray*}
{\cal B}_k:\quad {\cal D}\times \complex^n/2\pi\integer^n\times B^{2m}_{r_{k}}&\to& {\cal D}\times \complex^n/2\pi\integer^n\times B^{2m}_{r_{k-1}}\ ,\nonumber\\
(\tilde J,\tilde{\vartheta},\tilde P,\tilde Q)&\to&(\tilde I,\tilde{\varphi},\tilde p,\tilde q)
\end{eqnarray*}
verifies $\tilde I=\tilde J$ and
\begin{eqnarray}\label{quantitative estimate3}
|\tilde q- \tilde Q|&\leq& \frac{M^k_{1}}{1-\theta}|( \tilde P, \tilde Q)|^{k-1}\ ,\quad | \tilde p- \tilde P|\leq \frac{M^k_{2}}{1-\theta}|( \tilde P, \tilde Q)|^{k-1}\nonumber\\
|\tilde{\varphi}-\tilde{\vartheta}|&\leq&2^{k}M_0^k\,\max\left\{1,\ \left(\frac{M_1^k}{1-\theta}\right)^k|(\tilde P,\tilde Q)|^{k(k-2))}\right\}\,|(\tilde P,\tilde Q)|^k\ ,\nonumber\\
\end{eqnarray}
and, for any $k=2,\cdots,K$, the product $\p_k:={\cal B}_2\circ \cdots \circ{\cal B}_k$ puts $f$ in Birkhoff normal form up to order $k$:
$$f\circ \pi_{{k}}(\tilde J,\tilde{\vartheta},\tilde P,\tilde Q)=\sum_{i=1}^m\ \Omega_i(\tilde J)\,\frac{\tilde P_i^2+\tilde Q_i^2}{2}+{\cal P}_{{k}}(\tilde J,\tilde P,\tilde Q)+o_{{k}}(\tilde P,\tilde Q;\tilde J)\ .$$
\vskip.1in
\noindent
$ii)\quad$ 
The constants $M^k_j$, $M_{jh}^k$ are inductively defined as follows.
Let, for $K\geq 3$ and $k=3,\cdots,K$, 
$$o_{k-1}(\tilde p,\tilde q;\tilde I)=\sum_{|\a|_1+|\b|_1=k}p^k_{\a,\b}(\tilde I)\left(\frac{\tilde p+i\tilde q}{\sqrt{2}}\right)^{\a}\left(\frac{\tilde p-i\tilde q}{\sqrt{2}}\right)^{\b}+\tilde o_k(\tilde p,\tilde q;\tilde I)\qquad (i=\sqrt{-1})\ ,$$
with $\tilde o_k(\tilde p,\tilde q;\tilde I)/|(\tilde p,\tilde q)|^k\to 0$, as $(\tilde p,\tilde q)\to 0$; let
$$s_k(\tilde J,\tilde P,\tilde q):=\sum_{\a\neq \b}{2}\,{i}\,\frac{p^k_{\a,\b}(\tilde J)}{\O(\tilde J)\cdot(\a-\b)}\left(\frac{\tilde P+i\tilde q}{\sqrt{2}}\right)^{\a}\left(\frac{\tilde P-i\tilde q}{\sqrt{2}}\right)^{\b}\ .$$
Then, 
$$M^k_j:=\sup_{{\cal D}\times B^{2m}_1(0)}|\partial_{{\cal P}(j)}\,s_k|\ ,\qquad M^k_{jh}:=\sup_{{\cal D}\times B^{2m}_1(0)}\|\partial_{{\cal P}(j){\cal P}(h)}\,s_k\|\ ,\quad j,\ h=0,\, 1,\, 2\ ,$$
where 
$${\cal P}(0):=\tilde J\ ,\quad {\cal P}(1):=\tilde P\ ,\quad {\cal P}(2):=\tilde q\ .$$
\vskip.1in
\noindent
$iii)\quad$ 
The polynomials ${\cal P}_k$ are inductively defined as follows. Starting with  
${\cal P}_2\equiv 0$, and given ${\cal P}_{k-1}(\tilde I,\tilde p,\tilde q)$, $o_{k-1}(\tilde p,\tilde q;\tilde I)$, then, for $K\geq 3$ and $k=3,\cdots,K$,
$${\cal P}_k(\tilde J,\tilde P,\tilde Q)={\cal P}_{k-1}(\tilde J,\tilde P,\tilde Q)+{\cal Q}_k(\tilde J,\tilde P,\tilde Q)\ ,$$ 
where
\begin{eqnarray}\label{Qp}{\cal Q}_k(\tilde J,\tilde P,\tilde Q)=\left\{
\begin{array}{lrr}
0\quad \textrm{for odd}\quad k\\
\sum_{|\a|_1=k/2}p^k_{\a,\a}(\tilde J)\left(\frac{\tilde P_1^2+\tilde Q_1^2}{2}\right)^{\a_1}\cdots\left(\frac{\tilde P_m^2+\tilde Q_m^2}{2}\right)^{\a_m}\quad \textrm{for even}\quad k .\\
\end{array}\right.
\end{eqnarray}
\vskip.1in
\noindent
$iv)\quad$ The radii $r_k$ are inductively defined as follows. Starting with $r_2=r$ and given $r_{k-1}$, then, for $K\geq 3$ and $k=3,\cdots,K$, 
$$r_k=\frac{1-\theta}{\sqrt{2}}\,\r_k\ ,\quad \r_k:=\min\left\{\left(\frac{\theta}{M^k_{12}}\right)^{1/(k-2)},\ \left(\frac{\theta}{\sqrt{2}M^k_{1}}\right)^{1/(k-2)},\ \left(\frac{\theta}{\sqrt{2}M^k_{2}}\right)^{1/(k-2)},\ r_{k-1}\right\}\ .$$
\end{proposition}
\begin{remark}\label{Birkhoff algorithm}\rm
Observe that, rather than projecting the remainders $o_{k+1}$ with order $k+1$ of $f_k:=f\circ\cB_2\circ\cdots\circ\cB_k$ over the spaces $(p+iq)/\sqrt{2}$, $(p-iq)/\sqrt{2}$, a simple algorithm for the computation of $\cP_{k+1}$ comes from the identity
$$\cP_{k+1}=\frac{1}{(2\p)^m}\int_{\torus^n}o_{k+1}(p,q;I)|_{p_i=\sqrt{P_i^2+Q_i^2}\cos\varphi_i,\ q_i=\sqrt{P_i^2+Q_i^2}\sin\varphi_i}d\varphi_i$$
\end{remark}
Our proof of the previous Proposition is based on the following
\begin{lemma}\label{Birkhoff Lemma}
Let $0<\theta<1$, $r>0$, ${\cal D}\subseteq \complex^n$ and, for $J\in{\cal D}$, $\varphi\in \complex^n/2\pi\integer^n$, $(P,q)\in \complex^{2m}$,
$$S(J,P,\varphi,q):=J\varphi+Pq+s(J,P,q)$$
where 
$$s(J,P,q):=\sum_{|\a|_1+|\b|_1=k}\s_{\a,\b}(J)P^{\a}q^{\b}$$
is a polynomial with degree $k\geq 3$ in $(P,q)$, with analytic coefficients $J\to \s_{\a,\b}(J)$. Then, $S$ is the generating function of a (symplectic,) analytic transformation
$${\cal B}:\quad {\cal D}\times \complex^n/2\pi\integer^n\times B^{2m}_{r'}(0)\to {\cal D}\times \complex^n/2\pi\integer^n\times B^{2m}_{r}(0)\ ,$$
$$(J,\vartheta,P,Q)\to (I,\varphi,p,q)$$
with
$$r':=\frac{1-\theta}{\sqrt{2}}\,\r\ ,\quad \r:=\min\left\{\left(\frac{\theta}{M_{12}}\right)^{1/(k-2)},\ \left(\frac{\theta}{\sqrt{2}M_{1}}\right)^{1/(k-2)},\ \left(\frac{\theta}{\sqrt{2}M_{2}}\right)^{1/(k-2)},\ r\right\}\ .$$
such that $I=J$ and
\begin{eqnarray}\label{quantitative estimate2}
| q- Q|&\leq& \frac{M_1}{1-\theta}|( P, Q)|^{k-1}\ ,\quad | p- P|\leq \frac{M_2}{1-\theta}|( P, Q)|^{k-1}\nonumber\\
|{\varphi}-{\vartheta}|&\leq&2^{k}M_0\,\max\left\{1,\ \left(\frac{M_1}{1-\theta}\right)^k|(P,Q)|^{k(k-2))}\right\}\,|(P,Q)|^k ,\nonumber\\
\end{eqnarray}
where 
$$M_j:=\sup_{{\cal D}\times B^{2m}_1(0)}|\partial_{{\cal P}(j)}\,s|\ ,\quad M_{jh}:=\sup_{{\cal D}\times B^{2m}_1(0)}\|\partial^2_{{\cal P}(j){\cal P}(h)}\,s\|\ ,$$
if $\cP$ sends the set $\{0, 1, 2\}$ to the set $\{J, P, q\}$ as
$${\cal P}(0):=J\ ,\quad {\cal P}(1):=P\ ,\quad {\cal P}(2):=q$$
\end{lemma}
{\bf Proof.}
Observe preliminarly that, as $s(J,P,q)$, is a homogeneous polynomial in $(P,q)$ with degree $k$, then, $|\partial_{{\cal P}(j)}(J,P,q)|$, is a homogeneous function of $(P,q)$ with degree $k$ or $k-1$ for $j=0$, $j\neq 0$, respectively, and
\begin{eqnarray}\label{partial of s}|\partial_{{\cal P}(j)}\,s(J,P,q)|\leq 
\left\{\begin{array}{lrr}
M_0\,|(P,q)|^{k}\quad \textrm{if}\quad j=0\\
M_j\,|(P,q)|^{k-1}\quad \textrm{if}\quad j\neq 0\\
\end{array}
\right.
\end{eqnarray}
by the definitions of $M_j$. Similarly, $\|\partial^2_{{\cal P}(j){\cal P}(h)}\,s(I,P,Q)\|$ is a homogeneous function with degree $k$ (for $j=h=0$), or $k-1$ (for $j=0\neq h$), or $k-2$ (for $j,\ h\neq 0$), and
\begin{eqnarray}\label{derivatives of order two}
\|\partial^2_{{\cal P}(j){\cal P}(h)}\,s(I,P,Q)\|\leq 
\left\{\begin{array}{lrr}
M_{00}\,|(P,q)|^{k}\quad \textrm{if}\quad j=h=0\\
M_{0h}\,|(P,q)|^{k-1}\quad \textrm{if}\quad j=0\ ,\ h\neq 0\\
M_{jh}\,|(P,q)|^{k-2}\quad \textrm{if}\quad j,\ h\neq 0\ .\\
\end{array}
\right.
\end{eqnarray}
We construct ${\cal B}$ by its generating equations, which are are
\begin{eqnarray}\label{tilde y I varphi}
\left\{
\begin{array}{lrr}
 I= J\\
{\varphi}={\vartheta}-\partial_{ J}s( J, P, q)\\
 p= P+\partial_{ q}s( J, P, q)\\
\end{array}
\right.
\end{eqnarray}
where $ q$ is obtained by solving, with respect to $ q$, the implicit equation
\begin{eqnarray}\label{tilde x}
 q+\partial_{ P}s( J, P, q)= Q\ .
\end{eqnarray}
Let
$$\r:=\min\left\{\left(\frac{\theta}{M_{12}}\right)^{1/(k-2)},\ \left(\frac{\theta}{\sqrt{2}M_{1}}\right)^{1/(k-2)},\ \left(\frac{\theta}{\sqrt{2}M_{2}}\right)^{1/(k-2)},\ r\right\}\ .$$
By (\ref{derivatives of order two}) it follows, for $J\in{\cal D}$, $(P,q)\in B^m_{\r/\sqrt{2}}(0)\times B^m_{\r/\sqrt{2}}(0)(\subseteq B^{2m}_{\r}(0))$, 
$$\|\partial^2_{Pq}s(J,P,q)\|\leq M_{12}\,\r^{k-2}\leq \theta<1\ ,$$
which is enough to assert that the function $ q\to  q+\partial_{ P}s( J, P, q)$ is injective on $B^m_{\r/\sqrt{2}}(0)$, for any $(J,P)\in{\cal D}\times B^m_{\r/\sqrt{2}}(0)$. Now, using
$$|\partial_{ P}s|\leq M_1\r^{k-1}\leq \frac{\theta \r}{\sqrt{2}}\quad \textrm{for}\quad (J,P,q)\in{\cal D}\times B^m_{\r/\sqrt{2}}(0)\times B^m_{\r/\sqrt{2}}(0)$$
we also find that, for any $(J,P)\in{\cal D}\times B^m_{\r/\sqrt{2}}(0)$, the map $ q\to  q+\partial_{ P}s( J, P, q)$ is onto on $B^m_{r'}(0)$, with $r'=(1-\theta)\r/\sqrt{2}$. Let $\tilde q\in B^{m}_{\r/\sqrt{2}}(0)$ the unique solution of (\ref{tilde x}), for $J\in{\cal D}$, $(P,Q)\in B^{2m}_{r'}(0)(\subseteq B^{m}_{r'}(0)\times B^{m}_{r'}(0))$, and let $p= P+\partial_{ q}s( J, P, q)$.
Using
$$|\partial_{ q}s|\leq M_2\r^{k-1}\leq \frac{\theta \r}{\sqrt{2}}$$
we find $p\in B^m_{\r/{\sqrt{2}}}(0)$, namely, $(p,q)\in B^{2m}_{\r}(0)\subseteq B^{2m}_{r}(0)$. Taking also $\varphi={\vartheta}-\partial_{ J}s( J, P, q)$, we have constructed
$${\cal B}:\quad {\cal D}\times \complex^n/2\pi\integer^n\times B^{2m}_{r'}(0)\to {\cal D}\times \complex^n/2\pi\integer^n\times B^{2m}_{r}(0)\ .$$
In order to prove (\ref{quantitative estimate2}), using (\ref{partial of s}) and the estimate
\begin{eqnarray}\label{diff of deriv}
|\partial_{ P}s( J, P, q)-\partial_{ P}s( J, P, Q)|&=&\left|\int_{0}^1\sum_j\partial^2_{ Pq_j}s( J, P, Q+t( q- Q))( q_j- Q_j)dt\right|\nonumber\\
&\leq&\int_{0}^1\left|\sum_j\partial^2_{ Pq_j}s( J, P, Q+t( q- Q))( q_j- Q_j)\right|dt\nonumber\\
&\leq&\sup_{{\cal D}\times r(( P, Q),( P, q))}\|\partial^2_{Pq}s\|| q- Q|\ .\nonumber\\
\end{eqnarray}
where $r(P,Q)$ denotes the straight line from $P$ to $Q$,we get
\begin{eqnarray}\label{difference of x}
| q- Q|&=&|\partial_{ P}s( J, P, q)|\nonumber\\
&\leq&|\partial_{ P}s( J, P, Q)|+|\partial_{ P}s( J, P, q)-\partial_{ P}s( J, P, Q)|\nonumber\\
&\leq&M_1|( P, Q)|^{k-1}+\sup_{{\cal D}\times r(( P, Q),( P, q))}\|\partial^2_{Pq}s\|| q- Q|\ .\nonumber\\
\end{eqnarray}
As $( P, Q)$, $( P, q)\in B^{2m}_{\r}(0)$, then, $r(( P, Q),( P, q))\subseteq B^{2m}_{\r}(0)$, hence,
\begin{eqnarray}\label{diff of deriv1}
\sup_{{\cal D}\times r(( P, Q),( P, q))}\|\partial^2_{Pq}s(J,P,q)\|\leq M_{12}\r^{k-2}\leq \theta\ ,
\end{eqnarray}
giving so, by (\ref{difference of x}),
$$| q- Q|\leq M_1|( P, Q)|^{k-1}+\theta| q- Q|$$
namely,
\begin{eqnarray}\label{difference of x1}
| q- Q|\leq \frac{M_1}{1-\theta}|( P, Q)|^{k-1}\ .
\end{eqnarray}
The proof of
$$| p- P|\leq \frac{M_2}{1-\theta}|( P, Q)|^{k-1}$$
is quite similar and is omitted. Using now (\ref{difference of x1}), we obtain
\begin{eqnarray}\label{difference of theta}
|{\varphi}-{\vartheta}|&\leq&M_0|( P, q)|^k\nonumber\\
&=&M_0|(P,Q)+(0,q-Q)|^k\nonumber\\
&\leq&M_0(|(P,Q)|+|(0,q-Q)|)^k\nonumber\\
&\leq&2^{k}M_0\,\max\left\{|(P,Q)|^k,\ |q-Q|^k\right\}\nonumber\\
&\leq&2^{k}M_0\,\max\left\{|(P,Q)|^k,\ \left(\frac{M_1}{1-\theta}\right)^k|(P,Q)|^{k(k-1)}\right\}\nonumber\\
&=&2^{k}M_0|(P,Q)|^k\,\max\left\{1,\ \left(\frac{M_1}{1-\theta}\right)^k|(P,Q)|^{k(k-2))}\right\}
\end{eqnarray}
and the proof of (\ref{quantitative estimate2}) is complete. 
\vskip.1in
\noindent
{\bf Proof of Proposition \ref{Birkhoff Theorem}.}\ We proceed by induction on $K$. For $K=2$, $f$ is yet in Birkhoff normal form up to order $2$, and the Proposition is proved with  $r_2=r$, $\pi_2={\cal B}_2=\textrm{id}$, ${\cal P}_2\equiv 0$. Assuming, now, that Proposition \ref{Birkhoff Theorem} holds when $K-1$ replaces $K$, we want to prove it for $K$. Assume, then, that ${\cal D}$ is $(\a,K)$ non resonant for $\O$. Obviously, ${\cal D}$ is $(\a,K-1)$ non resonant. By the inductive hypothesis, we find 
\begin{eqnarray*}
\pi_{K-1}={\cal B}_2\circ \cdots \circ{\cal B}_{K-1}:\quad {\cal D}\times \complex^n/2\pi\integer^n\times B^{2m}_{r_{K-1}}&\to& {\cal D}\times \complex^n/2\pi\integer^n\times B^{2m}_{r}\nonumber\\
(\tilde I,\tilde{\varphi},\tilde p,\tilde q)&\to&(I,\varphi,p,q)=\pi_{K-1}(\tilde I,\tilde{\varphi},\tilde p,\tilde q)
\end{eqnarray*}
with ${\cal B}_2=\textrm{id}$ and $${\cal B}_k:\quad {\cal D}\times \complex^n/2\pi\integer^n\times B^{2m}_{r_{k}}\to {\cal D}\times \complex^n/2\pi\integer^n\times B^{2m}_{r_{k-1}}\ ,\quad k=3,\cdots,K-1\ ,$$ 
which puts $f$ into Birkhoff normal form up to order $K-1$:
\begin{eqnarray}\label{r-1}
f_{K-1}(\tilde I,\tilde{\varphi},\tilde p,\tilde q)&:=&f\circ \pi_{K-1}(\tilde I,\tilde{\varphi},\tilde p,\tilde q)=\sum_{i=1}^m\O_i(\tilde I)\,\frac{\tilde p_i^2+\tilde q_i^2}{2}+{\cal P}_{K-1}(\tilde I,\tilde p,\tilde q)\nonumber\\
&+&o_{K-1}(\tilde p,\tilde q;\tilde I\nonumber\\
\end{eqnarray}
where items $i)$, $\cdots$, $iv)$ hold, for $k\leq K-1$.
We prove that, defining 
\begin{eqnarray}\label{rK}
\r_K&:=&\min\left\{\left(\frac{\theta}{M^K_{12}}\right)^{1/(K-2)},\ \left(\frac{\theta}{\sqrt{2}M^K_{1}}\right)^{1/(K-2)},\ \left(\frac{\theta}{\sqrt{2}M^K_{2}}\right)^{1/(K-2)},\ r_{K-1}\right\}\ ,\nonumber\\
r_K&:=&\frac{1-\theta}{\sqrt{2}}\,\r_K\ ,\nonumber\\
\end{eqnarray}
we find a symplectic, analytic transformation 
\begin{eqnarray*}
{\cal B}_K:\quad {\cal D}\times \complex^n/2\pi\integer^n\times B^{2m}_{r_{K}}&\to& {\cal D}\times \complex^n/2\pi\integer^n\times B^{2m}_{r_{K-1}}\nonumber\\
(J,{\vartheta},P,Q)&\to&(\tilde I,\tilde{\varphi},\tilde p,\tilde q)={\cal B}_K(J,{\vartheta},P,Q)
\end{eqnarray*}
with $\tilde{\varphi}-\vartheta$, $\tilde p$, $\tilde q$ independent from $\vartheta$ and
verifying $\tilde I=J$ and
\begin{eqnarray}\label{Kderiv}
|\tilde q- Q|&\leq& \frac{M^K_{1}}{1-\theta}|(P, Q)|^{K-1}\ ,\quad | \tilde p- P|\leq \frac{M^K_{2}}{1-\theta}|(P, Q)|^{K-1}\nonumber\\
|\tilde{\varphi}-{\vartheta}|&\leq&2^{K}M^K_0|(P,Q)|^K\,\max\left\{1,\ \left(\frac{M^K_1}{1-\theta}\right)^K|(P,Q)|^{K(K-2))}\right\}\ ,\nonumber\\
\end{eqnarray}
which puts $f_{K-1}$ in Birkhoff normal form up to order $K$.
We construct ${\cal B}_K$ by means of a generating function $S_K( J, P,\tilde{\varphi},\tilde q)$ of the form
\begin{eqnarray}\label{generating function}
S_K( J, P,\tilde{\varphi},\tilde q)= J\cdot\tilde{\varphi}+ P\cdot \tilde q+s_K( J, P,\tilde q)
\end{eqnarray}
where $s_K( J, P,\tilde q)$ is a homogeneous polynomial in $( P,\tilde q)$ with degree $K$, which we write as:
\begin{eqnarray}\label{sab}
s_K( J, P,\tilde q)=\sum_{|\a|_1+|\b|_1=K}s^K_{\a,\b}( J)\left(\frac{ P+i\tilde q}{\sqrt{2}}\right)^{\a}\left(\frac{ P-i\tilde q}{\sqrt{2}}\right)^{\b}\ .
\end{eqnarray}
Splitting, in (\ref{r-1}), $o_{K-1}$ as
$$o_{K-1}(\tilde p,\tilde q;\tilde I)=\sum_{|a|_1+|\b|_1=K}p^K_{\a,\b}(\tilde I)\left(\frac{\tilde p+i\tilde q}{\sqrt{2}}\right)^{\a}\left(\frac{\tilde p-i\tilde q}{\sqrt{2}}\right)^{\b}+\tilde o_K(\tilde p,\tilde q;\tilde I)$$
where $\tilde o_K(\tilde p,\tilde q;\tilde I/|(\tilde p,\tilde q)|^{K}\to 0$ as $(\tilde p,\tilde q)\to 0$, and replacing the generating equations of ${\cal B}_K$
\begin{eqnarray*}
\left\{
\begin{array}{lrr}
\tilde I=  J\\
\tilde {\varphi}= {\vartheta}-\partial_{  J}s(  J,  P, \tilde q)\\
 Q=\tilde q+\partial_{  P}s(  J, P,\tilde  q)\\
\tilde p= P+\partial_{\tilde q}s( J, P,\tilde q)\\
\end{array}
\right.
\end{eqnarray*}
into the definition (\ref{r-1}) of $f_{K-1}$, we find that $f_{K-1}$ changes to 
\begin{eqnarray}\label{cal Hr}
f_K( J,{\vartheta}, P, Q)&=&\sum_{i=1}^m\,\frac{ P_i^2+ Q_i^2}{2}\,\O_i( J)+{\cal P}_{K-1}( J, P, Q)\nonumber\\
&+&\sum_{|\a|_1+|\b|_1=K}\left[\frac{i}{2}\,\O( J)\cdot(\a-\b)\,s^K_{\a,\b}( J)+p^K_{\a,\b}( J)\right]\nonumber\\
&\times&\,\left(\frac{ P+i Q}{\sqrt{2}}\right)^{\a}\left(\frac{ P-i Q}{\sqrt{2}}\right)^{\b}\nonumber\\
&+&o_K( P, Q;J) ,\nonumber\\
\end{eqnarray}
and this leads us to choice, in (\ref{sab}),
\begin{eqnarray}\label{choice} 
s^K_{\a,\b}( J)=\left\{
\begin{array}{lrr} 
0\quad \textrm{for}\quad \a=\b\\
\\
{2}\,{i}\,{p^K_{\a,\b}( J)}/({\O( J)\cdot(\a-\b)})\quad \textrm{for}\quad \a\neq \b\ .\\
\end{array}
\right.
\end{eqnarray}
The definition (\ref{choice}) is well put because $|\a-\b|_1\leq K$ (observe $|a_l-\b_l|\leq \max\{\a_l,\b_l\}$) and ${\cal D}$ is $(\a,K)$ non resonant for $\O$. The choice (\ref{choice}) allows us to kill, in the summand in (\ref{cal Hr}), all the terms with $\a\neq \b$, and $f_K$ is in Birkhoff normal form up to order $K$. In particular, when $K$ is odd, no term survives, and ${\cal P}_{K}\equiv {\cal P}_{K-1}$. For even values of $K$, by (\ref{cal Hr}), we find 
$$f_K( J,{\vartheta}, P, Q)=\sum_{i=1}^m\,\frac{ P_i^2+ Q_i^2}{2}\,\O_i( J)+{\cal P}_{K}( J, P, Q)+o_K(P, Q;J) ,\nonumber\\$$
where
$${\cal P}_K( J, P, Q)={\cal P}_{K-1}( J, P, Q)+{\cal Q}_K( J, P, Q)$$
with
$${\cal Q}_K( J, P, Q)=\sum_{|\a|_1=K/2}p^K_{\a,\a}( J)\left(\frac{ P_1^2+ Q_1^2}{2}\right)^{\a_1}\cdots\left(\frac{ P_m^2+ Q_m^2}{2}\right)^{\a_m}\ .$$
On the other hand, by Lemma \ref{Birkhoff Lemma}, the function 
$$S_K( J, P,\tilde{\varphi},\tilde q)= J\cdot\tilde{\varphi}+ P\cdot \tilde q+s_K( J, P,\tilde q)$$
with
$$s_K( J, P,\tilde q)=\sum_{\a\neq \b}{2}\,{i}\,\frac{p^K_{\a,\b}( J)}{\O( J)\cdot(\a-\b)}\left(\frac{ P+i\tilde q}{\sqrt{2}}\right)^{\a}\left(\frac{ P-i\tilde q}{\sqrt{2}}\right)^{\b}$$
generates an analytic (symplectic) transformation
\begin{eqnarray*}{\cal B}_K:\quad {\cal D}\times \complex^n/2\pi\integer^n\times B^{2m}_{r_K}(0)&\to& {\cal D}\times \complex^n/2\pi\integer^n\times B^{2m}_{r_{K-1}}(0)\nonumber\\
(J,\vartheta,P,Q)&\to& (\tilde I,\tilde{\varphi},\tilde p,\tilde q)={\cal B}_K(J,\vartheta,P,Q)\ ,
\end{eqnarray*}
with $r_K$ as in (\ref{rK}), with $\tilde I=J$ and $\tilde q$, $\tilde p$, $\tilde{\varphi}-\vartheta$ independent from $\vartheta$, such that (\ref{Kderiv}) holds.
This completes the proof.
\newpage
\section{Proof of Lemma \ref{diagonalization N bodies}}
\setcounter{equation}{0}
\label{Proof of Lemma on diagonalization}
\renewcommand{\theequation}{\ref{Proof of Lemma on diagonalization}.\arabic{equation}}
In this appendix, we prove the Lemma \ref{diagonalization N bodies}. For shortness, we will refer to the property (\ref{aijneq0})$\div$(\ref{tilde F}) for a given matrix $\cA$ with order $n$ as $(*)$--property.

\vskip.1in
\noi 
{\bf Proof.}\ 
We proceed by induction on $n$.
The assertion is trivially true for $j=2$, by direct computation
\beqano\l_{21}&=&\frac{1}{2}\left(a_{11}+a_{22}+\sqrt{(a_{11}-a_{22})^2+4a_{12}^2}\right)\nonumber\\
&=&a_{11}+O\left(\frac{a_{12}^2}{a_{11}-a_{22}}\right)\nonumber\\
&=&a_{11}+O\left(\frac{\d^{2n_{12}}}{1-\bar a_{22}\d^{n_{22}}}\right)\nonumber\\
&=&a_{11}+O\left(\d^{2n_{12}}\right)\nonumber\\
\eeqano
and, similarly,
$$\l_{22}=\frac{1}{2}\left(a_{11}+a_{22}-\sqrt{(a_{11}-a_{22})^2+4a_{12}^2}\right)=a_{22}+O\left(\d^{2n_{12}}\right)\ .$$
Assume, now, that the Lemma holds for $n-1$. Let $\cA$ a matrix with order $n$ with the $(*)$--property and let
$\cP(\l)$  its characteristic polynomial. We are interested to solve equation
$$\cP(\l,\d)=0$$
closely to any diagonal element $a_{jj}$ of $\cA$. 
We  use an Implicit Function Theorem argument.
We expand the determinant of $\cA-\l\id_n$ along the $j^{th}$ row, so to split $\cP(\l,\d)$ as
\beqano
\cP(\l)=f(\l,\d)+g(\l,\d)
\eeqano
with
$$f(\l,\d):=(a_{jj}(\d)-\l)\det[M_{jj}(\l,\d)]\ ,\quad g(\l,\d):=\sum_{k\neq j}(-1)^{k-j}a_{j,k}(\d)\det[M_{j,k}(\l,\d)]\ ,$$where  $M_{j,k}$ is the minor with order $n-1$ of $\cA-\l\id_n$ 
with place $(j,k)$.
In particular, if  $\tilde\l_1$, $\cdots$, $\tilde\l_{n-1}$ are the eigenvalues of $M_{jj}(\l,\d)$, then, $\det[M_{jj}(\l,\d)]$ is given by
$$\det[M_{jj}(\l,\d)]=\prod_{k=1}^{n-1}(\tilde\l_k(\d)-\l)\ .$$
But
$M_{jj}(\l,\d)$ has the $(*)$--property, so, by the inductive hypothesis, its eigenvalues $\tilde\l_1$, $\cdots$, $\tilde\l_{n-1}$, verify
$$|\tilde\l_k-a_{kk}(\d)|\leq C \d^{\tilde m_k}\ ,\quad k=1,\cdots,j-1\ ,\quad |\tilde\l_k-a_{k+1,k+1}(\d)|\leq C \d^{\tilde m_k}\ ,\quad k=j,\cdots,n-1\ ,$$ for suitable $\tilde m_k$, $1\leq k\leq n-1$. Let $c_n>0$  so small that
$$\frac{c_n(n-2)(1+c_n)^{n-2}}{(1-c_n)^{n-1}}<1\ ,$$
$\d$ so small that
$$\min_k |\tilde\l_k(\d)-a_{jj}(\d)|>0\ .$$
The function $\l\to f(\l,\d)$ vanishes for $\l=a_{jj}(\d)$ and, for any $\l$ in the complex ball
$$|\l-a_{jj}(\d)|\leq R(\d):=c_n\min_k |\tilde\l_k(\d)-a_{jj}(\d)|$$
it results
\beqano
|\partial_\l f(\l,\d)|&=&\left|\Big[\prod_{k=1}^{n-1}(\tilde\l_k-\l)+(a_{jj}-\l)\sum_{1\leq m\leq n-1}\prod_{k\neq m}(\tilde\l_k-\l)\Big]\right|\nonumber\\
&\geq&(1-c_n)^{n-1}\prod_{k=1}^{n-1}|\tilde\l_k(\d)-a_{jj}(\d)|\nonumber\\
&-&c_n\min_k |\tilde\l_k(\d)-a_{jj}(\d)|\sum_{1\leq m\leq n-1}\prod_{k\neq m}|\tilde\l_k(\d)-a_{jj}(\d)|(1+c)\nonumber\\
&\geq&(1-c_n)^{n-1}\left[1-\frac{c(n-2)(1+c)^{n-2}}{(1-c_n)^{n-1}}\right]\prod_{k=1}^{n-1}|\tilde\l_k(\d)-a_{jj}(\d)|\nonumber\\
&\geq&(1-c_n)^{n-1}\left[1-\frac{c(n-2)(1+c)^{n-2}}{(1-c_n)^{n-1}}\right]\prod_{k=1}^{j-1}||\tilde\l_k(\d)||a_{jj}(\d)|^{n-j+1}
\eeqano
having used the inequality
$$(1-c_n)|\tilde\l_k(\d)-a_{jj}(\d)|\leq |\tilde\l_k-\l| \leq (1+c_n)|\tilde\l_k(\d)-a_{jj}(\d)|\ .
$$
for $1\leq k\leq n-1$.
On the other hand, any minor $M_{jk}$ appearing in the perturbation $g(\l,\d)$ contains the coloumn
$$
c_{jk}^j\left(
\begin{array}{c}
a_{1j}\\
\vdots\\
a_{j-1,j}\\
a_{j+1,j}\\
\vdots\\
a_{n,j}
\end{array}
\right)
$$
which is of order $\d^{\min\{n_{j-1,j},\ n_{j,j+1}\}}$. 
The remining $n-2$ columns of $M_{jk}$, $c_{jk}^l$ have only at one place place, a coordinate of the kind $a_{mm}-\l$, with $m\neq k$, $j$, and the other coordinates are $a_{pm}$, with $p\neq k,j$ $|a_{pm}|\leq |a_{mm}|\leq \max\{|\tilde\l_m|,|a_{jj}|\}$
$$|a_{mm}-\l|\leq|a_{mm}-a_{jj}|+|a_{jj}-\l|\leq |a_{mm}-a_{jj}|+|\tilde\l_m-a_{jj}|\leq 2|\tilde\l_m-a_{jj}|\leq 4\max\{|\tilde\l_m|,|a_{jj}|\}$$
so that 
$$|\det[M_{jk}]|\leq 4^n(n-1)!|a_{jj}|^{n-j}\prod_{m=1}^{j-1}|\tilde\l_m|$$
we find that 
$$|g(\l,\d)|=\left|\sum_{k\neq j}(-1)^{k-j}a_{j,k}(\d)\det[M_{j,k}(\l,\d)]\right|\leq 4^nn!\d^{2\min\{n_{j-1,j},\ n_{j,j+1}\}}|a_{jj}|^{n-j}\prod_{m=1}^{j-1}|\tilde\l_m|^{n-j}$$
By Cauchy estimate we find
$$|\partial_\l g(\l,\d)|\leq 4^nn!\d^{2\min\{n_{j-1,j},\ n_{j,j+1}\}}|a_{jj}|^{n-j-1}\prod_{m=1}^{j-1}|\tilde\l_m|^{n-j}$$
so that the Implicit Function Theorem may be applied provided $\d$ is sufficiently small:
$$\sup_{|\l-a_{jj}|\leq r(\d)}\frac{|\partial_\l g(\l,\d)|}{|\partial_\l f(\l,\d)|}\leq \sup_{|\l-a_{jj}|\leq r(\d)}\frac{|g(\l,\d)|}{r(\d)|\partial_\l f(\l,\d)|}\leq\left(\frac{\d}{\bar\d}\right)^{2\min\{n_{j-1,j},\ n_{j,j+1}\}-2n_{jj}}\leq \frac{1}{2}$$
with $$r(\d)=|a_{jj}|/2\leq R(\d)$$
We find then a unique solution $\l_j(\d)$ of 
$$\cP(\l)=f(\l,\d)+g(\l,\d)=0$$
verifying
$$|\l_j(\d)-a_{jj}(\d)|\leq \sup_{|\l-a_{jj}|\leq r(\d)}\frac{|g(\l,\d)|}{|\partial_\l f(\l,\d)|}\leq C\d^{2\min\{n_{j-1,j},\ n_{j,j+1}\}-n_{jj}}$$
Let us now study the $1$--dimensional eigenspace ${\cal V}_j$ associated to the aigenvalue $\l_j$. Let
$$v_j=\left(
\begin{array}{c}
v_{1j}\\
v_{2j}\\
\vdots\\
v_{nj}
\end{array}
\right)$$
the unitary eigenvector associated to $\l_j$, so that $V=(v_{ij})$ is the unitary matrix which diagonaluzes $\cA$:
$$V^T\cA V=\textrm{\rm diag}\,(\l_1,\cdots,\l_n)\ ,\quad V^TV=\id_n\ .$$
Then, the vector $\hat v_j$ with dimension $n-1$ which is obtained by $v_j$ dropping its $j^{th}$ component is the unique  solution of
$$(M_{jj}-\l_j\id_{n-1})\hat v_j=-\hat a_j v_{jj}$$
where $\hat a_j$ is the $j^{th}$ coloumn of $\cA$ deprivated of its $j^{th}$ component. But, as noticed before, $M_{jj}$ is almost diagonal, so, we write $M_{jj}-\l_j\id_{n-1}$ as
$$M_{jj}-\l_j\id_{n-1}=(D_j-\l_j\id_{n-1})\Big[\id_{n-1}-B_j\Big]$$
where $D_j$ is the principal diagonal  of $M_{jj}$ and $B_j$ is the off--diagonal
$$B_j=(\l_j\id_{n-1}-D_j)^{-1}(M_{jj}-D_j)\ .$$
The (non zero) elements $b^j_{hk}$ of $B_j$ go to $0$ with $\d$ because they satisfy the following asymptotics
$$
b^{j}_{h,k}=\kappa^{j}_{h,k}(\d)\,\d^{m^{j}_{h,k}}$$
with
$$0<m^{j}_{h,k}:=\left\{
\begin{array}{lrr}
{n_{h,k}-n_{h,h}}\quad \textrm{if}\quad h\neq k=1,\cdots,j-1\\
\\
{n_{h,k+1}-n_{h,h}}\quad \textrm{if}\quad h=1,\cdots,j-1\ ,\quad k=j,\cdots,N-1\\
\\
{n_{h+1,k}-n_{j,j}}\quad \textrm{if}\quad h=j,\cdots,N-1\ ,\quad k=1,\cdots,j-1\\
\\
{n_{h+1,k+1}-n_{j,j}}\quad \textrm{if}\quad h\neq k=j,\cdots,N-1\ .
\end{array}
\right.
$$
by assumption. The matrix $\id_{n-1}-B_j$ is thus non singular for small $\d$ and
$$(\id_{n-1}-B_j)^{-1}=\id_{n-1}+\d^{m_j}\bar B_j(\d)$$
where
$$m_j=\min_{h,k}\,\{m^{j}_{h,k}\}\ .$$
It follows, 
\beqano
\hat v_j&=&-(M_{jj}-\l_j\id_{n-1})^{-1}\hat a_j v_{jj}\nonumber\\
&=&-\Big[\id_{n-1}-B_j\Big]^{-1}(D_j-\l_j\id_{n-1})^{-1}\hat a_j v_{jj}\nonumber\\
&=&-(\id_{n-1}+\d^{m_j}\bar B_j(\d))\left(
\begin{array}{c}
\frac{a_{j1}}{a_{11}-\l_j}\\
\vdots\\
\frac{a_{j,j-1}}{a_{j-1,j-1}-\l_j}\\
\frac{a_{j,j+1}}{a_{j+1,j+1}-\l_j}\\
\vdots\\
\frac{a_{j,n}}{a_{n,n}-\l_j}
\end{array}
\right)\,v_{jj}\nonumber\\
&=&\left(
\begin{array}{c}
\check v_{1j}\,\d^{\n_{1j}}\\
\vdots\\
\check v_{j-1,j}\,\d^{\n_{j-1,j}}\\
\check v_{j+1,j}\,\d^{\n_{j+1,j}}\\
\vdots\\
\check v_{n,j}\,\d^{\n_{n,j}}
\end{array}
\right)v_{jj}
\eeqano
where
\beqano\n_{k,j}&=&\arr{
n_{jk}-n_{kk}\qquad \textrm{for}\quad 1\leq k\leq j-1\\
n_{jk}-n_{jj}\qquad \textrm{for}\quad j+1\leq k\leq n\\
}\nonumber\\
\check v_{kj}(0)&=&\arr{
-\frac{a_{jk}(0)}{a_{kk}(0)}\qquad \textrm{for}\quad 1\leq k\leq j-1\\
\frac{a_{jk}(0)}{a_{jj}(0)}\qquad \textrm{for}\quad j+1\leq k\leq n\\
}
\eeqano
By normalization,
$$v_{jj}=\frac{1}{\sqrt{1+\sum_{k\neq j}\check v_{k,j}^2\,\d^{2\n_{k,j}}}}=1+\check v_{jj}\,\d^{\n_{j,j}}$$
where $\check v_{jj}$, $\n_{j,j}$ are determined expanding 
$$z\to (1+z)^{-1/2}=1-\frac{z}{2}+O(z^2)$$ as in (\ref{check v and nu}).
The proof is complete.

\newpage
\section{The General Cauchy Inequality}
\setcounter{equation}{0}
\label{The General Cauchy Inequality}
We state a Cauchy inequality for the operatorial norm
\begin{eqnarray*}|d_v\,F|_{B,A}=\max_{u\neq 0}\frac{
|d_vF(u)|_B}{|u|_A}\end{eqnarray*}
of the first derivative $d_v\,F$, as a linear operator from $A$ to $B$, of a given analytic map $F:\ A\to B$, where $A$, $B$ are complex Banach spaces, with norms $|\cdot|_A$, $|\cdot|_B$. The present form is due to P\"oschel \cite{POSCH93}, to whom we refer for the proof.
\begin{lemma}
Let $F$ be an analytic map from the open ball of radius $r$ around $v$ in $A$ into $B$, such that $|F|_B\leq M$ into this ball. Then, the inequality
$$|d_v\,F|_{B,A}\leq\frac{M}{r}$$
holds. 
\end{lemma}
\section{Quantitative Implicit Function Theorem}
\setcounter{equation}{0}
\label{Quantitative Implicit Function Theorem}
\begin{theorem}
Let $F=f+g:\quad C^1(D^n_R(0),\ \complex^n)$, where:
\begin{itemize}
\item[(i)] $f$ is a diffeomorphism of $D^n_R(0)$ such that $f(0)=0$ and Jacobian matrix $\partial f$ non degenerate on $D^n_R(0)$;
\item[(ii)]$\sup_{D_R(0)}\|\partial g\|\sup_{D_R(0)}\|(\partial f)^{-1}\|\leq \frac{1}{2}$;
\item[(iii)]$\frac{\sup_{D_R(0)}|g|\sup_{D_R(0)}\|(\partial f)^{-1}\|}{r}\leq \frac{1}{2}$, where $0<r\leq R$;
\end{itemize}
Then, there exists a unique $z_0\in 
B^n_r(0)$ such that $F(z_0)=0$.
\end{theorem}
\section{The Laplace Coefficients}
\setcounter{equation}{0}
\label{Laplace coefficients}
\renewcommand{\theequation}{\ref{Laplace coefficients}.\arabic{equation}}
The {\it Laplace coefficients} $b_{s,k}(\a)$, is defined as the $k^{th}$ Fourier coefficients of the function $t\to (1+\a^2-2\a\cos{t})^{-s}$ 
\begin{eqnarray}\label{laplace coefficients}
b_{s,k}(\a)=\frac{1}{2\pi}\int_0^{2\pi}\,\frac{\cos{kt}}{(1+\a^2-2\a\cos{t})^s}\,dt\ ,\ \a\in \complex\ ,\ |\a|\neq 1\ , \ 0<s\in \real\ ,\ k\in \integer\nonumber\\ 
\end{eqnarray}
\begin{lemma}\label{asymp}
The Laplace coefficients are analytic of $\a$, for $|\a|<1$, and verify
\begin{itemize}
\item[(i)]$b_{s,k}(-\a)=(-1)^k\,b_{s,k}(\a)$;
\item[(ii)]$b_{s,-k}(\a)=b_{s,k}(\a)$;
\item[(iii)]$b_{s,k}(1/\a)=\a^{2s}\,b_{s,k}(\a)$;
\item[(iv)]$b_{s,k+2}(\a)=\frac{k+1}{k+2-s}\left(\a+\frac{1}{\a}\right)\,b_{s,k+1}(\a)-\frac{k+s}{k+2-s}\,b_{s,k}(\a)$;
\item[(v)]if $k\geq 0$, $b_{s,k}(\a)=\a^k\,\b_{s,k}(\a)$ where $\b_{s,k}(\a)$ is an even function of $\a$,  verifying
\begin{eqnarray}\label{asymptotic1}
\b_{s,k}(\a)&=&\frac{s(s+1)\cdots(s+k-1)}{k!}+s\,\frac{s(s+1)\cdots(s+k)}{(k+1)!}\ \a^2+O(\a^4)\ ,\nonumber\\
\end{eqnarray}
where
$$\frac{s(s+1)\cdots(s+k-1)}{k!}:=1\quad \textrm{if}\quad k=0\ .$$
\end{itemize}
\end{lemma}
\vskip.1in
\noindent
Notice that, by {\sl (iv)}, all the $b_{s,k}(\a)'s$ with $|k|\geq 2$ may be expressed as linear functions of $b_{s,0}(\a)$, $b_{s,1}(\a)$.
\vskip.1in
\noindent
{\bf Proof.}\ 
Items {\it (i)$\div$(iv)} are imediate consequences of (\ref{laplace coefficients}); in particular, {\it(iv)} is found by integrating twice by parts. In order to prove {\it (v)}, from (\ref{laplace coefficients}),  
we introduce the {\it hypergeometric series}
\begin{eqnarray*}
\frac{1}{(1-w)^s}=\sum_{l\geq 0}\,\frac{s(s+1)\cdots(s+l-1)}{l!}\,w^l
\end{eqnarray*}
with
\begin{eqnarray*}
\frac{s(s+1)\cdots(s+l-1)}{l!}\equiv1\quad \textrm{for}\quad l=0\ .
\end{eqnarray*}
The hypergeometric series is uniformly convergent in every closed disk inside the set $\{|w|<1\}$,
therefore, we may expand, for $\{|\a|\leq r<1\}$,
\begin{eqnarray*}
\frac{1}{(1+\a^2-2\a\,\cos{t})^s}&=&\frac{1}{(1-\a\,e^{it})^s(1-\a\,e^{-it})^s}\nonumber\\
&=&\sum_{l,j}\frac{s(s+1)\cdots(s+l-1)}{l!}\frac{s(s+1)\cdots(s+j-1)}{j!}\a^{i+j}e^{i(l-j)t}
\end{eqnarray*}
hence, multiplying by $\cos{kt}$ and then integrating over $[0,2\pi]$, we find
\begin{eqnarray}\label{Laplace series} b_{s,k}(\a)
&=&\a^k\,\sum_{j}\frac{s(s+1)\cdots(s+j+k-1)}{(j+k)!}\frac{s(s+1)\cdots(s+j-1)}{j!}\,\a^{2j}\ .\nonumber\\
\end{eqnarray}
It follows, in particular, that the $b_{s,k}(\a)'s$ are analytic for $|\a|<1$, and {\it (v)} is obtained by truncation of (\ref{Laplace series}).
\newpage
\noi
{\Large \bf Index of Notations}
\addcontentsline{toc}{section}{Index of Notations}
\vskip.2in
\noi
The references denote the  number of the page of the first occurrence in the text.
\vskip.1in
\noi
${\eufm K}$, ${\cal K}$:\  {\bf Kolmogorov sets}, \pageref{cal K}, \pageref{plane problem}.
\vskip.1in
\noi
{\bf \large KAM}
\vskip.1in
\noi
{\bf Sets}\\
$\natural$, $\integer$, $\rational$, $\real$, $\complex$: usual number sets.\\
$B_r^p(I):$ $p$--dimensional real ball centered at $I$, with radius $r$, \pageref{A two times scale KAM theorem}\\
$D_r^p(I):$ $p$--dimensional complex disk centered at $I$, with radius $r$, \pageref{A two times scale KAM theorem}\\
${\cal I}_r:$ complex $r$--neighborhood of $\cI\subset\real^{p}$, \pageref{A two times scale KAM theorem}\\
$\cD^{\bar n,\hat n}: $ generalized Diophantine set, \pageref{KAM and degeneracy}\\
$\cD^{\bar n,\hat n}_{\g,\hat\g;\t}: $ $(\g,\hat\g;\t)$--generalized Diophantine set, \pageref{two gammas}\\
$\torus^p:$ real standard $n$--dimensional torus, \pageref{A two times scale KAM theorem}\\
$\torus^p_s:$ complex $s$--neighborhood of $\torus^p$, \pageref{A two times scale KAM theorem}\\
$\torus_{\complex}^p:$ complex standard $n$--dimensional torus, \pageref{A two times scale KAM theorem}\\
If ${\cal I}=\bar{\cal I}\times \hat{\cal I}$, $h:{\cal I}\to \complex$ is analytic, and  $\o:=\partial h$, $\bar\o$ means $\partial_{\bar I}h$, $\hat\o$ means $\partial_{\hat I}h$, when $I=(\bar I,\hat I)$ is the generic element of ${\cal I}=\bar{\cal I}\times \hat{\cal I}$, with $\bar I\in\bar{\cal I}$, $\hat I\in\hat{\cal I}$, \pageref{two scales KAM};\\
If ${\cal I}$ is as before, $\o:\quad{\cal I}\to{\cal O}$ is onto and $\n\in{\cal O}$, $\bar{\o}^{-1}(\n)$ means the projection over the $\bar I$--coordinate of $\o^{-1}(\n)$ and $\hat{\o}^{-1}(\n)$ means the projection over the $\hat I$--coordinate, \pageref{two scales KAM};
\vskip.1in
\noi
{\bf Differential Operators}\\
$D$ differential operator with respect to $(I,\varphi)$, \pageref{derivatives};\\
$\partial$, $\partial^2$ differential operators with respect to $I$, \pageref{two scales KAM}\\
\vskip.1in
\noi
{\bf Matrices}\\
$A^{[p,q]}$, $A_*^{[p,q]}$: submatrices of a matrix $A$, \pageref{A two times scale KAM theorem}
\vskip.1in
\noi
{\bf  Norms}\\
{\sl a) for numbers}\\
$|k|$: $1$--norm of $k\in \integer^p$, \pageref{A two times scale KAM theorem}\\
$|(I,\varphi)|_{\cal P}$: ${\cal P}$--norm of ${\cal I}_\r\times\torus^n_{s}$, \pageref{average}\\
\\
{\sl b) for analytic functions}\\
$\|f\|_{r,s}:$ Sup--Fourier norm of a real--analytic function $f$ on ${\cal I}_r\times\torus^p_s$,  ${\cal I}\subset \real^p$  compact, \pageref{A two times scale KAM theorem}\\
\\
{\sl c) for Lipschitz functions}\\
${\cal L}(f)$, ${\cal L}_+(f)$, ${\cal L}_-(f)$, ${\cal L}_{\|\cdot\|}(f)$, $\|f\|^{\rm Lip}_{\r,\cI}$ Lipschitz norms, \pageref{two scales KAM}\\
\\
{\sl d) for vector and matrix functions}\\
If  $\o:{\cal I}\to \real^n$, $|\o|$ means its operator  norm when $\o$ is seen as linear operator from $({\cal I},||_1)$ to $(\complex,||)$ (corresponds to $|\o|_{\infty}:=\max|\o_i|$, with $\o_i$ $i^{th}$ coordinate of $\o$), \pageref{two scales KAM};\\
If $U:{\cal I}\to {\rm Matr}(m\times n)$, $\|U\|$ means its operator  norm when $U$ is seen as linear operator from $({\cal I},||_1)$ to $(\complex^m,||_{\infty})$(corresponds to $|U|_{\infty}:=\max|U_{ij}|$, with $U_{ij}$ the entries of $U$), \pageref{two scales KAM};
\vskip.1in
\noi
{\bf \large Celestial Mechanics}
\vskip.1in
\noi
${\cal R}_{\rm x}$, ${\cal R}_{\rm z}$: elementary rotations, \pageref{momentum position}, \pageref{inversion};\\
Elliptic elements, , \pageref{momentum position}, \pageref{two body problem};\\
$\p_w$, $\a_w(u,v)$, \pageref{assumptions}, 125;\\
$\cA$, ${\cA_{\rm G}}$, \pageref{non resonance}.
\vskip.1in
\noi
{\bf Domains and maps}\\
$\phi_{\rm DP}$: plane Delaunay--Poincar\'e map, \pageref{DP1};\\
$\cC_*$: \pageref{assumptions};\\
$\cD_*$, $\Phi_*$: Deprit action--angle map, \pageref{D*};\\
$\cD_{\rm r}$, $\Phi_{\rm r}$, $\phi_{\rm r}$: regularization of (${\cal D}_*$, $\Phi_*$, $\phi_*$), \pageref{regularization};\\
$\Phi_{\rm BD}$: full reduction map from  cartesian variables to $\cD_{\rm r}$, \pageref{regularization};\\
$\phi_{\rm BD}=\Phi_{\rm BD}^{-1}$, \pageref{eufm R i};\\
$\cD_{\rm pr}$: domain of regularized partially reduced Deprit variables, \pageref{part red2};\\
$\phi_{\rm BD,pr}$: full reduction map, from $\cD_{\rm pr}$ to cartesian variables, \pageref{partial reduction prop}.\\
\newpage

\addcontentsline{toc}{section}{Bibliography}
\bibliographystyle{plain}
 
\end{document}